%% file: teukolskyvanishinglambda2.tex
\documentclass[10pt,reqno,a4paper]{article}

\usepackage{geometry}
\usepackage{float}

\usepackage{array} 
\usepackage{paralist} 
\usepackage{verbatim} 
\usepackage{amsfonts,graphicx,amssymb,amsmath,mathtools}

\usepackage{amsthm}
\usepackage{ifthen}

\usepackage{tikz} 
\usepackage{tikz-3dplot}
\usepackage{pgfplots}
\usetikzlibrary{pgfplots.fillbetween,positioning,intersections,decorations.markings}
\usetikzlibrary{patterns} 
 
\pgfdeclaredecoration{complete sines}{initial}
{
  \state{initial}[
  width=+0pt,
  next state=sine,
  persistent precomputation={\pgfmathsetmacro\matchinglength{
      \pgfdecoratedinputsegmentlength / int(\pgfdecoratedinputsegmentlength/\pgfdecorationsegmentlength)}
    \setlength{\pgfdecorationsegmentlength}{\matchinglength pt}
  }] {}
  \state{sine}[width=\pgfdecorationsegmentlength]{
    \pgfpathsine{\pgfpoint{0.25\pgfdecorationsegmentlength}{0.5\pgfdecorationsegmentamplitude}}
    \pgfpathcosine{\pgfpoint{0.25\pgfdecorationsegmentlength}{-0.5\pgfdecorationsegmentamplitude}}
    \pgfpathsine{\pgfpoint{0.25\pgfdecorationsegmentlength}{-0.5\pgfdecorationsegmentamplitude}}
    \pgfpathcosine{\pgfpoint{0.25\pgfdecorationsegmentlength}{0.5\pgfdecorationsegmentamplitude}}
  }
  \state{final}{}
}

\makeatletter
\tikzset{
  hatch distance/.store in=\hatchdistance,
  hatch distance=5pt,
  hatch thickness/.store in=\hatchthickness,
  hatch thickness=5pt
}
\pgfdeclarepatternformonly[\hatchdistance,\hatchthickness]{north east hatch}
{\pgfqpoint{-1pt}{-1pt}}
{\pgfqpoint{\hatchdistance}{\hatchdistance}}
{\pgfpoint{\hatchdistance-1pt}{\hatchdistance-1pt}}%
{
  \pgfsetcolor{\tikz@pattern@color}
  \pgfsetlinewidth{\hatchthickness}
  \pgfpathmoveto{\pgfqpoint{0pt}{0pt}}
  \pgfpathlineto{\pgfqpoint{\hatchdistance}{\hatchdistance}}
  \pgfusepath{stroke}
}
\makeatother

\usepackage{slashed}
\usepackage{xr-hyper}
\usepackage[dvipsnames]{xcolor}
\colorlet{rouge}{red!50!black}
\colorlet{vert}{green!40!black}
\usepackage[urlcolor=blue!80!black, bookmarksopen=true, colorlinks, linkcolor=vert, citecolor=rouge]{hyperref} 

\usepackage{tensor} 
\usepackage{cancel}

\usepackage{appendix}
\usepackage{enumitem}

\usepackage{bbm} 
\usepackage{bm}
\usepackage{mathrsfs}
\usepackage{zref-clever}
\zcsetup{cap=true}

\zcRefTypeSetup{equation}{
    Name-sg=,
    name-sg=,
    Name-pl=,
    name-pl=,
    Name-sg-ab=,
    name-sg-ab=,
    Name-pl-ab=,
    name-pl-ab=
}

\makeatletter
\renewcommand*\env@matrix[1][\arraystretch]{%
  \edef\arraystretch{#1}%
  \hskip -\arraycolsep
  \let\@ifnextchar\new@ifnextchar
  \array{*\c@MaxMatrixCols c}}
\makeatother

\DeclarePairedDelimiter{\bangle}{\langle}{\rangle}
\DeclarePairedDelimiter{\abs}{\lvert}{\rvert}
\DeclarePairedDelimiter{\norm}{\lVert}{\rVert}
\DeclarePairedDelimiter{\squareBrace}{[}{]}
\DeclarePairedDelimiter{\curlyBrace}{\{}{\}}
\DeclarePairedDelimiter{\PoissonB}{\{}{\}}
\DeclarePairedDelimiter{\evalAt}{.}{\vert}

\newcommand{\addtheorem}[2]{%
  \AddToHook{env/#1/begin}{%
    \zcsetup{countertype={theorem=#1}}%
  }%
  \zcRefTypeSetup{#1}{Name-sg=\textcolor{vert}{#2}}%
  \newtheorem{#1}[theorem]{#2}%
}

\newtheorem{theorem}{Theorem}[section]
\addtheorem{proposition}{Proposition}
\addtheorem{lemma}{Lemma}
\addtheorem{remark}{Remark}
\addtheorem{definition}{Definition}
\addtheorem{corollary}{Corollary}
\addtheorem{conjecture}{Conjecture}

\DeclareFontFamily{U}{mathx}{}
\DeclareFontShape{U}{mathx}{m}{n}{<-> mathx10}{}
\DeclareSymbolFont{mathx}{U}{mathx}{m}{n}
\DeclareMathAccent{\widehat}{0}{mathx}{"70}
\DeclareMathAccent{\widecheck}{0}{mathx}{"71}

\newcommand{\dBar}{\underline{d}}
\newcommand{\NT}{\operatorname{nt}}
\newcommand{\tauAux}{k}
\newcommand{\timelikeInf}{i^+}

\newcommand{\LHS}{left-hand side}
\newcommand{\RHS}{right-hand side}
\newcommand{\Real}{\mathbb{R}}
\newcommand{\Complex}{\mathbb{C}}
\newcommand{\Natural}{\mathbb{N}}

\newcommand{\Laplace}{\bigtriangleup}
\newcommand{\ImagUnit}{\mathbbm{i}}
\newcommand{\bOne}{\mathbbm{1}}
\newcommand{\Identity}{\operatorname{Id}}
\newcommand{\supp}{\operatorname{supp}\,}

\newcommand{\Divergence}{\operatorname{div}}
\newcommand{\Curl}{\operatorname{curl}}
\newcommand{\Hom}{\operatorname{Hom}}
\newcommand{\End}{\operatorname{End}}

\newcommand{\Manifold}{\mathcal{M}}
\newcommand{\LieDerivative}{\mathcal{L}}

\newcommand{\Ric}{\mathbf{Ric}}
\newcommand{\Riem}{\mathbf{R}}
\newcommand{\Weyl}{\mathbf{W}}

\newcommand{\ChristoffelTypeTwo}[3][]{
  \ifthenelse{\equal{#1}{}}
  {\tensor{\Gamma}{^{#2}_{#3}}} 
  {\tensor{\Gamma(#1)}{^{#2}_{#3}}}
}
\newcommand{\ChristoffelTypeOne}[3][]{
  \ifthenelse{\equal{#1}{}}
  {\tensor{\Gamma}{_{{#2#3}}}}
  {\tensor{\Gamma(#1)}{_{{#2#3}}}}
}
\newcommand{\Trace}{\operatorname{tr}}

\newcommand{\Sphere}{\mathbb{S}}
\newcommand{\DeformationTensor}[2][]{
  \ifthenelse{\equal{#1}{}}
  {\tensor[^{({#2})}]{\pi}{}}
  {\tensor[^{({#2})}]{\pi}{#1}}
}
\newcommand{\DeformationTensorLinearized}[2][]{
  \ifthenelse{\equal{#1}{}}
  {\tensor[^{({#2})}]{\widecheck{\pi}}{}}
  {\tensor[^{({#2})}]{\widecheck{\pi}}{#1}}
}
\newcommand{\DeformationTensorErr}[2][]{
  \ifthenelse{\equal{#1}{}}
  {\tensor[^{({#2})}]{\tilde{\pi}}{}}
  {\tensor[^{({#2})}]{\tilde{\pi}}{#1}}
}

\newcommand{\CovariantDeriv}{\mathbf{D}}

\newcommand{\ABar}{\underline{A}}
\newcommand{\BBar}{\underline{B}}

\newcommand{\XBar}{\underline{X}}
\newcommand{\XHat}{\widehat{X}}
\newcommand{\XHatBar}{\underline{\XHat}}
\newcommand{\XiBar}{\underline{\Xi}}
\newcommand{\HBar}{\underline{H}}
\newcommand{\omegaBar}{\underline{\omega}}
\newcommand{\CCOneFormJ}{\mathfrak{J}}
\newcommand{\LeftDual}[1]{\tensor[^*]{{#1}}{}}
\newcommand{\RightDual}[1]{{{#1}^{*}}}
\newcommand{\Metric}{\mathbf{g}}
\newcommand{\chiBar}{\underline{\chi}}
\newcommand{\etaBar}{\underline{\eta}}
\newcommand{\xiBar}{\underline{\xi}}
\newcommand{\alphaBar}{\underline{\alpha}}
\newcommand{\betaBar}{\underline{\beta}}
\newcommand{\chiTF}{\widehat{\chi}}
\newcommand{\chiBarTF}{\widehat{\chiBar}}
\newcommand{\aTrace}[1]{\tensor[^{(a)}]{\Trace{#1}}{}}
\newcommand{\SymTracelessTensorProd}{\widehat{\otimes}}
\newcommand{\ComplexDeriv}{\mathcal{D}}

\newcommand{\volForm}{\in}
\newcommand{\volFormHor}{\in}
\newcommand{\horProj}[1]{\tensor[^{(h)}]{#1}{}}
\newcommand{\horMetric}{g}
\newcommand{\realHorkTensor}[1]{\mathfrak{s}_{#1}}
\newcommand{\HorkTensor}[1]{\mathbf{O}_{#1}}
\newcommand{\ConformalComplexDeriv}{\tensor[^{(c)}]{\ComplexDeriv}{}}
\newcommand{\BHorizontalCurv}{\mathbf{B}}
\newcommand{\LaplaceHor}{\bigtriangleup}
\newcommand{\HorRiem}{\dot{\Riem}}
\newcommand{\HorLieDeriv}{\slashed{\LieDerivative}}
\newcommand{\HodgeOp}[1]{\slashed{\mathcal{D}}_{{#1}}}
\newcommand{\HodgeOpDual}[1]{\slashed{\mathcal{D}}^{*}_{{#1}}}

\newcommand{\ConformalInvDeriv}{\tensor[^{(c)}]{\nabla}{}}

\newcommand{\ConformalInvDualDeriv}{\tensor[^{*(c)}]{\nabla}{}}

\newcommand{\TeukOp}{\mathcal{L}}

\newcommand{\HorVFSpace}{\mathbf{O}}
\newcommand{\WaveOpHork}[1]{\dot{\Box}_{#1}}

 \newcommand{\HorCovDeriv}{\dot{\CovariantDeriv}}

\newcommand{\ein}{e^{(\operatorname{in})}}
\newcommand{\eout}{e^{(\operatorname{out})}}
\newcommand{\eglo}{e^{(\operatorname{glo})}}
\newcommand{\glo}{\operatorname{glo}}
\newcommand{\lambdaglo}{\lambda_{\glo}}
\newcommand{\chiglo}{\chi_{\glo}}

\newcommand{\frakWeightedDeriv}{\mathfrak{d}}
\newcommand{\frakWeightedDerivAngular}{\slashed{\frakWeightedDeriv}}

\newcommand{\TransformScalarWaveOp}[1][]{
  \ifthenelse{\equal{#1}{}}
  {\widehat{\Box}^{(0)}}
  {\widehat{\Box}^{(0)}_{#1}}
}
\newcommand{\ScalarWaveOp}[1][]{
  \ifthenelse{\equal{#1}{}}
  {\Box}
  {\Box_{#1}}
}
\newcommand{\VectorWaveOp}[1][]{
  \ifthenelse{\equal{#1}{}}
  {\Box^{(1)}}
  {\Box^{(1)}_{#1}}
}
\newcommand{\TensorWaveOp}[1][]{
  \ifthenelse{\equal{#1}{}}
  {\Box^{(2)}}
  {\Box^{(2)}_{#1}}
}
\newcommand{\ScalarWaveConjOp}[1][]{
  \ifthenelse{\equal{#1}{}}
  {\overline{\Box}^{(0)}}
  {\overline{\Box}^{(0)}_{#1}}
}
\newcommand{\ScalarWaveLaplaceOp}[1][]{
  \ifthenelse{\equal{#1}{}}
  {\widehat{\Box}^{(0)}}
  {\widehat{\Box}^{(0)}_{#1}}
}
\newcommand{\ReducedWaveOp}[1][]{
  \ifthenelse{\equal{#1}{}}
  {\widetilde{\Box}^{(0)}}
  {\widetilde{\Box}^{(0)}_{#1}}
}

\newcommand{\Lagrangian}{\mathcal{L}}
\newcommand{\KillT}{\mathbf{T}}
\newcommand{\KillPhi}{\mathbf{\Phi}}

\newcommand{\RedShiftY}{\mathbf{Y}}

\newcommand{\JCurrent}[1]{J^{#1}}
\newcommand{\JCurrentPert}[1]{\widetilde{J}^{#1}}

\newcommand{\JCurrentMod}[1]{\mathring{J}^{#1}}
\newcommand{\KCurrent}[1]{K^{#1}}
\newcommand{\KCurrentAngular}[1]{\slashed{K}^{#1}}
\newcommand{\KCurrentR}{\mathcal{A}}
\newcommand{\KCurrentLTwo}{\mathcal{P}}
\newcommand{\KCurrentPert}[1]{\widetilde{K}^{#1}}

\newcommand{\KCurrentSym}[1]{\mathbf{k}^{#1}}
\newcommand{\KCurrentPertSym}[1]{\tilde{\mathbf{k}}^{#1}}

\newcommand{\HprVFGen}{R}

\newcommand{\EnergyFlux}{\mathbf{E}}
\newcommand{\EnergyFluxFar}{\dot{\EnergyFlux}}
\newcommand{\EnergyFluxWeighted}{\dot{\EnergyFlux}_{p,R}}
\newcommand{\EnergyFluxWeightedOpt}[2]{\dot{\EnergyFlux}_{{#1},{#2}}}
\newcommand{\EnergyFluxCombined}{\EnergyFlux_p}
\newcommand{\EnergyFluxCombinedOpt}[1]{\EnergyFlux_{#1}}
\newcommand{\EnergyFluxInitAux}[1]{\underline{\EnergyFlux}_{#1}}

\newcommand{\SpacelikeFlux}{\mathbf{F}}
\newcommand{\SpacelikeFluxFar}{\dot{\mathbf{F}}}
\newcommand{\SpacelikeFluxWeighted}{\dot{\mathbf{F}}_p}
\newcommand{\SpacelikeFluxCombined}{\mathbf{F}_p}

\newcommand{\ForcingTermNorm}{\mathcal{N}}
\newcommand{\ForcingTermWeightedNorm}[2]{\mathcal{N}_{{#1},{#2}}}
\newcommand{\ForcingTermCombinedNorm}[1]{\mathcal{N}_{#1}}
\newcommand{\ForcingTermNormReg}[1]{\tensor[^{\operatorname{#1}}]{\ForcingTermNorm}{}}
\newcommand{\ForcingTermExternalNorm}[1]{\mathring{\mathcal{N}}_{#1}}
\newcommand{\WeightedBEFNorm}[1]{\dot{\mathbf{BEF}}_{#1}}
\newcommand{\CombinedBEFNorm}[1]{\mathbf{BEF}_{#1}}
\newcommand{\BEFNormAux}[1]{\widetilde{\mathbf{BEF}}_{#1}}
\newcommand{\CombinedEFNorm}[1]{\mathbf{EF}_{#1}}
\newcommand{\BulkNormWeighted}[1]{\mathbf{B}_{#1}}
\newcommand{\BulkNormWeightedAux}[1]{\widetilde{\mathbf{B}}_{#1}}
\newcommand{\BulkNormExternal}[1]{\mathring{\mathbf{B}}_{#1}}
\newcommand{\ExternalBEFNorm}[1]{\mathring{\mathbf{BEF}}_{#1}}

\newcommand{\EnergyHorizonDeg}{E_{\operatorname{deg}}}
\newcommand{\EnergyHorizonDegAux}{\dot{E}_{\operatorname{deg}}}

\newcommand{\Trapping}[1]{\tensor[^{(trap)}]{#1}{}}

\newcommand{\rTrapping}{\Trapping{r}}
\newcommand{\TrappedSet}{\mathbf{\Gamma}}

\newcommand{\trap}{\operatorname{trap}}
\newcommand{\nontrap}{\cancel{\trap}}

\newcommand{\WeylQ}[1]{\operatorname{\Op}_w\left(#1\right)}

\newcommand{\OpClass}{\Psi}
\newcommand{\TanSymClass}[1]{S_{\operatorname{tan}}^{#1}}
\newcommand{\TanOpClass}[1]{\Psi_{\operatorname{tan}}^{#1}}

\newcommand{\MixedSymClass}[2]{S^{#1}_{#2}}
\newcommand{\MixedOpClass}[2]{\Psi^{#1}_{#2}}
\newcommand{\Op}{\operatorname{Op}}
\newcommand{\Main}{\operatorname{princ}}
\newcommand{\Aux}{\bowtie}

\newcommand{\NablaAngular}{\nabla}

\newcommand{\EMTensor}{\mathbb{T}}
\newcommand{\Nonlinearity}{\mathcal{N}}

\newcommand{\MorNorm}{\mathbf{Mor}}
\newcommand{\MorNormTrap}{\mathring{\MorNorm}}
\newcommand{\MorrNorm}{\mathbf{Morr}}

\newcommand{\MorDualNorm}[1]{\MorNorm^{*}}
\newcommand{\MorrDualNorm}[1]{\MorrNorm^{*}}

\newcommand{\SdS}{Schwarzschild-de Sitter}
\newcommand{\KdS}{Kerr-de Sitter}

\newcommand{\Horizon}{\mathcal{H}}
\newcommand{\EventHorizon}{\mathcal{H}}
\newcommand{\CosmologicalHorizon}{\overline{\mathcal{H}}}
\newcommand{\EventHorizonFuture}{\mathcal{H}^+}
\newcommand{\CosmologicalHorizonFuture}{\overline{\mathcal{H}}^+}

\newcommand{\DomainOfIntegration}{\mathcal{D}}

\newcommand{\TAlmostKilling}{\widetilde{\KillT}}

\newcommand{\rAux}{\mathfrak{r}}

 \newcommand{\SubPOp}{\mathbf{S}}
 \newcommand{\SubPSym}{\mathbf{s}}

\newcommand{\ForcingTerm}{F}

\newcommand{\FreqAngular}{\eta}
\newcommand{\FreqPhi}{\FreqAngular_{\varphi}}
\newcommand{\FreqTheta}{\FreqAngular_{\theta}}

\newcommand{\PrinSymb}{p}
\newcommand{\RescaledPrinSymb}{{\abs*{q}^2 \PrinSymb}}

\newcommand{\SigmaStar}{\Sigma_{*}}

\newcommand{\HawkingVF}{\widehat{T}}
\newcommand{\HprVF}{\widehat{R}}
\newcommand{\CartarOp}{\mathcal{O}}
\newcommand{\MorawetzVF}{X}
\newcommand{\MorawetzSym}{\mathfrak{x}}
\newcommand{\MorawetzLagrangeCorr}{w}
\newcommand{\MorawetzLagrangeCorrSym}{\mathfrak{w}}
\newcommand{\MorawetzOneForm}{m}
\newcommand{\SquareDecomp}{\mathfrak{a}}
\newcommand{\SquareDecompOp}{\mathfrak{A}}
\newcommand{\SymClass}{S}

\newcommand{\rpVF}{X}
\newcommand{\rpLagrangianCorr}{w}
\newcommand{\rpOneForm}{m}
\newcommand{\rpKCurrent}{\KCurrent{\rpVF_p, \rpLagrangianCorr_p, \rpOneForm_p}[\psi]}
\newcommand{\rpJCurrent}{\JCurrent{\rpVF_p, \rpLagrangianCorr_p, \rpOneForm_p}[\psi]}
\newcommand{\rpBulkWeighted}[2]{\dot{\mathbf{B}}_{{#1};{#2}}}
\newcommand{\rpBulkCombined}[1]{\mathbf{B}_{#1}}

\newcommand{\bea}{\begin{eqnarray}}
\newcommand{\eea}{\end{eqnarray}}
\def\beaa{\begin{eqnarray*}}
\def\eeaa{\end{eqnarray*}}
\def\ba{\begin{array}}
\def\ea{\end{array}}
\def\be#1{\begin{equation} \label{#1}}
\def \eeq{\end{equation}}

\def\a{{\alpha}}

\def\b{{\beta}}
\def\be{{\beta}}
\def\ga{\gamma}

\def\de{\delta}
\def\De{\Delta}

\def\ka{\kappa}

\def\la{\lambda}
\def\La{\Lambda}

\def\Si{\Sigma}
\def\om{\omega}
\def\Om{\Omega}

\def\th{\theta}

\def\ka{\kappa}
\def\nab{\nabla}

\def\c{\cdot}

\def\CC{{\mathcal C}}
\def\MM{{\mathcal M}}

\def\LL{{\mathcal L}}

\def\HH{{\mathcal H}}

\def\DD{{\mathcal D}}

\def\RR{{\mathcal R}}

\def\HH{{\mathcal H}}

\def\D{{\bf D}}
\def\E{{\bf E}}

\def\O{{\bf O}}

\def\R{{\bf R}}

\def\T{{\bf T}}

\def\W{{\bf W}}
\def\g{{\bf g}}

\def\t{{\bf t}}

\def\RRR{{\Bbb R}}

\def\CCC{{\Bbb C}}
\def\f12{{\frac 1 2}}

\def\div{\mathrm{div}}
\def\curl{\mathrm{curl}\,}

\def\trch{{\mbox tr}\, \chi}

\def\chib{{\underline \chi}}

\def\etab{{\underline \eta}}
\def\omb{{\underline{\om}}}

\def\xib{{\underline \xi}}

\def\tr{\mbox{tr}}
\def\atr{\,^{(a)}\mbox{tr}}

\def\trchb{{\tr \,\chib}}

\def\atrch{\atr\chi}
\def\atrchb{\atr\chib}

\def\f12{\frac 1 2}
\def\bsplit{\begin{split}}

\newcommand{\init}{\mathrm{in}}

\newcommand{\half}{\frac{1}{2}}

\usepackage{mathtools}

\newcommand{\dr}{\partial}
\newcommand{\ffi}{\varphi}
\newcommand{\GO}[1]{O\left( #1 \right)}

\DeclareFontFamily{U}{mathx}{\hyphenchar\font45}
\DeclareFontShape{U}{mathx}{m}{n}{
      <5> <6> <7> <8> <9> <10>
      <10.95> <12> <14.4> <17.28> <20.74> <24.88>
      mathx10
      }{}
\DeclareSymbolFont{mathx}{U}{mathx}{m}{n}
\DeclareFontSubstitution{U}{mathx}{m}{n}
\DeclareMathAccent{\widecheck}{0}{mathx}{"71}

\def\dk{\mathfrak{d}}

\renewcommand{\tr}{\mathrm{tr}}
\renewcommand{\atr}{\,^{(a)}\mathrm{tr}}
\newcommand{\pth}[1]{\left( #1 \right)}

\newcommand{\e}{\varepsilon}

\renewcommand{\t}{\mathfrak{t}}

\newcommand{\q}{\mathfrak{q}}

\newcommand{\VFArthur}{\mathbf{Z}}
  \newcommand{\EnergyInit}[1]{\widehat{\mathbf{E}}_{#1}}
  \newcommand{\SigmaInit}{\widehat{\Sigma}_{\init}}
  
  \usepackage[sortcites,
backend=biber,
date=year,
style=alphabetic,
doi=false,
isbn=false,
  url=false,
  maxcitenames=4,
maxbibnames=4,
eprint=true]{biblatex} 

\DeclareNameAlias{author}{family-given}
\AtEveryBibitem{
  \clearfield{month}
  \clearfield{url}
  \clearfield{urlyear}
  \clearfield{urlmonth} 
  \ifentrytype{online}{
    \clearfield{year}}
  {
    \clearfield{eprint}}
}

\renewbibmacro*{volume+number+eid}{%
  \printfield{volume}%
  \setunit*{\addnbspace}%
  \printfield{number}%
  \setunit{\addcomma\space}%
  \printfield{eid}}
\DeclareFieldFormat[article]{volume}{\textbf{#1}}
\DeclareFieldFormat[article]{number}{\mkbibparens{#1}}
\DeclareFieldFormat*{title}{\textit{#1}}
\DeclareFieldFormat{journaltitle}{#1\isdot}
\renewbibmacro{in:}{}
\numberwithin{equation}{section}

\title{Teukolsky on slowly-rotating Kerr-de Sitter \\ in the vanishing $\Lambda$ limit}
\author{Allen Juntao Fang\thanks{University of M\"unster, 
    \href{mailto:allen.juntao.fang@uni-muenster.de}{allen.juntao.fang@uni-muenster.de}},
  J\'er\'emie  Szeftel\thanks{CNRS \& Laboratoire Jacques-Louis Lions, Sorbonne Universit\'e, \href{mailto:jeremie.szeftel@sorbonne-universite.fr}{jeremie.szeftel@sorbonne-universite.fr}} ,
  Arthur Touati\thanks{CNRS \& Institut Mathematiques de Bordeaux, Université de Bordeaux, \href{mailto:arthur.touati@math.u-bordeaux.fr}{arthur.touati@math.u-bordeaux.fr}}}
\date{}
\begin{document}

\maketitle

\begin{abstract}

 As a first step towards resolving a vanishing cosmological constant black hole stability conjecture,
  we prove energy, Morawetz and $r^p$-weighted estimates for solutions
  to the Teukolsky equations on a slowly-rotating Kerr-de Sitter
  background, which we derive using an extension of the non-integrable formalism of \cite{giorgiWaveEquationsEstimates2024}. The main feature of our estimates is their uniformity
  with respect to the cosmological constant $\La>0$ (thus allowed to
  tend to $0$), while they hold on the whole domain of outer
  communications, extending up to $r\sim \La^{-\half}$. As an
  application of our result, we recover well-known corresponding estimates for solutions
  to Teukolsky on a slowly-rotating Kerr background in the limit $\La\to 0$.
\end{abstract}

\tableofcontents

\section{Introduction}

The Einstein vacuum equations (EVE) are the governing equations for
Einstein's theory of relativity under the assumption of a vacuum
universe, and are given by
\begin{equation}
  \label{eq:EVE}
  \Ric(\Metric) - \Lambda \Metric = 0,
\end{equation}
where $(\Manifold, \Metric)$ is a Lorentzian manifold and $\Metric$
has signature $(-,+,+,+)$. In \zcref[noname]{eq:EVE}, $\Lambda$
denotes the cosmological constant, and $\Ric$ the Ricci tensor. While
in general, $\Lambda$ can take any real value, we will only be
interested in the present work in the case where $\Lambda\ge 0$.
As shown in \cite{choquet-bruhatGlobalAspectsCauchy1969,choquet-bruhatTheoremeDexistencePour1952},
the Einstein's vacuum equations \zcref[noname]{eq:EVE} actually possess a hidden hyperbolic
structure and admit a well-posed Cauchy problem. A natural question
that then arises is the existence and stability of special stationary
solutions. In the context of the evolution problem in general
relativity, the most interesting families of stationary solutions are the rotating black hole solutions, namely, the \KdS{} family in
the $\Lambda>0$ setting and the Kerr family in the $\Lambda=0$
setting.

\subsection{The Kerr-de Sitter and Kerr black holes}\label{section intro KdS et K}

The Kerr family $\g_{M,a}$ and the Kerr-de Sitter family
$\g_{M,a,\La}$ are the main families of stationary solutions of
interest in this article and represent rotating black hole solutions
to \zcref{eq:EVE} with $\Lambda=0$ and $\Lambda>0$ respectively. In the standard Boyer-Lindquist coordinates $(t,r,\th,\phi)$, the
metric $\g_{M,a,\La}$ reads
\begin{equation}\label{KdS first expression}
\begin{aligned}
  \g_{M,a,\La}  ={}& -\frac{\Delta}{\left(1+\gamma\right)^2\abs*{q}^2}(dt - a \sin^2\theta d\phi)^2
  + \abs*{q}^2 \left(
    \frac{1}{\De}dr^2 + \frac{1}{\ka}d\theta^2
  \right)
  + \frac{\kappa\sin^2\theta}{(1+\gamma)^2\abs*{q}^2}\left(a\,dt - (r^2+a^2)d\phi\right)^2,
\end{aligned}
\end{equation}
with
\begin{align*}
  \Delta & \vcentcolon= (r^2+a^2)\left(1-\frac{\Lambda}{3}r^2\right)-2Mr,& 
  q & \vcentcolon= r+ \ImagUnit a\cos\theta ,&  \kappa &\vcentcolon= 1+ \gamma\cos^2\theta,&  \gamma & \vcentcolon= \frac{\Lambda a^2}{3},
\end{align*}
where $M$ is the mass of the black hole and $|a|$ its angular momentum
per unit mass. The expression of the Kerr metric $\g_{M,a}$ in
Boyer-Lindquist coordinates can be deduced from \zcref[noname]{KdS
  first expression} by setting $\La=0$. A horizon is subextremal if
its surface gravity is positive, and we say that spacetime is
subextremal if all relevant Killing horizons are subextremal. In the
case of Kerr, since there is only one horizon, namely, the event
horizon, this condition reduces to the condition that the angular
momentum satisfies $\abs*{a}< M$. This is equivalent to the condition
that $\evalAt*{\Delta}_{\Lambda=0}$ has two distinct real roots. For
\KdS, the condition that both the event horizon and the cosmological
horizon are subextremal reduces to the parameter condition
\begin{equation*}
  -\left( 1 + \frac{\Lambda a^2}{3} \right)\left( \frac{a}{M} \right)^2
  + 12\left( 1 - \frac{\Lambda a^2}{3} \right)\Lambda a^2
  + \left( 1 - \frac{\Lambda a^2}{3} \right)^3
  - 9 \Lambda M^2 > 0,
\end{equation*}
which is clearly verified for $\abs*{a}, \Lambda \ll 1$.
At a more practical level, the subextremality condition of \KdS{} is
satisfied if $\Delta$ has four distinct real roots.  For a more
in-depth discussion of the behavior of subextremal waves on \KdS, we
refer the interested reader to
\cite{petersenStationarityFredholmTheory2024}.

Thanks to the explicit forms of the metrics in Boyer-Lindquist
coordinates in \zcref[noname]{KdS first expression}, we immediately
see that the convergence $\g_{M,a,\La}\longrightarrow \g_{M,a}$ as
$\La\to 0$ holds pointwise. However, due to the global geometry of
Kerr and Kerr-de Sitter, this limit is singular. This can be seen on
their Penrose diagrams, see \zcref{fig:penrose}.
\begin{figure}[H]
  \centering
  \begin{minipage}[t]{0.48\linewidth}
    \centering
    \input{Images/Intro-Penrose-KdS.tex}
  \end{minipage}%
  \begin{minipage}[t]{0.48\linewidth}
    \centering
    \input{Images/Intro-Penrose-Kerr.tex}
  \end{minipage}
  \caption{\textit{Penrose diagrams of subextremal \KdS{} (on the left) and subextremal Kerr (on the
    right). The stationary region (also known in the literature as the
    domain of outer communication) is shaded in gray in both
    figures, $r_{\HH,\La}$ and $r_{\overline{\HH},\La}$ are the largest positive roots of $\De$ when $\La>0$, and $r_{\HH,0}$ is the largest root of $\De$ when $\La=0$.}}
  \label{fig:penrose}
\end{figure}
Kerr is an asymptotically flat spacetime and in the region far from
the black hole, is dominated by its Minkowskian structure. \KdS{} on
the other hand, forms a cosmological horizon at
$r\sim \Lambda^{-\frac{1}{2}}$. As a result, the convergence
$\g_{M,a,\La}\longrightarrow \g_{M,a}$ cannot hold globally and is
only uniform on a region where $r\ll \La^{-\half}$.

\subsection{Black hole stability}\label{sec:black hole stability in the introduction}

\subsubsection{Kerr and Kerr-de Sitter stability}

Stationary black hole families such as Kerr and \KdS{} are particularly
interesting given their expected role in the final state conjecture as
physical final states for the universe. A basic requirement for the final state conjecture to hold is that the subextremal Kerr(-de Sitter) family is stable. Let $\mathcal{B}_\La$ denote either the subextremal Kerr-de Sitter ($\La>0$) or Kerr ($\La=0$) family of black hole solutions.
\begin{conjecture}[Stability of $\mathcal{B}_\La$]
  The family $\mathcal{B}_\La$ is a stable family of solutions to
  \zcref{eq:EVE} with $\La\geq 0$, i.e. the
  evolution of a small perturbation of the initial data generating
  $b_0\in \mathcal{B}_\La$ asymptotes in the appropriate sense
  to $b_{\varepsilon}\in \mathcal{B}_\La$.
\end{conjecture}
Black hole stability has been
the subject of intense mathematical research in recent decades. A
brief introduction to the stability of Kerr and \KdS{} is presented
below.

The first asymptotically flat solution of \zcref[noname]{eq:EVE} with $\La=0$ which
was proven to be nonlinearly stable was Minkowski in
the groundbreaking work
\cite{christodoulouGlobalNonlinearStability1993}. The Schwarzschild
family, the non-rotating subfamily of Kerr (i.e. $a=0$), was proven to be stable
under polarized axisymmetry in
\cite{klainermanGlobalNonlinearStability2020} and subsequently for a co-dimension 3 set of initial data in \cite{dafermosNonlinearStabilitySchwarzschild2021}.
Recently, the nonlinear stability of the slowly-rotating Kerr family
was proven in
\cite{klainermanConstructionGCMSpheres2022,klainermanEffectiveResultsUniformization2022,klainermanKerrStabilitySmall2023,giorgiWaveEquationsEstimates2024,shenKerrStabilityExternal2024}. A
common quantitative feature of the above stability results is the slow polynomial decay of
their perturbations, which comes from their asymptotically flat
structure.

The first proof of the nonlinear stability of a stationary solution to
\zcref[noname]{eq:EVE} with $\Lambda>0$ was the global stability of de
Sitter spacetimes in
\cite{friedrichExistence$n$geodesicallyComplete1986}. The first
proof of nonlinear stability of black hole spacetime solutions to
\zcref{eq:EVE} with $\Lambda>0$ was given in
\cite{hintzGlobalNonLinearStability2018} 
to prove the stability of the slowly-rotating \KdS{} family (see also
the alternative proof given in
\cite{fangLinearStabilitySlowlyRotating2026,fangNonlinearStabilitySlowlyRotating2026}). A
critical quantitative feature of the stability of the slowly-rotating
\KdS{} family (which is expected to extend to the full subextremal
family, see \cite{hintzConditionalNonlinearStability2025}) is that
perturbations decay exponentially back to the \KdS{} family. This can
be seen at the level of the linearized Einstein equations (or from more basic hyperbolic equations such as the scalar wave equation), and persists at the level of nonlinear perturbations.

\subsubsection{The vanishing $\La$ black hole stability conjecture}

While both slowly-rotating Kerr and slowly-rotating \KdS{} have been
shown to be qualitatively stable, as mentioned in the previous
sections, the quantitative nature of their stability differs
significantly. Perturbations of the slowly-rotating Kerr family only
decay back to the Kerr family at a polynomial rate, whereas
perturbations of the slowly-rotating \KdS{} family decay exponentially
back to the \KdS{} family. This is largely due to the difference in
the asymptotic structure of the spacetime. \KdS{} features a
cosmological horizon with positive surface gravity, which induces a
sink-source radial point structure at the level of the Hamiltonian
flow \cite{vasyMicrolocalAnalysisAsymptotically2013}. This structure
is annihilated in the $\Lambda\to 0$ limit, where the Hamiltonian flow
on Kerr is instead dominated by its scattering properties at null
infinity. This thus makes the vanishing $\Lambda$ limit a singular
limit similar to the vanishing viscosity problem from fluid dynamics. 
Motivated by this, and astrophysical observations suggesting that $\Lambda$ is extremely small, we formulate the following vanishing $\Lambda$ black hole stability conjecture.

\begin{conjecture}[Vanishing $\Lambda$ black hole stability]
  \label{conj:vanishingcosmologicalconstantlimit}
  Let $\g_{M,a}$ be a subextremal Kerr metric. For any
  suitably small perturbation
  $(\Sigma_0, g_0,
  k_0)$ of the initial data\footnote{Recall that the initial data for \zcref[noname]{eq:EVE} are a triplet $(\Si,g,k)$ formed by $\Si$ a 3-dimensional manifold, $g$ a Riemannian metric on $\Si$ and $k$ a symmetric 2-tensor on $\Si$, and moreover satisfying the so-called constraint equations.} of $\g_{M,a}$, and any family
  of initial data
  $(\Sigma_{\Lambda}, g_\Lambda,
  k_\Lambda)$ that converges in an appropriate sense to $(\Sigma_0, g_0,
  k_0)$ and is itself a small perturbation of the initial data of the subextremal \KdS{} metric $\g_{M,a,\La}$,
  the evolution of
  $(\Sigma_{\Lambda}, g_\Lambda,
  k_\Lambda)$ converges in an appropriate sense to the
  evolution of
  $(\Sigma_{0}, g_0,
  k_0)$ on regions where both $t,r\lesssim  \Lambda^{-\frac{1}{2}}$.
\end{conjecture}

\zcref{conj:vanishingcosmologicalconstantlimit} is the main motivation for this article. In particular, our main results outlined in \zcref{section rough statements} could be interpreted as the resolution of this conjecture in the linearized setting. Note that \zcref{conj:vanishingcosmologicalconstantlimit} describes the behavior of perturbations on $t\lesssim \Lambda^{-\frac{1}{2}}$. For $t\gtrsim \Lambda^{-\frac{1}{2}}$, we expect that perturbations of \KdS{} transition from the $\Lambda$-uniform polynomial decay of  \zcref{conj:vanishingcosmologicalconstantlimit} to the $\Lambda$-dependent exponential decay uncovered in \cite{hintzGlobalNonLinearStability2018}.

\begin{remark}
    We note that the specific singular limit of \KdS{} to Kerr is different from that taken in \cite{hintzModeStabilityShallow2021}, where $\Lambda=3$ and the mass converges to $0$. While the metric in that context formally converges to de Sitter, the Kerr analysis nevertheless plays a role in the near-zero frequency analysis of \cite{hintzModeStabilityShallow2021}.  
\end{remark}

\subsubsection{The Teukolsky equations}\label{section intro teuk}

The Teukolsky equations are wave equations for spin-weighted functions with the spin $s$ taking half-integer values. It encompasses scalar waves ($s=0$), Dirac ($s=\pm\frac{1}{2})$, Maxwell ($s=\pm1)$, and linearized gravity ($s=\pm2$). Given the primary objective of the article is to provide a first step towards \zcref{conj:vanishingcosmologicalconstantlimit} we focus in what follows on the case $s=2$, where the Teukolsky equations are a system of equations for the extreme null components of the
Weyl curvature tensor. 
In particular, the Teukolsky equations capture the
hidden wave-type structure embedded in the Einstein vacuum equations.

While the Teukolsky equations in the $\Lambda=0$
setting are well-studied
\cite{milletOptimalDecaySolutions2023,maUniformEnergyBound2020a,maSharpDecayTeukolsky2023,shlapentokh-rothmanBoundednessDecayTeukolsky2020,shlapentokh-rothmanBoundednessDecayTeukolsky2023,dafermosLinearStabilitySchwarzschild2019,dafermosBoundednessDecayTeukolsky2019}
and have played a critical role in stability proofs in the $\Lambda=0$
setting
\cite{dafermosNonlinearStabilitySchwarzschild2021,klainermanGlobalNonlinearStability2020,klainermanKerrStabilitySmall2023,giorgiWaveEquationsEstimates2024},
they have not been studied in great detail on $\Lambda>0$ black hole
solutions (see however the partial mode stability result in
\cite{casalsHiddenSpectralSymmetries2022}). This is since existing
proofs of stability for  $\Lambda>0$ black hole solutions generally rely on
generalized harmonic gauge, where all metric components are solved for
at the same time in a system of wave equations
\cite{fangNonlinearStabilitySlowlyRotating2026,hintzGlobalNonLinearStability2018}.

Moreover, the Teukolsky equations already capture the qualitative
difference in polynomial versus exponential decay rates for
perturbations of Kerr versus perturbations of \KdS. This can be
observed even from more basic equations such as the scalar wave
equation, as waves on Kerr decay at an inverse polynomial rate
\cite{blueSemilinearWaveEquations2003,anderssonHiddenSymmetriesDecay2015,
  dafermosDecaySolutionsWave2010, dafermosDecaySolutionsWave2016,
  tataruLocalEnergyEstimate2011}, while waves on Kerr-de Sitter decay
exponentially \cite{dyatlovAsymptoticsLinearWaves2015,
  vasyMicrolocalAnalysisAsymptotically2013,
  hintzGlobalAnalysisQuasilinear2016a,
  mavrogiannisQuasilinearWaveEquations2024,bonyDecayNonDecayLocal2008}.

\subsection{Rough statement of the results}\label{section rough statements}

Our main result on the Teukolsky equations on Kerr-de Sitter, which we derive in \zcref{sec:derivation-of-RW},  can be roughly stated as follows.

\begin{theorem}[Rough statement of \zcref{MAINTHEOREM}]{}\label{rough theo 1}
  Solutions to the Teukolsky equations on a slowly-rotating Kerr-de Sitter background satisfy a
  $\Lambda$-uniform integrated local energy decay estimate (consistent with
  polynomial decay) over a region extending up to $r\sim \La^{-\half}$.
\end{theorem}

The uniformity of the estimate with respect to $\La$ allows us to pass to the limit $\Lambda\to0$ and recover, as a corollary of \zcref[cap]{rough theo 1}, the well-known corresponding estimate for Kerr first derived in \cite{maUniformEnergyBound2020a,dafermosBoundednessDecayTeukolsky2019}.

\begin{corollary}[Rough statement of \zcref{theo comparaison}]\label{rough coro}
  Solutions to the Teukolsky equations on a slowly-rotating Kerr background satisfy an integrated local energy decay estimate (consistent with
  polynomial decay).
\end{corollary}

In view of the central role played the Teukolsky equations in the linearized Einstein vacuum equations, the $\Lambda$-uniform estimates in \zcref{rough theo 1} and their convergence in the $\Lambda\to 0$ limit in \zcref{rough coro} could be interpreted as a proof of  \zcref{conj:vanishingcosmologicalconstantlimit} in the linearized setting. Moreover, 
    similarly to how we deduce \zcref[cap]{rough coro} from \zcref[cap]{rough theo 1}, we expect \zcref[cap]{conj:vanishingcosmologicalconstantlimit} to give another proof of stability for Kerr in the vanishing $\La$ limit. 

\subsection{Sketch of the proof}

We now provide a sketch of proof of the main results in this paper, namely the derivation of the Teukolsky equations on Kerr-de Sitter, the uniform estimates with respect to $\La>0$, and the convergence $\La\to 0$.

\subsubsection{Derivation of Teukolsky and Regge-Wheeler}

As explained in \zcref{sec:black hole stability in the introduction}, the Teukolsky equations of linearized gravity are a system of wave equations for the extreme null components of the Weyl curvature tensor, defined with respect to a principal null frame. As in the Kerr case, there are two natural principal null frames in the Kerr-de Sitter spacetime, an ingoing one and an outgoing one, see \zcref{sec:KdS:principal-null-frames}. However, as opposed to the Kerr case, neither of those null frames are globally regular on the stationary region of Kerr-de Sitter: the ingoing frame is irregular at the cosmological horizon while the outgoing one is irregular at the event horizon. 

To circumvent this, we define a global principal null frame in \zcref{sec:global null frame} by conformally gluing the ingoing frame to the outgoing one. The derivation of the Teukolsky equations is thus entirely performed within the conformally invariant formalism first introduced in \cite{giorgiWaveEquationsEstimates2024} and recalled in \zcref{sec:main-eqns-with-conformally-invariant-derivatives}. Moreover, since the global principal null frame is non-integrable (a common feature for rotating black holes), the derivation of the Teukolsky equations as well as all the analysis carried out in this article are based on the horizontal non-integrable formalism of \cite{giorgiWaveEquationsEstimates2024}. We refer the reader to \zcref{prop:Teuk:Teuk-eq-for-A} for the precise form of the Teukolsky equations in this setting.

Finally, to overcome the presence of non-conservative first-order terms, we transform the Teukolsky equations via a physical space version of Chandrasekhar's transformation into a wave-transport system coupling the generalized Regge-Wheeler equation (a wave equation with no problematic first-order terms) to a system of transport equations. This procedure is standard in the $\La=0$ literature (see for instance \cite{dafermosLinearStabilitySchwarzschild2019,giorgiWaveEquationsEstimates2024}), but is new for Kerr-de Sitter. The generalized Regge-Wheeler system is given in \zcref{coro:Teukolsky-wave-transport}, while its derivation from Teukolsky is the content of \zcref[cap]{appendix:proposition T to RW}.

\subsubsection{The uniform estimates}\label{section intro uniform estimates}

Since the decay results of perturbations of \KdS{} are significantly
stronger than their Kerr equivalents, in order to prove uniform
results when $\La$ tends to 0, it is necessary to adapt the
comparatively more complicated methodology of Kerr to prove a
comparatively weaker but uniform in $\Lambda$ polynomial decay for
perturbations of \KdS. The main mechanism used throughout this article is thus an adaptation of the vectorfield multiplier method, where
different multipliers are chosen to overcome the three geometric
obstacles on black hole backgrounds: the event horizon, the trapped
set, and the far region (see \zcref{fig:obstacles}).  We present each of the issues in greater
detail, dividing them into issues that can and can not be handled
perturbatively in the vanishing $\Lambda$ limit.

To explain heuristically the divide between the perturbative and nonperturbative obstacles, observe that the
inverse metric for \KdS{} can be written as
\begin{equation}\label{inverse de KdS = Kerr + dS}
  \Metric_{M,a,\Lambda}^{-1} = \Metric_{M,a}^{-1} + O(r^2\Lambda),
\end{equation}
where the $O$ notation refers to the weight in $r$. As a result, for $r^2\Lambda\ll 1$, \KdS{} can be
viewed as a perturbation of Kerr and since the event horizon and the trapped
set in the slowly-rotating regime are geometric phenomena occuring
in the region $r\le 4M$, they can be handled perturbatively in
the vanishing $\Lambda$ limit. On the other hand, the cosmological
horizon is located at $r=r_{\overline{\mathcal{H}},\La}\sim \sqrt{\frac{3}{\Lambda}}$ so that the $O(r^2\Lambda)$ de
Sitter-type term in \zcref[noname]{inverse de KdS = Kerr + dS} competes directly with the Kerr contribution. Therefore, the main difficulty in achieving $\Lambda$-uniform estimates is in the far region, which cannot be treated perturbatively.

\begin{figure}[H]
  \centering
  \input{Images/intro-penrose-3-obstacles.tex}
  \caption{\textit{A qualitative depiction of the location of the main regions of interest: $A$ is the ergoregion, $B$ is the trapped set, and $C$ is the far region where $r^2\La\sim1$.}}
  \label{fig:obstacles}
\end{figure}

\paragraph{The perturbative obstacles.}

In the vanishing $\Lambda$ limit, both the event horizon and the
trapped null geodesics can be handled in a perturbative manner, we review here the common strategy to handle these phenomena.

In both Kerr and \KdS, the black hole induces an ergoregion
near the event horizon, i.e. a region where the Killing vectorfield
$\partial_t$ becomes spacelike, inducing a loss of ellipticity in energy estimates. However, in the case of subextremal Kerr and \KdS, the positive surface gravity
of the event horizon produces a redshift effect near the event
horizon. This in turn induces a locally exponentially damping effect on waves
near the event horizon, which, in the case of slowly-rotating black
holes, is sufficient to overcome the loss of ellipticity induced by
the ergoregion, via the use of redshift estimates first demonstrated
in \cite{dafermosRedshiftEffectRadiation2009}. These estimates reflect the underlying sink-source
radial point structure at the event horizon
\cite{vasyMicrolocalAnalysisAsymptotically2013} and capture the
physical redshift effect experienced by observers near the black hole.

A more substantial geometric obstacle to integrated local decay
estimates on Kerr and Kerr-de Sitter is the presence of trapped
null-bicharacteristics. For slowly rotating Kerr and Kerr-de Sitter, they remain in a spatial region
of size $O(a)$ around the photon sphere located at $r=3M$.  While these may
seem initially problematic, the trapped set on Kerr and \KdS{} is
normally hyperbolic. This unstable structure admits Morawetz (also
known as integrated local energy decay) estimates that lose some
regularity at the trapped set
\cite{dyatlovSpectralGapsNormally2016,blueSemilinearWaveEquations2003}. In
this article, we follow the approach pioneered in
\cite{tataruLocalEnergyEstimate2011} in the context of studying the scalar wave on Kerr, using a pseudo-differential
modification of the vectorfield method to capture the
frequency-dependent nature of trapping on \KdS.

\paragraph{The nonperturbative obstacle: the far region.} 
As opposed to the two previously mentioned difficulties, the far regions of Kerr and \KdS{} present different
geometric obstacles, since Kerr-de Sitter presents a cosmological
horizon while Kerr is asymptotically flat (see \zcref{fig:penrose}). This difference is clearly seen on the surface gravity $\kappa_{\CosmologicalHorizon,\La}$ of the cosmological horizon of Kerr-de Sitter, which is proportional to the Ricci coefficient evaluated at the cosmological horizon $-\underline{\om}_{|_{r = r_{\CosmologicalHorizon,\La} }}$ in the outgoing null frame, an exact expression of which can be deduced from \zcref{lemma:Kerr:outgoing-PG:Ric-and-curvature}:
\begin{align*}
    -\underline{\om}_{|_{r = r_{\CosmologicalHorizon,\La} }} & = \frac{1}{ r_{\CosmologicalHorizon,\La}^2 + a^2\cos^2\th} \pth{  r_{\CosmologicalHorizon,\La}   - M\pth{ 1   + \frac{2r_{\CosmologicalHorizon,\La}^2}{r_{\CosmologicalHorizon,\La}^2+a^2}  } + \frac{a^2}{r_{\CosmologicalHorizon,\La}} }.
\end{align*}
For a fixed $\La>0$ we can thus show that $\kappa_{\CosmologicalHorizon,\La}>0$ so that redshift estimates capturing the underlying sink-source radial point structure of the cosmological horizon are available. However, this cosmological redshift effect breaks down in the vanishing $\La$ limit since we have $\kappa_{\CosmologicalHorizon,\La}\sim \La^{\half}$ when $\La\to 0$, consistent with the absence of redshift at null infinity in Kerr. To capture the behavior of the wave at the cosmological horizon in a $\Lambda$-uniform manner, we instead adapt the $r^p$-weighted estimates (first shown in \cite{dafermosNewPhysicalSpaceApproach2010}) to \KdS. 

Qualitatively, this results in a $\Lambda$-parametrized family of $r$-weighted estimates capturing the full $H^1$ norm at fixed $\La>0$ but only tangential derivatives to the cosmological horizon uniformly over $\La>0$. Concretely, we prove the $r^p$-weighted estimates by constructing a vectorfield multiplier based on $r^p\partial_r$. In the Kerr case, the asymptotically flat structure of Kerr yields the positivity of the bulk and boundary terms forming the $r^p$-weighted estimates, while in the \KdS{} case, the $O\pth{r^2\Lambda}$ terms in \zcref[noname]{inverse de KdS = Kerr + dS} directly compete at leading order close to the cosmological horizon. Therefore, we must rely on a good sign rather than a perturbative argument using the smallness of $\La$ to close the estimates. For instance, we show that the bulk term $K_{\operatorname{KdS}}[\psi]$ decomposes as
\begin{equation*}
  K_{\operatorname{KdS}}[\psi] =  K_{\operatorname{Kerr}}[\psi] + K_{\Lambda}[\psi] + \text{lower order terms},
\end{equation*}
where the lower order terms are measured with respect to their weight in $r$, and 
where $K_{\Lambda}[\psi]$ is comparable to $K_{\operatorname{Kerr}}[\psi]$ but crucially satisfies the bound
\begin{equation*}
    K_{\Lambda}[\psi] \geq  (\Lambda r^2) r^{p-1}\left( \abs*{\widecheck{\nabla}_4\psi}^2 +  \abs*{r^{-1}\psi}^2\right).
\end{equation*}
We refer the reader to the beginning of \zcref{sec:rp} for a more in-depth discussion about the structure and sign of the terms
involved.

\subsubsection{The vanishing $\La$ limit for Teukolsky}

We comment here on the proof of \zcref{theo comparaison}, which is the main application of the uniform estimates with respect to the cosmological constant $\La$ proved in \zcref{MAINTHEOREM}. First, recall that integrated local energy decay estimates consistent with polynomial decay for solutions to the Teukolsky equations on Kerr were already derived, see the references in \zcref{section intro teuk}. The novelty of our work thus lies entirely in the method of proof of such estimates, namely we derive the estimates in the $\La=0$ case based on the corresponding uniform estimates in the $\La>0$ case discussed in \zcref{section intro uniform estimates} and by passing to the limit $\La\to 0$.

As explained in \zcref{section intro KdS et K}, the vanishing $\La$ limit is based on the pointwise convergence of the Kerr-de Sitter to the Kerr metric when $\La\to 0$. In order to compare the metrics and all the geometric structures attached to them, we define the Kerr manifold $\mathcal{M}_{\mathrm{tot},0}$ and the Kerr-de Sitter manifold $\mathcal{M}_{\mathrm{tot},\La}$ directly in the global coordinates of \zcref{sec:adaptedglobalcoordinates}, as represented on \zcref{fig:KdSonK}\footnote{Note that since $r_{\mathcal{H},\La}$ and $r_{\mathcal{H},0}$ slightly differ, the regions $\mathcal{M}_{\mathrm{tot},0}$ and $\mathcal{M}_{\mathrm{tot},\La}$ don't actually coincide near the event horizons. However since we clearly have $r_{\mathcal{H},\La} = r_{\mathcal{H},0} + O(\La)$ we choose not to distinguish between the two event horizons on \zcref{fig:KdSonK} for the sake of its readability. Moreover, the regions $\mathcal{M}_{\mathrm{tot},0}$ and $\mathcal{M}_{\mathrm{tot},\La}$ actually extend inside both the black hole and the cosmological regions (in the Kerr-de Sitter case), so that the boundary values for $r$ on \zcref{fig:KdSonK} actually differ by a small universal constant. See the beginning of \zcref{section preuve theorem comparaison} for the precise definition of $\mathcal{M}_{\mathrm{tot},0}$ and $\mathcal{M}_{\mathrm{tot},\La}$.}.

\begin{figure}[H]
  \centering
    \input{Images/intro-penrose-3.tex}
  \caption{\textit{Representation of $\mathcal{M}_{\mathrm{tot},\La}$ as a subset of $\mathcal{M}_{\mathrm{tot},0}$.}}
  \label{fig:KdSonK}
\end{figure}

Once the two geometric frameworks are defined respectively on $\mathcal{M}_{\mathrm{tot},0}$ and $\mathcal{M}_{\mathrm{tot},\La}$, the proof of \zcref{theo comparaison} proceeds in two steps. First, given a solution $A^0$ on $\mathcal{M}_{\mathrm{tot},0}$ to the Teukolsky equation of Kerr, we define initial data $A^\La$ for the Teukolsky equation of Kerr-de Sitter on $\Sigma_0\cap \mathcal{M}_{\mathrm{tot},\La}$ with the property that $A^\La\to A^0$ pointwise on $\Sigma_0$. However, the fact that Kerr-de Sitter and Kerr are not close when $r\sim r_{\overline{\mathcal{H}},\La}$ but rather somehow stay in a $O(1)$ neighborhood of each other (recall \eqref{inverse de KdS = Kerr + dS}) translates to the impossibility of proving that the Kerr-de Sitter initial energy of $A^\La$ converges to the Kerr initial energy of $A^0$. Nevertheless, we can prove that the latter uniformly bounds the former, which is enough for our purpose.

Finally, armed with the boundedness (with respect to $\La$) of the Kerr-de Sitter initial energy of $A^\La$ on $\Sigma_0\cap \mathcal{M}_{\mathrm{tot},\La}$, we can apply the uniform estimates of \zcref{MAINTHEOREM} and get boundedness of $A^\La$ on the whole region $\mathcal{M}_{\mathrm{tot},\La}$ in the topology associated to the energy, Morawetz and $r^p$-weighted norms of Kerr-de Sitter. Crucially using the fact that, by construction, every compact $K$ of $\mathcal{M}_{\mathrm{tot},0}$ belongs to $\mathcal{M}_{\mathrm{tot},\La}$ for $\La>0$ small enough depending on $K$, we can pass to the vanishing $\La$ limit and deduce the corresponding energy, Morawetz and $r^p$-weighted estimates in Kerr for $A^0$, as claimed in \zcref{theo comparaison}.

\begin{remark}
    This strategy of proof, namely the construction of appropriate initial data discussed in this section combined with uniform spacetime estimates discussed in \zcref{section intro uniform estimates}, is here applied to the Teukolsky equations but we expect it to extend to the full Einstein vacuum equations and thus prove \zcref{conj:vanishingcosmologicalconstantlimit}. Of course, the case of the Einstein vacuum equations presents many additional difficulties.
\end{remark}

\subsection{Organization of the article}

We give here an overview of the first sections of this article:
\begin{itemize}
\item \zcref[cap]{sec:nonintegrable-structures} introduces part of the
  non-integrable formalism from
  \cite{giorgiGeneralFormalismStability2020,giorgiWaveEquationsEstimates2024}. In
  particular, we define the complexified Weyl null components
  associated to a metric and a null frame, among which the component
  $A$, for which the Teukolsky equation will be derived.
\item \zcref[cap]{sec:KdS} introduces the Kerr-de Sitter spacetime and
  its main features. We define several important vector fields and
  null frames, as well as global coordinates and the spacetime regions
  we will work on. 
\item In \zcref[cap]{sec:derivation-of-RW}, we derive the main equation
  studied in this article, namely the Teukolsky equation on
  Kerr-de Sitter, as well as the Teukolsky wave-transport system,
  i.e. the generalized Regge-Wheeler equation satisfied by a suitable
  second-order quantity.
\item In \zcref[cap]{sec:main-theorem}, we first introduce all the
  energy, Morawetz and $r^p$-weighted norms considered here, and then
  state our main result, \zcref[cap]{MAINTHEOREM}, and its main application, \zcref[cap]{theo comparaison}.
\end{itemize}
The remaining sections are devoted to the proof of
\zcref[cap]{MAINTHEOREM} and \zcref[cap]{theo comparaison}, while the
two appendices prove technical results from \zcref[cap]{sec:KdS} and
\zcref[cap]{sec:derivation-of-RW}.  \zcref[cap]{sec:road-map} gives a
more precise roadmap to the article's main sections, i.e.
\zcref[cap]{sec:wave} to \zcref[cap]{section preuve theorem
  comparaison}.

\subsection{Acknowledgments}

A.J.F. acknowledges support from NSF award DMS-2303241 during the writing
of this paper, and through Germany’s Excellence
Strategy EXC 2044 390685587, Mathematics M\"{u}nster:
Dynamics–Geometry–Structure, from the Alexander von Humboldt
Foundation in the framework of the Alexander von Humboldt
Professorship endowed by the Federal Ministry of Education and
Research.  J.S. is supported by the ERC grant ERC-2023 AdG 101141855
BlaHSt.

\section{Non-integrable structures} 
\label{sec:nonintegrable-structures}

In this section, we review part of the non integrable formalism
introduced in
\cite{giorgiGeneralFormalismStability2020,giorgiWaveEquationsEstimates2024}
that we will need in order to define the Teukolsky equations.

\subsection{Null pairs and horizontal structures}
\label{sec:hor-struc}

Let $(\mathcal{M}, \Metric)$ be a $3+1$ Lorentzian metric, and let
$(e_3,e_4)$ be a fixed pair of null vectors i.e.,
\begin{equation*}
  \Metric(e_3,e_3)=\Metric(e_4,e_4)= 0,\qquad
  \Metric(e_3,e_4)=-2.
\end{equation*}
Given a fixed orientation on $\Manifold$, we will denote by
$\volForm$ the corresponding volume form. 

\begin{definition}
  A vectorfield $X$ is called $(e_3, e_4)$-\emph{horizontal}, or simply
  \emph{horizontal}, if
  \begin{equation*}
    \Metric\left(e_3, X\right) = \Metric\left(e_4, X\right) = 0.
  \end{equation*}
  We denote by $\HorkTensor{}(\Manifold)$ the set of horizontal
  vectorfields of $\Manifold$. We define the induced volume form on
  $\HorkTensor{}(\Manifold)$ by
  \begin{equation}
    \label{eq:hor-induced-vol-form}
    \volFormHor\left(X,Y\right)
    \vcentcolon= \frac{1}{2}\volFormHor\left(X, Y, e_3, e_4\right).
  \end{equation}
\end{definition}
\begin{remark}
  Given any null pair $(e_3,e_4)$, the horizontal vectorfields
  $\HorkTensor{}(\Manifold)$ define a sub-bundle of the tangent
  bundle $T\Manifold$. In the case of Schwarzschild-de Sitter, the
  horizontal sub-bundle is integrable, i.e. for any
  $X, Y\in\HorkTensor{}(\Manifold)$,
  $[X,Y]\in\HorkTensor{}(\Manifold)$. However, in rotating
    Kerr-de Sitter, the horizontal sub-bundle is nonintegrable.
\end{remark}

\begin{definition}
  \label{def:hor-proj}
  Given an arbitrary vectorfield $X$, we denote by $\horProj{X}$ its
  \emph{horizontal projection},
  \begin{equation}
    \label{eq:hor-proj}
    \horProj{X} \vcentcolon= X + \frac{1}{2}\Metric(X, e_3)e_4 + \frac{1}{2}\Metric(X, e_4)e_3. 
  \end{equation}
\end{definition}

A  $k$-covariant tensor-field $U$ is said to be horizontal, $U\in\HorkTensor{k}(\Manifold)$,  
if  for any $X_1,\ldots X_k$ we have 
$$U(X_1,\ldots X_k)=U( ^{(h)} X_1,\ldots  ^{(h)} X_k).$$

Given a horizontal structure defined by $(e_3, e_4)$, we can associate
to it a null frame by choosing horizontal vectorfields $e_1,e_2$ such
that $\horMetric(e_a,e_b)=\delta_{ab}$. By convention, we will say
that $(e_1,e_2)$ are positively oriented on
$\HorkTensor{}(\Manifold)$ if
\begin{equation}
  \label{eq:e1-e2-positiv-orientation}
  \volFormHor(e_1,e_2) = \frac{1}{2}\volFormHor(e_1,e_2,e_3,e_4) = 1.
\end{equation}

\begin{definition}
  \label{def:hor-struc:hor-tensor-field}
  We denote by $\realHorkTensor{0}=\realHorkTensor{0}(\Manifold)$ to
  the set of pairs of real scalar functions on $\Manifold$, by
  $\realHorkTensor{1}=\realHorkTensor{1}(\Manifold)$ the set of real
  horizontal $1$-forms on $\Manifold$, and by
  $\realHorkTensor{2}=\realHorkTensor{2}(\Manifold)$ the set of
  symmetric traceless horizontal real 2-tensors.
\end{definition}

\begin{definition}
We define the dual of a $\xi\in\realHorkTensor{1}$, and $U\in\realHorkTensor{2}$ by
$$\LeftDual{\xi}_a = \volFormHor_{ab}\xi_b,\qquad   \LeftDual{U}_{ab} = \volFormHor_{ac}U_{cb}.$$
\end{definition}

Note that given $\xi, \eta\in\realHorkTensor{1}$ and
$U\in\realHorkTensor{2}$, we have
$$\LeftDual{(\LeftDual{\xi})} = -\xi, \qquad  \LeftDual{U}_{ab} = - \RightDual{U}_{ba},\qquad \LeftDual{\xi}\cdot\eta=-\xi\cdot\LeftDual{\eta}.$$
Also, given $\xi,\eta\in \realHorkTensor{1}$, $U, V\in \realHorkTensor{2}$, we denote
  \begin{gather*}
    \xi\cdot\eta \vcentcolon= \delta^{ab}\xi_a\eta_b,\qquad
    \xi\wedge\eta\vcentcolon= \volFormHor^{ab}\xi_a\eta_b = \eta\cdot\LeftDual{\eta},\qquad
    (\xi\SymTracelessTensorProd\eta)_{ab} \vcentcolon= \frac{1}{2}\left( \xi_a\eta_b + \xi_b\eta_a - \delta_{ab}\xi\cdot\eta \right),\\
    (\xi\cdot U)_a \vcentcolon= \delta^{bc}\xi_bU_{ac}, \qquad U\wedge V \vcentcolon= \volFormHor^{ab}\tensor[]{U}{_a^c}V_{cb}.
  \end{gather*}

\begin{definition}
  \label{def:gamma-hor-chi-chiBar:def}
  For any horizontal vectorfields $X, Y,$ we define
  \begin{equation}
    \label{eq:hor-gamma:def}
    g(X,Y) = \Metric(X,Y)
  \end{equation}
  and
  \begin{equation}
    \label{eq:chi-chiBar:def}
      \chi(X,Y) = \Metric\left(\CovariantDeriv_Xe_4,Y\right),\qquad  \chiBar(X,Y) = \Metric\left(\CovariantDeriv_Xe_3,Y\right),
  \end{equation}
  where $\CovariantDeriv$ denotes the covariant derivative of $\Metric$.
\end{definition}

\begin{remark}
  Observe that $\chi$ and $\chiBar$ are symmetric if and only if the
  horizontal structure is integrable. Also, in the case where the
  horizontal structure is integrable, $g$ is the induced metric, and
  $\chi$, $\chiBar$ are the null second fundamental forms.
\end{remark}

We can consider $g$, $\chi$, $\chiBar$ as defined in 
\zcref{def:gamma-hor-chi-chiBar:def} as horizontal 2-tensor fields by
extending their definition to arbitrary vectorfield $X,Y$ by
\begin{equation*}
  g(X,Y) = \gamma\left(\horProj{X}, \horProj{Y}\right),\qquad  \chi(X,Y) = \chi\left(\horProj{X}, \horProj{Y}\right),\qquad
  \chiBar(X,Y) = \chiBar\left(\horProj{X}, \horProj{Y}\right).
\end{equation*}

Given $U\in \HorkTensor{2}(\Manifold)$, we define its trace $\Trace U$ and anti-trace $\aTrace{U}$ by
$$\Trace U \vcentcolon= \delta^{ab}U_{ab}, \qquad \aTrace{U} \vcentcolon= \volFormHor^{ab}U_{ab}.$$
Accordingly, we decompose $\chi, \chiBar$ as follows 
\begin{align*}
  \chi_{ab} =& \widehat{\chi}_{ab}
  + \frac{1}{2}\delta_{ab}\Trace \chi
  + \frac{1}{2}\volFormHor_{ab}\aTrace{\chi},\qquad
  \chiBar_{ab} = \widehat{\chiBar}_{ab}
  + \frac{1}{2}\delta_{ab}\Trace \chiBar
  + \frac{1}{2}\volFormHor_{ab}\aTrace{\chiBar},
\end{align*}
where $\widehat{\chi}$ and $\widehat{\chiBar}$ denotes the symmetric traceless part of $\chi$ and $\underline{\chi}$ respectively.

\subsection{Horizontal covariant derivatives}
\label{sec:hor-struc:cov-derivs}

We define the horizontal covariant operator $\nabla$ as follows. Given $X, Y\in\HorkTensor{}(\Manifold)$
 \begin{equation}
 \nabla_X Y \vcentcolon= \,^{(h)}(\CovariantDeriv_XY)=\CovariantDeriv_XY- \frac 1 2 \chiBar(X,Y)e_4 -  \frac 1 2 \chi(X,Y) e_3.
 \end{equation}
 In particular, for  all  $X,Y, Z\in\HorkTensor{}(\Manifold)$, we have $Z \horMetric (X,Y)=\horMetric(\nabla_Z X, Y)+ \horMetric(X, \nabla_ZY)$. In the integrable case, $\nabla$ coincides with the Levi-Civita connection
 of the metric induced on the integral surfaces of   $\HorkTensor{}(\Manifold)$.  
 Given $X$ horizontal, $\CovariantDeriv_4X$ and $\CovariantDeriv_3 X$ are in general not horizontal. We define $\nabla_4 X$ and $\nabla_3 X$  to be the horizontal projections
 of the former.  More precisely,
 \begin{align*}
 \nabla_4 X&\vcentcolon= \,^{(h)}(\CovariantDeriv_4 X)=\CovariantDeriv_4 X- \frac 1 2 \Metric(X, \CovariantDeriv_4 e_3 ) e_4- \frac 1 2  \Metric(X, \CovariantDeriv_4 e_4)  e_3 ,\\
 \nabla_3 X&\vcentcolon= \,^{(h)}(\CovariantDeriv_3 X)=\CovariantDeriv_3 X-   \frac 1 2 \Metric(X, \CovariantDeriv_3e_3) e_3 - \frac 1 2   \Metric(X, \CovariantDeriv_3 e_4 ) e_3. 
 \end{align*}
The definition can be easily extended to arbitrary  $\HorkTensor{k}(\Manifold)$ tensor-fields  $U$ 
\begin{align*}
 \nabla_4U(X_1,\ldots, X_k)&= e_4 (U(X_1,\ldots, X_k))- \sum_i U( X_1,\ldots, \nabla_4 X_i, \ldots X_k),\\
  \nabla_3 U(X_1,\ldots, X_k)&= e_3 (U(X_1,\ldots, X_k)) -\sum_i U( X_1,\ldots, \nabla_3 X_i, \ldots X_k).
 \end{align*}
 Given a horizontal $1$-form $\xi$, we define the frame independent operators
\begin{equation}
  \label{eq:div-curl-sym-grad:def}
  \Divergence \xi \vcentcolon= \delta^{ab}\nabla_b\xi_a,\qquad
  \Curl \xi  \vcentcolon= \volFormHor^{ab}\nabla_a\xi_b, \qquad
  (\nabla \SymTracelessTensorProd \xi)_{ba} \vcentcolon= \nabla_b\xi_a + \nabla_a\xi_b - \delta_{ab}(\Divergence \xi). 
\end{equation}
We also introduce the covariant derivative $\dot{\CovariantDeriv}$ acting on mixed tensors of the type $T_k(\Manifold)\otimes\HorkTensor{l}(\Manifold)$, i.e. tensors  of the form  $U_{\nu_1\ldots \nu_k,  a_1\ldots a_l}$, 
for which we define
\begin{align*}
  \dot{\CovariantDeriv}U_{\nu_1\ldots \nu_k,  a_1\ldots a_l}
  ={}& e_\mu(U_{\nu_1\ldots \nu_k,  a_1\ldots a_l}) -U_{\CovariantDeriv_\mu e_{\nu_1}\ldots \nu_k,  a_1\ldots a_l}-\ldots- U_{\nu_1\ldots \CovariantDeriv_\mu e_{\nu_k},  a_1\ldots a_l}\\
&- U_{\nu_1\ldots \nu_k, ^{(h)}(\CovariantDeriv_\mu e_{a_1})\ldots a_l}-  U_{\nu_1\ldots \nu_k,   a_1 \ldots ^{(h)}(\CovariantDeriv_\mu e_{a_l})}.
\end{align*}

\subsection{Ricci coefficients and Weyl components}
\label{sec:ricci-and-curvature}

We define the Ricci coefficients associated to the metric $\Metric$:
\begin{gather*}
  \chi_{ab} \vcentcolon= \Metric(\CovariantDeriv_a e_4, e_b), \qquad
  \chiBar_{ab} \vcentcolon= \Metric(\CovariantDeriv_a e_3, e_b ),\\
  \eta_a\vcentcolon= \frac{1}{2}\Metric(\CovariantDeriv_3 e_4, e_a),\qquad
  \etaBar_a\vcentcolon= \frac{1}{2}\Metric(\CovariantDeriv_4 e_3, e_a),\\
  \xi_a\vcentcolon=\frac{1}{2}\Metric(\CovariantDeriv_4e_4, e_a),\qquad
  \xiBar_a\vcentcolon=\frac{1}{2}\Metric(\CovariantDeriv_3e_3,e_a),\\
  \omega\vcentcolon=\frac{1}{4}\Metric(\CovariantDeriv_4e_4,e_3),\qquad
  \omegaBar\vcentcolon= \frac{1}{4}\Metric(\CovariantDeriv_3e_3,e_4),\\
  \zeta_a\vcentcolon=\frac{1}{2}\Metric(\CovariantDeriv_ae_4,e_3). 
\end{gather*}
They allow to express each covariant derivative of the frame vector fields:
\begin{equation}\label{eq:Ricci-formulas}
\begin{aligned}
\CovariantDeriv_a e_b & = \nabla_ae_b + \frac{1}{2}\chi_{ab}e_3 + \frac{1}{2}\chiBar_{ab}e_4,&
    \CovariantDeriv_a e_4 & = \chi_{ab}e_b -\zeta e_4,&
    \CovariantDeriv_a e_3 & = \chiBar_{ab}e_b -\zeta e_3,
\\ \CovariantDeriv_3 e_a & = \nabla_3 e_a + \eta_{a}e_3 + \xiBar_a e_4, &
    \CovariantDeriv_3 e_3 & = -2\omegaBar e_3 + 2\xiBar_b e_b,& 
    \CovariantDeriv_3 e_4 & = 2\omegaBar e_4 + 2\eta_b e_b,
\\ \CovariantDeriv_4 e_a & = \nabla_4 e_a + \etaBar_a e_4 + \xi_a e_3,&
    \CovariantDeriv_4 e_4 & = -2\omega e_4 + 2\xi_b e_b, &
    \CovariantDeriv_4 e_3 & = 2\omega e_3 + 2\etaBar_b e_b.
\end{aligned}
\end{equation}
We now define the following Weyl components,
\begin{align}
  \label{eq:curvature-components:def}
    \alpha_{ab} & \vcentcolon= \Weyl_{a4b4}, &
    \beta_a & \vcentcolon= \frac{1}{2}\Weyl_{a434}, &
    \rho & \vcentcolon= \frac{1}{4}\Weyl_{3434}, &
    \LeftDual{\rho} & \vcentcolon= \frac{1}{4}\LeftDual{\Weyl}_{3434}, &
    \betaBar_a  &\vcentcolon= \frac{1}{2}\Weyl_{a334}, &
    \alphaBar_{ab} & \vcentcolon= \Weyl_{a3b3},
\end{align}
where $\LeftDual{\Weyl}$ denotes the Hodge dual of the Weyl tensor
$\Weyl$. Recall the decomposition of the Riemann curvature in terms of the Weyl tensor, the Ricci tensor and the scalar curvature:
\begin{equation}
  \label{WeylRiemanngeneral}
  \Riem_{\alpha\beta\gamma\delta}=\Weyl_{\alpha\beta\gamma\delta}
  +\frac{1}{2} (\Metric_{\beta\delta}\Riem_{\alpha\gamma}+\Metric_{\alpha\gamma}\Riem_{\beta\delta}-\Metric_{\beta\gamma}\Riem_{\alpha\delta}-\Metric_{\alpha\delta}\Riem_{\beta\gamma})
  + \frac{\Riem}{6}\left(\Metric_{\alpha\delta}\Metric_{\beta\gamma} -\Metric_{\alpha\gamma}\Metric_{\beta\delta}\right)
  .
\end{equation}
In the specific case where the metric $\Metric$ is a solution to EVE with
$\Lambda>0$ as given in \zcref[noname]{eq:EVE}, this becomes
\begin{equation*}
  \Riem_{\alpha\beta\gamma\delta}
  = \Weyl_{\alpha\beta\gamma\delta}
  + \frac{\Lambda}{3}\left(
    \Metric_{\beta\delta}\Metric_{\alpha\gamma}
    -\Metric_{\beta\gamma}\Metric_{\alpha\delta}
  \right)
  .
\end{equation*}
In this case, we thus obtain
\begin{equation}\label{link Riemann Weyl}
\begin{gathered}
  \Riem_{a33b} =  -\alphaBar_{ab},\qquad \Riem_{a334} =  2\betaBar_{a},\qquad
  \Riem_{a44b} =  -\alpha_{ab},\\   
  \Riem_{a443} =  -2\beta_{a},\qquad
  \Riem_{a34b} =  \left(\rho + \frac{2\Lambda}{3}\right) \horMetric_{ab} - \sigma\volFormHor_{ab} ,\\
  \Riem_{a3cb} =  \volFormHor_{cb}\LeftDual{\betaBar}_{a},\qquad   \Riem_{a4cb} = -\volFormHor_{cb}\LeftDual{\beta}_{a},\qquad
  \Riem_{3434} = 4\rho - \frac{4\Lambda}{3},\\
  \Riem_{abcd} = \volFormHor_{ab}\volFormHor_{cd}\rho
                 + \frac{2\Lambda}{3}\left(
                 \horMetric_{bd}\horMetric_{ac}
                 - \horMetric_{bc}\horMetric_{ad}
                 \right),
  \qquad
  \Riem_{ab34} = 2\volFormHor_{ab}\LeftDual{\rho}.
\end{gathered}
\end{equation}
\begin{remark}
  \label{remark:Weyl-Riem-difference}
  Note that the only components that differ between $\Weyl$ and $\Riem$
  are $\Riem_{abcd}$, $\Riem_{3434}$, and $\Riem_{a34b}$.  
\end{remark}

\subsection{Properties of the horizontal covariant derivatives}

 In this section, we give various properties of the horizontal covariant derivatives, and introduce useful objects such as the Hodge operators. We start with the non-integrable analogue of the classical Gauss equation.

\begin{proposition}
  \label{prop:gauss-eqn:tensor-adapted}
  The following identities hold.
  \begin{enumerate}
  \item For $\psi\in \realHorkTensor{0}$:
    \begin{equation}
      \label{eq:gauss-eqn:scalar}
      \left[\nabla_a,\nabla_b\right]\psi
      =  \frac{1}{2}\left(
        \aTrace{\chi}\nabla_3
        + \aTrace{\chiBar}\nabla_4
      \right)\psi \volForm_{ab}.
    \end{equation}

   \item For $\psi\in \realHorkTensor{k}$ for $k=1,2$,
    \begin{equation}
      \label{eq:gauss-eqn:k-1-2}
      \left[\nabla_a,\nabla_b\right]\psi
      = \left(
        \frac{1}{2}\left(
          \aTrace{\chi}\nabla_3
          +\aTrace{\chiBar}\nabla_4
        \right)\psi
        + k\horProj{K}\LeftDual{\psi}
      \right)\volForm_{ab},
    \end{equation}
    where
    \begin{equation}
      \label{eq:h-K:def}
      \horProj{K}
      \vcentcolon= -\frac{1}{4}\Trace\chi\Trace\chiBar
      - \frac{1}{4}\aTrace{\chi}\aTrace{\chiBar}
      + \frac{1}{2}\chiTF\cdot\chiBarTF
      - \frac{1}{4}\Riem_{3434}.
    \end{equation}
  \end{enumerate}
\end{proposition}

\begin{proof}
  See Proposition 2.1.43 in \cite{giorgiWaveEquationsEstimates2024}.
\end{proof}

We now define the Hodge type operators we will make use of in what
follows. 
\begin{definition}
  \label{def:hodge-type-operators}
  Given an orthonormal basis of horizontal vectors $(e_1,e_2)$, we
  define the following \emph{Hodge type operators} (following
  \cite{christodoulouGlobalNonlinearStability1993}).
  \begin{gather*}
    \HodgeOp{1}: \realHorkTensor{1}\to \realHorkTensor{0}, \qquad
    \HodgeOp{1}\xi\vcentcolon=(\Divergence\xi, \Curl\xi),\qquad
    \HodgeOp{2}: \realHorkTensor{2}\to \realHorkTensor{1}, \qquad
    \HodgeOp{2}u_a\vcentcolon=\nabla^bu_{ab},\\
    \HodgeOpDual{1}: \realHorkTensor{0}\to \realHorkTensor{1},\qquad
    \HodgeOpDual{1}(f,f_{*})\vcentcolon= -\nabla_af + \volFormHor_{ab}\nabla_bf_{*},\qquad
    \HodgeOpDual{2}:\realHorkTensor{1}\to\realHorkTensor{2},\qquad
    \HodgeOpDual{2}\xi\vcentcolon= \frac{1}{2}\nabla\SymTracelessTensorProd\xi.
  \end{gather*}
\end{definition}

These operators will be used for elliptic estimates, themselves based on the following Bochner type identities for the horizontal Laplacian (see Proposition 2.1.48 in \cite{giorgiWaveEquationsEstimates2024}). 
\begin{lemma}
  \label{lemma:laplacian-bochner}
  We have the following identities regarding the horizontal Laplacian
  $\LaplaceHor_k \vcentcolon= \nabla^a\nabla_a$.
  \begin{enumerate}
  \item Let $\psi$ be a scalar function. Then
    \begin{equation}
      \label{eq:laplacian-bochner:s0}
      \abs*{\LaplaceHor_0}^2
      = \abs*{\nabla^2\psi}^2
      + \horProj{K}\abs*{\nabla\psi}^2
      + \operatorname{Err}_0[\LaplaceHor_0\psi]
      + \Divergence_0[\LaplaceHor_0\psi],
    \end{equation}
    where
    \begin{equation*}      
      \begin{split}
        \operatorname{Err}_0[\LaplaceHor_0\psi]
        ={}& -\frac{1}{2}\nabla\psi\cdot \left( \aTrace{\chi}\nabla_3 + \aTrace{\underline{\chi}}\nabla_4 \right)\LeftDual{\nabla}\psi,\\
        \Divergence_0[\LaplaceHor_0\psi]
        ={}& \nabla_a\left( \nabla^a\psi\cdot\LaplaceHor_0\psi - \frac{1}{2}\nabla^a\abs*{\nabla\psi}^2 \right).
      \end{split}
    \end{equation*}
  \item Let $\psi\in \realHorkTensor{1}(\Manifold)$. Then
    \begin{equation}
      \label{eq:laplacian-bochner:s1}
      \abs*{\LaplaceHor_1\psi}^2
      = \abs*{\nabla^2\psi}^2
      + \horProj{K}\left(\abs*{\nabla\psi}^2 - 2\horProj{K}\abs*{\psi}^2\right)
      + \operatorname{Err}_1[\LaplaceHor_1\psi]
      + \Divergence_1[\LaplaceHor_1\psi],
    \end{equation}
    where
    \begin{equation*}
      \begin{split}
        \operatorname{Err}_1[\LaplaceHor_1\psi]
        ={}& -\frac{1}{2}\nabla\psi\cdot \left( \aTrace{\chi}\nabla_3 + \aTrace{\underline{\chi}}\nabla_4 \right)\LeftDual{\nabla}\psi
             + \frac{1}{2}\abs*{\left(\aTrace{\chi}\nabla_3 + \aTrace{\underline{\chi}}\nabla_4\right)\psi}^2\\
            & - \frac{3}{2}\horProj{K}\left( \aTrace{\chi}\nabla_3 + \aTrace{\underline{\chi}}\nabla_4 \right)\psi\cdot\LeftDual{\psi}
             + \horProj{K}\Divergence[\HodgeOp{1}\psi]
             ,\\
        \Divergence_1[\LaplaceHor_1\psi]
        ={}& \nabla_a\left( \nabla^a\psi\cdot\LaplaceHor_1\psi
             - \nabla_c\psi\cdot\nabla^c\nabla^a\psi
             \right).
      \end{split}
    \end{equation*}
  \item Let $\psi\in \realHorkTensor{2}(\Manifold)$. Then
    \begin{equation}
      \label{eq:laplacian-bochner:s2}
      \abs*{\LaplaceHor_2\psi}^2
      = \abs*{\nabla^2\psi}^2
      + \horProj{K}\left(\abs*{\nabla\psi}^2 - 6\horProj{K}\abs*{\psi}^2\right)
      + \operatorname{Err}_2[\LaplaceHor_2\psi]
      + \Divergence_2[\LaplaceHor_2\psi],
    \end{equation}
    where
    \begin{equation*}
      \begin{split}
        \operatorname{Err}_2[\LaplaceHor_2\psi]
        ={}& -\frac{1}{2}\nabla\psi\cdot \left( \aTrace{\chi}\nabla_3 + \aTrace{\underline{\chi}}\nabla_4 \right)\LeftDual{\nabla}\psi
             + \frac{1}{2}\abs*{\left(\aTrace{\chi}\nabla_3 + \aTrace{\underline{\chi}}\nabla_4\right)\psi}^2\\
            & - 3\horProj{K}\left( \aTrace{\chi}\nabla_3 + \aTrace{\underline{\chi}}\nabla_4 \right)\psi\cdot\LeftDual{\psi}
             + 2\horProj{K}\Divergence[\HodgeOp{2}\psi]
             ,\\
        \Divergence_2[\LaplaceHor_2\psi]
        ={}& \nabla_a\left( \nabla^a\psi\cdot\LaplaceHor_1\psi
             - \nabla_c\psi\cdot\nabla^c\nabla^a\psi
             \right).
      \end{split}
    \end{equation*}
  \end{enumerate}
\end{lemma}

\subsection{The wave operator acting on horizontal tensors}

We can now define the main operator we will work with. 
\begin{definition}
  \label{def:wave-op}
  We define the wave operator acting on $\psi\in \realHorkTensor{k}(\Manifold)$ by
  \begin{equation}
    \label{eq:horizontal-k-tensor-wave:def}
    \WaveOpHork{k}\psi \vcentcolon= \Metric^{\mu\nu}\HorCovDeriv_{\mu}\HorCovDeriv_{\nu}\psi.
  \end{equation}
\end{definition}

\begin{proposition}
  \label{prop:hor-cov-deriv:hor-vf-curvature-formula}
  For a tensor $\psi\in \HorVFSpace(\Manifold)$, we have the curvature
  formula
  \begin{equation}
    \label{eq:horizontal-structure:horizontal-curvature-formula}
    \left(
      \dot{\CovariantDeriv}_\mu\dot{\CovariantDeriv}_{\nu}
      - \dot{\CovariantDeriv}_{\nu}\dot{\CovariantDeriv}_\mu
    \right) \Psi_{a} = \dot{\Riem}_{ab\mu\nu}\Psi^{b},
  \end{equation}
  where, with connection coefficients given by
  $(\Lambda_{\alpha})_{\beta\gamma} = \Metric\left(\CovariantDeriv_{\alpha}e_{\gamma},e_{\beta}\right)$,
  \begin{equation}
    \label{eq:horizontal-structure:B-horizontal-curvature:def}
    \begin{split}
      \dot{\Riem}_{ab\mu\nu}
      & \vcentcolon= \Riem_{ab\mu\nu} + \frac{1}{2}\BHorizontalCurv_{ab\mu\nu},\\
      \BHorizontalCurv_{ab\mu\nu}
      & \vcentcolon= (\Lambda_{\mu})_{3a}(\Lambda_{\nu})_{b4}
        + (\Lambda_{\mu})_{4a}(\Lambda_{\nu})_{b3}
        - (\Lambda_{\nu})_{3a}(\Lambda_{\mu})_{b4}
        - (\Lambda_{\nu})_{4a}(\Lambda_{\mu})_{b3}.
    \end{split}
  \end{equation}
  More generally, for a mixed tensor $\Psi\in T_{1}(\Manifold)\otimes \HorVFSpace_{1}(\Manifold)$, we have
  \begin{equation*}
    \left(
      \dot{\CovariantDeriv}_{\mu}\dot{\CovariantDeriv}_{\nu}
      - \dot{\CovariantDeriv}_{\nu}\dot{\CovariantDeriv}_{\mu}
    \right)\Psi_{\lambda a}
    = \tensor{\Riem}{_{\lambda}^{\sigma}_{\mu\nu}}\Psi_{\sigma a}
    + \tensor{\dot{\Riem}}{_{a}^{b}_{\mu\nu}}\Psi_{\lambda b},
  \end{equation*}
  with an immediate generalization to tensors $\Psi\in T_k (\Manifold)\otimes\HorVFSpace_{l}(\Manifold)$.
\end{proposition}

\begin{proof}
See Proposition 2.1.27 in \cite{giorgiWaveEquationsEstimates2024}.
\end{proof}

\begin{proposition}
  \label{prop:B-hor:components}
  The components of $\BHorizontalCurv$ are given by
  \begin{align*}
    \BHorizontalCurv_{abc3}
    &=  - \Trace\chiBar\left(\delta_{ca}\eta_{b} - \delta_{ca}\eta_{a}\right)
      -\aTrace{\chiBar}\left(\volFormHor_{ca}\eta_{b}-\volFormHor_{cb}\eta_{a}\right)
      +2\left(
      -\chiBarTF_{ca}\eta_{b}  +\chiBarTF_{cb}\eta_{a} -\chi_{ca}\xiBar_{b} + \chi_{cb}\xiBar_{a}
      \right),\\
    \BHorizontalCurv_{abc4}
    &=  - \Trace\chi\left(\delta_{ca}\etaBar_{b} - \delta_{ca}\etaBar_{a}\right)
      -\aTrace{\chi}\left(\volFormHor_{ca}\etaBar_{b}-\volFormHor_{cb}\etaBar_{a}\right)
      +2\left(
      -\chiTF_{ca}\etaBar_{b}  +\chiTF_{cb}\etaBar_{a} -\chiBar_{ca}\xi_{b} + \chiBar_{cb}\xi_{a}
      \right),\\
    \BHorizontalCurv_{abcd} &= \left(
                              -\frac{1}{2}\Trace\chi \Trace\chiBar
                              - \frac{1}{2}\aTrace{\chi}\aTrace{\chiBar}
                              + \chiTF\cdot\chiBarTF
                              \right)\volFormHor_{ab}\volFormHor_{cd}.
  \end{align*}
\end{proposition}

\begin{proof}
  See Proposition 2.2.4 in
  \cite{giorgiWaveEquationsEstimates2024}.
\end{proof}

We also recall the definition of the horizontal Lie
derivative from Section 2.2.8 of
\cite{giorgiWaveEquationsEstimates2024} for the sake of
completeness. 

\begin{definition}[Horizontal Lie derivatives]
  \label{def:hor-Lie-derivative:hor-tensor}
  Given vectorfields $X,Y$, the horizontal Lie derivative $\HorLieDeriv$ is given by
  \begin{equation}
    \label{eq:hor-Lie-derivative:hor-tensor:vectorfield}
    \HorLieDeriv_XY \vcentcolon= \LieDerivative_XY
    + \frac{1}{2}\Metric(\LieDerivative_XY,e_3)e_4
    + \frac{1}{2}\Metric(\LieDerivative_XY, e_4)e_3.
  \end{equation}
  Given a horizontal covariant $k$-tensor $U$, the horizontal Lie
  derivative $\HorLieDeriv_XU$ is defined to be the projection of
  $\LieDerivative_XU$ to the horizontal space. Thus, for horizontal
  indices $A = a_1\cdots a_k$,
  \begin{equation}
    \label{eq:hor-Lie-derivative:hor-tensor:hor-tensors}
    \left(\HorLieDeriv_XU\right)_A\vcentcolon= \nabla_XU_A
    + \CovariantDeriv_{a_1}X^bU_{b\cdots a_k}
    + \cdots +
    \CovariantDeriv_{a_k}X^bU_{a_1\cdots b}.
  \end{equation}
\end{definition}

\subsection{Complex notations}
\label{sec:complexnotationsRicciandcurvature}

In this section, we recall the complex notations introduced first in \cite{giorgiWaveEquationsEstimates2024}.
Recall the set of real horizontal $k$-tensors
$\realHorkTensor{k} = \realHorkTensor{k}\left(\mathcal{M},
  \Real\right)$ on $\mathcal{M}$ from \zcref[cap]{def:hor-struc:hor-tensor-field}.
We now define the complexified version of horizontal tensors on $\Manifold$.

\begin{definition}
  \label{def:complex-horizontal-tensors}
  We denote by
  $\realHorkTensor{k}(\Complex) = \realHorkTensor{k}(\Manifold,
  \Complex)$ the set of complex anti-self dual $k$-tensors on
  $\Manifold$. In particular,
  \begin{itemize}
  \item $a+\ImagUnit b\in \realHorkTensor{0}(\Complex)$ is a complex
    scalar function on $\Manifold$ if
    $a,b\in \realHorkTensor{0}(\Real)$.
  \item
    $F = f + \ImagUnit\LeftDual{f}\in \realHorkTensor{1}(\Complex)$ is
    a complex anti-self dual $1$-tensor on $\Manifold$ if
    $f\in \realHorkTensor{1}(\Real)$.
  \item $U = u+\ImagUnit\LeftDual{u}\in \realHorkTensor{2}(\Complex)$
    is a complex anti-self dual symmetric traceless 2-tensor on
    $\Manifold$ if $u\in \realHorkTensor{2}(\Real)$.
  \end{itemize}
\end{definition}

\begin{definition}
  \label{def:complex-curvature-and-ricci:def}
  We define the following complexified Weyl components
  \begin{equation}
    \label{eq:complex-curvature:def}
    A:=\alpha+i\LeftDual\alpha, \quad B:=\beta+i\LeftDual\beta, \quad P:=\rho+i\LeftDual\rho,\quad \BBar:=\betaBar+i\LeftDual\betaBar, \quad \ABar:=\alphaBar+i\LeftDual\alphaBar,
  \end{equation}    
  with $A, \ABar\in\realHorkTensor{2}(\mathbb{C})$, $B, \BBar\in\realHorkTensor{1}(\mathbb{C})$, $P\in\realHorkTensor{0}(\mathbb{C})$, and the following complexified Ricci coefficients
  \begin{equation}
    \label{eq:complex-ricci:def}
    \begin{gathered}
      X=\chi+i\LeftDual\chi, \quad \XBar=\chiBar+i\LeftDual\chiBar, \quad H=\eta+i\LeftDual\eta, \quad \HBar=\etaBar+i\LeftDual\etaBar, \quad Z=\zeta+i\LeftDual\zeta, \\ 
      \Xi=\xi+i\LeftDual\xi, \quad \XiBar=\xiBar+i\LeftDual\xiBar,
    \end{gathered}
  \end{equation}
  with $\widehat{X}, \widehat{\XBar}\in\realHorkTensor{2}(\mathbb{C})$,
  $H, \HBar, Z, \Xi, \XiBar\in\realHorkTensor{1}(\mathbb{C})$, where
  $\widehat{X}, \widehat{\XBar}$, as well as $\Trace X, \Trace\XBar$ are
  given by
  $$\Trace X = \Trace\chi -i\aTrace{\chi}, \quad \widehat{X}=\widehat{\chi}+i\LeftDual{\widehat{\chi}}, \quad \Trace\widehat{\XBar} = \Trace\chiBar -i\aTrace{\chiBar}, \quad \widehat{\XBar}=\widehat{\chiBar}+i\LeftDual{\widehat{\chiBar}}.$$
\end{definition}

\begin{definition}
  \label{def:complex-derivatives}
  We define derivatives of complex quantities as follows.
  \begin{enumerate}
  \item For two scalar real-valued functions $a$ and $b$, we define
    \begin{equation*}
      \ComplexDeriv(a+\ImagUnit b)
      \vcentcolon= \left(\nabla + \ImagUnit\LeftDual{\nabla}\right)(a+\ImagUnit b).
    \end{equation*}
  \item For a real-valued $1$-form $f$, we define
    \begin{equation*}
      \ComplexDeriv\cdot (f+\ImagUnit\LeftDual{f})
      \vcentcolon= \left(\nabla + \ImagUnit\LeftDual{\nabla}\right)\cdot (f + \ImagUnit\LeftDual{f}),
    \end{equation*}
    and
    \begin{equation*}
      \ComplexDeriv\SymTracelessTensorProd(f+\ImagUnit\LeftDual{f})
      \vcentcolon= \left(\nabla + \ImagUnit\LeftDual{\nabla}\right)\SymTracelessTensorProd(f+ \ImagUnit\LeftDual{f}).
    \end{equation*}
  \item For a symmetric traceless 2-tensor $u$, we define
    \begin{equation*}
      \ComplexDeriv\cdot (u+\ImagUnit\LeftDual{u})
      \vcentcolon= (\nabla + \ImagUnit\LeftDual{\nabla})\cdot (u+\ImagUnit \LeftDual{u}). 
    \end{equation*}
  \end{enumerate}
\end{definition}

\subsection{Conformally invariant derivatives}
\label{sec:main-eqns-with-conformally-invariant-derivatives}

In this section, we recall conformally invariant derivatives introduced first in \cite{giorgiWaveEquationsEstimates2024}. Their use is motivated by the fact that we will often consider families of horizontal tensors depending on conformally equivalent null frames. For a given null pair $(e_4,e_3)$, define its conformal class $[(e_4,e_3)]$ to be the set of all null pairs $(e_4',e_3')$ such that
\begin{align}\label{eq:e3-e4-rescale-conformal-transformation}
    e_{4}' = \lambda e_{4},\qquad
  e_{3}' = \lambda^{-1} e_{3},
\end{align}
for some scalar function $\la>0$. Crucially, note that the horizontal structure associated to a null pair $(e_4,e_3)$ only depends on its conformal class since it is invariant under \zcref[noname]{eq:e3-e4-rescale-conformal-transformation}. A common choice of positively oriented horizontal vectorfields $e_1,e_2$ is thus made for all elements of $[(e_4,e_3)]$.

\begin{definition}\label{def:s-conformally-invariants}
    Let $(e_4,e_3)$ be a null pair and denote by $\O(\MM)$ its horizontal structure. For $k\geq 0$ and $s\in\RRR$, we say that the family of elements of $\O_k(\MM)$
    \begin{equation}\label{bourbaki}
        U = \pth{U^{(e'_4,e'_3)}}_{(e'_4,e'_3)\in[(e_4,e_3)]}
    \end{equation}
    is $s$-conformally invariant if $U^{(\la e_4,\la^{-1}e_3)} = \la^s U^{(e_4,e_3)}$ for all scalar functions $\la>0$.
\end{definition}

Natural and important examples of families of horizontal tensors enjoying conformal invariance properties are (some of) the Ricci coefficients and the Weyl components. Indeed one can easily check that
\begin{gather*}
  \label{eq:conformally-invariant-ricci-curvature}
  \Trace\chiBar ' = \lambda^{-1}\Trace\chiBar,\qquad
  \aTrace{\chiBar'} = \lambda^{-1}\aTrace{\chi}, \qquad
  \Trace\chi' = \lambda\Trace\chi,\qquad
  \aTrace{\chi'} = \lambda\aTrace{\chi},\\
  \xi' = \lambda^{2}\xi,\qquad
  \eta' = \eta, \qquad
  \etaBar' = \etaBar,\qquad
  \xiBar' = \lambda^{-2}\xiBar,\\
  \alpha' = \lambda^{2}\alpha,\qquad
  \beta' = \lambda \beta,\qquad
  \rho' = \rho,\qquad
  \LeftDual{\rho}' = \LeftDual{\rho},\qquad
  \betaBar' = \lambda^{-1}\betaBar,\qquad
  \alphaBar' = \lambda^{-2}\alphaBar,
\end{gather*}
where the primed quantities are associated to $(\la e_4,\la^{-1}e_3)$ and the unprimed quantities are associated to $(e_4,e_3)$. The Ricci coefficients $\om$, $\underline{\om}$ and $\zeta$ don't define $s$-conformally invariant maps since they transform as
\begin{equation}
  \label{eq:conformally-not-invariant-ricci-curvature}
  \omegaBar' = \lambda^{-1}\left(
    \omegaBar + \frac{1}{2}e_{3}(\log \lambda)
  \right),\qquad
  \omega' = \lambda\left(
    \omega - \frac{1}{2}e_{4}(\log \lambda)
  \right),\qquad
  \zeta' = \zeta - \nabla(\log\lambda).
\end{equation}

Once a null pair $(e_4,e_3)$ is fixed, the notion of horizontal covariant derivatives of tensors can be naturally extended to families of the form \zcref[noname]{bourbaki}, by simply setting $\pth{\nab_\mu U}^{(e'_4,e'_3)}= \nab_{e'_\mu}U^{(e_4',e_3')}$ for $\mu=1,2,3,4$. Observe that if $U$ is $s$-conformally invariant in the sense of \zcref[cap]{def:s-conformally-invariants}, then the family $\nab_\mu U$ does not enjoy any conformal invariance property. We correct for this by defining the conformally invariant derivatives.

\begin{lemma}\label{lemma:cnab}
    Let $(e_4,e_3)$ be a null pair, $\O(\MM)$ its horizontal structure and $U$ as in \zcref[noname]{bourbaki} being $s$-conformally invariant. The families of horizontal tensors
    \begin{align*}
        \,^{(c)}\nab_3 U &  \vcentcolon = \nab_3 U - 2 s \underline{\om}U,
        \\ \,^{(c)}\nab_4 U &  \vcentcolon = \nab_4 U + 2 s {\om}U,
        \\ \,^{(c)}\nab_a U &  \vcentcolon = \nab_a U + s \zeta_a U,
    \end{align*}
    are respectively $(s-1)$-, $(s+1)$- and $s$-conformally invariant.
\end{lemma}

\begin{proof}
    Straightforward computation using \zcref[noname]{eq:conformally-not-invariant-ricci-curvature}.
\end{proof}

In most of this article, we will be interested in a specific conformal
class $[(e_4,e_3)]$ (see already \zcref[cap]{sec:global null frame})
and it will always be clear which particular representative
$(e'_4,e'_3)\in[(e_4,e_3)]$ is being considered, and thus which
element of $\pth{U^{(e'_4,e'_3)}}_{(e'_4,e'_3)\in[(e_4,e_3)]}$ is
being considered. Therefore, to avoid cumbersome notations, we will
write $U$ to denote any element of this family as an abuse of
notation.

\begin{definition}
  \label{def:conformally-invariant-derivatives}
  Let $(e_4,e_3)$ be a null pair, and consider $s$-conformally invariant tensors. We define the following conformally invariant complex angular
  derivatives.
  \begin{enumerate}
  \item For $a+\ImagUnit b\in \realHorkTensor{0}(\Complex)$, we define
    \begin{equation*}
      \ConformalComplexDeriv(a+\ImagUnit b) \vcentcolon= \left(\ConformalInvDeriv + \ImagUnit\ConformalInvDualDeriv\right)(a+\ImagUnit b).
    \end{equation*}
  \item For
    $f + \ImagUnit\LeftDual{f}\in \realHorkTensor{1}(\Complex)$, we
    define
    \begin{align*}
      \ConformalComplexDeriv\cdot(f+\ImagUnit\LeftDual{f})
      &\vcentcolon= \left(\ConformalInvDeriv + \ImagUnit\ConformalInvDualDeriv \right)\cdot (f+\ImagUnit\LeftDual{f}),\\
      \ConformalComplexDeriv\SymTracelessTensorProd(f+\ImagUnit\LeftDual{f})
      &\vcentcolon= \left(\ConformalInvDeriv + \ImagUnit\ConformalInvDualDeriv \right)\SymTracelessTensorProd (f+\ImagUnit\LeftDual{f}).
    \end{align*}
  \item For
    $u+\ImagUnit \LeftDual{u}\in \realHorkTensor{2}(\Complex)$, we
    define
    \begin{equation*}
      \ConformalComplexDeriv\cdot(u+\ImagUnit\LeftDual{u})
      \vcentcolon= \left(\ConformalInvDeriv + \ImagUnit\ConformalInvDualDeriv \right)\cdot (u+\ImagUnit\LeftDual{u}),\\
    \end{equation*}
  \item We use the notation
    \begin{equation*}
      \overline{\ConformalComplexDeriv} \vcentcolon= \ConformalInvDeriv - \ImagUnit\ConformalInvDualDeriv.
    \end{equation*}
  \end{enumerate}
\end{definition}


\section{The Kerr-de Sitter spacetime}
\label{sec:KdS}

In this section, we introduce the main features of the Kerr-de
Sitter spacetimes.

\subsection{The \KdS{} metric}
\label{sec:BL-coordinates}

We first introduce basic properties of the \KdS{} metric. 

\subsubsection{\KdS{} metric in Boyer-Lindquist coordinates}

In Boyer-Lindquist coordinates, the \KdS{} metric takes the form\footnote{
As the Boyer-Lindquist coordinates $(t,r,\theta,\phi)$ are regular
only for $\theta\neq 0,\pi$, we supplement them with $(t,r,x,y)$
coordinates where
$$x\vcentcolon =\sin\theta\cos\varphi, \qquad y\vcentcolon =\sin\theta\sin\varphi,$$
with $(t,r,x,y)$ regular for $\theta\neq\frac{\pi}{2}$.}
\begin{equation}
  \label{eq:KdS-metric:metric:BL}
  \begin{split}
  \Metric_{M,a,\Lambda} =& -\frac{\Delta}{\left(1+\gamma\right)^2\abs*{q}^2}(dt - a \sin^2\theta d\phi)^2
  + \abs*{q}^2 \left(
    \Delta^{-1}dr^2 +  \kappa^{-1}d\theta^2
  \right) + \frac{\kappa\sin^2\theta}{(1+\gamma)^2\abs*{q}^2}\left(a\,dt - (r^2+a^2)d\phi\right)^2,
 \end{split}
\end{equation}
and the inverse metric takes the form
\begin{equation}
  \label{eq:KdS-metric:inverse:BL}
  \abs*{q}^2\Metric_{M,a,\Lambda}^{-1} = 
    \Delta \,\partial_r^2
    +  \frac{(1+\gamma)^2}{\kappa \sin^2\theta}\left(a\sin^2\theta\,\partial_t + \partial_\phi\right)^2
    + \kappa\,\partial_\theta^2
    - \frac{(1+\gamma)^2}{\Delta}\left(\left(r^2+a^2\right)\partial_t + a\partial_\phi \right)^2
    ,
\end{equation}
where
\begin{align}
  \label{eq:metric-aux-qtys}
    q &\vcentcolon= r+ \ImagUnit a\cos\theta,& \Delta & \vcentcolon= (r^2+a^2)\left(1-\frac{\Lambda}{3}r^2\right)-2Mr,&
  \kappa &\vcentcolon= 1+ \gamma\cos^2\theta,&  \gamma & \vcentcolon= \frac{\Lambda a^2}{3}.
\end{align}
For convenience, we will also define
\begin{equation}
  \label{eq:Upsilon:def}
  \Upsilon \vcentcolon= \frac{\Delta}{r^2+a^2}.
\end{equation}
From the expression for the \KdS{} metric in
\zcref[noname]{eq:KdS-metric:metric:BL}, it is immediately clear that \KdS{}
admits two Killings vectors,
\begin{equation}
  \label{eq:KillT-KillPhi:def}
  \KillT\vcentcolon= \partial_t, \qquad
  \KillPhi\vcentcolon= \partial_{\phi}. 
\end{equation}
These span the Killing vectorfields of \KdS. In \SdS, the space of
Killing vectorfields is spanned by $\KillT$, and the rotational
vectorfields.

\subsubsection{Killing Horizons}
\label{sec:horizons}

\KdS{} features both a black hole region and a cosmological region,
which are separated from the stationary region, the main region of
interest, by Killing horizons.
\begin{definition}
  Subextremal \KdS{} black holes have two Killing horizons, the
  \emph{event horizon} and the \emph{cosmological horizon}, which in
  Boyer-Lindquist coordinates are $r$-constant null hypersurfaces
  $\{r = r_{\EventHorizon}\}$ and $\{r = r_{\CosmologicalHorizon}\}$
  respectively, where $r_{\EventHorizon}$ and
  $r_{\CosmologicalHorizon}$ are the largest positive roots of
  $\Delta$ (recall the definition of $\Delta$ in
  \zcref[noname]{eq:metric-aux-qtys}).
\end{definition}
In \SdS, the horizons are located at
\begin{equation*}
  r_{\EventHorizon} = \frac{2}{\sqrt{\Lambda}}\cos\left(\frac{1}{3}\cos^{-1}\left(3\sqrt{\Lambda}M\right)
    + \frac{\pi}{3}\right),\qquad
  r_{\CosmologicalHorizon} = \frac{2}{\sqrt{\Lambda}}\cos\left(\frac{1}{3}\cos^{-1}\left(3\sqrt{\Lambda}M\right)
    -  \frac{\pi}{3}\right).
\end{equation*}
We can Taylor expand these to be %
\begin{equation*}
  r_{\EventHorizon} = 2M + \frac{8M^3\Lambda}{3} + O(\Lambda^2),\qquad
  r_{\CosmologicalHorizon} = \sqrt{\frac{3}{\Lambda}} - M - \frac{\sqrt{3}}{2}M^2\sqrt{\Lambda} + O(\Lambda).   
\end{equation*}
Since we are only considering the slowly-rotating \KdS{} family, the
locations of the event and cosmological horizons in the \KdS{} black
holes we consider will simply be a small perturbation of the \SdS{}
locations.

\subsubsection{Inverse metric}
\label{sec:inverse-metric}

In this section, we discuss some properties of the \KdS{} inverse
metric which will be heavily used in what follows. We first introduce
the vectorfields
\begin{align}
  \HawkingVF\vcentcolon={}& (1+\gamma)\left(\KillT + \frac{a}{r^2+a^2}\KillPhi\right),    \label{eq:HawkingVF-def}\\
  \HprVF\vcentcolon={}& \frac{\Delta}{r^2+a^2}\partial_r, \label{eq:HprVF-def}
\end{align}
observing that $\HprVF$ is regular at the horizon, and that
$\HawkingVF$ is the \KdS{} analogue of the Hawking vectorfield on
Kerr. A simple computation leads to
\begin{align*}
    \Metric\pth{\HawkingVF, \HawkingVF}
      = -\frac{\De\abs*{q}^2}{\left(r^2+a^2\right)^2},
\end{align*}
showing that $\HawkingVF$ is thus timelike for
    $r_{\EventHorizon}<r<r_{\CosmologicalHorizon}$, and null
    on the horizons
    $r=r_{\EventHorizon}, r=r_{\CosmologicalHorizon}$.

Using the vectorfields $\HawkingVF$ and $\HprVF$ introduced in
\zcref[noname]{eq:HawkingVF-def} and \zcref[noname]{eq:HprVF-def} respectively, we can
rewrite the inverse \KdS{} metric in \zcref[noname]{eq:KdS-metric:inverse:BL}
in a simplified format.
\begin{lemma}
  \label{lemma:KdS-metric:inverse:THat-RHat}
  The inverse \KdS{} metric can be written in the form
  \begin{equation*}
    \abs*{q}^2\Metric^{\alpha\beta}
    = \Delta\partial_r^{\alpha}\partial_r^{\beta}
    + \frac{1}{\Delta}\mathcal{R}^{\alpha\beta},
  \end{equation*}
  with
  \begin{equation*}
    \mathcal{R}^{\alpha\beta} = -(r^2+a^2)^2\HawkingVF^{\alpha}\HawkingVF^{\beta}
    + \Delta \CartarOp^{\alpha\beta},
  \end{equation*}
  where
  \begin{equation}
    \label{eq:CartarOp:def}
    \abs*{q}^2\mathcal{O}^{\alpha\beta} = \kappa \partial_{\theta}^{\alpha}\partial_{\theta}^{\beta}
    + \frac{(1+\gamma)^2}{\kappa\sin^2\theta}\left(a\sin^2\theta\,\partial_t^{\alpha} + \partial_\phi^{\alpha}\right)\left(a\sin^2\theta\,\partial_t^{\beta} + \partial_\phi^{\beta}\right).
  \end{equation}
  As a result, inverse \KdS{} metric in \zcref[noname]{eq:KdS-metric:inverse:BL} can be expressed as
  \begin{equation}
    \label{eq:KdS-metric:inverse:THat-RHat}
    \abs*{q}^2\Metric^{\alpha\beta}
    = - \frac{\left(r^2+a^2\right)^2}{\Delta}\HawkingVF^{\alpha}\HawkingVF^{\beta}
    + \Delta \partial_r^{\alpha}\partial_r^{\beta}
    + \abs*{q}^2\CartarOp^{\alpha\beta}.
  \end{equation}  
\end{lemma}
\begin{proof}
  This is a straightforward computation.
\end{proof}

\subsection{Principal null frames in \KdS}
\label{sec:KdS:principal-null-frames}

In this section, we set up the principal outgoing and ingoing null
frames in \KdS{} that will be used throughout the paper.

\begin{definition}
  The following are the principal null directions in \KdS.
  \begin{enumerate}
  \item We denote the \emph{principal outgoing null pair} $(\eout_3, \eout_4)$, for which $e_4$ is geodesic, by 
    \begin{equation}
      \label{eq:e3-e4-outgoing:def}
      \begin{split}
        \eout_4 &= \frac{(1+\gamma)(r^2+a^2)}{\Delta}\partial_t
      + \partial_r
      + \frac{a(1+\gamma)}{\Delta} \partial_{\phi},\\
      \eout_3 &= \frac{(1+\gamma)(r^2+a^2)}{\abs*{q}^2}\partial_t
      - \frac{\Delta}{\abs*{q}^2}\partial_r
      + \frac{a(1+\gamma)}{\abs*{q}^2}\partial_{\phi}.
      \end{split}      
    \end{equation}
  \item We denote the \emph{principal ingoing null pair} $(\ein_3,\ein_4)$, for which $\ein_3$ is geodesic, by  
    \begin{equation}
      \label{eq:e3-e4-ingoing:def}
      \begin{split}
        \ein_3 &= \frac{(1+\gamma)(r^2+a^2)}{\Delta}\partial_t - \partial_r + \frac{a(1+\gamma)}{\Delta} \partial_{\phi},\\
        \ein_4 &= \frac{(1+\gamma)(r^2+a^2)}{\abs*{q}^2}\partial_t + \frac{\Delta}{\abs*{q}^2}\partial_r + \frac{a(1+\gamma)}{\abs*{q}^2}\partial_{\phi}.
      \end{split}      
    \end{equation}
  \item In both cases, a principal null frame can be obtained by adding
    the vectorfields
    \begin{equation}
      \label{eq:e1-e2:def}
      e_1 = \frac{\sqrt{\kappa}}{\abs*{q}}\partial_{\theta},\qquad
      e_2 = \frac{1+\gamma}{\sqrt{\kappa}}\left(\frac{a\sin\theta}{\abs*{q}}\partial_t + \frac{1}{\abs*{q}\sin\theta}\partial_\phi\right). 
    \end{equation}
  \end{enumerate}
\end{definition}
\begin{remark}
  Observe that in the case $\Lambda=0$, these choices reduce to the
  principal ingoing and outgoing null frames chosen in (3.3.1) and
  (3.3.12) of \cite{giorgiWaveEquationsEstimates2024} respectively.
\end{remark}

With $(e_1,e_2)$ as defined above, we have the following relationship
between the Killing vectorfields $\KillT, \KillPhi$ and the 
vectorfield $\HawkingVF$ defined in \zcref[noname]{eq:HawkingVF-def}. %
  \begin{align}
    \label{eq:KillT-KillPhi-HawkingVF:relationship}
      \KillT &= \frac{r^2+a^2}{(1+\gamma)\abs*{q}^2}\HawkingVF
               - \frac{a\sqrt{\kappa}\sin\theta}{\abs*{q}(1+\gamma)}e_2, &
      \KillPhi &= - \frac{a(r^2+a^2)\sin^2\theta}{(1+\gamma)\abs*{q}^2}\HawkingVF
                 +\frac{(r^2+a^2)\sqrt{\kappa}\sin\theta}{\abs*{q}(1+\gamma)}e_2.
  \end{align}

\subsubsection{Canonical complex one-form \texorpdfstring{$\CCOneFormJ$}{J} }
\label{sec:canonical-complex-one-form-J}

\begin{definition} 
  \label{def:canonical-complex-one-form-J}
  We define the following complex one-form $\CCOneFormJ$ in \KdS{} by $\CCOneFormJ\vcentcolon= j + \ImagUnit \LeftDual{j}$
  where the real one-form $j$ is defined by $j_1 = 0$ and $j_2 = \frac{\sqrt{\kappa}\sin\theta}{\abs*{q}}$. Hence,
  \begin{equation*}
    \CCOneFormJ_1 = \frac{\ImagUnit\sqrt{\kappa}\sin\theta}{\abs*{q}},\qquad
    \CCOneFormJ_2 = \frac{\sqrt{\kappa}\sin\theta}{\abs*{q}}. 
  \end{equation*}
\end{definition}

\begin{lemma}
  \label{lemma:canonical-complex-one-form-j-basic-identities}
  The following identities hold true.
  \begin{enumerate}
  \item We have
    \begin{equation*}
      \LeftDual{\CCOneFormJ} = -\ImagUnit \CCOneFormJ,\qquad
      \CCOneFormJ\cdot \overline{\CCOneFormJ} = \frac{2\kappa(\sin\theta)^2}{\abs*{q}^2}.
    \end{equation*}

  \item In both outgoing and ingoing null frames we have
    \begin{equation*}
      \nabla_4\CCOneFormJ + \frac{1}{2}\Trace X \CCOneFormJ = 0, \qquad
      \nabla_3 \CCOneFormJ + \frac{1}{2}\Trace \XBar \CCOneFormJ = 0. 
    \end{equation*}
  \item The complex one-form $\CCOneFormJ$ verifies
    \begin{equation*}
      \overline{\ComplexDeriv}\cdot \CCOneFormJ
      = \frac{4\ImagUnit(r^2+a^2)\cos\theta}{\abs*{q}^4} %
      + \frac{4}{3}\ImagUnit \Lambda\cos\theta\left(1 - \frac{r^2(r^2+a^2)}{\abs*{q}^4}\right)
      , \qquad
      \ComplexDeriv\SymTracelessTensorProd \CCOneFormJ = 0.
    \end{equation*}
  \item We have
    \begin{equation*}
      \ComplexDeriv(q) = -a\CCOneFormJ, \qquad
      \ComplexDeriv(\overline{q}) = a\CCOneFormJ. 
    \end{equation*}
  \end{enumerate}
\end{lemma}
\begin{proof}
  The proof follows by a direct computation. 
\end{proof}
\begin{remark}
  There are some differences in the algebra of the $\CCOneFormJ$
  between the Kerr case and the \KdS{} case. We still have $\Divergence j = 0$ and $\nabla\SymTracelessTensorProd j_{12} = 0$ but now 
  \begin{align*}
\nabla_aj_b & = \frac{\cos\theta}{\abs*{q}^4}\left(
                                             r^2+a^2 + \gamma\left(\abs*{q}^2\cos^2\theta - r^2\sin^2\theta\right)
                                             \right)\volFormHor_{ab}.
  \end{align*}
\end{remark}


\subsubsection{Properties of the outgoing principal null frame in \KdS}
\label{sec:KdS:principal-outgoing-null-frame}

In this section, let $(e_3,e_4) = (\eout_3, \eout_4)$ denote the principal outgoing null frame as defined in \zcref[noname]{eq:e3-e4-outgoing:def}, which is regular towards the future for all $r>r_{\EventHorizon}$.

\begin{lemma}
  \label{lemma:Kerr:outgoing-PG:Ric-and-curvature}
  The outgoing principal null frame enjoys the following properties.
  \begin{enumerate}
  \item The nonvanishing real Ricci coefficients and null curvature components are given by
    \begin{equation*}
    \begin{gathered}
      \Trace\chiBar = -\frac{2r \Delta}{\abs*{q}^{4}},\qquad
      \aTrace{\chiBar} = 2\frac{a \Delta \cos\theta}{\abs*{q}^{4}},\qquad
      \Trace\chi = \frac{2r }{\abs*{q}^{2}},\qquad
      \aTrace{\chi} = 2\frac{a  \cos\theta}{\abs*{q}^{2}},\\    
      (\eta_1, \eta_2) = \left( -\frac{a^{2}\sqrt{\kappa}\cos\theta\sin\theta}{\abs*{q}^3}, \frac{a r \sqrt{\kappa}\sin\theta}{\abs*{q}^3} \right),\qquad
      (\etaBar_1, \etaBar_2) = \left( -\frac{a^{2}\sqrt{\kappa}\cos\theta\sin\theta}{\abs*{q}^3}, -\frac{a r \sqrt{\kappa}\sin\theta}{\abs*{q}^3} \right),\\
      (\zeta_1, \zeta_2) =  \left( \frac{a^{2}\sqrt{\kappa}\cos\theta\sin\theta}{\abs*{q}^3}, \frac{a r \sqrt{\kappa}\sin\theta}{\abs*{q}^3} \right),\\
      \omegaBar = \frac{1}{2}\partial_r\left(\frac{\Delta}{\abs*{q}^2}\right),\qquad
      \rho = -2\frac{M r \left(r^2 - 3 a^{2}\cos^2\theta\right)}{\abs*{q}^{6}}, \qquad
      \LeftDual{\rho} = -2aM\cos\theta\frac{a^2\cos^2\theta - 3r^2} {\abs*{q}^6}.
    \end{gathered}
  \end{equation*}

  \item The complex null curvature components with respect to the frame are given by
    \begin{equation*}
      A = B = \BBar = \ABar=0, \qquad P = -\frac{2M}{q^3}.
    \end{equation*}
  \item The vanishing complex Ricci coefficients in \KdS{} are
    \begin{equation*}
      \XHat = \XHatBar = \Xi = \XiBar = \omega = 0. 
    \end{equation*}
  \item The non-vanishing complex Ricci coefficients are
    \begin{gather*}
      \Trace X = \frac{2}{q}\qquad
      \Trace \XBar = -\frac{2\Delta q}{\abs*{q}^4},\qquad  \omegaBar = \frac{1}{2}\partial_{r}\left(\frac{\Delta}{\abs*{q}^{2}}\right),\qquad
      \HBar = -\frac{a\overline{q}}{\abs*{q}^2}\CCOneFormJ,\qquad
      H = \frac{a q}{\abs*{q}^{2}}\CCOneFormJ,\qquad
      Z = \frac{a\overline{q}}{\abs*{q}^{2}}\CCOneFormJ,
    \end{gather*}
    with the one-form $\CCOneFormJ$ as defined in
    \zcref[cap]{def:canonical-complex-one-form-J}.

  \item Recalling that $(\Lambda_{\alpha})_{\beta\gamma} = \Metric\left(\CovariantDeriv_{\alpha}e_{\gamma},e_{\beta}\right)$ we have 
    \begin{equation}
      \label{eq:Kerr:outgoing:Lambda-Ricci}
      \begin{split}
        &(\Lambda_3)_{12} = - \frac{a\Delta\cos\theta}{\abs*{q}^4},\qquad
      (\Lambda_4)_{12} = - \frac{a\cos\theta}{\abs*{q}^2}, \qquad (\Lambda_1)_{12} = 0, \\
      & (\Lambda_2)_{12} = - \frac{(r^2+a^2+\gamma(\abs*{q}^2\cos^2\theta - r^2\sin^2\theta))\cot\theta}{\abs*{q}^3\sqrt{\kappa}}.
      \end{split}      
    \end{equation}
  \end{enumerate}
\end{lemma}

\begin{lemma}
  \label{lemma:T-R-Z:principal-outgoing:expression}
  The following identities hold true in the principal outgoing frame.
  \begin{equation}
    \label{eq:T-R-Z:principal-outgoing:expression:e4-e3}
    e_4 = \frac{(1+\gamma)(r^2+a^2)}{\Delta}\left(\HawkingVF + \frac{1}{1+\gamma}\HprVF\right),\qquad
    e_3 = \frac{r^2+a^2}{\abs*{q}^2}\left(\HawkingVF  - \HprVF\right).
  \end{equation}
  In addition, 
  \begin{equation}
    \label{eq:T-R-Z:principal-outgoing:expression:THat-RHat}
    \begin{split}
      \HawkingVF
      =&{} \frac{1}{2(1+\gamma)}\left(\frac{\Delta}{r^2+a^2}e_4^{\operatorname{(out)}}
                    + \frac{\abs*{q}^2}{r^2+a^2}e_3^{\operatorname{(out)}}\right),\qquad
      \HprVF ={} \frac{1}{2(1+\gamma)}\left(\frac{\Delta}{r^2+a^2}e_4^{\operatorname{(out)}}
                    - \frac{\abs*{q}^2}{r^2+a^2}e_3^{\operatorname{(out)}}\right),
    \end{split}
  \end{equation}
and
  \begin{equation}
    \label{eq:T-R-Z:principal-outgoing:expression:T-Z}
    \begin{split}
      \KillT ={}& \frac{1}{2(1+\gamma)}\left(\frac{\Delta}{\abs*{q}^2}e_4^{\operatorname{(out)}} + e_3^{\operatorname{(out)}}\right)
      - \frac{a\sqrt{\kappa}\sin\theta}{\abs*{q}(1+\gamma)}e_2,\\
      \KillPhi ={}& \frac{\left(r^2+a^2\right)\sqrt{\kappa}\sin\theta}{(1+\gamma)\abs*{q}}e_2
                  - \frac{1}{2(1+\gamma)}a\sin^2\theta\left(\frac{\Delta}{\abs*{q}^2}e_4^{\operatorname{(out)}} +  e_3^{\operatorname{(out)}}\right).
    \end{split}
  \end{equation}
\end{lemma}

\begin{remark}
We remark that in \KdS, we have the following asymptotic behavior for
$\horProj{K}$ as defined in \zcref[noname]{eq:h-K:def}
\begin{equation}
  \label{eq:K:asymptotic-behavior}
  \horProj{K} = r^{-2} + O(ar^{-3}).
\end{equation}
\end{remark}

\subsubsection{Properties of the ingoing principal null frame on \KdS}
\label{sec:KdS:principal-ingoing-null-frame}

In this section, let $(e_3,e_4)=(\ein_3,\ein_4)$ denote the principal ingoing null
vectors as defined in \zcref[noname]{eq:e3-e4-ingoing:def} which is regular
towards the future for all $r<r_{\CosmologicalHorizon}$.

\begin{lemma}
  \label{lemma:Kerr:ingoing-PG:Ric-and-curvature}
  The ingoing principal null frame enjoys the following properties.
  \begin{enumerate}
  \item The nonvanishing real Ricci coefficients and null curvature components are given by
    \begin{equation*}
      \begin{gathered}
        \Trace\chi = \frac{2r \Delta}{\abs*{q}^4},\qquad
        \aTrace{\chi} = 2\frac{a \Delta \cos\theta}{\abs*{q}^4},\qquad
        \Trace\chiBar = -\frac{2r }{\abs*{q}^2},\qquad
        \aTrace{\chiBar} = 2\frac{a  \cos\theta}{\abs*{q}^2},\\    
         (\eta_1, \eta_2) = \left( -\frac{a^2\sqrt{\kappa}\cos\theta\sin\theta}{\abs*{q}^3}, \frac{a r \sqrt{\kappa}\sin\theta}{\abs*{q}^3} \right),\qquad
         (\etaBar_1, \etaBar_2) = \left( -\frac{a^2\sqrt{\kappa}\cos\theta\sin\theta}{\abs*{q}^3}, -\frac{a r \sqrt{\kappa}\sin\theta}{\abs*{q}^3} \right),\\
        (\zeta_1,\zeta_2) =  \left(- \frac{a^2\sqrt{\kappa}\cos\theta\sin\theta}{\abs*{q}^3}, \frac{a r \sqrt{\kappa}\sin\theta}{\abs*{q}^3} \right),\\
        \omega = -\frac{1}{2}\partial_r\left(\frac{\Delta}{\abs*{q}^2}\right),\qquad
        \rho = -2\frac{M r \left(r^2 - 3 a^2\cos^2\theta\right)}{\abs*{q}^6}, \qquad
        \LeftDual{\rho} = -2aM\cos\theta\frac{a^2\cos^2\theta - 3r^2} {\abs*{q}^6}.
      \end{gathered}
    \end{equation*}
  \item The complex null curvature  components with respect to the frame are given by
    \begin{equation*}
      A = B = \BBar = \ABar=0, \qquad P = -\frac{2M}{q^3}.
    \end{equation*}
  \item The vanishing complex Ricci coefficients in \KdS{} are
    \begin{equation*}
      \XHat = \XHatBar = \Xi = \XiBar = \omegaBar = 0. 
    \end{equation*}
  \item The non-vanishing complex Ricci coefficients are
    \begin{gather*}
      \Trace X = \frac{2\Delta\overline{q}}{\abs*{q}^{4}}\qquad
      \Trace \XBar = -\frac{2}{\overline{q}},\qquad  \omega = -\frac{1}{2}\partial_{r}\left(\frac{\Delta}{\abs*{q}^{2}}\right),\\
      \HBar = -\frac{a\overline{q}}{\abs*{q}^2}\CCOneFormJ,\qquad
      H = \frac{a q}{\abs*{q}^{2}}\CCOneFormJ,\qquad
      Z = \frac{a q}{\abs*{q}^{2}}\CCOneFormJ,
    \end{gather*}
    with the one-form $\CCOneFormJ$ as defined in
    \zcref[cap]{def:canonical-complex-one-form-J}.

  \item Recalling that $(\Lambda_{\alpha})_{\beta\gamma} = \Metric\left(\CovariantDeriv_{\alpha}e_{\gamma},e_{\beta}\right)$ we have 
    \begin{equation}
      \label{eq:Kerr:ingoing:Lambda-Ricci}
      \begin{split}
        & (\Lambda_3)_{12} =  \frac{a\cos\theta}{\abs*{q}^2},\qquad
           (\Lambda_4)_{12} =  \frac{a\Delta\cos\theta}{\abs*{q}^4}, \qquad (\Lambda_1)_{12} = 0,\\
        & (\Lambda_2)_{12} = - \frac{(r^2+a^2+\gamma(\abs*{q}^2\cos^2\theta - r^2\sin^2\theta))\cot\theta}{\abs*{q}^3\sqrt{\kappa}}.
      \end{split}      
    \end{equation}
  \end{enumerate}
\end{lemma}

\begin{lemma}
  \label{lemma:T-R-Z:principal-ingoing:expression}
  The following identities hold true in the principal ingoing frame.
  \begin{equation}
    \label{eq:T-R-Z:principal-ingoing:expression:e4-e3}
    e_4^{(\operatorname{in})} = \frac{r^2+a^2}{\abs*{q}^2}\left(\HawkingVF + \frac{1}{1+\gamma}\HprVF\right),\qquad
    e_3^{(\operatorname{in})} = \frac{(1+\gamma)(r^2+a^2)}{\Delta}\left(\HawkingVF  - \HprVF\right).
  \end{equation}
  In addition, 
  \begin{align}
    \label{eq:T-R-Z:principal-ingoing:expression:THat-RHat}
      \HawkingVF
      &= \frac{1}{2(1+\gamma)}\left(
         \frac{\abs*{q}^2}{r^2+a^2}e_4^{\operatorname{(in)}}
         + \frac{\Delta}{r^2+a^2}e_3^{\operatorname{(in)}}
         \right),&
      \HprVF & =  \frac{1}{2(1+\gamma)}\left(
                \frac{\abs*{q}^2}{r^2+a^2}e_4^{\operatorname{(in)}}
                - \frac{\Delta}{r^2+a^2}e_3^{\operatorname{(in)}}
                \right),
\end{align}
  and
  \begin{equation}
    \label{eq:T-R-Z:principal-ingoing:expression:T-Z}
    \begin{split}
      \KillT ={}& \frac{1}{2(1+\gamma)}\left(\frac{\Delta}{\abs*{q}^2}e_3^{\operatorname{(in)}} + e_4^{\operatorname{(in)}}\right)
                  - \frac{a\sqrt{\kappa}\sin\theta}{\abs*{q}(1+\gamma)}e_2,\\
      \KillPhi ={}& \frac{\left(r^2+a^2\right)\sqrt{\kappa}\sin\theta}{(1+\gamma)\abs*{q}}e_2
                    - \frac{1}{2(1+\gamma)}a\sin^2\theta\left(e_4^{\operatorname{(in)}} +  \frac{\Delta}{\abs*{q}^2}e_3^{\operatorname{(in)}}\right).
    \end{split}
  \end{equation}
\end{lemma}

\subsubsection{Global principal null frame}
\label{sec:global null frame}

In \KdS, neither the principal outgoing nor the principal ingoing null
frame are defined globally. To circumvent this, we introduce the
following principal global frame. It will be used in particular to write the Teukolsky wave-transport
system globally (see \zcref[cap]{sec:derivation-of-RW}).

\begin{definition}
  \label{def:global null frame}
  Choose $r_{\EventHorizon}<r_0<r_{\CosmologicalHorizon}$ large compared to $M$, and fix $\chiglo=\chiglo(r)\in C^{\infty}(\mathcal{M}; [0,1])$ a cutoff function such that
  \begin{equation*}
    \chiglo(r) = 0 \text{ for }r\leq r_0- M,\qquad
    \chiglo(r) = 1 \text{ for }r \geq r_0 + M
    .
  \end{equation*}
  We define $\left( \eglo_4, \eglo_3 \right)$, the \emph{principal global null
    frame of \KdS}, by
  \begin{equation*}
    \left(\eglo_4,\eglo_3\right) \vcentcolon =\left(\lambdaglo(r, \theta) \eout_4 ,\lambdaglo(r, \theta)^{-1}\eout_3\right),
  \end{equation*}  
  where $\lambdaglo$ is defined by
  \begin{equation}
    \label{eq:lambdaglo:def}
    \lambdaglo\vcentcolon= (1-\chiglo(r))\frac{\Delta}{\abs*{q}^2} + \chiglo(r).
  \end{equation}
  In particular $\left(\eglo_4,\eglo_3\right) = \left(\ein_4, \ein_3\right)$ for $r\le r_0 - M$ and $\left(\eglo_4,\eglo_3\right) = \left(\eout_4, \eout_3\right)$ for $r\ge r_0 - M$.
\end{definition}

\begin{remark}
  The following formulas hold for the principal global null frame
  \begin{gather}
    \label{eq:global-frame:vanishing-qtys}
    A = \underline{A}=0,\quad
    B = \BBar = 0, \quad 
    \XHat = \XHatBar =0,\quad \Xi = \XiBar  = 0,\\
    \label{eq:global-frame:non-vanishing-qtys}
    P = -\frac{2M}{q^3},\quad
    \HBar = -\frac{a\overline{q}}{\abs*{q}^2}\CCOneFormJ,\quad
    H = \frac{a q}{\abs*{q}^{2}}\CCOneFormJ
    .
  \end{gather}
  Indeed \zcref[noname]{eq:global-frame:vanishing-qtys} and
  \zcref[noname]{eq:global-frame:non-vanishing-qtys} follow from the values of
  the Ricci coefficients and Weyl components on \KdS{} in 
  \zcref[cap]{lemma:Kerr:outgoing-PG:Ric-and-curvature} and the conformal
  invariance of the above quantities in the sense of
  \zcref[cap]{def:s-conformally-invariants}.
\end{remark}

From now on, we will always use
$(e_1,e_2,e_3,e_4) = \left( e_1,e_2,\eglo_3, \eglo_4 \right)$. Note
that \zcref[cap]{lemma:T-R-Z:principal-outgoing:expression} implies
the following:
\begin{equation}
  \label{eq:T-R-Z:principal-global:expression:e4-e3}
  \lambdaglo^{-1}e_4 = \frac{(1+\gamma)(r^2+a^2)}{\Delta}\left(\HawkingVF + \frac{1}{1+\gamma}\HprVF\right),\qquad
  \lambdaglo e_3 = \frac{r^2+a^2}{\abs*{q}^2}\left(\HawkingVF  - \HprVF\right),
\end{equation}
\begin{align}
  \label{eq:T-R-Z:principal-global:expression:THat-RHat}
    \HawkingVF &
    = \frac{1}{2(1+\gamma)}\left(\frac{\Delta}{r^2+a^2}\lambdaglo^{-1}e_4
       + \frac{\abs*{q}^2}{r^2+a^2}\lambdaglo e_3\right),&
    \HprVF & = \frac{1}{2(1+\gamma)}\left(\frac{\Delta}{r^2+a^2}\lambdaglo^{-1}e_4
              - \frac{\abs*{q}^2}{r^2+a^2}\lambdaglo e_3\right),
\end{align}
and
\begin{equation}
  \label{eq:T-R-Z:principal-global:expression:T-Z}
  \begin{split}
    \KillT ={}& \frac{1}{2(1+\gamma)}\left(\frac{\Delta}{\abs*{q}^2}\lambdaglo^{-1}e_4 + \lambdaglo e_3\right)
                - \frac{a\sqrt{\kappa}\sin\theta}{\abs*{q}(1+\gamma)}e_2,\\
    \KillPhi ={}& \frac{\left(r^2+a^2\right)\sqrt{\kappa}\sin\theta}{(1+\gamma)\abs*{q}}e_2
                  - \frac{1}{2(1+\gamma)}a\sin^2\theta\left(\frac{\Delta}{\abs*{q}^2}\lambdaglo^{-1} e_4 +  \lambdaglo e_3\right).
  \end{split}
\end{equation}

\subsection{Global coordinates}
\label{sec:adaptedglobalcoordinates}

The Boyer-Lindquist coordinates of \zcref[cap]{sec:BL-coordinates} are
only regular between the two horizons of the Kerr-de Sitter metric. To
circumvent this, we introduce new global coordinates that will be used
throughout the article.

\begin{definition}\label{definition de k et de h}
  Consider $\chiglo$ the cutoff introduced in \zcref[cap]{def:global
    null frame}. We define $k$ and $h$ to be any function satisfying
  \begin{align*}
    k'(r) &  =(1+\ga) \pth{1-2\chi_{\mathrm{glo}}(r)}\pth{ \frac{r^2+a^2}{\De} - \frac{M^2}{r^2}},&
     h'(r) & = (1-2\chi_{\mathrm{glo}}(r))\frac{a(1+\ga)}{\De},
  \end{align*}
  and set $\tau  = t + k(r)$ and $\ffi  = \phi + h(r)$ where $t$ and $\phi$ are the Boyer-Lindquist $(t, \phi)$ coordinates.
\end{definition}

\begin{lemma}\label{lem:action-global-frame}
  The action of the global principal null frame of
  \zcref[cap]{def:global null frame} on the $(\tau,r,\th,\ffi)$
  coordinates is given by:
  \begin{align*}
    e_4(r) &  = \lambdaglo, & e_4(\tau) & =\lambdaglo\left(\frac{(1+\gamma)(r^2+a^2)}{\Delta} + \tauAux'(r)\right),
    \\ e_4(\theta)& =0, &  e_4(\varphi) & =\lambdaglo\left(\frac{a(1+\gamma)}{\Delta} + h'(r)\right),
\\    e_3(r) & = -\frac{\Delta}{\lambdaglo\abs*{q}^2} ,& e_3(\tau) & =  \lambdaglo^{-1}\left(\frac{(1+\gamma)(r^2+a^2)}{\abs*{q}^2} - \frac{\Delta}{\abs*{q}^2}\tauAux'(r)\right),
    \\ e_3(\theta)& =0 ,& e_3(\varphi) & = \lambdaglo^{-1}\left( \frac{a(1+\gamma)}{\abs*{q}^2} - \frac{\Delta}{\abs*{q}^2}h'(r)\right),
    \end{align*}
\begin{align*}
    e_1(r) & = 0, & e_1(\tau) &= 0, &  e_1(\theta) & = \frac{\sqrt{\kappa}}{\abs*{q}}, & e_1(\varphi) & = 0,
\\    e_2(r) & = 0, & e_2(\tau) & = \frac{a(1+\gamma)\sin\theta}{\abs*{q}\sqrt{\kappa}}, &  e_2(\theta) & = 0, & e_2(\varphi) &= \frac{1+\gamma}{\abs*{q}\sqrt{\kappa}\sin\theta}.
  \end{align*}
\end{lemma}

\begin{proof}
    Straightforward computations using the definition of global principal null frame, the action of the principal outgoing null frame on the Boyer-Lindquist coordinates (which can be deduced from \zcref[cap]{eq:e3-e4-outgoing:def}) and the definition of $(\tau,\ffi)$.
\end{proof}

\begin{remark}
    Note already that \zcref[cap]{lem:action-global-frame} and \zcref[noname]{eq:T-R-Z:principal-global:expression:T-Z} imply
    \begin{equation}
      \label{eq:tau-foliation:properties:normalization}
      \KillT(\tau) = 1,\qquad
      \nabla(\tau) = a\Re \CCOneFormJ.
    \end{equation}
\end{remark}

\begin{proposition}\label{prop:properties-of-global-coordinates}
  If $a$, $\de_\HH$ and $\La$ are small enough compared to $M$, then the following holds on $r_\HH(1-\de_\HH) \leq r \leq r_{\overline{\HH}}(1+\de_\HH)$:
  \begin{align}\label{eq:tau-foliation:properties:NSigma-uniformly-timelike}
    \g(\D\tau,\D\tau) & \leq -  \frac{M^2}{8r^2},
  \end{align}
  and
  \begin{align}\label{eq:tau-foliation:properties:ingoing}
    e_4(\tau) & > 0,& e_3(\tau) & > 0,& | \nab \tau|^2 \leq \frac{8}{9}e_4(\tau) e_3(\tau).
  \end{align}
  Moreover, if $r\geq r_0+M$ then
  \begin{align}\label{eq:tau-foliation:properties:asymptotic-behavior}
    &\frac{M^2}{r^2} \lesssim e_4(\tau) \lesssim \frac{M^2}{r^2} ,& 1 \lesssim e_3(\tau) \lesssim 1.
  \end{align}
  Finally, the Kerr-de Sitter metric is regular in the coordinates $(\tau,r,\th,\ffi)$.
\end{proposition}

\begin{proof}
 See \zcref[cap]{appendix:properties-of-global-coordinates}. 
\end{proof}

We will thus use $(\tau, r, \theta, \varphi)$ as a global coordinate system\footnote{The $(\tau, r, \theta, \varphi)$ coordinate system is global up to needing different coordinates to cover the north and south poles of the spheres.}. \zcref[cap]{prop:properties-of-global-coordinates} shows that the constant $\tau$ hypersurfaces, denoted $\Si(\tau)$, are spacelike hypersurfaces but that they are also asymptotically null. To include an external statement, we also need a uniformly timelike coordinate.

\begin{definition}\label{definition de ell}
    Consider $\underline{\chi}$ a smooth cutoff function such that
    \begin{align*}
        \underline{\chi}(r) = 0 \text{ for $r\leq 2r_0$}, \qquad \underline{\chi}(r) = 1 \text{ for $r\geq 3r_0$},
    \end{align*}
    and $0<\underline{\chi}<1$ on $(2r_0,3r_0)$. Define $\ell(r)  \vcentcolon = \int_{2r_0}^r \underline{\chi} $ and $\underline{t}  \vcentcolon  = \tau + \ell(r)$.
\end{definition}

The next proposition gathers properties of $\underline{t}$.

\begin{proposition}\label{prop:properties-underline-t}
  If $r_0$ is large enough compared to $M$, the following holds. First, we have the following asymptotics
  \begin{align}\label{eq:action-frame-underline-t}
      e_4(\underline{t}) & =   1 +  \GO{\frac{M^2}{r^2}}  ,& e_3(\underline{t}) & =  2  - \frac{\De}{|q|^2} +  \GO{\frac{M^2}{r^2}} ,& |\nab\underline{t}| & = \GO{\frac{a}{r}}, 
  \end{align}
  and
  \begin{align}\label{eq:asympt-gtt}
      \g(\D\underline{t},\D\underline{t})   =  -2 + \frac{\De}{|q|^2} + \GO{\frac{M^2}{r^2}}.
  \end{align}
Moreover, we have
\begin{align}\label{eq:underline-t-timelike}
    \g(\D\underline{t},\D\underline{t}) & <0,
\end{align}
on $r_\HH(1-\de_\HH) \leq r \leq r_{\overline{\HH}}(1+\de_\HH)$ and the coordinate $\underline{t}$ is regular on this region.
\end{proposition}

\begin{proof}
See \zcref[cap]{appendix:proof lem underline t}.
\end{proof}

\subsection{Spacetime regions}\label{sec:spacetime-regions}

Based on the coordinates $(\tau,r,\th,\ffi)$ and $\underline{t}$
defined in \zcref[cap]{sec:adaptedglobalcoordinates}, we now define
the main manifold considered in this article, along with some
important regions and hypersurfaces.

\begin{definition}\label{def:main-regions}
  We consider the following manifold and regions of it:
  \begin{align*}
    \MM_{\mathrm{tot}} & \vcentcolon = \pth{\mathbb{R}_\tau \times \left(
                     r_{\EventHorizon}(1-\delta_{\Horizon}),
                     r_{\CosmologicalHorizon}(1+\delta_{\Horizon}) \right)_r\times
                     \Sphere^2}\cap \{\underline{t}\geq 0\},
    \\    \MM & \vcentcolon = \MM_{\mathrm{tot}} \cap \{\tau \geq 0\},
    \\ \MM_e & \vcentcolon = \MM_{\mathrm{tot}} \cap \{\tau < 0\},
  \end{align*}
  endowed with the \KdS{} metric. For $\tau'\in\mathbb{R}$, consider
  the following hypersurface
  \begin{align*}
    \Si(\tau') & \vcentcolon = \MM_{\mathrm{tot}} \cap \{\tau=\tau'\} .
  \end{align*}
  Finally, consider the initial hypersurface 
  \begin{equation}
    \label{eq:true-init-hypersurface:def}
    \widehat{\Si}_\init  \vcentcolon = \MM_{\mathrm{tot}} \cap \{ \underline{t} = 0\}.
  \end{equation}
\end{definition}

\begin{remark}
  Note that the region $\MM_e$ is an external region in the sense that on $\MM_e$ we have
  \begin{align}\label{eq:external-properties}
    2r_0 & <r \leq r_{\overline{\HH}}(1+\de_\HH) , & \tau_{\overline{\HH}}&\leq \tau < 0 ,& |\tau|& \leq r, 
  \end{align}
  where \[\tau_{\overline{\HH}}\vcentcolon = -r_{\overline{\HH}}(1+\de_\HH)+3r_0 - \int_{2r_0}^{3r_0}\underline{\chi}.\] Moreover, the family $(\Si(\tau))_{\tau\geq 0}$ is a foliation by spacelike hypersurfaces of $\MM$, while \[(\Si(\tau))_{-r_{\overline{\HH}}(1+\de_\HH) + c_0 \leq \tau < 0},\] is a foliation by spacelike hypersurfaces of $\MM_e$.
\end{remark}

In the following, we define important regions of $\MM$ in view of the properties of  \KdS.

\begin{definition}
  Let $0<\delta_{\trap}, \delta_{\operatorname{red}}\ll1$
  be sufficiently small\footnote{We will require that
    \begin{equation*}
      \frac{\abs*{a}}{M} \ll \delta_{\operatorname{red}}, \delta_{\operatorname{trap}}\ll 1.
    \end{equation*}} to be fixed later. 
  \begin{enumerate}
  \item We define
    $\Manifold_{\operatorname{trap}}$
    the \emph{trapping region} of $\Manifold$ to be the set 
    
    \begin{equation}
      \label{eq:M-trap:def}
      \Manifold_{\operatorname{trap}}
      \vcentcolon= \Manifold\cap\curlyBrace*{\left|1-\frac{3M}{r}\right| \le \delta_{\operatorname{trap}}}, 
    \end{equation} 
    and    $\Manifold_{\cancel{\operatorname{trap}}}$
    the \emph{non-trapping region} of $\Manifold$ to be the complement of the
    trapping region
    $\Manifold_{\operatorname{trap}}$.
  \item We define the black hole \emph{redshift region}
    $\Manifold_{\operatorname{red}}$ to be the region
    \begin{equation}
      \label{eq:M-redshift:event-horizon:def}
      \Manifold_{\operatorname{red}} \vcentcolon= \Manifold \cap \curlyBrace*{r\le r_{\HH}\left(1+\delta_{\operatorname{red}}\right)},
    \end{equation}
    such that the connected component of the ergoregion containing the
    event horizon is contained in $\Manifold_{\operatorname{red}}$.
  \end{enumerate}
\end{definition}

\begin{figure}[ht]
  \centering
  \input{Images/Global-foliation.tex}
  \caption{\textit{Penrose diagram of $\MM_{\mathrm{tot}}$: $\MM$ is in gray while $\MM_e$ is in brown.}}
  \label{fig:global-foliation}
\end{figure}

Finally, we introduce the boundaries of the various spacetime regions considered in this article.

\begin{definition}
  The boundaries of $\MM_{\mathrm{tot}}$, $\MM$ and $\MM_e$ are given by
  \begin{align*}
    \dr\MM_{\mathrm{tot}} & = \mathcal{A} \cup \widehat{\Si}_\init \cup \overline{\Si}_* ,
    \\ \dr\MM & = \mathcal{A} \cup \Si(0) \cup \SigmaStar,
    \\ \dr \MM_e & = \pth{ \Si(0)\cap\{r\geq 2r_0\} } \cup  \pth{ \widehat{\Si}_\init\cap\{r\geq 2r_0\} }  \cup \Si_{*,e},
  \end{align*}
  where 
  \begin{align*}
    \mathcal{A} & \vcentcolon= \curlyBrace{r=r_{\EventHorizon}(1-\delta_{\Horizon})}\cap \MM_{\mathrm{tot}},
    \\ \overline{\Si}_* &\vcentcolon=
                          \curlyBrace{r=r_{\CosmologicalHorizon}(1+\delta_{\Horizon})}\cap\MM_{\mathrm{tot}}, 
    \\\SigmaStar & \vcentcolon = \overline{\Si}_* \cap\{\tau\geq 0\},
    \\ \Si_{*,e} & \vcentcolon = \overline{\Si}_* \cap\{\tau\leq 0\}.
  \end{align*}
  See \zcref{fig:global-foliation} for an overview of the regions of $\MM_{\mathrm{tot}}$.
\end{definition}

\section{The Teukolsky wave-transport system in \KdS}
\label{sec:derivation-of-RW}

In this section, we derive the main system studied in this article, namely the Teukolsky wave-transport system in \KdS. We start by the general derivation of the Teukolsky equation, first derived in \cite{khanalRotatingBlackHole1983} in the Newman-Penrose formalism using mode decomposition.

\begin{proposition}
  \label{prop:Teuk:Teuk-eq-for-A}
  Consider a solution of the Einstein vacuum equations with cosmological constant $\La$
  and a null frame. Then the $2$-conformally invariant complexified
  Weyl component $A\in \realHorkTensor{2}(\Complex)$ introduced in
  \zcref[noname]{eq:complex-curvature:def} satisfies
  \begin{align}\label{eq:teuk-with-Q}
    \mathcal{L}(A) & = \mathcal{Q}\pth{\Xi,\underline{\Xi},\widehat{X},\widehat{\underline{X}},A,\underline{A},B,\underline{B}},
  \end{align}
  where
  \begin{equation}
    \label{eq:Teuk:A:def}
    \begin{split}
      \TeukOp(A)\vcentcolon={}&
                                - \ConformalInvDeriv_{4}\ConformalInvDeriv_{3} A
                                + \frac{1}{4}\ConformalComplexDeriv\SymTracelessTensorProd\left(\overline{\ConformalComplexDeriv}\cdot A\right)\\
                              &  + \left(-\frac{1}{2}\Trace X - 2\overline{\Trace X}\right)\ConformalInvDeriv_{3}A  - \frac{1}{2}\Trace\XBar\ConformalInvDeriv_{4}A
                                + \left(4H + \HBar + \overline{\HBar}\right)\cdot \ConformalInvDeriv A\\
                              & + \left(-\overline{\Trace X}\Trace \XBar + 2\overline{P} - \frac{2\Lambda}{3} +  2 H\cdot \overline{\underline{H}} \right)A,
    \end{split}
  \end{equation}
  and where $\mathcal{Q}$ is a quadratic expression in its variables and their first derivatives.
\end{proposition}

\begin{proof}
  See proof in \zcref[cap]{appendix:prop:Teuk:Teuk-eq-for-A}.
\end{proof}

In the context of the nonlinear stability of Kerr-de Sitter, i.e when
the metric and the null frame are perturbations of the metric and the
principal null frame of Kerr-de Sitter, the right hand side
$\mathcal{Q}\pth{\Xi,\underline{\Xi},\widehat{X},\widehat{\underline{X}},A,\underline{A},B,\underline{B}}$
of \zcref[noname]{eq:teuk-with-Q} vanishes quadratically. Therefore
the linearisation of \zcref[noname]{eq:teuk-with-Q} is given by the
Teukolsky equation
\begin{align}\label{eq:teuk}
  \mathcal{L}(A) & = 0,
\end{align}
where $\mathcal{L}$ is as defined in
\zcref[noname]{eq:Teuk:A:def}. From now on, we will consider generic
solutions to \zcref[noname]{eq:teuk}. To derive estimates for such
solutions, it will be useful to derive the following generalized
Regge-Wheeler equation.

\begin{proposition}
  \label{proposition T to RW}
  Let $A\in \mathfrak{s}_2(\mathbb{C})$ and assume that $A$ is
  $2$-conformally invariant with respect to the global principal null pair in Kerr-de Sitter and moreover is solving \zcref[noname]{eq:teuk}. The tensor $\mathfrak{q}\in \realHorkTensor{2}(\Complex)$ defined by
  \begin{align}\label{definition of q frak}
    \q & \vcentcolon= \frac{q}{\bar{q}} r^2 \,^{(c)}\nab_3 \,^{(c)} \nab_3 \pth{ \frac{\bar{q}^4}{r^2}A} ,
  \end{align}
  is $0$-conformally invariant and satisfies the equation
  \begin{align}\label{gRW}
    \dot{\Box}_2\q -  \ImagUnit\frac{4a(1+\ga)\cos\th}{|q|^2} \nab_{\T}\q -  \widetilde{V}  \q ={}   q\bar{q}^3 \bigg(& \,^{(c)}\nab_\VFArthur\,^{(c)}\nab_3 A + W^{[3]} \,^{(c)}\nab_3 A 
    \\ & + W^{[4]} \,^{(c)}\nab_4 A  + W^{[h]} \cdot \,^{(c)}\nab A  + W^{[0]} A \bigg),\nonumber
  \end{align}
  where the potential $\widetilde{V}$ is real and satisfies
  \begin{align*}
    \widetilde{V}   =  V + V_0,\qquad
    V =  \frac{4\De}{(r^2+a^2)|q|^2} + 2\La,\qquad V_0 = \GO{ a r^{-4}},
  \end{align*}
  and the vector field $\VFArthur$ is given by
  \begin{equation*}
    \VFArthur  = - \lambdaglo^{-1} (1+\ga)   \frac{8a\De }{r^2|q|^4}  \pth{   a  \T   +  \mathbf{\Phi} },
  \end{equation*}
  with $\lambdaglo$ is defined as in \zcref[cap]{def:global null frame},
  and where the lower order coefficients satisfy
  \begin{align*}
    q\bar{q}^3 W^{[4]},\quad q\bar{q}^3 W^{[3]}, \quad q\bar{q}^3 W^{[h]} = \GO{a} \quad \text{and} \quad q\bar{q}^3W^{[0]}=\GO{a r^{-1}}.
  \end{align*}
\end{proposition}
\begin{proof}
  See proof in \zcref[cap]{appendix:proposition T to RW}.
\end{proof}

The following is an immediate corollary of \zcref[cap]{proposition T to RW}.
\begin{corollary}
  \label{coro:Teukolsky-wave-transport}
  Let $A\in \realHorkTensor{2}(\Complex)$ be a solution to
  $\mathcal{L}(A)=0$ as in \zcref[cap]{proposition T to RW}, $\mathfrak{q}$ as defined in \zcref[noname]{definition of q frak},
  and define $\Psi\in \realHorkTensor{1}(\Complex)$ by
  \begin{equation}
    \label{eq:Teukolsky:Psi:def}
    \Psi \vcentcolon= \ConformalInvDeriv_3\left(\frac{\overline{q}^4}{r^2} A\right).
  \end{equation}
  Then $\left(\mathfrak{q}, \Psi, A\right)$ satisfy the following
  so-called Teukolsky wave-transport system
  \begin{align}
    \label{eq:full-RW-system:wave}
    \left( \WaveOpHork{2}-V \right)\mathfrak{q}
    ={}& \ImagUnit\frac{4a(1+\gamma)\cos\theta}{\abs*{q}^{2}}\nabla_{\KillT}\mathfrak{q}
         + N_0 + N_L,\\
    \label{eq:full-RW-system:transport}
    \begin{split}
      \ConformalInvDeriv_3\Psi ={}& \frac{\overline{q}}{q r^2}\mathfrak{q},\\
      \ConformalInvDeriv_3\left(\frac{\overline{q}^4}{r^2}A\right)
      ={}& \Psi,
    \end{split}
  \end{align}
  where 
  \begin{align}
    \label{eq:full-RW-system:N0:def}
    N_0 = V_0\mathfrak{q},\qquad
    V_0=O(ar^{-4}),\qquad
    \Im (V_0)=0,
  \end{align}
  and
  \begin{equation}
    \label{eq:full-RW-system:NL:def}
    N_L = -\frac{8(1+\gamma)\overline{q}^2\Delta}{\lambdaglo r^2\abs*{q}^2}\left( a^2\nabla_{\KillT}+a\nabla_{\KillPhi} \right)\nabla_3A
    + O(a r^{-3})\dk A + O(ar^{-4})A
    .
  \end{equation}
\end{corollary}


\section{Main theorems}
\label{sec:main-theorem}

In this section, we provide a precise statement of our main results in
\zcref[cap]{MAINTHEOREM} and \zcref[cap]{theo comparaison}. We begin by introducing the main
domains of integration and norms that will be used in what follows.

\subsection{Integration conventions}
\label{sec:domain-of-integration}

In this section, we introduce some notations and conventions related to integrations on $\MM_{\mathrm{tot}}$ (recall the notations introduced in \zcref[cap]{sec:spacetime-regions}). We start with the region $\MM$ where $\tau\geq 0$:
\begin{itemize}
    \item If $0\leq \tau_1< \tau_2$, $\mathcal{M}(\tau_1,\tau_2)$ will denote $\mathcal{M}\cap \curlyBrace*{\tau\in [\tau_1, \tau_2]}$, and similarly for $\mathcal{M}_{\trap}(\tau_1,\tau_2)$, $\mathcal{M}_{\cancel{\trap}}(\tau_1,\tau_2)$, $\mathcal{M}_{\operatorname{red}}(\tau_1,\tau_2)$.
    \item If $0\leq \tau_1< \tau_2$ and $I$ is an interval in $\Real$, then $\mathcal{M}_{r\in I}(\tau_1,\tau_2)$ will denote $\mathcal{M}\cap \curlyBrace*{\tau\in [\tau_1, \tau_2]}\cap \{r\in I\}$, and $\Sigma_{r\in I}(\tau)$ will denote $\Sigma(\tau)\cap \curlyBrace*{r\in I}$.
    \item If $0\leq \tau_1< \tau_2$, we define $\mathcal{A}(\tau_1,\tau_2)=\mathcal{A}\cap \{\tau\in [\tau_1, \tau_2]\}$ and $\SigmaStar(\tau_1,\tau_2)=\SigmaStar\cap \{\tau\in [\tau_1, \tau_2]\}$. In particular, we have
    \begin{align*}
        \dr\mathcal{M}(\tau_1,\tau_2) = \mathcal{A}(\tau_1,\tau_2) \cup \Sigma(\tau_1) \cup \Sigma(\tau_2) \cup \SigmaStar(\tau_1,\tau_2).
    \end{align*}
\end{itemize}
\begin{figure}[h]
  \centering
  \input{Images/Internal-Boundary.tex}
  \caption{\textit{A Penrose diagram depicting $\dr\mathcal{M}(\tau_1,\tau_2)$.}}
  \label{fig:internal-region}
\end{figure}

We now consider the external region $\MM_e$ where $\tau<0$:
\begin{itemize}
    \item If $\tau_1\in[\tau_{\overline{\HH}},0]$, $\mathcal{M}_e(\tau_{\overline{\HH}},\tau_1)$ will denote $\MM_e\cap\{\tau\in[\tau_{\overline{\HH}},\tau_1]\}$.
    \item If $\tau_1\in[\tau_{\overline{\HH}},0]$, we define $\widehat{\Si}_{\init}(\tau_1) = \widehat{\Si}_{\init}\cap\mathcal{M}_e(\tau_{\overline{\HH}},\tau_1)$ and $\SigmaStar(\tau_1)=\SigmaStar\cap \{\tau\in [\tau_{ \overline{\HH} }, \tau_1]\}$. In particular, we have
    \begin{align*}
        \dr\mathcal{M}_e(\tau_{\overline{\HH}},\tau_1) = \Si(\tau_1)\cup \SigmaStar(\tau_1) \cup \widehat{\Si}_{\init}(\tau_1).
    \end{align*}
\end{itemize}
\begin{figure}[ht]
  \centering
  \input{Images/External-boundary.tex}
  \caption{\textit{A Penrose diagram depicting $\dr\mathcal{M}_e(\tau_{\CosmologicalHorizon},\tau)$.}}
  \label{fig:external-region}
\end{figure}

Our convention for integrals on $\mathcal{A}$, $\SigmaStar$, and
$\Sigma(\tau)$ is that for $f$ a scalar function on $\mathcal{M}$
\begin{equation*}
  \begin{split}
    \int_{\mathcal{A}(\tau_1,\tau_2)}f ={}& \int_{\Sphere^2}\int_{\tau_1}^{\tau_2}f\,r^2d \tau d\mathring{\gamma},\\
    \int_{\SigmaStar(\tau_1,\tau_2)}f ={}& \int_{\Sphere^2}\int_{\tau_1}^{\tau_2}f\,r^2d \tau d\mathring{\gamma},\\
    \int_{\Sigma(\tau)}f ={}& \int_{\Sphere^2}\int_{r_{\EventHorizon}(1-\delta_{\Horizon})}^{r_{\CosmologicalHorizon}(1+\delta_{\Horizon})}f\,r^2drd\mathring{\gamma}.
  \end{split}  
\end{equation*}
Note that the volume integration does not coincide with the natural
volume form of $\mathcal{A}$, $\SigmaStar$, and $\Sigma(\tau)$ as
Riemannian manifolds. Indeed since, $\mathcal{A}$ and $\SigmaStar$ are
almost null, and $\Sigma(\tau)$ is asymptotically null, we use the
volume form consistent with null hypersurfaces to avoid any
degeneration.  For consistency with these
integration conventions, we will use
\begin{equation}
  \label{eq:hypersurface-normal-def}
  N_{\Sigma(\tau)} = -\Metric^{\alpha\beta}\partial_{\alpha}\tau\partial_{\beta},\qquad
  N_{\SigmaStar} = -\Metric^{\alpha\beta}\partial_{\alpha}r\partial_{\beta},\qquad
  N_{\mathcal{A}} = \Metric^{\alpha\beta}\partial_{\alpha}r\partial_{\beta},
\end{equation}
as normals to the respective spacelike hypersurfaces\footnote{Note that these normals satisfy
  \begin{equation*}
    N_{\Sigma(\tau)}(r) = 1 + O(r^{-2}),\qquad
    N_{\SigmaStar}(\tau) = 1 + O(r^{-2}),\qquad
    N_{\mathcal{A}}(\tau) = \frac{(r^2+a^2)(1+\gamma)}{\abs*{q}^2} +O(\delta_{\Horizon}).
  \end{equation*}
}. 
These conventions are naturally extended to the external region $\MM_e$, and for integrals on $\widehat{\Si}_{\init}\cap\{r\geq 2r_0\}$ we also use the volume $r^2dr d\mathring{\gamma}$ and consider
\begin{align}
  \label{eq:initial-data-hypersurface-normal:def}
    N_{\widehat{\Si}_{\init}} & = - \g^{\a\b}\dr_\a \underline{t} \dr_\b,
\end{align}
to be the normal to $\widehat{\Si}_{\init}\cap\{r\geq 2r_0\}$.

\subsection{Norms}
\label{sec:main-quantities}

In this section, we present the main norms that will appear in
the statement and proof of the main theorem.  Recall that
$(e_1,e_2,e_3, e_4) = \left( e_1,e_2,\eglo_3, \eglo_4 \right)$ denotes
the global principal null frame defined in \zcref[cap]{def:global
  null frame}.

\subsubsection{Energy-Morawetz norms}
\label{sec:Energy-Morawetz-norms}

Throughout this section, we will use $\psi\in \realHorkTensor{2}$ to
denote any symmetric traceless real-valued horizontal 2-tensor. Note
that all norms extend naturally to
$\psi\in \realHorkTensor{2}(\Complex)$ by combining
the norms of the real and imaginary parts of $\psi$. Moreover, until the last part of this section, we work in $\MM$ and consider $0\leq \tau_1\leq \tau_2$.

\paragraph{Morawetz quantities.}

We first define the Morawetz norms on various regions of $\mathcal{M}$.
\begin{equation}
  \label{eq:eq:Morawetz-bulk-norms:trapping-nontrapping-norm-def}
  \begin{split}
    \MorNorm_{\operatorname{trap}}[\psi](\tau_1,\tau_2)
    \vcentcolon={}& \int_{\Manifold_{\operatorname{trap}}(\tau_1,\tau_2)}\left(\frac{1}{r^2}\abs*{\nabla_{\HprVF} \psi}^2
                    + \frac{1}{r^3}\abs*{\psi}^2\right)
                    ,\\
    \MorNorm_{\cancel{\operatorname{trap}}}[\psi](\tau_1,\tau_2)
    \vcentcolon={}& \int_{\Manifold_{\cancel{\operatorname{trap}}}(\tau_1,\tau_2)}
                    \left( \frac{1}{r^2}\abs*{\nabla_{\HprVF} \psi}^2
                    + \frac{1}{r^3}\abs*{\psi}^2
                    + \frac{1}{r}\abs*{\nabla\psi}^2
                    + \frac{1}{r^2}\abs*{\HawkingVF \psi}^2 \right),\\
    \MorNorm_{\operatorname{red}}[\psi](\tau_1,\tau_2)
    \vcentcolon={}& \int_{\Manifold_{\operatorname{red}}(\tau_1,\tau_2)}
                    \left( \abs*{\nabla_{3}\psi}^2
                    + \abs*{\nabla_{4}\psi}^2
                    + \abs*{\nabla\psi}^2
                    + \abs*{\psi}^2 \right).
  \end{split}  
\end{equation}

We now define the combined Morawetz norms 
\begin{equation}
  \label{eq:Morawetz-bulk-norms:def}
  \begin{split}
    \MorNorm[\psi](\tau_1,\tau_2)
    \vcentcolon={}& \MorNorm_{\operatorname{trap}}[\psi](\tau_1,\tau_2)
                    + \MorNorm_{\cancel{\operatorname{trap}}}[\psi](\tau_1,\tau_2)
                    ,\\
    \MorrNorm[\psi](\tau_1,\tau_2)
    \vcentcolon={}&\MorNorm[\psi](\tau_1,\tau_2)
                    + \MorNorm_{\operatorname{red}}[\psi](\tau_1,\tau_2)
                    + \int_{{\Manifold_{r>4M}(\tau_1,\tau_2)}}r^{-1-\delta}\abs*{\nabla_3\psi}^2,
  \end{split}
\end{equation}
where $\delta>0$ is a sufficiently small positive quantity will be fixed later.

Observe that we have the following equivalent formulation for
$\MorrNorm[\psi](\tau_1,\tau_2)$:
\begin{equation*}
  \begin{split}
    \MorrNorm[\psi](\tau_1,\tau_2)
    \simeq{}& \MorNorm_{\operatorname{trap}}[\psi](\tau_1,\tau_2)
              + \MorNorm_{\operatorname{red}}[\psi](\tau_1,\tau_2)\\
            &  + \int_{\Manifold_{\cancel{\operatorname{trap}}}(\tau_1,\tau_2)}
              \pth{r^{-2}\left(\frac{\Delta}{r^2+a^2}\abs*{\nabla_4\psi}^2 + r^{-1}\abs*{\psi}^2\right)
              + r^{-1}\abs*{\nabla\psi}^2
              + r^{-1-\delta}\abs*{\nabla_3\psi}^2}.
  \end{split}  
\end{equation*}

\paragraph{Weighted bulk quantities.}

Define, for $0<p<2$, and any $R\ge 4M$, 
\begin{equation}
  \label{eq:rp:weighted-bulk:def}
  \begin{split}
    \rpBulkWeighted{p}{R}[\psi](\tau_1,\tau_2)
    \vcentcolon={}& \int_{\Manifold_{r\ge R}(\tau_1,\tau_2)}
                    r^{p-1}\left(p\abs*{\widecheck{\nabla}_4\psi}^2
                    + (2-p)\abs*{\nabla\psi}^2
                    + r^{-2}\abs*{\psi}^2\right),\\
    \rpBulkCombined{p}[\psi](\tau_1,\tau_2)
    \vcentcolon={}& \MorrNorm[\psi](\tau_1,\tau_2)
                    + \rpBulkWeighted{p}{4M}[\psi](\tau_1,\tau_2),
  \end{split}
\end{equation}
where
\begin{equation}
  \label{eq:nable-4-check:def}
  \widecheck{\nabla}_4\psi \vcentcolon= r^{-1}\nabla_4(r\psi). 
\end{equation}
Observe that for $\delta\le p \le 2-\delta$, we have
\begin{equation*}
  \rpBulkCombined{p}[\psi](\tau_1,\tau_2)
  \simeq \MorrNorm[\psi](\tau_1,\tau_2)
  + \int_{\Manifold_{r\ge 4M}(\tau_1,\tau_2)}r^{p-1}
  \left(
    \abs*{\widecheck{\nabla}_4\psi}^2
    + \abs*{\nabla\psi}^2
    + r^{-2}\abs*{\nabla_3\psi}^2
    + r^{-2}\abs*{\psi}^2
  \right),
\end{equation*}
which follows from
\begin{equation*}
  \int_{\Manifold_{r\ge 4M}(\tau_1,\tau_2)}r^{p-3}\abs*{\nabla_3\psi}^2
  \lesssim \int_{\Manifold_{r\ge 4M}(\tau_1,\tau_2)}r^{-1-\delta}\abs*{\nabla_3\psi}^2.
\end{equation*}

\paragraph{Basic energy-flux quantity.}

The basic energy-flux quantity on a hypersurface $\Sigma(\tau)$ is defined by %
\begin{equation}
  \label{eq:energy-flux:basic:def}
  \EnergyFlux[\psi](\tau)
  \vcentcolon= \int_{\Sigma(\tau)}
  \left(
    \abs*{\nabla_4\psi}^2
    + \frac{1}{r^2}\abs*{\nabla_3\psi}^2
    + \abs*{\nabla\psi}^2
    + r^{-2}\abs*{\psi}^2
  \right).
\end{equation}

\paragraph{Horizon degenerate energy norm.}
For $\psi\in \realHorkTensor{2}$, we define the following energy norm along
$\Sigma(\tau)$
\begin{equation}
  \label{eq:horizon-degenerate-energy}
  \begin{split}
    \EnergyFlux_{\operatorname{deg}}[\psi](\tau)
    \vcentcolon={}& \int_{\Sigma_{r\le 4M}(\tau)}
                    \left(
                    \abs*{\nabla_{4}\psi}^2
                    + \frac{\abs*{\Delta}}{r^4}\abs*{\nabla_{3}\psi}^2
                    + \abs*{\nabla\psi}^2
                    + r^{-2}\abs*{\psi}^2
                    \right)\\
                  & + \int_{\Sigma_{r\ge 4M}(\tau)}
                    \left(
                    \frac{1}{r^2}\abs*{\nabla_{3}\psi}^2
                    + \frac{\abs*{\Delta}}{r^2}\abs*{\nabla_{4}\psi}^2
                    + \abs*{\nabla\psi}^2
                    + r^{-2}\abs*{\psi}^2
                    \right).
  \end{split}      
\end{equation}

\paragraph{Weighted energy-flux type quantities.}

Define
\begin{equation}
  \label{eq:energy-flux:weighted:def}
  \EnergyFluxWeighted[\psi](\tau) \vcentcolon=
  \begin{cases}\displaystyle
    \int_{\Sigma_{r\ge R}(\tau)} r^p\left(\abs*{\widecheck{\nabla}_4\psi}^2+r^{-2}\abs*{\psi}^2\right) &\text{for } p\le 1-\delta,\\
    \displaystyle
    \int_{\Sigma_{r\ge R}(\tau)} r^p\left(\abs*{\widecheck{\nabla}_4\psi}^2+r^{-p-1-\delta}\abs*{\psi}^2\right) &\text{for } p> 1-\delta,
  \end{cases}
\end{equation}
and
\begin{equation}
  \label{eq:energy-flux:combined:def}
  \EnergyFluxCombined[\psi](\tau)\vcentcolon= \EnergyFlux[\psi](\tau) + \EnergyFluxWeightedOpt{p}{4M}[\psi](\tau).
\end{equation}

\paragraph{Flux quantities.}

The basic flux quantities along the spacelike hypersurfaces
$\mathcal{A}$ and $\SigmaStar$ are given by
\begin{equation}
  \label{eq:flux-spacelike:basic-def}
  \begin{split}
    \mathbf{F}_{\mathcal{A}}[\psi](\tau_1,\tau_2)
    &\vcentcolon= \int_{\mathcal{A}(\tau_1,\tau_2)}
      \left(
      \abs*{\nabla_{4}\psi}^2
      + \delta_{\Horizon}
      \abs*{\nabla_{3}\psi}^2
      + \abs{\nabla \psi}^2
      + r^{-2}\abs*{\psi}^2
      \right),\\
    \mathbf{F}_{\SigmaStar}[\psi](\tau_1,\tau_2)
    & \vcentcolon= \int_{\SigmaStar(\tau_1,\tau_2)}\left(
      \delta_{\Horizon}\abs*{\nabla_{4}\psi}^2
      +\abs*{\nabla_{3}\psi}^2
      +\abs*{\nabla\psi}^2
      + r^{-2}\abs*{\psi}^2
      \right),\\
    \mathbf{F}[\psi](\tau_1,\tau_2)
    &\vcentcolon={}  \mathbf{F}_{\mathcal{A}}[\psi](\tau_1,\tau_2) +  \mathbf{F}_{\SigmaStar}[\psi](\tau_1,\tau_2)
      .
  \end{split}
\end{equation}

\paragraph{Weighted flux quantities.}

For $\delta\le p\le 2-\delta$ where $0<\delta\ll 1$ is chosen later, we define
\begin{equation}
  \label{eq:flux-rp-weighted:def}
  \begin{split}
    \SpacelikeFluxWeighted[\psi](\tau_1,\tau_2)
    \vcentcolon={}& \int_{\SigmaStar(\tau_1,\tau_2)}r^p\left(
                    \delta_{\Horizon}\abs*{\nabla_{4}\psi}^2
                     + \abs*{\nabla\psi}^2
                    + r^{-2}\abs*{\psi}^2
                    \right)
                    ,\\
    \mathbf{F}_p[\psi](\tau_1,\tau_2)
    \vcentcolon={}& \mathbf{F}[\psi](\tau_1,\tau_2)
                    + \dot{\mathbf{F}}_p[\psi](\tau_1,\tau_2).
  \end{split}  
\end{equation}

\paragraph{Forcing term norms.}

We define the following basic norm for the forcing terms %
\begin{equation}
  \label{eq:forcing-term-basic-norm:def}
  \begin{split}
    \ForcingTermNorm[\psi,N](\tau_1,\tau_2)
    \vcentcolon={}& \int_{\Manifold(\tau_1,\tau_2)}\left(\abs*{\nabla_{\HprVF}\psi} + r^{-1}\abs*{\psi}\right)\abs*{N}  + \abs*{\int_{\Manifold_{\operatorname{trap}}(\tau_1,\tau_2)}\nabla_{\TAlmostKilling}\psi\cdot N}
                    + \int_{\Manifold_{\cancel{\operatorname{trap}}}(\tau_1,\tau_2)}\abs*{D\psi}\abs*{N}\\
                  & + \int_{\Manifold(\tau_1,\tau_2)}\abs*{N}^2
                    + \sup_{\tau\in [\tau_1,\tau_2]}\int_{\Sigma(\tau)}\abs*{N}^2
                    + \int_{\SigmaStar(\tau_1,\tau_2)}\abs*{N}^2,
  \end{split}  
\end{equation}
where $\TAlmostKilling$ is the vectorfield constructed in
\zcref[noname]{eq:TALmostkilling:def}.

We define the following weighted norms for the forcing term $N$. %
\begin{equation}
  \label{eq:forcing-term-weighted-norms:def}
  \begin{split}
    \ForcingTermWeightedNorm{p}{R}[\psi, N](\tau_1,\tau_2)
    \vcentcolon={}& \int_{\Manifold_{r\ge R}(\tau_1,\tau_2)}\left( \abs*{r^p\widecheck{\nabla}_4\psi} + \abs*{r^{p-1}\psi} \right)\cdot \abs*{N},
    \\
    \ForcingTermCombinedNorm{p}[\psi, N](\tau_1,\tau_2)
    \vcentcolon={}& \ForcingTermNorm[\psi, N](\tau_1,\tau_2)
                    + \ForcingTermWeightedNorm{p}{4M}[\psi,N](\tau_1,\tau_2).
  \end{split}  
\end{equation}

\paragraph{Combined norms.}

We denote the combined norm
\begin{align}
  \label{eq:combined-BEF-norms:def}
  \CombinedBEFNorm{p}[\psi](\tau_1,\tau_2)&\vcentcolon=
  \sup_{\tau\in [\tau_1,\tau_2]}\EnergyFluxCombined[\psi](\tau)
  + \rpBulkCombined{p}[\psi](\tau_1,\tau_2)
  + \SpacelikeFluxCombined[\psi](\tau_1,\tau_2),\\
  \label{eq:combined-EF-norms:def}
  \CombinedEFNorm{p}[\psi](\tau_1,\tau_2)
  &\vcentcolon= \sup_{\tau\in[\tau_1,\tau_2]}\EnergyFluxCombined[\psi](\tau)
  + \SpacelikeFluxCombined[\psi](\tau_1,\tau_2).
\end{align}
We also define the auxiliary combined norm
\begin{equation}
  \label{eq:combined-wBEF-norm:def}
  \begin{aligned}
    \WeightedBEFNorm{p,R}[\psi](\tau_1,\tau_2)
    &\vcentcolon= \sup_{\tau\in[\tau_1,\tau_2]}\EnergyFluxWeighted[\psi](\tau)
      + \SpacelikeFluxWeighted[\psi](\tau_1,\tau_2)
      + \rpBulkWeighted{p}{R}[\psi](\tau_1,\tau_2)
 + \int_{\Manifold_{r\ge R}(\tau_1,\tau_2)}r^{-1-\delta}\abs*{\nabla_3\psi}^2.
  \end{aligned}
\end{equation}

\paragraph{Higher-order norms.}

Consider the following set of weighted derivatives, %
\begin{align*}
  \frakWeightedDerivAngular  \vcentcolon= \curlyBrace*{r\nab},\qquad
  \frakWeightedDeriv  \vcentcolon= \curlyBrace*{\nabla_{3}, r\nabla_{4}, \frakWeightedDerivAngular}.
\end{align*}
For all norms $Q[\psi]$ defined above, we define the respective
higher-order norms by
\begin{equation}\label{HON}
  Q^s[\psi] \vcentcolon = \sum_{k\le s}Q[\frakWeightedDeriv^k\psi].
\end{equation}

\paragraph{Norms in the external region.} In the external region $\MM_e$, we define for $\tau\in [\tau_{\CosmologicalHorizon}, 0]$,
\begin{align}
  \label{eq:external-bulk-norm:def}
  \BulkNormExternal{p}[\psi](\tau_{\CosmologicalHorizon}, \tau)
  &\vcentcolon={} \int_{\Manifold_e(\tau_{\CosmologicalHorizon},\tau)}r^{p-1}
  \left(
    \abs*{\widecheck{\nabla}_4\psi}^2
    + \abs*{\nabla\psi}^2
    + r^{-2}\abs*{\nabla_3\psi}^2
    + r^{-2}\abs*{\psi}^2
  \right),\\
  \label{eq:external-BEF-norm:def}
  \ExternalBEFNorm{p}[\psi](\tau_{\CosmologicalHorizon}, \tau)
  &\vcentcolon={} \BulkNormExternal{p}[\psi](\tau_{\CosmologicalHorizon}, \tau)
  + \sup_{\tau'\in [\tau_{\CosmologicalHorizon},\tau]}\EnergyFluxCombinedOpt{p}[\psi](\tau')
  + \SpacelikeFluxFar_p[\psi](\tau_{\CosmologicalHorizon},\tau),
\end{align}
where $\EnergyFluxCombinedOpt{p}[\psi](\tau')$ is defined by \zcref[noname]{eq:energy-flux:combined:def} (note that if $\tau \leq 0$, then the restriction $r\geq 4M$ is always satisfied) and $\SpacelikeFluxFar_p[\psi](\tau_{\CosmologicalHorizon},\tau)$ by \zcref[cap]{eq:flux-rp-weighted:def}. Moreover, we define
\begin{align}
    \label{eq:forcing-term-external-norm:def}
  \ForcingTermExternalNorm{p}[\psi, N](\tau_{\CosmologicalHorizon}, \tau)
  &\vcentcolon= \int_{\Manifold_e(\tau_{\CosmologicalHorizon},\tau)}r^{p-1}\left(r\abs*{\widecheck{\nabla}_4\psi} + \abs*{\psi}\right)\cdot\abs*{N}.
\end{align}
In addition, we introduce a norm for the data on $\widehat{\Si}_\init\cap\{r\geq 2r_0\}$, based on the following set of weighted derivatives:
\begin{align*}
    \widehat{\dk} & \vcentcolon= \curlyBrace*{r\nabla_{3}, r\nabla_{4},\frakWeightedDerivAngular}.
\end{align*}
For $\tau\in[\tau_{\overline{\HH}},0]$, $\de\leq p\leq 2-\de$ and $s\geq 0$ we define
\begin{align*}
    \widehat{\mathbf{E}}^s_p[\psi](\tau) & \vcentcolon = \int_{\SigmaInit(\tau)}r^{p+2}\left| \widehat{\dk}^{\leq s+2}\psi \right|^2.
\end{align*}
Recalling \zcref[cap]{definition de ell}, observe that on $\widehat{\Si}_\init$, the restriction $\tau\leq \tau_1$ defining $\SigmaInit(\tau_1)$
(with $\tau_1\leq 0$) is equivalent to $r\geq \ell^{-1}(-\tau_1)$, where $\ell^{-1}: [\tau_{\CosmologicalHorizon},0] \longrightarrow \left[2r_0,r_{\overline{\HH}}(1+\de_\HH)\right]$.

\subsubsection{Combined norms for \texorpdfstring{$(\mathfrak{q}, A)$}{q, A}}
\label{sec:combined-norms-q-A}

We now introduce the norms adapted to the Teukolsky wave-transport
system in \zcref[cap]{coro:Teukolsky-wave-transport}, similarly starting with the region $\{\tau\geq 0\}$.

\paragraph{Norms for \texorpdfstring{$A$}{A}.}
For $0\leq \tau_1\leq \tau_2$, we define modified bulk and energy-flux norms for $A\in\realHorkTensor{2}(\Complex)$, solving the Teukolsky equation.
  \begin{equation}
    \label{eq:A-norms}
    \begin{aligned}
      \BulkNormWeighted{p}[A](\tau_1,\tau_2)
      \vcentcolon ={}& \int_{\mathcal{M}(\tau_1,\tau_2)}r^{p+1}\left(
            r^4\abs*{\nabla_3\nabla_3A}^2
            + r^4\abs*{\nabla_4\nabla_3A}^2
            + r^4\abs*{\nabla\nabla_3A}^2\right.\\
          &\left. \hspace{3cm}+r^2\abs*{\nabla_3A}^2
            + r^2\abs*{\nabla_4A}^2
            + r^2\abs*{\nabla A}^2
            + \abs*{A}^2
            \right),
      \\
      \EnergyFluxCombined[A](\tau)
      \vcentcolon={}& \int_{\Sigma(\tau)}r^{p+2}\Big(
            r^4\chi_{\NT}^2\abs*{\nabla_4\nabla_3A}^2
            + r^2\abs*{\nabla_{\HprVF}\nabla_3A}^2
             + r^2\abs*{\nabla_3A}^2
            + r^2\abs*{\nabla_4A}^2
            + \abs*{A}^2
            \Big),\\
      \SpacelikeFluxCombined[A](\tau_1,\tau_2)
      \vcentcolon={}& \int_{\mathcal{A}(\tau_1,\tau_2)\cup\SigmaStar(\tau_1,\tau_2)}r^{p+2}\Big(
            r^4\chi_{\NT}^2\abs*{\nabla_4\nabla_3A}^2
            + r^2\abs*{\nabla_{\HprVF}\nabla_3A}^2 
          + r^2\abs*{\nabla_3A}^2
            + r^2\abs*{\nabla_4A}^2
            + \abs*{A}^2
            \Big),
    \end{aligned}    
  \end{equation}
  where $\chi_{\NT}=\chi_{\NT}(r)\in C^{\infty}$ is a smooth
  cutoff function such that
  \begin{equation}
    \label{eq:chi-nontrapping-asymptotic:def}
    \chi_{\NT}(r) = 0 \text{ for } \abs*{r-3M}< \delta_{\operatorname{trap}},\qquad
    \chi_{\NT}(r) = 1 \text{ for }   r > 4M.
  \end{equation}
  We also define the combined norms
  \begin{equation}
    \label{eq:BEF-A-norm:def}
    \CombinedBEFNorm{p}[A](\tau_1,\tau_2)
    \vcentcolon= \BulkNormWeighted{p}[A](\tau_1,\tau_2)
    + \sup_{\tau\in [\tau_1,\tau_2]}\EnergyFluxCombined[A](\tau)
    + \SpacelikeFluxCombined[A](\tau_1,\tau_2).
  \end{equation}
Finally, we also define an auxiliary bulk norm for $A$ that does not
control its angular derivatives
  \begin{equation}
    \label{eq:bulk-A-aux-norm:def}
    \BulkNormWeightedAux{p}[A](\tau_1,\tau_2)
    \vcentcolon=\int_{\mathcal{M}(\tau_1,\tau_2)}r^{p+1}\left(
      r^4\abs*{\nabla_3\nabla_3A}^2
      + r^4\abs*{\nabla_4\nabla_3A}^2
      + r^2\abs*{\nabla_3A}^2
      + r^2\abs*{\nabla_4A}^2
      + \abs*{A}^2
    \right),
  \end{equation}
  as well as the corresponding auxiliary combined norm
  \begin{equation}
    \label{eq:BEF-A-aux-norm:def}
    \BEFNormAux{p}[A](\tau_1,\tau_2)
    \vcentcolon=\BulkNormWeightedAux{p}[A](\tau_1,\tau_2)
    + \sup_{\tau\in [\tau_1,\tau_2]}\EnergyFluxCombined[A](\tau)
    + \SpacelikeFluxCombined[A](\tau_1,\tau_2).
  \end{equation}

\paragraph{Combined norms for $(\mathfrak{q}, A)$.}

We define 
\begin{equation}
  \label{eq:q-A-combined norms:def}
  \begin{split}
    \BulkNormWeighted{p}[\mathfrak{q},A](\tau_1,\tau_2)\vcentcolon={}&  \BulkNormWeighted{p}[\mathfrak{q}](\tau_1,\tau_2) + \BulkNormWeighted{p}[A](\tau_1,\tau_2),\\
    \EnergyFluxCombined[\mathfrak{q}, A](\tau_1,\tau_2)\vcentcolon={}& \EnergyFluxCombined[\mathfrak{q}](\tau_1,\tau_2) + \EnergyFluxCombined[A](\tau_1,\tau_2),\\
    \SpacelikeFluxCombined[\mathfrak{q}, A](\tau_1,\tau_2)\vcentcolon={}& \SpacelikeFluxCombined[\mathfrak{q}](\tau_1,\tau_2) + \SpacelikeFluxCombined[A](\tau_1,\tau_2),
  \end{split}  
\end{equation}
where $\BulkNormWeighted{p}[\mathfrak{q}](\tau_1,\tau_2)$,
$\EnergyFluxCombined[\mathfrak{q}](\tau_1,\tau_2)$, and
$\SpacelikeFluxCombined[\mathfrak{q}](\tau_1,\tau_2)$ refer to the
norms defined in \zcref[cap]{sec:Energy-Morawetz-norms}. Finally, we set
\begin{equation}
    \CombinedBEFNorm{p}[\mathfrak{q}, A](\tau_1,\tau_2) \vcentcolon = \BulkNormWeighted{p}[\mathfrak{q},A](\tau_1,\tau_2) + \EnergyFluxCombined[\mathfrak{q}, A](\tau_1,\tau_2) + \SpacelikeFluxCombined[\mathfrak{q}, A](\tau_1,\tau_2).
\end{equation}

\paragraph{Combined norms for $(\mathfrak{q}, A)$ in the external region.} The norms defined in \zcref[noname]{eq:A-norms} for $A$ can be naturally extended to the external region, by simply replacing $\MM(\tau_1,\tau_2)$ by $\MM_e(\tau_{\overline{\HH}},\tau)$ (for $\tau\leq 0$) which defines
\begin{align*}
    \BulkNormExternal{p}[A](\tau_{\CosmologicalHorizon}, \tau) & \vcentcolon = \int_{\MM_e(\tau_{\overline{\HH}},\tau)}r^{p+1}\left(
            r^4\abs*{\nabla_3\nabla_3A}^2
            + r^4\abs*{\nabla_4\nabla_3A}^2
            + r^4\abs*{\nabla\nabla_3A}^2\right.\\
          &\left. \hspace{4cm}+r^2\abs*{\nabla_3A}^2
            + r^2\abs*{\nabla_4A}^2
            + r^2\abs*{\nabla A}^2
            + \abs*{A}^2
            \right) ,
    \\ \SpacelikeFlux_p[A](\tau_{\CosmologicalHorizon}, \tau)
      \vcentcolon={}& \int_{\SigmaStar(\tau_{\overline{\HH}},\tau)}r^{p+2}\Big(
            r^4\abs*{\nabla_4\nabla_3A}^2
            + r^2\abs*{\nabla_{\HprVF}\nabla_3A}^2 + r^2\abs*{\nabla_3A}^2
            + r^2\abs*{\nabla_4A}^2
            + \abs*{A}^2
            \Big),\\
  \EnergyFluxInitAux{p}[A](\tau)
      \vcentcolon={}& \int_{\SigmaInit(\tau)}r^{p+2}\Big(
            r^4\abs*{\nabla_4\nabla_3A}^2
            + r^2\abs*{\nabla_{\HprVF}\nabla_3A}^2
             + r^2\abs*{\nabla_3A}^2
            + r^2\abs*{\nabla_4A}^2
            + \abs*{A}^2
            \Big),
\end{align*}
while $\mathbf{E}_p[A](\tau)$ stays as in\footnote{We remark that
  $\chi_{\operatorname{nt}}=1$ in the external region.}
\zcref[noname]{eq:A-norms}. We will also define the following
auxiliary initial energy norm on $\SigmaInit(\tau)$, 
\begin{equation}
  \label{eq:external:combined-q-A-norms:energy}
  \EnergyFluxInitAux{p}^s[\mathfrak{q}, A](\tau)
  \vcentcolon={} \EnergyFluxCombined[\mathfrak{q}](\tau) + \EnergyFluxCombined[A](\tau),\\
\end{equation}

We define
\begin{align}
    \ExternalBEFNorm{p}[A](\tau_{\CosmologicalHorizon}, \tau) \vcentcolon = \BulkNormExternal{p}[A](\tau_{\CosmologicalHorizon}, \tau) + \mathbf{E}_p[A](\tau) + \SpacelikeFlux_p[A](\tau_{\CosmologicalHorizon}, \tau),
\end{align}
and finally
\begin{equation}
  \label{eq:external-BEF-norm}  
  \ExternalBEFNorm{p}[\mathfrak{q}, A](\tau_{\overline{\HH}},\tau) \vcentcolon = \ExternalBEFNorm{p}[A](\tau_{\CosmologicalHorizon}, \tau) +  \ExternalBEFNorm{p}[\psi](\tau_{\CosmologicalHorizon}, \tau),
\end{equation}
where $\ExternalBEFNorm{p}[\psi](\tau_{\CosmologicalHorizon}, \tau)$ has been defined in \zcref[noname]{eq:external-BEF-norm:def}.

\begin{remark}
    Higher-order norms for $A$ or in the external regions are denoted and defined as in \zcref[noname]{HON}.
\end{remark}

\subsection{Precise statement of the main theorem}
\label{subsec:main-theorem}

The precise statement of our main theorem is as follows.

\begin{theorem}
  \label{MAINTHEOREM}
  Let $(\mathfrak{q}, A)$ be a solution to the Teukolsky
  wave-transport system in \zcref[noname]{eq:full-RW-system:wave} and
  \zcref[noname]{eq:full-RW-system:transport} on $\MM_{\mathrm{tot}}$ endowed with
  $\Metric_{M,a,\Lambda}$, a slowly-rotating Kerr-de Sitter black
  hole metric. Let $s\geq 0$, $\de\in(0,1)$ and $p\in[\de,2-\de]$.  There exists $\La_0>0$ such that the following estimates hold uniformly over $\La\in (0,\La_0]$: for all $0\leq \tau_1< \tau_2$ we have
      \begin{equation}\label{eq:main-q-A-combined-estimate}
        \CombinedBEFNorm{p}^s[\mathfrak{q}, A](\tau_1,\tau_2) \lesssim  \EnergyFluxCombined^s[\mathfrak{q}, A](\tau_1),
      \end{equation}
      and for all $\tau_{\overline{\HH}}\leq \tau \leq 0$ we have
      \begin{equation}\label{eq:main-external-estimate}
          \ExternalBEFNorm{p}^s[\mathfrak{q}, A](\tau_{\overline{\HH}},\tau) \lesssim \widehat{\mathbf{E}}^{s+1}_p[A](\tau).
      \end{equation}
\end{theorem}

Although it is not explicit from our notations for the various norms, we recall that the integration domains hidden in \zcref[noname]{eq:main-q-A-combined-estimate} and \zcref[noname]{eq:main-external-estimate} extend to $r\sim \La^{-\half}$ and in particular to regions where Kerr-de Sitter is not an approximation of Kerr, so that the uniformity of these estimates over $\La\in(0,\La_0]$ is far from trivial, and is the main result of this article.

\subsection{The vanishing $\La$ limit for Teukolsky}

In this section we state a consequence of \zcref[cap]{MAINTHEOREM} in the vanishing $\La$ limit where Kerr-de Sitter converges pointwise to Kerr. To this end, we define the Kerr manifold in the coordinates from \zcref[cap]{sec:adaptedglobalcoordinates} by
    \begin{align*}
        \MM_{\mathrm{tot},0} &  =  \pth{\mathbb{R}_\tau \times \left( r_{\EventHorizon,0}(1-\delta_{\Horizon}), +\infty \right)_r\times \Sphere^2}\cap \{\underline{t}\geq 0\},
    \end{align*}
    where $r_{\EventHorizon,0}=M+\sqrt{M^2-a^2}$ is the location of
    the event horizon in Kerr.  On this manifold, endowed with a
    slowly-rotating Kerr metric $\g_{M,a}$, we consider the global
    principal null frame $(e^0_\mu)_{\mu=1,2,3,4}$, obtained by
    setting $\La=0$ in the corresponding global principal null frame
    introduced for Kerr-de Sitter in \zcref[cap]{sec:global null
      frame}. We can thus consider the same Teukolsky wave transport
    system \zcref[noname]{eq:full-RW-system:wave} and
    \zcref[noname]{eq:full-RW-system:transport} on Kerr. For a
    solution $(\q^0,A^0)$ of this system, the norms are defined as in
    \zcref[cap]{sec:main-quantities} but on the slightly smaller
    manifold defined by
    \begin{align*}
        \MM'_{\mathrm{tot},0} &  =  \pth{\mathbb{R}_\tau \times \left( r_{\EventHorizon,0}\pth{1-\frac{\delta_{\Horizon}}{2}}, +\infty \right)_r\times \Sphere^2}\cap \{\underline{t}\geq 0\},
    \end{align*}
    and are denoted similarly but with the addition of a subscript
    $0$. Moreover, we define
    $\widehat{\Si}_{\init,0}=\MM_{\mathrm{tot},0}\cap\{\underline{t}=0\}$
    and add a $0$ exponent or subscript to quantities (such as
    weighted derivatives or horizontal tensorial norms) defined with
    the Kerr metric.

    The uniformity of the estimates proved in \zcref{MAINTHEOREM} and the fact that they hold on the region $r\lesssim \La^{-\half}$ allow us to pass to the limit $\La\to 0$ and get the following corollary.

\begin{corollary}\label{theo comparaison}
    Let $(\q^0,A^0)$ be a solution to the Teukolsky wave-transport system of Kerr on $\MM_{\mathrm{tot},0}$. Let $s\geq 6$, $0<\de<1$ and $p\in[\de,2-\de)$ and assume that
    \begin{align}\label{assumption comparaison}
        \int_{\widehat{\Si}_{\init,0}} r^{4-\de} \left| \pth{ \widehat{\dk}^0}^{\leq s+2} A^0 \right|^2_0 < +\infty.
    \end{align}
   Then for all $0\leq \tau_1< \tau_2$ we have 
    \begin{align}\label{internal estimate pour Kerr}
        \mathbf{BE}^{s-6,0}_p\left[ \q^0,A^0\right](\tau_1,\tau_2)\lesssim \mathbf{E}^{s-6,0}_p\left[ \q^0,A^0\right](\tau_1),
    \end{align}
    and for all $\tau\leq 0$ we have
    \begin{align}\label{external estimate pour Kerr}
        \mathring{\mathbf{BE}}^{s-6,0}_p \left[ \q^0,A^0\right](-\infty,\tau) \lesssim \widehat{\mathbf{E}}^{s-5}_p\left[A^0\right](\tau).
    \end{align}
\end{corollary}

\begin{remark}
    The proof of \zcref[cap]{theo comparaison} doesn't yield the control of fluxes on the left and right boundaries of $\MM'_{\mathrm{tot},0}$. However, they can easily be recovered \textit{a posteriori} from \zcref[noname]{internal estimate pour Kerr} and \zcref[noname]{external estimate pour Kerr}.
\end{remark}

\begin{remark}
    Note that in order to pass to the limit $\La\to0$, it would have been enough to prove uniform estimates up to $r\sim \La^{-\eta}$ for some $0<\eta<\half$, a region where Kerr-de Sitter converges uniformly to Kerr. \zcref{theo comparaison} should thus be considered as a limited application of \zcref{MAINTHEOREM}.
\end{remark}

\subsection{Organisation of the rest of the article}
\label{sec:road-map}

The remaining sections of this article are devoted to the proof of \zcref[cap]{MAINTHEOREM} and \zcref[cap]{theo comparaison}. Here, we give a brief overview of their content.

\zcref[cap]{sec:wave} contains generic computations on the wave
operator $\WaveOpHork{k} - V$, appearing in the generalized
Regge-Wheeler equation \zcref[noname]{eq:full-RW-system:wave}. The proof of \zcref[cap]{MAINTHEOREM} is spread over
\zcref[cap]{sec:Killing-energy-estimate} to
\zcref[cap]{sec:external-estimates}. More precisely, we prove the
internal estimate \zcref[noname]{eq:main-q-A-combined-estimate} from
\zcref[cap]{sec:Killing-energy-estimate} to
\zcref[cap]{sec:transport}:
\begin{itemize}
\item In \zcref[cap]{sec:Killing-energy-estimate} we control the
  degenerate energies $\mathbf{E}_{\mathrm{deg}}[\psi]$ by using an almost
  Killing timelike vectorfield as our vectorfield multiplier.
\item \zcref[cap]{sec:setup-freq-analysis} contains a brief
  introduction to the pseudo-differential framework used in
  \zcref[cap]{sec:Morawetz} to prove a Morawetz estimate which
  captures the degeneracy at the trapped set in \KdS.
\item \zcref[cap]{sec:rp} proves a $r^p$-weighted estimate in a region
  where $r>r_0$, where $r_0$ is large but not $\Lambda$-dependent,
  that in particular includes the cosmological horizon. 
\item \zcref[cap]{sec:redshift} handles the ergoregion near the black
  hole event horizon via a redshift estimate.
\item \zcref[cap]{sec:combined-estimates-on-KdS} and
  \zcref[cap]{sec:transport} conclude the proof of
  \zcref[noname]{eq:main-q-A-combined-estimate} by first stating a
  complete estimate for the solution $\q$ to
  \zcref[noname]{eq:full-RW-system:wave} and then estimating $A$ via
  the full Teukolsky wave-transport system
  \zcref[noname]{eq:full-RW-system:wave} and
  \zcref[noname]{eq:full-RW-system:transport}.
\end{itemize}
The external estimate \zcref[noname]{eq:main-external-estimate} is
proved in \zcref[cap]{sec:external-estimates}, which thus concludes
the proof of \zcref[cap]{MAINTHEOREM}. Finally, \zcref[cap]{theo comparaison} is proved in \zcref[cap]{section preuve theorem comparaison}.

\section{Wave operator and energy-momentum tensor}
\label{sec:wave}

To study the tensorial wave equation
\zcref[noname]{eq:full-RW-system:wave}, we collect material related to
the conformal wave operator $\WaveOpHork{2}$ in this section. To this
end, we consider more generally the model equation
\begin{equation}
  \label{eq:model-wave-equation}
  \WaveOpHork{k}\psi - V\psi= \ForcingTerm,
\end{equation}
where $V$ is a real potential, and $\ForcingTerm$ represents some
suitable forcing term.

\subsection{Energy-momentum tensor and divergence theorem \texorpdfstring{of \zcref[noname]{eq:model-wave-equation}}{}}
\label{sec:EM-Tensor}

\begin{definition}
  \label{def:EM-Tensor}
  The \emph{energy-momentum tensor} associated to
  \zcref[noname]{eq:model-wave-equation} is given by
  \begin{equation}
    \label{eq:EM-Tensor:def}
    \EMTensor_{\mu\nu}[\psi]
    \vcentcolon= \HorCovDeriv_{\mu}\psi\cdot\HorCovDeriv_{\nu}\psi
    - \frac{1}{2}\Metric_{\mu\nu}\left(
      \HorCovDeriv_{\lambda}\psi \cdot \HorCovDeriv^{\lambda}\psi
      + V\psi\cdot\psi
    \right)
    = \HorCovDeriv_{\mu}\psi\cdot\HorCovDeriv_{\nu}\psi
    - \frac{1}{2}\Metric_{\mu\nu}\mathcal{L}[\psi],
  \end{equation}
  where the Lagrangian $\mathcal{L}[\psi]$ is given by
  \begin{equation}
    \label{eq:model-wave-equation:Lagrangian}
    \mathcal{L}[\psi] = \Metric^{\mu\nu}\HorCovDeriv_{\mu}\psi\cdot\HorCovDeriv_{\nu}\psi
    + V \psi\cdot\psi.
  \end{equation}
\end{definition}

It is easy to compute the null components of the energy-momentum tensor.
\begin{lemma}
  \label{lemma:EM-tensor:null-components}
  In any null frame, the energy-momentum tensor has the following null components:
  \begin{equation}
    \label{eq:EM-tensor:null-components}
    \begin{gathered}
      \EMTensor[\psi]_{33} = \abs*{\nabla_3\psi}^2,\qquad
      \EMTensor[\psi]_{44} = \abs*{\nabla_4\psi}^2,\qquad
      \EMTensor[\psi]_{34} = \abs*{\nabla\psi}^2 + V \abs*{\psi}^2,\\
      \EMTensor[\psi]_{3a} = \nabla_3\psi\cdot\nabla_a\psi,\qquad
      \EMTensor[\psi]_{4a} = \nabla_4\psi\cdot\nabla_a\psi,\qquad
      \delta^{bc}\EMTensor[\psi]_{bc} = \nabla_3\psi\cdot\nabla_4\psi - V\abs*{\psi}^2.
    \end{gathered} 
  \end{equation}
\end{lemma}

The energy-momentum tensor has the following important divergence
property.
\begin{lemma}
  \label{lemma:EM-Tensor:divergence-prop}
  Given a solution $\psi\in \realHorkTensor{k}$ of equation
  \zcref[noname]{eq:model-wave-equation}, we have
  \begin{equation}
    \label{eq:EM-Tensor:divergence-prop}
    \CovariantDeriv^{\nu}\EMTensor_{\mu\nu}
    = \HorCovDeriv_{\mu}\psi\cdot\left(\WaveOpHork{k}\psi - V\psi\right)
    + \HorCovDeriv^{\nu}\psi^A\HorRiem_{AB\nu\mu}\psi^{B}
    - \frac{1}{2}\CovariantDeriv_{\mu}V\abs*{\psi}^2.
  \end{equation}
\end{lemma}

\begin{proof}
  See the proof of Lemma 4.7.1 in \cite{giorgiWaveEquationsEstimates2024}.
\end{proof}

We will heavily use the following applications of the divergence
property in \zcref[cap]{lemma:EM-Tensor:divergence-prop}.

\begin{proposition}
  \label{prop:div-thm:general}
  Let $\psi\in \realHorkTensor{k}$ be a solution to
  \zcref[noname]{eq:model-wave-equation}, $\MorawetzVF$ be a vectorfield,
  $\MorawetzLagrangeCorr$ be a scalar, and $\MorawetzOneForm$ be a
  one-form, and define
  \begin{equation}
    \label{eq:J-current:def}
    \JCurrent{\MorawetzVF, \MorawetzLagrangeCorr, \MorawetzOneForm}_{\mu}[\psi]
    \vcentcolon= \EMTensor_{\mu\nu}\MorawetzVF^{\nu}
    + \MorawetzLagrangeCorr\psi\cdot\HorCovDeriv_{\mu}\psi
    - \frac{1}{2}\abs*{\psi}^2\partial_{\mu}\MorawetzLagrangeCorr
    + \frac{1}{2}\abs*{\psi}^2\MorawetzOneForm_{\mu}.
  \end{equation}
  Then,
  \begin{equation}
    \label{eq:div-thm:general}
    \CovariantDeriv\cdot\JCurrent{\MorawetzVF, \MorawetzLagrangeCorr, \MorawetzOneForm}[\psi]
    = \KCurrent{\MorawetzVF, \MorawetzLagrangeCorr, \MorawetzOneForm}[\psi]    
    + \left(\MorawetzVF(\psi)+ \MorawetzLagrangeCorr\psi\right)\cdot \left(\WaveOpHork{k}\psi - V\psi\right),
  \end{equation}
  where
  \begin{equation}
    \label{eq:K-current:def}
    \begin{split}
      \KCurrent{\MorawetzVF, \MorawetzLagrangeCorr, \MorawetzOneForm}[\psi]
      \vcentcolon={}& \EMTensor[\psi]\cdot\DeformationTensor{\MorawetzVF}
            + \MorawetzVF^{\mu}\HorCovDeriv^{\nu}\psi^a\HorRiem_{ab\nu\mu}\psi^b
            - \frac{1}{2}\MorawetzVF(V)\abs*{\psi}^2
           + \MorawetzLagrangeCorr\mathcal{L}[\psi]
            - \frac{1}{2}\abs*{\psi}^2\ScalarWaveOp[\Metric]\MorawetzLagrangeCorr
            + \frac{1}{2}\Divergence\left(\abs*{\psi}^2\MorawetzOneForm\right).
    \end{split}    
  \end{equation}
  We emphasize that in our notation we have $\DeformationTensor[_{\alpha\beta}]{X} = \CovariantDeriv_{(\alpha}X_{\beta)} = -\frac{1}{2}\LieDerivative_Xg$,
  where the parentheses in the indices denote the usual symmetrization $X_{(\alpha\beta)} = \frac{1}{2}\left(X_{\alpha\beta} + X_{\beta\alpha}\right)$.
\end{proposition}
\begin{proof}
  See the proof of Proposition 4.7.2 in \cite{giorgiWaveEquationsEstimates2024}. 
\end{proof}

\begin{remark}
  It will be convenient to introduce the notation
  \begin{equation}
    \label{eq:KCurrent:without-Riem}
    \KCurrent{\MorawetzVF,\MorawetzLagrangeCorr,\MorawetzOneForm}_0[\psi]
    \vcentcolon= \KCurrent{\MorawetzVF, \MorawetzLagrangeCorr, \MorawetzOneForm}[\psi]
    - \MorawetzVF^{\mu}\HorCovDeriv^{\nu}\psi^a\HorRiem_{ab\nu\mu}\psi^b,
  \end{equation}
  to be analogous to the standard $K$-current for scalar wave
  equations.
\end{remark}

\subsection{Basic deformation tensors}
\label{sec:basic-deformation-tensors}

In this section, we compute the deformation tensors of $e_3$ and $e_4$
in the principal outgoing and principal ingoing null frames.
\begin{lemma}
  \label{lemma:e3-e4-deformation-tensors}
  Let $(e_3, e_4)$ be a null frame. Then, 
  \begin{gather}
    \label{eq:e3-deformation-tensors}
    \DeformationTensor[_{44}]{e_3}
    = -4\omega, \qquad
    \DeformationTensor[_{34}]{e_3}
    = 2\omegaBar, \qquad
    \DeformationTensor[_{33}]{e_3}
    =0, \\
    \DeformationTensor[_{4a}]{e_3}
    = \etaBar_a-\zeta_a,\qquad
    \DeformationTensor[_{3a}]{e_3}
    = \xiBar_a \qquad
    \DeformationTensor[_{ab}]{e_3}
    = \frac{1}{2}\left(\Trace \chiBar\delta_{ab} + 2\chiBarTF_{ab}\right),
  \end{gather}
  and
  \begin{gather}
    \label{eq:e4-deformation-tensors}
    \DeformationTensor[_{33}]{e_4}
    = -4\omegaBar, \qquad
    \DeformationTensor[_{34}]{e_4}
    = 2\omega, \qquad
    \DeformationTensor[_{44}]{e_4}
    =0, \\
    \DeformationTensor[_{3a}]{e_4}
    = \eta_a+\zeta_a,\qquad
    \DeformationTensor[_{4a}]{e_4}
    = \xi_a \qquad
    \DeformationTensor[_{ab}]{e_4}
    = \frac{1}{2}\left(\Trace \chi\delta_{ab} + 2\chiTF_{ab}\right).
  \end{gather}
  In particular, if $(e_3,e_4)=\left(e_3^{(out)}, e_4^{(out)}\right)$ denotes
  the principal outgoing null frame we have 
  \begin{gather}
    \label{eq:e3-deformation-tensors:principal-outgoing}
    \DeformationTensor[_{44}]{e_3}
    = 0, \qquad
    \DeformationTensor[_{34}]{e_3}
    = \partial_r\left(\frac{\Delta}{\abs*{q}^2}\right), \qquad
    \DeformationTensor[_{33}]{e_3}
    =0, \\
    \DeformationTensor[_{4a}]{e_3}
    = -\Re \frac{a\overline{q}}{\abs*{q}^2}\CCOneFormJ_a,\qquad
    \DeformationTensor[_{3a}]{e_3}
    =  0, \qquad
    \DeformationTensor[_{ab}]{e_3}
    = -\frac{r\Delta}{\abs*{q}^4}\delta_{ab} ,
  \end{gather}
  and
  \begin{gather}
    \label{eq:e4-deformation-tensors:principal-outgoing}
    \DeformationTensor[_{33}]{e_4}
    = -2\partial_r\left(\frac{\Delta}{\abs*{q}^2}\right), \qquad
    \DeformationTensor[_{34}]{e_4}
    = 0, \qquad
    \DeformationTensor[_{44}]{e_4}
    =0, \\
    \DeformationTensor[_{4a}]{e_4}
    = 0,\qquad
    \DeformationTensor[_{3a}]{e_4}
    = \frac{ar}{\abs*{q}^2}\Re\left(\CCOneFormJ_a\right), \qquad
    \DeformationTensor[_{ab}]{e_4}
    = \frac{r}{\abs*{q}^2}\delta_{ab}.
  \end{gather}

  In particular, if $(e_3,e_4)=\left(e_3^{(in)}, e_4^{(in)}\right)$ denotes
  the principal ingoing null frame we have 
  \begin{gather}
    \label{eq:e3-deformation-tensors:principal-ingoing}
    \DeformationTensor[_{44}]{e_3}
    = 2\partial_r\left(\frac{\Delta}{\abs*{q}^2}\right), \qquad
    \DeformationTensor[_{34}]{e_3}
    = 0, \qquad
    \DeformationTensor[_{33}]{e_3}
    =0, \\
    \DeformationTensor[_{4a}]{e_3}
    = -\frac{ar}{\abs*{q}^2}\Re \CCOneFormJ_a,\qquad
    \DeformationTensor[_{3a}]{e_3}
    =  0, \qquad
    \DeformationTensor[_{ab}]{e_3}
    = -\frac{r}{\abs*{q}^2}\delta_{ab} ,
  \end{gather}
  and
  \begin{gather}
    \label{eq:e4-deformation-tensors:principal-ingoing}
    \DeformationTensor[_{33}]{e_4}
    = 0, \qquad
    \DeformationTensor[_{34}]{e_4}
    = -\partial_r\left(\frac{\Delta}{\abs*{q}^2}\right), \qquad
    \DeformationTensor[_{44}]{e_4}
    =0, \\
    \DeformationTensor[_{4a}]{e_4}
    = 0,\qquad
    \DeformationTensor[_{3a}]{e_4}
    = \frac{ar}{\abs*{q}^2}\Re\left(\CCOneFormJ_a\right), \qquad
    \DeformationTensor[_{ab}]{e_4}
    = \frac{r\Delta}{\abs*{q}^4}\delta_{ab}.
  \end{gather}
\end{lemma}

\begin{proof}
  Denoting $(e_4,e_3) = (\eout_4,\eout_3)$, we can calculate that
  \begin{align*}
    \DeformationTensor[_{44}]{e_3}
    &=  \Metric\left(\CovariantDeriv_4e_3,e_4\right)
      = -4\omega,\\
    \DeformationTensor[_{34}]{e_3}
    &= \frac{1}{2}\Metric\left(\CovariantDeriv_3e_3,e_4\right)
      = 4\omegaBar,\\
    \DeformationTensor[_{33}]{e_3}
    &= \Metric\left(\CovariantDeriv_3e_3,e_3\right)
      =0,\\
    \DeformationTensor[_{4a}]{e_3}
    &= \frac{1}{2}\left( \Metric\left(\CovariantDeriv_4e_3,e_a\right)
      + \Metric\left(\CovariantDeriv_ae_3,e_4\right) \right)
      = \etaBar_a-\zeta_a,
    \\
    \DeformationTensor[_{3a}]{e_3}
    &= \frac{1}{2}\left( \Metric\left(\CovariantDeriv_3e_3,e_a\right)
      + \Metric\left(\CovariantDeriv_ae_3,e_3\right) \right)
      =\xiBar_a,
    \\
    \DeformationTensor[_{ab}]{e_3}
    &= \frac{1}{2}\left( \Metric\left(\CovariantDeriv_ae_3, e_b\right)
      + \Metric\left(\CovariantDeriv_be_3, e_a\right) \right)
      = \frac{1}{2}\left( \chiBar_{ab} + \chiBar_{ba} \right).
  \end{align*}
Denoting $(e_4,e_3) = (\ein_4,\ein_3)$, we can similarly calculate
  \begin{align*}
    \DeformationTensor[_{33}]{e_4}
    &= \Metric\left(\CovariantDeriv_3e_4,e_3\right) = -4\omega,\\
    \DeformationTensor[_{34}]{e_4}
    &= \frac{1}{2}\Metric\left(\CovariantDeriv_4e_4,e_3\right) = 4\omega,\\
    \DeformationTensor[_{44}]{e_4}
    &= \Metric\left(\CovariantDeriv_4e_4,e_4\right)
      = 0,\\
    \DeformationTensor[_{3a}]{e_4}
    &= \frac{1}{2}\left(\Metric\left(\CovariantDeriv_3e_4,e_a\right)
      + \Metric\left(\CovariantDeriv_ae_4,e_3\right)\right)
      = \eta_a + \zeta_a,\\
    \DeformationTensor[_{4a}]{e_4}
    &=\frac{1}{2}\left(\Metric\left(\CovariantDeriv_4e_4,e_a\right)
      + \Metric\left(\CovariantDeriv_ae_4,e_4\right)\right)
      = \xi_a,\\
    \DeformationTensor[_{ab}]{e_4}
    &= \frac{1}{2}\left(\Metric\left(\CovariantDeriv_ae_4, e_b\right)
      + \Metric\left(\CovariantDeriv_be_4, e_a\right)\right)
      = \frac{1}{2}\left(\chi_{ab} + \chi_{ba}\right)
      = \frac{1}{2}\left(\Trace \chi\delta_{ab} + 2\chiTF_{ab}\right).
  \end{align*}
  Plugging in the exact values for the Ricci coefficients from in the
  principal outgoing and principal ingoing null frame from
  \zcref[cap]{lemma:Kerr:outgoing-PG:Ric-and-curvature} and
  \zcref[cap]{lemma:Kerr:ingoing-PG:Ric-and-curvature} respectively concludes
  the proof of \zcref[cap]{lemma:e3-e4-deformation-tensors}.
\end{proof}

\subsection{Decomposition of the wave operator in frames}

We start by decomposing the wave operator in terms of
null frames.
\begin{lemma}
  \label{coro:wave-operator-null-frame-decomp:wave-2}
  The wave operator for $\psi\in \realHorkTensor{k}$ is given by
  \begin{equation}
    \label{eq:wave-operator-null-frame-decomp:wave-2}
    \begin{split}
      \WaveOpHork{2}\psi
      ={}& - \nabla_4\nabla_3\psi - \frac{1}{2}\Trace\chiBar\nabla_4\psi
           + \left(2\omega - \frac{1}{2}\Trace\chi\right)\nabla_3\psi
           + \Laplace_2\psi
           + 2\etaBar\cdot\nabla\psi
          + 2\ImagUnit \left(\LeftDual{\rho} - \eta\wedge \etaBar\right)\psi,
    \end{split}
  \end{equation}
    where we recall that $\Laplace_k = \nabla^a\nabla_a$ denotes the horizontal Laplacian
  acting on $\realHorkTensor{k}$.
\end{lemma}
\begin{proof}
  See proof of Lemma 4.7.5 in
  \cite{giorgiWaveEquationsEstimates2024}.
\end{proof}

It will also be convenient to have the following expressions for the wave operator.
\begin{lemma}
  \label{lemma:wave-using-THat-RHat}
  In exact Kerr-de Sitter we have 
  \begin{equation}
    \label{eq:wave-using-THat-RHat}
    \begin{split}
      \abs*{q}^2\WaveOpHork{k}\psi
      ={}& \frac{(r^2+a^2)^2}{\Delta}\left(-\nabla_{\HawkingVF}\nabla_{\HawkingVF}\psi + \nabla_{\HprVF}\nabla_{\HprVF}\psi\right)
           + 2r\nabla_{\HprVF}\psi
           + \abs*{q}^2\LaplaceHor_k\psi
           + \abs*{q}^2\left(\eta + \etaBar\right)\cdot\nabla\psi.
    \end{split}    
  \end{equation}
\end{lemma}
\begin{proof}
  The proof of \zcref[noname]{eq:wave-using-THat-RHat} follows identically as
  in Section C.7 of \cite{giorgiWaveEquationsEstimates2024}. 
\end{proof}

\subsection{General commutation formulas}\label{sec:general-commutation-formulas}

We review some basic commutation formulas.
\begin{lemma}
  \label{lemma:commutation-formula:basic}
  Let $U_A=U_{a_1\cdots a_k}$ be a general $k$-horizontal tensor field. We have
    \begin{align}
      \label{eq:commutation-formula:nabla3-nablab}
        [\nabla_{3},\nabla_{b}]U_A &
        = -\chiBar_{bc}\nabla_{c}U_{A}
             + \left(\eta_{b}-\zeta_{b}\right)\nabla_3U_{A}
             + \xiBar_{b}\nabla_{4}U_{A}
            + \sum_{i=1}^{k}\left(\volFormHor_{a_{i}c}\LeftDual{\betaBar}_{b} + \frac{1}{2}\BHorizontalCurv_{a_{i}c3b}\right)\tensor[]{U}{_{a_{1}\cdots}^{b}_{{\cdots a_{k}}}},\\
      \label{eq:commutation-formula:nabla4-nablab}
        [\nabla_{4},\nabla_{b}]U_A &
        = -\chi_{bc}\nabla_{c}U_{A}
             + \left(\etaBar_{b}+\zeta_{b}\right)\nabla_4U_{A}
             + \xi_{b}\nabla_{3}U_{A}
            + \sum_{i=1}^{k}\left(\volFormHor_{a_{i}c}\LeftDual{\beta}_{b} + \frac{1}{2}\BHorizontalCurv_{a_{i}c4b}\right)\tensor[]{U}{_{a_{1}\cdots}^{b}_{{\cdots a_{k}}}},\\
      \label{eq:commutation-formula:nabla4-nabla3}
        [\nabla_4,\nabla_3]U_{A} &
        = 2\left(\etaBar_b-\eta_b\right)\nabla_{b}U_{A}
             + 2\omega\nabla_{3}U_{A}
             - 2\omegaBar\nabla_{4}U_{A}
            + \sum_{i=1}^{k}\left(
             -\volFormHor_{a_{i}b}\LeftDual{\rho}
             + \frac{1}{2}\BHorizontalCurv_{a_{i}b43}
             \right) \tensor[]{U}{_{a_{1}\cdots}^{b}_{{\cdots a_{k}}}}.
    \end{align}
\end{lemma}

\begin{proof}
  The proof is identical to that of Lemma 2.2.7 in
  \cite{giorgiWaveEquationsEstimates2024} (see \zcref[cap]{remark:Weyl-Riem-difference}).
\end{proof}

If we plug in the values for $\BHorizontalCurv$ from
\zcref[cap]{prop:B-hor:components} and use the identities in
\zcref[noname]{eq:global-frame:vanishing-qtys}, we have the following
corollary of \zcref[cap]{lemma:commutation-formula:basic}.

\begin{corollary}
  \label{corollary:commutation-formula:B-applied}
  Let $U_A=U_{a_1\cdots a_k}\in \realHorkTensor{k}(\Manifold)$. Then
  the following hold in the global frame as defined in
  \zcref[cap]{def:global null frame}:
    \begin{equation}
      \label{eq:commutation-formula:B-applied:e3eb}
      \begin{split}
        \left[ \nabla_3,\nabla_b \right]U_A
        ={}& - \frac{1}{2}\left(\Trace \underline{\chi} \nabla_bU_A + \aTrace{\underline{\chi}}\LeftDual{\nabla}_bU_A\right)
             + \left( \eta_b-\zeta_b \right)\nabla_3 U_A \\
           & + \frac{1}{2}\sum_{i=1}^k\left( \delta_{a_ib}\Trace\underline{\chi}
             + \volFormHor_{ba_i}\aTrace{\underline{\chi}} \right)\eta_c\tensor[]{U}{_{a_1\cdots}^c_{\cdots a_k}}
           -\frac{1}{2}\sum_{i=1}^k\eta_{a_i} \left( \Trace\underline{\chi}U_{a_1\cdots b\cdots a_k}
             +\aTrace{\underline{\chi}}U_{a_1\cdots b\cdots a_k} \right),
      \end{split}
    \end{equation}
    \begin{equation}
      \label{eq:commutation-formula:B-applied:e4eb}
      \begin{split}
        \left[ \nabla_3,\nabla_b \right]U_A
        ={}& - \frac{1}{2}\left(\Trace \chi \nabla_bU_A + \aTrace{\chi}\LeftDual{\nabla}_bU_A\right)
             + \left( \underline{\eta}_b+\zeta_b \right)\nabla_4 U_A\\
           &+ \frac{1}{2}\sum_{i=1}^k\left( \delta_{a_ib}\Trace\chi
             + \volFormHor_{ba_i}\aTrace{\chi} \right)\underline{\eta}_c\tensor[]{U}{_{a_1\cdots}^c_{\cdots a_k}}
           -\frac{1}{2}\sum_{i=1}^k\underline{\eta}_{a_i} \left( \Trace\underline{\chi}U_{a_1\cdots b\cdots a_k}
             +\aTrace{\chi}U_{a_1\cdots b\cdots a_k} \right),
      \end{split}
    \end{equation}
    \begin{equation}
      \label{eq:commutation-formula:B-applied:e3e4}
      \begin{split}
        \left[ \nabla_4,\nabla_3 \right]U_A
        ={}& 2\left(\underline{\eta}_b-\eta_b\right)\nabla_bU_A
             + 2\omega\nabla_3U_A
             - 2\underline{\omega}\nabla_3U_A
           + 2\sum_{i=1}^k\left( \eta_{a_i}\underline{\eta}_b
             - \underline{\eta}_{a_i}\eta_b
             - \volFormHor_{a_ib}\LeftDual{\rho} \right)\tensor[]{U}{_{a_1\cdots}^c_{\cdots a_k}}.
      \end{split}
    \end{equation}
\end{corollary}

We specialize to the case of $\realHorkTensor{0}$,
$\realHorkTensor{1}$, and $\realHorkTensor{2}$ in the following
corollary.
\begin{corollary}
  \label{coro:commutation-formula:realHorkTensor}
  The following commutation formulas hold in the global frame defined
  in \zcref[cap]{def:global null frame}. Given $f\in \realHorkTensor{0}$, $F\in \realHorkTensor{1}$ and $u\in \realHorkTensor{2}$ we have
    \begin{align*}
      [\nabla_3, \nabla_a]f
      ={}& -\frac{1}{2}\left(
           \Trace\chiBar \nabla_af + \aTrace{\chiBar}\LeftDual{\nabla}_af
           \right)
           + \left(\eta_a-\zeta_a\right)\nabla_3f,\\
      [\nabla_4, \nabla_a]f
      ={}& -\frac{1}{2}\left(
           \Trace\chi\nabla_af + \aTrace{\chi}\LeftDual{\nabla}_af
           \right)
           + \left(\etaBar_a+\zeta_a\right)\nabla_4f,\\
      [\nabla_4,\nabla_3]f
      ={}& 2\left(\etaBar-\eta\right)\cdot \nabla f
           + 2\omega \nabla_3f
           - 2\omegaBar\nabla_4f,
    \end{align*}
    \begin{align*}
      [\nabla_3,\nabla_a]F_b
      ={}& -\frac{1}{2}\Trace\chiBar\left(
           \nabla_aF_b+\eta_bF_a-\delta_{ab}\eta\cdot F
           \right)
           -\frac{1}{2}\aTrace{\chiBar}\left(
           \LeftDual{\nabla}_aF_b + \eta_b\LeftDual{u}_a - \volFormHor_{ab}\eta\cdot F
           \right)
          + \left(\eta_a-\zeta_a\right)\nabla_{3}F_b,\\
      [\nabla_4,\nabla_a]F_b
      ={}& -\frac{1}{2}\Trace\chi\left(
           \nabla_aF_b + \etaBar_bF_a - \delta_{ab}\etaBar\cdot F
           \right)
           -\frac{1}{2}\aTrace{\chi}\left(
           \LeftDual{\nabla_aF_b + \etaBar_b\LeftDual{F}_a-\volFormHor_{ab}\etaBar\cdot F}
           \right)
          + \left(
           \etaBar_a + \zeta_a 
           \right)\nabla_4F_b,\\
    \end{align*}
    \begin{align*}
      \squareBrace*{\nabla_3, \nabla_a}u_{bc}
      ={}& - \frac{1}{2}\Trace\chiBar\left(
           \nabla_a u_{bc}
           + \eta_bu_{ac}
           + \eta_{c}u_{ab}
           - \delta_{ab}\left(\eta\cdot u\right)_c
           - \delta_{ac}\left(\eta\cdot u\right)_b
           \right)\\
         & - \frac{1}{2}\aTrace{\chiBar}\left(
           \LeftDual{\nabla}_au_{bc}
           + \eta_b\LeftDual{u}_{ac}
           + \eta_c\LeftDual{u}_{ab}
           - \volFormHor_{ab}\left(\eta\cdot u\right)_c
           - \volFormHor_{ac}\left(\eta\cdot u\right)_b
           \right)
          + \left(\eta-\zeta\right)_a\nabla_3u_{bc},\\
      \squareBrace*{\nabla_4,\nabla_a}u_{bc}
      ={}& - \frac{1}{2}\Trace\chi\left(
           \nabla_au_{hc}
           + \etaBar_bu_{ac}
           + \etaBar_cu_{ab}
           - \delta_{ab}\left(\etaBar \cdot u\right)_c
           - \delta_{ac}\left(\etaBar\cdot u\right)_b
           \right)\\
         &- \frac{1}{2}\aTrace{\chi}\left(
           \LeftDual{\nabla}_au_{bc}
           + \etaBar_b\LeftDual{u}_{ac}
           + \etaBar_c\LeftDual{u}_{ab}
           - \volFormHor_{ab}\left(\etaBar\cdot u\right)_c
           - \volFormHor_{ac}\left(\etaBar \cdot u\right)_b
           \right)
          + \left(\etaBar + \zeta\right)_a\ConformalInvDeriv_4u_{bc},\\
      \squareBrace*{\nabla_4,\nabla_3}u_{ab}
      ={}& 2\omega\nabla_3u_{ab}
           - 2\omegaBar\nabla_4 u_{ab}
           + 2\left(\etaBar-\eta\right)_c\nabla_cu_{ab}
           + 4\eta\SymTracelessTensorProd(\etaBar\cdot u)
           - 4 \etaBar\SymTracelessTensorProd(\eta\cdot u)
           - 4\LeftDual{\rho}\LeftDual{u}_{ab}.
    \end{align*}
\end{corollary}

\begin{proof}
  The proof is identical to that of Corollary 2.2.9 from
  \cite{giorgiWaveEquationsEstimates2024} using also
  \zcref[noname]{eq:global-frame:vanishing-qtys}.
\end{proof}

We also immediately have the following corollary.
\begin{corollary}
  \label{coro:commutation-formula:div-curl-realHorkTensor}
  The following commutation formulas hold in the global frame defined
  in \zcref[cap]{def:global null frame}. Given $F\in \realHorkTensor{1}$ and $U\in \realHorkTensor{2}$ we have
    \begin{align}
      \left[\nabla_3,\nabla\cdot \right]F
      ={}& -\frac{1}{2}\Trace\chiBar\left(\nabla\cdot F - \eta\cdot F\right)
           + \frac{1}{2}\aTrace{\chi}\left(
           \nabla\cdot\LeftDual{F}
           -\eta\cdot\LeftDual{F}
           \right)
           + \left(\eta-\zeta\right)\cdot \nabla_3F,
            \label{eq:commutation-formula:1-tensor:3-div}\\
      \left[\nabla_4,\nabla\cdot \right]F
      ={}& -\frac{1}{2}\Trace\chi\left(\nabla\cdot F - \etaBar\cdot F\right)
           + \frac{1}{2}\aTrace{\chi}\left(
           \nabla\cdot\LeftDual{F}
           -\etaBar\cdot\LeftDual{F}
           \right)
           + \left(\etaBar+\zeta\right)\cdot \nabla_4F
           ,       \label{eq:commutation-formula:1-tensor:4-div} \\    
      \left[\nabla_3, \nabla\SymTracelessTensorProd\right]F
      ={}& - \frac{1}{2}\Trace\chiBar\left(\nabla\SymTracelessTensorProd F + \eta\SymTracelessTensorProd F\right)
           - \frac{1}{2}\aTrace{\chiBar}\LeftDual{\left(\nabla\SymTracelessTensorProd F + \etaBar\SymTracelessTensorProd F\right)}
           + \left(\eta-\zeta\right)\cdot \nabla_3 F
           ,       \label{eq:commutation-formula:1-tensor:3-curl}\\
      \left[\nabla_4, \nabla\SymTracelessTensorProd\right]F
      ={}& - \frac{1}{2}\Trace\chi\left(\nabla\SymTracelessTensorProd F + \etaBar\SymTracelessTensorProd F\right)
           - \frac{1}{2}\aTrace{\chi}\LeftDual{\left(\nabla\SymTracelessTensorProd F + \etaBar\SymTracelessTensorProd F\right)}
           + \left(\etaBar+\zeta\right)\cdot \nabla_4 F
           ,       \label{eq:commutation-formula:1-tensor:4-curl}
      \\ \left[\nabla_3, \nabla\cdot\right]U
      ={}& - \frac{1}{2}\Trace\chiBar\left(
           \nabla\cdot U - 2\eta\cdot U
           \right)
           + \frac{1}{2}\aTrace{\chiBar}\left(\nabla\cdot \LeftDual{U}  - 2\eta \cdot \LeftDual{U}\right)
           + \left(\eta - \zeta\right)\cdot \nabla_3U
           , \label{eq:commutation-formula:2-tensor:3-div}\\
      \left[\nabla_4, \nabla\cdot \right]U
      ={}& - \frac{1}{2}\Trace\chi\left(\nabla\cdot U - 2\etaBar\cdot U\right)
           + \frac{1}{2}\aTrace{\chi}\left(\nabla\cdot \LeftDual{U} - 2\etaBar\cdot \LeftDual{U}\right)
           + \left(\etaBar + \zeta\right)\cdot \nabla_4U. \label{eq:commutation-formula:2-tensor:4-div}
    \end{align}
\end{corollary}

The commutators $\left[ \,^{(c)}\nab_3,\,^{(c)}\nab_a \right]$ and $\left[ \,^{(c)}\nab_4,\,^{(c)}\nab_a \right]$ will respectively require the null structure formulas for $\nab_3\zeta + 2 \nab \underline{\om}$ and $\nab_4\zeta + 2 \nab \om$. However, the curvature components appearing in those null structure equations are $\R_{a334}$ and $\R_{a443}$ which coincide with the Weyl tensor in the $\La\neq 0$ case (see \zcref[noname]{link Riemann Weyl}). Therefore, every commutator involving one horizontal conformal derivative can be directly taken from Lemma 4.2.2 in \cite{giorgiWaveEquationsEstimates2024}. However, the commutator $\left[ \,^{(c)}\nab_3,\,^{(c)}\nab_4 \right]$ will require the null structure for $\nabla_4\underline{\om} + \nabla_3\omega$, namely
\begin{align*}
\nabla_4\underline{\om} + \nabla_3\omega & = 4 \om \underline{\om}  + \xi \cdot \underline{\xi} + \zeta\cdot (\eta- \underline{\eta})  - \eta\cdot \underline{\eta}  + \frac{1}{4} \R_{3434}.
\end{align*}
However from \zcref[noname]{link Riemann Weyl} we find $\frac{1}{4} \R_{3434}=\rho-\frac{\La}{3}$. Therefore, the commutation formulas for $\left[ \,^{(c)}\nab_3,\,^{(c)}\nab_4 \right]$ applied to any kind of tensors are given by the formulas from Lemma 4.2.2 in \cite{giorgiWaveEquationsEstimates2024} with $\rho$ replaced by $\rho-\frac{\La}{3}$. We thus have the following corollary for conformally invariant derivatives. 

\begin{corollary}
  \label{lemma:commutation-formula:conformally-invariant-derivaties}
  The following commutation formulas involving the conformally
  invariant derivatives hold in the global frame defined
  in \zcref[cap]{def:global null frame}. Given $h\in \realHorkTensor{0}(\Complex)$, $F\in\realHorkTensor{1}(\Complex)$ and $U\in\realHorkTensor{2}(\Complex)$ all $s$-conformally invariant we have  
    \begin{align}
      \left[\ConformalInvDeriv_4,\ConformalComplexDeriv\right]h
      ={}& -\frac{1}{2}\Trace X\ConformalComplexDeriv h
           + \HBar \ConformalInvDeriv_{4}h  + s\frac{1}{2}\Trace X\HBar h ,
           \label{eq:commutation-formula:conformally-invariant-derivaties:scalar:4-A}\\
      \left[\ConformalInvDeriv_3,\ConformalComplexDeriv\right]h
      ={}& -\frac{1}{2}\Trace \XBar\ConformalComplexDeriv h
           + H \ConformalInvDeriv_{3}h
            - s \frac{1}{2}\Trace \XBar H h , 
           \label{eq:commutation-formula:conformally-invariant-derivaties:scalar:3-A}\\
      \left[\ConformalInvDeriv_3,\ConformalInvDeriv_4\right]h
      ={}& 2 \left(\eta-\etaBar\right)\cdot \ConformalInvDeriv h
           + 2s \left(\rho - \frac{\Lambda}{3} -\eta\cdot \etaBar \right)h, \label{eq:commutation-formula:conformally-invariant-derivaties:scalar:3-4}
       \\ \squareBrace*{\ConformalInvDeriv_4, \overline{\ConformalComplexDeriv\cdot}} F
        ={}& -\frac{1}{2}\overline{\Trace X}\left(
             \overline{\ConformalComplexDeriv}\cdot F
             - (s+1)\overline{\HBar}\cdot F
             \right)
             + \overline{\HBar}\cdot \ConformalInvDeriv_4F, \label{eq:commutation-formula:conformally-invariant-derivatives:1-tensor:4-div} \\
        \squareBrace*{\ConformalInvDeriv_3, \overline{\ConformalComplexDeriv\cdot}} F
        ={}& -\frac{1}{2}\overline{\Trace \XBar}\left(
             \overline{\ConformalComplexDeriv}\cdot F
             + (s-1)\overline{\HBar}\cdot F
             \right)
             + \overline{H}\cdot \ConformalInvDeriv_3F,
      \\ \left[\ConformalInvDeriv_4, \overline{\ConformalComplexDeriv\cdot}\right]U
      ={}& -\frac{1}{2}\overline{\Trace X}\left(
           \overline{\ConformalComplexDeriv}\cdot U
           - (s+2)\overline{\HBar}\cdot U
           \right)
           + \overline{\HBar}\cdot \ConformalInvDeriv_4U ,  \label{eq:eq:commutation-formula:conformally-invariant-derivaties:tensor:4-div}\\
      \left[\ConformalInvDeriv_3,\overline{\ConformalComplexDeriv}\cdot\right]U
      ={}& - \frac{1}{2}\overline{\Trace \XBar}\left(
           \overline{\ConformalComplexDeriv}\cdot U
           + (s-2)\overline{H}\cdot U
           \right)
           + \overline{H}\cdot \ConformalInvDeriv_3U ,
           \label{eq:eq:commutation-formula:conformally-invariant-derivaties:tensor:3-div}\\
           \left[\ConformalInvDeriv_{3},\ConformalInvDeriv_{4}\right] U
      ={}& 2\left(\eta-\etaBar\right)\cdot \ConformalInvDeriv U
           + 2 s\left(
           \rho- \frac{\Lambda}{3} -\eta\cdot \etaBar
           \right)U
           + 4\ImagUnit \left(-\LeftDual{\rho} + \eta\wedge \etaBar\right)U.
           \label{eq:commutation-formula:conformally-invariant-derivaties:tensor:3-4}
    \end{align}
\end{corollary}

\subsection{Commutation properties of the wave operator and the horizontal Laplacian}
\label{sec:wave-operator-basic-commutation-properties}

In this section, we review some basic commutation properties of the
horizontal wave operator. 

\begin{lemma}
  \label{lemma:wave-operator-basic-commutation}
  Let $\WaveOpHork{k}$ denote the wave operator on
  $\realHorkTensor{k}$. Then,
  \begin{equation}
    \label{eq:wave-operator-basic-commutation:r-e4}
    \begin{split}
      \squareBrace*{r\nabla_4,\WaveOpHork{2}}\psi
      ={}& -\nabla_4\nabla_4\psi
           - r\left(\frac{1}{2}\Trace\chi - 2\omega\right)\WaveOpHork{2}\psi
           - r\left(\frac{1}{2}\Trace \chi + 2\omega\right)\Laplace_2\psi
         \\& + O(r^{-2})\frakWeightedDeriv^{\le 1}\psi
           + O(r^{-3})\frakWeightedDeriv^{\le 2}\psi
           + r \nabla_3\left(\xi\cdot\HorCovDeriv_a\psi\right),
    \end{split}
  \end{equation}
  and
  \begin{align}
    \label{eq:wave-operator-basic-commutation:Lie-Killing}
      \left( \HorLieDeriv_{\KillT} \WaveOpHork{2} - \WaveOpHork{2}\HorLieDeriv_{\KillT} \right)\psi_{ab}
      &= 0 , &
      \left( \HorLieDeriv_{\KillPhi} \WaveOpHork{2} - \WaveOpHork{2}\HorLieDeriv_{\KillPhi} \right)\psi_{ab} & =  0, 
  \end{align}
  and
  \begin{equation}
    \label{eq:wave-operator-basic-commutation:Hodge-op}
    \begin{split}
      \abs*{q}\HodgeOp{2}\WaveOpHork{2}\psi
      - \WaveOpHork{1}\abs*{q}\HodgeOp{2}\psi
      ={}& 3\horProj{K}\abs*{q}\HodgeOp{2}\psi
           - \frac{2a(1+\gamma)\cos\theta}{\abs*{q}}\LeftDual{\HodgeOp{2}}\HorLieDeriv_{\KillT}\psi
           - \abs*{q}(\eta+\etaBar)\cdot\WaveOpHork{2}\psi\\
          & + O(ar^{-2})\frakWeightedDeriv^{\le 1}\psi
            + O(ar^{-3})\frakWeightedDeriv^{\le 2}\psi,%
    \end{split} 
  \end{equation}
  where $\horProj{K}$ is as defined in \zcref[noname]{eq:h-K:def}.
\end{lemma}

\begin{proof}
  For \zcref[noname]{eq:wave-operator-basic-commutation:r-e4}, see the proof
  for Lemma 4.7.11 in \cite{giorgiWaveEquationsEstimates2024},
  observing that the only null structure equation or null Bianchi
  equation used was the equation for $\nabla_4\Trace \chi$, which does
  not change in the presence of a cosmological constant since it involves the curvature coefficient $\R_{a44b}$ (see \zcref[noname]{link Riemann Weyl}). For \zcref[noname]{eq:wave-operator-basic-commutation:Lie-Killing}, see the
  proof of Corollary 4.3.4 in \cite{giorgiWaveEquationsEstimates2024}. For \zcref[noname]{eq:wave-operator-basic-commutation:Hodge-op}, see the
  proof of Lemma 4.7.13 in \cite{giorgiWaveEquationsEstimates2024}. We
  note that the only change is that for \KdS, 
  \begin{equation*}
    \aTrace{\chi}e_3 + \aTrace{\chi}e_4 + 2(\eta+\etaBar)\cdot\LeftDual{\nabla}
    = \frac{4a(1+\gamma)\cos\theta}{\abs*{q}^2}\KillT. 
  \end{equation*}
  This concludes the proof of \zcref[cap]{lemma:wave-operator-basic-commutation}.
\end{proof}

We also give the following basic commutators with the horizontal
Laplacian.
\begin{lemma}
  \label{lemma:commutation:horizontal-laplacian}
  The following commutations formulas hold true for a $2$-tensor
  $\psi\in \realHorkTensor{2}$:
  \begin{equation*}
    \begin{split}
      \left[\nabla_3,\abs*{q}^2\Laplace_2\right]\psi
      ={}& \left(\eta-\zeta\right)\cdot\abs*{q}^2\nabla_3\nabla\psi
           + \left(\eta-\zeta\right)\cdot\abs*{q}^2\nabla\nabla_3\psi
           + \Divergence\left(\eta-\zeta\right)\abs*{q}^2\nabla_3\psi\\
      &- \frac{1}{2}\abs*{q}^2\left(\nabla\Trace\chiBar\cdot\nabla\psi + \nabla\aTrace{\chiBar}\cdot\LeftDual{\nabla}\psi\right)
      + O(ar^{-4})\frakWeightedDeriv^{\le 1}\psi,\\
      \left[\nabla_4,\abs*{q}^2\Laplace_2\right]
      ={}& \left(\etaBar+\zeta\right)\cdot\abs*{q}^2\nabla_4\nabla\psi
           + \left(\etaBar+\zeta\right)\cdot\abs*{q}^2\nabla\nabla_4\psi
           + \Divergence\left(\etaBar+\zeta\right)\abs*{q}^2\nabla_4\psi\\
      &- \frac{1}{2}\abs*{q}^2\left(\nabla\Trace\chi\cdot\nabla\psi + \nabla\aTrace{\chi}\cdot\LeftDual{\nabla}\psi\right)
      + O(ar^{-4})\frakWeightedDeriv^{\le 1}\psi
      + r^2\nabla_3(\xi\cdot\HorCovDeriv_a \psi).
    \end{split}
  \end{equation*}
\end{lemma}
\begin{proof}
  The proof is identical to that of Lemma 4.7.10 in
  \cite{giorgiWaveEquationsEstimates2024}.
\end{proof}

\section{Energy estimate}
\label{sec:Killing-energy-estimate}

The main goal of this section will be to prove the following energy estimate for solutions of the equation
\begin{equation}
  \label{eq:model-problem-gRW}
  \WaveOpHork{2}\psi - V\psi = \ImagUnit\frac{4a(1+\gamma)\cos\theta}{\abs*{q}^{2}}\nabla_{\KillT}\psi
  + N
  ,\qquad V =  \frac{4\Delta}{(r^2+a^2)|q|^2} + 2\Lambda.
\end{equation}

\begin{proposition}
  \label{prop:Killing-energy-estimate}
  For $\psi\in \realHorkTensor{2}$ a solution of \zcref[noname]{eq:model-problem-gRW}, 
  \begin{equation}
    \label{eq:Killing-energy-estimate}
    \begin{split}
      &\EnergyFlux_{\operatorname{deg}}\left[\psi\right](\tau_2) + \int_{\mathcal{A}(\tau_1,\tau_2)}\abs*{\nabla_{4}\psi}^2 + \int_{\SigmaStar(\tau_1,\tau_2)}\abs*{\nabla_{3}\psi}^2    \\
      \lesssim{}& \EnergyFlux_{\operatorname{deg}}[\psi](\tau_1)
                  + \delta_{\Horizon}\left(
                  \EnergyFlux_{r\le r_{\EventHorizon}}[\psi](\tau_2)
                  + \EnergyFlux_{r\ge r_{\CosmologicalHorizon}}[\psi](\tau_2)
                  + \SpacelikeFlux_{\mathcal{A}}[\psi](\tau_1, \tau_2)
                  + \SpacelikeFlux_{\SigmaStar}[\psi](\tau_1,\tau_2)
                  \right)\\
      &+ \frac{a}{M}\MorNorm[\psi](\tau_1,\tau_2)
        + \abs*{\int_{\Manifold(\tau_1,\tau_2)}\nabla_{\TAlmostKilling}\psi\cdot N}
        + \int_{\Manifold(\tau_1,\tau_2)}\abs*{N}^2,
    \end{split}
  \end{equation}
  where $\TAlmostKilling$ is defined in \zcref[cap]{def:TAlmostKilling}, and where the implicit constant is independent of $\Lambda$. 
\end{proposition}

\subsection{The almost Killing vectorfield \texorpdfstring{$\TAlmostKilling$}{T-delta}}
\label{sec:TFixer:construction}

There are no globally timelike Killing vectorfields on \KdS{}
backgrounds. However, we can define a vectorfield which is Killing
outside of two disconnected components that avoid the horizons as well
as a neighborhood of the trapped set\footnote{Recall that for slowly
  rotating Kerr(-de Sitter), the trapped set is an $O(a)$ neighborhood
  of $\curlyBrace*{r=3M}$.}, which has an $O(a)$ deformation tensor and is timelike
up to the horizons, where it becomes null.

\begin{definition}
  \label{def:TAlmostKilling}
  We define the vectorfield
  \begin{equation}
    \label{eq:TALmostkilling:def}
    \TAlmostKilling 
    \vcentcolon= \KillT + \chi_{\delta_{\trap}}\KillPhi,
  \end{equation}
  where
  \begin{equation}
    \label{eq:TAlmostKilling:chi-delta:def}
    \chi_{\delta_{\trap}}\vcentcolon=\frac{a}{r^2+a^2}\chi_0\left(\delta_{\trap}^{-1}\frac{r-3M}{r}\right), 
  \end{equation}
  where $0<\delta_{\trap}\ll1$, and $\chi_0$ is a smooth bump
  function satisfying
  \begin{equation*}
    \chi_0(x) =0,\text{ for }\abs*{x}\le 1,\qquad 
      \chi_0(x) =1 \text{ for }  \abs*{x}\ge 2.
      \end{equation*}
\end{definition}

\begin{remark}
  Observe that $\TAlmostKilling = \KillT$ at the trapped
  set\footnote{Recall the definition of the trapped set from
    \zcref[cap]{lemma:trapping:KdS}.}, and $\TAlmostKilling=\HawkingVF$
  away from it, and is thus timelike in the entirety of the trapped
  region $\mathcal{M}_{\trap}$.
\end{remark}

The vectorfield $\TAlmostKilling$ enjoys the following simple property.
\begin{lemma}
  \label{lemma:TAlmostKilling:K-current:base}
  For $\TAlmostKilling$, we have that
  \begin{equation}
    \label{eq:TAlmostKilling:Deformation-tensor}
    \LieDerivative_{\TAlmostKilling}\left(\abs*{q}^2\Metric^{\alpha\beta}\right)
    = -2\delta_{\trap}\left(\partial_r\chi_{\delta_{\trap}}\right)\partial_{\phi}^{\alpha}\partial_r^{\beta}, \qquad
    \DeformationTensor{\TAlmostKilling}^{\alpha\beta} = \frac{\delta_{\trap}\partial_r\chi_{\delta_{\trap}}}{\abs*{q}^2}\partial_{\phi}^{\alpha}\partial_r^{\beta},
  \end{equation}
  where
  \begin{equation*}
    \partial_r\chi_{\delta_{\trap}} = -\frac{2ar}{\left(r^2+a^2\right)^2}\chi_0
    + O\left(a\delta_{\trap}^{-1}\chi_0'\right).
  \end{equation*}
  In addition, we have that 
  \begin{equation}
    \label{eq:TAlmostKilling:K-current:base}
    \begin{split}
      \abs*{q}^2\EMTensor\cdot\DeformationTensor{\TAlmostKilling_{\delta_{\trap}}}
      ={}& \left(r^2+a^2\right)\partial_r\chi_{\delta_{\trap}}\nabla_{\phi}\psi\cdot\nabla_{\HprVF}\psi\\
      ={}& \left(-\frac{2ar}{r^2+a^2}\chi_0 + O\left(a\delta_{\trap}^{-1}\right)\chi_0'\right)\nabla_{\phi}\psi\cdot\nabla_{\HprVF}\psi.
    \end{split}    
  \end{equation}
\end{lemma}
\begin{proof}
  Observe from \zcref[noname]{eq:CartarOp:def} that $\LieDerivative_{\TAlmostKilling}\abs*{q}^2 \mathcal{O}^{\alpha\beta}
    = \LieDerivative_{\TAlmostKilling}\frac{(r^2+a^2)^2}{\Delta}\HawkingVF^{\alpha}\HawkingVF^{\beta}
    = 0$. As a result, using \zcref[cap]{lemma:KdS-metric:inverse:THat-RHat}, we observe that
  \begin{align*}
    \LieDerivative_{\TAlmostKilling}\left(\abs*{q}^2\Metric^{\alpha\beta}\right)
    ={}& \LieDerivative_{\TAlmostKilling}\left(\Delta\partial_r^{\alpha}\partial_r^{\beta}
         + \abs*{q}^2 \mathcal{O}^{\alpha\beta}
         - \frac{(r^2+a^2)^2}{\Delta}\HawkingVF^{\alpha}\HawkingVF^{\beta}
         \right)\\
    ={}& \Delta[\TAlmostKilling, \partial_r]^{\alpha}\partial_r^{\beta}
         + \Delta\partial_r^{\alpha}[\TAlmostKilling, \partial_r]^{\beta}\\
    ={}& -2\Delta(\partial_r\chi_{\delta_{\trap}})\partial_{\phi}^{\alpha}\partial_r^{\beta}.
  \end{align*}
  The result in \zcref[noname]{eq:TAlmostKilling:K-current:base} then follows
  directly.
\end{proof}

\subsection{Properties of $\TAlmostKilling$}

In this section, we gather properties of the almost Killing vectorfield
$\TAlmostKilling$ defined in \zcref[cap]{def:TAlmostKilling}, starting with the main application of the divergence theorem.

\begin{lemma}
  \label{lemma:Killing-energy-estimate:basic-div-thm}
  Let $\TAlmostKilling$ be the almost Killing vectorfield defined in
\zcref[cap]{def:TAlmostKilling}. We have 
  \begin{equation}
    \label{eq:Killing-energy-estimate:basic-div-thm} 
    \CovariantDeriv\cdot\JCurrent{\TAlmostKilling, 0, 0}[\psi]
    = \EMTensor\cdot \DeformationTensor{\TAlmostKilling}
    + A_{\nu}\LeftDual{\psi}\cdot \HorCovDeriv^{\nu}\psi
    + \left(\TAlmostKilling-\KillT\right)^{\mu}\HorCovDeriv^{\nu}\psi^a\HorRiem_{ab\nu\mu}\psi^b
    + \nabla_{\TAlmostKilling}\psi\cdot\left(\WaveOpHork{2}\psi-V\psi\right),
  \end{equation}
  where
  \begin{equation}
    \label{eq:Killing-energy-estimate:A-def}
    A_{\mu} \vcentcolon= \volFormHor^{ab}\HorRiem_{ab\mu\nu}\KillT^{\nu}.
  \end{equation}
\end{lemma}
\begin{proof}
  Observe that from \zcref[noname]{eq:div-thm:general}, we have that
  \begin{align*}
      \CovariantDeriv\cdot\JCurrent{\TAlmostKilling, 0, 0}[\psi] &
      = \KCurrent{\TAlmostKilling, 0, 0}[\psi]
           + \nabla_{\TAlmostKilling}\psi\cdot\left(\WaveOpHork{2}\psi-V\psi\right),&
           \KCurrent{\TAlmostKilling, 0, 0}[\psi] &
    = \EMTensor\cdot \DeformationTensor{\TAlmostKilling}
    + \TAlmostKilling^{\mu}\HorCovDeriv^{\nu}\psi^a\cdot\HorRiem_{ab\nu\mu}\psi^b.   
  \end{align*}

   We can rewrite this as
  \begin{equation*}
    \KCurrent{\TAlmostKilling, 0, 0}[\psi]
    = \EMTensor\cdot \DeformationTensor{\TAlmostKilling}
    + \KillT^{\mu}\HorCovDeriv^{\nu}\psi^a\HorRiem_{ab\nu\mu}\psi^b
    + \left(\TAlmostKilling-\KillT\right)^{\mu}\HorCovDeriv^{\nu}\psi^a\cdot\HorRiem_{ab\nu\mu}\psi^b.
  \end{equation*}
  Since $\HorRiem_{ab\nu\mu}$ is antisymmetric with respect to the
  $(a,b)$ indices, we can rewrite
  \begin{equation*}
    \KCurrent{\TAlmostKilling, 0, 0}[\psi]
    = \EMTensor\cdot \DeformationTensor{\TAlmostKilling}
    + \frac{1}{2}\KillT^{\mu} \volFormHor^{ab}\HorRiem_{ab\nu\mu}\LeftDual{\psi} \cdot \HorCovDeriv^{\nu}\psi
    + \left(\TAlmostKilling-\KillT\right)^{\mu}\HorCovDeriv^{\nu}\psi^a\cdot \HorRiem_{ab\nu\mu}\psi^b.
  \end{equation*}
  We now make use of the spacetime 1-form $A_{\mu}$ from
  \zcref[noname]{eq:Killing-energy-estimate:A-def} to see that
  \begin{equation*}
    \KCurrent{\TAlmostKilling, 0, 0}[\psi]
    = \EMTensor\cdot \DeformationTensor{\TAlmostKilling}
    + A_{\nu}\LeftDual{\psi}\cdot \HorCovDeriv^{\nu}\psi
    + \left(\TAlmostKilling-\KillT\right)^{\mu}\HorCovDeriv^{\nu}\psi^a\HorRiem_{ab\nu\mu}\psi^b.
  \end{equation*}
  This concludes the proof of \zcref[cap]{lemma:Killing-energy-estimate:basic-div-thm}. 
\end{proof}
Next, we introduce some auxiliary lemmas which will be used in the
proof of \zcref[cap]{prop:Killing-energy-estimate}. Critical to
the proof of \zcref[cap]{prop:Killing-energy-estimate} will be properties of the one-form $A_{\mu}$. We begin by computing the values of $A$.
\begin{lemma}
  \label{lemma:Killing-energy:A-values}
  Let $A$ be the spacetime $1$-form defined in
  \zcref[noname]{eq:Killing-energy-estimate:A-def}.  Then, for a given null 
  frame $(e_3,e_4,e_1,e_2)$ we have that
  \begin{equation}
    \label{eq:Killing-energy:A-values}
    \begin{split}
      A_4 ={}& -4\LeftDual{\rho}\KillT^3
               - 4\left(\etaBar\wedge \eta\right)\KillT^3
               + \Trace \chi \left(\horProj{\KillT}\wedge \etaBar\right)
               - \aTrace{\chi}\left(\etaBar\cdot \horProj{\KillT}\right).\\
      A_3 ={}& 4\LeftDual{\rho}\KillT^4
               + 4\left(\etaBar\wedge \eta\right)\KillT^4
               + \Trace \chiBar \left(\horProj{\KillT}\wedge \eta\right)
               - \aTrace{\chiBar}\left(\eta\cdot \horProj{\KillT}\right),\\
      A_e ={}& \left(-\Trace\chiBar \LeftDual{\eta}_e
               + \aTrace{\chiBar}\eta_e
               \right)\KillT^3
               + \left(
               -\Trace\chi\LeftDual{\etaBar}_e
               + \aTrace{\chi}\etaBar_e
               \right)\KillT^4\\
             & - \frac{1}{2}\left(
               4\rho - \frac{4\Lambda}{3}
               + \Trace\chi\Trace\chiBar
               + \aTrace{\chi}\aTrace{\chiBar}
               \right)\LeftDual{\horProj{\KillT}}_e.
    \end{split}
  \end{equation}
\end{lemma}
\begin{proof}
  We begin by rewriting
  \begin{equation*}
    A_{\mu} = \volFormHor^{ab}\HorRiem_{ab\mu3}\KillT^3
    + \volFormHor^{ab}\HorRiem_{ab\mu4}\KillT^4
    + \volFormHor^{ab}\HorRiem_{ab\mu c}\KillT^c.
  \end{equation*}
  We can then compute each of the components of $A_{\mu}$
  individually. In \KdS, we have that $\horProj{\KillT}^a = \KillT^a$
  where we recall that $a,b,c,d,e$ are used to denote the horizontal
  indices. Then we can compute that
  \begin{align*}
    A_4 ={}& \volFormHor^{ab}\HorRiem_{ab43}\KillT^3
             + \volFormHor^{ab}\HorRiem_{ab4c}\KillT^c\\
    ={}& \volFormHor^{ab}\left(
         -2\volFormHor_{ab}\LeftDual{\rho}
         - 2\left(
         \etaBar_a\eta_b - \eta_a\etaBar_b
         \right)
         \right)\KillT^3
            + \frac{1}{2}\volFormHor^{ab}\left(
             \Trace\chi\left(\delta_{ac}\eta_b - \delta_{bc}\eta_a\right)
             + \aTrace{\chi}\left(\volFormHor_{ca}\etaBar_b - \volFormHor_{cb}\etaBar_a\right)\KillT^c
             \right)\\
    ={}& -4\LeftDual{\rho}\KillT^3
         - 4\left(\etaBar\wedge \eta\right)\KillT^3
         + \frac{1}{2}\volFormHor^{ab}\left(
         \Trace\chi\left(\etaBar_b\KillT_a - \etaBar_a\KillT_b\right)
         + \aTrace{\chi}\left(
         -\etaBar_b\horProj{\KillT}_a + \horProj{\KillT}_b\etaBar_a
         \right)
         \right)\\
    ={}& -4\LeftDual{\rho}\KillT^3
         - 4\left(\etaBar\wedge \eta\right)\KillT^3
         + \Trace\chi\left(\horProj{\KillT}\wedge \etaBar\right)
         - \aTrace{\chi}\left(\etaBar\cdot \horProj{\KillT}\right),
  \end{align*}
  and
  \begin{align*}
    A_3={}& \volFormHor^{ab}\HorRiem_{ab34}\KillT^4 + \volFormHor^{ab}\HorRiem_{ab3c}\KillT^c = 4\LeftDual{\rho}\KillT^4
         +  4\left(\etaBar\wedge \eta\right)\KillT^4
         -\aTrace{\chiBar}\left(\eta\cdot\horProj{\KillT}\right),
  \end{align*}
  and
  \begin{align*}
    A_e ={}& \volFormHor^{ab}\HorRiem_{abe3}\KillT^3
             + \volFormHor^{ab}\HorRiem_{abe4}\KillT^4
             + \volFormHor^{ab}\HorRiem_{abed}\KillT^{d}\\
    ={}& \frac{1}{2}\volFormHor^{ab}\left(
         -\Trace\chiBar\left(\delta_{ea}\eta_b-\delta_{eb}\eta_a\right)
         -\aTrace{\chiBar}\left(\volFormHor_{ea}\eta_b-\volFormHor_{eb}\eta_a\right)
         \right)\KillT^3\\
           &+ \frac{1}{2}\volFormHor^{ab}\left(
             -\Trace\chi\left(\delta_{ea}\etaBar_b-\delta_{eb}\etaBar_a\right)
             - \aTrace{\chi}\left(\volFormHor_{ea}\etaBar_b-\volFormHor_{eb}\etaBar_a\right)             
             \right)\KillT^4\\
           &+ \frac{1}{2}\volFormHor^{ab}\left(
             -2\volFormHor_{ab}\volFormHor_{ed}\rho
             + \frac{2\Lambda}{3}\left(\delta_{ae}\delta_{bd}- \delta_{be}\delta_{ad}\right)
             - \frac{1}{2}\left(
             \Trace\chi\Trace\chiBar
             + \aTrace{\chi}\aTrace{\chiBar}
             \right)\volFormHor_{ab}\volFormHor_{ed}
             \right)\KillT^d\\
    ={}& \left(\chiBar\LeftDual{\eta}+ \aTrace{\chiBar}\eta\right)\KillT^3
         + \left(-\Trace\chi\LeftDual{\etaBar}+ \aTrace{\chi}\eta\right)\KillT^4
        - \frac{1}{2}\left(
         4\rho - \frac{4\Lambda}{3}
         + \Trace\chi\Trace\chiBar
         + \aTrace{\chi}\aTrace{\chiBar}
         \right)\LeftDual{(\horProj{\KillT})}.
  \end{align*}
  This concludes the proof of
  \zcref[cap]{lemma:Killing-energy-estimate:basic-div-thm}.
\end{proof}

A critical property of $A_{\mu}$ is that it satisfies a divergence
property.

\begin{corollary}
  \label{coro:Killing-energy:A-div-form}
  We have that for $A$ defined as in \zcref[noname]{eq:Killing-energy-estimate:A-def},
  \begin{equation}
    \label{eq:Killing-energy:A-div-form}
    A_{\mu} = -\CovariantDeriv_{\mu}\left(\Im \left(\frac{2M}{(1+\gamma)q^2} + \frac{2q\Lambda}{3(1+\gamma)}\right)\right).
  \end{equation}
\end{corollary}
\begin{proof}
  Since $A$ is $0$-conformally invariant, it suffices to verify
  \zcref[noname]{eq:Killing-energy:A-div-form} in the principal
  ingoing frame $(e_3,e_4)=(\ein_3,\ein_4)$.  In \KdS, we have that
  \begin{equation*}
    \ein_4\cdot\KillT = \frac{1}{2(1+\gamma)}, \qquad
    \ein_3\cdot\KillT = \frac{\Delta}{2(1+\gamma)\abs*{q}^2},\qquad
    \KillT^b = - \frac{a}{1+\gamma}\Re(\CCOneFormJ)_b.
  \end{equation*}
  Substituting these values into the identities of
  \zcref[cap]{lemma:Killing-energy:A-values}, we have that
  \begin{align*}
    A_4
    ={}& -2 \LeftDual{\rho}\frac{\Delta}{(1+\gamma)\abs*{q}^2}
         - 2\left(\etaBar\wedge \eta\right)\frac{\Delta}{(1+\gamma)\abs*{q}^2}
         - \frac{a\sqrt{\kappa}}{(1+\gamma)}\Trace\chi\left(
         \Re(\CCOneFormJ)\wedge \etaBar
         \right)
         + \frac{a\sqrt{\kappa}}{(1+\gamma)}\aTrace{\chi}\left(\eta\cdot \Re(\CCOneFormJ)\right),\\
    A_3
    ={}& \frac{2}{1+\gamma}\LeftDual{\rho}
         + \frac{2}{1+\gamma}\left(\etaBar\wedge\eta\right)
         - \frac{a\sqrt{\kappa}}{(1+\gamma)}\Trace\chiBar\left(\Re(\CCOneFormJ) \wedge \eta\right)
         + \frac{a\sqrt{\kappa}}{(1+\gamma)}\aTrace{\chiBar}\left(\eta\cdot \Re(\CCOneFormJ)\right), \\
    A_b={}& \frac{\Delta}{2(1+\gamma)\abs*{q}^2}\left(
            -\Trace\chiBar\LeftDual{\eta}_b
            + \aTrace{\chiBar}\eta_b
            \right)
            + \frac{1}{2(1+\gamma)}\left(
            - \Trace\chi \LeftDual{\etaBar}_b
            + \aTrace{\chi}\etaBar_b
            \right)\\
       &+ \frac{a\sqrt{\kappa}}{2(1+\gamma)}\left(
         4\rho - \frac{4\Lambda}{3}
         + \Trace\chi\Trace\chiBar
         + \aTrace{\chi}\aTrace{\chiBar}
         \right)\LeftDual{\Re(\CCOneFormJ)}_b.
  \end{align*}
  Now observe that
  \begin{align*}
    2(1+\gamma)A_b
    ={}& \left(
         -\Re(\Trace\XBar)\Im(H_b)
         -\Im(\Trace\XBar)\Re(H_b)
         \right) \frac{\Delta}{\abs*{q}^2}
         + \left(
         -\Re(\Trace X)\Im\HBar_b
         -\Im(\Trace X)\Re\HBar_b
         \right)\\
       &+ a\sqrt{\kappa}\left(
         4\Re P
         + \Re(\Trace X)\Re(\Trace \XBar)
         + \Im(\Trace X)\Im(\Trace \XBar)
         \right)\LeftDual{\Re(\CCOneFormJ)}_b\\
    ={}& -\Im(\Trace\XBar H_b)\frac{\Delta}{\abs*{q}^2}
         - \Im (\Trace X \HBar_b)
         + a\Re\left(4P - \frac{4\Lambda}{3} + \Trace X\overline{\Trace\XBar}\right)\Re\LeftDual{\CCOneFormJ}_b.
  \end{align*}
  We can now substitute in the exact values of the Ricci coefficients
  and the null curvature components in the principal ingoing frame
  from \zcref[cap]{lemma:Kerr:ingoing-PG:Ric-and-curvature}, so that
  \begin{align*}
    2(1+\gamma)A_b
    ={}& \Im\left(\frac{2}{\overline{q}}\frac{a q}{\abs{q}^2}\CCOneFormJ_b\right)\frac{\Delta}{\abs*{q}^2}
         + \Im\left(\frac{2\Delta\overline{q}}{\abs*{q}^4}\frac{a\overline{q}}{\abs*{q}^2}\CCOneFormJ_b\right)
         - a\Re\left(\frac{8M}{q^3} + \frac{4\Delta}{\abs*{q}^2q^2}\right)\LeftDual{\Re \CCOneFormJ}_b\\
    ={}& \frac{2a\Delta}{\abs*{q}^2}\left(
         \Im\left(
         \frac{1}{\overline{q}^2}
         + \frac{1}{q^2}
         \right)\Im\CCOneFormJ_b
         - \Re\left(\frac{2}{q^2}\right)\LeftDual{\Re\CCOneFormJ}_b
         \right)
         - a\Re\left(\frac{8M}{q^3}\right)\LeftDual{\Re\CCOneFormJ_b}
         - a \frac{4\Lambda}{3}\LeftDual{\Re\CCOneFormJ_b}\\
    ={}& - a\Re\left(\frac{8M}{q^3}\right)\LeftDual{\Re\CCOneFormJ_b}
         - a \frac{4\Lambda}{3}\LeftDual{\Re\CCOneFormJ_b}.
  \end{align*}
  We infer from $\nabla\cos\theta = \Re \ImagUnit\CCOneFormJ = -\LeftDual{\Re \CCOneFormJ}$ that
  \begin{align*}
    A_b ={}& a \Re\left(\frac{4M}{q^3}\right)\nabla_b\cos\theta
    ={}\Im \frac{4M}{q^3}\nabla_bq
    ={} -\nabla_b\left(\Im \frac{2M}{q}\right).
  \end{align*}
  Next, we consider $A_4$.
  \begin{align*}
    A_4={}& - 2 \LeftDual{\rho}\frac{\Delta}{\abs*{q}^2}
            - 2\left(\etaBar \cdot \LeftDual{\eta}\right)\frac{\Delta}{\abs*{q}^2}
            - a\Trace \chi \left(\Re{\CCOneFormJ}\cdot \LeftDual{\etaBar}\right)
            + a \aTrace{\chi}\left(\etaBar\cdot \Re\CCOneFormJ\right)\\
    ={}& -2\LeftDual{\rho}
         - 2\left(\Re\HBar\cdot\Im H\right)\frac{\Delta}{\abs*{q}^2}
         - a\Im\left(\Trace X\HBar\right)\cdot \Re \CCOneFormJ.
  \end{align*}
  Recall that $\Im \CCOneFormJ = \LeftDual{\Re \CCOneFormJ}$. As a result, we have that
  \begin{align*}
    A_4={}& -2\LeftDual{\rho}\frac{\Delta}{\abs*{q}^2}
            + 2a^2\frac{\Delta}{\abs*{q}^6}\Re\left(\overline{q}\CCOneFormJ
            \right)\cdot \Im\left(q\CCOneFormJ\right)
            - a\Im \left(
            \frac{2}{q}\frac{\Delta}{\abs*{q}^2}\left(-\frac{a\overline{q}}{\abs*{q}^2}\CCOneFormJ\right)
            \right)\cdot\Re\CCOneFormJ 
    = -2\LeftDual{\rho}\frac{\Delta}{\abs*{q}^2}.
  \end{align*}
  Since in the principal ingoing frame,
  $e_4(q) = \frac{\Delta}{\abs*{q}^2}$, we can now plug in the value
  of $\LeftDual{\rho}$ on Kerr-de Sitter from 
  \zcref[cap]{lemma:Kerr:ingoing-PG:Ric-and-curvature} into the expression
  above to see that
  \begin{align*}
    A_4={}& 4\Im\frac{M}{q^3}\frac{\Delta}{\abs*{q}^2}
    ={} 4\Im\left(\frac{M}{q^3}\right)e_4(q)
    ={} -e_4\left(\Im \frac{2M}{q^2}\right).
  \end{align*}
  Following a similar computation, we obtain $A_3 = -e_3\left(\Im \frac{2M}{q^2}\right)$. Since $e_3(\theta) = e_4(\theta) = 0$, we finally get
  \begin{equation*}
    A_4 = e_4\left(\frac{2M}{(1+\gamma)q^2} + \frac{2q\Lambda}{3(1+\gamma)}\right),\qquad
    A_3 = e_3\left(\frac{2M}{(1+\gamma)q^2} + \frac{2q\Lambda}{3(1+\gamma)}\right),
  \end{equation*}
  as desired.
\end{proof}

\begin{lemma}
  \label{lemma:energy:error-terms}
  Let $\psi\in \realHorkTensor{2}$ be a solution to \zcref[noname]{eq:model-problem-gRW}. 
  Consider the modified current
  \begin{equation*}
    \JCurrentMod{\TAlmostKilling, 0, 0, \mathring{w}}_{\mu}[\psi]
    = \JCurrent{\TAlmostKilling, 0, 0}_{\mu}[\psi]
    + \mathring{w}\LeftDual{\psi}\cdot \HorCovDeriv_{\mu}\psi, \qquad
    \mathring{w} = \Im \left(\frac{M}{(1+\gamma)q^2} + \frac{2a\Lambda}{3(1+\gamma)}\right),
  \end{equation*}
  where the vectorfield $\TAlmostKilling$ is given by Definition
  \ref{def:TAlmostKilling}. Then, 
  \begin{equation*}
    \begin{split}
      &\abs*{\CovariantDeriv\cdot\left(\JCurrentMod{\TAlmostKilling, 0, 0, \mathring{w}}[\psi] + \KillT\frac{4a\cos\theta}{\abs*{q}^2}\abs*{\psi}^2\right)
        - \left(\nabla_{\TAlmostKilling}\psi + \mathring{w}\LeftDual{\psi}\right)\cdot N 
        }\\
      \lesssim{}& \bOne_{\mathcal{M}_{\cancel{\text{trap}}}}\left(
                  \delta_{\trap}^{-1}\frac{a}{r^3}\abs*{\nabla_{\HprVF}\psi}\abs*{\nabla_{\KillPhi}\psi}
                  + \frac{a}{r^4}\abs*{\nabla_{\KillT}\psi}\abs*{\nabla_{\KillPhi}\psi}
                  + \frac{a}{r^3}\left(
                  \abs*{\nabla_{\HprVF}\psi}
                  + \abs*{\nabla_{\KillT}\psi}
                  + \abs*{\nabla \psi}
                  \right)\abs*{\psi}
                  \right).
    \end{split}
  \end{equation*}
\end{lemma}
\begin{proof}  
  Combining \zcref[noname]{eq:Killing-energy-estimate:basic-div-thm} and
  \zcref[noname]{eq:Killing-energy:A-div-form}, we have that
  \begin{equation}
    \label{eq:energy:div-thm-with-div-structure}
    \begin{split}
      \CovariantDeriv\cdot\JCurrent{\TAlmostKilling, 0, 0}[\psi]
      ={}& \EMTensor\cdot \DeformationTensor{\TAlmostKilling}
           - \CovariantDeriv_{\nu}\left(\Im\left(\frac{M}{\left(1+\gamma\right)q^2} + \frac{2q\Lambda}{3(1+\gamma)}\right)\right)\LeftDual{\psi}\cdot \HorCovDeriv^{\nu}\psi\\
         & + \left(\TAlmostKilling-\KillT\right)^{\mu}\HorCovDeriv^{\nu}\psi^a\HorRiem_{ab\nu\mu}\psi^b
           + \nabla_{\TAlmostKilling}\psi\cdot \left(\WaveOpHork{2}\psi - V\psi\right).
    \end{split}
  \end{equation}
  Due to the divergence structure present in the second term on the
  \RHS{} of \zcref[noname]{eq:energy:div-thm-with-div-structure}, we can
  modify the divergence identity in
  \zcref[noname]{eq:energy:div-thm-with-div-structure} by a lower-order term
  to cancel the second term. To this end, we consider the modified
  current
  \begin{equation}
    \label{eq:energy:JCurrentMod:def}
    \JCurrentMod{\TAlmostKilling, 0, 0, \mathring{w}}_{\mu}[\psi]
    = \JCurrent{\TAlmostKilling, 0, 0}_{\mu}[\psi]
    + \mathring{w}\LeftDual{\psi}\cdot \HorCovDeriv_{\mu}\psi,
  \end{equation}
  where $\mathring{w}\vcentcolon= \mathring{w}(r, \cos\theta)$ is a function
  that we will specify later. Then observe that
  \begin{align*}
    \CovariantDeriv^{\mu}\left(\mathring{w}\LeftDual{\psi}\cdot \HorCovDeriv_{\mu}\psi\right)
    ={}& \mathring{w}\LeftDual{\psi}\cdot \HorCovDeriv^{\mu}\HorCovDeriv_{\mu}\psi
         + \mathring{w}\LeftDual{\HorCovDeriv}^{\mu}\psi\cdot \HorCovDeriv_{\mu}\psi
         + \CovariantDeriv^{\mu}(\mathring{w})\LeftDual{\psi}\cdot \HorCovDeriv_{\mu}\psi\\
    ={}& \mathring{w}\LeftDual{\psi}\cdot \WaveOpHork{2}\psi
         + \CovariantDeriv^\mu(\mathring{w})\LeftDual{\psi}\cdot\HorCovDeriv_{\mu}\psi.
  \end{align*}
  Using the form of the model equation in
  \zcref[noname]{eq:model-problem-gRW}, we have that
  \begin{align*}
    \CovariantDeriv^{\mu}\left(\mathring{w}\LeftDual{\psi}\cdot \HorCovDeriv_{\mu}\psi\right)
    ={}& \mathring{w}\LeftDual{\psi}\cdot \left(
         V\psi
         - \frac{4a(1+\gamma)\cos\theta}{\abs*{q}^2}\LeftDual{\nabla_{\KillT}\psi}
         + N
         \right)
         + \CovariantDeriv^{\mu}(\mathring{w})\LeftDual{\psi}\dot \HorCovDeriv_{\mu}\psi\\
    ={}& \CovariantDeriv^{\mu}(\mathring{w})\LeftDual{\psi}\cdot \HorCovDeriv_{\mu}\psi
         - \CovariantDeriv_{\mu}\left(\KillT^{\mu}\left(\mathring{w}\frac{4a(1+\gamma)\cos\theta}{\abs*{q}^2}\abs*{\psi}^2\right)\right)
         + \mathring{w}\LeftDual{\psi}\cdot N. 
  \end{align*}
  We now specify our choice of $\mathring{w}$. We choose $\mathring{w}$
  \begin{equation}
    \label{eq:energy:w-choice}
    \mathring{w} \vcentcolon= \Im \left(\frac{M}{(1+\gamma)q^2} + \frac{2a\Lambda}{3(1+\gamma)}\right).
  \end{equation}
  Plugging this into the definition of
  $\JCurrentMod{\TAlmostKilling, 0, 0, \mathring{w}}[\psi]$ in
  \zcref[noname]{eq:energy:JCurrentMod:def}, we have that
  \begin{equation}
    \label{eq:energy:error-terms:div-thm}
    \begin{split}    
      \CovariantDeriv\cdot\left(\JCurrentMod{\TAlmostKilling, 0, 0, \mathring{w}}[\psi] + \KillT\left(\mathring{w}\frac{4a(1+\gamma)\cos\theta}{\abs*{q}^2}\abs*{\psi}^2\right)\right)
      ={}& \EMTensor\cdot \DeformationTensor{\TAlmostKilling}
           + \left(\TAlmostKilling-\KillT\right)^{\mu}\HorCovDeriv^{\nu}\psi^a\HorRiem_{ab\nu\mu}\psi^b\\
         & + \nabla_{\TAlmostKilling}\psi\cdot \left(\WaveOpHork{2}\psi - V\psi\right)
           + \mathring{w}\LeftDual{\psi}\cdot N.
    \end{split}
  \end{equation}
  We now control each of the terms on the \RHS{} of
  \zcref[noname]{eq:energy:error-terms:div-thm} individually.  Using
  \zcref[noname]{eq:TAlmostKilling:K-current:base}, we can control
  \begin{align*}
    \abs*{\EMTensor\cdot\DeformationTensor{\TAlmostKilling}}
    \lesssim{}& \left(\frac{4 a r}{\abs*{q}^2\left(r^2+a^2\right)}\abs*{\chi_0} + O\left(a\delta_{\trap}^{-1}\right)\abs*{\chi_0'}
                \right)\abs*{\nabla_{\phi}\psi}\abs*{\nabla_{\HprVF}\psi}\\
    \lesssim{}& \bOne_{\mathcal{M}_{\slashed{\text{trap}}}}\delta_{\trap}^{-1}\frac{a}{r^3}\abs*{\nabla_{\phi}\psi}\abs*{\nabla_{\HprVF}\psi}.
  \end{align*}
  Recalling the definition of $\TAlmostKilling$ in 
  \zcref{def:TAlmostKilling}, we have that
  \begin{align}
    \left(\TAlmostKilling- \KillT\right)^{\mu}\HorCovDeriv^{\nu}\psi^a\HorRiem_{ab\nu\mu}\psi^b\notag
    ={}&\frac{1}{2}\chi_{\delta_{\trap}}(r)\KillPhi^{\mu}\volFormHor^{ab}\HorRiem_{ab\nu\mu}\LeftDual{\psi}\cdot \HorCovDeriv^{\nu}\psi\notag\\
    ={}& -\frac{1}{4}\chi_{\delta_{\trap}}(r)\KillPhi^{\mu}\volFormHor^{ab}\left(
          \HorRiem_{ab3\mu}\LeftDual{\psi}\cdot\nabla_4\psi
         + \HorRiem_{ab4\mu}\LeftDual{\psi}\cdot\nabla_3\psi
         + \HorRiem_{abc\mu}\LeftDual{\psi}\cdot\nabla^c\psi
         \right).\label{eq:energy:error-terms:conclusion:aux1}
  \end{align}
  Recall that in the principal ingoing null frame, we have $\KillPhi^3 = O\left(a r^{-2}\Delta\right)$, $\KillPhi^4= O(a)$ and $\KillPhi^c = O(r^2)\Re\CCOneFormJ^c$. Recalling the definition of $\HorRiem$ and $\BHorizontalCurv$
  \zcref[noname]{eq:horizontal-structure:B-horizontal-curvature:def} and the
  values of $\BHorizontalCurv$ in \zcref[cap]{prop:B-hor:components},
  we have that
  \begin{align*}
    \KillPhi^{\mu}\volFormHor^{ab}\HorRiem_{ab3\mu}
    ={}& \KillPhi^4\volFormHor^{ab}\HorRiem_{ab34}
         + \KillPhi^c\volFormHor^{ab}\HorRiem_{ab3c}\\
    ={}& O(a)\left(\LeftDual{\rho} + \eta\etaBar\right)
         + O\left(r\right)\chiBar\eta,
    \\
    ={}& O\left(ar^{-2}\right),\\
    \KillPhi^{\mu}\volFormHor^{ab}\HorRiem_{ab4\mu}
    ={}& \KillPhi^3\volFormHor^{ab}\HorRiem_{ab43}
         + \KillPhi^c\volFormHor^{ab}\HorRiem_{ab4c}\\
    ={}& O\left(ar^{-2}\Delta\right)\left(\LeftDual{\rho} + \eta\etaBar\right)
         + O(r)\chi\etaBar\\
    ={}& O\left(ar^{-4}\Delta\right),\\
    \KillPhi^{\mu}\volFormHor^{ab}\HorRiem_{abc\mu}
    ={}& \KillPhi^4\volFormHor^{ab}\HorRiem_{abc4}
         + \KillPhi^3\volFormHor^{ab}\HorRiem_{abc3}
         + \KillPhi^d\volFormHor^{ab}\HorRiem_{abcd}\\
    ={}&O(a)\chi\etaBar
         + O\left(ar^{-2}\Delta\right)\chiBar\eta
         + O(r)\left(\rho + \Lambda + \chi\chiBar\right)
    \\
    ={}& O\left(r^{-1} + r\Lambda\right)\\
    ={}& O\left(r^{-1}\right),
  \end{align*}
  where the last line follows from the fact that on $\Manifold$,
  $\Lambda\lesssim r^{-2}$.
  Recalling the definition of $\chi_{\delta_{\trap}}$ from
  \zcref[noname]{eq:TAlmostKilling:chi-delta:def}, we then have that
  \begin{align}
    & \abs*{\left(\TAlmostKilling - \KillT\right)^{\mu}\HorCovDeriv^\nu\psi^a\HorRiem_{ab\nu\mu}\psi^b}\notag \\
    \lesssim{}& \frac{a}{r^2}\bOne_{\mathcal{M}_{\slashed{\text{trap}}}}\left(
                \abs*{\KillPhi^{\mu}\volFormHor^{ab}\HorRiem_{ab3\mu}}\abs*{\LeftDual{\psi}\cdot \nabla_4\psi}
                + \abs*{\KillPhi^\mu\volFormHor^{ab}\HorRiem_{ab4\mu}}\abs*{\LeftDual{\psi}\cdot \nabla_3\psi}
                + \abs*{\KillPhi^{\mu}\volFormHor^{ab}\HorRiem_{abc\mu}}\abs*{\LeftDual{\psi}\cdot \nabla^c\psi}
                \right)\notag \\
    \lesssim{}& \frac{a}{r^3}\bOne_{\mathcal{M}_{\slashed{\text{trap}}}}\left(
                \abs*{\nabla_4\psi}
                + r^{-2}\abs*{\Delta \nabla_3\psi}
                + \abs*{\nabla \psi}\right)\abs*{\psi}\notag \\
    \lesssim{}& \frac{a}{r^3}\bOne_{\mathcal{M}_{\slashed{\text{trap}}}}\left(
                \abs*{\nabla_{\HprVFGen}\psi}
                + \abs*{\nabla_T\psi}
                + \abs*{\nabla \psi}
                \right)\abs*{\psi}.\label{eq:energy:error-terms:conclusion:aux2}
  \end{align}
  Then, using the equation \zcref[noname]{eq:model-problem-gRW},
  \begin{align}
    \nabla_{\TAlmostKilling}\psi\cdot \left(\WaveOpHork{2}\psi - V\psi\right)
    ={}& -\frac{4a(1+\gamma)\cos\theta}{\abs*{q}^2}\nabla_{\TAlmostKilling}\psi \cdot \LeftDual{\nabla}_{\KillT}\psi
         + \nabla_{\TAlmostKilling}\psi\cdot N\notag\\
    ={}& -\frac{4a(1+\gamma)\cos\theta}{\abs*{q}^2}\chi_{\delta_{\trap}}\nabla_{\KillPhi}\psi\cdot \LeftDual{\nabla}_{\KillT}\psi
         + \nabla_{\TAlmostKilling}\psi\cdot N, \label{eq:energy:error-terms:conclusion:aux3}
  \end{align}
  where we have used the crucial cancellation $\nabla_{\KillT}\psi\cdot \LeftDual{\nabla}_{\KillT}\psi = 0$. Combining \zcref[noname]{eq:energy:error-terms:div-thm},
  \zcref[noname]{eq:energy:error-terms:conclusion:aux1},
  \zcref[noname]{eq:energy:error-terms:conclusion:aux2}, and
  \zcref[noname]{eq:energy:error-terms:conclusion:aux3} concludes the proof of
  \zcref[cap]{lemma:energy:error-terms}.
\end{proof}
We also have the following lemma regarding the boundary terms that
appear in proving the energy estimate in 
\zcref[cap]{prop:Killing-energy-estimate}.

\begin{lemma}
  \label{lemma:energy:boundary}
  The following hold true for $c_0>0$ sufficiently small and for
  $\frac{a}{M}\ll 1$,
  \begin{align}    
      \int_{\Sigma(\tau)}\EMTensor\left(\HawkingVF, N_{\Sigma}\right)
      \ge{}&c_0 \EnergyFlux_{\operatorname{deg}}[\psi](\tau)
             - O(\delta_{\Horizon})E_{r\le r_{\EventHorizonFuture}}[\psi](\tau)
             - O(\delta_{\Horizon})\EnergyFlux_{r\ge r_{\CosmologicalHorizonFuture}}[\psi](\tau),
             \label{eq:energy:boundary:Sigma}
    \\
      \int_{\SigmaStar(\tau_1,\tau_2)}\EMTensor\left(\HawkingVF, N_{\SigmaStar}\right)
    \gtrsim{}&  \int_{\SigmaStar(\tau_1,\tau_2)}\abs*{\nabla_{3}\psi}^2
               -\delta_{\Horizon}\SpacelikeFlux_{\SigmaStar}[\psi](\tau_1,\tau_2),
               \label{eq:energy:boundary:SigmaStar}
    \\
      \int_{\mathcal{A}(\tau_1,\tau_2)}\EMTensor\left(\HawkingVF, N_{\mathcal{A}}\right)
    \gtrsim{}& \int_{\mathcal{A}(\tau_1,\tau_2)}\abs*{\nabla_{4}\psi}^2 
               - \delta_{\Horizon}\SpacelikeFlux_{\mathcal{A}}[\psi](\tau_1,\tau_2).
               \label{eq:energy:boundary:ACal}
  \end{align}
\end{lemma}

\begin{proof}
  Given our convention for $N_{\Sigma}$ from
  \zcref[noname]{eq:hypersurface-normal-def}, we have $\EMTensor\left(e_\mu,N_{\Sigma}\right)
    = \frac{1}{2}e_3(\tau)\EMTensor_{\mu 4}
         + \frac{1}{2}e_4(\tau)\EMTensor_{\mu 3}
         - \nabla^a(\tau)\EMTensor_{\mu a}$ for $\mu\in\{3,4\}$.
  Recalling from \zcref[noname]{eq:EM-tensor:null-components} the components
  of the energy-momentum 
  tensor, we infer that
  \begin{align*}
    \EMTensor\left(e_4,N_{\Sigma}\right)
    \ge{}& \frac{1}{2}e_3(\tau)\abs*{\nabla_4\psi}^2
    + \frac{1}{2}e_4(\tau)\left(\abs*{\nabla\psi}^2 + V\abs*{\psi}^2\right)
    - \abs*{\nabla\tau}\abs*{\nabla_4\psi}\abs*{\nabla\psi},\\
    \EMTensor\left(e_3,N_{\Sigma}\right)
    \ge{}& \frac{1}{2}e_3(\tau)\left(\abs*{\nabla\psi}^2 + V\abs*{\psi}^2\right)
           + \frac{1}{2}e_4(\tau)\abs*{\nabla_3\psi}^2
           - \abs*{\nabla\tau}\abs*{\nabla_3\psi}\abs*{\nabla\psi}.
  \end{align*}
  From the properties of $\tau$ in \zcref[cap]{prop:properties-of-global-coordinates} we obtain
  \begin{align*}
    \abs*{\nabla\tau}\abs*{\nabla_4\psi}\abs*{\nabla\psi}
    \le{}& \sqrt{\frac{8}{9}e_4(\tau)e_3(\tau)}\abs*{\nabla_4\psi}\abs*{\nabla\psi}\\
    \le{}& \sqrt{\frac{8}{9}}\left(
           \frac{1}{2}e_3(\tau)\abs*{\nabla_4\psi}^2
           + \frac{1}{2}e_4(\tau)\abs*{\nabla\psi}^2
           \right),\\
    \abs*{\nabla\tau}\abs*{\nabla_3\psi}\abs*{\nabla\psi}
    \le{}& \sqrt{\frac{8}{9}e_4(\tau)e_3(\tau)}\abs*{\nabla_3\psi}\abs*{\nabla\psi}\\
    \le{}& \sqrt{\frac{8}{9}}\left(
           \frac{1}{2}e_3(\tau)\abs*{\nabla\psi}^2
           + \frac{1}{2}e_4(\tau)\abs*{\nabla_3\psi}^2
           \right).
  \end{align*}
  Thus, we have that
  \begin{align*}
    \EMTensor(e_4,N_{\Sigma})
    \ge{}& \left(1-\sqrt{\frac{8}{9}}\right)
           \left(
           \frac{1}{2}e_3(\tau)\abs*{\nabla_4\psi}^2
           + \frac{1}{2}e_4(\tau)\abs*{\nabla\psi}^2
           \right)
           + \frac{1}{2}e_4(\tau)V\abs*{\psi}^2,\\
    \EMTensor(e_3,N_{\Sigma})
    \ge{}& \left(1-\sqrt{\frac{8}{9}}\right)
           \left(
           \frac{1}{2}e_3(\tau)\abs*{\nabla\psi}^2
           + \frac{1}{2}e_4(\tau)\abs*{\nabla_3\psi}^2
           \right)
           + \frac{1}{2}e_3(\tau)V\abs*{\psi}^2. 
  \end{align*}
  Again from the properties of $\tau$ in \zcref[cap]{prop:properties-of-global-coordinates} and also from $V\ge - O(\delta_{\Horizon})r^{-2}\abs*{\psi}^2
    \left( \bOne_{r\le r_{\EventHorizonFuture}}+\bOne_{r\ge r_{\CosmologicalHorizonFuture}} \right)$ (using the particular form of $V$ in
  \zcref[noname]{eq:model-problem-gRW}) we can infer the existence of
  a constant $c_0>0$ such that on $\Manifold$,
  \begin{align*}
    \EMTensor(e_4,N_{\Sigma})
    \ge{}& c_0\left(
           \frac{1}{2}e_3(\tau)\abs*{\nabla_4\psi}^2
           + \frac{1}{2}e_4(\tau)\abs*{\nabla\psi}^2
           \right)
           - O(\delta_{\Horizon})r^{-2}\left(\bOne_{r\le r_{\EventHorizonFuture}} + \bOne_{r\le r_{\CosmologicalHorizon}}\right)\abs*{\psi}^2,\\
    \EMTensor(e_3,N_{\Sigma})
    \ge{}& c_0\left(
           \frac{1}{2}e_3(\tau)\abs*{\nabla\psi}^2
           + \frac{1}{2}e_4(\tau)\abs*{\nabla_3\psi}^2
           \right)
           - O(\delta_{\Horizon})r^{-2}\left(\bOne_{r\le r_{\EventHorizonFuture}} + \bOne_{r\le r_{\CosmologicalHorizon}}\right)\abs*{\psi}^2. 
  \end{align*}
  Using the decomposition of $\HawkingVF$ in
  \zcref[noname]{eq:T-R-Z:principal-global:expression:THat-RHat}, we
  have that
  \begin{align}
    \EMTensor\left(\HawkingVF, N_{\Sigma}\right)
    ={}&\frac{\lambdaglo}{2(1+\gamma)}\frac{\abs*{q}^2}{r^2+a^2}\EMTensor(e_3,N_{\Sigma})
         + \frac{1}{2\lambdaglo(1+\gamma)}\frac{\Delta}{r^2+a^2}\EMTensor(e_4,N_{\Sigma})\notag\\
    \ge{}& c_0\left(
           \frac{\lambdaglo}{r^2}\abs*{\nabla_3\psi}^2
           + \frac{\abs*{\Delta}}{\lambdaglo r^2}\abs*{\nabla_4\psi}^2
           + \abs*{\nabla\psi}^2
           \right)\notag \\
       & - O(\delta_{\Horizon})\left(
         \bOne_{r\ge r_{\CosmologicalHorizonFuture}}\abs*{\nabla_4\psi}^2
         + \bOne_{r\le r_{\EventHorizon}}\abs*{\nabla_3\psi}^2
         + \left(\bOne_{r\le r_{\EventHorizon}}
         + \bOne_{r\le r_{\CosmologicalHorizon}}\right)r^{-2} \abs*{\psi}^2
         \right),
         \label{eq:Killing-energy:energy-flux:cosmological}
  \end{align}
  where $\lambdaglo$ is given by \zcref[noname]{eq:lambdaglo:def}.
  Using the Poincare inequality in \zcref[cap]{lemma:Poincare}, and the
  null decomposition of $\KillT$ in
  \zcref[noname]{eq:T-R-Z:principal-outgoing:expression:T-Z}, we have
  that
  \begin{equation*}
    \frac{1}{r^2}\int_S\abs*{\psi}^2
    \lesssim \int_S\abs*{\nabla\psi}^2
    + \frac{\lambdaglo}{ r^2}\int_S\abs*{\nabla_3\psi}^2
    + \frac{\abs*{\Delta}}{\lambdaglo\abs*{q}^2} \abs*{\nabla_4\psi}^2.
  \end{equation*}
  Combined with
  \zcref[noname]{eq:Killing-energy:energy-flux:cosmological}, this
  gives us that there exists some $c_0> 0$ such that
  \begin{equation}
    \label{eq:Killing-energy:energy-flux:cosmological:conclusion}
    \begin{split}
      \int_{\Sigma(\tau)}\EMTensor\left(\HawkingVF, N_{\Sigma}\right)
    &= \frac{\lambdaglo}{2(1+\gamma)}\int_{\Sigma(\tau)}\frac{\abs*{q}^2}{r^2+a^2}\EMTensor\left(e_3, N_{\Sigma}\right)
      + \frac{1}{2\lambdaglo(1+\gamma)}\int_{\Sigma(\tau)}\frac{\Delta}{r^2+a^2}\EMTensor\left(e_4, N_{\Sigma}\right)\notag\\
    &\ge c_0\int_{\Sigma(\tau)}\left(
      \frac{\abs*{\Delta}}{\lambdaglo r^2}\abs*{\nabla_4\psi}^2
      + \frac{\lambdaglo}{r^2}\abs*{\nabla_3\psi}^2
      + \abs*{\nabla\psi}^2
      + r^{-2}\abs*{\psi}^2
      \right)
      - O(\delta_{\Horizon})\EnergyFlux_{r\ge r_{\CosmologicalHorizon}}[\psi]\\
    &\ge c_0\EnergyFlux_{\operatorname{deg}}[\psi](\tau)
      - O(\delta_{\Horizon})\EnergyFlux_{r\ge r_{\CosmologicalHorizon}}[\psi](\tau)
      - O(\delta_{\Horizon})\EnergyFlux_{r\le r_{\EventHorizonFuture}}[\psi](\tau)
      .
    \end{split}          
  \end{equation}
  This yields \zcref[noname]{eq:energy:boundary:Sigma}. We now move onto proving
  \zcref[noname]{eq:energy:boundary:SigmaStar}. Recalling our
  convention for $N_{\SigmaStar}$ from
  \zcref[noname]{eq:hypersurface-normal-def},
  \begin{align*}
    \EMTensor\left(e_4, N_{\SigmaStar}\right)
    ={}& \frac{1}{2}\EMTensor\left(e_4, e_3\right) - \frac{\Delta}{2\abs*{q}^2}\EMTensor(e_4,e_4)
    = \frac{1}{2}\left(\abs*{\nabla\psi}^2 + V\abs*{\psi}^2\right)
         - \frac{\Delta}{2\abs*{q}^2}\abs*{\nabla_4\psi}^2,
\\    \EMTensor\left(e_3, N_{\SigmaStar}\right)
    ={}& -\frac{\Delta}{2\abs*{q}^2}\EMTensor\left(e_3, e_4\right) + \frac{1}{2}\EMTensor(e_3,e_3)
    = -\frac{\Delta}{2\abs*{q}^2}\left(\abs*{\nabla\psi}^2 + V\abs*{\psi}^2\right)
         + \frac{1}{2}\abs*{\nabla_3\psi}^2.
  \end{align*}
  As a result, we have that
  \begin{align*}
    \EMTensor\left(\HawkingVF, N_{\SigmaStar}\right)
    ={}& \frac{1}{2(1+\gamma)}\frac{\abs*{q}^2}{r^2+a^2}\EMTensor\left(e_3,N_{\SigmaStar}\right)
         + \frac{1}{2(1+\gamma)}\frac{\Delta}{r^2+a^2}\EMTensor\left(e_4,N_{\SigmaStar}\right)\\
    ={}& \frac{1}{4(1+\gamma)}\frac{\abs*{q}^2}{r^2+a^2}\abs*{\nabla_3\psi}^2
         - \frac{1}{4(1+\gamma)}\frac{\Delta^2}{\abs*{q}^2\left(r^2+a^2\right)}\abs*{\nabla_4\psi}^2.    
  \end{align*}
  As a result, we have that
  \begin{equation*}
    \int_{\SigmaStar}\EMTensor\left(\HawkingVF,N_{\SigmaStar}\right)
    \gtrsim \int_{\SigmaStar(\tau_1,\tau_2)}\abs*{\nabla_{3}\psi}^2
    -\delta_{\Horizon}\SpacelikeFlux_{\SigmaStar}[\psi](\tau_1,\tau_2), 
  \end{equation*}
  as desired, proving \zcref[noname]{eq:energy:boundary:SigmaStar}. We now move on to proving
  \zcref[noname]{eq:energy:boundary:ACal}. Given our convention for
  $N_{\mathcal{A}}$ from \zcref[noname]{eq:hypersurface-normal-def},
  we can compute that
  \begin{align*}
    \EMTensor\left(e_4, N_{\mathcal{A}}\right)
    ={}&  -\frac{\Delta}{2\abs*{q}^2}\EMTensor\left(e_3, e_4\right)
         + \frac{1}{2}\EMTensor(e_4,e_4)
    ={} -\frac{\Delta}{2\abs*{q}^2}\left(\abs*{\nabla\psi}^2 + V\abs*{\psi}^2\right)
         + \frac{1}{2} \abs*{\nabla_4\psi}^2,
    \\ \EMTensor\left(e_3, N_{\mathcal{A}}\right)
    ={}&  \frac{1}{2}\EMTensor\left(e_3, e_4\right)
         -\frac{\Delta}{2\abs*{q}^2}\EMTensor(e_3,e_3)
    ={} \frac{1}{2}\left(\abs*{\nabla\psi}^2 + V\abs*{\psi}^2\right)
         - \frac{\Delta}{2\abs*{q}^2}\abs*{\nabla_3\psi}^2.
  \end{align*}
  As a result, we have that
  \begin{align*}
    \EMTensor\left(\HawkingVF, N_{\mathcal{A}}\right)
    ={}& \frac{1}{2(1+\gamma)}\frac{\abs*{q}^2}{r^2+a^2}\EMTensor\left(e_4,N_{\mathcal{A}}\right)
         + \frac{1}{2(1+\gamma)}\frac{\Delta}{r^2+a^2}\EMTensor\left(e_3,N_{\mathcal{A}}\right)\\
    ={}& \frac{1}{4(1+\gamma)}\frac{\abs*{q}^2}{r^2+a^2}\abs*{\nabla_4\psi}^2
         - \frac{1}{4(1+\gamma)}\frac{\Delta^2}{\abs*{q}^2\left(r^2+a^2\right)}\abs*{\nabla_3\psi}^2.
  \end{align*}  
  As a result, we have that
  \begin{align*}
    \int_{\mathcal{A}}\EMTensor\left(\HawkingVF, N_{\mathcal{A}}\right)
    \gtrsim \int_{\mathcal{A}(\tau_1,\tau_2)}\abs*{\nabla_{4}\psi}^2 
    - \delta_{\Horizon}\SpacelikeFlux_{\mathcal{A}}[\psi](\tau_1,\tau_2),
  \end{align*}
  as desired, proving \zcref[cap]{lemma:energy:boundary}.
\end{proof}

\subsection{Proof of \zcref[cap]{prop:Killing-energy-estimate}}

We are now ready to prove \zcref[cap]{prop:Killing-energy-estimate}.
\begin{proof}[Proof of Proposition \ref{prop:Killing-energy-estimate}]

  Recall that $\JCurrent{\TAlmostKilling, 0, 0, \mathring{w}}[\psi]$
  satisfies \zcref[noname]{eq:energy:error-terms:div-thm}. Integrating
  \zcref[noname]{eq:energy:error-terms:div-thm} on $\Manifold(\tau_1,\tau_2)$
  and applying \zcref[cap]{lemma:energy:error-terms}, we have that
  \begin{align*}
    &\int_{\Sigma(\tau_2)}\left(
      \JCurrent{\TAlmostKilling, 0, 0, \mathring{w}}[\psi]
      + \mathring{w}\frac{4a(1+\gamma)\cos\theta}{\abs*{q}^2}\abs*{\psi}^2\KillT
      \right)\cdot N_{\Sigma}\\
    & + \int_{\SigmaStar(\tau_1,\tau_2)}\left(
      \JCurrent{\TAlmostKilling, 0, 0, \mathring{w}}[\psi]
      + \mathring{w}\frac{4a(1+\gamma)\cos\theta}{\abs*{q}^2}\abs*{\psi}^2\KillT
      \right)\cdot N_{\SigmaStar}\\
    & + \int_{\mathcal{A}(\tau_1,\tau_2)}\left(
      \JCurrent{\TAlmostKilling, 0, 0, \mathring{w}}[\psi]
      + \mathring{w}\frac{4a(1+\gamma)\cos\theta}{\abs*{q}^2}\abs*{\psi}^2\KillT
      \right)\cdot N_{\mathcal{A}}
    \\
    \lesssim{}&  \int_{\Sigma(\tau_1)}\left(
                \JCurrent{\TAlmostKilling, 0, 0, \mathring{w}}[\psi]
                + \KillT\mathring{w}\frac{4a(1+\gamma)\cos\theta}{\abs*{q}^2}\abs*{\psi}
                \right)\cdot N_{\Sigma}
                + \frac{a}{M}\MorNorm[\psi](\tau_1,\tau_2)\\
              & + \abs*{\int_{\Manifold(\tau_1,\tau_2)}\nabla_{\TAlmostKilling}\psi\cdot N}
                + \abs*{\int_{\Manifold(\tau_1,\tau_2)}\mathring{w}\LeftDual{\psi}\cdot N}.
  \end{align*}
  To estimate the boundary terms, observe that for any vectorfield $n$,
  \begin{align}
    \label{eq:energy:J-T-tilde-q:aux}
    \begin{split}
      \left(
        \JCurrent{\TAlmostKilling, 0, 0, \mathring{w}}[\psi]
        + \KillT\mathring{w}\frac{4a(1+\gamma)\cos\theta}{\abs*{q}^2}\abs*{\psi}^2
        \right)\cdot n
      ={}& \EMTensor\left(\TAlmostKilling, n\right)
           + \mathring{w}\LeftDual{\psi}\cdot \nabla_n\psi
           + \Metric(\KillT, n)\mathring{w}\frac{4a(1+\gamma)\cos\theta}{\abs*{q}^2}\abs*{\psi}^2\\
      ={}& \EMTensor\left(\TAlmostKilling, n\right)
           - O\left(aMr^{-3} + a\Lambda\right)\abs*{\psi}\abs*{\nabla_n\psi}\\
      & - O\left(a^2 M r^{-5} + a^2r^{-2}\Lambda\right)\abs*{\Metric(\KillT, n)}\abs*{\psi}^2,         
    \end{split}
  \end{align}
  where the last line follows from the fact $\abs*{\mathring{w}} \lesssim aMr^{-3} + a\Lambda$ (recall the definition of $\mathring{w}$ in
  \zcref[noname]{eq:energy:w-choice}).

  Recall that from the definition of $\TAlmostKilling$ in
  ~\zcref[noname]{eq:TALmostkilling:def} that
  $\TAlmostKilling=\HawkingVF$ on $\SigmaStar$ and on
  $\mathcal{A}$. As a result, using \zcref[cap]{lemma:energy:boundary},
  we have that 
  \begin{align*}
    \int_{\SigmaStar(\tau_1,\tau_2)}\EMTensor\left(\TAlmostKilling, N_{\SigmaStar}\right)
    \ge{}&\int_{\SigmaStar(\tau_1,\tau_2)}\abs*{\nabla_{3}\psi}^2
           -\delta_{\Horizon}\SpacelikeFlux_{\SigmaStar}[\psi](\tau_1,\tau_2),\\
    \int_{\mathcal{A}(\tau_1,\tau_2)}\EMTensor\left(\TAlmostKilling, N_{\mathcal{A}}\right)
    \gtrsim{}&\int_{\mathcal{A}(\tau_1,\tau_2)}\abs*{\nabla_{4}\psi}^2 
               - \delta_{\Horizon}\SpacelikeFlux_{\mathcal{A}}[\psi](\tau_1,\tau_2).
  \end{align*}
  As a result of \zcref[noname]{eq:energy:J-T-tilde-q:aux}, we have then that
  for $\abs*{a}\ll M$, 
  \begin{align*}
    \int_{\SigmaStar(\tau_1,\tau_2)}\left(
    \JCurrent{\TAlmostKilling, 0, 0, \mathring{w}}[\psi]
    + \KillT\mathring{w}\frac{4a(1+\gamma)\cos\theta}{\abs*{q}^2}\abs*{\psi}^2
    \right)\cdot N_{\SigmaStar}
    \gtrsim{}& \int_{\SigmaStar(\tau_1,\tau_2)}\abs*{\nabla_{3}\psi}^2
               -\delta_{\Horizon}\SpacelikeFlux_{\SigmaStar}[\psi](\tau_1,\tau_2),\\
    \int_{\mathcal{A}(\tau_1,\tau_2)}\left(
    \JCurrent{\TAlmostKilling, 0, 0, \mathring{w}}[\psi]
    + \KillT\mathring{w}\frac{4a(1+\gamma)\cos\theta}{\abs*{q}^2}\abs*{\psi}^2
    \right)\cdot N_{\mathcal{A}}
    \gtrsim{}& \int_{\mathcal{A}(\tau_1,\tau_2)}\abs*{\nabla_{4}\psi}^2 
               -\delta_{\Horizon}\SpacelikeFlux_{\mathcal{A}}[\psi](\tau_1,\tau_2).
  \end{align*}
  We now turn our attention to the $\Sigma_{\tau}$ boundary term. From
  the definition of $\TAlmostKilling$ in \zcref[noname]{eq:TALmostkilling:def},
  \begin{equation*}
    \JCurrent{\TAlmostKilling,0,0}[\psi]\cdot N_{\Sigma}
    = \EMTensor\left(\HawkingVF, N_{\Sigma}\right)
    + \frac{a}{r^2+a^2}\left(\chi_0\left(\delta_{\trap}^{-1}\frac{r-3M}{r}\right)-1\right)\EMTensor\left(\KillPhi, N_{\Sigma}\right).
  \end{equation*}
  Thus,
  \begin{align*}
    \int_{\Sigma(\tau)}\JCurrent{\TAlmostKilling,0,0}[\psi]\cdot N_{\Sigma}
    \ge{}& c_0\EnergyFlux_{\operatorname{deg}}[\psi](\tau)
    - O(a)\int_{\Sigma(\tau)}\abs*{\EMTensor\left(\KillPhi, N_{\Sigma}\right)}\bOne_{\frac{r-3M}{r}\le 2\delta_{\trap}}
    \\& - O(\delta_{\Horizon}) \pth{ \EnergyFlux_{r\le r_{\EventHorizonFuture}}[\psi]
    + \EnergyFlux_{r\ge r_{\CosmologicalHorizonFuture}}[\psi] }. 
  \end{align*}
  As a result, we have that for $a$ sufficiently small, and for
  $\delta$ sufficiently small, 
  \begin{equation*}
    \int_{\Sigma(\tau)}\left(\JCurrent{\TAlmostKilling,0,0,\mathring{w}}[\psi]
      + \KillT \mathring{w}\frac{4a(1+\gamma)\cos\theta}{\abs*{q}^2}\abs*{\psi}^2\right)\cdot N_{\Sigma}
    \ge \frac{c_0}{2}\EnergyFlux_{\operatorname{deg}}[\psi]
    - O(\delta_{\Horizon})\EnergyFlux_{r\le r_{\EventHorizonFuture}}[\psi]
    - O(\delta_{\Horizon})\EnergyFlux_{r\ge r_{\CosmologicalHorizonFuture}}[\psi].
  \end{equation*}
  Thus,
  \begin{equation}
    \label{eq:energy:expanded-final-estimate}
    \begin{split}
      &\EnergyFlux_{\operatorname{deg}}[\psi](\tau_2)
      + \int_{\mathcal{A}(\tau_1,\tau_2)}\abs*{\nabla_{4}\psi}^2 
      + \int_{\SigmaStar(\tau_1,\tau_2)}\abs*{\nabla_{3}\psi}^2\\
      \lesssim {}& \EnergyFlux_{\operatorname{deg}}[\psi](\tau_1)
                   + \delta_{\Horizon}\left(
                   \EnergyFlux_{r\le r_{\EventHorizonFuture}}[\psi](\tau_2)
                   + \EnergyFlux_{r\ge r_{\CosmologicalHorizonFuture}}[\psi](\tau_2)
                   + \SpacelikeFlux_{\mathcal{A}}[\psi](\tau_1,\tau_2)
                   + \SpacelikeFlux_{\SigmaStar}[\psi](\tau_1,\tau_2)
                   \right)\\
                 & + \frac{\abs*{a}}{M}\MorNorm[\psi](\tau_1,\tau_2)
                   + \abs*{\int_{\Manifold(\tau_1,\tau_2)}\nabla_{\TAlmostKilling}\psi\cdot N}
                   + \abs*{\int_{\Manifold(\tau_1,\tau_2)}\mathring{w}\LeftDual{\psi}\cdot N}.
    \end{split}
  \end{equation}  
  Then using the fact that
  $\abs*{\mathring{w}}\lesssim a M r^{-3} + a\Lambda$ and
  Cauchy-Schwarz to control $\mathring{w}\LeftDual{\psi}\cdot N$ we obtain
  \begin{equation*}
    \begin{split}
      &\EnergyFlux_{\operatorname{deg}}[\psi](\tau_2)
        + \int_{\mathcal{A}(\tau_1,\tau_2)}\abs*{\nabla_{4}\psi}^2 
        + \int_{\SigmaStar(\tau_1,\tau_2)}\abs*{\nabla_{3}\psi}^2\\
      \lesssim {}& \EnergyFlux_{\operatorname{deg}}[\psi](\tau_1)
                   + \delta_{\Horizon}\left(
                   \EnergyFlux_{r\le r_{\EventHorizonFuture}}[\psi](\tau_2)
                   + \EnergyFlux_{r\ge r_{\CosmologicalHorizonFuture}}[\psi](\tau_2)
                   + \SpacelikeFlux_{\mathcal{A}}[\psi](\tau_1,\tau_2)
                   + \SpacelikeFlux_{\SigmaStar}[\psi](\tau_1,\tau_2)
                   \right)\\
      & + \frac{\abs*{a}}{M}\MorNorm[\psi](\tau_1,\tau_2)
        + \abs*{\int_{\Manifold(\tau_1,\tau_2)}\nabla_{\TAlmostKilling}\psi\cdot N}
        + \int_{\Manifold(\tau_1,\tau_2)}\abs*{N}^2,
    \end{split}
  \end{equation*}  
  as desired.   
  This concludes the proof of \zcref[cap]{prop:Killing-energy-estimate}.
\end{proof}

\section{Setup for frequency analysis}
\label{sec:setup-freq-analysis}

In order to capture the frequency-dependent behavior of the trapped
set in \KdS, we use pseudo-differential multipliers to derive the
relevant Morawetz estimates. We briefly review some basic
pseudodifferential theory that will be used in what follows. For a
more thorough reference, we refer the reader to
\cite{hormanderAnalysisLinearPartial2007,hintzIntroductionMicrolocalAnalysis2025}

\subsection{Pseudo-differential operators on vector bundles}
\label{sec:calculus}

We briefly recall some auxiliary notions.
\begin{definition}[Left and right projections]
  \label{def:left-right-projection}
  We define the \emph{left projection} and the \emph{right projection} by
  \begin{align}
    \label{def:left-projection}
      \pi_L: \mathcal{M}\times \mathcal{M}&\to \mathcal{M}, \quad \pi_L(p,q)= p,
    \\
    \label{def:right-projection}
      \pi_R: \mathcal{M}\times \mathcal{M}&\to \mathcal{M},\quad \pi_R(p,q)= q.  
  \end{align}
\end{definition}

\begin{definition}[Right density bundle]
  \label{def:right-density-bundle}
  We define the \emph{right density bundle} over $\mathcal{M}$,
  denoted $\Omega_R = \Omega_R\mathcal{M}$ by
  \begin{equation}
    \label{eq:right-density-bundle:def}
    \Omega_R\vcentcolon=\pi_R^{*}(\Omega \mathcal{M}).
  \end{equation}
\end{definition}

\begin{theorem}[Schwartz kernel theorem] %
  Every continuous linear map
  $C_0^{\infty}(\mathcal{M}, E\otimes
  \Omega^{\frac{1}{2}}\mathcal{M})\to \mathcal{D}'(\mathcal{N},
  F\otimes \Omega^{\frac{1}{2}}\mathcal{N})$ has a kernel
  $K\in \mathcal{D}'(Y\times X, \Hom(E,F)\otimes
  \Omega^{\frac{1}{2}}(Y\times X))$, where $\Hom(E,F)_{y,x}$ denotes
  the space of linear maps from $E_x\to F_y$.
\end{theorem}

We first review the construction of a pseudo-differential operator on
a tensor bundle (such as the $(e_3,e_4)$-horizontal tensor
bundle)\footnote{Recall that the $(e_3,e_4)$-horizontal tensor bundle
  $\realHorkTensor{2}(\Manifold)$ is not integrable for $a\neq0$. See \zcref[cap]{sec:nonintegrable-structures}.}. These are
natural applications of the standard definition of pseudo-differential
operators acting on sections of vector bundles (see Chapter 18 of
\cite{hormanderAnalysisLinearPartial2007}, or Section 5.7 of
\cite{hintzIntroductionMicrolocalAnalysis2025}).

\begin{definition}[Definition of pseudo-differential operators on
  tensor bundles. See Definition 18.1.32 in
  \cite{hormanderAnalysisLinearPartial2007}.]
  \label{def:PsiDO:pseudo-diff-on-tensor-bundle}
  Let $\pi_E:E\to \mathcal{M}$ and $\pi_F:F\to \mathcal{M}$ be complex
  $C^{\infty}$ tensor bundles over a $C^{\infty}$ manifold $\mathcal{M}$ of rank
  $N_E$, $N_F$ respectively. Then $A \in \OpClass^m(\mathcal{M}; E,F),$ the
  space of \emph{pseudo-differential operators of order $m$ from
    sections of $E$ to sections of $F$} is the space of continuous
  linear maps $A: C_0^{\infty}(\mathcal{M},E) \to C^{\infty}(\mathcal{M}, F)$ such that
  \begin{enumerate}
  \item If $\phi,\psi\in C^{\infty}(\mathcal{M})$ have disjoint
    supports, then there exists $K\in C^{\infty}(M^2;\pi^{*}_LF\otimes \pi^{*}_R(E^{*}\otimes \mathcal{M}) )$
    such that
    \begin{equation*}
      \phi A\psi = K
    \end{equation*}
  \item Let $U\subset \mathcal{M}$ be an open set and let
    $G: U\mapsto G(U)\subset \Real^n$ be a diffeomorphism. For every
    open $Y\subset \mathcal{M}$, where $E$ and $F$ are trivialized by
  \begin{equation*}
    \phi_E: \pi_E^{-1}(Y)\to G(Y) \times \Complex^{N_E}, \qquad
    \phi_F: \pi_F^{-1}(Y) \to G(Y) \times  \Complex^{N_F},
  \end{equation*}
  if $\psi\in C_0^{\infty}(Y)$, then there exists an $N_F\times N_E$
  matrix of order-$m$ pseudo-differential operators
  $A_{ij}\in \SymClass^m_c(G(Y))$ such that on $Y$
  \begin{equation*}
    \psi A(\psi u)(x)_i = \sum_{j=1}^{N_{E}}G^{*}\left(A_{ij}\left((G^{-1})^{*}u\right)_j\right),
    \qquad
    u\in C_0^{\infty}(\mathcal{M},E),\quad i=1,\cdots, N_F.
  \end{equation*}
  \end{enumerate}
  If $A$ is a pseudo-differential operator of order $m$ from sections
  of $E$ to sections of $F$, then we write $A \in \OpClass^m(\mathcal{M}; E,F)$.
\end{definition}

In practice, we will only use pseudo-differential operators with
Schwartz kernels supported near the trapped set of \KdS. To this end,
it is convenient to recall the definition of properly supported
pseudo-differential operators.
\begin{definition}[Properly supported pseudo-differential operators]
  A pseudo-differential operator $A\in \OpClass^m(\Manifold)$ with
  Schwartz kernel
  $K\in \mathcal{D}'(\mathcal{M}\times \mathcal{M}; \Omega_R)$ is said
  to be \emph{properly supported} if both $\pi_L: \supp K \to \mathcal{M}$ and
  $\pi_R:\supp K\to \mathcal{M}$ are proper. 
\end{definition}

We can also naturally extend the definition of a symbol to
pseudo-differential operators acting on sections of tensor bundles.
\begin{definition}
  \label{def:PsiDO:symbol-tensor-bundle}
  Let $\pi:T^{*}\Manifold\to \Manifold$ denote the projection.
  Given a rank $k$ vector bundle $G\to \Manifold$, 
  we denote by $S^m(T^{*}\Manifold, \pi^{*}G)  \subset C^{\infty}(T^{*}\Manifold, \pi^{*}G)$ the space of all smooth functions which in local coordinates and in
  a trivialization of $G$ are matrices of symbols on $\Real^n$ of
  order $m$.
\end{definition}

In other words, $a\in S^m(T^{*}\Manifold, \pi^{*}G)$ if
  $a$ has the following property. For each coordinate chart
  $\phi:U\to \phi(U)\subset \Real^n$ on $\mathcal{M}$ on which $G$ is
  trivial with trivialization
  $\tau: \pi^{-1}(U)\to \phi(U)\times
  \Real^k$, setting
  \begin{equation*}
    b(x,\xi) := (\tau^{-1})^{*}(a\vert_{\pi^{-1}(U)})(x,\xi)
    = a(\tau^{-1}(x,\xi))\in C^\infty(\phi(U)\times \Real^k),
  \end{equation*}
  then for any
  $\chi\in C^{\infty}_0(F(U))\subset C^{\infty}_0(\Real^d)$, we have
  that $\chi(x)b(x,\xi)\in S^m(\Real^n;\Real^k)$. 

\begin{remark}
  Invariantly, if $A$ is an order-$m$ pseudo-differential operators
  acting between sections of vector bundles $E\to \mathcal{M}$ and
  $F\to \mathcal{M}$, then
  $\sigma^m(A)\in (S^m/S^{m-1})(T^{*}\mathcal{M}; \pi^{*}\Hom(E,F))$,
  where $\pi:T^{*}\mathcal{M}\to \mathcal{M}$.
\end{remark}

For the coordinate-inclined, symbols can be characterized and
described by coordinate charts. 
\begin{lemma}[Symbols and charts]
  \label{lemma:PsiDO:symbols-and-charts}
  In the notation of \zcref[cap]{def:PsiDO:symbol-tensor-bundle},
  for $\chi\in C^{\infty}_0(\phi(U))$, $b\in S^m(\Real^n;\Real^k)$,
  \begin{equation*}
    \tau^{*}(\chi b)\in S^m(T^{*}\mathcal{M};\pi^{*}G). 
  \end{equation*}
  Moreover, let $a\in C^{\infty}(G)$. Then, $a\in S^m(T^{*}\mathcal{M};\pi^{*}G)$ if and only
  if there exist symbols $b_i\in S^m(\Real^d; \Real^k)$ and
  $\chi_i\in C^{\infty}_0(\phi_i(U_i))$ such that
  \begin{equation*}
    a = \sum_i\left( \tau_i^{-1} \right)^{*}(\chi_ib_i). 
  \end{equation*}
\end{lemma}

It will be convenient for computational purposes to factor out the
half density bundle.

\begin{definition}[Density bundles]
  \label{def:density-bundle}
  Let $\sigma\in \Real$. The \emph{$\alpha$-density bundle on
    $\mathcal{M}$, denoted by $\Omega^{\alpha}\mathcal{M}$} is the vector bundle consisting of
  sections that, expressed in local coordinates $\{x_{j}\}_{j=1}^{n}$,
  are functions $u$ such that given two choices of local coordinates,
  $\{x_{j}\}_{j=1}^{n}$ and $\{x'_{j}\}_{j=1}^{n}$, satisfy the
  transformation law
  \begin{equation}
    \label{eq:density-bundle}    
    u'\abs*{dx'}^{\alpha} = u\abs*{dx}^{\alpha}.
  \end{equation}
  An equivalent formulation is the definition
  \begin{equation*}
    \Omega^{\alpha}\mathcal{M}\vcentcolon = \curlyBrace*{\xi: \wedge^{n}\mathcal{M}\to \Real: \omega(\mu v) = \abs*{\mu}^{\alpha}\omega(v), v\in \wedge^n\mathcal{M}, \mu\in\Real}.
  \end{equation*}
  
\end{definition}

\begin{proposition}[Compositions and adjoints]
  \label{prop:PsiDO-over-bundle:composition-and-adjoint}
  Let $\mathcal{M}$ be a smooth manifold and let $E,F,G$ denote three
  vector bundles over $\mathcal{M}$.
  \begin{enumerate}
  \item Let $A\in \Psi^m(\mathcal{M}; F,G)$ and
    $B\in \Psi^{m'}(\mathcal{M}; E,F)$ with at least one $A$ and $B$ properly
    supported. Then $A\circ B\in \Psi^{m+m'}(\mathcal{M};E,G)$ and
    $\sigma^{m+m'}(A\circ B) = \sigma^m(A)\circ \sigma^{m'}(B)$.
  \item Any
    $A\in \Psi^m(\mathcal{M};
    E\otimes\Omega^{\frac{1}{2}}\mathcal{M},F\otimes\Omega^{\frac{1}{2}}\mathcal{M})$
    has an adjoint
    $A^{*}\in \Psi^m(\mathcal{M}; F^{*}\otimes
    \Omega^{\frac{1}{2}}\mathcal{M}, E^{*}\otimes
    \Omega^{\frac{1}{2}}\mathcal{M})$ with symbol $a^{*}$, and moreover, if $A$
    is properly supported, then $A^{*}$ is also properly supported.
  \end{enumerate}
\end{proposition}

For the purpose of constructing the desired pseudo-differential
multipliers to prove a Morawetz estimate, we will be primarily
interested in principally scalar pseudo-differential operators. In practice, we will work with principally scalar operators.
\begin{definition}[Principally scalar operators]
  \label{def:principally-scalar-operators}
  Let $m\in\Real$, $\mathcal{M}$ be a smooth manifold, and $E$ a
  vector bundle over $\mathcal{M}$, and
  $A\in \Psi^m(\mathcal{M};E) = \Psi^m(\mathcal{M};E, E)$. Then we say
  that $A$ is \emph{principally scalar} if its principal symbol is
  multiplication by scalars on the fibers of E i.e, if there exists a
  symbol $a\in S^m(T^{*}\mathcal{M})$ such that $\sigma^m(A)(x,\xi) = a(x,\xi)\Identity_{E_x}$.
\end{definition}

\begin{remark}
  Both the covariant derivative $\HorCovDeriv_X$, and the horizontal
  wave operator $\WaveOpHork{2}$ are principally scalar differential
  operators. 
\end{remark}

\subsection{Mixed pseudo-differential operators}
\label{sec:mixed-PDO-class}

In proving our Morawetz estimates, it will be convenient to deal with
a class of pseudo-differential multipliers which are
pseudo-differential in space but differential in time. This class of
operators will allow for some flexibility to account for the
frequency-dependent behavior of trapping for \KdS, while also allowing
the multiplier to produce local energy decay estimates. These types of
pseudo-differential operators have been used previously to prove
Morawetz estimates. We refer the reader to
\cite{maEnergyMorawetzEstimatesWave2024,
  lindbladLocalEnergyEstimate2020, tataruLocalEnergyEstimate2011} for
in-depth introductions to such operators. Below we review their basic
construction and definitions, as well as some basic properties of
these mixed pseudo-differential operators.

\begin{definition}[$x_0$-tangential symbols on $\Real^d$]
  For $m\in \Real$, let
  $\TanSymClass{m}\subset C^{\infty}(\Real^{d-1}\times \Real^{d-1})$, the
  set of $x_0$-tangential symbols of order $m$, consist of functions
  $a(x, \xi): \Real^{d-1}\times \Real^{d-1}\to
  \Complex$ such that for all multi-indices $\alpha,\beta$,
  \begin{equation*}
    \abs*{\partial_x^{\alpha}\partial_{\xi}^{\beta}a(x,\xi)}\le C_{\alpha,\beta}\bangle*{\xi}^{m-\abs*{\beta}},
  \end{equation*}
  for all $x\in \Real_{x_0}\times \Real^{d-1}$. 
\end{definition}
\begin{remark}
  When compared to the mixed symbols introduced in
  \cite{maEnergyMorawetzEstimatesWave2024}, it should be noted that
  our symbols, much like those used in
  \cite{tataruLocalEnergyEstimate2011} are independent of both $x_0$
  and $\xi_0$, rather than just $\xi_0$. This has to do with the
  stationarity of Kerr(-de Sitter), and substantially simplifies our
  analysis. 
\end{remark}

\begin{definition}[$x_0$-mixed symbols on $\Real^d$]
  \label{def:mixed-sym-class}
  For $m\in \Real$, $n\in \Natural$, we define the class
  $\MixedSymClass{m}{n}(\Real^d)\vcentcolon=\MixedSymClass{m}{n}[x_0](\Real^d)\subset
  C^{\infty}(\Real^d\times \Real^d)$ of \emph{mixed symbols of order
    $(m,n)$}, such that, decomposing $\xi = (\xi_0,\xi')$,
  \begin{equation*}
    a(x,\xi) = \sum_{j=0}^na_{m-j}(x,\xi')\xi_0^j, \qquad a_{m-j}\in \TanSymClass{m-j}(\Real^d).
  \end{equation*}
\end{definition}

We next define the mixed symbols on a manifold.
\begin{definition}[Mixed symbols on a vector bundle]
  \label{def:mixed-symbols-bundle}
  Let $\pi:T^{*}\Manifold\to \Manifold$ denote the projection.
  Given a rank $k$ vector bundle $G\to \Manifold$, 
  we denote by
  \begin{equation*}
    \MixedSymClass{m}{n}(T^{*}\Manifold, \pi^{*}G)  \subset C^{\infty}(T^{*}\Manifold, \pi^{*}G)
  \end{equation*}
  the space of all smooth functions which in local coordinates and in
  a trivialization of $G$ are matrices of symbols in
  $\MixedSymClass{m}{n}(\Real^d)$. 
\end{definition}

In other words,
  $a\in \MixedSymClass{m}{n}(T^{*}\Manifold, \pi^{*}G)$ if $a$ has the
  following property. For each coordinate chart
  $F:U\mapsto F(U) = \Real_{x_0}\times \widetilde{F}(U)\subset
  \Real\times \Real^{d-1}$ on $\mathcal{M}$ on which $G$ is trivial
  with trivialization
  $\tau: \pi^{-1}(U)\mapsto U\times \Real^k\mapsto
  \Real_{x_0}\times\widetilde{F}(U)\times
  \Real^k$, setting
  \begin{equation*}
    b(x,\xi) \vcentcolon= (\tau^{-1})^{*}(a\vert_{\pi^{-1}(U)})(x,\xi)
    = a(\tau^{-1}(x,\xi))\in C^\infty(F(U)\times \Real^k),
  \end{equation*}
  then for any
  $\chi\in C^{\infty}_0(F(U))\subset C^{\infty}_0(\Real^d)$, we have
  that $\chi(x)b(x,\xi)\in \MixedSymClass{m}{n}(\Real^n;\Real^k)$.

\subsection{Weyl calculus}
\label{sec:Weyl-calculus}

It will be convenient to quantize our symbols by using the Weyl
quantization as opposed to the standard quantization. While the
properties of the Weyl calculus presented in this subsection are
standard, we review them briefly for the sake of clarity. For a more
in-depth reference on the Weyl calculus, we refer the reader to
Section 18.5 of \cite{hormanderAnalysisLinearPartial2007}.

\begin{definition}[Weyl quantization]
  \label{def:Weyl-quantization} Given a symbol $s\in S^m(\Real^n)$,
  the \emph{Weyl quantization of $s$} is the pseudo-differential operator of
  order $m$ given by
  \begin{equation*}
    \WeylQ{\phi}(x) = \frac{1}{(2\pi)^n}\iint
    s\left(\frac{x+y}{2}, \xi\right)e^{{\ImagUnit(x-y)\cdot
        \xi}}\phi(y)\,dyd\xi.
  \end{equation*}
\end{definition}

We can also define the Weyl quantization of a mixed symbol.
\begin{definition}[Weyl quantization of mixed symbol]
  Let $(m,n)\in \Real\times\Natural$, and $a\in \MixedSymClass{m}{n}(\Real^d)$
  with
  \begin{equation*}
    a(x,\xi) = \sum_{j=0}^na_{m-j}(x,\xi')\xi_0^j,\qquad a_{m-j}\in \TanSymClass{m-j}(\Real^d). 
  \end{equation*}
  Then $\WeylQ{a}$, the \emph{Weyl quantization of $a$}, is given by
  \begin{equation*}
    \WeylQ{a}=\sum_{j=0}^n\sum_{k=0}^j2^{-k}\binom{j}{k}\WeylQ{D_{x^0}^ka_{m-j}}D^{j-k}_{x^0},
  \end{equation*}
  where $\WeylQ{D_{x^0}^ka_{m-j}}$ is the Weyl quantization in
  $\Real^{n-1}$ of the $x^0$-tangential symbol $D_{x^0}^ka_{m-j}$.
\end{definition}

The Weyl quantization enjoys several convenient properties.
\begin{lemma}[Proposition 5.15 in \cite{maEnergyMorawetzEstimatesWave2024}]
  \label{lemma:Weyl-quantization:props}
  The following properties hold.
  \begin{enumerate}
  \item The Weyl quantization of a real symbol is self-adjoint with
    respect to the complex inner product.
  \item If $s\in \MixedSymClass{m_1}{n_1}(\Real^k)$,
    $q\in \MixedSymClass{m_2}{n_2}(\Real^k)$ then
    \begin{align*}
      \WeylQ{sq} &= \frac{1}{2}\left(\WeylQ{s}\WeylQ{q} + \WeylQ{q}\WeylQ{s}\right)\in \MixedOpClass{m_1+m_2-2}{n_1+n_2}(\Real^k),
      \\ \WeylQ{s}\WeylQ{q} &= \WeylQ{sq} + \frac{1}{2}\WeylQ{\PoissonB{s,q}}\mod \MixedOpClass{m_1+m_2-2}{n_1+n_2}(\Real^k).
    \end{align*}
  \item If $s\in S^m(\Real^k)$, $q\in S^n(\Real^k)$, then
    \begin{align*}
      \left[\WeylQ{p},\WeylQ{s}\right]  = \frac{1}{2\ImagUnit}\WeylQ{\PoissonB*{p,s}}   \in \MixedOpClass{m_1+m_2-1}{n_1+n_2}(\Real^k),
      \\\left[\WeylQ{p},\WeylQ{s}\right] - \frac{1}{2\ImagUnit}\WeylQ{\PoissonB*{p,s}}   \in \MixedOpClass{m_1+m_2-3}{n_1+n_2}(\Real^k).
    \end{align*}
  \end{enumerate}
\end{lemma}

In practice, we will use the Weyl quantization to construct
principally scalar mixed-pseudodifferential operators acting on
sections of a vector bundle. We recall that there is no canonical way
of quantizing symbols in $S^m(T^{*}\mathcal{M}; \pi^{*}\Hom(E,F))$ to
pseudodifferential operators in $\Psi^m(\mathcal{M};E, F)$. However,
the following construction will be useful in what follows. We start with the (again useful, but not canonical) Weyl quantization
of symbols in $\SymClass^m(T^{*}\mathcal{M})$.
\begin{definition}[Weyl quantization of pseudodifferential symbols on manifolds]
  \label{def:WeylQ:manifold}
  Let $(m, n)\in\Real\times \Natural$, $\mathcal{M}$ a smooth
  $d$-dimensional manifold, $E\to \mathcal{M}$ a vector bundle over
  $\mathcal{M}$, and $a\in S^m(T^{*}\Manifold)$. Then, fix a locally
  finite open cover $\bigcup_iU_i= \mathcal{M}$ by coordinate charts
  $\phi_i:U_i\to \phi_i(U_i)\subset \Real^d$, where $\overline{U}_i$ are
  compact, and let $\{\chi_i\}_i,\chi_i\in C^{\infty}_0(U_i)$ be a
  partition of unity subordinate to $\{U_i\}$. Then, for
  $\widetilde{\chi}_i\in C^{\infty}_0(U_i)$, with
  $\widetilde{\chi}_i=1$ in a neighborhood of $\supp\chi_i$, we define
  \begin{equation*}
    \WeylQ{a} := \sum_iF^{*}_i\WeylQ{(\phi_i^{-1})^{*}(\chi_ia)}\widetilde{\chi}_i(\phi_i^{-1})^{*},
  \end{equation*}
  where $\WeylQ{\cdot}: S^m(T^{*}\Real^d)\to \Psi^m(\Real^d) $ is as
  defined before.
\end{definition}

This naturally extends to the definition for the Weyl
quantization of mixed symbols
$\MixedSymClass{m}{n}(T^{*}\mathcal{M};\pi^{*}\Hom(E,F))$.
\begin{definition}[Weyl quantization of mixed-pseudodifferential
  operators on sections of a vector bundle]
  Let $(m, n)\in\Real\times \Natural$, $\mathcal{M}$ a smooth
  $d$-dimensional manifold with a global coordinate $x_0$,
  $E\to \mathcal{M}$ a vector bundle over $\mathcal{M}$, and
  $a\in \MixedSymClass{m}{n}(\Real^d)$. Then, fix a locally finite
  open cover $\bigcup_iU_i= \mathcal{M}$ by coordinate charts
  $\phi_i:U_i\to \phi_i(U_i) = \Real_{x_0}\times \widetilde{\phi}_i\subset
  \Real_{x_0}\times\Real^{d-1}$, where $\overline{U}_i$ are compact,
  and let $\{\chi_i\}_i,\chi_i\in C^{\infty}_0(U_i)$ be a partition of
  unity subordinate to $\{U_i\}$. Moreover, let
  $\tau_i: \pi^{-1}(U_i)\to \Real_{x_0}\times \widetilde{\phi}_i\times
  \Real^k$.

  Then, for
  $\widetilde{\chi}_i\in C^{\infty}_0(U_i)$, with
  $\widetilde{\chi}_i=1$ in a neighborhood of $\supp\chi_i$,
  given $a\in \MixedSymClass{m}{n}(T^{*}\mathcal{M};\pi^{*}\Hom(E,F))$,
  we define
  \begin{equation*}
    \WeylQ{a} := \sum_i\tau^{*}_i\WeylQ{(\tau_i^{-1})^{*}(\chi_ia)}\widetilde{\chi}_i(\tau_i^{-1})^{*}.
  \end{equation*}
\end{definition}

\begin{proposition}[Properties of the quantization map over bundles]
  \label{prop:Weyl-quantization:props:bundle}
  Let $\mathcal{M}$ be a $d$-dimensional manifold, and $E$ and $F$ be
  $k_E$ and $k_F$ rank respectively vector bundles over $\mathcal{M}$.
  The map $\WeylQ{\cdot}$ has the following basic properties.
  \begin{enumerate}
  \item $\WeylQ{\cdot}$ is a continuous, linear maps that takes
    values in the subspace of $\Psi^m(\mathcal{M})$ of properly
    supported operators;
  \item  $\WeylQ{1} = \Identity$;
  \item $\WeylQ{\cdot}: S^m(T^{*}\mathcal{M}; \pi^{*}\Hom(E,F))\to \MixedOpClass{m}{n}(\mathcal{M}; E, F)$
    is surjective modulo $\MixedSymClass{-\infty}{n}(\mathcal{M})$. That is,
    \begin{equation*}
      \MixedOpClass{m}{n}(\mathcal{M}; E, F) = \WeylQ{\MixedSymClass{m}{n}(T^{*}\mathcal{M}; \pi^{*}\Hom(E,F)} + \MixedOpClass{-\infty}{n}(\mathcal{M}; E, F).
    \end{equation*}    
  \end{enumerate}
  Moreover, the Weyl quantization also satisfies the following properties. 
  \begin{enumerate}
  \item The Weyl quantization of a real symbol is self-adjoint with
    respect to the complex inner product.
  \item If
    $s\in \MixedSymClass{m_1}{n_1}(T^{*}\mathcal{M};
    \pi^{*}\Hom(E,F))$,
    $q\in \MixedSymClass{m_2}{n_2}(T^{*}\mathcal{M};
    \pi^{*}\Hom(E,F))$, at least one of which is principally scalar,
    then
    \begin{align*}
      \WeylQ{sq} & = \frac{1}{2}\left(\WeylQ{s}\WeylQ{q} + \WeylQ{q}\WeylQ{s}\right)\mod \MixedOpClass{m_1+m_2-2}{n_1+n_2}(\mathcal{M}, E, F),
      \\ \WeylQ{s}\WeylQ{q} & = \WeylQ{sq} + \frac{1}{2}\WeylQ{\PoissonB{s,q}}\mod \MixedOpClass{m_1+m_2-2}{n_1+n_2}(\mathcal{M}, E, F).
    \end{align*}
  \item If
    $s\in \MixedSymClass{m_1}{n_1}(T^{*}\mathcal{M};
    \pi^{*}\Hom(E,F))$,
    $q\in \MixedSymClass{m_2}{n_2}(T^{*}\mathcal{M};
    \pi^{*}\Hom(E,F))$, at least one of which is principally scalar,
    then $\left[\WeylQ{q},\WeylQ{s}\right]\in\MixedOpClass{m_1+m_2-1}{n_1+n_2}(\mathcal{M}, E, F)$ and 
    \begin{align*}
    \left[\WeylQ{q},\WeylQ{s}\right] & - \frac{1}{2\ImagUnit}\WeylQ{\PoissonB*{q,s}} \in \MixedOpClass{m_1+m_2-3}{n_1+n_2}(\mathcal{M}, E, F).
    \end{align*}
  \end{enumerate}
\end{proposition}
\begin{proof}
  The proof follows from generalizing \zcref[cap]{lemma:Weyl-quantization:props}.
\end{proof}

The Weyl calculus also enjoys the following Garding-type inequality.
\begin{lemma}[]
  \label{thm:fefferman-phong-tataru}
  Let $a_j, b$ be real mixed-symbols
  in $\MixedSymClass{1}{1}(T^{*}\mathcal{M};\pi^{*}\End(E))$ that
  have kernels supported in $\mathcal{M}_{\operatorname{trap}}$ such
  that $\abs*{b}\le \sum \abs*{a_j}$, then
  \begin{equation}
    \label{eq:fefferman-phong-tataru}
    \norm*{\WeylQ{b}(x,D)\psi}_{L^2(\mathcal{M}(\tau_1,\tau_2))}
    \lesssim \sum\norm*{\WeylQ{a_j}(x,D)\psi}_{L^2(\mathcal{M}(\tau_1,\tau_2))}
    + \norm*{\psi}_{L^2(\mathcal{M}(\tau_1,\tau_2))}.
  \end{equation}
\end{lemma}

\subsection{Commutation with the wave operator}
\label{sec:pseudo-diff:wave-commutation}

In this section, we present the main equivalent of the divergence
theorem in \zcref[cap]{sec:EM-Tensor} that we will use when
considering pseudo-differential multipliers. We recall that throughout
this section we will use the terms skew-symmetric and symmetric to be
defined with respect to the inner product structure in 
\zcref[cap]{def:hor-inner-product}. For clarity of the exposition, we define
$P_i\in
\TanSymClass{i}(T^{*}\mathcal{M};\pi^{*}\End(\realHorkTensor{2}(\mathcal{M})))$ for 
$i=0,1,2$ as the symmetric operators such that the following holds
\begin{equation*}
  \WaveOpHork{2} = \sum_{i=0}^2P_i\left(\frac{1}{\ImagUnit}\HorCovDeriv_t\right)^{2-i}.
\end{equation*}
We now introduce the natural inner product associated to
$\realHorkTensor{2}$. 
\begin{definition}
  \label{def:hor-inner-product}
  Let $\phi, \psi\in\realHorkTensor{2}$ be a horizontal 2-tensor. Then
  we define the inner product
  \begin{equation*}
    \bangle*{\phi, \psi}_{L^2(\Manifold(\tau_1,\tau_2))} \vcentcolon= \int_{\Manifold(\tau_1,\tau_2)}\phi\cdot \overline{\psi},
  \end{equation*}
  where $\phi\cdot \overline{\psi}\vcentcolon= \phi_{ab}\overline{\psi}^{ab}$.
\end{definition}

\begin{lemma}
  \label{lemma:PsiDO-divergence-theorem}
  Let $\psi\in \realHorkTensor{2}(\Manifold)$ such that
  $\supp\psi\subset \Manifold_{\operatorname{trap}}$ and let
  \begin{equation*}
    \widetilde{\MorawetzVF}\vcentcolon=
    \widetilde{\MorawetzVF}_1
    + \widetilde{\MorawetzVF}_0\HorCovDeriv_t, \qquad \widetilde{\MorawetzVF}_i\in \TanSymClass{i}(T^{*}\Manifold;
  \pi^{*}\End(\realHorkTensor{2}(\mathcal{M}))),
  \end{equation*}
  be a first-order anti-symmetric properly supported
  pseudodifferential operator on $\realHorkTensor{2}(\Manifold)$. Moreover, let
  \begin{equation*}
    \widetilde{\MorawetzLagrangeCorr}\vcentcolon=
    \widetilde{\MorawetzLagrangeCorr}_0
    + \widetilde{\MorawetzLagrangeCorr}_{-1}\HorCovDeriv_t, \qquad \widetilde{\MorawetzLagrangeCorr}_i\in \TanSymClass{i}(T^{*}\Manifold;
  \pi^{*}\End(\realHorkTensor{2}(\mathcal{M}))),
  \end{equation*}
  be a zero-order symmetric properly supported pseudodifferential operator on
  $\realHorkTensor{2}(\Manifold)$. Then,
  \begin{equation}
    \label{eq:frequency-setup:wave-commutation-formula}
    \bangle*{\WaveOpHork{2}\psi, \left(\widetilde{\MorawetzVF}+\widetilde{\MorawetzLagrangeCorr}\right)\psi}_{L^2(\Manifold(\tau_1,\tau_2))}
    = \KCurrentPert{\widetilde{\MorawetzVF}, \widetilde{\MorawetzLagrangeCorr}}[\psi](\tau_1,\tau_2)
    + \JCurrentPert{\widetilde{\MorawetzVF}, \widetilde{\MorawetzLagrangeCorr}}[\psi](\tau_2)
    - \JCurrentPert{\widetilde{\MorawetzVF}, \widetilde{\MorawetzLagrangeCorr}}[\psi](\tau_1),
  \end{equation}
  where
  \begin{align}
    2\KCurrentPert{\widetilde{\MorawetzVF}, \widetilde{\MorawetzLagrangeCorr}}[\psi](\tau_1,\tau_2)
    \vcentcolon={}& \bangle*{\left[\widetilde{\MorawetzVF}, \WaveOpHork{2}\right]\psi,\psi}_{L^2(\Manifold(\tau_1,\tau_2))}
          + \bangle*{\left(\widetilde{\MorawetzLagrangeCorr}\WaveOpHork{2} + \WaveOpHork{2}\widetilde{\MorawetzLagrangeCorr}\right)\psi, \psi}_{L^2(\Manifold(\tau_1,\tau_2))},
          \label{eq:KCurrentPert-def}\\
    \JCurrentPert{\widetilde{\MorawetzVF}, \widetilde{\MorawetzLagrangeCorr}}[\psi](\tau)
    \vcentcolon=&\Re\bangle*{n_{\Sigma(\tau)}\psi, \widetilde{\MorawetzVF}\psi}_{L^2(\Sigma(\tau))}
        +\Re\bangle*{n_{\Sigma(\tau)}\psi, \widetilde{\MorawetzLagrangeCorr}\psi}_{L^2(\Sigma(\tau))}\notag \\
       &+ \Re\bangle*{\WaveOpHork{2}\psi,\Metric^{-1}_{M,a,\Lambda}(n_{\Sigma(\tau)},d\tau)\widetilde{\MorawetzVF}_0\psi}_{L^2(\Sigma(\tau))}
        ,
        \label{eq:JCurrentPert-def}
  \end{align}
  where $n_{\Sigma(\tau)}$ here denotes the standard unit normal to $\Sigma(\tau)$.
  Critically,
  $\JCurrentPert{\widetilde{\MorawetzVF},
    \widetilde{\MorawetzLagrangeCorr}}[\psi](\tau)$ satisfies the following inequality
  \begin{equation}
    \label{eq:JCurrentPert:basic-control}
    \abs*{\JCurrentPert{\widetilde{\MorawetzVF}, \widetilde{\MorawetzLagrangeCorr}}[\psi](\tau)}
    \lesssim \EnergyHorizonDeg[\psi](\tau).
  \end{equation}
\end{lemma}
\begin{proof}
  Recall that the wave operator $\WaveOpHork{2}$ is a symmetric
  operator as an operator on $\realHorkTensor{2}(\Manifold)$, %
  Then, we have that 
  \begin{align*}
    2\Re\bangle*{\WaveOpHork{2}\psi, \widetilde{\MorawetzVF}_1\psi}_{L^2(\Manifold(\tau_1,\tau_2))}
    ={}& \bangle*{\psi, \WaveOpHork{2}^{*}\widetilde{\MorawetzVF}_1\psi}_{L^2(\Manifold(\tau_1,\tau_2))}
         + \bangle*{\widetilde{\MorawetzVF}_1^{*}\WaveOpHork{2}\psi, \psi}_{L^2(\Manifold(\tau_1,\tau_2))}
         + \evalAt*{\JCurrentPert{\widetilde{\MorawetzVF}_1,0}[\psi](\tau)}_{\tau=\tau_1}^{\tau=\tau_2}
    \\
    ={}& \bangle*{\psi, \WaveOpHork{2}\widetilde{\MorawetzVF}_1\psi}_{L^2(\Manifold(\tau_1,\tau_2))}
         - \bangle*{\widetilde{\MorawetzVF}_1\WaveOpHork{2}\psi, \psi}_{L^2(\Manifold(\tau_1,\tau_2))}
         + \evalAt*{\JCurrentPert{\widetilde{\MorawetzVF}_1,0}[\psi](\tau)}_{\tau=\tau_1}^{\tau=\tau_2}
    \\
    ={}& \Re\bangle*{\left[\WaveOpHork{2},\widetilde{\MorawetzVF}_1\right]\psi,\psi}_{L^2(\Manifold(\tau_1,\tau_2))}
         + \evalAt*{\JCurrentPert{\widetilde{\MorawetzVF}_1,0}[\psi](\tau)}_{\tau=\tau_1}^{\tau=\tau_2},
  \end{align*}
  where $\abs*{\JCurrentPert{\widetilde{\MorawetzVF}_1,0}[\psi](\tau)}  \lesssim{} \EnergyHorizonDeg[\psi](\tau)$. Similarly, we also have that
  \begin{align*}
    & 2\Re\bangle*{\WaveOpHork{2}\psi, \widetilde{\MorawetzVF}_0\HorCovDeriv_t  \psi}_{L^2(\Manifold(\tau_1,\tau_2))}\\
    ={}&\bangle*{\WaveOpHork{2}\psi, \widetilde{\MorawetzVF}_0\HorCovDeriv_t \psi}_{L^2(\Manifold(\tau_1,\tau_2))}
         + \bangle*{\widetilde{\MorawetzVF}_0\HorCovDeriv_t \psi, \WaveOpHork{2}\psi}_{L^2(\Manifold(\tau_1,\tau_2))}
         + \evalAt*{\JCurrentPert{\widetilde{\MorawetzVF}_0\HorCovDeriv_t,0}[\psi](\tau)}_{\tau=\tau_1}^{\tau=\tau_2}
    \\
    ={}& - \bangle*{\widetilde{\MorawetzVF}_0\HorCovDeriv_t \WaveOpHork{2}\psi, \psi}_{L^2(\Manifold(\tau_1,\tau_2))}
         + \bangle*{\WaveOpHork{2}\left(\widetilde{\MorawetzVF}_0\HorCovDeriv_t \psi\right), \psi}_{L^2(\Manifold(\tau_1,\tau_2))}
         + \evalAt*{\JCurrentPert{\widetilde{\MorawetzVF}_0\HorCovDeriv_t,0}[\psi](\tau)}_{\tau=\tau_1}^{\tau=\tau_2}
    \\
    ={}& \Re\bangle*{\left[\WaveOpHork{2},\widetilde{\MorawetzVF}_0\HorCovDeriv_t\right]\psi,\psi}_{L^2(\Manifold(\tau_1,\tau_2))}
         + \evalAt*{\JCurrentPert{\widetilde{\MorawetzVF}_0\HorCovDeriv_t,0}[\psi](\tau)}_{\tau=\tau_1}^{\tau=\tau_2},
  \end{align*}
  where $\left|\JCurrentPert{\widetilde{\MorawetzVF}_0\HorCovDeriv_t,0}[\psi](\tau) \right|   \lesssim{} \EnergyHorizonDeg[\psi](\tau)$. Now let us consider the result of multiplying by the Lagrangian correction.
  We have that
  \begin{align*}
    & 2\Re\bangle*{\WaveOpHork{2}\psi, \widetilde{\MorawetzLagrangeCorr}_0\psi}_{L^2(\Manifold(\tau_1,\tau_2))}\\
    ={}& \bangle*{\psi, \WaveOpHork{2}^{*}\widetilde{\MorawetzLagrangeCorr}_0\psi}_{L^2(\Manifold(\tau_1,\tau_2))}
         + \bangle*{\widetilde{\MorawetzLagrangeCorr}_0^{*}\WaveOpHork{2}\psi, \psi}_{L^2(\Manifold(\tau_1,\tau_2))}
         + \evalAt*{\JCurrentPert{0, \MorawetzLagrangeCorr_0}[\psi](\tau)}_{\tau=\tau_1}^{\tau=\tau_2}
    \\
    ={}& \bangle*{\psi, \WaveOpHork{2}\widetilde{\MorawetzLagrangeCorr}_0\psi}_{L^2(\Manifold(\tau_1,\tau_2))}
         + \bangle*{\widetilde{\MorawetzLagrangeCorr}_0\WaveOpHork{2}\psi, \psi}_{L^2(\Manifold(\tau_1,\tau_2))}
         + \evalAt*{\JCurrentPert{0, \MorawetzLagrangeCorr_0}[\psi](\tau)}_{\tau=\tau_1}^{\tau=\tau_2}
    \\
    ={}& \Re\bangle*{\left(\WaveOpHork{2}\widetilde{\MorawetzLagrangeCorr}_0 + \widetilde{\MorawetzLagrangeCorr}_0\WaveOpHork{2}\right)\psi,\psi}_{L^2(\Manifold(\tau_1,\tau_2))}
         + \evalAt*{\JCurrentPert{0, \MorawetzLagrangeCorr_0}[\psi](\tau)}_{\tau=\tau_1}^{\tau=\tau_2},
  \end{align*}
  where $\abs*{\JCurrentPert{0, \MorawetzLagrangeCorr_0}[\psi](\tau)}
  \lesssim{} \EnergyHorizonDeg[\psi](\tau)$. Similarly, we also have that
  \begin{align*}
    &2\Re\bangle*{\WaveOpHork{2}\psi, \left(\widetilde{\MorawetzLagrangeCorr}_{-1}\HorCovDeriv_t + \HorCovDeriv_t\widetilde{\MorawetzLagrangeCorr}_{-1}\right)\psi}_{L^2(\Manifold(\tau_1,\tau_2))}\\
    ={}& \bangle*{\psi, \WaveOpHork{2}^{*}\widetilde{\MorawetzLagrangeCorr}_{-1}\HorCovDeriv_t\psi}_{L^2(\Manifold(\tau_1,\tau_2))}
         + \bangle*{\left(\widetilde{\MorawetzLagrangeCorr}_{-1}\HorCovDeriv_t\right)^{*}\WaveOpHork{2}\psi, \psi}_{L^2(\Manifold(\tau_1,\tau_2))}\\
    ={}& \bangle*{\psi, \WaveOpHork{2}\left(\widetilde{\MorawetzLagrangeCorr}_{-1}\HorCovDeriv_t\psi\right)}_{L^2(\Manifold(\tau_1,\tau_2))}
         + \bangle*{\left(\widetilde{\MorawetzLagrangeCorr}_{-1}\HorCovDeriv_t\right)\WaveOpHork{2}\psi, \psi}_{L^2(\Manifold(\tau_1,\tau_2))}
         + \evalAt*{\JCurrentPert{0, \widetilde{\MorawetzLagrangeCorr}_{-1}\HorCovDeriv_t}[\psi](\tau)}_{\tau=\tau_1}^{\tau=\tau_2}
    \\
    ={}& \Re\bangle*{\left(\WaveOpHork{2}\widetilde{\MorawetzLagrangeCorr}_{-1}\HorCovDeriv_t + \widetilde{\MorawetzLagrangeCorr}_{-1}\HorCovDeriv_t \WaveOpHork{2}\right)\psi,\psi}_{L^2(\Manifold(\tau_1,\tau_2))}
         + \evalAt*{\JCurrentPert{0, \widetilde{\MorawetzLagrangeCorr}_{-1}\HorCovDeriv_t}[\psi](\tau)}_{\tau=\tau_1}^{\tau=\tau_2},
  \end{align*}
  where $\abs*{\JCurrentPert{0, \widetilde{\MorawetzLagrangeCorr}_{-1}\HorCovDeriv_t}[\psi](\tau)}
  \lesssim{} \EnergyHorizonDeg[\psi](\tau)$. We then conclude by using the definitions of
  $\KCurrentPert{\widetilde{\MorawetzVF},\widetilde{\MorawetzLagrangeCorr}}[\psi]$,
  $\JCurrentPert{\widetilde{\MorawetzVF},\widetilde{\MorawetzLagrangeCorr}}[\psi]$
  in \zcref[noname]{eq:KCurrentPert-def} and \zcref[noname]{eq:JCurrentPert-def}. 
\end{proof}

\subsection{Trapped set of \KdS}
\label{sec:trapped-set}

In this subsection, we discuss the well-known properties of
\textit{trapped null geodesics} in \KdS{} which remain in a compact
spatial region for all time (see for example Proposition 3.1 of
\cite{dyatlovAsymptoticsLinearWaves2015} and Subsection 6.4 of
\cite{vasyMicrolocalAnalysisAsymptotically2013}). These null geodesics represent a
fundamental high-frequency geometric obstacle to decay.  To analyze
the dynamics of the trapped set $\TrappedSet_{M,a,\Lambda}$ in frequency space, we
consider null-bicharacteristics. rather than null-geodesics, as
null-geodesics are just the physical projection of the integral curves
of null-bicharacteristics.

It is instructive to first consider the trapped null geodesics in
\SdS, where we can write out the trapped set explicitly, and make some
fundamental observations. Since we are only interested in
slowly-rotating Kerr(-de Sitter), the trapped null geodesics in \KdS{}
can be considered perturbations of the trapped null geodesics in \SdS. We
will make this notion more rigorous in what follows. On \SdS, the trapped set can be located entirely physically.
\begin{lemma}
  \label{lemma:trapping:SdS}
  For $\Metric_{M,0,\Lambda}$ a \SdS{} (or Schwarzschild) background, the
  trapped set is given by
  \begin{equation}
    \label{eq:trapping:SdS}
    \TrappedSet_{M, 0, \Lambda} = \curlyBrace*{(t, r, \omega;\sigma, \xi, \eta): r=3M, \xi=0, \PrinSymb_{M,0,\Lambda}=0},
  \end{equation}
  where $\PrinSymb_{M, 0, \Lambda}$ is the principal symbol of the scalar wave
  operator $\ScalarWaveOp[\Metric_{M, 0, \Lambda}]$. Moreover, the trapped set is
  unstable in the sense that
  \begin{equation}
    \label{eq:trapping:SdS:hyperbolic-trapping-condition}
    \Delta_{M, 0, \Lambda} > 0,\quad
    \PrinSymb_{M, 0, \Lambda} =0,\quad
    \pm (r-3M) > 0, \quad
    H_{\PrinSymb_{M, 0, \Lambda}}r = 0 \implies \pm H_{\PrinSymb_{M, 0, \Lambda}}^2r > 0,
  \end{equation}
  where  we denote by $H_{\PrinSymb_{M, 0, \Lambda}}$ the
  Hamiltonian vectorfield of $\Metric_{M, 0, \Lambda}$.
\end{lemma}
\begin{proof}
  See \cite{fangLinearStabilitySlowlyRotating2026}.
\end{proof}

\begin{remark}
  The physical projection of $\TrappedSet_{M, 0, \Lambda}$ is
  exactly the photon sphere, $r=3M$. 
\end{remark}

We now move onto the trapped set in the case of \KdS. In this case,
the trapped set exhibits frequency-dependent behavior.
\begin{lemma}
  \label{lemma:trapping:KdS}
  For $\Metric_{M, a, \Lambda}$ a \KdS{} background, the trapped set
  \begin{equation}
    \label{eq:trapping:KdS:trapped-set}
    \TrappedSet_{M, a, \Lambda} = \curlyBrace*{(t, r, \omega; \sigma, \xi, \FreqAngular): r=\rTrapping_{M,a,\Lambda}(\sigma, \FreqAngular), \xi=0, \PrinSymb_{M,a,\Lambda}=0},
  \end{equation}
  where $\PrinSymb_{M, a, \Lambda}$ is the principal symbol of the scalar wave
  operator $\ScalarWaveOp[\Metric_{M, a, \Lambda}]$, and $\rTrapping_{M, a, \Lambda}(\sigma,
  \FreqAngular)$ is a function satisfying the following properties. 
  \begin{enumerate}
  \item $\rTrapping_{M, a, \Lambda}(\sigma,\FreqAngular)$ lies in an $O(a)$ neighborhood
    of $r=3M$ for all $\sigma, \FreqAngular$.
  \item $\rTrapping_{M, a, \Lambda}(\sigma, \FreqAngular)$ is smooth in $\sigma,
    \FreqAngular$, as well as the black hole parameters $(M, a, \Lambda)$.  
  \end{enumerate}
  Moreover, the trapped set is unstable in the sense that
  \begin{equation}
    \label{eq:trapping:KdS:hyperbolic-trapping-condition}
    \begin{gathered}
      \Delta_{M, a, \Lambda}>0,
    \quad \PrinSymb_{M, a, \Lambda}=0,
    \quad \pm (r-\rTrapping_{M, a, \Lambda})>0,\quad
    H_{\PrinSymb_{M, a, \Lambda}} r=0 \implies \pm H^2_{\PrinSymb_{M, a, \Lambda}} r > 0,
    \end{gathered}    
  \end{equation}
  where we denote by $H_{\PrinSymb_{M, a, \Lambda}}$ the Hamiltonian vectorfield of
  $\Metric_{M, a, \Lambda}$. 
\end{lemma}

\begin{proof}
  See \cite{fangLinearStabilitySlowlyRotating2026}.
\end{proof}

\subsection{Microlocal Morawetz norms}

The main technicality with the Morawetz bulk quantity will be the
introduction of microlocally weighted norms in order to capture the
frequency-dependent behavior of the trapped set in Kerr-de Sitter. Since we are working in the slowly-rotating regime, the trapped set is
localized to a $O(a)$ neighborhood of $r=3M$, where there are no geometric features
apart from trapping, in particular, no superradiance. As such, it will
be convenient to just work with the Boyer-Lindquist coordinates on
phase space $(t,r,\omega;\sigma,\xi,\eta)$. We remark that while we
will write things out explicitly in
$(\theta,\phi;\eta_{\theta},\eta_{\phi})$ coordinates on the sphere,
all the norms extend to cover the entire sphere.

To this end, let us write out the principal symbol of the scalar wave on \KdS{} as %
\begin{equation}
  \label{eq:trapping-norm:sigmai-def}
  \PrinSymb_{M,a,\Lambda} = \Metric^{-1}_{M,a,\Lambda}(dt, dt)\left(\sigma-\sigma_1\right)\left(\sigma-\sigma_2\right),
\end{equation}
where $\sigma_i = \sigma_i(r,\phi;\xi,\eta)$ are distinct smooth
1-homogeneous real symbols with respect to the spatial frequency
variables. Define then
\begin{equation}
  \label{eq:trapping-norm:elli-def}
  \ell_i(r,\varphi;\xi,\eta)
  = r - r_i(\sigma_i,\eta),
\end{equation}
where
\begin{equation}
  \label{eq:trapping-norm:ri-def}
  r_i(\sigma_i,\eta) = \rTrapping_{M,a,\Lambda}(\sigma_i, \eta),
\end{equation}
and where $\sigma_i$ are as defined in
(\ref{eq:trapping-norm:sigmai-def}), and $\rTrapping_{M,a,\Lambda}$ is
as constructed in \zcref[cap]{lemma:trapping:KdS}. We can now introduce the frequency-dependent Morawetz norms we will use
\begin{equation}
  \label{eq:MorNormTrap:def}
  \begin{split}
    \MorNormTrap[\psi](\tau_1,\tau_2)
    \vcentcolon={}& \int_{\Manifold_{\operatorname{trap}}(\tau_1,\tau_2)}\frac{1}{r^2}\sum_{i\neq j}\abs*{\mathring{\chi}\left(1-\frac{r_i(\sigma_i,\eta)}{r}\right)\left(D_t-\sigma_j(D,x)\right)(\mathring{\chi}\psi)}^2 + \MorNorm[\psi](\tau_1,\tau_2)
  .
  \end{split}  
\end{equation}


\section{Morawetz Estimates}
\label{sec:Morawetz}

In this section, we present the proof for the Morawetz estimate for
the model problem in \zcref[noname]{eq:model-problem-gRW}.
We state the main result of this section.

\begin{proposition}[Morawetz estimate]
  \label{prop:Morawetz:KdS:main}
  Let $\psi\in\realHorkTensor{2}$ be a  solution to the model problem in
  \zcref[noname]{eq:model-problem-gRW}.  Then the following estimates hold in the region
  $\Manifold(\tau_1,\tau_2)$, $1\le \tau_1<\tau_2\le \tau_{*}$.
  \begin{equation}
    \label{eq:Morawetz:KdS:main}
    \begin{split}
      \MorNorm[\psi](\tau_1,\tau_2)
      \lesssim{}& \sup_{\tau\in[\tau_1,\tau_2]}\EnergyHorizonDeg[\psi]
                  + \int_{\mathcal{A}(\tau_1,\tau_2)}\abs*{\nabla_{4}\psi}^2
                  + \int_{\SigmaStar(\tau_1,\tau_2)}\abs*{\nabla_{3}\psi}^2\\ 
                &+ \delta_{\Horizon} \pth{\SpacelikeFlux_{\mathcal{A}}[\psi](\tau_1,\tau_2)
                  + \SpacelikeFlux_{\SigmaStar}[\psi](\tau_1,\tau_2)}
                +  \int_{\Manifold(\tau_1,\tau_2)}\left(\abs*{\nabla_{\HprVF}\psi}+r^{-1}\abs*{\psi} + |N|\right)\abs*{N}.
    \end{split}    
  \end{equation}
\end{proposition}

The key to proving \zcref[cap]{prop:Morawetz:KdS:main} is to generate
a positive bulk term.  This is proven in three main steps. We provide a brief
overview below.

\paragraph{Step 1: The non-trapping region.}

We first choose physical space multipliers
$\MorawetzVF_{M,a,\Lambda}$, $\MorawetzLagrangeCorr_{M,a,\Lambda}$,
  and $\MorawetzOneForm_{M,a,\Lambda}$ such that for any sphere
$S(\tau,r)\subset \Manifold_{\cancel{\operatorname{trap}}}$,
\begin{equation}
  \label{eq:Morawetz:non-trapping:rough}
  \begin{split}
    &\int_{S(\tau, r)}\KCurrent{\MorawetzVF_{M,a,\Lambda}, \MorawetzLagrangeCorr_{M,a,\Lambda}, \MorawetzOneForm_{M,a,\Lambda}}[\psi]\\
    \gtrsim{} & \left(1+O(a)+O(\sqrt{\Lambda})\right)\int_{S(\tau, r)}\pth{
           \frac{M}{r^2}\abs*{\HprVF\psi}^2
           + \frac{M}{r^2}\abs*{\HawkingVF \psi}^2
           + r^{-1}\abs*{\nabla \psi}^2      
           + r^{-3}\abs*{\psi}^2}.
  \end{split} 
\end{equation}
The estimate \zcref[noname]{eq:Morawetz:non-trapping:rough} is proved in
\zcref[cap]{prop:Morawetz:outside-trapping:KCurrent-coercive}. Its proof of
   is perturbative in
  nature. In fact, up to terms which vanish at the event and
  cosmological horizons (such as the exact form of
  $\Delta_{M,a,\Lambda}$, and other terms depending on
  $\Delta_{M,a,\Lambda}$), the proof of Morawetz is in fact
  perturbative from Schwarzschild, since we are in the
  slowly-rotating, $\Lambda$-small regime.

\paragraph{Step 2: The trapping region.}

We next consider the bulk terms in the trapping region. Unlike the previous step, the analysis done in this step
is not perturbative off of Schwarzschild but rather perturbative off of Kerr as the main
geometric structure used is that of normally hyperbolic trapping,
which holds uniformly for slowly-rotating Kerr(-de Sitter).

To deal with the frequency-dependent nature of trapping, the main goal
will be to find $\MorawetzVF'$ and $\MorawetzLagrangeCorr'$
pseudo-differential $O(a)$-modifiers of the physical space multipliers
$\MorawetzVF_{M,a,\Lambda}$, $\MorawetzLagrangeCorr_{M,a,\Lambda}$
chosen in the previous step so as to guarantee that for $\psi$
supported near $r=3M$ we have a positive bulk term when applying the
divergence theorem in \zcref[cap]{lemma:PsiDO-divergence-theorem}. To
this end, we fix symbols
$\MorawetzSym\in \MixedSymClass{1}{1}(\Manifold,
\realHorkTensor{2}),
\MorawetzLagrangeCorrSym\in\MixedSymClass{0}{1}(\Manifold,
\realHorkTensor{2})$ such that
\begin{equation*}
  \abs*{q}^2\left(\frac{1}{2\ImagUnit}H_p\MorawetzSym + p\MorawetzLagrangeCorrSym\right) > 0,
\end{equation*}
which will ensure that the top order terms in the final Morawetz
estimate are non-negative. This is proven in
\zcref[cap]{lemma:ILED-KdS:Bulk}. Afterwards, we show that various
error terms are controllable in
\zcref[cap]{lemma:ILED-near:LoT-control:zero-order-mixed-term,lemma:ILED-near:LoT-control:exchange-degeneracy-trick,lemma:ILED-near:KdS:extra-terms}.

\paragraph{Step 3: Combining the bulk terms}

We next combine the pointwise positive bulk term in the non-trapping
region proved in
\zcref[cap]{prop:Morawetz:outside-trapping:KCurrent-coercive} with the
symbolic positivity in \zcref[cap]{lemma:ILED-KdS:Bulk} to prove that 
\begin{equation}
  \label{eq:Morawetz:cutoff-PsiDO:bulk:rough}
  \begin{split}
    & a \KCurrentPert{\widetilde{\MorawetzVF}, \widetilde{\MorawetzLagrangeCorr}}_{\Main}[\psi](\tau_1,\tau_2)
      + \int_{\Manifold(\tau_1,\tau_2)}\KCurrent{\MorawetzVF_{M,a,\Lambda}, \MorawetzLagrangeCorr_{M,a,\Lambda}, \MorawetzOneForm_{M,a,\Lambda}}[\psi]\\
    \gtrsim{}&  \MorNormTrap[\psi](\tau_1,\tau_2)
               - O(a)\left(
               \norm*{\psi}_{L^2_{\operatorname{comp}}(\Manifold(\tau_1,\tau_2))}^2
               + \int_{\Manifold(\tau_1,\tau_2)}\abs*{N}^2
               + \sup_{[\tau_1,\tau_2]}\EnergyHorizonDeg[\psi](\tau)
               \right),               
  \end{split}    
\end{equation}
where
$\KCurrentPert{\widetilde{\MorawetzVF},
  \widetilde{\MorawetzLagrangeCorr}}_{\Main}[\psi](\tau_1,\tau_2)$,
defined in \zcref[noname]{eq:KCurrent-princ:def} represents the
principal terms in the divergence theorem
\zcref[cap]{lemma:PsiDO-divergence-theorem} applied with appropriately
cutoff quantizations of $\MorawetzSym$ and $\MorawetzLagrangeCorrSym$.
The estimate in \zcref[noname]{eq:Morawetz:cutoff-PsiDO:bulk:rough} is
proven in \zcref[cap]{lemma:Morawetz:cutoff-PsiDO:bulk}, and uses the
fact that since the pseudodifferential perturbations are $O(a)$, the
cutoff errors generated, for $a$ sufficiently small, can be controlled
by the pointwise positivity proven in
\zcref[cap]{prop:Morawetz:outside-trapping:KCurrent-coercive}. 

We emphasize that this approach maintains the heuristic philosophy in
\cite{tataruLocalEnergyEstimate2011} that the only geometric behavior
on Kerr(-de Sitter) that is frequency-based is the trapped null geodesics.

\subsection{Morawetz estimates: non-trapping}
\label{sec:Morawetz:SdS}

In this section, we first consider the region of our spacetime outside
of a neighborhood of the trapped set. Since we are only considering
the slowly-rotating Kerr-de Sitter family, the desired estimate on
this region can be proven by considering just a single
(non-pseudodifferential) vectorfield multiplier. 
We refer the reader to Section 10.1 of
\cite{klainermanGlobalNonlinearStability2020} for a proof of a similar
estimate to solutions of equations of the form
\zcref[noname]{eq:model-problem-gRW} on Schwarzschild, and to Section 7.2 of
\cite{giorgiWaveEquationsEstimates2024} for a discussion of similar
conditional Morawetz estimates on Kerr.

The main difficulty of proving a Morawetz estimate is to adapt the
energy estimates to the degeneracy of the trapped set of null
bi-characteristics. The main result of this section is as follows.
\begin{proposition}
  \label{prop:Morawetz:outside-trapping:KCurrent-coercive}
  Let $\Metric_{M,a,\Lambda}$ be a subextremal, slowly-rotating \KdS{}
  metric. Then, there exists a choice of
  $\left(\MorawetzVF_{M,a,\Lambda},
    \MorawetzLagrangeCorr_{M,a,\Lambda},
    \MorawetzOneForm_{M,a,\Lambda}\right)$ such that
  \begin{enumerate}
  \item $\MorawetzVF_{M, a, \Lambda} = \mathcal{F}_{M,a,\Lambda}\partial_r
    = f_{M,a,\Lambda}(r)\HprVF$ is a smooth vectorfield,
  \item $\MorawetzLagrangeCorr_{M, a, \Lambda}=\MorawetzLagrangeCorr_{M,
      a, \Lambda}(r)$ is a smooth function such that
    \begin{equation} 
      \label{eq:ILED:SdS:q:falloff}
      \abs*{\MorawetzLagrangeCorr_{M, a, \Lambda}}\lesssim \abs*{\Delta} r^{-3},\qquad
      \abs*{\partial_r\MorawetzLagrangeCorr_{M, a, \Lambda}}\lesssim r^{-2},
    \end{equation}
  \item $\MorawetzOneForm_{M, a, \Lambda}=\MorawetzOneForm_{M, a, \Lambda}(r)$ is a smooth, spherically symmetric one-form,
  \end{enumerate}
  with the property that for any sphere $S(\tau,r)\subset \Manifold_{\cancel{\operatorname{trap}}}$,
  \begin{equation}
    \label{eq:Morawetz:outside-trapping:KCurent-coercive}
    \begin{split}
      &\int_{S(\tau, r)}\KCurrent{\MorawetzVF_{M,a,\Lambda}, \MorawetzLagrangeCorr_{M,a,\Lambda}, \MorawetzOneForm_{M,a,\Lambda}}[\psi]
    \gtrsim{} \int_{S(\tau, r)}
    \left(
      \frac{M}{r^2}\abs*{\nabla_{\HprVF}\psi}^2
      + \frac{(r-3M)^2}{r^2}\left(
        \frac{M}{r^2}\abs*{\nabla_{\HawkingVF} \psi}^2
        + r^{-1}\abs*{\nabla \psi}^2
      \right)
      + r^{-3}\abs*{\psi}^2
    \right).
    \end{split}    
  \end{equation}
\end{proposition}

To begin, we will choose $\left( \MorawetzVF, \MorawetzLagrangeCorr \right)$ so that
the principal terms in
$\KCurrent{\MorawetzVF, \MorawetzLagrangeCorr, 0}[\psi]$ have the
desired coercivity. To this end, we first consider
\begin{equation*}
  \begin{split}
    \KCurrent{\MorawetzVF, \MorawetzLagrangeCorr, \MorawetzOneForm}_0[\psi]
    \vcentcolon={}& \EMTensor[\psi]\cdot\DeformationTensor{\MorawetzVF}
          - \frac{1}{2}\MorawetzVF(V)\abs*{\psi}^2
          + \MorawetzLagrangeCorr\mathcal{L}[\psi]
          - \frac{1}{2}\abs*{\psi}^2\ScalarWaveOp[\Metric]\MorawetzLagrangeCorr
          + \frac{1}{2}\Divergence\left(\abs*{\psi}^2\MorawetzOneForm\right).
  \end{split}    
\end{equation*}

\begin{lemma}
  \label{lemma:ILED:SdS:Bulk:principal-terms}
  There exists a choice of
  $(\MorawetzVF,\MorawetzLagrangeCorr) =
  (\MorawetzVF_{M,a,\Lambda},\MorawetzLagrangeCorr_*)$
  satisfying the conditions in \zcref[cap]{prop:Morawetz:outside-trapping:KCurrent-coercive} such that
  \begin{equation*}
    \begin{split}
      \abs*{q}^2\KCurrent{\MorawetzVF_{M, a, \Lambda},\MorawetzLagrangeCorr_*,0}_0[\psi]
      ={}& \KCurrentR\abs*{\Delta\nabla_r\psi}^2
           + \KCurrentAngular{\MorawetzVF_{M, a, \Lambda},\MorawetzLagrangeCorr_*,0}[\psi]
           + \KCurrentLTwo\abs*{\psi}^2
           + \mathfrak{E}^{\MorawetzVF_{M,a,\Lambda}, \MorawetzLagrangeCorr_{M,a,\Lambda},0}[\psi].
    \end{split}      
  \end{equation*}
  where 
  \begin{equation}
    \label{eq:Morawetz:SdS:bulk-coefficient-conditions}
    \begin{split}
      \KCurrentAngular{\MorawetzVF_{M, a, \Lambda},\MorawetzLagrangeCorr_*,0}[\psi]
      &\gtrsim \frac{1}{r}\frac{(r-3M)^2}{r^2}\CartarOp^{\alpha\beta}\HorCovDeriv_{\alpha}\psi\cdot\HorCovDeriv_{\beta}\psi,\\
      \KCurrentR
      &\gtrsim \frac{1}{r^2(r^2+a^2)},\\
      \KCurrentLTwo
      &= \KCurrentLTwo_0 + \KCurrentLTwo_1,\\ 
      \mathfrak{E}^{\MorawetzVF_{M,a,\Lambda}, \MorawetzLagrangeCorr_{M,a,\Lambda},0}[\psi]
      &\lesssim O(ar^{-3})\left(\abs*{\HawkingVF\psi}^2 + \abs*{\nabla\psi}^2\right),
    \end{split}      
  \end{equation}
  with
  \begin{align*}
    2\KCurrentLTwo_0&= \frac{M\left(54M^2-46Mr+9r^2\right)}{r^4}
                      + \frac{5M\Lambda}{3} - \frac{4r^3\Lambda^2}{9} + O(a^2r^{-3}),\\
    \KCurrentLTwo_1&=\frac{4(r-3M)\Upsilon\left(r-4M + r^3\frac{\Lambda}{3}\right)}{r^3}
                     + \frac{2(r-3M)(r^3\Lambda - 3M)}{3r}\Lambda + O(a^2r^{-3}).
  \end{align*}
\end{lemma}

\begin{proof}  
  We first choose 
  \begin{equation}
        \label{eq:ILED:SdS:X-def}
    \MorawetzVF_{M,a,\Lambda} = f_{M, a, \Lambda}(r)\HprVF = \mathcal{F}_{M,a,\Lambda}(r)\partial_r,
  \end{equation}
  where $f_{M, a, \Lambda}(r)$ is bounded near $r=3M$. In
  particular, we will choose
  \begin{equation}
    \label{eq:ILED:SdS:X-def:f-def}
    f_{M, a, \Lambda}(r)= \frac{r-3M}{r}.
  \end{equation}
  This is equivalent to
  $\mathcal{F}_{M,a,\Lambda}(r) =
  \frac{r-3M}{r}\frac{\Delta}{r^2+a^2}$.
  We can first
  calculate that for
  $\MorawetzVF_{M,a,\Lambda} = \mathcal{F}(r)\partial_r$ as defined
  above in \zcref[noname]{eq:ILED:SdS:X-def:f-def},
  \begin{align*}   
    2\abs*{q}^2\DeformationTensor[]{\MorawetzVF_{M,a,\Lambda}}^{\alpha\beta}
    ={}& \left(2\Delta\partial_r\mathcal{F}_{M,a,\Lambda} - \mathcal{F}_{M,a,\Lambda}\partial_r\Delta\right)\partial_r^{\alpha}\partial_r^{\beta}
         - \mathcal{F}\partial_r\left(\frac{1}{\Delta}\mathcal{R}^{\alpha\beta}\right)
         + \MorawetzVF_{M,a,\Lambda}(\abs*{q}^2)\Metric^{\alpha\beta}.
  \end{align*}
  Thus, we have that
  \begin{equation}
    \label{eq:ILED:EMT-X}
    \begin{split}
      2\abs*{q}^2\EMTensor[\psi] \cdot \DeformationTensor[]{\MorawetzVF_{M,a,\Lambda}}
      ={}& \left(2\Delta\partial_r\mathcal{F}_{M,a,\Lambda} - \mathcal{F}_{M, a, \Lambda}\partial_r\Delta\right)\abs*{\partial_r\psi}^2
           - \mathcal{F}_{M,a,\Lambda}\partial_r\left(\frac{1}{\Delta}\mathcal{R}^{\alpha\beta}\right)\HorCovDeriv_{\alpha}\psi \cdot\HorCovDeriv_{\beta}\psi\\
          & + \MorawetzVF(\abs*{q}^2)\left(\Lagrangian[\psi] - V\abs*{\psi}^2\right)
           - \Lagrangian[\psi]\abs*{q}^2\nabla_{\Metric}\cdot \MorawetzVF.
    \end{split}  
  \end{equation}
  We will fix the Lagrangian correction in two steps. We first use a
  Lagrangian correction to remove the last $\Lagrangian[\psi]$ term that
  appears in \zcref[noname]{eq:ILED:EMT-X}.
  Consider 
  \begin{equation}
    \label{eq:ILED:q0:def}
    \MorawetzLagrangeCorr_0
    = - \frac{1}{2}\abs*{q}^2\Divergence\left(\abs*{q}^{-2}\MorawetzVF_{M,a,\Lambda}\right).
  \end{equation}
  Using from \zcref{eq:model-wave-equation:Lagrangian} the fact that $\abs*{q}^2\Lagrangian[\psi]
    = \Delta\abs*{\nabla_r\psi}^2
    + \frac{1}{\Delta}\mathcal{R}^{\alpha\beta}\HorCovDeriv_{\alpha}\psi\cdot\HorCovDeriv_{\beta}\psi
    + \abs*{q}^2V\abs*{\psi}^2$ we see that
  \begin{equation*}
    \begin{split}
          \abs*{q}^2\KCurrent{\MorawetzVF_{M,a,\Lambda},\MorawetzLagrangeCorr_0,0}_0[\psi]
    ={}& \left(\Delta\partial_r\mathcal{F}_{M,a,\Lambda} - \frac{1}{2}\mathcal{F}_{M, a, \Lambda}\partial_r\Delta\right)\abs*{\partial_r\psi}^2
           - \frac{1}{2}\mathcal{F}_{M, a, \Lambda}\partial_r\left(\frac{1}{\Delta}\mathcal{R}^{\alpha\beta}\right)\HorCovDeriv_{\alpha}\psi \cdot\HorCovDeriv_{\beta}\psi\\
          & - \frac{1}{2}\left( \MorawetzVF_{M, a, \Lambda}(\abs*{q}^2)V + \abs*{q}^2\MorawetzVF_{M, a, \Lambda}(V) + \abs*{q}^2\ScalarWaveOp[\Metric]\MorawetzLagrangeCorr_0\right)\abs*{\psi}^2.
    \end{split}
  \end{equation*}
  We will now use an additional Lagrangian corrector to remove some of
  components of the $\mathcal{R}$ term. To this end, consider
  \begin{equation}
    \label{eq:Morawetz:w1:def}
    \MorawetzLagrangeCorr_1
    = \frac{(r-3M)(3r^2(r-3M)+a^4r\Lambda + a^2(3M+3r+r^3\Lambda)}{3r(r^2+a^2)^2},
  \end{equation}
  so that
  \begin{equation*}    
    \left( \mathcal{F}_{M, a, \Lambda}\partial_r\left(\frac{1}{\Delta}\mathcal{R}^{\alpha\beta}\right) - 2\MorawetzLagrangeCorr_1\frac{\mathcal{R}^{\alpha\beta}}{\Delta} \right)\HawkingVF_{\alpha}\HawkingVF_{\beta} = 0.
  \end{equation*}
  We then have that 
  \begin{equation*}
    \begin{split}
      \abs*{q}^2\KCurrent{\MorawetzVF_{M,a,\Lambda}, \MorawetzLagrangeCorr_0 + \MorawetzLagrangeCorr_1,0}_0[\psi]
      ={}& \left(\Delta\partial_r\mathcal{F}_{M,a,\Lambda} - \frac{1}{2}\mathcal{F}_{M, a, \Lambda}\partial_r\Delta + \Delta\MorawetzLagrangeCorr_1\right)\abs*{\partial_r\psi}^2\\
         & - \frac{1}{2}\left(\mathcal{F}_{M, a, \Lambda}\partial_r\left(\frac{1}{\Delta}\mathcal{R}^{\alpha\beta}\right) - 2\MorawetzLagrangeCorr_1\frac{\mathcal{R}^{\alpha\beta}}{\Delta}\right)\HorCovDeriv_{\alpha}\psi \cdot\HorCovDeriv_{\beta}\psi\\
          & - \frac{1}{2}\left( \MorawetzVF_{M, a, \Lambda}(\abs*{q}^2)V + \abs*{q}^2\MorawetzVF_{M, a, \Lambda}(V) + \abs*{q}^2\ScalarWaveOp[\Metric]\left(\MorawetzLagrangeCorr_0+\MorawetzLagrangeCorr_1\right) - 2\abs*{q}^2\MorawetzLagrangeCorr_1V\right)\abs*{\psi}^2.
    \end{split}
  \end{equation*}
  We define
  \begin{equation*}
    \begin{split}
      \mathcal{A} ={}& \Delta\partial_r\mathcal{F}_{M,a,\Lambda} - \frac{1}{2}\mathcal{F}_{M, a, \Lambda}\partial_r\Delta + \Delta\MorawetzLagrangeCorr_1,\\
      \mathcal{U}^{\alpha\beta}={}& - \frac{1}{2}\left(\mathcal{F}_{M, a, \Lambda}\partial_r\left(\frac{1}{\Delta}\mathcal{R}^{\alpha\beta}\right) - 2\MorawetzLagrangeCorr_1\frac{\mathcal{R}^{\alpha\beta}}{\Delta}\right),\\
    2\mathcal{P}_0 ={}& -\abs*{q}^2\ScalarWaveOp[\Metric](\MorawetzLagrangeCorr_0+\MorawetzLagrangeCorr_1),\\
    2\mathcal{P}_1 ={}& -\left(\MorawetzVF_{M,a,\Lambda}(\abs*{q}^2)V
                       + \abs*{q}^2\MorawetzVF_{M,a,\Lambda}(V)
                       - 2\abs*{q}^2\MorawetzLagrangeCorr_1V
                       \right).
    \end{split}
  \end{equation*}
  Then, a (lengthy) computation shows that
  \begin{equation}
    \label{eq:Morawetz:K-prelim-computations}
    \begin{split}
      \mathcal{A} ={}& \frac{3M\Delta^2}{r^2(r^2+a^2)}                       
                       ,\\
      \mathcal{U}^{\alpha\beta}={}& \frac{(r-3M)^2}{r^3}\mathcal{O}^{\alpha\beta}
                                    + O(ar^{-3})\left(
                                    \HawkingVF^{\alpha}\KillPhi^{\beta}
                                    \HawkingVF^{\beta}\KillPhi^{\alpha}
                                    \right)
                                    ,\\
      2\mathcal{P}_0 ={}& \frac{M(54M^2-46 Mr + 9r^2)}{r^4} + \frac{5M\Lambda}{3} - \frac{4r^3\Lambda^2}{9} + O(a^2r^{-3}),\\
    \mathcal{P}_1={}& \frac{4(r-2M)(4-3M)(r-4M)}{r^4} + \frac{2M(r-3M)}{3r}\Lambda + \frac{2}{9}r^2(r-3M)\Lambda^2 + O(a^2r^{-3}).
    \end{split}
  \end{equation}

  Defining
  \begin{equation}
    \label{eq:Morawetz:w-star-explicit}
    \MorawetzLagrangeCorr_{*}
    = \MorawetzLagrangeCorr_0 +
    \MorawetzLagrangeCorr_1
    = \frac{3a^2M + r^2(2r-3M)}{2r^2(r^2+a^2)^2}\Delta
    ,
  \end{equation}
  and
  $\KCurrentAngular{\MorawetzVF_{M,a,\Lambda},
    \MorawetzLagrangeCorr_{*}, 0}[\psi]$ and
  $\mathfrak{E}^{\MorawetzVF_{M,a,\Lambda},\MorawetzLagrangeCorr_{*},0}[\psi]$
  such that
  \begin{align*}
    \mathcal{U}^{\alpha\beta}\HorCovDeriv_{\alpha}\psi\cdot\HorCovDeriv_{\beta}\psi
    & = \KCurrentAngular{\MorawetzVF_{M,a,\Lambda}, \MorawetzLagrangeCorr_{*}, 0}[\psi]
    + \mathfrak{E}^{\MorawetzVF_{M,a,\Lambda},\MorawetzLagrangeCorr_{*},0}[\psi],\\
    \KCurrentAngular{\MorawetzVF_{M,a,\Lambda}, \MorawetzLagrangeCorr_{*}, 0}[\psi]
    &= \frac{(r-3M)^2}{r^3}
      \CartarOp^{\alpha\beta}\HorCovDeriv_{\alpha}\psi\cdot\HorCovDeriv_{\beta}\psi,
  \end{align*}
  concludes the proof of
  \zcref[cap]{lemma:ILED:SdS:Bulk:principal-terms}.
\end{proof}

There are a few remaining issues to handle. While
$\KCurrent{\MorawetzVF_{M,a,\Lambda},
  \MorawetzLagrangeCorr_*,0}_0[\psi]$ is principally positive (up to
$O(a)$ terms), it does not control all derivatives. Notably, it lacks
control over the $\HawkingVF$ derivative. Second, we would like to
have a bulk term that is positive at both the principal level and at
the level of the $L^2$ terms.  We also have to show how to control the
mixed terms in
$\KCurrent{\MorawetzVF,\MorawetzLagrangeCorr,\MorawetzOneForm}[\psi]-\KCurrent{\MorawetzVF,\MorawetzLagrangeCorr,\MorawetzOneForm}_0[\psi]$.

We first address the issue of the lower-order positivity. To do this,
we will choose the appropriate $\MorawetzOneForm_{M,a,\Lambda}$
term. To this end, observe that
\begin{equation*}
  \KCurrent{\MorawetzVF,\MorawetzLagrangeCorr,\MorawetzOneForm}[\psi]
  = \KCurrent{\MorawetzVF,\MorawetzLagrangeCorr,0}[\psi]
  + \frac{1}{2}\Divergence_{M,a,\Lambda}\left(\abs*{\psi}^2\MorawetzOneForm\right).
\end{equation*}
If we pick $\MorawetzOneForm = v(r)\partial_r$ as stipulated in 
\zcref[cap]{prop:Morawetz:outside-trapping:KCurrent-coercive}, then we can compute
\begin{equation*}
  \frac{1}{2}\Divergence_{M,a,\Lambda}\left(\abs*{\psi}^2\MorawetzOneForm\right)
  = v(r)\psi\cdot\nabla_r\psi
  + \frac{1}{2}\left(\partial_rv + \frac{2r}{\abs*{q}^2}v\right)\abs*{\psi}^2.
\end{equation*}
As a result, using \zcref[cap]{lemma:ILED:SdS:Bulk:principal-terms}, we
have that for $m=v(r)\partial_r$, and
$(\MorawetzVF,\MorawetzLagrangeCorr) =
\left(\MorawetzVF_{M,a,\Lambda},\MorawetzLagrangeCorr_*\right)$
as chosen in \zcref[cap]{lemma:ILED:SdS:Bulk:principal-terms},
\begin{equation}
  \label{eq:ILED:SdS:KCurrent-decomp-into-Hardy}
  \begin{split}
    \abs*{q}^2\KCurrent{\MorawetzVF_{M,a,\Lambda},\MorawetzLagrangeCorr_{*},\MorawetzOneForm}_0[\psi]
  ={}& \delta\mathcal{A}\abs*{\partial_r\psi}^2
       + \delta\KCurrentAngular{\MorawetzVF_{M,a,\Lambda},\MorawetzLagrangeCorr_{*},0}[\psi]
       + \KCurrent{\MorawetzVF_{M,a,\Lambda}, \MorawetzLagrangeCorr_{*}, 0}_{\delta, v}[\psi]\\
     & + (1-\delta)\left(\KCurrentAngular{\MorawetzVF_{M,a,\lambda},\MorawetzLagrangeCorr_{*},0}\abs*{\nabla\psi}^2
       - \frac{(r-3M)^2}{r^3}\abs*{\psi}^2
       \right) + \mathfrak{E}^{\MorawetzVF_{M,a,\Lambda}, \MorawetzLagrangeCorr_{*},0}[\psi]
    ,
  \end{split}  
\end{equation}
where $\mathfrak{E}^{\MorawetzVF_{M,a,\Lambda}, \MorawetzLagrangeCorr_{*},0}[\psi]=\mathfrak{E}^{\MorawetzVF_{M,a,\Lambda}, \MorawetzLagrangeCorr_{M,a,\Lambda},0}[\psi]$ and 
\begin{equation}
  \label{eq:ILED:SdS:KCurrent-with-Hardy:def}
  \begin{split}
    \KCurrent{\MorawetzVF_{M,a,\Lambda}, \MorawetzLagrangeCorr_{*}, 0}_{\delta, v}[\psi]
    ={}& (1-\delta)\KCurrentR\abs*{\partial_r\psi}^2
         + \left(\KCurrentLTwo + \left(1-\delta\right)
         \frac{(r-3M)^2}{r^3}
         \right)\abs*{\psi}^2\\
       & + \frac{\abs*{q}^2}{2}\left(
         2 v(r)\psi\cdot \partial_r\psi
         + \left(
         \partial_rv + \frac{2}{r}v
         \right)\abs*{\psi}^2
         \right).
  \end{split}  
\end{equation}
Observe that the first two terms on the \RHS{} of
\zcref[noname]{eq:ILED:SdS:KCurrent-decomp-into-Hardy} integrate to something
positively directly using \zcref[cap]{lemma:ILED:SdS:Bulk:principal-terms}.  We first argue that the
second line in the right-hand side of
\zcref[noname]{eq:ILED:SdS:KCurrent-decomp-into-Hardy} integrates to a
positive term.  To this end, first recall the following Poincare
estimate.
\begin{lemma}[Poincare estimate]
  \label{lemma:Poincare}
  Let $S$ be a sphere of radius $r$. Then for
  $\psi\in \realHorkTensor{2}$, we have
  \begin{equation*}
    \int_S\abs*{\nabla \psi}^2 \ge \frac{2}{r^2}\int_S\abs*{\psi}^2
    -O(a)\int_S\left(\abs*{\nabla\psi}^2 + r^{-2}\abs*{\nabla_{\KillT}\psi}^2 + r^{-4}\abs*{\psi}^2\right).
  \end{equation*}
\end{lemma}
\begin{proof}
  The proof given of Lemma 7.2.3 of \cite{giorgiWaveEquationsEstimates2024} extends to
  our setting as the proof is independent of the presence of a
  cosmological constant. 
\end{proof}
Thus we see that to show
\zcref{prop:Morawetz:outside-trapping:KCurrent-coercive}, it suffices to show that
$\KCurrent{\MorawetzVF, \MorawetzLagrangeCorr, 0}_{\delta, v}[\psi]$
is positive. To this end, we have the following application of the
Hardy inequality.
\begin{lemma}
  \label{lemma:ILED:SDS:KCurrent-with-Hardy}
  There exists a function $v(r)$ such that
  $v(r) = O\left(\Delta r^{-\frac{9}{2}}\sqrt{M}\right)$, and
  $\delta>0$ sufficiently small such that for all
  $r\ge r_{\EventHorizonFuture}(1-\delta_{\Horizon})$, there exists
  some $c_{\delta}>0$, $c(\delta)=O(\delta)$ such that
  \begin{equation}
    \label{eq:ILED:SdS:KCurrent-with-Hardy}
    \KCurrent{\MorawetzVF, \MorawetzLagrangeCorr, 0}_{\delta, v}[\psi]
    \ge c_{\delta}\left(\frac{M\Delta^2}{r^4}\abs*{\partial_r\psi}^2 + r^{-1}\abs*{\psi}^2\right).
  \end{equation}
  \end{lemma}
\begin{proof}
  Observe that we can rewrite \zcref[noname]{eq:ILED:SdS:KCurrent-with-Hardy:def} to see that
  \begin{align*}
    \KCurrent{\MorawetzVF_{M,a,\Lambda}, \MorawetzLagrangeCorr_{*}, 0}_{\delta, v}[\psi]
    ={}& \left(1-\delta\right)\KCurrentR\abs*{\partial_r\psi + \frac{\abs*{q}^2}{2(1-\delta)\KCurrentR}v(r)\psi}^2
         - \frac{\abs*{q}^4}{4(1-\delta)\KCurrentR}v^2\abs*{\psi}^2\\
       & +\left(
         \KCurrentLTwo
         + (1-\delta)\frac{(r-3M)^2}{r^3}
         \right)\abs*{\psi}^2
         + \frac{\abs*{q}^2}{2}\left(\partial_rv + \frac{2v}{r}\right)\abs*{\psi}^2.
  \end{align*}
  As a result,
  \begin{equation*}
    \KCurrent{\MorawetzVF_{M,a,\Lambda}, \MorawetzLagrangeCorr_{*}, 0}_{\delta, v}[\psi]
    \ge{} \left(
      \KCurrentLTwo
      + (1-\delta)\frac{(r-3M)^2}{r^3}
      + \frac{\abs*{q}^2}{2}\left(\partial_rv + \frac{2v}{r}\right)
      - \frac{\abs*{q}^4}{4(1-\delta)\KCurrentR}v^2
    \right)\abs*{\psi}^2.
  \end{equation*}
  Thus, it suffices to show that for $\frac{a}{M}\ll 1$ and $\delta>0$
  sufficiently small, there exists a function $v(r)$ such that
  \begin{equation*}
    \abs*{q}^2E\vcentcolon= \KCurrentLTwo
    + (1-\delta)\frac{(r-3M)^2}{r^3}
      + \frac{\abs*{q}^2}{2}\left(\partial_rv + \frac{2v}{r}\right)
      - \frac{\abs*{q}^4}{4(1-\delta)\KCurrentR}v^2> 0.
  \end{equation*}
  By continuity, it suffices to find such a $v(r)$ such that $E>0$ for
  $\delta=0$ and $a=0$. Recall that we have the following exact computations on \SdS: 
  \begin{align*}
    \KCurrentLTwo &= \KCurrentLTwo_0 + \KCurrentLTwo_1,\\
    2\KCurrentLTwo_0&= \frac{M\left(54M^2-46Mr+9r^2\right)}{r^4}                      
                      + \frac{5M\Lambda}{3} - \frac{4r^3\Lambda^2}{9},\\
    \KCurrentLTwo_1&=\frac{4(r-2M)(r-3M)(r-4M)}{r^4}
                     + \frac{2M(r-3M)}{3r}\Lambda
                     + \frac{2}{9}r^2(r-3M)\Lambda^2,\\
    \KCurrentR &= 3M\Upsilon^2.
  \end{align*}
  Let us now denote $\widetilde{v} = r^2v$. Then
  \begin{equation*}
    E = r^{-2}\KCurrentLTwo + \frac{1}{r^3}\left(1-\frac{3M}{r}\right)^2
    + \frac{1}{2r^2}\partial_r\widetilde{v} - \frac{1}{12Mr^2\Upsilon^2}\widetilde{v}^2.
  \end{equation*}
  Using the values for $\KCurrentLTwo$ in
  \zcref[noname]{eq:Morawetz:SdS:bulk-coefficient-conditions}, we have that
  \begin{align*}
    E ={}& \frac{8r^3 - 63M r^2 + 162M^2r - 138M^3}{2r^6}
           + \frac{1}{r^3}\left(1-\frac{3M}{r}\right)^2
           + \frac{M(3r-4M)}{2r^3}\Lambda
           - \frac{2}{3}M\Lambda^2
       \\&   + \frac{1}{2r^2}\partial_r\widetilde{v}
           - \frac{1}{12Mr^2\Upsilon^2}\widetilde{v}^2\\
    ={}& \frac{10r^3 - 75 M r^2 + 180M^2r - 138 M^3}{2r^6}
         + \frac{M(3r-4M)}{2r^3}\Lambda
         - \frac{2}{3}M\Lambda^2
         + \frac{1}{2r^2}\partial_r\widetilde{v}
         - \frac{1}{12Mr^2\Upsilon^2}\widetilde{v}^2.
  \end{align*}  
  Thus, we have that
  \begin{align*}
    E = \frac{10r^3 - 75 M r^2 + 180M^2r - 138 M^3}{2r^6}
         + \frac{M(3r-4M)}{2r^3}\Lambda
         - \frac{2}{3}M\Lambda^2    
         + \frac{1}{2r^2}\partial_r\widetilde{v}
         - \frac{1}{12Mr^2\Upsilon^2}\widetilde{v}^2.
  \end{align*}
  We then have, denoting $\widetilde{v} = \Upsilon\mathfrak{v}$,
  \begin{equation}
    \label{eq:ILED:SDS:KCurrent-with-Hardy:E-aux}
    \begin{split}
      E ={}& \frac{10r^3 - 75 M r^2 + 180M^2r - 138 M^3}{2r^6}
             - \frac{1}{12Mr^2}\mathfrak{v}^2\\
           & + \frac{M(3r-4M)}{2r^3}\Lambda
             - \frac{2}{3}M\Lambda^2    
             + \frac{1}{2r^2}\left(\frac{2M}{r^2} - \frac{2\Lambda}{3}r\right)\mathfrak{v}
             + \frac{1}{2r^2}\left(1-\frac{2M}{r} - \frac{\Lambda}{3}r^2\right)\partial_r\mathfrak{v}
             .
    \end{split}
  \end{equation}
  Observe that the first line of
  \zcref[noname]{eq:ILED:SDS:KCurrent-with-Hardy:E-aux} is $O(r^{-3})$, while
  the second line is $O(r^{-3}\sqrt{\Lambda})$. Thus, for $\Lambda$
  sufficiently small, to prove that $E>0$, it suffices to choose
  $\mathfrak{v}$ such that the first line is positive. To this end,
  define
  \begin{equation*}
    E_0\vcentcolon= \frac{10r^3 - 75 M r^2 + 180M^2r - 138 M^3}{2r^6}
    - \frac{1}{12Mr^2}\mathfrak{v}^2.
  \end{equation*}
  We now consider $\mathfrak{v} = A\sqrt{\frac{M}{r}}$, so that
  \begin{equation*}
    E_0 ={} \frac{10r^3 - 75 M r^2 + 180M^2r - 138 M^3}{2r^6}
    - \frac{A}{12r^3}\left(-18M^{\frac{3}{2}}r^{-\frac{3}{2}} + 3\sqrt{M}r^{-\frac{1}{2}} + A\right).
  \end{equation*}
  Then, let us define $x=\frac{r}{2M}$, so that 
  \begin{align*}
    2r^2E_0 ={}&  \frac{40x^3 - 150 x^2 + 180x - 69}{8Mx^4}
                 + \frac{1 }{24 M x}\left(-2A^2 +9\sqrt{2}x^{- \frac{3}{2}} -3\sqrt{2}x^{-\frac{1}{2}}\right). 
  \end{align*}
  We can easily check that for $A=\frac{1}{\sqrt{2}}$, this is
  positive and strictly increasing for all $x\ge 1$. Moreover, exactly
  at $x=1$, we see that $r^3E_0 = \frac{1}{3}$.
    Since $\Upsilon(2M)<0$, we have that in fact, $E>0$ in the entirety
  of the stationary region of \SdS{}, and by perturbation, in the
  entirety of the stationary region for all slowly-rotating \KdS, as
  desired.
  It is easy to then verify that
  $v = O(\Upsilon r^{-\frac{5}{2}} M^{\frac{1}{2}}) = O(\Delta
  r^{-\frac{9}{2}} M^{\frac{1}{2}})$.
\end{proof}

We are now ready to conclude the proof of 
\zcref[cap]{prop:Morawetz:outside-trapping:KCurrent-coercive}.
\begin{proof}[Proof of Proposition \ref{prop:Morawetz:outside-trapping:KCurrent-coercive}]
  Having proven \zcref[cap]{lemma:ILED:SDS:KCurrent-with-Hardy}, we now
  have that
  \begin{equation}
    \begin{split}
      \abs*{q}^2\KCurrent{\MorawetzVF_{M,a,\Lambda},\MorawetzLagrangeCorr_*,\MorawetzOneForm_{M,a,\Lambda}}[\psi]
      ={}& \delta\mathcal{A}\abs*{\partial_r\psi}^2
           + \delta\KCurrentAngular{\MorawetzVF_{M,a,\Lambda},\MorawetzLagrangeCorr_*,0}[\psi]
           +  \KCurrent{\MorawetzVF_{M,a,\Lambda}, \MorawetzLagrangeCorr_*, 0}_{\delta, v}[\psi]
      \\
         & + (1-\delta)\left(\KCurrentAngular{\MorawetzVF_{M,a,\Lambda},\MorawetzLagrangeCorr_*,0}[\psi]
           - \frac{(r-3M)^2}{r^3}\abs*{\psi}^2\right)\\
           & + \mathfrak{E}^{\MorawetzVF_{M,a,\Lambda}, \MorawetzLagrangeCorr_{*},0}[\psi]
             + O(ar^{-3})\abs*{\psi}^2,
    \end{split}  
  \end{equation}
  where the integral of all but the last two terms on the right-hand side is positive. It remains to gain coerciveness of the $\HawkingVF$ derivatives. To
  do so, we borrow positivity from the angular and radial derivatives
  in
  $\KCurrent{\MorawetzVF_{M,a,\Lambda},\MorawetzLagrangeCorr_*,\MorawetzOneForm_{M,a,\Lambda}}[\psi]$
  via the Lagrangian correction. Consider
  \begin{equation}
    \label{eq:Morawetz:outside-trapping:q1-def}
    \MorawetzLagrangeCorr_{\delta_1} \vcentcolon= -\delta_1\frac{M\Delta(r-3M)^2}{r^2(r^2+a^2)^2}.
  \end{equation}
  Then
  \begin{equation*}
    \begin{split}
      \abs*{q}^2\KCurrent{0,\MorawetzLagrangeCorr_{\delta_1},0}_0[\psi]
      ={}& \delta_1\frac{M(r-3M)^2}{r^2}\abs*{\HawkingVF\psi}^2
           - \delta_1\frac{M\Delta^2(r-3M)^2}{r^2(r^2+a^2)^2}\abs*{\partial_r\psi}^2\\
         & - \delta_1\frac{M\Delta (r-3M)^2}{r^2(r^2+a^2)^2}\abs*{q}^2\mathcal{O}^{\alpha\beta}\HorCovDeriv_{\alpha}\psi\cdot\HorCovDeriv_{\beta}\psi
           -\delta_1\abs*{q}^2\frac{M\Delta(r-3M)^2}{r^4(r^2+a^2)}V\abs*{\psi}^2
      \\
         &- \frac{\delta_1}{2}\abs*{q}^2\nabla^\alpha\partial_\alpha \frac{M\Delta(r-3M)^2}{r^4(r^2+a^2)}|\psi|^2.
    \end{split}      
  \end{equation*}
    Then, defining
  \begin{equation}
    \label{eq:Morawetz:outside-trapping:q-def}
    \MorawetzLagrangeCorr_{M, a, \Lambda}\vcentcolon= \MorawetzLagrangeCorr_* + \MorawetzLagrangeCorr_{\delta_1},
  \end{equation}
  we have that
  \begin{equation}
    \label{eq:ILED-near:sum-of-squares:SdS}
    \begin{split}
      \abs*{q}^2\KCurrent{\MorawetzVF_{M, a, \Lambda}, \MorawetzLagrangeCorr_{M, a, \Lambda},\MorawetzOneForm_{M,a,\Lambda}}[\psi]
      ={}& \alpha_{M, a, \Lambda}^2\abs*{\HawkingVF\psi}^2
           + \beta_{M, a, \Lambda}^2 \abs*{\HprVF\psi}^2
           + \KCurrentAngular{\MorawetzVF_{M,a,\Lambda},\MorawetzLagrangeCorr_{M,a,\Lambda},\MorawetzOneForm_{M,a,\Lambda}}[\psi] 
      \\
         & + (1-\delta)\left(\KCurrentAngular{\MorawetzVF_{M,a,\Lambda},\MorawetzLagrangeCorr_{*},0}[\psi] - \frac{(r-3M)^2}{r^3}\abs*{\psi}^2\right)\\
         & + \KCurrent{\MorawetzVF_{M,a,\Lambda},\MorawetzLagrangeCorr_0,0}_{\delta,v}[\psi]
           + \mathfrak{E}^{\MorawetzVF_{M,a,\Lambda},\MorawetzLagrangeCorr_*,0}[\psi],
    \end{split}        
  \end{equation}
  where
  \begin{equation}
    \label{eq:Morawetz:alpha-beta-decomp}
    \begin{split}
      \alpha_{M, a, \Lambda}^2
      &=  \delta_1\frac{M(r-3M)^2}{r^2},\\
      \beta_{M, a, \Lambda}^2
      &=  \delta\left( \frac{r^2+a^2}{\Delta} \right)^2\KCurrentR
        - \delta_1\frac{M(r-3M)^2}{r^2}, \\
      \KCurrentAngular{\MorawetzVF_{M,a,\Lambda},\MorawetzLagrangeCorr_{M,a,\Lambda},\MorawetzOneForm_{M,a,\Lambda}}[\psi] 
      &= \delta\KCurrentAngular{\MorawetzVF_{M,a,\Lambda},\MorawetzLagrangeCorr_*,0}[\psi]
        -\delta_1\frac{M\Delta(r-3M)^2}{r^2(r^2+a^2)^2}\abs*{q}^2\mathcal{O}^{\alpha\beta}\HorCovDeriv_{\alpha}\psi\cdot\HorCovDeriv_{\beta}\psi.
    \end{split}
  \end{equation}
  Then, recalling the exact value of $\mathcal{A}$ from
  \zcref[noname]{eq:Morawetz:K-prelim-computations}, we observe that
  \begin{align*}
    \mathcal{A}_{\delta_1}
    \vcentcolon={}& \delta\mathcal{A} - \delta_1\frac{M\Delta^2(r-3M)^2}{r^2(r^2+a^2)^2}
    ={} \frac{\Delta^2}{r^2(r^2+a^2)}\left(
         3M\delta
         - \delta_1\frac{M(r-3M)^2}{r^2+a^2}
         \right),
  \end{align*}
  and that for $\delta_1$ sufficiently small with respect to $\delta$, there exists some $c_0$ such that 
  \begin{equation*}
    3M\delta - \delta_1\frac{M(r-3M)^2}{r^2+a^2}> c_0 M.
  \end{equation*}
  As a result for $\delta_1$ sufficiently small with respect to $\delta$,
  \begin{equation*}
    \mathcal{A}_{\delta_1}\abs*{\partial_r\psi}^2 > \delta M \abs*{\HprVF \psi}^2.
  \end{equation*}
  Moreover, for $\delta_1$ sufficiently small (with respect to $\delta$), we have that
  \begin{equation*}
    \alpha_{M,a,\Lambda}^2, \beta^2_{M,a,\Lambda} > 0, \qquad
    \KCurrentAngular{\MorawetzVF_{M,a,\Lambda},\MorawetzLagrangeCorr_{M,a,\Lambda},\MorawetzOneForm_{M,a,\Lambda}}[\psi] \gtrsim \CartarOp^{\alpha\beta}\HorCovDeriv_{\alpha}\psi\cdot\HorCovDeriv_{\beta}\psi.
  \end{equation*}

  Moreover, observe that for $a$ sufficiently small, we have from
  \zcref[cap]{lemma:ILED:SdS:Bulk:principal-terms}
  on
  $\Manifold_{\cancel{\operatorname{trap}}}(\tau_1,\tau_2)$,
  \begin{equation*}
    \abs*{\mathfrak{E}^{\MorawetzVF_{M,a,\Lambda},\MorawetzLagrangeCorr_{*},0}[\psi]} < \frac{1}{2}
    \left(\alpha_{M,a,\Lambda}^2 \abs*{\HawkingVF\psi}^2 + \KCurrentAngular{\MorawetzVF_{M,a,\Lambda},\MorawetzLagrangeCorr_{M,a,\Lambda},0}[\psi]\right).
  \end{equation*}
  Thus,
  $\mathfrak{E}^{\MorawetzVF_{M,a,\Lambda},\MorawetzLagrangeCorr_{*},0}[\psi]$
  can be absorbed by
  $\alpha_{M,a,\Lambda}^2\abs*{\HawkingVF\psi}^2 +
  \KCurrentAngular{\MorawetzLagrangeCorr_{M,a,\Lambda},
    \MorawetzLagrangeCorr_{M,a,\Lambda},\MorawetzOneForm_{M,a,\Lambda}}[\psi]$
  on $\Manifold_{\cancel{\operatorname{trap}}}(\tau_1,\tau_2)$. We also observe that from
  \zcref[cap]{lemma:ILED:SDS:KCurrent-with-Hardy}, we have that for
  $\delta_1$ sufficiently small relative to $\delta$,
  \begin{equation*}
    \delta_1\abs*{q}^2\left(
      \frac{M\Delta(r-3M)^2}{r^4(r^2+a^2)}V\abs*{\psi}^2
      + \frac{1}{2}\nabla^{\alpha}\partial_{\alpha}\frac{M\Delta(r-3M)^2}{r^4(r^2+a^2)}
    \right)\abs*{\psi}^2
    = \delta_1O(r^{-2})\abs*{\psi}^2
    \lesssim \KCurrent{\MorawetzVF_{M,a,\Lambda},\MorawetzLagrangeCorr_0,0}_{\delta,v}[\psi],
  \end{equation*}
  where we used the fact that $V = O(r^{-2})$, recalling that
  $\Lambda\lesssim r^{-2}$. Finally, integrating and applying the Poincar\'{e} inequality from
  \zcref[cap]{lemma:Poincare} concludes the proof of
  \zcref[cap]{prop:Morawetz:outside-trapping:KCurrent-coercive}.
\end{proof}

The key property of the choice of multiplier in
$(\MorawetzVF_{M,a,\Lambda}, \MorawetzLagrangeCorr_{M,a,\Lambda},
\MorawetzOneForm_{M,a,\Lambda})$ is its bulk property in
\zcref[noname]{eq:Morawetz:outside-trapping:KCurent-coercive}. However, to
close the proof of the Morawetz estimate, we will also need the
following, easily shown properties of the boundary terms that appear
in the estimate.
\begin{corollary}
  \label{coro:Morawetz:outside-trapping:boundary-terms}
  The
  $(\MorawetzVF_{M,a,\Lambda}, \MorawetzLagrangeCorr_{M,a,\Lambda},
  \MorawetzOneForm_{M,a,\Lambda})$ of \zcref[cap]{prop:Morawetz:outside-trapping:KCurrent-coercive} can be chosen
  so that
  \begin{align}   
    \abs*{\int_{\Sigma(\tau)} \JCurrent{\MorawetzVF_{M,a,\Lambda}, \MorawetzLagrangeCorr_{M,a,\Lambda}, \MorawetzOneForm_{M,a,\Lambda}}[\psi]\cdot N_{\Sigma}}
    \lesssim{}& \EnergyHorizonDeg[\psi](\tau),
                \label{eq:Morawetz:outside-trapping:Sigma-energy-bound}
    \\
    \abs*{\int_{\mathcal{A}(\tau_1,\tau_2)}\JCurrent{\MorawetzVF_{M,a,\Lambda}, \MorawetzLagrangeCorr_{M,a,\Lambda}, \MorawetzOneForm_{M,a,\Lambda}}[\psi]\cdot N_{\mathcal{A}}}
    \lesssim{}&\delta_{\Horizon}\SpacelikeFlux_{\mathcal{A}}[\psi](\tau_1,\tau_2)
                + \int_{\mathcal{A}(\tau_1,\tau_2)}\abs*{\nabla_{4}\psi}^2 
                ,
           \label{eq:Morawetz:outside-trapping:ACal-energy-bound}
    \\
    \abs*{\int_{\SigmaStar(\tau_1,\tau_2)}\JCurrent{\MorawetzVF_{M,a,\Lambda}, \MorawetzLagrangeCorr_{M,a,\Lambda}, \MorawetzOneForm_{M,a,\Lambda}}[\psi]\cdot N_{\SigmaStar}}
    \lesssim{}& \delta_{\Horizon}\SpacelikeFlux_{\SigmaStar}[\psi](\tau_1,\tau_2)
               + \int_{\SigmaStar(\tau_1,\tau_2)}\abs*{\nabla_{3}\psi}^2
                .
           \label{eq:Morawetz:outside-trapping:SigmaStar-energy-bound}                                     
  \end{align}
\end{corollary}
\begin{proof}
  The estimates follow from the definitions of
  $\EnergyHorizonDeg[\psi](\tau)$ in
  \zcref[noname]{eq:horizon-degenerate-energy}, and the explicit choices in
  \zcref[noname]{eq:ILED:SdS:X-def}, \zcref[noname]{eq:ILED:SdS:X-def:f-def},
  \zcref[noname]{eq:Morawetz:outside-trapping:q-def}, and
  \zcref[cap]{lemma:ILED:SDS:KCurrent-with-Hardy}.
\end{proof}

\subsection{Morawetz estimate: trapping}
\label{sec:Morawetz:KdS}

In this section, we will fill in the necessary analysis in the
trapping region of \KdS{} to prove
\zcref[cap]{prop:Morawetz:KdS:main}. To prove a Morawetz estimate on the trapping region of \KdS, the main
change that is needed is a finer analysis of the trapping behavior. To
accommodate this, the main goal is to perturb
$\left( \MorawetzVF_{M,a,\Lambda}, \MorawetzLagrangeCorr_{M,a,\Lambda}
\right)$ by
$\left( \widetilde{\MorawetzVF}, \widetilde{\MorawetzLagrangeCorr}
\right)$ such that there exists a sum of squares expression
\begin{equation*}
  H_p\MorawetzSym + p\MorawetzLagrangeCorrSym = \sum_{j=1}^k\SquareDecomp_j^2,
\end{equation*}
where $\MorawetzSym$ is the symbol of
$\MorawetzVF_{M,a,\Lambda} + \widetilde{\MorawetzVF}$ and
$\MorawetzLagrangeCorrSym$ is the symbol of
$\MorawetzLagrangeCorr_{M,a,\Lambda}+\widetilde{\MorawetzLagrangeCorr}$,
and $p$ denotes the principal symbols of the scalar wave
operator. Combined with the pseudo-differential divergence theorem in
\zcref[cap]{lemma:PsiDO-divergence-theorem}, this will show that we
have a bulk that is positive at the principal order. To extend this
positivity to lower-order terms, and to correct for the fact that the
Regge-Wheeler equation is not exactly the scalar wave, we need to gain
control of lower-order terms. This is nontrivial due to the
degeneration of some of the top-order derivatives in the bulk term at
the trapped set. However, we will take advantage of the principal-level bulk by using the fact that $\partial_r(r-\rTrapping)\in \OpClass^{-1}$ to control the lower-order terms (see
\zcref[cap]{sec:Morawetz:lower-order-bulk} where this idea is used in
full several times).

We now first prove several critical lemmas and then show how these
combine later afterwards to give a Morawetz estimate for generalized
Regge-Wheeler on exact \KdS.

\subsubsection{Principal level bulk terms}
\label{sec:Morawetz:principal-order-bulk}

We begin with an analysis of the principal level bulk terms that arise
in the trapping region.  The analysis is similar to how the Morawetz
estimate was proven in
\cite{fangLinearStabilitySlowlyRotating2026}. However, here we more
carefully keep track of powers of $r$, and note that while the
perturbative argument is the same, in
\cite{fangLinearStabilitySlowlyRotating2026}, the perturbation was
done off a vectorfield that would prove the Morawetz estimate for
solutions to the wave equation on \SdS, while in the present setting,
the perturbation is essentially off of a vectorfield that would prove
the Morawetz estimate for axially-symmetric solutions to the wave
equation on \KdS.

\begin{lemma}
  \label{lemma:ILED-KdS:Bulk}
  
  There exists
  \begin{enumerate}
  \item a pseudo-differential operator
    $\widetilde{\MorawetzVF}\in \MixedOpClass{1}{1}(\mathcal{M};
    \realHorkTensor{2}(\mathcal{M}))$ such that
    \begin{equation}
      \label{eq:ILED-KdS:X-pert-def}
      \widetilde{\MorawetzVF} = \widetilde{\MorawetzVF}_0\HorCovDeriv_t +  \widetilde{\MorawetzVF}_1,
    \end{equation}
    where
    $\widetilde{\MorawetzVF}_i\in \TanOpClass{i}(\mathcal{M};
    \realHorkTensor{2}(\mathcal{M}))$ is an anti-symmetric operator
    with symbol $\widetilde{\MorawetzSym}_i$,
  \item a pseudo-differential operator
    $\widetilde{\MorawetzLagrangeCorr}\in
    \MixedOpClass{0}{1}(\mathcal{M}; \realHorkTensor{2}(\mathcal{M}))$
    such that
    \begin{equation}
      \label{eq:ILED-KdS:q-pert-def}
      \widetilde{\MorawetzLagrangeCorr}
      = \widetilde{\MorawetzLagrangeCorr}_0
      + \widetilde{\MorawetzLagrangeCorr}_{-1}\HorCovDeriv_{t},
    \end{equation}
    where
    $\widetilde{\MorawetzLagrangeCorr}_i\in
    \TanOpClass{i}(\mathcal{M}, \realHorkTensor{2}(\mathcal{M}))$ is a self-adjoint
    pseudo-differential operator with symbol
    $\widetilde{\MorawetzLagrangeCorrSym}_i$,
  \end{enumerate}
  such that
  \begin{equation}
    \label{eq:Morawetz:KdS:bulk:square-sum-decomposition}
    \abs*{q}^2\left(
      \frac{1}{2\ImagUnit}H_p\MorawetzSym + \PrinSymb \MorawetzLagrangeCorrSym
    \right)
    = \sum_{j=1}^{7}\SquareDecomp_j^2,
  \end{equation}
  where
  $\SquareDecomp_j\in \MixedSymClass{1}{1}(T^{*}\mathcal{M};
  \pi^{*}(\End(\realHorkTensor{2}(\mathcal{M})))$, are principally scalar, and
  \begin{gather*}
    \MorawetzSym = \MorawetzSym_{M, a, \Lambda} + a\widetilde{\MorawetzSym}, \qquad
    \MorawetzSym_{M, a, \Lambda} = \ImagUnit \mathcal{F}_{M, a, \Lambda}(r)\xi, \\
    \MorawetzLagrangeCorrSym = \MorawetzLagrangeCorrSym_{M, a, \Lambda} + a\widetilde{\MorawetzLagrangeCorrSym},\qquad
    \MorawetzLagrangeCorrSym = \MorawetzLagrangeCorrSym_{M, a, \Lambda} + a\widetilde{\MorawetzLagrangeCorr},\qquad
    \MorawetzLagrangeCorrSym_{M, a, \Lambda} = \MorawetzLagrangeCorr_{M, a, \Lambda} - \frac{1}{2}\CovariantDeriv_{\Metric_{M, a, \Lambda}}\cdot \MorawetzVF_{M, a, \Lambda}.
  \end{gather*}
  Moreover, the decomposition in
  \zcref[noname]{eq:Morawetz:KdS:bulk:square-sum-decomposition} extends the
  decomposition in \zcref[noname]{eq:ILED-near:sum-of-squares:SdS} in the
  sense that if $a=0$, then
  $\widetilde{\MorawetzLagrangeCorrSym} = \widetilde{\MorawetzSym} = 0$, and 
  \begin{equation}    
    \sum_{j=1}^7\SquareDecomp_j^2 = \mathfrak{k}^{\alpha\beta}\zeta_{\alpha}\zeta_{\beta},
  \end{equation}
  where $\zeta = \{\sigma,\xi,\eta\}\in T^{*}\Manifold$, and
  $\mathfrak{k}^{\alpha\beta}\HorCovDeriv_{\alpha}\psi\HorCovDeriv_{\beta}\psi
  =
  \KCurrent{\MorawetzVF_{M,a,\Lambda},\MorawetzLagrangeCorr_{M,a,\Lambda},\MorawetzOneForm_{M,a,\Lambda}}[\psi] + O(\abs*{\psi}^2)$.
\end{lemma}
\begin{remark}
  Observe that  the specific forms of
  $\widetilde{\MorawetzVF}$, $\widetilde{\MorawetzLagrangeCorr}$ in
  \zcref[noname]{eq:ILED-KdS:X-pert-def} and \zcref[noname]{eq:ILED-KdS:q-pert-def}
  are not unique due to their pseudodifferential nature.
\end{remark}
\begin{remark}
  The definitions of $\MorawetzSym_{M,a,\Lambda}$,
  $\MorawetzLagrangeCorrSym_{M,a,\Lambda}$ are chosen so that
  \begin{equation*}
    \WeylQ{\MorawetzSym_{M,a,\Lambda} + \MorawetzLagrangeCorrSym_{M,a,\Lambda}} = \MorawetzVF_{M,a,\Lambda} + \MorawetzLagrangeCorr_{M,a,\Lambda}. 
  \end{equation*}
\end{remark}
We remark that the main difficulty in the ensuing proof is the
limitation that $\widetilde{\mathfrak{x}}$ must be a
mixed-pseudodifferential operator, i.e. that
$\widetilde{\mathfrak{x}}\in \MixedSymClass{1}{1}(T^{*}\mathcal{M};
\pi^{*}(\End(\realHorkTensor{2}(\mathcal{M})))$. This is a necessary
restriction in our method to allow for the application of the
divergence theorem in \zcref[cap]{lemma:PsiDO-divergence-theorem}. A
natural choice for $\widetilde{\mathfrak{x}}$ is
\begin{equation*}
  \MorawetzSym'\vcentcolon=\ImagUnit \Delta\abs*{q}^{-2}\left(r-\rTrapping_{M,a,\Lambda}(\sigma,\FreqPhi)\right)\xi
  = \ImagUnit\frac{r-\rTrapping_{M,a,\Lambda}(\sigma,\FreqPhi)}{2\abs*{q}^2}H_{\abs*{q}^2\PrinSymb}r.
\end{equation*}
However,
$\MorawetzSym'\not\in \MixedSymClass{1}{1}(T^{*}\mathcal{M};
\pi^{*}(\End(\realHorkTensor{2}(\mathcal{M})))$. To make a good choice
for $\widetilde{\MorawetzSym}$, we instead divide
$\MorawetzSym'-\MorawetzSym_{M,a,\Lambda}$ by $p$ using the Malgrange
preparation theorem and use the remainder as
$\widetilde{\MorawetzSym}$, taking advantage of the fact that
$\widetilde{\MorawetzSym}+\MorawetzSym_{M,a,\Lambda}$ and
$\MorawetzSym'$ are equal up to a multiple of $p$.

To show that $H_p\left( \MorawetzSym_{M,a,\Lambda}+\widetilde{\MorawetzSym} \right)$ is positive, it suffices to
divide $H_p\widetilde{\MorawetzSym}$ by $p$, again using the Malgrange
representation theorem to write that
\begin{equation*}
  H_p\left( \MorawetzSym_{M,a,\Lambda}+\widetilde{\MorawetzSym} \right)
  ={} \gamma_2+\gamma_1\sigma + qp.
\end{equation*}
Using the fact that by definition,
$\MorawetzSym_{M,a,\Lambda}+\widetilde{\MorawetzSym}=\MorawetzSym'$ if
$p=0$, we can then solve an explicit system of equations for
$\gamma_1,\gamma_2$ by considering
$\sigma=\sigma_1,\sigma_2$. Finally, since we only need a good
multiplier when $p=0$, we will show that for some sufficiently good
choice of $\mathfrak{b}$, there exists a sum of square decomposition
\begin{equation*}
  \gamma_2+\gamma_1\sigma + \mathfrak{b}p = \sum_{i=1}^k\SquareDecomp_i^2.
\end{equation*}

\begin{proof}[Proof of Lemma \ref{lemma:ILED-KdS:Bulk}]
For the sake of simplifying some of our ensuing calculations, define
  \begin{equation*}
    \MorawetzLagrangeCorrSym'_{M, a, \Lambda} = \MorawetzLagrangeCorrSym_{M, a, \Lambda}
    - 2\PoissonB{\ln q, \MorawetzSym_{M, a, \Lambda}}, \qquad
    \tilde{\MorawetzLagrangeCorrSym}'_{M, a, \Lambda} =
    \tilde{\MorawetzLagrangeCorrSym}_{M, a, \Lambda} -2\PoissonB{\ln q, \tilde{\MorawetzSym}}, 
  \end{equation*}
  so that
  \begin{equation*}
    \abs*{q}^2\left(
      \frac{1}{2\ImagUnit}H_\PrinSymb
      (\MorawetzSym_{M,a,\Lambda}+a\tilde{\MorawetzSym})
      +(\MorawetzLagrangeCorrSym_{M,a,\Lambda}+a\tilde{\MorawetzLagrangeCorrSym})\PrinSymb  
    \right)
    = \frac{1}{2\ImagUnit}H_{\abs*{q}^2\PrinSymb} (\MorawetzSym_{M,a,\Lambda}+a\tilde{\MorawetzSym}) +
    \left(\tilde{\MorawetzLagrangeCorrSym}'
      +a\tilde{\MorawetzLagrangeCorrSym}'\right)(\abs*{q}^2\PrinSymb).
  \end{equation*}
  We first choose $\tilde{\MorawetzSym}$ so that
  $H_{\RescaledPrinSymb}
  (\MorawetzSym_{M,a,\Lambda}+a\tilde{\MorawetzSym})$ vanishes at the
  trapped set $\TrappedSet$. The most immediate choice for this is the
  symbol
  \begin{equation*}
    \MorawetzSym'\vcentcolon=\ImagUnit \Delta\abs*{q}^{-2}\left(r-\rTrapping_{M,a,\Lambda}(\sigma,\FreqPhi)\right)\xi
    = \ImagUnit\frac{r-\rTrapping_{M,a,\Lambda}(\sigma,\FreqPhi)}{2\abs*{q}^2}H_{\abs*{q}^2\PrinSymb}r.
  \end{equation*}
  This symbol is clearly well defined and smooth in a neighborhood of
  the trapped set. Moreover, we can calculate that on the
  characteristic set $\PrinSymb=0$, we have
  \begin{equation*}
    2 H_{\RescaledPrinSymb}\MorawetzSym' =
    \left(\frac{1}{\abs*{q}^2} - \frac{2(r-\rTrapping_{M,a,\Lambda})\partial_r\abs*{q}}{\abs*{q}^3}\right)(H_{\RescaledPrinSymb} r)^2
    + \frac{r-\rTrapping_{M,a,\Lambda}(\sigma,\FreqPhi)}{\abs*{q}^2}H_{\RescaledPrinSymb}^2r.
  \end{equation*}
  We can compute that for
  $\PrinSymb=0$, we have that
  \begin{equation*}
    H_{\RescaledPrinSymb}^2r
    = 2\Delta\partial_r\left(\frac{(1+\gamma)^2}{\Delta}\left((r^2+a^2)\sigma +a\FreqPhi\right)\right).
  \end{equation*}
  Since $\rTrapping_{M,a,\Lambda}$ is the unique minimum of
  $\frac{(1+\gamma)^2}{\Delta}\left((r^2+a^2)\sigma
    +a\FreqPhi\right)$, and we are in a $\delta_{\operatorname{trap}}$ neighborhood of
  $r=3M$, there exist positive homogeneous symbols
  $\alpha,\beta\in \Psi^0_{\hom}(r,\sigma,\eta_{\varphi})$ such that
  on $\PrinSymb=0$, near $r=3M$,
  \begin{equation}
    \label{eq:ILED-near:sum-of-squares:KdS-ILED-sym-characteristic}
    H_{\RescaledPrinSymb}\MorawetzSym' =
    \alpha^2(r,\sigma, \FreqPhi) (r-\rTrapping_{M,a,\Lambda})^2
    + \beta^2(r,\sigma, \FreqPhi)\xi^2.
  \end{equation}  
  Unfortunately, the problem with $\MorawetzSym'$ is that it is not a
  polynomial in $\sigma$ unless $a=0$. Thus, for $a\neq 0$,
  $\MorawetzSym'$ cannot be directly used in conjunction with our
  integration-by-parts or divergence theorem method to produce a
  Morawetz estimate. We will overcome this difficulty with the aid of
  the Malgrange preparation theorem.
  Observe that we defined $\MorawetzSym$ so that it is smooth in $a$,
  and so that $\MorawetzSym'-\MorawetzSym_{M,a,\Lambda} \in a\MixedSymClass{2}{0}(T^{*}\Manifold)$. Thus the Malgrange preparation theorem gives us the existence of
  homogeneous
  $\tilde{\MorawetzSym}_i\in \TanSymClass{i}(T^{*}\mathcal{M})$, $i=0,1$
  and homogeneous  $\rAux \in S^{-1}(T^{*}\mathcal{M})$
  such that
  \begin{equation}
    \label{eq:ILED-near:x-div-theorem}
    \frac{1}{\ImagUnit}\left(\MorawetzSym' - \MorawetzSym_{M,a,\Lambda}\right)
    = a\left(\tilde{\MorawetzSym}_1
      + \tilde{\MorawetzSym}_0\sigma\right)
    + a \rAux\PrinSymb.
  \end{equation}
  Now, we define
  \begin{equation*}
    \frac{1}{\ImagUnit}\tilde{\MorawetzSym} = \tilde{\MorawetzSym}_1+\tilde{\MorawetzSym}_0\sigma,
  \end{equation*}
  so that on $\PrinSymb=0$ we have $\MorawetzSym = \MorawetzSym_{M,a,\Lambda} + a\tilde{\MorawetzSym} = \MorawetzSym'$. Thus $\MorawetzSym$ is a symbol which is a polynomial in $\sigma$
  and moreover vanishes at the trapped set $\TrappedSet$.
  
  \textit{A priori}, $H_{\RescaledPrinSymb} \tilde{\MorawetzSym}$ is a
  third degree polynomial in $\sigma$. Applying the Malgrange
  preparation theorem again yields that there exist some
  $\gamma_1\in \TanSymClass{1}(T^{*}\mathcal{M}), \gamma_2 \in
  \TanSymClass{2}(T^{*}\mathcal{M})$,
  $f_0\in \TanSymClass{0}(T^{*}\mathcal{M}), f_{-1}\in \sigma
  \TanSymClass{-1}(T^{*}\mathcal{M})$ such that
  \begin{equation*}
    \frac{1}{2\ImagUnit\abs*{q}^2} H_{\RescaledPrinSymb}(\MorawetzSym_{M,a,\Lambda}+a\tilde{\MorawetzSym})
    + \MorawetzLagrangeCorrSym_{M,a,\Lambda}' (\RescaledPrinSymb)
    = \gamma_2+\gamma_1\sigma  + \left(
      e_{M,a,\Lambda} + a(f_0+f_{-1}\sigma)\right)(\sigma-\sigma_1)(\sigma-\sigma_2),
  \end{equation*}
  observing that
  \begin{equation*}
    e_{M,a,\Lambda} \vcentcolon= \delta_1\alpha_{M,a,\Lambda}^2
  \end{equation*}
  is the coefficient for $\sigma^2$ in the expression for
  $\frac{1}{2\ImagUnit}H_{\PrinSymb}\MorawetzSym_{M,a,\Lambda} +
  \MorawetzLagrangeCorrSym_{M,a,\Lambda}\PrinSymb$ (see \zcref[noname]{eq:ILED-near:sum-of-squares:SdS}).                 
  It now remains to demonstrate that $\gamma_2+\gamma_1\sigma
  +e_{M,a,\Lambda}(\sigma-\sigma_1)(\sigma-\sigma_2)$ can be expressed as a
  sum of squares up to some error in
  $a\PrinSymb\MixedSymClass{0}{1}(T^{*}\mathcal{M})$. If this were true, we could write
  \begin{equation}
    \label{eq:ILED-near:sum-of-squares:aux-1}
    \gamma_2+\gamma_1\sigma + e_{M,a,\Lambda}(\sigma-\sigma_1)(\sigma-\sigma_2)
    = \sum \SquareDecomp_j^2 + a(g_0+g_{-1}\sigma)(\sigma-\sigma_1)(\sigma-\sigma_2). 
  \end{equation}
  We could then define $\tilde{\MorawetzLagrangeCorrSym}$ such that
  \begin{equation*}
    \tilde{\MorawetzLagrangeCorrSym}' = -2\left(f_0+g_0+(f_{-1}+g_{-1})\sigma\right),
  \end{equation*}
  so that the
  $a\PrinSymb\MixedSymClass{0}{1}(T^{*}\mathcal{M})$ terms are all
  canceled. We now return to showing
  \zcref[noname]{eq:ILED-near:sum-of-squares:aux-1}. Recall that on
  $\PrinSymb=0$ we have $H_{\RescaledPrinSymb}(\MorawetzSym_{M,a,\Lambda} + a\tilde{\MorawetzSym}) = H_{\RescaledPrinSymb}\MorawetzSym'$. As a result of
  \zcref[noname]{eq:ILED-near:sum-of-squares:KdS-ILED-sym-characteristic}, we
  now have that if $\sigma=\sigma_i$, which in particular implies that for $\PrinSymb=0$,
  \begin{equation*}
    \gamma_2+\gamma_1\sigma =\alpha^2(r, \sigma, \FreqPhi)(r-\rTrapping_{M,a,\Lambda})^2 + \beta^2(r,\sigma,\FreqPhi)\xi^2.
  \end{equation*}
  We can solve for $\gamma_2, \gamma_1$ explicitly now by considering
  the two-dimensional system of equations
  \begin{equation*}
    \begin{split}
      \gamma_2+\gamma_1\sigma_i &= \frac{1}{4}\alpha_i^2(\sigma_1-\sigma_2)^2 + \beta_i^2\xi^2,\\
      \alpha_i &= \frac{2|\sigma_i|}{\sigma_1-\sigma_2}\alpha(r,
                 \sigma_i, \FreqPhi)(r-\rTrapping_{M,a,\Lambda}(\sigma_i, \FreqPhi))\in \TanSymClass{0}(T^{*}\mathcal{M}),\\
      \beta_i &= \beta(r, \sigma_i, \FreqPhi) \in \TanSymClass{0}(T^{*}\mathcal{M}).
    \end{split}
  \end{equation*}
  Solving the system yields 
  \begin{equation}
    \begin{split}
      \label{eq:ILED-near:sum-of-squares:gamma-def}
      \gamma_2 &= \frac{1}{4}(\sigma_1-\sigma_2)(\alpha_2^2\sigma_1 - \alpha_1^2\sigma_2) + \frac{\sigma_1\beta_2^2 - \sigma_2\beta_1^2}{\sigma_1-\sigma_2}\xi^2,  \\
      \gamma_1 &= \frac{1}{4}(\sigma_1-\sigma_2)(\alpha_1^2-\alpha_2^2)+\frac{\beta_1^2 - \beta_2^2}{\sigma_1-\sigma_2}\xi^2. 
    \end{split} 
  \end{equation}
  We first add together the first two terms in $\gamma_i$ to see that
  \begin{equation}
    \label{eq:ILED-near:sum-of-squares:gamma-sum-first-terms}
    \begin{split}
      (\sigma_1-\sigma_2)\left(\alpha_2^2\sigma_1 -\alpha_1^2\sigma_2 + \sigma(\alpha_1^2-\alpha_2^2)\right)
      ={}&(1-{\delta_1})(\alpha_1(\sigma-\sigma_2)-\alpha_2(\sigma-\sigma_1))^2\\
         &+ {\delta_1}\left(\alpha_1(\sigma-\sigma_2)+\alpha_2(\sigma-\sigma_1)\right)^2
         - 4\mathfrak{e}(\sigma-\sigma_1)(\sigma-\sigma_2),
    \end{split}
  \end{equation}
  where $\mathfrak{e} = \frac{(\alpha_1-\alpha_2)^2}{4}+ {\delta_1}\alpha_1\alpha_2$. Recall that when $a = 0$ we have
  $\alpha_1=\alpha_2=\alpha_{M,0,\Lambda}$, $\sigma_2=-\sigma_1$,
  $\beta_1=\beta_2=\beta_{M,0,\Lambda}$, and that
  $\MorawetzLagrangeCorrSym_{M,0,\Lambda} =
  \delta_1\alpha_{M,0,\Lambda}^2$.  This implies that
  \begin{equation}
    \label{eq:ILED-near:sum-of-squares:e-diff}
    \mathfrak{e}-e_{M,a,\Lambda}\in a\MixedSymClass{0}{1}(T^{*}\mathcal{M}),
  \end{equation}
  as desired.  We now add together the second terms in the $\gamma_i$
  given in \zcref[noname]{eq:ILED-near:sum-of-squares:gamma-def}
  \begin{equation}
    \label{eq:ILED-near:sum-of-squares:gamma-sum-second-terms}
    \begin{split}
      \frac{\sigma_1\beta_2^2-\sigma_2\beta_1^2}{\sigma_1-\sigma_2}
      + \sigma\frac{\beta_1^2-\beta_2^2}{\sigma_1-\sigma_2}
      ={}&
           \frac{1}{2}\left(\beta_1^2+\beta_2^2 -Ca\right) +
           \frac{(Ca-\beta_2^2+\beta_1^2)(\sigma-\sigma_2)^2}{2(\sigma_1-\sigma_2)^2}\\
         &+
           \frac{(Ca-\beta_1^2+\beta_2^2)(\sigma-\sigma_1)^2}{2(\sigma_1-\sigma_2)^2}
           + O(a)\PrinSymb.
    \end{split}
  \end{equation}
  Summing \zcref[noname]{eq:ILED-near:sum-of-squares:gamma-sum-first-terms}
  and \zcref[noname]{eq:ILED-near:sum-of-squares:gamma-sum-second-terms}
  together, we have that
  \begin{equation*}
    \begin{split}
      \frac{1}{2\ImagUnit}H_{\RescaledPrinSymb}(
        \MorawetzSym_{M,a,\Lambda}+a\tilde{\MorawetzSym}) + \RescaledPrinSymb\MorawetzLagrangeCorrSym'
      ={}& \frac{1-{\delta_1}}{4}\left(
           \alpha_1(\sigma-\sigma_2)-\alpha_2(\sigma-\sigma_1)
           \right)^2
           + \frac{{\delta_1}}{4}\left(\alpha_1(\sigma-\sigma_2) + \alpha_2(\sigma-\sigma_1)\right)^2\\
      &+ \frac{1}{2}\left(\beta_1^2+\beta_2^2 - Ca\right)\xi^2
        + \frac{(Ca-\beta_2^2+\beta_1^2)(\sigma-\sigma_2)^2}{2(\sigma_1-\sigma_2)^2}\xi^2\\
      &+ \frac{(Ca-\beta_1^2+\beta_2^2)(\sigma-\sigma_1)^2}{2(\sigma_1-\sigma_2)^2}\xi^2
        +a(\sigma-\sigma_1)(\sigma-\sigma_2)\MixedSymClass{0}{1}(T^{*}\mathcal{M}).
    \end{split}
  \end{equation*}
  We then pick
  \begin{align*}
    \SquareDecomp_1^2 &= \frac{\delta_1}{4}\left(\alpha_1(\sigma-\sigma_2)+\alpha_2(\sigma-\sigma_1)\right)^2,\\
    \SquareDecomp_2^2 &= \frac{1}{2}\left(\beta_1^2 + \beta_2^2 - Ca\right)\xi^2,\\
    \SquareDecomp_{3}^2 &= \frac{\FreqTheta^2}{|\FreqAngular|^2 + \Delta\xi^2}\frac{(1-\delta_1)}{4} \left(\alpha_1(\sigma-\sigma_2)-\alpha_2(\sigma-\sigma_1)\right)^2, \\ 
    \SquareDecomp_{4}^2 &= \frac{\FreqPhi^2}{|\FreqAngular|^2 + \Delta\xi^2}\frac{(1-\delta_1)}{4} \left(\alpha_1(\sigma-\sigma_2)-\alpha_2(\sigma-\sigma_1)\right)^2,\\
    \SquareDecomp_{5}^2 &= \frac{\Delta\xi^2}{|\FreqAngular|^2 + \Delta\xi^2}\frac{(1-\delta_1)}{4}\left(\alpha_1(\sigma-\sigma_2)-\alpha_2(\sigma-\sigma_1)\right)^2,\\
    \SquareDecomp_6^2 &= \frac{\left(Ca-\beta_2^2 + \beta_1^2\right)\left(\sigma-\sigma_2\right)^2}{2\left(\sigma_1-\sigma_2\right)^2}\xi^2,\\
    \SquareDecomp_7^2 &= \frac{\left(Ca-\beta_1^2 + \beta_2^2\right)\left(\sigma-\sigma_1\right)^2}{2\left(\sigma_1-\sigma_2\right)^2}\xi^2,
  \end{align*}
  concluding the proof of \zcref[cap]{lemma:ILED-KdS:Bulk}.
\end{proof}

\begin{remark}
    We observe that in particular, when $a = 0$,
  $\SquareDecomp_i$ are all symbols of differential operators, as
  $\evalAt*{\curlyBrace*{\SquareDecomp_i}_{i=3}^7}_{\eta_{\phi}=0}=0$,
  and similarly, when $a=0$, $\SquareDecomp_i$ are all symbols of
  differential operators. Moreover, we see that when $a\eta_{\phi}=0$,
  $\sum_{i=1}^7\SquareDecomp_i^2$ agrees precisely with the expansion
  in \zcref[cap]{prop:Morawetz:outside-trapping:KCurrent-coercive}.
\end{remark}

\begin{corollary}
  \label{coro:Morawetz:elliptic-equivalence}
  The family of symbols $\{\SquareDecomp_j\}_{j=1,7}$ is elliptically
  equivalent with the family of symbols
  $\curlyBrace*{\ell_1(\sigma-\sigma_2), \ell_2(\sigma-\sigma_1),
    \xi^2}$ where $\ell_i$, $\sigma_j$, are as defined in
  \zcref[noname]{eq:trapping-norm:elli-def} and
  \zcref[noname]{eq:trapping-norm:sigmai-def} respectively, in the
  sense that there is a representation of the form
  \begin{equation*}
    \SquareDecomp = \mathfrak{M} \mathfrak{b},\qquad
    \mathfrak{b} =
    \begin{pmatrix*}
      \ell_1(\sigma-\sigma_2)\\
      \ell_2(\sigma-\sigma_1)\\
      \xi
    \end{pmatrix*},
  \end{equation*}
  where the symbol valued matrix $\mathfrak{M}$ has maximum rank 3 everywhere. 
\end{corollary}

\begin{proof}
  We observe that
  \begin{equation*}
    \frac{2\abs*{\sigma_i}}{\sigma_1-\sigma_2}\alpha(r, \sigma_i, \eta_{\varphi})\ell_i
    = \alpha_i.
  \end{equation*}
  As a result, it suffices to show that there exists some
  $\widetilde{\mathfrak{M}}$ such that
  \begin{equation*}
    \mathfrak{a}=\widetilde{\mathfrak{M}}
    \begin{pmatrix*}
      \alpha_1(\sigma-\sigma_2)\\
      \alpha_2(\sigma-\sigma_1)\\
      \xi
    \end{pmatrix*}.
  \end{equation*}
  To this end, consider
  \begin{equation*}
    \widetilde{\mathfrak{M}}
    =
    \begin{pmatrix}
      \frac{\sqrt{\delta_1}}{2} & \frac{\sqrt{\delta_1}}{2}& 0\\
      0& 0 & \frac{1}{\sqrt{2}}\sqrt{\beta_1^2+\beta_2^2-Ca}    \\
      \frac{\FreqTheta}{\sqrt{|\FreqAngular|^2 + \Delta\xi^2}}\frac{\sqrt{1-\delta_1}}{2} & -\frac{\FreqTheta}{\sqrt{|\FreqAngular|^2 + \Delta\xi^2}}\frac{\sqrt{1-\delta_1}}{2} & 0\\
      \frac{\FreqPhi}{\sqrt{|\FreqAngular|^2 + \Delta\xi^2}}\frac{\sqrt{1-\delta_1}}{2} & -\frac{\FreqTheta}{\sqrt{|\FreqAngular|^2 + \Delta\xi^2}}\frac{\sqrt{1-\delta_1}}{2} & 0\\
      \frac{\sqrt{\Delta}\xi}{\sqrt{|\FreqAngular|^2 + \Delta\xi^2}}\frac{\sqrt{1-\delta_1}}{2} & -\frac{\FreqTheta}{\sqrt{|\FreqAngular|^2 + \Delta\xi^2}}\frac{\sqrt{1-\delta_1}}{2} & 0\\
      0 & 0 & \frac{\sqrt{Ca-\beta_2^2 + \beta_1^2}(\sigma-\sigma_2)}{\sqrt{2}(\sigma_1-\sigma_2)}\\
      0 & 0 & \frac{\sqrt{Ca-\beta_1^2 + \beta_2^2}(\sigma-\sigma_1)}{\sqrt{2}(\sigma_1-\sigma_2)}
    \end{pmatrix},
  \end{equation*}
  which is of rank 3 everywhere. 
\end{proof}

\subsubsection{Lower-order bulk terms}
\label{sec:Morawetz:lower-order-bulk}

We now move on to controlling the lower-order bulk terms that arise in
Morawetz. These arise from a few different areas. First, lower-order
bulk terms will arise from the fact that since we are now dealing with
$\WaveOpHork{2}$, the wave operator acting on sections of
$\realHorkTensor{2}$, instead of the scalar wave operator
$\ScalarWaveOp$, we have lower-order bulk terms arising from the
commutation of covariant derivatives. We will also have lower-order
terms arising from the fact that we are using pseudo-differential
multipliers. Fortunately, a key benefit will be that all of the
relevant lower-order terms will be of $O(a)$.

We begin with a basic observation concerning skew-adjoint operators.
\begin{lemma}
  \label{lemma:Morawetz:lower-order:skew-adjoint}
  Let
  $P\in \MixedOpClass{1}{1}(\Manifold;\realHorkTensor{2}(\Manifold))$
  be a skew-Hermitian operator supported on
  $\Manifold_{\operatorname{trap}}(\tau_1,\tau_2)$. Then,
  \begin{equation}
    \label{eq:Morawetz:lower-order:skew-adjoint}
    \Re\bangle*{P\psi, \psi}_{L^2(\Manifold(\tau_1,\tau_2))}
    \lesssim \sup_{\tau\in[\tau_1,\tau_2]}\EnergyHorizonDeg[\psi](\tau).
  \end{equation}  
\end{lemma}
\begin{proof}
  We decompose $P = P_{-1}D_t^2 + P_0D_t + P_1$ where $P_i \in \TanOpClass{i}(\Manifold;\realHorkTensor{2}(\Manifold))$
  are individually skew-Hermitian. Then, since $P_1$ is
  skew-Hermitian, we have that
  \begin{align*}
     2\Re\bangle*{P_1\psi, \psi}_{L^2(\Manifold(\tau_1,\tau_2))}
    ={}& \bangle*{P_1\psi, \psi}_{L^2(\Manifold(\tau_1,\tau_2))}
         + \bangle*{\psi, P_1\psi}_{L^2(\Manifold(\tau_1,\tau_2))}\\
    ={}& \bangle*{P_1\psi, \psi}_{L^2(\Manifold(\tau_1,\tau_2))}
         - \bangle*{P_1\psi, \psi}_{L^2(\Manifold(\tau_1,\tau_2))}\\
    ={}&0.
  \end{align*}
  We can also compute that
  \begin{align*}
     2\Re\bangle*{P_{-1}D_t^2\psi, \psi}_{L^2(\Manifold(\tau_1,\tau_2))}
    ={}& \bangle*{P_{-1}D_t^2\psi, \psi}_{L^2(\Manifold(\tau_1,\tau_2))}
         + \bangle*{\psi, P_{-1}D_t^2\psi}_{L^2(\Manifold(\tau_1,\tau_2))}\\
    ={}& -\bangle*{\psi, P_{-1}D_t^2\psi}_{L^2(\Manifold(\tau_1,\tau_2))}
         + \bangle*{P_{-1}D_t^2\psi, \psi}_{L^2(\Manifold(\tau_1,\tau_2))}\\
    & +O\left( \norm*{D_t\psi}_{L^2(\Sigma(\tau_2)\bigcap \Manifold_{\operatorname{trap}})}^2 
      + \norm*{D_t\psi}_{L^2(\Sigma(\tau_1)\bigcap \Manifold_{\operatorname{trap}})}^2
      \right)
    \\
    &+ O\left( \norm*{\psi}_{H^1(\Sigma(\tau_2)\bigcap \Manifold_{\operatorname{trap}})}^2 
      + \norm*{\psi}_{H^1(\Sigma(\tau_1)\bigcap \Manifold_{\operatorname{trap}})}^2
      \right)\\
    ={}& O\left( \sup_{\tau\in [\tau_1,\tau_2]}\EnergyHorizonDeg[\psi](\tau)\right).
  \end{align*}
  On the other hand, we observe that,
  \begin{align*}
    & 2\Re\bangle*{P_0D_t\psi, \psi}_{L^2(\Manifold(\tau_1,\tau_2))}\\
    ={}& \bangle*{P_0D_t\psi, \psi}_{L^2(\Manifold(\tau_1,\tau_2))}
         + \bangle*{\psi, P_0D_t\psi}_{L^2(\Manifold(\tau_1,\tau_2))}\\
    ={}& -\bangle*{\psi, P_0D_t\psi}_{L^2(\Manifold(\tau_1,\tau_2))}
         + \bangle*{P_0D_t\psi, \psi}_{L^2(\Manifold(\tau_1,\tau_2))}
         +O\left( \norm*{\psi}_{L^2(\Sigma(\tau_2)\bigcap \Manifold_{\operatorname{trap}})}^2 
         + \norm*{\psi}_{L^2(\Sigma(\tau_1)\bigcap \Manifold_{\operatorname{trap}})}^2\right)
    \\
    ={}& O\left( \sup_{\tau\in [\tau_1,\tau_2]}\EnergyHorizonDeg[\psi](\tau)\right),
  \end{align*}
  where we used the fact that
  $P_0\in\TanOpClass{0}(\Manifold;\realHorkTensor{2}(\Manifold))$ and that
  $\EnergyHorizonDeg[\psi](\tau)$ controls the $L^2(\Sigma(\tau))$
  norm in $\Manifold_{\operatorname{trap}}(\tau_1,\tau_2)$.
\end{proof}

The following lemma allows us to control error
terms arising from using pseudo-differential multipliers. Recall the definition of the $H^{-1}_{\operatorname{comp}}$
norm as relevant for us:
  \begin{equation*}
    \norm*{\psi}_{H^{-1}_{\operatorname{comp}}(\Manifold(\tau_1,\tau_2))}^2
    \vcentcolon = \norm*{\WeylQ{\bangle*{\xi + \abs*{\eta}}^{-1}}\psi}^2_{L^2(\Manifold(\tau_1,\tau_2))},
  \end{equation*}
  for $\psi$ supported on $\Manifold_{\operatorname{trap}}(\tau_1,\tau_2)$.

\begin{lemma}
  \label{lemma:ILED-near:LoT-control:zero-order-mixed-term}
  We have that
  \begin{equation}
    \label{eq:ILED-near:LoT-control:zero-order-mixed-term}
    \norm*{\HorCovDeriv_{\KillT}\psi}_{H^{-1}_{\operatorname{comp}}(\Manifold(\tau_1,\tau_2))}^2
    \lesssim \norm*{\psi}_{L^2_{\operatorname{comp}}(\Manifold(\tau_1,\tau_2))}^2
    + \int_{\Manifold(\tau_1,\tau_2)}\abs*{N}^2
    + \sup_{[\tau_1,\tau_2]}\EnergyHorizonDeg[\psi](\tau).
  \end{equation}
\end{lemma}
\begin{proof}
  Consider some compactly supported self-adjoint operator
  $Q \in \TanOpClass{-1}(\Manifold;\realHorkTensor{2}(\Manifold))$. We
  use $Q^2$ as a multiplier to see that
  \begin{equation}
    \label{eq:ILED-near:LoT-control:zero-order-mixed-term:aux1}
    \begin{split}
      &2\Re\bangle*{\WaveOpHork{2}\psi - V\psi, (g^{tt})^{-1}Q^2\psi}_{L^2(\Manifold(\tau_1,\tau_2))}\\
    ={}& \norm*{Q\HorCovDeriv_t\psi}_{L^2(\Manifold(\tau_1,\tau_2))}^2
    + O\left(\norm*{Q\HorCovDeriv_t\psi}_{L^2(\Manifold(\tau_1,\tau_2))}\norm*{\psi}_{L^2_{\operatorname{comp}}(\Manifold(\tau_1,\tau_2))}\right)\\
      & + O\left(\norm*{\psi}^2_{L^2_{\operatorname{comp}}(\Manifold(\tau_1,\tau_2))}
      + \sup_{[\tau_1,\tau_2]}\EnergyHorizonDeg[\psi](\tau)
    \right).
    \end{split}    
  \end{equation}
  Then, using the fact that $\psi$ is a solution to 
  \zcref[noname]{eq:model-problem-gRW},  we have that
  \begin{equation}
    \label{eq:ILED-near:LoT-control:zero-order-mixed-term:aux2}
    \begin{split}
      & 2\Re\bangle*{\WaveOpHork{2}\psi - V\psi, \frac{1}{\Metric_{M,a,\Lambda}^{-1}(dt, dt)}Q^2\psi}_{L^2(\Manifold(\tau_1,\tau_2))}\\
      ={}& 2\Re\bangle*{-\frac{4a\cos\theta}{\abs*{q}^2}\LeftDual{\nabla}_{\KillT}\psi, \frac{1}{\Metric_{M,a,\Lambda}^{-1}(dt, dt)}Q^2\psi}_{L^2(\Manifold(\tau_1,\tau_2))}\\
         &  + 2\Re\bangle*{N, \frac{1}{\Metric_{M,a,\Lambda}^{-1}(dt, dt)}Q^2\psi}_{L^2(\Manifold(\tau_1,\tau_2))}\\
      ={}& O \left( \norm*{\HorCovDeriv_{\KillT}\psi}_{H^{-1}_{\operatorname{comp}}(\Manifold(\tau_1,\tau_2))}\norm*{\psi}_{L^2_{\operatorname{comp}}(\Manifold(\tau_1,\tau_2))}
           + \int_{\Manifold(\tau_1,\tau_2)}\abs*{N}^2
           + \norm*{\psi}_{L^2_{\operatorname{comp}}(\Manifold(\tau_1,\tau_2)}^2 \right).
    \end{split}
  \end{equation}  
  Combining
  \zcref[noname]{eq:ILED-near:LoT-control:zero-order-mixed-term:aux1} and
  \zcref[noname]{eq:ILED-near:LoT-control:zero-order-mixed-term:aux2}, we have
  that
  \begin{equation*}
    \norm*{Q\HorCovDeriv_t\psi}_{L^2(\Manifold(\tau_1,\tau_2))}^2
    \lesssim \norm*{\psi}_{L^2_{\operatorname{comp}}(\Manifold(\tau_1,\tau_2))}^2
    + \varepsilon \norm*{\HorCovDeriv_{\KillT}\psi}_{H^{-1}_{\operatorname{comp}}(\Manifold(\tau_1,\tau_2))}^2
    + \int_{\Manifold(\tau_1,\tau_2)}\abs*{N}^2
    + \sup_{[\tau_1,\tau_2]}\EnergyHorizonDeg[\psi](\tau),
  \end{equation*}
  which allows us to close since $Q$ is arbitrary and we can absorb
  the
  $\varepsilon
  \norm*{\HorCovDeriv_{\KillT}\psi}_{H^{-1}_{\operatorname{comp}}(\Manifold(\tau_1,\tau_2))}^2$
  term from the right-hand side onto the left-hand side for
  $\varepsilon$ sufficiently small.
\end{proof}

We are now ready to control the main subprincipal error term. 
\begin{lemma}
  \label{lemma:ILED-near:LoT-control:exchange-degeneracy-trick}
  Let
  $\SubPOp\in
  \MixedOpClass{1}{1}(\Manifold;\realHorkTensor{2}(\Manifold))$ be a
  first order symmetric pseudo-differential operator supported on
  $\Manifold_{\operatorname{trap}}(\tau_1,\tau_2)$ such that
  \begin{equation*}
    \SubPOp = \SubPOp_1 + \frac{1}{2}\left( \SubPOp_0D_t + D_t\SubPOp_0 \right), \qquad
    \SubPOp_i\in \TanOpClass{i}(\Manifold;\realHorkTensor{2}(\Manifold)).
  \end{equation*}
  Then 
  \begin{equation}
    \label{eq:ILED-near:LoT-control:exchange-degeneracy-trick}
    \Re\bangle*{\SubPOp \psi, \psi}_{L^2(\Manifold_{\operatorname{trap}}(\tau_1,\tau_2))}
    \lesssim \MorNormTrap[\psi](\tau_1,\tau_2)
    + \sup_{\tau\in[\tau_1,\tau_2]}\EnergyHorizonDeg[\psi](\tau).
  \end{equation}
\end{lemma}
Before proving
\zcref{lemma:ILED-near:LoT-control:exchange-degeneracy-trick}, we
first summarize briefly the main idea which will be used several times
throughout the proof. Recall from \zcref{lemma:ILED-KdS:Bulk} that
symbolically, we control
\begin{equation*}
  (r-\rTrapping(\sigma_1))\left( \sigma-\sigma_2 \right)
  + (r-\rTrapping(\sigma_2))\left( \sigma-\sigma_1 \right)
  + \xi^2,
\end{equation*}
where we recall that $\xi$ is the radial frequency.  The main game in
the proof of
\zcref{lemma:ILED-near:LoT-control:exchange-degeneracy-trick} is to
use the following integration by parts trick
\begin{equation*}
  \bangle*{\partial_r(r-\rTrapping(\sigma_2))(\sigma-\sigma_1)\psi, \psi}
  = \bangle*{(r-\rTrapping(\sigma_2))(\sigma-\sigma_1)\psi, \partial_r\psi}
  + l.o.t,
\end{equation*}
and the decomposition
\begin{equation*}
  \sigma = -\frac{\sigma_2(\sigma-\sigma_1)}{\sigma_1-\sigma_2}
  + \frac{\sigma_1(\sigma-\sigma_2)}{\sigma_1-\sigma_2}, \qquad
  1 = \frac{\sigma-\sigma_2}{\sigma_1-\sigma_2} - \frac{\sigma -\sigma_1}{\sigma_1-\sigma_2}
\end{equation*}
to force out the appearance of $\sigma-\sigma_2$ and
$\sigma-\sigma_1$.

\begin{proof}[Proof of Lemma \ref{lemma:ILED-near:LoT-control:exchange-degeneracy-trick}]
  Observe that for $a$ sufficiently
  small with respect to $M$, we have that
  $\sigma_1-\sigma_2$ is an elliptic operator.  We now define
  \begin{equation*}
    \rTrapping_i \vcentcolon= \rTrapping_{M,a,\Lambda}(\sigma_i, \FreqPhi),
  \end{equation*}
  and  write
  \begin{equation*}
    \SubPSym = \SubPSym_0\sigma + \SubPSym_1,\qquad
    \SubPSym_i\in \TanSymClass{i}(T^{*}\mathcal{M}; \pi^{*}\End(\realHorkTensor{2}(\mathcal{M}))), 
  \end{equation*}
  where $\SubPSym$ is the symbol of $\SubPOp$.  We will handle the two
  terms separately. We first handle $\SubPSym_0\sigma$. For this term,
  we define
  \begin{equation*}
    \SubPSym_0^{(2)} \vcentcolon= -\frac{\sigma_2(\sigma-\sigma_1)}{\sigma_1-\sigma_2}\SubPSym_0,\qquad
    \SubPSym_0^{(1)} \vcentcolon= \frac{\sigma_1(\sigma-\sigma_2)}{\sigma_1-\sigma_2}\SubPSym_0,
  \end{equation*}
  so that $\SubPSym_0\sigma = \SubPSym_0^{(1)} + \SubPSym_0^{(2)}$. We now compute
  \begin{align*}
    \WeylQ{\SubPSym_0^{(i)}}
    ={}& \frac{1}{2}\left(
         \WeylQ{\SubPSym_0^{(i)}}(\partial_r(r-\rTrapping_i))
         + (\partial_r(r-\rTrapping_i))\WeylQ{\SubPSym_0^{(i)}}
         \right) \mod \MixedOpClass{0}{1}(\Manifold;\realHorkTensor{2}(\Manifold)) \notag \\
    ={}& - \frac{1}{2}\left(
         \WeylQ{\SubPSym_0^{(i)}}(r-\rTrapping_i)\partial_r
         + (r-\rTrapping_i)\partial_r\WeylQ{\SubPSym_0^{(i)}}
         \right)\\
       & + \WeylQ{\ImagUnit\SubPSym_0^{(i)}\left(r-\rTrapping_i\right)\xi}
         \mod \MixedOpClass{0}{1}(\Manifold;\realHorkTensor{2}(\Manifold)) \notag \\
    ={}&\WeylQ{\ImagUnit\SubPSym_0^{(i)}\left(r-\rTrapping_i\right)\xi}
         - \frac{1}{2}\left(
         \WeylQ{\SubPSym_0^{(i)}}(r-\rTrapping_i)\partial_r
         + (r-\rTrapping_i)\WeylQ{\SubPSym_0^{(i)}}\partial_r
         \right)\notag \\
       &  -\frac{1}{2} (r-\rTrapping_i)\left[\partial_r,\WeylQ{\SubPSym_0^{(i)}}\right]
         \mod \MixedOpClass{0}{1}\notag(\Manifold;\realHorkTensor{2}(\Manifold)) \\
    ={}&\WeylQ{\ImagUnit\SubPSym_0^{(i)}\left(r-\rTrapping_i\right)\xi}
         - \WeylQ{\SubPSym_0^{(i)}\rTrapping_i}\partial_r\\
       &- \frac{1}{2}(r-\rTrapping_i)\left[\partial_r,\WeylQ{\SubPSym_0^{(i)}}\right]
         \mod \MixedOpClass{0}{1}(\Manifold;\realHorkTensor{2}(\Manifold)).
  \end{align*}
  Now, integrating by parts, we have that
  \begin{equation}
    \label{eq:ILED-near:LoT-control:exchange-degeneracy-trick:IbP:S0}
    \begin{split}
      &\bangle*{\SubPOp_0D_t \psi, \psi}_{L^2(\Manifold_{\operatorname{trap}}(\tau_1,\tau_2))}\\
      ={}& \sum_{i=1,2}\bangle*{\WeylQ{\ImagUnit(r-\rTrapping_i)\SubPSym_0^{(i)}\xi} \psi, \psi}_{L^2(\Manifold_{\operatorname{trap}}(\tau_1,\tau_2))}\\
      & - \sum_{i=1,2}\bangle*{\WeylQ{(r-\rTrapping_i)\SubPSym_0^{(i)}}\partial_r\psi,  \psi}_{L^2(\Manifold_{\operatorname{trap}}(\tau_1,\tau_2))} \\
      &- \frac{1}{2}\sum_{i=1,2} \bangle*{(r-\rTrapping_i)\squareBrace*{\partial_r, \WeylQ{\SubPSym_0^{(i)}}}\psi, \psi}_{L^2(\Manifold_{\operatorname{trap}}(\tau_1,\tau_2))}
        + O(\MorNormTrap[\psi](\tau_1,\tau_2)), 
    \end{split}
  \end{equation}
  where the fact that the lower-order error terms are
  $O(\MorNormTrap[\psi](\tau_1,\tau_2))$ is given by
  \zcref[cap]{lemma:ILED-near:LoT-control:zero-order-mixed-term}. We treat each of the terms on the right-hand side of
  \zcref[noname]{eq:ILED-near:LoT-control:exchange-degeneracy-trick:IbP:S0}
  individually. First, observe that since
  $\ImagUnit\SubPSym_0^{(i)}\left(r-\rTrapping_i\right)\xi$ is a
  purely imaginary symbol, we have that its quantization
  $\WeylQ{\ImagUnit\SubPSym_0^{(i)}\left(r-\rTrapping_i\right)\xi}
  \in\MixedOpClass{2}{1}(\Manifold;\realHorkTensor{2}(\Manifold))$ is
  an antisymmetric operator. As a result, we can directly use 
  \zcref[cap]{lemma:Morawetz:lower-order:skew-adjoint} to control it,
  \begin{equation*}
    \sum_{i=1,2}\bangle*{\WeylQ{\ImagUnit(r-\rTrapping_i)\SubPSym_0^{(i)}\xi} \psi, \psi}_{L^2(\Manifold_{\operatorname{trap}}(\tau_1,\tau_2))}
    \lesssim \MorNormTrap[\psi](\tau_1,\tau_2) + \sup_{\tau\in[\tau_1,\tau_2]}\EnergyHorizonDeg[\psi](\tau).
  \end{equation*}
  Then since $\SubPSym_0^{{(i)}}$ are real symbols,
  $\WeylQ{(r-\rTrapping_i)\left(\SubPSym_0^{(i)})\right)}$ are
  symmetric operators. As a result, using 
  \zcref[cap]{lemma:ILED-near:LoT-control:zero-order-mixed-term} and
  Cauchy-Schwarz,
  \begin{align*}
   & \abs*{ \Re\bangle*{\WeylQ{(r-\rTrapping_i)\SubPSym_0^{(i)}} \partial_r \psi, \psi}_{L^2(\Manifold_{\operatorname{trap}}(\tau_1,\tau_2))}}\\
    \le {}& \abs*{ \Re\bangle*{\partial_r \psi, \WeylQ{(r-\rTrapping_i)\SubPSym_0^{(i)}} \psi}_{L^2(\Manifold_{\operatorname{trap}}(\tau_1,\tau_2))}}
            + \sup_{\tau\in[\tau_1,\tau_2]}\EnergyHorizonDeg[\psi](\tau).
  \end{align*}
  Then, using \zcref[cap]{lemma:ILED-KdS:Bulk,thm:fefferman-phong-tataru}, we have that
  \begin{equation*}
    \abs*{ \Re\bangle*{\partial_r \psi, \WeylQ{(r-\rTrapping_i)\SubPSym_0^{(i)}} \psi}_{L^2(\Manifold_{\operatorname{trap}}(\tau_1,\tau_2))}}
    \lesssim\MorNormTrap[\psi](\tau_1,\tau_2).
  \end{equation*}  
  Similarly, we have that using \zcref[cap]{lemma:ILED-KdS:Bulk,thm:fefferman-phong-tataru}
  \begin{equation*}
    \norm*{\WeylQ{(r-\rTrapping_i)\SubPSym_0^{(i)}\xi}\psi}_{L^2(\Manifold_{\operatorname{trap}}(\tau_1,\tau_2))}^2
    \lesssim \MorNormTrap[\psi](\tau_1,\tau_2).
  \end{equation*}
  To handle the commutator term in
  \zcref[noname]{eq:ILED-near:LoT-control:exchange-degeneracy-trick:IbP:S0},
  we first observe that
  \begin{align*}
    \PoissonB*{\xi, \SubPSym_0^{(2)}} &= -(\sigma-\sigma_1)\partial_r\left(\frac{\sigma_2\SubPSym_0}{\sigma_1-\sigma_2}\right)
                                        +\frac{ \SubPSym_0\sigma_2\partial_r\sigma_1}{\sigma_1-\sigma_2},&
    \PoissonB*{\xi, \SubPSym_0^{(1)}} &= (\sigma-\sigma_2)\partial_r\left(\frac{\sigma_1\SubPSym_0}{\sigma_1-\sigma_2}\right)
                                        -\frac{ \SubPSym_0\sigma_1\partial_r\sigma_2}{\sigma_1-\sigma_2},
  \end{align*}
  Thus, we can decompose
  \begin{equation}
    \label{eq:ILED-near:LoT-control:exchange-degeneracy-trick:IbP:S0:commutator:decomposition}
    \begin{split}
    &\sum_{i=1,2} \bangle*{(r-\rTrapping_i)\squareBrace*{\partial_r, \WeylQ{\SubPSym_0^{(i)}}}\psi, \psi}_{L^2(\Manifold_{\operatorname{trap}}(\tau_1,\tau_2))}\\
    ={}& \bangle*{(r-\rTrapping_1)\WeylQ{
      (\sigma-\sigma_2)\partial_r\left(\frac{\sigma_1\SubPSym_0}{\sigma_1-\sigma_2}\right)
      }\psi,\psi}_{L^2(\Manifold_{\operatorname{trap}}(\tau_1,\tau_2))} \\
    &  - \bangle*{(r-\rTrapping_1)\WeylQ{
      \frac{ \SubPSym_0\sigma_1\partial_r\sigma_2}{\sigma_1-\sigma_2}
      }\psi,\psi}_{L^2(\Manifold_{\operatorname{trap}}(\tau_1,\tau_2))}\\
    &- \bangle*{(r-\rTrapping_2)\WeylQ{
      (\sigma-\sigma_1)\partial_r\left(\frac{\sigma_2\SubPSym_0}{\sigma_1-\sigma_2}\right)
      }\psi,\psi}_{L^2(\Manifold_{\operatorname{trap}}(\tau_1,\tau_2))} \\
    &  + \bangle*{(r-\rTrapping_2)\WeylQ{
      \frac{ \SubPSym_0\sigma_2\partial_r\sigma_1}{\sigma_1-\sigma_2}
      }\psi,\psi}_{L^2(\Manifold_{\operatorname{trap}}(\tau_1,\tau_2))}\\
    & +\bangle*{\MixedOpClass{0}{1}(\Manifold;\realHorkTensor{2}(\Manifold))\psi,\psi}_{L^2(\Manifold_{\operatorname{trap}}(\tau_1,\tau_2))}. 
  \end{split}
  \end{equation}
We can immediately bound the error term by using
  \zcref[cap]{lemma:ILED-near:LoT-control:zero-order-mixed-term} so that
  \begin{equation*}
    \bangle*{\MixedOpClass{0}{1}(\Manifold;\realHorkTensor{2}(\Manifold))\psi,\psi}_{L^2(\Manifold_{\operatorname{trap}}(\tau_1,\tau_2))}
    \lesssim \MorNormTrap[\psi](\tau_1,\tau_2).
  \end{equation*}
  We then see that the first and third line of the right-hand side of
  \zcref[noname]{eq:ILED-near:LoT-control:exchange-degeneracy-trick:IbP:S0:commutator:decomposition}
  are controlled by the Morawetz norm. Thus, using Cauchy-Schwarz and
  \zcref[cap]{lemma:ILED-near:LoT-control:zero-order-mixed-term}, we
  have that
  \begin{align*}
    \MorNormTrap[\psi](\tau_1,\tau_2)\gtrsim{}
    &\abs*{\bangle*{(r-\rTrapping_1)\WeylQ{
      (\sigma-\sigma_2)\partial_r\left(\frac{\sigma_1\SubPSym_0}{\sigma_1-\sigma_2}\right)
      }\psi,\psi}_{L^2(\Manifold_{\operatorname{trap}}(\tau_1,\tau_2))}}\\
    &+ \abs*{\bangle*{(r-\rTrapping_2)\WeylQ{
      (\sigma-\sigma_1)\partial_r\left(\frac{\sigma_2\SubPSym_0}{\sigma_1-\sigma_2}\right)
      }\psi,\psi}_{L^2(\Manifold_{\operatorname{trap}}(\tau_1,\tau_2))}}.
  \end{align*}
  It remains to handle the second and fourth line of the right-hand
  side of
  \zcref[noname]{eq:ILED-near:LoT-control:exchange-degeneracy-trick:IbP:S0:commutator:decomposition}. Without
  loss of generality, we handle
  $(\sigma_1-\sigma_2)^{-1}{\SubPSym_0\sigma_2\partial_r\sigma_1}$ since
  $(\sigma_1-\sigma_2)^{-1}{\SubPSym_0\sigma_1\partial_r\sigma_2}$ is
  handled identically. To this end, define
  \begin{equation*}
    \tilde{\SubPSym}_{0,2} \vcentcolon= \frac{\SubPSym_0\sigma_2\partial_r\sigma_1}{\sigma_1-\sigma_2}\in  \TanSymClass{1}(T^{*}\mathcal{M};\pi^{*}\End(\realHorkTensor{2}(\mathcal{M})),\qquad
    \tilde{\SubPSym}_{0,2}^{(2)} \vcentcolon= -\frac{\sigma-\sigma_1}{\sigma_1-\sigma_2} \tilde{\SubPSym}_{0,2},\qquad
    \tilde{\SubPSym}_{0,2}^{(1)} \vcentcolon= \frac{\sigma-\sigma_2}{\sigma_1-\sigma_2} \tilde{\SubPSym}_{0,2},
  \end{equation*}
  so that $\sum_{i=1,2}\tilde{\SubPSym}_{0,2}^{(i)} = \tilde{\SubPSym}_{0,2}$. Then we have that
  \begin{equation*}
    \begin{split}
      &\bangle*{(r-\rTrapping_2)\WeylQ{\tilde{\SubPSym}_{0,2}}\psi,\psi}_{L^2(\Manifold_{\operatorname{trap}}(\tau_1,\tau_2))}\\
    ={}& \sum_{i=1,2}\bangle*{\WeylQ{(r-\rTrapping_2)\tilde{\SubPSym}_{0,2}^{(i)}}\psi ,\psi }_{L^2(\Manifold_{\operatorname{trap}}(\tau_1,\tau_2))}
    \mod \MorNormTrap[\psi](\tau_1,\tau_2)
    .
    \end{split}    
  \end{equation*}
  Then observe that by construction,
  \begin{equation*}
    \norm*{\WeylQ{(r-\rTrapping_2)\tilde{\SubPSym}_{0,2}^{(2)}}\psi}_{L^2(\Manifold_{\operatorname{trap}}(\tau_1,\tau_2))}
    \lesssim \MorNormTrap[\psi](\tau_1,\tau_2).
  \end{equation*}
  Thus, it suffices to control
  $\bangle*{\WeylQ{(r-\rTrapping_2)\tilde{\SubPSym}_{0,2}^{(1)}}\psi,\psi}_{L^2(\Manifold_{\operatorname{trap}}(\tau_1,\tau_2))}$.   
  To this end, we use integration by parts to write that up to error terms controlled by $\MorNormTrap[\psi](\tau_1,\tau_2)$,
  \begin{align}
    & \bangle*{\WeylQ{(r-\rTrapping_2)\tilde{\SubPSym}_{0,2}^{(1)}}\psi,\psi}_{L^2(\Manifold_{\operatorname{trap}}(\tau_1,\tau_2))}\notag\\
    ={}& \bangle*{\WeylQ{\ImagUnit(r-\rTrapping_1)(r-\rTrapping_2)\tilde{\SubPSym}_{0,2}^{(1)}\xi} \psi,\psi}_{L^2(\Manifold_{\operatorname{trap}}(\tau_1,\tau_2))} \notag\\
       &- \bangle*{\WeylQ{(r-\rTrapping_1)(r-\rTrapping_2)\tilde{\SubPSym}_{0,2}^{(1)}}\partial_r \psi,\psi}_{L^2(\Manifold_{\operatorname{trap}}(\tau_1,\tau_2))}\notag\\
    &-\frac{1}{2} \bangle*{\WeylQ{(r-\rTrapping_1)(r-\rTrapping_2)}\squareBrace*{\partial_r,\WeylQ{\tilde{\SubPSym}_{0,2}^{(1)}}}\psi,\psi}_{L^2(\Manifold_{\operatorname{trap}}(\tau_1,\tau_2))}
      .
      \label{eq:ILED-near:LoT-control:exchange-degeneracy-trick:IbP:S0:rep1}
  \end{align}
  Since
  $\ImagUnit(r-\rTrapping_1)(r-\rTrapping_2)\tilde{\SubPSym}_{0,2}^{(1)}\xi$
  is a purely imaginary symbol, we have that its quantization
  $\left(\ImagUnit(r-\rTrapping_1)(r-\rTrapping_2)\tilde{\SubPSym}_{0,2}^{(1)}\xi\right)^w$
  is an antisymmetric operator. As a result, using 
  \zcref[cap]{lemma:Morawetz:lower-order:skew-adjoint}, we have that
  \begin{equation*}
    \abs*{\bangle*{\WeylQ{\ImagUnit(r-\rTrapping_1)(r-\rTrapping_2)\tilde{\SubPSym}_{0,2}^{(1)}\xi}
        \psi,\psi}_{L^2(\Manifold_{\operatorname{trap}}(\tau_1,\tau_2))}}
    \lesssim \sup_{\tau\in[\tau_1,\tau_2]}\EnergyHorizonDeg[\psi](\tau)
    .
  \end{equation*}
  On the other hand, since
  $(r-\rTrapping_1)(r-\rTrapping_2)\tilde{\SubPSym}_{0,2}^{(1)}$ is a real symbol, we have that
  \begin{align*}
    & \bangle*{\WeylQ{(r-\rTrapping_1)(r-\rTrapping_2)\tilde{\SubPSym}_{0,2}^{(1)}}\partial_r \psi,\psi}_{L^2(\Manifold_{\operatorname{trap}}(\tau_1,\tau_2))}\\
    ={}& \bangle*{\partial_r \psi,\WeylQ{(r-\rTrapping_1)(r-\rTrapping_2)\tilde{\SubPSym}_{0,2}^{(1)}}\psi}_{L^2(\Manifold_{\operatorname{trap}}(\tau_1,\tau_2))}
         + O\left(\sup_{\tau\in[\tau_1,\tau_2]}\EnergyHorizonDeg[\psi](\tau)\right)\\
    \lesssim{}& \MorNormTrap[\psi](\tau_1,\tau_2) + \sup_{\tau\in[\tau_1,\tau_2]}\EnergyHorizonDeg[\psi](\tau).
  \end{align*}
  Observing again that
  $(r-\rTrapping_1)(r-\rTrapping_2)\Op\left(\tilde{\SubPSym}_{0,2}^{(1)}\right)$
  is Hermitian up to a
  $\MixedOpClass{0}{1}(\Manifold;\realHorkTensor{2}(\Manifold))$ %
  term, we control the first three terms on the right-hand side of
  \zcref[noname]{eq:ILED-near:LoT-control:exchange-degeneracy-trick:IbP:S0:rep1}
  by a combination of integration by parts, Cauchy-Schwarz, and 
  \zcref[cap]{lemma:ILED-near:LoT-control:zero-order-mixed-term} after taking
  $\varepsilon_{\SubPOp}$ and $\delta_{\operatorname{trap}}$
  sufficiently small. To handle the final commutator term, we again
  observe that
  \begin{equation*}
    \PoissonB*{\xi, \tilde{\SubPSym}_{0,2}^{(1)}}
    = (\sigma-\sigma_2)\partial_r\left(\frac{\tilde{\SubPSym}_{0,2}}{\sigma_1-\sigma_2}\right)
    - \tilde{\SubPSym}_{0,2}\frac{\partial_r\sigma_2}{\sigma_1-\sigma_2}.
  \end{equation*}
  Again we can thus decompose
  \begin{align}
    \label{eq:ILED-near:LoT-control:exchange-degeneracy-trick:IbP:S0:commutator2:decomposition}
    \begin{split}
      &\bangle*{(r-\rTrapping_1)(r-\rTrapping_2)\squareBrace*{\partial_r,\WeylQ{\tilde{\SubPSym}_{0,2}^{(1)}}}\psi,\psi}_{L^2(\Manifold_{\operatorname{trap}}(\tau_1,\tau_2))}\\
      ={}& \bangle*{(r-\rTrapping_1)(r-\rTrapping_2)\WeylQ{ (\sigma-\sigma_2)\partial_r\left(\frac{\tilde{\SubPSym}_{0,2}}{\sigma_1-\sigma_2}\right)}\psi,\psi}_{L^2(\Manifold_{\operatorname{trap}}(\tau_1,\tau_2))}\\         
      &- \bangle*{(r-\rTrapping_1)(r-\rTrapping_2)\WeylQ{\tilde{\SubPSym}_{0,2}\frac{\partial_r\sigma_2}{\sigma_1-\sigma_2}}\psi,\psi}_{L^2(\Manifold_{\operatorname{trap}}(\tau_1,\tau_2))}
        + O(\MorNormTrap[\psi](\tau_1,\tau_2))
        .
    \end{split}
  \end{align}
  To control the first term on the right-hand side of
  \zcref[noname]{eq:ILED-near:LoT-control:exchange-degeneracy-trick:IbP:S0:commutator2:decomposition},
  we again have from the Weyl calculus that
  \begin{align*}
    &\bangle*{(r-\rTrapping_1)(r-\rTrapping_2)\WeylQ{ (\sigma-\sigma_2)\partial_r\left(\frac{\tilde{\SubPSym}_{0,2}}{\sigma_1-\sigma_2}\right)}\psi,\psi}_{L^2(\Manifold_{\operatorname{trap}}(\tau_1,\tau_2))}\\
    ={}& \bangle*{\WeylQ{(r-\rTrapping_1)(r-\rTrapping_2) (\sigma-\sigma_2)\partial_r\left(\frac{\tilde{\SubPSym}_{0,2}}{\sigma_1-\sigma_2}\right)}\psi,\psi }_{L^2(\Manifold_{\operatorname{trap}}(\tau_1,\tau_2))}
         \mod \MorNormTrap[\psi](\tau_1,\tau_2). 
  \end{align*}
  Thus, using \zcref[cap]{thm:fefferman-phong-tataru} and
  Cauchy-Schwarz to contain the main term on the \RHS{} and \zcref[cap]{lemma:ILED-near:LoT-control:zero-order-mixed-term} to control
  the error term, we have immediately that
  \begin{align*}
    &\abs*{  \bangle*{(r-\rTrapping_1)(r-\rTrapping_2)\WeylQ{ (\sigma-\sigma_2)\partial_r\left(\frac{\tilde{\SubPSym}_{0,2}}{\sigma_1-\sigma_2}\right)}\psi,\psi}_{L^2(\Manifold_{\operatorname{trap}}(\tau_1,\tau_2))}}\\
    \lesssim{}& \MorNormTrap[\psi](\tau_1,\tau_2)
                +  \sup_{\tau\in[\tau_1,\tau_2]}\EnergyHorizonDeg[\psi](\tau).
  \end{align*}
  We now control the second term on the right-hand side of
  \zcref[noname]{eq:ILED-near:LoT-control:exchange-degeneracy-trick:IbP:S0:commutator2:decomposition}.
  We partition one final time to define
  \begin{gather*}
      \tilde{\SubPSym}_{0,2,1}\vcentcolon= \tilde{\SubPSym}_{0,2}\frac{\partial_r\sigma_2}{\sigma_1-\sigma_2}\in \TanSymClass{1}(T^{*}\mathcal{M};\pi^{*}\End(\realHorkTensor{2}(\mathcal{M}))),\\
    \tilde{\SubPSym}_{0,2,1}^{(2)}\vcentcolon= -\frac{\sigma-\sigma_1}{\sigma_1-\sigma_2}\tilde{\SubPSym}_{0,2,1},\qquad
    \tilde{\SubPSym}_{0,2,1}^{(1)}\vcentcolon= \frac{\sigma-\sigma_2}{\sigma_1-\sigma_2}\tilde{\SubPSym}_{0,2,1},
  \end{gather*}
  so that $\tilde{\SubPSym}_{0,2,1} = \tilde{\SubPSym}_{0,2,1}^{(2)} + \tilde{\SubPSym}_{0,2,1}^{(1)}$. Then we have that
  \begin{align*}
    &\bangle*{(r-\rTrapping_2)(r-\rTrapping_1)\WeylQ{\tilde{\SubPSym}_{0,2,1}}\psi,\psi}_{L^2(\Manifold_{\operatorname{trap}}(\tau_1,\tau_2))}\\
    ={}& \sum_{i=1,2}\bangle*{\WeylQ{(r-\rTrapping_2)(r-\rTrapping_1)\tilde{\SubPSym}_{0,2,1}^{(i)}}\psi,\psi}_{L^2(\Manifold_{\operatorname{trap}}(\tau_1,\tau_2))}\\
    & + O\left(\MorNormTrap[\psi](\tau_1,\tau_2) + \sup_{\tau\in[\tau_1,\tau_2]}\EnergyHorizonDeg[\psi](\tau)\right),
  \end{align*}
  where again the error term being
  $O\left(\MorNormTrap[\psi](\tau_1,\tau_2) +
  \sup_{\tau\in[\tau_1,\tau_2]}\EnergyHorizonDeg[\psi](\tau)\right)$ is
  given by \zcref[cap]{lemma:ILED-near:LoT-control:zero-order-mixed-term}. But by construction, 
  \begin{align*}
    &\norm*{\WeylQ{(r-\rTrapping_2)(r-\rTrapping_1)\tilde{\SubPSym}_{0,2,1}^{(1)}}\psi}^2_{L^2(\Manifold_{\operatorname{trap}}(\tau_1,\tau_2))}\\
    &+ \norm*{\WeylQ{(r-\rTrapping_2)(r-\rTrapping_1)\tilde{\SubPSym}_{0,2,1}^{(2)}}\psi}^2_{L^2(\Manifold_{\operatorname{trap}}(\tau_1,\tau_2))}\\
    \lesssim{}& \MorNormTrap[\psi](\tau_1,\tau_2).
  \end{align*}
  As a result, using Cauchy-Schwarz and 
  \zcref[cap]{lemma:ILED-near:LoT-control:zero-order-mixed-term}, we have that
  \begin{equation*}
    \begin{split}
      &\sum_{i=1,2}\abs*{\bangle*{(r-\rTrapping_2)(r-\rTrapping_1)\WeylQ{\tilde{\SubPSym}_{0,2,1}^{(i)}}\psi,\psi}_{L^2(\Manifold_{\operatorname{trap}}(\tau_1,\tau_2))}}\\
        \lesssim{}& \MorNormTrap[\psi](\tau_1,\tau_2) +\sup_{\tau\in[\tau_1,\tau_2]}\EnergyHorizonDeg[\psi](\tau). 
    \end{split}    
  \end{equation*}    
  We now move on to handling
  $\Re\bangle*{\SubPOp_1 \psi,\psi}_{L^2(\Manifold_{\operatorname{trap}}(\tau_1,\tau_2))}$. The
  main idea will be the same as when handling
  $\SubPOp_0\partial_t$. To this end, we define
  \begin{equation*}
    \SubPSym_1^{(2)}\vcentcolon= -\frac{\sigma-\sigma_1}{\sigma_1-\sigma_2}\SubPSym,\qquad
    \SubPSym_1^{(1)}\vcentcolon= \frac{\sigma-\sigma_2}{\sigma_1-\sigma_2}\SubPSym,
  \end{equation*}
  so that $\SubPSym = \SubPSym_1^{(1)}+\SubPSym_1^{(2)}$. Again using the Weyl calculus, we have that 
  \begin{align}
    \label{eq:ILED-near:LoT-control:exchange-degeneracy-trick:IbP:S1}
    \begin{split}
      \bangle*{\WeylQ{\SubPSym_1^{(i)}}\psi, \psi}_{L^2(\Manifold_{\operatorname{trap}}(\tau_1,\tau_2))}
    ={}& \bangle*{\WeylQ{\ImagUnit(r-\rTrapping_i)\SubPSym_1^{(i)}\xi} \psi, \psi}_{L^2(\Manifold_{\operatorname{trap}}(\tau_1,\tau_2))}\\
    & + \bangle*{\WeylQ{(r-\rTrapping_i)\SubPSym_{1}^{(i)}}\psi, \partial_r \psi}_{L^2(\Manifold_{\operatorname{trap}}(\tau_1,\tau_2))} \\
    &+ \bangle*{(r-\rTrapping_i)\squareBrace*{\partial_r, \WeylQ{\SubPSym_1^{(i)}}}\psi, \psi}_{L^2(\Manifold_{\operatorname{trap}}(\tau_1,\tau_2))}\\
    &+ O\left(\MorNormTrap[\psi](\tau_1,\tau_2) + \sup_{\tau\in[\tau_1,\tau_2]}\EnergyHorizonDeg[\psi](\tau)\right), 
    \end{split}
  \end{align}
  where the error terms being
  $O\left(\MorNormTrap[\psi](\tau_1,\tau_2) +
    \sup_{\tau\in[\tau_1,\tau_2]}\EnergyHorizonDeg[\psi](\tau)\right)$
  comes from 
  \zcref[cap]{lemma:ILED-near:LoT-control:zero-order-mixed-term}. We treat each of the terms on the right-hand side of
  \zcref[noname]{eq:ILED-near:LoT-control:exchange-degeneracy-trick:IbP:S1}
  individually. Observing again that
  $\WeylQ{\ImagUnit(r-\rTrapping_i)\SubPSym_1^{(i)}\xi}$ is an
  antisymmetric operator, we have from 
  \zcref[cap]{lemma:Morawetz:lower-order:skew-adjoint} that
  \begin{equation*}
    \abs*{\bangle*{\WeylQ{\ImagUnit(r-\rTrapping_i)\SubPSym_1^{(i)}\xi} \psi, \psi}_{L^2(\Manifold_{\operatorname{trap}}(\tau_1,\tau_2))}}
    \lesssim \MorNormTrap[\psi](\tau_1,\tau_2) + \sup_{\tau\in[\tau_1,\tau_2]}\EnergyHorizonDeg[\psi](\tau).
  \end{equation*}
  To control the second term on the \RHS, we use 
  \zcref{thm:fefferman-phong-tataru} and Cauchy-Schwarz to directly bound
  \begin{equation*}
    \abs*{\bangle*{\WeylQ{(r-\rTrapping_i)\SubPSym_{1}^{(i)}}\psi, \partial_r \psi}_{L^2(\Manifold_{\operatorname{trap}}(\tau_1,\tau_2))}}
    \lesssim \MorNormTrap[\psi](\tau_1,\tau_2). 
  \end{equation*}
  To handle the commutator term in
  \zcref[noname]{eq:ILED-near:LoT-control:exchange-degeneracy-trick:IbP:S1},
  we proceed as previously and observe that (without loss of
  generality, we will just consider the $i=1$ case)
  \begin{align*}
    &\bangle*{(r-\rTrapping_1)\squareBrace*{\partial_r, \WeylQ{\SubPSym_1^{(1)}}}\psi, \psi}_{L^2(\Manifold_{\operatorname{trap}}(\tau_1,\tau_2))}\\
    ={}& \bangle*{(r-\rTrapping_1)\left((\sigma-\sigma_2)\partial_r\WeylQ{\frac{\SubPSym_1}{\sigma_1-\sigma_2}\right)}\psi, \psi}_{L^2(\Manifold_{\operatorname{trap}}(\tau_1,\tau_2))}\\
       &- \bangle*{(r-\rTrapping_1)\WeylQ{\frac{\SubPSym_1\partial_r\sigma_2}{\sigma_1-\sigma_2}}\psi, \psi}_{L^2(\Manifold_{\operatorname{trap}}(\tau_1,\tau_2))}
         + O\left(\MorNormTrap[\psi](\tau_1,\tau_2) + \sup_{\tau\in [\tau_1,\tau_2]}\EnergyHorizonDeg[\psi](\tau)\right),
  \end{align*}
  where the error being
  $O\left(\MorNormTrap[\psi](\tau_1,\tau_2) + \sup_{\tau\in
      [\tau_1,\tau_2]}\EnergyHorizonDeg[\psi](\tau)\right)$ is given
  by \zcref[cap]{lemma:ILED-near:LoT-control:zero-order-mixed-term}.
  To handle the second term, we can define
  \begin{equation*}
    \tilde{\SubPSym}_{1,1}\vcentcolon= \frac{\SubPSym_1\partial_r\sigma_2}{\sigma_1-\sigma_2}\in \TanSymClass{1}(T^*\Manifold;\pi^*\End(\realHorkTensor{2}(\Manifold))),\qquad
    \tilde{\SubPSym}_{1,1}^{(2)} \vcentcolon= -\frac{\sigma-\sigma_1}{\sigma_1-\sigma_2}\tilde{\SubPSym}_{1,1},\qquad
    \tilde{\SubPSym}_{1,1}^{(1)} \vcentcolon= \frac{\sigma-\sigma_2}{\sigma_1-\sigma_2}\tilde{\SubPSym}_{1,1}
  \end{equation*}
  so that $\tilde{\SubPSym}_{1,1} = \sum_{i=1,2}\tilde{\SubPSym}_{1,1}^{(i)}$ and as a result, 
  \begin{align*}
    \bangle*{(r-\rTrapping_1)\WeylQ{\tilde{\SubPSym}_{1,1}}\psi, \psi}_{L^2(\Manifold_{\operatorname{trap}}(\tau_1,\tau_2))}
    ={}& \bangle*{\WeylQ{(r-\rTrapping_1)\tilde{\SubPSym}_{1,1}^{(1)}}\psi, \psi}_{L^2(\Manifold_{\operatorname{trap}}(\tau_1,\tau_2))}\\
       &+ \bangle*{\WeylQ{(r-\rTrapping_1)\tilde{\SubPSym}_{1,1}^{(2)}}\psi, \psi}_{L^2(\Manifold_{\operatorname{trap}}(\tau_1,\tau_2))}\\
       & + O\left(\MorNormTrap[\psi](\tau_1,\tau_2) + \sup_{\tau\in [\tau_1,\tau_2]}\EnergyHorizonDeg[\psi](\tau)\right),        
  \end{align*}
  where it is clear by construction that the first term on the \RHS{}
  respects the degeneracy in the Morawetz norm and thus can be
  controlled by using Cauchy-Schwarz and
  \zcref{lemma:ILED-near:LoT-control:zero-order-mixed-term}. To control
  the second term on the right-hand side, we repeat the integration by
  parts argument above to see that
  \begin{align}
    &\bangle*{\WeylQ{(r-\rTrapping_1)\tilde{\SubPSym}_{1,1}^{(2)}}\psi,\psi}_{L^2(\Manifold_{\operatorname{trap}}(\tau_1,\tau_2))}\notag\\
    ={}& \bangle*{\WeylQ{\ImagUnit(r-\rTrapping_2)(r-\rTrapping_1)\tilde{\SubPSym}_{1,1}^{(2)}\xi} \psi,\psi}_{L^2(\Manifold_{\operatorname{trap}}(\tau_1,\tau_2))} \notag\\
       &- \bangle*{\WeylQ{(r-\rTrapping_2)(r-\rTrapping_1)\tilde{\SubPSym}_{1,1}^{(2)}}\psi,\partial_r\psi}_{L^2(\Manifold_{\operatorname{trap}}(\tau_1,\tau_2))}\notag\\
       & -\frac{1}{2} \bangle*{(r-\rTrapping_2)\squareBrace*{\partial_r,\WeylQ{(r-\rTrapping_1)\tilde{\SubPSym}_{1,1}^{(2)}}}\psi,\psi}_{L^2(\Manifold_{\operatorname{trap}}(\tau_1,\tau_2))}. \label{eq:ILED-near:LoT-control:exchange-degeneracy-trick:IbP:S1:rep1}
  \end{align}
  Again, since
  $\WeylQ{\ImagUnit(r-\rTrapping_2)(r-\rTrapping_1)\tilde{\SubPSym}_{1,1}^{(2)}\xi}$
  is anti-symmetric, we have from 
  \zcref[cap]{lemma:Morawetz:lower-order:skew-adjoint} that
  \begin{equation*}
    \abs*{\bangle*{\WeylQ{\ImagUnit(r-\rTrapping_2)(r-\rTrapping_1)\tilde{\SubPSym}_{1,1}^{(2)}\xi} \psi,\psi}_{L^2(\Manifold_{\operatorname{trap}}(\tau_1,\tau_2))}}
    \lesssim \sup_{\tau\in[\tau_1,\tau_2]}\EnergyHorizonDeg[\psi](\tau).
  \end{equation*}
  To handle the second term on the \RHS{} of
  \zcref[noname]{eq:ILED-near:LoT-control:exchange-degeneracy-trick:IbP:S1:rep1},
  we again use Cauchy-Schwarz and 
  \zcref[cap]{thm:fefferman-phong-tataru} to see that
  \begin{equation*}
    \abs*{\bangle*{\WeylQ{(r-\rTrapping_2)(r-\rTrapping_1)\tilde{\SubPSym}_{1,1}^{(2)}}\psi,\partial_r\psi}_{L^2(\Manifold_{\operatorname{trap}}(\tau_1,\tau_2))}}
    \lesssim \MorNormTrap[\psi](\tau_1,\tau_2) + \sup_{\tau\in[\tau_1,\tau_2]}\EnergyHorizonDeg[\psi](\tau). 
  \end{equation*}
  To handle the last term on the right-hand side of
  \zcref[noname]{eq:ILED-near:LoT-control:exchange-degeneracy-trick:IbP:S1:rep1},
  we see from the Weyl calculus and
  \zcref{lemma:ILED-near:LoT-control:zero-order-mixed-term} that
  \begin{align*}
    &\bangle*{(r-\rTrapping_2)(r-\rTrapping_1)\squareBrace*{\partial_r,\WeylQ{\tilde{\SubPSym}_{1,1}^{(2)}}}\psi,\psi}_{L^2(\Manifold_{\operatorname{trap}}(\tau_1,\tau_2))}\\
    = {}& -\bangle*{\WeylQ{(r-\rTrapping_2)(r-\rTrapping_1)
          (\sigma-\sigma_1)\partial_r\left(\frac{\tilde{\SubPSym}_{1,1}}{\sigma_1-\sigma_2}\right)
          }\psi,\psi}_{L^2(\Manifold_{\operatorname{trap}}(\tau_1,\tau_2))}\\
    & + \bangle*{\WeylQ{((r-\rTrapping_2)(r-\rTrapping_1)
      \frac{\tilde{\SubPSym}_{1,1}\partial_r\sigma_1}{\sigma_1-\sigma_2}
      }\psi,\psi}_{L^2(\Manifold_{\operatorname{trap}}(\tau_1,\tau_2))}\\
    &  + O\left(\MorNormTrap[\psi](\tau_1,\tau_2) + \sup_{\tau\in[\tau_1,\tau_2]}\EnergyHorizonDeg[\psi](\tau)
      \right),
  \end{align*}
  where the first term on the right-hand side is handled directly as
  before using Cauchy-Schwarz, and 
  \zcref[cap]{lemma:ILED-near:LoT-control:zero-order-mixed-term}. The second
  term on the right-hand side is handled by a final decomposition
  \begin{equation*}    
    \tilde{\SubPSym}_{1,1,2}\vcentcolon= \tilde{\SubPSym}_{1,1}\frac{\partial_r\sigma_1}{\sigma_1-\sigma_2}\in S^1,\qquad
    \tilde{\SubPSym}_{1,1,2}^{(2)}\vcentcolon= -\frac{\sigma-\sigma_1}{\sigma_1-\sigma_2}\tilde{\SubPSym}_{1,1,2},\qquad
    \tilde{\SubPSym}_{1,1,2}^{(1)}\vcentcolon= \frac{\sigma-\sigma_2}{\sigma_1-\sigma_2}\tilde{\SubPSym}_{1,1,2}.  
  \end{equation*}
  Then it is apparent that
  \begin{equation*}
    \sum_{i=1,2}\norm*{\WeylQ{(r-\rTrapping_1)(r-\rTrapping_2)\tilde{\SubPSym}_{1,1,2}^{(i)}}\psi}_{L^2(\Manifold_{\operatorname{trap}}(\tau_1,\tau_2))}
    \lesssim \MorNormTrap[\psi](\tau_1,\tau_2)
    + \sup_{\tau\in [\tau_1,\tau_2]}\EnergyHorizonDeg[\psi](\tau),
  \end{equation*}
  and we conclude the proof of 
  \zcref[cap]{lemma:ILED-near:LoT-control:exchange-degeneracy-trick} by using
  Cauchy-Schwarz and \zcref[cap]{lemma:ILED-near:LoT-control:zero-order-mixed-term}.
\end{proof}

We now turn to the application of the
\zcref[cap]{lemma:ILED-near:LoT-control:exchange-degeneracy-trick} and
the preceding lemmas in the section.

\begin{lemma}
  \label{lemma:ILED-near:KdS:extra-terms}
  Consider
  $\left(\MorawetzVF_{M,a,\Lambda},
    \MorawetzLagrangeCorr_{M,a,\Lambda},
    \MorawetzOneForm_{M,a,\Lambda}\right)$ as chosen in the proof
  of \zcref[cap]{prop:Morawetz:outside-trapping:KCurrent-coercive} and define $\widecheck{K}[\psi]\vcentcolon={} \MorawetzVF_{M,a,\Lambda}^{\mu}\HorCovDeriv^{\nu}\psi^a\HorRiem_{ab\nu\mu}\psi^b$. We have
  \begin{equation*}
    \abs*{\int_{\Manifold(\tau_1,\tau_2)}\widecheck{K}[\psi]}\lesssim O(a)\left(\MorNormTrap[\psi](\tau_1,\tau_2) + \sup_{\tau\in [\tau_1,\tau_2]}\EnergyHorizonDeg[\psi](\tau)\right).
  \end{equation*}
\end{lemma}
\begin{proof}
  We first consider the term generated by
  $\MorawetzVF_{M,a,\Lambda}^{\mu}\HorCovDeriv^{\nu}\psi^a\HorRiem_{ab\nu\mu}\psi^{b}$. Using
  the fact that $\HorRiem_{ab\nu\mu}$ is antisymmetric with respect to
  the $(a,b)$ indices, we rewrite
  \begin{equation*}
    \MorawetzVF_{M,a,\Lambda}^{\mu}\HorCovDeriv^{\nu}\psi^a\HorRiem_{ab\nu\mu}\psi^b
    = \frac{1}{2}\MorawetzVF_{M,a,\Lambda}^{\mu}\volFormHor^{ab}\HorRiem_{ab\nu\mu}\LeftDual{\psi}\cdot \HorCovDeriv^{\nu}\psi.
  \end{equation*}
  We first observe that, for $\Metric_{M,a,\Lambda}$ a \KdS{} metric,
  we have that $\HorRiem_{ab\nu r} \lesssim \frac{a}{r^4}$. Using that $\MorawetzVF_{M,a,\Lambda} = \mathcal{F}(r)\partial_r$, and denoting $A_{\nu}\vcentcolon= \frac{1}{r}\mathcal{F}(r)\volFormHor^{ab}\HorRiem_{ab\nu r}$ (which satisfies $A_{\nu} \lesssim \frac{a}{r^4}$) we have 
  \begin{equation*}
    \frac{1}{2}\MorawetzVF_{M,a,\Lambda}^{\mu}\volFormHor^{ab}\HorRiem_{ab\nu\mu}\LeftDual{\psi}\cdot \HorCovDeriv^{\nu}\psi
    = A_{\nu}\LeftDual{\psi}\cdot \HorCovDeriv^{\nu}\psi.
  \end{equation*}
  Now, let us define the differential operator $B\in \Psi^1_1$ such that
  \begin{equation*}
    \bangle*{B\psi,\psi}_{L^2(\Manifold(\tau_1,\tau_2))} = \bangle*{A_{\nu}\LeftDual{\psi},\HorCovDeriv^{\nu}\psi}_{L^2(\Manifold(\tau_1,\tau_2))},
  \end{equation*}
  that is $B:\psi \mapsto \left(\HorCovDeriv^{\nu}\right)^{*}\left(A_\nu\LeftDual{\psi}\right)$. Then notice that since $\HorCovDeriv$ is principally skew-symmetric
  and the dual operator $\psi\mapsto \LeftDual{\psi}$ is
  skew-symmetric, $B$ is a first-order, principally symmetric operator
  which we can decompose as $B = B_1 + B_0$ where $B_1$ is a first-order symmetric differential operator, and
  $B_0$ is the lower-order remainder term satisfying
  \begin{align*}
    \norm*{B_1\psi}_{L^2(\Manifold(\tau_1,\tau_2)}
    \lesssim{}& a\norm*{r^{-2}\HorCovDeriv\psi}_{L^2(\Manifold(\tau_1,\tau_2))}
                + a\norm*{r^{-2}\psi}_{L^2(\Manifold(\tau_1,\tau_2))},\\
    \norm*{B_0\psi}_{L^2(\Manifold(\tau_1,\tau_2))}
    \lesssim{}& a\norm*{r^{-2}\psi}_{L^2(\Manifold(\tau_1,\tau_2))},
  \end{align*}
  We now show how to control
  $\bangle*{B\psi,\psi}_{L^2(\Manifold(\tau_1,\tau_2))}$. The main
  idea will be to apply the previous lemmas, although we note that
  some caution must be exercised since cutoffs need to be introduced. Define the cutoffs
  \begin{equation*}
    \mathring{\chi}^2 + \breve{\chi}^2 = 1,\qquad \supp{\mathring{\chi}}\in [3M-\delta_{\operatorname{trap}}, 3M+\delta_{\operatorname{trap}}].
  \end{equation*}
  that localize to and away from the trapped set respectively. We have
  \begin{align*}
    \bangle*{B\psi,\psi}_{L^2(\Manifold(\tau_1,\tau_2))}
    ={}& \bangle*{B\psi, \left(\mathring{\chi}^2 + \breve{\chi}^2\right)\psi}_{L^2(\Manifold(\tau_1,\tau_2))}\\
    ={}& \bangle*{B\left(\mathring{\chi}\psi\right), \mathring{\chi}\psi}_{L^2(\Manifold(\tau_1,\tau_2))}
         +  \bangle*{B\left(\breve{\chi}\psi\right), \breve{\chi}\psi}_{L^2(\Manifold(\tau_1,\tau_2))}
         + O\left(a\norm*{r^{-2}\psi}^2_{L^2(\Manifold(\tau_1,\tau_2))}\right).
  \end{align*}
  On the support of $\breve{\chi}$, $\MorNorm[\cdot]$ is coercive. As a result, we have that
  \begin{align*}
    \bangle*{B\psi,\psi}_{L^2(\Manifold(\tau_1,\tau_2))}
    ={}& \bangle*{B\left(\mathring{\chi}\psi\right), \mathring{\chi}\psi}_{L^2(\Manifold(\tau_1,\tau_2))}     
         + O\left(a\MorNorm[\psi](\tau_1,\tau_2)\right),\\
    ={}& \bangle*{B_1\left(\mathring{\chi}\psi\right), \mathring{\chi}\psi}_{L^2(\Manifold(\tau_1,\tau_2))}
         + O\left(a\MorNorm[\psi](\tau_1,\tau_2)\right).
  \end{align*}
  Then, we use
  \zcref[cap]{lemma:ILED-near:LoT-control:exchange-degeneracy-trick} to
  control
  \begin{equation*}
    \abs*{\bangle*{B\left(\mathring{\chi}\psi\right), \mathring{\chi}\psi}_{L^2(\Manifold(\tau_1,\tau_2))}}
    \lesssim a\left(\MorNormTrap[\psi](\tau_1,\tau_2)
      + \sup_{\tau\in [\tau_1,\tau_2]}\EnergyHorizonDeg[\psi](\tau)\right),
  \end{equation*}
  as desired, concluding the proof of \zcref{lemma:ILED-near:KdS:extra-terms}.
\end{proof}

\subsection{Combining the trapping and non-trapping bulk terms}
\label{sec:combined-trapping-nontrapping-bulk}

As constructed in \zcref[cap]{lemma:ILED-KdS:Bulk}, both
$\widetilde{\MorawetzSym}$ and $\widetilde{\MorawetzLagrangeCorrSym}$
are homogeneous symbols and can be made smooth simply by a standard
truncation away from low frequencies. Since they are also only defined
near $r=3M$, we will also truncate the symbols in physical
space. Define
\begin{equation}
  \label{eq:Morawetz:chi-trapping-cutoff:def}
  \mathring{\chi} =
  \begin{cases}
    1 & \abs*{r-3M} < \delta_{\operatorname{trap}}\\
    0 & \abs*{r-3M} > 2\delta_{\operatorname{trap}},
  \end{cases}
\end{equation}
and define $\breve{\chi}$ so that $1 = \mathring{\chi}^2 + \breve{\chi}^2$. We truncate our pseudo-differential multipliers:
\begin{equation}
  \label{eq:Morawetz:morawetz-VF-LagrangeCorr-with-cutoff:def}
  \widetilde{\MorawetzVF} = \mathring{\chi}\WeylQ{\widetilde{\MorawetzSym}}\mathring{\chi}, \qquad
  \widetilde{\MorawetzLagrangeCorr} = \mathring{\chi}\WeylQ{\widetilde{\MorawetzLagrangeCorr}}\mathring{\chi}.
\end{equation}
From the Weyl calculus, we can only say that
there exist
$K^w_i\in \TanOpClass{i}(\Manifold;\realHorkTensor{2}(\Manifold))$
such that
\begin{equation*}
  \left(\left[\ScalarWaveOp[\Metric], \widetilde{\MorawetzVF}\right]
    + \ScalarWaveOp[\Metric]\widetilde{\MorawetzLagrangeCorr}
    + \widetilde{\MorawetzLagrangeCorr}\ScalarWaveOp[\Metric]
  \right)
  = K^w_2 + 2K_1^w\HorCovDeriv_t + K_0^w\HorCovDeriv_t^2 + K^w_{-1}\HorCovDeriv_t^3. 
\end{equation*}
In order to eliminate the $K^w_{-1}\HorCovDeriv_t^3$ term, we
slightly alter the choice of pseudo-differential Lagrangian corrector
so that
\begin{equation*}
  \widetilde{\MorawetzLagrangeCorr}
  = \mathring{\chi}\widetilde{\MorawetzLagrangeCorrSym}\mathring{\chi}
  - \widetilde{\MorawetzLagrangeCorrSym}_{\Aux}\HorCovDeriv_t,
\end{equation*}
where the operator $\widetilde{\MorawetzLagrangeCorrSym}_{\Aux}$ is chosen so that
\begin{equation*}
  \Metric_{M,a,\Lambda}^{-1}(dt, dt)\WeylQ{\widetilde{\MorawetzLagrangeCorrSym}_{\Aux}}
  + \WeylQ{\widetilde{\MorawetzLagrangeCorrSym}_{\Aux}}\Metric_{M,a,\Lambda}^{-1}(dt, dt)
  = K^w_{-1}.
\end{equation*}
With this choice of $\widetilde{\MorawetzLagrangeCorr}$, we can
enforce that actually $K^w_{-1}=0$, so that
\begin{equation*}
  \left(\left[\ScalarWaveOp[\Metric], \widetilde{\MorawetzVF}\right]
    + \ScalarWaveOp[\Metric]\widetilde{\MorawetzLagrangeCorr}
    + \widetilde{\MorawetzLagrangeCorr}\ScalarWaveOp[\Metric]
  \right)
  = K^w_2 + 2K_1^w\HorCovDeriv_t + K_0^w\HorCovDeriv_t^2. 
\end{equation*}
We now have that the principal symbol of the bilinear product in
$\KCurrentPert{\widetilde{\MorawetzVF},
  \widetilde{\MorawetzLagrangeCorr}}[\psi]$
will be equal to
\begin{equation*}
  \KCurrentPertSym{\widetilde{\MorawetzVF}, \widetilde{\MorawetzLagrangeCorr}} \vcentcolon=
  \mathring{\chi}^2\left(\frac{1}{2\ImagUnit}\PoissonB*{\PrinSymb, \widetilde{\MorawetzSym}}
    + \PrinSymb\widetilde{\MorawetzLagrangeCorrSym} \right)
  + \frac{1}{\ImagUnit}\mathring{\chi}\PoissonB{\PrinSymb, \mathring{\chi}}.
\end{equation*}
We now split
\begin{equation}
  \label{eq:KCurrent-princ-aux:decomposition}
  \KCurrentPert{\widetilde{\MorawetzVF},\widetilde{\MorawetzLagrangeCorr}}[\psi]
  = \KCurrentPert{\widetilde{\MorawetzVF},\widetilde{\MorawetzLagrangeCorr}}_{\Main}[\psi]
  + \KCurrentPert{\widetilde{\MorawetzVF},\widetilde{\MorawetzLagrangeCorr}}_{\Aux}[\psi],
\end{equation}
where
\begin{equation}
  \label{eq:KCurrent-princ:def}
  \begin{split}
    \KCurrentPert{\widetilde{\MorawetzVF},\widetilde{\MorawetzLagrangeCorr}}_{\Main}[\psi]
    \vcentcolon={}& \bangle*{\widetilde{K}^w_{2,\Main}\psi, \psi}_{L^2(\Manifold(\tau_1,\tau_2))}
                    + 2\Re\bangle*{\widetilde{K}^w_{1,\Main}\psi, \HorCovDeriv_t\psi}_{L^2(\Manifold(\tau_1,\tau_2))}\\
                  &+ \bangle*{\widetilde{K}^w_{0,\Main}\HorCovDeriv_t\psi, \HorCovDeriv_t\psi}_{L^2(\Manifold(\tau_1,\tau_2))},    
  \end{split}
\end{equation}
where $\widetilde{K}^w_{2,\Main}$, $i=0,1,2$ are defined so that
\begin{equation*}
  \sum_{i=0}^2
  \begin{pmatrix}
    2\\ i
  \end{pmatrix}
  \widetilde{K}^w_{i,\Main}\HorCovDeriv_t^i
  = \mathring{\chi}\WeylQ{\frac{1}{2\ImagUnit}\PoissonB*{\PrinSymb,\widetilde{\MorawetzSym}} + \PrinSymb\widetilde{\MorawetzLagrangeCorrSym}}\mathring{\chi}.
\end{equation*}
We similarly define $K^w_{i,\Aux}, i=0,1,2$ so that
\begin{equation*}
  \begin{split}
    \KCurrentPert{\widetilde{\MorawetzVF},\widetilde{\MorawetzLagrangeCorr}}_{\Aux}[\psi](\tau_1,\tau_2)
  ={}& \bangle*{K_{2,\Aux}^w\psi,\psi}_{L^2(\Manifold(\tau_1,\tau_2))}
  + 2\Re \bangle*{K_{1,\Aux}^w\HorCovDeriv_t\psi,\psi}_{L^2(\Manifold(\tau_1,\tau_2))}\\
  &+ \bangle*{K_{2,\Aux}^w\HorCovDeriv_t\psi,\HorCovDeriv_t\psi}_{L^2(\Manifold(\tau_1,\tau_2))}
  - \bangle*{\WeylQ{\frac{1}{\ImagUnit}\mathring{\chi}\widetilde{\MorawetzVF}H_{\PrinSymb}\mathring{\chi}}\psi,\psi}_{L^2(\Manifold(\tau_1,\tau_2))}.
  \end{split}  
\end{equation*}
Then, observe that the norm $\MorNorm[\psi](\tau_1,\tau_2)$ controls
the full spacetime $H^1$ norm outside an $O(a)$ neighborhood of the
trapped set, and that for $a$ sufficiently small,
$\KCurrentPert{\widetilde{\MorawetzVF},\widetilde{\MorawetzLagrangeCorr}}_{\Aux,\DomainOfIntegration}[\psi]$
has principal symbols with support away from
$\Manifold_{\operatorname{trap}}$. As a result, we can write that
\begin{equation}
  \label{eq:Morawetz:K-ext-bound}
  \abs*{\KCurrentPert{\widetilde{\MorawetzVF},\widetilde{\MorawetzLagrangeCorr}}_{\Aux}[\psi](\tau_1,\tau_2)}
  \lesssim \MorNorm[\psi](\tau_1,\tau_2)
  + \norm*{\HorCovDeriv_t\psi}_{H^{-1}_{\operatorname{comp}}(\Manifold)(\tau_1,\tau_2))}^2.
\end{equation}
We now prove the following lemma on the cutoff pseudo-differential operators.
\begin{lemma}
  \label{lemma:Morawetz:cutoff-PsiDO:bulk}  
  Let
  $(\MorawetzVF_{M,a,\Lambda}, \MorawetzLagrangeCorr_{M,a,\Lambda},
  \MorawetzOneForm_{M,a,\Lambda})$ be as chosen in
   \zcref[cap]{prop:Morawetz:outside-trapping:KCurrent-coercive},
  and let $\widetilde{\MorawetzVF}$ and
  $\widetilde{\MorawetzLagrangeCorr}$ be as defined in
  \zcref[noname]{eq:Morawetz:morawetz-VF-LagrangeCorr-with-cutoff:def}, where
  $\widetilde{\MorawetzSym}$ and
  $\widetilde{\MorawetzLagrangeCorrSym}$ are as chosen in 
  \zcref[cap]{lemma:ILED-KdS:Bulk}. Then,
  \begin{equation}
    \label{eq:Morawetz:cutoff-PsiDO:bulk}
    \begin{split}
      & a \KCurrentPert{\widetilde{\MorawetzVF}, \widetilde{\MorawetzLagrangeCorr}}_{\Main}[\psi](\tau_1,\tau_2)
        + \int_{\Manifold(\tau_1,\tau_2)}\KCurrent{\MorawetzVF_{M,a,\Lambda}, \MorawetzLagrangeCorr_{M,a,\Lambda}, \MorawetzOneForm_{M,a,\Lambda}}[\psi]\\
      \gtrsim{}&  \MorNormTrap[\psi](\tau_1,\tau_2)
                 - O(a)\left(
                 \norm*{\psi}_{L^2_{\operatorname{comp}}(\Manifold(\tau_1,\tau_2))}^2
                 + \int_{\Manifold(\tau_1,\tau_2)}\abs*{N}^2
                 + \sup_{\tau\in[\tau_1,\tau_2]}\EnergyHorizonDeg[\psi](\tau)
                 \right).
    \end{split}    
  \end{equation}
\end{lemma}
\begin{proof}
  We first decompose
  \begin{equation*}
    \begin{split}
      \int_{\Manifold(\tau_1,\tau_2)}\KCurrent{\MorawetzVF_{M,a,\Lambda}, \MorawetzLagrangeCorr_{M,a,\Lambda}, \MorawetzOneForm_{M,a,\Lambda}}[\psi]
    ={}& \int_{\Manifold(\tau_1,\tau_2)}\mathring{\chi}^2\KCurrent{\MorawetzVF_{M,a,\Lambda}, \MorawetzLagrangeCorr_{M,a,\Lambda}, \MorawetzOneForm_{M,a,\Lambda}}[\psi]\\
    & + \int_{\Manifold(\tau_1,\tau_2)}\breve{\chi}^2\KCurrent{\MorawetzVF_{M,a,\Lambda}, \MorawetzLagrangeCorr_{M,a,\Lambda}, \MorawetzOneForm_{M,a,\Lambda}}[\psi].
    \end{split}    
  \end{equation*}
  Then, since
  $\KCurrent{\MorawetzVF_{M,a,\Lambda},
    \MorawetzLagrangeCorr_{M,a,\Lambda},
    \MorawetzOneForm_{M,a,\Lambda}}[\psi]$ is pointwise positive
  outside of a neighborhood of $r=3M$, then
  \begin{equation}
    \label{eq:Morawetz:PsiDO-bulk:outer:representation}
    \begin{split}
      \int_{\Manifold(\tau_1,\tau_2)}\breve{\chi}^2\KCurrent{\MorawetzVF_{M,a,\Lambda}, \MorawetzLagrangeCorr_{M,a,\Lambda}, \MorawetzOneForm_{M,a,\Lambda}}[\psi]
      &\gtrsim \int_{\Manifold(\tau_1,\tau_2)}\breve{\chi}^2\left(r^{-2}\abs*{\HorCovDeriv \psi}^2 + r^{-4}\abs*{\psi}^2\right)\\
      &\gtrsim\MorNorm_{\nontrap}[\psi](\tau_1,\tau_2) . 
    \end{split}
  \end{equation}
  On the other hand, we see that
  \begin{align*}
    &a\KCurrentPert{\widetilde{\MorawetzVF},\widetilde{\MorawetzLagrangeCorr}}_{\Main}[\psi](\tau_1,\tau_2)
    + \int_{\Manifold(\tau_1,\tau_2)}\mathring{\chi}^2\KCurrent{\MorawetzVF_{M,a,\Lambda}, \MorawetzLagrangeCorr_{M,a,\Lambda}, \MorawetzOneForm_{M,a,\Lambda}}[\psi]\\
    ={}& \Re\int_{\Manifold(\tau_1,\tau_2)}\mathring{\chi}^2\left(\KCurrentSym{\MorawetzVF_{M,a,\Lambda}, \MorawetzLagrangeCorr_{M,a,\Lambda}, \MorawetzOneForm_{M,a,\Lambda}}_{(2)}\right)^{\alpha\beta}\HorCovDeriv_{\alpha}\psi\cdot \overline{\HorCovDeriv_{\beta}\psi}
    + \mathring{\chi}^2\KCurrentSym{\MorawetzVF_{M,a,\Lambda}, \MorawetzLagrangeCorr_{M,a,\Lambda}, \MorawetzOneForm_{M,a,\Lambda}}_{(0)}\abs*{\psi}^2\\
    &+ a\Re \int_{\Manifold(\tau_1,\tau_2)}K^w_{2,\Main}\psi\cdot \overline{\psi} + 2 K^w_{1,\Main}\HorCovDeriv_t\psi\cdot \overline{\psi} + K^w_{0,\Main}\HorCovDeriv_t\psi\cdot\overline{\HorCovDeriv_t\psi},
  \end{align*}
  where for $a$ sufficiently small, we have that
  \begin{equation*}
    \left(\KCurrentSym{\MorawetzVF_{M,a,\Lambda}, \MorawetzLagrangeCorr_{M,a,\Lambda}, \MorawetzOneForm_{M,a,\Lambda}}_{(2)}\right)^{\alpha\beta}\eta_{\alpha}\eta_{\beta} = \frac{1}{2\ImagUnit}H_{\PrinSymb}\MorawetzVF_{M,a,\Lambda} + \PrinSymb \MorawetzLagrangeCorr_{M,a,\Lambda}, \qquad
    \KCurrentSym{\MorawetzVF_{M,a,\Lambda}, \MorawetzLagrangeCorr_{M,a,\Lambda}, \MorawetzOneForm_{M,a,\Lambda}}_{(0)} > 0.
  \end{equation*}
  Moreover, recall that by construction,
  \begin{equation*}
    K^w_{2,\Main} + 2K^w_{1,\Main}\HorCovDeriv_t + K^w_{0,\Main}\HorCovDeriv_t^2
    = \mathring{\chi}\WeylQ{\frac{1}{2\ImagUnit}H_{\PrinSymb}\widetilde{\MorawetzSym} + \PrinSymb\widetilde{\MorawetzLagrangeCorrSym}}\mathring{\chi}.
  \end{equation*}
  Observe that by \zcref[cap]{lemma:ILED-KdS:Bulk}, we have that
  $\SquareDecomp_k(a)$ are in general symbols of pseudodifferential
  operators. However, $\evalAt*{\SquareDecomp_k}_{a=0}$ are
  symbols of differential operators (in particular, they agree with
  the principal symbol of the decomposition in \zcref[cap]{prop:Morawetz:outside-trapping:KCurrent-coercive}). As a
  result, we can decompose $\SquareDecomp_k =
    \evalAt*{\SquareDecomp_k}_{a=0}
    + \SquareDecomp_k
    - \evalAt*{\SquareDecomp_k}_{a=0}$
  where $\evalAt*{\SquareDecomp_k}_{a=0}$ is the symbol of
  a differential operator, and $\SquareDecomp_k - \evalAt*{\SquareDecomp_k}_{a\eta_{\phi}=0}
    \in a \MixedSymClass{2}{2}\left(T^{*}\Manifold;\pi^{*}\End(\realHorkTensor{2}(\Manifold))\right)$. Then, define the operators
  \begin{equation*}
    \SquareDecompOp_k
    = \mathring{\chi}\evalAt*{\SquareDecomp_k}_{a=0}(x, D)
    + \WeylQ{\SquareDecomp_k(a) - \evalAt*{\SquareDecomp_k}_{a=0}}\mathring{\chi}.
  \end{equation*}
  Then from the Weyl calculus, we have that
  \begin{equation}
    \label{eq:Morawetz:PsiDO-bulk:inner:representation}
    \begin{split}
      &a\KCurrentPert{\widetilde{\MorawetzVF},\widetilde{\MorawetzLagrangeCorr}}_{\Main}[\psi](\tau_1,\tau_2)
        + \int_{\Manifold(\tau_1,\tau_2)}\mathring{\chi}^2\KCurrent{\MorawetzVF_{M,a,\Lambda}, \MorawetzLagrangeCorr_{M,a,\Lambda}, \MorawetzOneForm_{M,a,\Lambda}}[\psi] \\
      ={}& \int_{\Manifold(\tau_1,\tau_2)}\sum_k\abs*{\SquareDecompOp_k\psi}^2
           + \KCurrentSym{\MorawetzVF_{M,a,\Lambda}, \MorawetzLagrangeCorr_{M,a,\Lambda}, \MorawetzOneForm_{M,a,\Lambda}}_{(0)}\mathring{\chi}^2\abs*{\psi}^2\\
      & + \Re \int_{\Manifold(\tau_1,\tau_2)}R^w_2\psi\cdot \overline{\psi} + 2R^w_1\HorCovDeriv_t\psi\cdot\overline{\psi} + R^w_0\HorCovDeriv_t\psi\cdot\overline{\HorCovDeriv_t\psi},      
    \end{split}
  \end{equation}
  where the remainder terms $R^w_i$ have symbols
  $r_j\in a
  \TanSymClass{j-2}(\Manifold;\realHorkTensor{2}(\Manifold))$ from the
  Weyl calculus.  Combining
  \zcref[noname]{eq:Morawetz:PsiDO-bulk:inner:representation,eq:Morawetz:PsiDO-bulk:outer:representation} and
  using
  \zcref[cap]{lemma:ILED-near:LoT-control:exchange-degeneracy-trick,lemma:ILED-near:LoT-control:zero-order-mixed-term}
  to control the remainder terms concludes the proof of
  \zcref[cap]{lemma:Morawetz:cutoff-PsiDO:bulk}.
\end{proof}

\subsection{Concluding the proof of \zcref[cap]{prop:Morawetz:KdS:main}}

We first prove a series of auxiliary lemmas that help control terms in
the Morawetz estimate that arise from the terms on the \RHS{} of
\zcref[noname]{eq:model-problem-gRW}.

\paragraph{Controlling the \RHS}

On the \RHS{} of \zcref[noname]{eq:model-problem-gRW}, we observe that we have two
terms: a forcing term, and a first-order derivative of $\psi$. The
forcing term is generally not a problem, since any norms are absorbed
into $\ForcingTermCombinedNorm{p}[\psi,N](\tau_1,\tau_2)$. However,
the presence of first-order derivative of $\psi$ on the \RHS{} may
appear concerning, as this type of term generally leads to
principal-order terms appearing when multiplying the equation with a
vectorfield. Fortunately, this term causes no such complications here
thanks to its anti-symmetric nature. This in essence allows such terms
to be controlled by a combination of lower-order terms and boundary
terms. Since the derivative on the \RHS{} is $O(a)$, this means that
we will in effect be able to absorb all resulting terms except for
those involving the forcing term, back onto the left-hand side. In
what follows, we show exactly how this is done, beginning with a lemma highlighting how this antisymmetry is used.
\begin{lemma}
  \label{lemma:Morawetz:RHS:anti-symmetry:simplification}
  We have that
  \begin{equation}
    \label{eq:Morawetz:RHS:anti-symmetry:simplification}
    \begin{split}
      \left(\nabla_{\MorawetzVF_{M,a,\Lambda}}\psi + \MorawetzLagrangeCorr_{M,a,\Lambda}\psi\right)\cdot\nabla_{\KillT}(\LeftDual{\psi})
      ={}&  \frac{1}{2}\left(\nabla_r\left(\mathcal{F}_{M,a,\Lambda}\psi\cdot\nabla_{\KillT}\LeftDual{\psi}\right)
           - \nabla_{\KillT}\left(\psi\cdot\nabla_{\MorawetzVF_{M,a,\Lambda}}\psi\right)\right)\\
         & - \frac{1}{2}\partial_r(\mathcal{F}_{M,a,\Lambda})\psi\cdot\nabla_{\KillT}\LeftDual{\psi}
           - \mathcal{F}_{M,a,\Lambda}\LeftDual{\rho}\abs*{\psi}^2.
    \end{split}
  \end{equation}
\end{lemma}
\begin{proof}
  We observe that
  \begin{align*}
    \nabla_{\MorawetzVF_{M,a,\Lambda}}\psi\cdot\nabla_{\KillT}\LeftDual{\psi}
    ={}& \nabla_r\left(\mathcal{F}_{M,a,\Lambda} \psi\cdot\nabla_{\KillT}\LeftDual{\psi}\right)
         - \partial_r\left(\mathcal{F}_{M,a,\Lambda}\right)\psi\cdot\nabla_{\KillT}\LeftDual{\psi}
         - \mathcal{F}_{M,a,\Lambda}\nabla_r\nabla_{\KillT}\LeftDual{\psi}.
  \end{align*}
  Then, observe that
  \begin{align*}
    \mathcal{F}_{M,a,\Lambda}\nabla_r\nabla_{\KillT}\LeftDual{\psi}
    ={}& \mathcal{F}_{M,a,\Lambda}\nabla_{\KillT}\nabla_r\LeftDual{\psi}
         + \mathcal{F}_{M,a,\Lambda}\squareBrace*{\nabla_r, \nabla_{\KillT}}\LeftDual{\psi}\\
    ={}& \partial_t\left(\psi\cdot\nabla_{\MorawetzVF_{M,a,\Lambda}}\LeftDual{\psi}\right)
         - \nabla_{\KillT}\psi\cdot\nabla_{\MorawetzVF_{M,a,\Lambda}}\LeftDual{\psi}
         -2 \mathcal{F}_{M,a,\Lambda}\psi\cdot\LeftDual{\rho}\LeftDual{\left(\LeftDual{\psi}\right)}\\
    ={}& \partial_t\left(\psi\cdot\nabla_{\MorawetzVF_{M,a,\Lambda}}\LeftDual{\psi}\right)
         - \nabla_{\KillT}\psi\cdot\nabla_{\MorawetzVF_{M,a,\Lambda}}\LeftDual{\psi}
         + 2 \mathcal{F}_{M,a,\Lambda}\LeftDual{\rho}\abs*{\psi}^2.
  \end{align*}
  As a result, we have that
  \begin{align*}
    2\nabla_{\MorawetzVF_{M,a,\Lambda}}\psi\cdot\nabla_{\KillT}\LeftDual{\psi}    
    ={}& \nabla_{\MorawetzVF_{M,a,\Lambda}}\psi\cdot\nabla_{\KillT}\LeftDual{\psi}
         + \nabla_{\KillT}\psi\cdot\nabla_{\MorawetzVF_{M,a,\Lambda}}\LeftDual{\psi}\\
       & + \nabla_r\left(\mathcal{F}_{M,a,\Lambda}\psi\cdot\nabla_{\KillT}\LeftDual{\psi}\right)
         - \nabla_{\KillT}\left(\psi\cdot\nabla_{\MorawetzVF_{M,a,\Lambda}}\psi\right)\\
       &  - \partial_r(\mathcal{F}_{M,a,\Lambda})\psi\cdot\nabla_{\KillT}\LeftDual{\psi}
         - 2\mathcal{F}_{M,a,\Lambda}\LeftDual{\rho}\abs*{\psi}^2. 
  \end{align*}
  Then, recalling the antisymmetry of the dual operator and $\nabla_r\psi\cdot\nabla_t\LeftDual{\psi} = -\nabla_t\psi\cdot\nabla_r\LeftDual{\psi}$
  we have that
  \begin{equation*}
    \begin{split}
      2\nabla_{\MorawetzVF_{M,a,\Lambda}}\psi\cdot\nabla_{\KillT}\LeftDual{\psi}    
      ={}& \nabla_r\left(\mathcal{F}_{M,a,\Lambda}\psi\cdot\nabla_{\KillT}\LeftDual{\psi}\right)
           - \nabla_{\KillT}\left(\psi\cdot\nabla_{\MorawetzVF_{M,a,\Lambda}}\psi\right)\\
         &  - \partial_r(\mathcal{F}_{M,a,\Lambda})\psi\cdot\nabla_{\KillT}\LeftDual{\psi}
           - 2\mathcal{F}_{M,a,\Lambda}\LeftDual{\rho}\abs*{\psi}^2,
    \end{split}
  \end{equation*}
  as desired.  
\end{proof}

We can then use \zcref[cap]{lemma:Morawetz:RHS:anti-symmetry:simplification} to prove the
following lemma which will help us control the error outside the
trapping region.
\begin{lemma}
  \label{lemma:Morawetz:RHS:estimate}
  Let $\psi$ be a solution to the model problem in
  \zcref[noname]{eq:model-problem-gRW}.  Let
  $\MorawetzVF_{M,a,\Lambda}$ and
  $\MorawetzLagrangeCorr_{M,a,\Lambda}$ be as chosen in \zcref[cap]{lemma:ILED:SdS:Bulk:principal-terms}. Then, for $0<\delta_2\ll 1$
  sufficiently small,
  \begin{equation}
    \label{eq:Morawetz:RHS:estimate:main-terms}
    \begin{split}
      &\left(\HorCovDeriv_{\MorawetzVF_{M,a,\Lambda}}\psi + \MorawetzLagrangeCorr_{M,a,\Lambda}\psi\right)\cdot\left(\WaveOpHork{2}\psi-V\psi\right)\\
    \ge{}& -\delta_2\frac{M}{r}\abs*{\HorCovDeriv_{\HprVF}\psi}^2
    - \delta_2\left(1-\frac{3M}{r}\right)^2\frac{M}{r^2}\abs*{\HorCovDeriv_{\HawkingVF}\psi}^2\\
    & + O(1)\left(\abs*{\HorCovDeriv_{\HprVF}\psi}+ r^{-1}\abs*{\psi}\right)\abs*{N}
    + O\left(a^3r^{-6}\right)\abs*{\HorCovDeriv_{\KillPhi}\psi}^2
    + O\left(ar^{-4}\right)\abs*{\psi}^2\\
    & - \CovariantDeriv_{\mu}\left(\frac{2a\cos\theta}{\abs*{q}^2}\partial_r^{\mu}\mathcal{F}_{M,a,\Lambda}\psi\cdot\HorCovDeriv_{\KillT}\LeftDual{\psi}\right)
    + \partial_t\left(\frac{2a\cos\theta}{\abs*{q}^2}\mathcal{F}_{M,a,\Lambda}\psi\cdot\HorCovDeriv_r\LeftDual{\psi}\right)
    \end{split}    
  \end{equation}
  and for $0<\delta_3\ll 1$ sufficiently small,
  \begin{equation}
    \label{eq:Morawetz:RHS:estimate:curvature-terms}
    \begin{split}
      &\abs*{\left(
      \left(\LeftDual{\rho} + \etaBar\wedge \eta\right)\frac{r^2+a^2}{\Delta}
      + \frac{2a^3r\cos\theta\sin^2\theta}{\abs*{q}^6}
      \right)\mathcal{F}_{M,a,\Lambda}\HorCovDeriv_{\HawkingVF}\psi\cdot\LeftDual{\psi}}
      + \abs*{\left(
      \frac{2a^2r\cos\theta \mathcal{F}_{M,a,\Lambda}}{(r^2+a^2)\abs*{q}^4}\HorCovDeriv_{\KillPhi}\psi\cdot\LeftDual{\psi}
      \right)}\\
      \lesssim{}& \delta_3\left(1-\frac{3M}{r}\right)^2\left(\frac{M}{r^2}\abs*{\HorCovDeriv_{\HawkingVF}\psi}^2 + \frac{1}{r^4}\abs*{\HorCovDeriv_{\KillPhi}\psi}^2\right)
      + O\left(a^2r^{-6}\right)\abs*{\psi}^2. 
    \end{split}
  \end{equation}
\end{lemma}
\begin{proof}
  Since $\psi$ solves \zcref[noname]{eq:model-problem-gRW}, we have that
  \begin{equation*}
    \left(\HorCovDeriv_{\MorawetzVF_{M,a,\Lambda}}\psi + \MorawetzLagrangeCorr_{M,a,\Lambda}\psi\right)\cdot\left(\WaveOpHork{2}\psi-V\psi\right)
    = \left(\HorCovDeriv_{\MorawetzVF_{M,a,\Lambda}}\psi + \MorawetzLagrangeCorr_{M,a,\Lambda}\psi\right)\cdot\left(-\frac{4a\cos\theta}{\abs*{q}^2} \LeftDual{\nabla_{\KillT}\psi} + N\right).
  \end{equation*}
  We first consider the term $-\frac{4a\cos\theta}{\abs*{q}^2}\left(\HorCovDeriv_{\MorawetzVF_{M,a,\Lambda}}\psi + \MorawetzLagrangeCorr_{M,a,\Lambda}\psi\right)\cdot \LeftDual{\nabla_{\KillT}\psi}$. Using \zcref[cap]{lemma:Morawetz:RHS:anti-symmetry:simplification}, we have that
  \begin{align*}
    -\frac{4a\cos\theta}{\abs*{q}^2}\left(\HorCovDeriv_{\MorawetzVF_{M,a,\Lambda}}\psi + \MorawetzLagrangeCorr_{M,a,\Lambda}\psi\right)\cdot \LeftDual{\nabla_{\KillT}\psi}
    ={}& -\frac{2a\cos\theta}{\abs*{q}^2}\left(\nabla_r\left(\mathcal{F}_{M,a,\Lambda}\psi\cdot\nabla_{\KillT}\LeftDual{\psi}\right)
         - \nabla_{\KillT}\left(\psi\cdot\nabla_{\MorawetzVF_{M,a,\Lambda}}\psi\right)\right)\\
       & + \frac{2a\cos\theta}{\abs*{q}^2}\left(
         \partial_r(\mathcal{F}_{M,a,\Lambda})\psi\cdot\nabla_{\KillT}\LeftDual{\psi}
         + 2\mathcal{F}_{M,a,\Lambda}\LeftDual{\rho}\abs*{\psi}^2\right)\\
    & - \frac{4a\cos\theta}{\abs*{q}^2}\MorawetzLagrangeCorr_{M,a,\Lambda}\psi\cdot\LeftDual{\nabla}_{\KillT}\psi.
  \end{align*}
  Recall from the choices made in \zcref[cap]{prop:Morawetz:outside-trapping:KCurrent-coercive} that $\mathcal{F}_{M,a,\Lambda} = \frac{(r-3M)\Delta}{r(r^2+a^2)}$ and $\MorawetzLagrangeCorr_{M,a,\Lambda} = \MorawetzLagrangeCorr_{*} + \MorawetzLagrangeCorr_{\delta_1}$ where
  \begin{equation*}
    \begin{split}
      \MorawetzLagrangeCorr_{*} ={}& \MorawetzLagrangeCorr_0 + \MorawetzLagrangeCorr_1,\\
      \MorawetzLagrangeCorr_0 ={}& - \frac{1}{2}\abs*{q}^2\Divergence\left(\abs*{q}^{-2}\MorawetzVF_{M,a,\Lambda}\right),
    \end{split}    
  \end{equation*}
  which in particular gives $\partial_r\left(\mathcal{F}_{M,a,\Lambda}\right) - 2\MorawetzLagrangeCorr_0 ={} 0$. As a result,
  \begin{align*}
    \left(\nabla_{\MorawetzVF_{M,a,\Lambda}}\psi + \MorawetzLagrangeCorr_{M,a,\Lambda}\psi\right)\cdot\nabla_{\KillT}\LeftDual{\psi}
    ={}& \frac{1}{2}\nabla_r\left(\mathcal{F}_{M,a,\Lambda}\psi\cdot\nabla_{\KillT}\LeftDual{\psi}\right)
         - \frac{1}{2}\nabla_{\KillT}\left(\psi\cdot\nabla_{\MorawetzVF_{M,a,\Lambda}}\psi\right)\\
       & + \left(\MorawetzLagrangeCorr_1 + \MorawetzLagrangeCorr_{\delta_1}\right)\psi\cdot\nabla_{\KillT}\LeftDual{\psi}
         - 2\mathcal{F}_{M,a,\Lambda}\LeftDual{\rho}\abs*{\psi}^2.
  \end{align*}
  Recall from
  \zcref[noname]{eq:Morawetz:w1:def,eq:Morawetz:outside-trapping:q1-def}
  that
  \begin{equation*}
    \begin{split}
      \MorawetzLagrangeCorr_1 + \MorawetzLagrangeCorr_{\delta_1}
      ={}& \frac{(r-3M)(3r^2(r-3M)+a^4r\Lambda + a^2(3M+3r+r^3\Lambda)}{3r(r^2+a^2)^2}
           -\delta_1\frac{M\Delta(r-3M)^2}{r^2(r^2+a^2)^2}\\
      ={}& \left(1-\frac{3M}{r}\right)O(r^{-1}).  
    \end{split}
  \end{equation*}
  As a result, for any $\delta_2$, we can estimate
  \begin{align*}
    &\abs*{\frac{4a\cos\theta}{\abs*{q}^2}\left(\MorawetzLagrangeCorr_1+\MorawetzLagrangeCorr_{\delta_1}\right)\psi\cdot\LeftDual{\nabla}_{\KillT}\psi}\\
    \lesssim{}& \frac{a}{\abs*{q}^2r}\left(
                \delta_2r\left(1-\frac{3M}{r}\right)^2\abs*{\nabla_{\KillT}\psi}^2 + \delta_2^{-1}r^{-1}\abs*{\psi}^2
                \right)
                + O(ar^{-4})\abs*{\psi}^2\\
    \lesssim{}& \frac{a}{\abs*{q}^2}\frac{(r-3M)^2}{r^3}\left(
                \delta_2r\abs*{\nabla_{\HawkingVF}\psi}^2
                + \delta_2\frac{2ar}{r^2+a^2}\abs*{\nabla_{\HawkingVF}\psi\cdot\nabla_{\KillPhi}}
                + \delta_2\frac{a^2r}{\left(r^2+a^2\right)^2}\abs*{\nabla_{\KillPhi}\psi}^2 
                \right)
                + O(ar^{-4})\abs*{\psi}^2.
  \end{align*}
  As a result, we have that
  \begin{equation}
    \label{eq:Morawetz:RHS:estimate:main-terms:aux-1}
    \begin{split}
      \abs*{\frac{4a\cos\theta}{\abs*{q}^2}\left(\MorawetzLagrangeCorr_1+\MorawetzLagrangeCorr_{\delta_1}\right)\cdot\LeftDual{\nabla}_{\KillT}\psi}
      \ge{}&  - \delta_2\left(1-\frac{3M}{r}\right)^2\frac{M}{r^2}\abs*{\nabla_{\HawkingVF}\psi}^2
             + O(a^3r^{-6})\abs*{\nabla_{\KillPhi}\psi}^2
             + O(ar^{-4})\abs*{\psi}^2\\
           & - \CovariantDeriv\cdot\left(\frac{2a\cos\theta}{\abs*{q}^2}\mathcal{F}_{M,a,\Lambda}\psi\cdot\nabla_{\KillT}\LeftDual{\psi}\partial_r\right)
            + \partial_t\left(\frac{2a\cos\theta}{\abs*{q}^2}\mathcal{F}\psi\cdot\nabla_r\LeftDual{\psi}\right).  
    \end{split}
  \end{equation}
  As a result, recalling that
  $\MorawetzVF = \mathcal{F}_{M,a,\Lambda}\partial_r = f_{M,a,\Lambda}\HprVF$, and that
  $\MorawetzLagrangeCorr_{M,a,\Lambda} = O(r^{-1})$, we have that
  \begin{equation}
    \label{eq:Morawetz:RHS:estimate:main-terms:aux-2}
    \abs*{\left(\nabla_{\MorawetzVF}\psi + \MorawetzLagrangeCorr\psi\right)\cdot N}
    \lesssim \left(\abs*{\nabla_{\HprVF}\psi} + r^{-1}\abs*{\psi}\right)\abs*{N}.
  \end{equation}
  Combining \zcref[noname]{eq:Morawetz:RHS:estimate:main-terms:aux-1} and
  \zcref[noname]{eq:Morawetz:RHS:estimate:main-terms:aux-2}
  yields~\zcref[noname]{eq:Morawetz:RHS:estimate:main-terms}, as desired. We now prove
  \zcref[noname]{eq:Morawetz:RHS:estimate:curvature-terms}. Recall from
  \zcref[cap]{lemma:Kerr:ingoing-PG:Ric-and-curvature,lemma:Kerr:outgoing-PG:Ric-and-curvature},
  and
  \zcref[noname]{eq:ILED:SdS:X-def}, that
  \begin{equation*}
    \abs*{\LeftDual{\rho}}\lesssim \frac{aM}{r^4},\qquad
    \abs*{\eta} + \abs*{\etaBar} \lesssim \frac{a}{r^2}, \qquad
    \abs*{\mathcal{F}_{M,a,\Lambda}} \lesssim \frac{\abs*{\Delta}}{r^2}\frac{\abs*{r-3M}}{r}.
  \end{equation*}
  As a result, we have that
  \begin{equation*}
    \begin{split}
      &\abs*{\left(
      \left(\LeftDual{\rho} + \etaBar\wedge \eta\right)\frac{r^2+a^2}{\Delta}
      + \frac{2a^3r\cos\theta\sin^2\theta}{\abs*{q}^6}
      \right)\mathcal{F}_{M,a,\Lambda}\nabla_{\HawkingVF}\psi\cdot\LeftDual{\psi}}
      + \abs*{\frac{2a^2r\cos\theta\mathcal{F}_{M,a,\Lambda}}{(r^2+a^2)\abs*{q}^4}\nabla_{\KillPhi}\psi\cdot\LeftDual{\psi}}\\
      \lesssim{}& \frac{aM}{r^4}\abs*{1-\frac{3M}{r}}\abs*{\nabla_{\HawkingVF}\psi}\abs*{\psi}
      + \frac{aM}{r^5}\abs*{1-\frac{3M}{r}}\abs*{\nabla_{\KillPhi}\psi}\abs*{\psi}\\
    \lesssim{}& \delta_3\left(1-\frac{3M}{r}\right)^2
                  \left(\frac{M}{r^2}\abs*{\nabla_{\HawkingVF}\psi}^2 + \frac{1}{r^4}\abs*{\nabla_{\KillPhi}\psi}^2\right)
      + O(a^2r^{-6})\abs*{\psi}^2,
    \end{split}
  \end{equation*}
  as desired, proving \zcref[noname]{eq:Morawetz:RHS:estimate:curvature-terms}.
\end{proof}

We also need a similar estimate within the trapping region
adapted to pseudo-differential operators.
\begin{lemma}
  \label{lemma:Morawetz:RHS:PsiDO}
  Let $\widetilde{\MorawetzVF}, \widetilde{\MorawetzLagrangeCorr}$ be
  as chosen in \zcref[cap]{prop:Morawetz:outside-trapping:KCurrent-coercive}, and let
  $\psi$ be a solution to the model problem in \zcref[noname]{eq:model-problem-gRW}.  In
  addition, consider $\dot{\chi}=\dot{\chi}(r)$ such that
  \begin{equation*}
    \dot{\chi} =
    \begin{cases}
      1 & \abs*{r-3M}<\delta_{\operatorname{trap}},\\
      0 & \abs*{r-3M}>\delta_{\operatorname{trap}},
    \end{cases}
  \end{equation*}
  where 
  $\TrappedSet_{M,a,\Lambda}\subset \{\abs*{r-3M}<\delta_{\operatorname{trap}}\}$. 
  Then,
  \begin{equation}
   \label{eq:Morawetz:RHS:PsiDO:main-terms}
    2\Re\abs*{\bangle*{\left(\dot{\chi}\MorawetzVF_{M,a,\Lambda}+ \widetilde{\MorawetzVF} + \dot{\chi}\MorawetzLagrangeCorr_{M,a,\Lambda} + \widetilde{\MorawetzLagrangeCorr}\right)\psi,\WaveOpHork{2}\psi - V\psi}_{L^2(\Manifold_{\operatorname{trap}}(\tau_1,\tau_2))}}
    \lesssim \ForcingTermNorm[\psi,N]
    + a\MorNorm[\psi](\tau_1,\tau_2),
  \end{equation}
  and
  \begin{equation}
    \label{eq:Morawetz:RHS:PsiDO:curvature-terms}
    \begin{split}
      &\abs*{\bangle*{\dot{\chi}\left(
      \left(\LeftDual{\rho} + \etaBar\wedge \eta\right)\frac{r^2+a^2}{\Delta}
      + \frac{2a^3r\cos\theta\sin^2\theta}{\abs*{q}^6}
      \right)\mathcal{F}_{M,a,\Lambda}\HorCovDeriv_{\HawkingVF}\psi,\LeftDual{\psi}}_{L^2(\Manifold_{\trap}(\tau_1,\tau_2))}}\\
      &+ \abs*{\bangle*{\dot{\chi}
      \frac{2a^2r\cos\theta \mathcal{F}_{M,a,\Lambda}}{(r^2+a^2)\abs*{q}^4}\HorCovDeriv_{\KillPhi}\psi,\LeftDual{\psi}}_{L^2(\Manifold_{\trap}(\tau_1,\tau_2))}}\\
      \lesssim{}& a\MorNormTrap[\psi](\tau_1,\tau_2)
                  + a\sup_{\tau\in[\tau_1,\tau_2]}\EnergyHorizonDeg(\tau).
    \end{split}
  \end{equation}
\end{lemma}
\begin{proof}
  As we will see, we will only need extremely rough properties of the
  multipliers to prove this lemma. Thus, let us define
  \begin{equation*}
    \MorawetzVF^{\Aux}\vcentcolon=\dot{\chi}\MorawetzVF_{M,a,\Lambda}+ a\widetilde{\MorawetzVF},\qquad
    \MorawetzLagrangeCorr^{\Aux}\vcentcolon=\dot{\chi}\MorawetzLagrangeCorr_{M,a,\Lambda} + a\widetilde{\MorawetzLagrangeCorr}. 
  \end{equation*}    
  Since $\psi$ solves \zcref[noname]{eq:model-problem-gRW}, we have that
  \begin{equation*}
    \bangle*{\left(\MorawetzVF^{\Aux} + \MorawetzLagrangeCorr^{\Aux}\right)\psi,\WaveOpHork{2}\psi - V\psi}_{L^2(\Manifold_{\operatorname{trap}}(\tau_1,\tau_2))}
    = \bangle*{\left(\MorawetzVF^{\Aux} + \MorawetzLagrangeCorr^{\Aux}\right)\psi,-\frac{4a\cos\theta}{\abs*{q}^2}\LeftDual{\nabla}_{\KillT}\psi + N}_{L^2(\Manifold_{\operatorname{trap}}(\tau_1,\tau_2))}.
  \end{equation*}
  Then,
  \begin{equation*}
    \abs*{\bangle*{\left(\MorawetzVF^{\Aux} + \MorawetzLagrangeCorr^{\Aux}\right)\psi,N}_{L^2(\Manifold_{\operatorname{trap}}(\tau_1,\tau_2))}}
    \lesssim \ForcingTermNorm[\psi, N](\tau_1,\tau_2) 
  \end{equation*}
  directly from the fact that
  $\MorawetzVF^{\Aux} + \MorawetzLagrangeCorr^{\Aux}$ is a first-order operator. Next, we observe that
  $\MorawetzVF^{\Aux} + \MorawetzLagrangeCorr^{\Aux}$ is
  principally skew-symmetric while
  $-\frac{4a\cos\theta}{\abs*{q}^2}\LeftDual{\nabla}_{\KillT}\psi$ is
  principally symmetric in the sense that there exists a decomposition
  \begin{equation*}
    \MorawetzVF^{\Aux} + \MorawetzLagrangeCorr^{\Aux}
    = \widetilde{M}_1 + \widetilde{M}_0 \in \MixedOpClass{1}{1},\qquad
    -\frac{4a\cos\theta}{\abs*{q}^2}\LeftDual{\nabla}_{\KillT}
    = \widetilde{F}_1 + \widetilde{F}_0,
  \end{equation*}
  where $\widetilde{M}_i,\widetilde{F}_i$ are order-$i$ (mixed
  pseudo)-differential operators, such that $\widetilde{M}_1$ is
  skew-symmetric\footnote{Here, we recall that we define symmetry and
    skew-symmetry with respect to the actions on the
    $L^2_{comp}(\Manifold)$ inner product.} and $\widetilde{F}_1$ is
  symmetric. Moreover,
  \begin{align*}
    2\Re\bangle*{\widetilde{M}_1\psi, \widetilde{F}_1\psi}_{L^2(\Manifold_{\operatorname{trap}}(\tau_1,\tau_2))}
    ={}& \bangle*{\squareBrace*{\widetilde{F}_1,\widetilde{M}_1}\psi,\psi}_{L^2(\Manifold_{\operatorname{trap}}(\tau_1,\tau_2))},
  \end{align*}
  then observe that using the Weyl calculus, we have that $\squareBrace*{\widetilde{F}_1,\widetilde{M}_1}\in \MixedOpClass{1}{2}$
  and that moreover $\squareBrace*{\widetilde{F}_1,\widetilde{M}_1}$
  itself is principally skew-symmetric. Then let us write $\squareBrace*{\widetilde{F}_1,\widetilde{M}_1}
    = \widetilde{P} + \widetilde{Q}$
    where $\widetilde{P}\in a\MixedOpClass{1}{2}$ is skew-symmetric, and
  $\widetilde{Q}\in a\MixedOpClass{0}{2}$.  Then, using
  \zcref[cap]{lemma:Morawetz:lower-order:skew-adjoint} and
  \zcref[cap]{lemma:ILED-near:LoT-control:zero-order-mixed-term}, we
  can bound
  \begin{equation*}
    \abs*{\Re\bangle*{\widetilde{P}\psi,\psi}_{L^2(\Manifold_{\operatorname{trap}}(\tau_1,\tau_2))}}
    \lesssim a\norm*{\psi}^2_{L^2_{\operatorname{comp}}(\Manifold(\tau_1,\tau_2))}
    + a\norm*{N}^2_{L^2_{\operatorname{comp}}(\Manifold(\tau_1,\tau_2))}
    + a\sup_{\tau\in [\tau_1,\tau_2]}\EnergyHorizonDeg[\psi](\tau),
  \end{equation*}
  where we used $a$ small to control the
  $-\frac{4a\cos\theta}{\abs*{q}^2}\LeftDual{\nabla}_{\KillT}\psi$
  term that would have appeared on the right-hand
  side\footnote{Another skew-symmetry argument would also control this
    term, but for simplification, we use the smallness of $a$ to
    absorb it.}. The remaining terms involved are all of order $a$, and can be
  absorbed using a combination of
  \zcref[cap]{lemma:Morawetz:lower-order:skew-adjoint} and
  \zcref[cap]{lemma:ILED-near:LoT-control:zero-order-mixed-term} to
  prove \zcref[noname]{eq:Morawetz:RHS:PsiDO:main-terms}. To prove \zcref[noname]{eq:Morawetz:RHS:PsiDO:curvature-terms}, we first observe that $\bangle*{\psi,\LeftDual{\phi}}_{L^2(\Manifold(\tau_1,\tau_2))}
    = \bangle*{\LeftDual{\phi}, \psi}_{L^2(\Manifold(\tau_1,\tau_2))}$. Then, we observe that both
\begin{align*}
    \dot{\chi}\left( \left(\LeftDual{\rho} + \etaBar\wedge
      \eta\right)\frac{r^2+a^2}{\Delta} +
    \frac{2a^3r\cos\theta\sin^2\theta}{\abs*{q}^6}
  \right)\mathcal{F}_{M,a,\Lambda}\LeftDual{\HorCovDeriv}_{\HawkingVF} \quad \text{and} \quad \dot{\chi} \frac{2a^2r\cos\theta
    \mathcal{F}_{M,a,\Lambda}}{(r^2+a^2)\abs*{q}^4}\LeftDual{\HorCovDeriv}_{\KillPhi}
\end{align*}
 are principally symmetric operators in $a\MixedOpClass{1}{1}$. Thus,
  from \zcref[cap]{lemma:Morawetz:lower-order:skew-adjoint,lemma:ILED-near:LoT-control:exchange-degeneracy-trick}, we can
  conclude \zcref[noname]{eq:Morawetz:RHS:PsiDO:curvature-terms} directly. 
\end{proof}

We are now ready to prove \zcref[cap]{prop:Morawetz:KdS:main}.
\begin{proof}[Proof of Proposition \ref{prop:Morawetz:KdS:main}]
  Adding the divergence theorem in \zcref[noname]{eq:div-thm:general} with $(\MorawetzVF, \MorawetzLagrangeCorr, \MorawetzOneForm) =
    \left(\MorawetzVF_{M,a,\Lambda},
      \MorawetzLagrangeCorr_{M,a,\Lambda},
      \MorawetzOneForm_{M,a,\Lambda}\right)$
  as chosen in
  \zcref[cap]{prop:Morawetz:outside-trapping:KCurrent-coercive}
  to the commutation formula in
  \zcref[noname]{eq:frequency-setup:wave-commutation-formula}, we have that
  \begin{equation}
    \label{eq:Morawetz:KdS:main-combined-div-thm}
    \begin{split}
      &J_{\operatorname{Mor}}(\tau_1,\tau_2)
      + \int_{\Manifold(\tau_1,\tau_2)}\KCurrent{\MorawetzVF_{M,a,\Lambda},\MorawetzLagrangeCorr_{M,a,\Lambda},\MorawetzOneForm_{M,a,\Lambda}}[\psi]
      + a\KCurrentPert{\widetilde{\MorawetzVF},\widetilde{\MorawetzLagrangeCorr}}[\psi](\tau_1,\tau_2) 
      ={} \bangle*{\left(\WaveOpHork{2}-V\right)\psi, \left(\MorawetzVF + \MorawetzLagrangeCorr\right)\psi}_{L^2(\Manifold(\tau_1,\tau_2))},
    \end{split}
  \end{equation}
  where $\MorawetzVF = \MorawetzVF_{M,a,\Lambda} + a\widetilde{\MorawetzVF}$, 
    $\MorawetzLagrangeCorr = \MorawetzLagrangeCorr_{M,a,\Lambda} + a\widetilde{\MorawetzLagrangeCorr}$ and
  \begin{equation*}
    \begin{split}
      J_{\operatorname{Mor}}(\tau_1,\tau_2)
      ={}& \int_{\Sigma(\tau_2)}\JCurrent{\MorawetzVF_{M,a,\Lambda},\MorawetzLagrangeCorr_{M,a,\Lambda},\MorawetzOneForm_{M,a,\Lambda}}[\psi]\cdot N_{\Sigma}
           - \int_{\Sigma(\tau_1)}\JCurrent{\MorawetzVF_{M,a,\Lambda},\MorawetzLagrangeCorr_{M,a,\Lambda},\MorawetzOneForm_{M,a,\Lambda}}[\psi]\cdot N_{\Sigma}\\
         & + \int_{\SigmaStar}\JCurrent{\MorawetzVF_{M,a,\Lambda},\MorawetzLagrangeCorr_{M,a,\Lambda},\MorawetzOneForm_{M,a,\Lambda}}[\psi]\cdot N_{\SigmaStar}
           + \evalAt*{a\JCurrentPert{\widetilde{\MorawetzVF},\widetilde{\MorawetzLagrangeCorr}}[\psi](\tau)}^{\tau=\tau_2}_{\tau=\tau_1}.
    \end{split}    
  \end{equation*}
  From the bounds in \zcref[noname]{eq:JCurrentPert:basic-control}, we
  have that
  \begin{equation}
    \label{eq:Morawetz:KdS:boundary-bound-1}
    \abs*{\evalAt*{a\JCurrentPert{\widetilde{\MorawetzVF},\widetilde{\MorawetzLagrangeCorr}}[\psi](\tau)}^{\tau=\tau_2}_{\tau=\tau_1}}
    \lesssim a\sup_{\tau\in[\tau_1,\tau_2]}\EnergyHorizonDeg[\psi](\tau).
  \end{equation}
  Combining \zcref{eq:Morawetz:KdS:boundary-bound-1} with \zcref[cap]{coro:Morawetz:outside-trapping:boundary-terms}, we now have that
  \begin{equation*}
    \begin{split}
      \abs*{J_{\operatorname{Mor}}(\tau_1,\tau_2)}
      \lesssim{}& \sup_{\tau\in[\tau_1,\tau_2]}\EnergyHorizonDeg[\psi](\tau)
                  + \delta_{\Horizon}\SpacelikeFlux_{\mathcal{A}}[\psi](\tau_1,\tau_2)
                  + \delta_{\Horizon}\SpacelikeFlux_{\SigmaStar}[\psi](\tau_1,\tau_2)\\
                & + \int_{\mathcal{A}(\tau_1,\tau_2)}\abs*{\nabla_{4}\psi}^2 
                  + \int_{\SigmaStar(\tau_1,\tau_2)}\abs*{\nabla_{3}\psi}^2
                  .
    \end{split}    
  \end{equation*}
  With regards to the bulk terms, we have
  \begin{equation*}
    \begin{split}
      &\int_{\Manifold(\tau_1,\tau_2)}\KCurrent{\MorawetzVF_{M,a,\Lambda},\MorawetzLagrangeCorr_{M,a,\Lambda},\MorawetzOneForm_{M,a,\Lambda}}[\psi]
    + a\KCurrentPert{\widetilde{\MorawetzVF},\widetilde{\MorawetzLagrangeCorr}}[\psi](\tau_1,\tau_2)\\
      ={}& a\KCurrentPert{\widetilde{\MorawetzVF},\widetilde{\MorawetzLagrangeCorr}}_{\operatorname{princ}}[\psi](\tau_1,\tau_2)
           + a\KCurrentPert{\widetilde{\MorawetzVF},\widetilde{\MorawetzLagrangeCorr}}_{\bowtie}[\psi](\tau_1,\tau_2)
    + \int_{\Manifold(\tau_2,\tau_1)}\KCurrent{\MorawetzVF_{M,a,\Lambda},\MorawetzLagrangeCorr_{M,a,\Lambda},\MorawetzOneForm_{M,a,\Lambda}}[\psi]
    .
    \end{split}    
  \end{equation*}
  Then, using \zcref[noname]{eq:Morawetz:K-ext-bound} and
  \zcref[noname]{eq:Morawetz:cutoff-PsiDO:bulk}, we have that for $a$
  sufficiently small, there exists some $\delta_*$ such that 
  \begin{equation*}
    \begin{split}
      &\int_{\Manifold(\tau_1,\tau_2)}\KCurrent{\MorawetzVF_{M,a,\Lambda},\MorawetzLagrangeCorr_{M,a,\Lambda},\MorawetzOneForm_{M,a,\Lambda}}[\psi]
        + a\KCurrentPert{\widetilde{\MorawetzVF},\widetilde{\MorawetzLagrangeCorr}}[\psi](\tau_1,\tau_2)\\
      \ge{}& \delta_{*}\MorNormTrap[\psi](\tau_1,\tau_2)
             - O(a)\left( \norm*{\psi}_{L^2_{\operatorname{comp}}(\Manifold(\tau_1,\tau_2))}^2
             + \int_{\Manifold(\tau_1,\tau_2)}\abs*{N}^2
             + \sup_{[\tau_1,\tau_2]}\EnergyHorizonDeg[\psi](\tau) \right)
             .
    \end{split}    
  \end{equation*}
  Using \zcref[cap]{lemma:Morawetz:RHS:estimate} in
  $\Manifold_{\cancel{\operatorname{trap}}}(\tau_1,\tau_2)$ and 
  \zcref[cap]{lemma:Morawetz:RHS:PsiDO} in
  $\Manifold_{\operatorname{trap}}(\tau_1,\tau_2)$, as well as 
  \zcref[cap]{lemma:ILED-near:LoT-control:zero-order-mixed-term} to control
  the error terms, we now have that
  \begin{equation*}
    \begin{split}
      \left( \delta_{*}-\delta_2 - \delta_3\right)\MorNormTrap[\psi](\tau_1,\tau_2)
      \lesssim{}& \sup_{\tau\in [\tau_1,\tau_2]}\EnergyHorizonDeg[\psi](\tau)
      + \delta_{\Horizon}\SpacelikeFlux_{\mathcal{A}}[\psi](\tau_1,\tau_2)
                  + \delta_{\Horizon}\SpacelikeFlux_{\SigmaStar}[\psi](\tau_1,\tau_2)\\
      & + \int_{\mathcal{A}(\tau_1,\tau_2)}\abs*{\nabla_{4}\psi}^2 
                  + \int_{\SigmaStar(\tau_1,\tau_2)}\abs*{\nabla_{3}\psi}^2
      + a \MorNorm[\psi](\tau_1,\tau_2)
      + \ForcingTermNorm[\psi,N](\tau_1,\tau_2),
    \end{split}
  \end{equation*}
  where we used Cauchy Schwarz and the fact that
  \begin{equation*}
    \begin{split}
      &\int_{\Manifold_{\cancel{\operatorname{trap}}}(\tau_1,\tau_2)}\left( - \CovariantDeriv_{\mu}\left(\frac{2a\cos\theta}{\abs*{q}^2}\partial_r^{\mu}\mathcal{F}_{M,a,\Lambda}\psi\cdot\HorCovDeriv_{\KillT}\LeftDual{\psi}\right)
        + \partial_t\left(\frac{2a\cos\theta}{\abs*{q}^2}\mathcal{F}_{M,a,\Lambda}\psi\cdot\HorCovDeriv_r\LeftDual{\psi}\right) \right)\\
      \lesssim{} & \sup_{\tau\in [\tau_1,\tau_2]}\EnergyHorizonDeg[\psi](\tau)
                   + \delta_{\Horizon}\SpacelikeFlux_{\mathcal{A}}[\psi](\tau_1,\tau_2)
                   + \delta_{\Horizon}\SpacelikeFlux_{\SigmaStar}[\psi](\tau_1,\tau_2)\\
      & + \int_{\mathcal{A}(\tau_1,\tau_2)}\abs*{\nabla_{4}\psi}^2 
        + \int_{\SigmaStar(\tau_1,\tau_2)}\abs*{\nabla_{3}\psi}^2
        .
    \end{split}    
  \end{equation*}
  Then, for $\delta_2,\delta_3$, and $a$ sufficiently
  small, we have that 
  \begin{equation*}
    \begin{split}
      \MorNormTrap[\psi](\tau_1,\tau_2)
      \lesssim{}& \sup_{\tau\in [\tau_1,\tau_2]}\EnergyHorizonDeg[\psi](\tau)
                  + \delta_{\Horizon}\SpacelikeFlux_{\mathcal{A}}[\psi](\tau_1,\tau_2)
                  + \delta_{\Horizon}\SpacelikeFlux_{\SigmaStar}[\psi](\tau_1,\tau_2)\\
                & + \int_{\mathcal{A}(\tau_1,\tau_2)}\abs*{\nabla_{4}\psi}^2 
                  + \int_{\SigmaStar(\tau_1,\tau_2)}\abs*{\nabla_{3}\psi}^2
                  + \ForcingTermNorm[\psi,N](\tau_1,\tau_2),
    \end{split}
  \end{equation*}
  as desired, concluding the proof of \zcref{prop:Morawetz:KdS:main}.
\end{proof}


\section{$r^p$-weighted estimates}
\label{sec:rp}

In this section, we will prove the $r^p$-weighted estimates for
solutions to the model problem in
\zcref[noname]{eq:model-problem-gRW}, which play an equivalent role to
the redshift estimate near the event horizon for \KdS{} and which
degenerate to become the standard $r^p$-weighted estimate for Kerr in
the vanishing $\Lambda$ limit.

\begin{proposition}
  \label{prop:rp:Kerr}
  Let $\psi\in \realHorkTensor{2}$ be a solution to
  \zcref[noname]{eq:model-problem-gRW}, and $R\gg M$ be sufficiently
  large. Then we have that for $\delta\le p \le 2-\delta$,
  \begin{equation}
    \label{eq:rp:Kerr}
    \begin{split}
      \WeightedBEFNorm{p,\ge R}^s[\psi](\tau_1,\tau_2)
    \lesssim{}& \EnergyFluxWeightedOpt{p}{\frac{R}{2}}^s[\psi](\tau_1)
    + R^{p+2}\MorNorm_{\frac{R}{2}\le r \le R }^s[\psi](\tau_1,\tau_2)\\
              &  + \ForcingTermWeightedNorm{p}{\frac{R}{2}}^s[\psi, N](\tau_1,\tau_2)
                + \EnergyHorizonDeg^s[\psi](\tau_2). 
    \end{split}    
  \end{equation}
\end{proposition}

We briefly remark on the differences when comparing the $r^p$-weighted
estimates proven here for \KdS{} with the $r^p$-weighted estimates
typically proven for solutions to wave-type equations on
asymptotically flat spacetimes (see in particular Section 10.3 of
\cite{giorgiWaveEquationsEstimates2024} for the Kerr setting). As the
region where $\Lambda r^2\sim 1$ (the region close to the cosmological
horizon of \KdS) is where Kerr and \KdS{} display truly non-uniform
geometric features, it is precisely here where one expects to
encounter the greatest difficulty in proving a uniform
estimate. Indeed, in the ensuing estimates, we see that the
contribution of $\Lambda$ to the estimates is non-pertubative in the
sense that there is a top-order term proportional to
$\Lambda r^2$. Namely, the following de Sitter-type terms will appear: in the bulk term (see \zcref[noname]{eq:rp:bulk:specific}), we will have
    \begin{equation*}
      \Lambda r^2 r^{p-1}\left( \abs*{\widecheck{\nabla}_4\psi}^2 + \abs*{r^{-1}\psi}^2 \right),
    \end{equation*}   
    while in the flux term over $\SigmaStar$  (see \zcref[cap]{prop:rp:main-boundary-estimates}), we will have
  \begin{equation*}
    r^{p-2}(\Lambda r^2)\abs*{\psi}^2 +  \frac{-\Delta}{\abs*{q}^2}r^p\abs*{\widecheck{\nabla}_4\psi}^2.
  \end{equation*}
The fact that all of these terms are positive on their respective
regions of integration plays a key role in allowing us to achieve
uniform estimates in the far region of the spacetime. 

\subsection{Preliminaries}

Throughout this section, we will be working in a section of \KdS{}
where $r^2\Lambda\sim 1$,
thus the global principal null frame $(e_3,e_4)$ will coincide with the outgoing principal null frame
$(\eout_3,\eout_4)$. We first collect some properties of our $\tau$-foliation that will be used in what follows.

\begin{lemma}
  \label{lemma:rp:tau-foliation-properties}
  We have that
  \begin{align}
    N_{\Sigma} ={}& e_4 + \frac{1}{2}r^{-2}\lambda e_3 + Y^be_b,\qquad
                    Y = O(ar^{-1}),
                    \label{eq:n-Sigma:null-decomposition}\\
    \Metric(N_{\Sigma}, N_{\Sigma}) ={}& -2r^{-2}\lambda + \abs*{Y}^2
                                         \lesssim -\frac{M^2}{r^2} \label{eq:n-Sigma-nSigma-metric-value}.
  \end{align}
  with $\lambda$ satisfying $M^2\lesssim \lambda \lesssim M^2$, $D\lambda = O(R^{-1})$ and $e_4(\lambda) = O(M^2r^{-3})$ and with $Y$ satisfying $e_4(r^2\abs*{Y}^2) = O(M^2r^{-3})$.
\end{lemma}
\begin{proof}
  The lemma follows in view of the construction of $\tau$ in
  \zcref[cap]{prop:properties-of-global-coordinates}, and observing that
  $N_{\Sigma}$ as defined in \zcref{eq:hypersurface-normal-def} is
  normalized so that $\Metric(N_{\Sigma}, e_3)=-2$.
\end{proof}

\begin{definition}
  \label{def:rp:nu-Sigma:def}
  In addition, we define the vectorfield $\nu_{\Sigma}$
  \begin{equation}
    \label{eq:rp:nu-Sigma:def}
    \nu_{\Sigma}\vcentcolon= e_4 - \frac{1}{2}r^{-2}\lambda e_3,
  \end{equation}
  which is clearly tangent to $\Sigma(\tau)$ and satisfies
  \begin{equation}
    \label{eq:rp:nu-Sigma:abs}
    \Metric(\nu_{\Sigma}, \nu_{\Sigma}) = 2\lambda r^{-2}.
  \end{equation}
\end{definition}
We similarly define a tangent vectorfield to $\SigmaStar$.
\begin{definition}
  \label{def:rp:nu-SigmaStar:def}
  Define $\nu_{\SigmaStar}$ by
  \begin{equation}
    \label{eq:rp:nu-SigmaStar:def}
    \nu_{\SigmaStar}
    \vcentcolon=  \frac{1}{2}e_4(r)e_3
    - \frac{1}{2}e_3(r)e_4,
  \end{equation}
  which is clearly tangent to $\SigmaStar$ and moreover satisfies,
  \begin{equation}
    \label{eq:rp:nu-SigmaStar:abs}
    \Metric(\nu_{\SigmaStar}, \nu_{\SigmaStar}) = e_4(r)e_3(r) = \frac{-\Delta}{\abs*{q}^2},
  \end{equation}
  where we recall that on $\SigmaStar$, $-\Delta>0$.
\end{definition}

We first prove the following lemma computing some divergences.
\begin{lemma}
  \label{lemma:div-Sigma}
  We denote by $\Divergence_{\Sigma}$ the divergence corresponding to
  the renormalized volume form\footnote{Recall the convention we
    introduced in \zcref{eq:hypersurface-normal-def}}
  $\frac{1}{\sqrt{-\Metric(N_{\Sigma},N_{\Sigma})}}\volForm$ on
  $\Sigma(\tau)$. Then we have that
  \begin{equation}
    \label{eq:div-Sigma}
    \Divergence_\Sigma(\nu_\Sigma)=\Trace\chi
    +O(M^2r^{-3}).
  \end{equation}
\end{lemma}
\begin{proof}
  We define the auxiliary vectorfields $\widetilde{e}_a \vcentcolon= e_a - \frac{Y^a}{\abs*{N_{\Sigma}}^2_{\Metric}}N_{\Sigma}$ for $a=1,2$, so that
  $ \left( \nu_{\Sigma}, \widetilde{e}_1, \widetilde{e}_2 \right)$
  form a basis of the tangent space of $\Sigma(\tau)$. We then have that
  \begin{align}
    \gamma^{\Sigma}_{ab}
    \vcentcolon={}& \Metric(\widetilde{e}_a,\widetilde{e}_b)={} \delta_{ab} + \frac{Y^aY^b}{\abs*{N_{\Sigma}}_{\Metric}^2}, \label{eq:rp:div-computation:gamma-ab}\\
    \label{eq:rp:div-computation:gamma-inv-ab}
    (\gamma^{\Sigma})^{ab} ={}& \frac{\abs*{N_{\Sigma}}_{\Metric}^2}{\abs*{N_{\Sigma}}_{\Metric}^2 + \abs*{Y}^2}
                                \left( \delta^{ab} + \frac{\LeftDual{Y}^a\LeftDual{Y}^b}{\abs*{N_{\Sigma}}_{\Metric}^2} \right). 
  \end{align}
  Denoting by $\underline{\Divergence}_{\Sigma}$ the induced divergence on
  $\Sigma$ by $\Metric$ and using
  $\Metric(\nu_{\Sigma}, \widetilde{e}_a)=0$ we have that
  \begin{align*}
    \underline{\Divergence}_{\Sigma}(\nu_{\Sigma})
    ={}& \frac{1}{\abs*{\nu_{\Sigma}}_{\Metric}^2}\Metric(\CovariantDeriv_{\nu_{\Sigma}}\nu_{\Sigma}, \nu_{\Sigma})
         + (\gamma^{\Sigma})^{ab}\Metric(\CovariantDeriv_{\widetilde{e}_a}\nu_{\Sigma}, \widetilde{e}_b)
    ={} \frac{1}{2\abs*{\nu_{\Sigma}}_{\Metric}^2}\nu_{\Sigma}(\abs*{\nu_{\Sigma}}_{\Metric}^2)
         + (\gamma^{\Sigma})^{ab}\mathcal{K}_{ab},
  \end{align*}
  where
  \begin{align*}
    \mathcal{K}_{ab} &\vcentcolon=
    \Metric \left( \CovariantDeriv_{e_a-\frac{Y^a}{\abs*{N_{\Sigma}}_{\Metric}^2}N_{\Sigma}}\left(
        e_4 - \frac{1}{2}r^{-2}\lambda e_3\right),  e_b-\frac{Y^b}{\abs*{N_{\Sigma}}_{\Metric}^2}N_{\Sigma}
    \right)
  \\& = \Metric\left(\D_{e_a -\frac{Y^a}{\abs*{N_\Sigma}_{\Metric}^2}N_\Sigma}e_4, e_b -\frac{Y^b}{\abs*{N_\Sigma}_{\Metric}^2}N_\Sigma\right)%
  - \frac{1}{2}r^{-2}\la\Metric\left(\D_{e_a -\frac{Y^a}{\abs*{N_\Sigma}_{\Metric}^2}N_\Sigma}e_3, e_b -\frac{Y^b}{\abs*{N_\Sigma}_{\Metric}^2}N_\Sigma\right)\\
&+\frac{1}{2}\left(e_a -\frac{Y^a}{\abs*{N_\Sigma}_{\Metric}^2}N_\Sigma\right)(r^{-2}\la)\frac{Y^b}{\abs*{N_\Sigma}_{\Metric}^2}\Metric(e_3, N_\Sigma).
  \end{align*}
  Now, using \zcref{lemma:rp:tau-foliation-properties}, we have that
  \begin{align*}
    \mathcal{K}_{ab} 
    ={}& \chi_{ab}
         - \frac{Y^aY^c}{\abs*{N_\Sigma}_{\Metric}^2}\chi_{cb}
         -\frac{Y^bY^c}{\abs*{N_\Sigma}_{\Metric}^2}\chi_{ac}
         +\frac{Y^aY^b}{\abs*{N_\Sigma}_{\Metric}^4}Y^cY^d\chi_{cd}\\
       & + O(r^{-1})\eta
         + O(r^{-2})\underline{\omega}
         + O(r)\xi
         + O(r^{-2})\underline{\chi}
         + O(1)\omega
         + O(r^{-1})\underline{\eta}
    \\
       & -\left(e_a -\frac{Y^a}{\abs*{N_\Sigma}_{\Metric}^2}N_\Sigma\right)(r^{-2}\la)\frac{Y^b}{\abs*{N_\Sigma}_{\Metric}^2}
         +O(M^2r^{-3})
         .
  \end{align*}
  Using now the form of $\gamma^{ab}$ from
  \zcref{eq:rp:div-computation:gamma-inv-ab}, and the fact that
  $\gamma^{ab}= O(1)$, we can now compute that
  \begin{align*}
    \gamma^{ab}\mathcal{K}_{ab}
    ={}& \frac{\abs*{N_{\Sigma}}^2}{\abs*{N_{\Sigma}}^2 + \abs*{Y}^2} \left(
         \delta^{ab}+ \frac{\LeftDual{Y}^a\LeftDual{Y}^b}{\abs*{N_{\Sigma}}_{\Metric}^2} 
         \right)\\
       & \times \left(
         \chi_{ab}
         - \frac{Y^aY^c}{\abs*{N_\Sigma}_{\Metric}^2}\chi_{cb}
         -\frac{Y^bY^c}{\abs*{N_\Sigma}_{\Metric}^2}\chi_{ac}
         +\frac{Y^aY^b}{\abs*{N_\Sigma}_{\Metric}^4}Y^cY^d\chi_{cd}
         -\left(e_a -\frac{Y^a}{\abs*{N_\Sigma}_{\Metric}^2}N_\Sigma\right)(r^{-2}\la)\frac{Y^b}{\abs*{N_\Sigma}_{\Metric}^2}\right)\\
       & + O(r^{-1})\eta
         + O(r^{-2})\underline{\omega}
         + O(r)\xi
         + O(r^{-2})\underline{\chi}
         + O(1)\omega
         + O(r^{-1})\underline{\eta}
         + O(m^2r^{-3})\\
    ={}& \frac{\abs*{N_{\Sigma}}^2}{\abs*{N_{\Sigma}}^2 + \abs*{Y}^2} \left(
         \Trace\chi
         - \frac{\abs*{Y}^2}{2\abs*{N_\Sigma}_{\Metric}^2}\Trace\chi
         + \frac{\abs*{Y}^4}{2\abs*{N_{\Sigma}}_{\Metric}^4}\Trace \chi
         -\frac{1}{\abs*{N_{\Sigma}}_{\Metric}^2}\left(Y - \frac{\abs*{Y}^2}{\abs*{N_{\Sigma}}_{\Metric}^2}N_{\Sigma}\right)(r^{-2}\la)
         \right)
    \\
       & + O(r^{-1})\eta
         + O(r^{-2})\underline{\omega}
         + O(r)\xi
         + O(r^{-2})\underline{\chi}
         + O(1)\omega
         + O(r^{-1})\underline{\eta}
         + O(M^2r^{-3}).
  \end{align*}
  Now, recalling that
  $\abs*{N_{\Sigma}}_{\Metric}^2 = -2\lambda r^{-2} + \abs*{Y}^2$, we
  have that
  \begin{align*}
     \Trace\chi
         - \frac{\abs*{Y}^2}{2\abs*{N_\Sigma}_{\Metric}^2}\Trace\chi
    + \frac{\abs*{Y}^4}{2\abs*{N_{\Sigma}}_{\Metric}^4}\Trace \chi
    ={}& %
         \Trace\chi
         + \frac{\abs*{Y}^2  \left( \abs*{Y}^2- \abs*{N_{\Sigma}}_{\Metric}^2 \right)}{2\abs*{N_{\Sigma}}_{\Metric}^4}\Trace \chi\\
    ={}& \Trace\chi
         + \frac{\lambda \abs*{Y}^2}{r^2\abs*{N_{\Sigma}}_{\Metric}^4}\Trace \chi\\
    ={}& \Trace\chi
         + \frac{2\lambda \abs*{Y}^2}{r^3\abs*{N_{\Sigma}}_{\Metric}^4}
         + O(M^2r^{-3}).
  \end{align*}
  We can also compute that
  \begin{align*}
    -\frac{1}{\abs*{N_{\Sigma}}_{\Metric}^2} \left( Y - \frac{\abs*{Y}^2}{\abs*{N_{\Sigma}}_{\Metric}^2}N_{\Sigma} \right)\left( r^{-2}\lambda \right)
    ={}& -\frac{1}{\abs*{N_{\Sigma}}_{\Metric}^2} \left(
         Y - \frac{\abs*{Y}^2}{\abs*{N_{\Sigma}}_{\Metric}^2}\left( N_{\Sigma}-e_4 \right)
         \right)\left( r^{-2}\lambda \right)
         + \frac{\abs*{Y}^2}{\abs*{N_{\Sigma}}_{\Metric}^4}e_4\left( r^{-2}\lambda \right)\\
    ={}& -\frac{2\abs*{Y}^2}{r^3\abs*{N_{\Sigma}}_{\Metric}^4}\lambda
         + \frac{\abs*{Y}^2}{r^2\abs*{N_{\Sigma}}_{\Metric}^4}e_4\left(\lambda \right)
          -\frac{1}{\abs*{N_{\Sigma}}_{\Metric}^2} \left(
         Y - \frac{\abs*{Y}^2}{\abs*{N_{\Sigma}}_{\Metric}^2}\left( N_{\Sigma}-e_4 \right)
         \right)\left( r^{-2}\lambda \right).
  \end{align*}
  As a result, we have that
  \begin{align*}
    (\gamma^{\Sigma})^{ab}\mathcal{K}_{ab}
    ={}& \Trace\chi
         + O(1)e_4(\lambda)
         + O(r^2)Y(r^{-2}\lambda)
         + O(1)e_3(r^{-2}\lambda)\\
         &+ O(r^{-1})\eta
         + O(r^{-2})\underline{\omega}
         + O(r)\xi
         + O(r^{-2})\underline{\chi}
         + O(1)\omega
         + O(r^{-1})\underline{\eta}
           + O(M^2r^{-3})\\
    ={}& \Trace\chi
         + O(1)e_4(\lambda)
         + O(1)Y(\lambda)
         + O(r^{-2})e_3(\lambda)\\
       & + O(r)\xi
         + O(r^{-1})\eta
         + O(r^{-2})\underline{\omega}         
         + O(r^{-2})\underline{\chi}
         + O(1)\omega
         + O(r^{-1})\underline{\eta}
           + O(M^2r^{-3})
           .
  \end{align*}
  Then, using the control of $\lambda$ and $Y$ in \zcref{lemma:rp:tau-foliation-properties}, and the
  values for the Ricci coefficients from
  \zcref{lemma:Kerr:outgoing-PG:Ric-and-curvature}, we obtain $(\gamma^{\Sigma})^{ab}\mathcal{K}_{ab}
    = \Trace\chi + O(M^2r^{-3})$. As a result, we have that
  \begin{align*}
    \underline{\Divergence}_{\Sigma}(\nu_{\Sigma})
    ={}& \frac{1}{\abs*{\nu_{\Sigma}}_{\Metric}^2}\nu_{\Sigma}(\abs*{\nu_{\Sigma}}_{\Metric}^2)
         + \Trace\chi + O(M^2r^{-3}). 
  \end{align*}
  We now recall the relationship between
  $\Divergence_{\Sigma}(\nu_{\Sigma})$ and
  $\underline{\Divergence}_{\Sigma}(\nu_{\Sigma})$ is given by
  \begin{align*}
    \Divergence_{\Sigma}(\nu_{\Sigma})
    ={}& \sqrt{-\abs*{N_{\Sigma}}_{\Metric}^2}\,\underline{\Divergence}_{\Sigma}\left(\frac{\nu_{\Sigma}}{\sqrt{-\abs*{N_{\Sigma}}_{\Metric}^2}}\right)
    ={} \underline{\Divergence}_{\Sigma}(\nu_{\Sigma})
         - \frac{1}{2\abs*{N_{\Sigma}}_{\Metric}^2}\nu_{\Sigma} \left( \abs*{N_{\Sigma}}_{\Metric}^2 \right).
  \end{align*}
  Therefore, using \zcref{eq:rp:nu-Sigma:abs} and
  \zcref{eq:n-Sigma-nSigma-metric-value}, we have that
  \begin{align*}
    \Divergence_{\Sigma}(\nu_{\Sigma})
    ={}& \Trace\chi
         + \frac{1}{2}\left(e_4 - \frac{1}{2}r^{-2}\lambda e_3\right)\ln \left( 1 - \frac{r^2\abs*{Y}^2}{2\lambda} \right)
         + O(M^2r^{-3}) \\
    ={}& \Trace\chi
         + O(1)e_4(\lambda)
         + O(1)e_4(r^2\abs*{Y}^2) + O(r^{-2})e_3(\lambda)
         + O(r^{-2})e_3(r^2\abs*{Y}^2)
        + O(M^2r^{-3}) .
  \end{align*}
  Using again \zcref{lemma:rp:tau-foliation-properties}, we deduce that $\Divergence_{\Sigma}(\nu_{\Sigma})
    = \Trace\chi
    + O(M^2r^{-3})$ as desired.
\end{proof}

We have a similar lemma concerning the divergence of $\nu_{\SigmaStar}$. 
\begin{lemma}
  \label{lemma:rp:div-SigmaStar}
  We have that on $\SigmaStar$,
  \begin{equation}
    \label{eq:rp:div-SigmaStar}
    \Divergence_{\SigmaStar}\nu_{\SigmaStar}
    = \frac{1}{2}e_4(r)\Trace\underline{\chi}
    - \frac{1}{2}e_3(r)\Trace\chi. 
  \end{equation}
\end{lemma}
\begin{proof}
  Observe that $\left( \nu_{\SigmaStar}, e_1,e_2 \right)$ form a basis
  of the tangent space of $\Sigma(\tau)$ with
  $\Metric(\nu_{\SigmaStar}, e_a)=0$. Denoting by
  $\Divergence_{\SigmaStar}$ the induced divergence on $\SigmaStar$ by
  $\Metric$, we have that
  \begin{align*}
    \underline{\Divergence}_{\SigmaStar}(\nu_{\SigmaStar})
    ={}& \frac{1}{\abs*{\nu_{\SigmaStar}}_{\Metric}^2}\Metric(\CovariantDeriv_{\nu_{\SigmaStar}}\nu_{\SigmaStar}, \nu_{\SigmaStar})
         + \delta^{ab}\Metric(\CovariantDeriv_a\nu_{\SigmaStar}, e_b)\\
    ={}& \frac{1}{2\abs*{\nu_{\SigmaStar}}_{\Metric}^2}\nu_{\SigmaStar}(\abs*{\nu_{\SigmaStar}}_{\Metric}^2)
         + \frac{1}{2}e_4(r)\Trace \underline{\chi}
         - \frac{1}{2}e_3(r)\Trace \chi
         + O(r^{-2})\\
    ={}& \frac{1}{2}e_4(r)\Trace \underline{\chi}
         - \frac{1}{2}e_3(r)\Trace \chi
         .
  \end{align*}
  Next, $\Metric(N_{\SigmaStar}, N_{\SigmaStar}) = \frac{\Delta}{\abs*{q}^2}$ (see \zcref{eq:hypersurface-normal-def}) implies $\nu_{\SigmaStar}(\abs*{N_{\SigmaStar}}_{\Metric}^2) =0$. We can thus compute
  \begin{align*}
    \Divergence_{\SigmaStar}(\nu_{\SigmaStar})
    ={}& \sqrt{\abs*{N_{\SigmaStar}}_{\Metric}^2}\,\underline{\Divergence}_{\SigmaStar}\left(\frac{\nu_{\SigmaStar}}{\sqrt{\abs*{N_{\SigmaStar}}_{\Metric}^2}}\right)
    ={} \underline{\Divergence}_{\SigmaStar}(\nu_{\SigmaStar})
         - \frac{1}{\abs*{N_{\SigmaStar}}_{\Metric}^2}\nu_{\SigmaStar} \left( \abs*{N_{\SigmaStar}}_{\Metric}^2 \right)\\
    ={}& \frac{1}{2}e_4(r)\Trace \underline{\chi}
         - \frac{1}{2}e_3(r)\Trace \chi,
  \end{align*}
  as desired, proving \zcref{lemma:rp:div-SigmaStar}.
\end{proof}

We first prove the preliminary $r$-weighted estimate for $\nabla_3\psi$ in
the region $r\ge R$.
\begin{proposition}
  \label{prop:rp:preliminary-nabla-3-estimate}
  Let
  $M\ll R <
  r_{\CosmologicalHorizon}\left(1+\delta_{\Horizon}\right)$. Then we
  have that 
  \begin{equation*}
    \label{eq:rp:preliminary-nabla-3-estimate}
    \begin{split}
      \int_{\Manifold_{\ge R}(\tau_1,\tau_2)}r^{-1-\delta}\abs*{\nabla_3\psi}^2
      \lesssim{}& \int_{\Manifold_{\ge R}(\tau_1,\tau_2)}r^{-1-\delta}\left(
                  \abs*{\nabla_4\psi}^2 + \abs*{\nabla\psi}^2 + r^{-2}\abs*{\psi}^2
                  \right)\\
                & + \EnergyFlux_{\ge \frac{R}{2}}[\psi](\tau_1)
                  + R\MorNorm_{\frac{R}{2}\le r\le R}[\psi](\tau_1, \tau_2)
                  + \ForcingTermCombinedNorm{}[\psi, N](\tau_1,\tau_2)\\
                & + \EnergyHorizonDeg[\psi](\tau_2)
                  + \int_{\Sigma_{\ge \frac{R}{2}}(\tau_2)}a^2r^{-2-\delta}\abs*{\nabla_4\psi}^2
      \\         
                & + \int_{\SigmaStar(\tau_1,\tau_2)}\left( \frac{a}{\abs*{q}^2}r^{-\delta}\abs*{\nabla\psi}^2  + \frac{\Delta^2}{\abs*{q}^4}r^{-\delta}\abs*{\widecheck{\nabla}_4\psi}^2 \right),
    \end{split}
  \end{equation*}
  where $R$ is a fixed sufficiently large constant.
\end{proposition}

To prove \zcref[cap]{prop:rp:preliminary-nabla-3-estimate}, we
first prove the following auxiliary lemma.
\begin{lemma}
  \label{lemma:rp:preliminary-nabla-3-estimate:aux}
  In the region $r\ge R$, we have the following identity
  \begin{equation}
    \label{eq:rp:preliminary-nabla-3-estimate:aux:bulk}
    (1+\gamma)\KCurrent{f_{-\delta}\KillT, 0, 0}[\psi]
    = \frac{1}{4}\delta r^{-1-\delta}\abs*{\nabla_3\psi}^2
    - \frac{1}{4}\frac{\Delta^2}{\abs*{q}^4}\delta r^{-1-\delta}\abs*{\nabla_4\psi}^2
    + O(R^{-1})r^{-1-\delta}\left(\abs*{D\psi}^2 + r^{-2}\abs*{\psi}^2\right).
  \end{equation}
  Moreover,
  \begin{equation}
    \label{eq:rp:preliminary-nabla-3-estimate:aux:boundary:Sigma}
    \begin{split}
      \frac{1}{2}f_{-\delta}\lambda \left(
        \frac{1}{2}r^{-1}\nabla_3\psi + 2\lambda^{-1}r Y^b\nabla_b\psi
        \right)^2
      &\le
              f_{-\delta}\left(
             \frac{1}{r^2}\abs*{\nabla_{3}\psi}^2
             + \frac{\abs*{\Delta}}{r^2}\abs*{\nabla_{4}\psi}^2
             + \abs*{\nabla\psi}^2
             + r^{-2}\abs*{\psi}^2
             \right)\\
      & +  \frac{a^2}{r^2}f_{-\delta}\abs*{\nabla_4\psi}^2
        + \JCurrent{f_{-\delta}\KillT, 0, 0}[\psi]\cdot N_{\Sigma}
             ,
    \end{split}    
  \end{equation}
  and
  \begin{equation}
    \label{eq:rp:preliminary-nabla-3-estimate:aux:boundary:SigmaStar}    
    \JCurrent{f_{-\delta}\KillT, 0, 0}[\psi]\cdot N_{\SigmaStar}
    \ge{} c_0\abs*{\nabla_3\psi}^2
    - C_0\left(
      \frac{a}{\abs*{q}^2}f_{-\delta}\abs*{\nabla\psi}^2
      + \frac{\Delta^2}{\abs*{q}^4}f_{-\delta}\abs*{\widecheck{\nabla}_4\psi}^2
    \right).
  \end{equation}
\end{lemma}

\begin{proof}
  We consider the vectorfield multiplier $X=f(r)\KillT$. Since $\KillT$
  is Killing, we have $\DeformationTensor[_{\alpha\beta}]{X}
    = \CovariantDeriv_{(\alpha}f\KillT_{\beta)}$. Then recalling the expression for $\KillT$ in the principal outgoing
  frame in \zcref[noname]{eq:T-R-Z:principal-outgoing:expression:T-Z}, we have
  that
  \begin{align*}
    \EMTensor\cdot\DeformationTensor[]{X}    
    ={}& \EMTensor^{3\beta}\KillT_{\beta}e_3f
         +  \EMTensor^{4\beta}\KillT_{\beta}e_4f\\
    ={}& - \frac{1}{2}\left(
         \EMTensor(e_4,\KillT)e_3f
         +\EMTensor(e_3,\KillT)e_4f
         \right)\\
    ={}& -\frac{1}{2}\EMTensor\left(e_4,\frac{1}{2(1+\gamma)}\left(\frac{\Delta}{\abs*{q}^2}e_4 + e_3
         - \frac{a\sqrt{\kappa}\sin\theta}{\abs*{q}}e_2\right)\right)e_3f\\
       & - \frac{1}{2}\EMTensor\left(e_3,\frac{1}{2(1+\gamma)}\left(\frac{\Delta}{\abs*{q}^2}e_4 + e_3
         - \frac{a\sqrt{\kappa}\sin\theta}{\abs*{q}}e_2\right)\right)e_4f\\
    ={}& -\frac{1}{4} \frac{\Delta^2}{(1+\gamma)\abs*{q}^4}\partial_rf\abs*{\nabla_4\psi}^2
         + \frac{1}{4} \frac{\partial_rf}{1+\gamma} \abs*{\nabla_3\psi}^2\\
       & - \frac{a\sqrt{\kappa}\sin\theta}{2\abs*{q}}\left(
         -\frac{\Delta}{(1+\gamma)\abs*{q}^2}\partial_rf\nabla_4\psi\cdot\nabla_2\psi
         +\partial_rf\nabla_3\psi\cdot\nabla_2\psi
         \right) 
         .         
  \end{align*}
  As a result, using $\KCurrent{f_{-\delta}\KillT, 0, 0}[\psi]
    = \EMTensor\cdot\DeformationTensor[]{f_{-\delta}\KillT}
    - \frac{1}{2}f_{-\delta}\KillT(V)\abs*{\psi}^2$ we can calculate that in fact
  \begin{equation*}
    \begin{split}
      (1+\gamma)\KCurrent{f_{-\delta}\KillT, 0, 0}[\psi]
    ={}& \frac{1}{4} \frac{\Delta^2}{\abs*{q}^4}\partial_rf_{-\delta}\abs*{\nabla_4\psi}^2
    - \frac{1}{4} \partial_rf_{-\delta} \abs*{\nabla_3\psi}^2\\
    & - \frac{a\sin\theta}{2\abs*{q}}\left(
      -\frac{\delta}{\abs*{q}^2}\partial_rf_{-\delta}\nabla_4\psi\cdot\nabla_2\psi
      +\partial_rf_{-\delta}\nabla_3\psi\cdot\nabla_2\psi
    \right) .
    \end{split}    
  \end{equation*}
  Since we have that $f_{-\delta} = r^{-\delta}$ for $r\ge R$, we have
  that on $r\ge R$,
  \begin{equation*}
    \begin{split}
      (1+\gamma)\KCurrent{f_{-\delta}\KillT, 0, 0}[\psi]
      ={}& -\frac{s}{4}r^{-\delta-1} \frac{\Delta^2}{\abs*{q}^4}\abs*{\nabla_4\psi}^2
           + \frac{\delta}{4}r^{-\delta-1}\abs*{\nabla_3\psi}^2\\
         & + \frac{a\delta r^{-\delta-1}\sin\theta}{2\abs*{q}}\left(
           -\frac{\Delta}{\abs*{q}^2}\nabla_4\psi\cdot\nabla_2\psi
           +\nabla_3\psi\cdot\nabla_2\psi
           \right) ,
    \end{split}    
  \end{equation*}
  from which the inequality in
  \zcref[noname]{eq:rp:preliminary-nabla-3-estimate:aux:bulk} follows
  immediately using Cauchy-Schwarz. We now show
  \zcref[noname]{eq:rp:preliminary-nabla-3-estimate:aux:boundary:Sigma}. To
  begin, we can calculate explicitly that
  \begin{align*}
    \JCurrent{f_{-\delta}\KillT, 0, 0}[\psi]\cdot N_{\Sigma}
    ={}& f_{-\delta}\EMTensor\left(\frac{\Delta}{\abs*{q}^2}e_4+e_3 - a\frac{\sqrt{\kappa}\sin\theta}{\abs*{q}}e_2, e_4 + \frac{1}{2}r^{-2}\lambda e_3 + Y^be_b\right),\\
    ={}&
         \begin{aligned}[t]
           f_{-\delta}&\left(
         \frac{\Delta}{\abs*{q}^2}\EMTensor_{44}[\psi]
         + \frac{1}{2}r^{-2}\lambda\EMTensor_{33}[\psi]
         - a\frac{\sqrt{\kappa}\sin\theta}{\abs*{q}}Y^b\EMTensor_{2b}[\psi]
         + \left(\frac{\Delta\lambda}{2r^2\abs*{q}^2} + 1\right)\EMTensor_{34}[\psi]\right.\\
         & \left.\quad + Y^b\EMTensor_{3b}[\psi]
           + \frac{\Delta}{\abs*{q}^2}Y^b\EMTensor_{4b}[\psi]
           - a\frac{\sqrt{\kappa}\sin\theta}{\abs*{q}}\EMTensor_{24}[\psi]
         \right).
         \end{aligned}         
  \end{align*}
  We can expand this as
  \begin{equation}
    \label{eq:rp:preliminary-nabla-3-estimate:aux:first-J-current}
    \begin{split}
      \frac{1}{f_{-\delta}}\JCurrent{f_{-\delta}\KillT, 0, 0}[\psi]\cdot N_{\Sigma}
      ={}& \frac{\Delta}{\abs*{q}^2}\left(
           \nabla_4\psi + \frac{1}{2}Y^b\nabla_b\psi
           \right)^2
           - \frac{1}{4}\abs*{Y^b\nabla_b\psi}^2
           + \frac{1}{2}\lambda \left(
           \frac{1}{2}r^{-1}\nabla_3\psi + 2\lambda^{-1}r Y^b\nabla_b\psi
           \right)^2\\
         & - \frac{2r^2}{\lambda^2}\abs*{Y^b\nabla_b\psi}^2
           + \frac{1}{4}\lambda r^{-2}\abs*{\nabla_3\psi}^2
           - a \frac{\sqrt{\kappa}\sin\theta}{\abs*{q}}Y^b\EMTensor_{2b}[\psi]
           + \left(\frac{\Delta\lambda}{2r^2\abs*{q}^2} + 1\right)\EMTensor_{34}[\psi]\\
      & - a\frac{\sqrt{\kappa}\sin\theta}{\abs*{q}}\EMTensor_{24}[\psi]
      .
    \end{split}    
  \end{equation}
  We can compute that
  \begin{equation*}
    \begin{split}
      a \frac{\sqrt{\kappa}\sin\theta}{\abs*{q}}Y^b\EMTensor_{2b}
    ={}& a \frac{\sqrt{\kappa}\sin\theta}{\abs*{q}}\left(
      Y^1\nabla_1\psi\cdot\nabla_2\psi
      + Y^2\left(
        \nabla_3\psi\cdot\nabla_4\psi
        - \frac{1}{2}\abs*{\nabla_1\psi}^2
        + \frac{1}{2}\abs*{\nabla_2\psi}^2
      \right)
         \right)\\
      \lesssim{}& a^2 r^{-2}\abs*{\nabla\psi}^2
                  + a^3r^{-2}\abs*{\nabla_4\psi}^2 
                  + a^3r^{-2}\abs*{\nabla_3\psi}^2
                  .
    \end{split}    
  \end{equation*}
  Similarly,
  \begin{equation*}
    \frac{1}{4}\abs*{Y^b\nabla_b\psi}^2
    + \frac{r^2}{2\lambda^2}\abs*{Y^b\nabla_b\psi}^2
    + O\left(a r^{-2}\right)\abs*{\nabla\psi}^2
    \le r^{-2}\abs*{\nabla\psi}^2,
  \end{equation*}
  and
  \begin{equation*}
    \abs*{a\frac{\sqrt{\kappa}\sin\theta}{\abs*{q}}\EMTensor_{24}[\psi]}
    \lesssim a^2r^{-2}\abs*{\nabla_4\psi}^2
    + \abs*{\nabla\psi}^2.
  \end{equation*}
  Thus, we have that 
  \begin{equation*}
    \begin{split}
      \frac{1}{2}f_{-\delta}\lambda \left(
        \frac{1}{2}r^{-1}\nabla_3\psi + 2\lambda^{-1}r Y^b\nabla_b\psi
        \right)^2
      \le{}& 
              f_{-\delta}\left(
             \frac{1}{r^2}\abs*{\nabla_{3}\psi}^2
             + \frac{\abs*{\Delta}}{r^2}\abs*{\nabla_{4}\psi}^2
             + \abs*{\nabla\psi}^2
             + r^{-2}\abs*{\psi}^2
             \right)\\
      & +  \frac{a^2}{r^2}f_{-\delta}\abs*{\nabla_4\psi}^2
        + \JCurrent{f_{-\delta}\KillT, 0, 0}[\psi]\cdot N_{\Sigma}
             ,
    \end{split}    
  \end{equation*}
  as desired, proving
  \zcref[noname]{eq:rp:preliminary-nabla-3-estimate:aux:boundary:Sigma}. Similarly, we can compute that
  \begin{align*}
    \JCurrent{f_{-\delta}\KillT, 0, 0}[\psi]\cdot N_{\SigmaStar}
    = {}& f_{-\delta}\EMTensor\left(\frac{\Delta}{\abs*{q}^2}e_4+e_3 - a\frac{\sqrt{\kappa}\sin\theta}{\abs*{q}}e_2, \frac{1}{2}e_3 + \frac{-\Delta}{2\abs*{q}^2}e_4\right)\\
    ={}& f_{-\delta}\left(
         \frac{-\Delta^2}{2\abs*{q}^4}\EMTensor_{44}[\psi]
         + \frac{1}{2}\EMTensor_{33}[\psi]
         - a\frac{\sqrt{\kappa}\sin\theta}{2\abs*{q}}\left(
         \frac{-\Delta}{\abs*{q}^2}\EMTensor_{24}[\psi]
         + \EMTensor_{23}[\psi]
         \right)\right).
  \end{align*}
  Then, we can observe that there exists some $C>0$ such that 
  \begin{align*}
    \abs*{a\frac{\sqrt{\kappa}\sin\theta}{2\abs*{q}}\left(\frac{-\Delta}{\abs*{q}^2}\EMTensor_{24}[\psi] + \EMTensor_{23}[\psi]\right)}
    \lesssim{}&  \frac{-a\Delta\sqrt{\kappa}}{\abs*{q}^2}\left(\abs*{q}^{-2}\abs*{\nabla_2\psi}^2 + \abs*{\nabla_4\psi}^2\right)
     + a\sqrt{\kappa}\left(\abs*{q}^{-2}\abs*{\nabla_2\psi}^2 +  \abs*{\nabla_3\psi}^2\right).
  \end{align*}
  Then, for $a$ sufficiently small, we have %
  that there exists some $c_0>0$ potentially small and $C_0>0$
  potentially large such that
  \begin{equation*}
    \JCurrent{f_{-\delta}\KillT, 0, 0}[\psi]\cdot N_{\SigmaStar}
    \ge %
    c_0\abs*{\nabla_3\psi}^2
    - C_0\left(
      \frac{a}{\abs*{q}^2}f_{-\delta}\abs*{\nabla\psi}^2
      + \frac{\Delta^2}{\abs*{q}^4}f_{-\delta}\abs*{\widecheck{\nabla}_4\psi}^2
    \right),
  \end{equation*}
  as desired.   
\end{proof}

We are now ready to prove \zcref[cap]{prop:rp:preliminary-nabla-3-estimate}.
\begin{proof}[Proof of Proposition \ref{prop:rp:preliminary-nabla-3-estimate}]
  Using the divergence theorem in \zcref[noname]{eq:div-thm:general}, we have that
  \begin{equation*}
    \begin{split}
      &\int_{\Sigma(\tau_2)}\JCurrent{f_{-\delta}\KillT, 0, 0}[\psi]\cdot N_{\Sigma}
      + \int_{\SigmaStar}\JCurrent{f_{-\delta}\KillT, 0, 0}[\psi]\cdot N_{\SigmaStar}\\
        \lesssim{}& \int_{\Sigma(\tau_1)}\JCurrent{f_{-\delta}\KillT, 0, 0}[\psi]\cdot N_{\Sigma}
                    + \int_{\Manifold_{r\ge \frac{R}{2}}(\tau_1,\tau_2)}\KCurrent{f_{-\delta}\KillT, 0, 0}[\psi]
                    + \ForcingTermNorm[\psi, N](\tau_1,\tau_2).
    \end{split}
  \end{equation*}
  Then, we see from \zcref[noname]{eq:rp:preliminary-nabla-3-estimate:aux:bulk} that
  \begin{equation*}
    \begin{split}
      \int_{\Manifold_{r\ge R}(\tau_1,\tau_2)}r^{-1-\delta}\abs*{\nabla_3\psi}^2
    \lesssim{}& \int_{\Manifold_{r\ge \frac{R}{2}}(\tau_1,\tau_2)}\KCurrent{f_{-\delta}\KillT, 0, 0}[\psi]    
    + R\MorNorm_{\frac{R}{2}\le r\le R}[\psi](\tau_1,\tau_2)\\
    & + \int_{\Manifold_{r\ge R}(\tau_1,\tau_2)}r^{-1-\delta}\left(\abs*{\nabla_4\psi}^2 + \abs*{\nabla\psi}^2 + r^{-2}\abs*{\psi}^2\right).
    \end{split}    
  \end{equation*}
  Combining this with the inequalities in
  \zcref[noname]{eq:rp:preliminary-nabla-3-estimate:aux:boundary:Sigma} and
  \zcref[noname]{eq:rp:preliminary-nabla-3-estimate:aux:boundary:SigmaStar}
  conclude the proof of \zcref[cap]{prop:rp:preliminary-nabla-3-estimate}.
\end{proof}

\subsection{Proof of \zcref[cap]{prop:rp:Kerr} for \texorpdfstring{$s=0$}{s=0}}

To prove \zcref[cap]{prop:rp:Kerr}, we consider the choices
\begin{equation}
  \label{eq:rp-Kerr:multiplier-choice:general}
  \begin{gathered}
    \rpVF = f(r)\left(e_4 + \frac{1}{2r^2}\lambda e_3\right),\quad
    \rpLagrangianCorr = \frac{r}{\abs*{q}^2}f(r), \quad
    \rpOneForm = \frac{r}{\abs*{q}^2}\partial_rf(r)e_4 + \rpOneForm_{\Lambda},\quad
    \rpOneForm_{\Lambda}= \frac{2-p}{3}\Lambda f(r)e_4,
  \end{gathered}
\end{equation}
where $\lambda$ is as defined in
\zcref[noname]{eq:n-Sigma:null-decomposition}, and where $f_p$ is a
non-negative function that is equal to $r^p$ for $r\ge R$, and zero
for $r\le \frac{R}{2}$, where $R$ is a fixed and sufficiently large
constant.  We remark that $\rpOneForm_{\Lambda}$ is a top order term
in $\rpOneForm$, i.e. it is nonperturbative.

\begin{remark}
  In view of \zcref[noname]{eq:n-Sigma:null-decomposition}, we have that $\rpVF = f\left(N_{\Sigma} - Y\right)$.
\end{remark}

We again recall that throughout this section, the global principal
null frame $(e_3,e_4)$ coincides with the principal outgoing
null frame.

\paragraph{Bulk Terms.}

We begin with the following lemma regarding the bulk terms that
result when applying the divergence theorem with the multiplier
choices in \zcref[noname]{eq:rp:multiplier-choice:specific}.
\begin{lemma}
  \label{lemma:rp:bulk:basic-computations}
  The following hold true.
  \begin{enumerate}
  \item The deformation tensor of the vectorfield
    $X_1 = f(r)e_4$ is given by
    \begin{gather*}
      \DeformationTensor{X_1}_{44} = 0,\qquad
      \DeformationTensor{X_1}_{43} = \left(e_4f\right)\Metric_{34},\qquad
      \DeformationTensor{X_1}_{33} = -4f\partial_r\left(\frac{\Delta}{\abs*{q}^2}\right) - 4e_3f\\
      \DeformationTensor{X_1}_{4a} = 0,\qquad
      \DeformationTensor{X_1}_{3a} = 0,\qquad
      \DeformationTensor{X_1}_{ab} = 2f\left( \frac{1}{2}\frac{2r}{\abs*{q}^2}\Metric_{ab}\right),
    \end{gather*}

  \item The deformation tensor of the vectorfield
    $\DeformationTensor{X_2} =
    \frac{1}{2}r^{-2}\lambda fe_3$ is given by
    \begin{gather*}
    \DeformationTensor{X_2}_{33} = 0,\qquad
    \DeformationTensor{X_2}_{44} = -4r^{-2}\lambda f\omega - 2e_4\left(r^{-2}\lambda f\right),\\
    \DeformationTensor{X_2}_{3a} = r^{-2}f\lambda\xiBar_a,\qquad
    \DeformationTensor{X_2}_{4a} = r^{-2}f\lambda\left(\etaBar_a + \zeta_a\right),\qquad
    \DeformationTensor{X_2}_{3a} = r^{-2}f\lambda\xiBar_a,\\
    \DeformationTensor{X_2}_{ab} = r^{-2}\lambda f\left(\chiBarTF_{ab} + \frac{1}{2}\Trace\chiBar \Metric_{ab}\right),\qquad
    \DeformationTensor{X_2}_{34} = -2r^{-2}\lambda f\omegaBar - 2e_3\left(r^{-2}\lambda f\right),
  \end{gather*}

  \item For $\rpLagrangianCorr = \frac{r}{\abs*{q}^2}f$,
    \begin{equation}
      \label{eq:rp:bulk:Lagrangian-corrector}
      \begin{split}
       \ScalarWaveOp[\Metric]\rpLagrangianCorr
        ={}& \left(r^{-1} - \frac{r}{3}\Lambda\right)\partial_{rr}f
             - \frac{2\Lambda}{3}\partial_rf
             + \frac{2\Lambda}{3r}f\\
           & + \left(O(r^{-2}) + a^2O(\Lambda r^{-1})\right)\partial_{rr}f
             +\left( O(r^{-3}) + a^2O(r^{-2}\Lambda)\right)\partial_rf\\
           & + \left(O(r^{-4}) + a^2O(r^{-3}) + a^2O(r^{-1}\Lambda)\right)f.
      \end{split}
    \end{equation}
  \item For
    $\rpOneForm = \frac{r}{\abs*{q}^2}\partial_rfe_4 + 
    \rpOneForm_{\Lambda}$, where $\rpOneForm_{\Lambda} = \frac{2-p}{3}\Lambda f e_4$ we have that
    \begin{equation}
      \label{eq:rp:bulk:one-form}
      \begin{split}
        \Divergence\left(\abs*{\psi}^2\rpOneForm\right)
        ={}&2r^{-1}\partial_rf \nabla_4\psi\cdot\psi
             + \left(
             r^{-2}\partial_rf
             + r^{-1}\partial_{rr}f
             \right)\abs*{\psi}^2\\
           &+ (2-p)\Lambda f\nabla_4\psi\cdot\psi
             +(2-p)\Lambda\left(
             \partial_rf + 2r^{-1}f
             \right)\abs*{\psi}^2\\         
           & + a^2O(r^{-3})\partial_rf\left(\abs*{\nabla_4\psi}^2 + \abs*{\psi}^2\right)\\
           & + a^2O\left(r^{-2}\right)\partial_{rr}f\abs*{\psi}^2
             + a^2O(r^{-2})f \Lambda \abs*{\psi}^2.
      \end{split}      
    \end{equation}
  \end{enumerate}
\end{lemma}
\begin{proof}
  Using the fact that $\DeformationTensor{X_1}_{\mu\nu}
    = \Metric\left(\CovariantDeriv_{\mu}X_1, e_{\nu}\right)
    + \Metric\left(\CovariantDeriv_{\nu}X_1, e_{\mu}\right)$ we have that
  \begin{gather*}
    \DeformationTensor{X_1}_{44} = 0,\qquad
    \DeformationTensor{X_1}_{43} = -2e_4f + 4f\omega,\qquad
    \DeformationTensor{X_1}_{33} = -8f\omegaBar - 4e_3f\\
    \DeformationTensor{X_1}_{4a} = 2f\xi_a,\qquad
    \DeformationTensor{X_1}_{3a} = 2f\left(\eta + \zeta\right)_a,\qquad
    \DeformationTensor{X_1}_{ab} = 2f\left(\chiTF_{ab} + \frac{1}{2}\Trace\chi\Metric_{ab}\right).
  \end{gather*}
  Plugging in the exact values in the principal outgoing
  frame from \zcref[cap]{lemma:Kerr:outgoing-PG:Ric-and-curvature}, we
  get
  \begin{gather*}
    \DeformationTensor{X_1}_{44} = 0,\qquad
    \DeformationTensor{X_1}_{43} = -2e_4f,\qquad
    \DeformationTensor{X_1}_{33} = -4f\partial_r\left(\frac{\Delta}{\abs*{q}^2}\right) - 4e_3f\\
    \DeformationTensor{X_1}_{4a} = 0,\qquad
    \DeformationTensor{X_1}_{3a} = 0,\qquad
    \DeformationTensor{X_1}_{ab} = 2f\left( \frac{1}{2}\frac{2r}{\abs*{q}^2}\Metric_{ab}\right),
  \end{gather*}
  as desired. Similarly, we have that
  \begin{gather*}
    \DeformationTensor{X_2}_{33} = 0,\qquad
    \DeformationTensor{X_2}_{44} = -4r^{-2}\lambda f\omega - 2e_4\left(r^{-2}\lambda f\right),\\
    \DeformationTensor{X_2}_{3a} = r^{-2}f\lambda\xiBar_a,\qquad
    \DeformationTensor{X_2}_{4a} = r^{-2}f\lambda\left(\etaBar_a + \zeta_a\right),\qquad
    \DeformationTensor{X_2}_{3a} = r^{-2}f\lambda\xiBar_a,\\
    \DeformationTensor{X_2}_{ab} = r^{-2}\lambda f\left(\chiBarTF_{ab} + \frac{1}{2}\Trace\chiBar \Metric_{ab}\right),\qquad
    \DeformationTensor{X_2}_{34} = -2r^{-2}\lambda f\omegaBar - 2e_3\left(r^{-2}\lambda f\right).
  \end{gather*}
  Plugging in the exact values in the principal outgoing
  frame from \zcref{lemma:Kerr:outgoing-PG:Ric-and-curvature}, we get
  \begin{gather*} 
    \DeformationTensor{X_2}_{33} = 0,\qquad
    \DeformationTensor{X_2}_{44}
    =  - 2e_4\left(r^{-2}\lambda f\right),\\
    \DeformationTensor{X_2}_{3a}
    = 0,\qquad
    \DeformationTensor{X_2}_{4a} = 0,\qquad
    \DeformationTensor{X_2}_{3a} = 0,\\
    \DeformationTensor{X_2}_{ab}
    = - \lambda f \frac{\Delta}{r\abs*{q}^4} \Metric_{ab},\qquad
    \DeformationTensor{X_2}_{34}
    = -2r^{-2}\lambda f\partial_r\left(\frac{\Delta}{\abs*{q}^2}\right)
    - 2e_3\left(r^{-2}\lambda f\right),
  \end{gather*}
  as desired. We now move on to showing \zcref[noname]{eq:rp:bulk:Lagrangian-corrector} observe that
  \begin{align*}
    \ScalarWaveOp[\Metric]\rpLagrangianCorr
    ={}& \frac{1}{\sqrt{\abs*{\Metric}}}\partial_{\alpha}\left(\sqrt{\abs*{\Metric}} \Metric^{\alpha\beta}\partial_{\beta}\right)\rpLagrangianCorr
    ={} \frac{1}{\sqrt{\abs*{\Metric}}}\partial_{r}\left(\sqrt{\abs*{\Metric}} \Metric^{rr}\partial_{r}\right)\rpLagrangianCorr
         + \frac{1}{\sqrt{\abs*{\Metric}}}\partial_{\theta}\left(\sqrt{\abs*{\Metric}} \Metric^{\theta\theta}\partial_{\theta}\right)\rpLagrangianCorr.
  \end{align*}
  We can calculate that for the choice of $\rpLagrangianCorr$ made in
  \zcref[noname]{eq:rp-Kerr:multiplier-choice:general},
  \begin{align*}
    \frac{1}{\sqrt{\abs*{\Metric}}}\partial_{r}\left(\sqrt{\abs*{\Metric}} \Metric^{rr}\partial_{r}\right)\rpLagrangianCorr
    ={}&\frac{r\Delta}{\abs*{q}^4}\partial_{rr}f
         + 2\frac{\Delta}{\abs*{q}^2}\partial_r\left(\frac{r}{\abs*{q}^2}\right)\partial_rf
         + \frac{\Delta}{\abs*{q}^2}\partial_{rr}\left(\frac{r}{\abs*{q}^2}\right)f\\
       & + \frac{\partial_r\Delta}{\abs*{q}^2}\left(
         \partial_r\left(\frac{r}{\abs*{q}^2}\right)f
         + \frac{r}{\abs*{q}^2}\partial_rf
         \right)\\
    ={}& \frac{r\Delta}{\abs*{q}^4}\partial_{rr}f
         + \left(
         \frac{2\Delta}{\abs*{q}^2}\partial_r\left(\frac{r}{\abs*{q}^2}\right)
         + \frac{2\partial_r\Delta}{\abs*{q}^2}\frac{r}{\abs*{q}^2}
         \right)\partial_rf\\
       & + \left(
         \frac{\Delta}{\abs*{q}^2}\partial_{rr}\left(\frac{r}{\abs*{q}^2}\right)
         + \frac{2\partial_r\Delta}{\abs*{q}^2}\partial_r\left(\frac{r}{\abs*{q}^2}\right)
         \right)f\\
         \frac{1}{\sqrt{\abs*{\Metric}}}\partial_{\theta}\left(\sqrt{\abs*{\Metric}} \Metric^{\theta\theta}\partial_{\theta}\rpLagrangianCorr\right)
    ={}& \frac{2r\Delta}{\abs*{q}^2\sin\theta}\partial_{\theta}\left(\sin\theta \partial_{\theta}\abs*{q}^{-2}\right)f.     
  \end{align*}
  Observing that $\partial_r\abs*{q} = \frac{r}{\abs*{q}}$ and $\partial_r\left(\frac{r}{\abs*{q}^2}\right)
    = \frac{-r^2+a^2\cos^2\theta}{\abs*{q}^4}$ we obtain
  \begin{align*}
    \ScalarWaveOp[\Metric]\rpLagrangianCorr
    ={}& \frac{r\Delta}{\abs*{q}^4}\partial_{rr}f
    + \left(
      \frac{2\Delta}{\abs*{q}^2}\partial_r\left(\frac{r}{\abs*{q}^2}\right)
      + \frac{\partial_r\Delta}{\abs*{q}^2}\frac{r}{\abs*{q}^2}
         \right)\partial_rf\\
       & + \left(
         \frac{\Delta}{\abs*{q}^2}\partial_{rr}\left(\frac{r}{\abs*{q}^2}\right)
         + \frac{\partial_r\Delta}{\abs*{q}^2}\partial_r\left(\frac{r}{\abs*{q}^2}\right)
         \right)f
         + \frac{r\Delta}{\abs*{q}^2\sin\theta}\partial_{\theta}\left(\sin\theta \partial_{\theta}\abs*{q}^{-2}\right)f.
  \end{align*}
  We can then compute that 
  \begin{align*}
    \frac{r\Delta}{\abs*{q}^4}
    ={}& \frac{r^3}{\abs*{q}^4}
    - \frac{r^5}{3\abs*{q}^4}\Lambda
    + O\left(r^{-2}\right)
    + O\left(a^2 \Lambda r^{-1}\right),\\
    \frac{2\Delta}{\abs*{q}^2}\partial_r\left(\frac{r}{\abs*{q}^2}\right)
    + \frac{\partial_r\Delta}{\abs*{q}^2}
    ={}& - \frac{2r^5}{3\abs*{q}^5}\Lambda
         + O\left(r^{-3}\right)
         + O\left(a^2r^{-2}\Lambda\right),\\
    \frac{\Delta}{\abs*{q}^2}\partial_{rr}\left(\frac{r}{\abs*{q}^2}\right)
    + \frac{\partial_r\Delta}{\abs*{q}^2}\partial_r\left(\frac{r}{\abs*{q}^2}\right)
    ={}& \frac{2r^{5}}{3\abs*{q}^6}\Lambda
         + O\left(r^{-4}\right)
         + O\left(a^2r^{-3}\Lambda\right).
  \end{align*}
  We can also compute that $\frac{r\Delta}{\abs*{q}^2\sin\theta}\partial_{\theta}\left(\sin\theta \partial_{\theta}\abs*{q}^{-2}\right)
    = a^2O(r^{-3})f + a^2O(r^{-1})\Lambda f$. As a result, we have that
  \begin{align*}
    \ScalarWaveOp[\Metric]\rpLagrangianCorr
    ={}& \left(\frac{r^3}{\abs*{q}^4}
         - \frac{r^5}{3\abs*{q}^4}\Lambda\right)\partial_{rr}f
         - \frac{2r^5}{3\abs*{q}^5}\Lambda\partial_rf
         +   \frac{2r^5}{3\abs*{q}^6}\Lambda f\\
       & + \left(O(r^{-2}) + a^2O(\Lambda r^{-1})\right)\partial_{rr}f
         +\left( O(r^{-3}) + a^2O(r^{-2}\Lambda)\right)\partial_rf \\
       & + \left(O(r^{-4}) + a^2\left(O(r^{-3})+O(r^{-1}\Lambda)\right)\right)f\\
    ={}& \left(r^{-1} - \frac{r}{3}\Lambda\right)\partial_{rr}f
         - \frac{2\Lambda}{3}\partial_rf
         + 2\frac{\Lambda}{3r}f\\
       & + \left(O(r^{-2}) + a^2O(\Lambda r^{-1})\right)\partial_{rr}f
         +\left( O(r^{-3}) + a^2O(r^{-2}\Lambda)\right)\partial_rf\\
       & + \left(O(r^{-4}) + a^2\left(O(r^{-3})+O(r^{-1}\Lambda)\right)\right)f
         .
  \end{align*}
  Finally, we move on to proving \zcref[noname]{eq:rp:bulk:one-form}. Observe
  that $\Divergence e_4 = \frac{2r}{\abs*{q}^2}
    = 2r^{-1}
    + a^2O\left(r^{-3}\right)$. As a result, we have that
  \begin{align}
    \Divergence\left(\abs*{\psi}^2 \frac{r}{\abs*{q}^2}\partial_rf e_4\right)
    ={}& \frac{2r}{\abs*{q}^2}\partial_rf \nabla_4\psi\cdot \psi
    + \left(\partial_r\left(
        \frac{r}{\abs*{q}^2}
      \right)\partial_rf
      + \frac{r}{\abs*{q}^2}\partial_{rr}f
         + \frac{r}{\abs*{q}^2}\partial_rf\Divergence e_4
         \right)\abs*{\psi}^2\notag \\
    ={}& 2r^{-1}\partial_rf \nabla_4\psi\cdot\psi
         + \left(
         r^{-2}\partial_rf
         + r^{-1}\partial_{rr}f
         \right)\abs*{\psi}^2\notag\\
       & + a^2O(r^{-3})\partial_rf\left(\abs*{\nabla_4\psi}^2 + \abs*{\psi}^2\right)
         + a^2O\left(r^{-2}\right)\partial_{rr}f\abs*{\psi}^2. \label{eq:rp:bulk:one-form:1}
  \end{align}
  We also have
  \begin{align}
    \Divergence\left(\abs*{\psi}^2(2-p)\Lambda r^{p}e_4\right)
    ={}& (2-p)\Lambda f\nabla_4\psi\cdot\psi
         +(2-p)\Lambda\left(
         \partial_rf + 2r^{-1}f
         \right)\abs*{\psi}^2
         + a^2O(r^{-2})f \Lambda \abs*{\psi}^2. \label{eq:rp:bulk:one-form:2}
  \end{align}
  Combining \zcref[noname]{eq:rp:bulk:one-form:1} and
  \zcref[noname]{eq:rp:bulk:one-form:2} yields \zcref[noname]{eq:rp:bulk:one-form} as
  desired.
\end{proof}

We have the following proposition.
\begin{proposition}
  \label{prop:rp:bulk}
  Fix $\delta>0$ and $\delta\le p \le 2-\delta$. Then
  the following statements hold.
  \begin{enumerate}
  \item Let $\left(\rpVF,\rpLagrangianCorr,\rpOneForm\right)$ be as
    chosen in \zcref[noname]{eq:rp-Kerr:multiplier-choice:general}. Then,
    \begin{equation}
      \label{eq:rp:bulk:general}
      \KCurrent{\rpVF, \rpLagrangianCorr, \rpOneForm}[\psi]
      = \frac{1}{2}\partial_rf\abs*{\widecheck{\nabla}_4\psi}^2
      + \frac{1}{2}\left(\frac{2}{r}f - \partial_rf\right)\left(
        \abs*{\nabla\psi}^2 + V\abs*{\psi}^2
      \right)
      + \KCurrent{\rpVF, \rpLagrangianCorr, \rpOneForm}_{(\Lambda)}[\psi] + Err,
    \end{equation}
    where %
    \begin{equation*}
      \KCurrent{\rpVF, \rpLagrangianCorr, \rpOneForm}_{(\Lambda)}[\psi]
      = \frac{\Lambda}{3}\left(2rf-r^2\partial_rf\right)\abs*{\nabla_4\psi}^2 
          + \frac{\Lambda}{3}\left(r\partial_{rr}f + 2\partial_rf - 2fr^{-1}\right)\abs*{\psi}^2,
    \end{equation*}
    and
    \begin{align*}
      Err ={}& O\left(r^{-3}f + r^{-2}\partial_rf\right)\abs*{\nabla_3\psi}^2
               + \left(O\left(r^{-3}\right) + O\left(r^{-1}\right)\Lambda\right)f\abs*{\nabla_4\psi}^2\\
             & + \left(O\left(r^{-3}\right)f + \Lambda O\left(r^{-1}\right)f  + O\left(r^{-2}\right)\partial_rf + \Lambda O\left(r^{-1}\right)\right)\abs*{\nabla\psi}^2\\
           &  + \left(\left(O(r^{-2}) + a^2O(\Lambda r^{-1})\right)\partial_{rr}f
             +\left( O(r^{-3}) + a^2O(r^{-2}\Lambda)\right)\partial_rf\right.\\
           & \left.+ \left(O(r^{-4}) + O\left(r^{-2}\right)\Lambda
             + a^2\left(O(r^{-3}) + O(r^{-1}\Lambda)\right)\right)f\right)\abs*{\psi}^2.
    \end{align*}
  \item In the specific case where $f=f_p$ and 
    \begin{equation}
      \label{eq:rp:multiplier-choice:specific}
      \rpVF_p \vcentcolon= \rpVF = f_p\left(e_4 + \frac{1}{2}r^{-2}\lambda e_3\right),\qquad
      \rpLagrangianCorr_p \vcentcolon= \frac{r}{\abs*{q}^2}f_p,\qquad
      \rpOneForm_p \vcentcolon= \frac{r}{\abs*{q}^2}\partial_rf_p e_4 + \rpOneForm_{\Lambda},
    \end{equation}
    we have that
    \begin{equation}
      \label{eq:rp:bulk:specific}
      \rpKCurrent
      = \frac{p}{2}r^{p-1}\abs*{\widecheck{\nabla}_4\psi}^2
      + \frac{1}{2}r^{p-1}\left(2 - p\right)\left(
        \abs*{\nabla\psi}^2 + V\abs*{\psi}^2
      \right)
      + \KCurrent{\rpVF_p, \rpLagrangianCorr_p, \rpOneForm_p}_{(\Lambda)}[\psi] + Err,
    \end{equation}
    where 
    \begin{equation*}
      \KCurrent{\rpVF_p, \rpLagrangianCorr_p, \rpOneForm_p}_{(\Lambda)}[\psi]
      > \pth{\Lambda r^2} r^{p-1}\pth{\abs*{\widecheck{\nabla}_4\psi}^2
      + \abs*{r^{-1}\psi}^2},
    \end{equation*}
    and
    \begin{equation*}
      \begin{split}
        Err ={}& \left(O\left(r^{p-4}\right)
                 + O\left(r^{p-2}\right)\Lambda\right)\abs*{\psi}^2
                 + O\left(r^{p-4}\right)\abs*{\nabla_3\psi}^2
                 + \left(O\left(r^{p-2}\right) + O(r^p)\Lambda\right)\abs*{\nabla_4\psi}^2\\
               & + O(r^{p-2} + O(r^p)\Lambda)\abs*{\nabla\psi}^2.
      \end{split}      
    \end{equation*}
  \end{enumerate}
\end{proposition}
\begin{proof}

  We have the computation
  \begin{equation*}
    \begin{split}    
      2\KCurrent{fe_4,0,0}[\psi]
      ={}& \left( -f\partial_r\left(\frac{\Delta}{\abs*{q}^2}\right) + \frac{\Delta}{\abs*{q}^2}\partial_rf \right)\nabla_4\psi\cdot \nabla_4\psi
           + \left(2\frac{fr}{\abs*{q}^2} - \partial_rf\right)\left(\abs*{\NablaAngular\psi}^2
           + V\abs*{\psi}^2
           \right)\notag \\
         & - \left(2\frac{fr}{\abs*{q}^2} V+f\partial_rV\right)\abs*{\psi}^2
           - 2\frac{fr}{\abs*{q}^2}
           \left(
           \HorCovDeriv^{\mu}\psi\cdot  \HorCovDeriv_{\mu}\psi
           + V\abs*{\psi}^2
           \right)
           .
    \end{split}
  \end{equation*}
  Similarly, we can compute that
  \begin{align*}
    2\KCurrent{\frac{1}{2}r^{-2}f\lambda e_3,0,0}[\psi]
    ={}& -\frac{1}{2}e_4\left(r^{-2}\lambda f\right)\nabla_3\psi\cdot \nabla_3\psi
         - \left(r^{-2}\lambda f \partial_r\left(\frac{\Delta}{\abs*{q}^2}\right) - \frac{\Delta}{\abs{q}^2}\partial_r\left(r^{-2}\lambda f\right)\right)\left(\abs*{\NablaAngular \psi}^2 + V\abs*{\psi}^2\right)\notag\\
       & - \frac{\lambda \Delta f}{r\abs*{q}^4}\abs*{\NablaAngular\psi}^2
         - \frac{\lambda \Delta f}{r\abs*{q}^4}\left(\HorCovDeriv^{\mu}\psi\cdot\HorCovDeriv_{\mu}\psi + V\abs*{\psi}^2\right)
         - \frac{1}{2}r^{-2}f\lambda e_3(V)\abs*{\psi}^2.
  \end{align*}
  Since we have that $1\lesssim\lambda\lesssim 1$, we have
  that
  \begin{align*}
    &\KCurrent{\frac{1}{2}r^{-2}f\lambda e_3,0,0}[\psi]\\
    ={}& O\left(r^{-3}f + r^{-2}\partial_rf\right)\abs*{\nabla_3\psi}^2
         + \left(
         O\left(r^{-5}\right)f
         + O\left(r^{-4}\right)\partial_rf
         + \left(O\left(r^{-3}\right)f
         + O\left(r^{-2}\right)\partial_rf\right)\Lambda\right)\abs*{\psi}^2\\
       & +\left(
         O\left(r^{-3}\right)f
         + O\left(r^{-2}\right)\partial_rf
         + O\left(r^{-1}\right)\Lambda f
         + O(1)\Lambda\partial_rf
         \right)\abs*{\nabla\psi}^2
    + O\left(r^{-3}+ r^{-1}\Lambda\right)f\abs*{\nabla_4\psi}^2.
  \end{align*}
  As a result, with $\rpLagrangianCorr$ as chosen in
  \zcref[noname]{eq:rp-Kerr:multiplier-choice:general}, we have that
  \begin{align*}
    2\KCurrent{X,\rpLagrangianCorr,0}[\psi]
    ={}& \left(
         -f\partial_r\left(\frac{\Delta}{\abs*{q}^2}\right)
         + \frac{\Delta}{\abs*{q}^2}\partial_rf \right)\nabla_4\psi\cdot \nabla_4\psi
         + \left(2\frac{fr}{\abs*{q}^2} - \partial_rf\right)\left(\abs*{\NablaAngular\psi}^2
         + V\abs*{\psi}^2
         \right)\\
    & -\ScalarWaveOp[\Metric]\rpLagrangianCorr\abs*{\psi}^2
      - 2\left(\frac{fr}{\abs*{q}^2}V + \frac{1}{2}f\partial_rV\right)\abs*{\psi}^2
      +Err\\
    ={}& \left(
         -f\partial_r\left(\frac{\Delta}{\abs*{q}^2}\right)
         + \frac{\Delta}{\abs*{q}^2}\partial_rf \right)\nabla_4\psi\cdot \nabla_4\psi
         + \left(2\frac{fr}{\abs*{q}^2} - \partial_rf\right)\left(\abs*{\NablaAngular\psi}^2
         + V\abs*{\psi}^2
         \right)\\
    & + \left(
      - \left(2\frac{fr}{\abs*{q}^2}V + f\partial_rV\right)
      - r^{-1}\partial_{rr}f
      + \frac{\Lambda}{3}\left(r\partial_{rr}f + 2\partial_rf - 2fr^{-1}\right)\right)\abs*{\psi}^2
      + Err,
  \end{align*}
  where
  \begin{align*}
    Err ={}& \left(\left(O(r^{-2}) + a^2O(\Lambda r^{-1})\right)\partial_{rr}f
             +\left( O(r^{-3}) + a^2O(r^{-2}\Lambda)\right)\partial_rf\right.\\
           & \left.+ \left(O(r^{-4}) + O\left(r^{-2}\right)\Lambda
             + a^2\left(O(r^{-3}) + O(r^{-1}\Lambda)\right)\right)f\right)\abs*{\psi}^2.
  \end{align*}
  Observe also that
  \begin{align*}
    - \partial_r\left(\frac{\Delta}{\abs*{q^2}}\right)
    &= -\frac{2M}{r^2} + \frac{2r}{3} \Lambda    
    + a^2O(r^{-3}) + a^2O(r^{-1}\Lambda),
    \\\frac{\Delta}{\abs*{q}^2} - 1
    &= - \frac{2M}{r} - \frac{\Lambda}{3}r^2
    + a^2O(r^{-2}) + a^2O(1)\Lambda.
  \end{align*}
  As a result, we have that
  \begin{equation*}
    \begin{split}
      2\KCurrent{X,\rpLagrangianCorr,0}[\psi]
      ={}& \left(\partial_rf + O(r^{-1} + r\Lambda)\partial_rf + O(r^{-2} + \Lambda)f\right) \abs*{\nabla_4\psi}^2\\
         & +\left(\frac{2}{r}f - \partial_rf + O(r^{-2}+\Lambda)f + O(r^{-1}+r\Lambda)\partial_rf\right)\left(\abs*{\nabla\psi}^2 + V\abs*{\psi}^2\right)\\
         &   + \left(      - \left(2\frac{fr}{\abs*{q}^2}V + f\partial_rV\right)
           - r^{-1}\partial_{rr}f\right)\abs*{\psi}^2
           + \frac{\Lambda}{3}\left(2rf-r^2\partial_rf\right)\abs*{\nabla_4\psi}^2 \\
         & + \frac{\Lambda}{3}\left(r\partial_{rr}f + 2\partial_rf - 2fr^{-1}\right)\abs*{\psi}^2
           + Err,
    \end{split}
  \end{equation*}
  which yields \zcref[noname]{eq:rp:bulk:general}. We now prove \zcref[noname]{eq:rp:bulk:specific}. In the case where
  $\rpVF = r^pe_4$, we have that,
  \begin{align*}
    \left(\frac{\Delta}{\abs*{q}^2}-1\right)\partial_rf-f \partial_r\left(\frac{\Delta}{\abs*{q}^2}\right)
    ={}& \left(2-p\right)\Lambda r^{p+1} + O\left(r^{p-2}\right) + O\left(r^p\right)\Lambda .
    \\ \frac{2r}{\abs*{q}^2}V + \partial_rV
    ={}& \frac{4\Lambda}{3r} + O\left(r^{-4}\right) + a^2O\left(r^{-2}\Lambda\right).
  \end{align*}
  where we also used the expression for $V$ from \zcref[noname]{eq:model-problem-gRW}. In the case of $f=f_p=r^p$, we then have that
  \begin{align*}
    \Divergence\left(\abs*{\psi}^2(2-p)\Lambda f_p e_4\right)
    =&{} (2-p)\Lambda\left(
       2r^p\nabla_4\psi\cdot\psi
       + r^{p-1}\abs*{\psi}^2
       + \left(p+1\right)r^{p-1}\abs*{\psi}^2
       \right)\\
    ={}& (2-p)\Lambda r^{p+1}\abs*{\widecheck{\nabla}_4\psi}^2
         - (2-p)\Lambda r^{p+1}\abs*{\nabla_4\psi}^2
         + (2-p)(p+1)\Lambda r^{p+1}\abs*{r^{-1}\psi}^2. 
  \end{align*}
  We can also check $\frac{\Lambda}{3}\left(r\partial_{rr}f + 2\partial_rf - 2fr^{-1}\right)
    = r^{p-1}\left(
      p^2 + p - 2\right)\frac{\Lambda}{3}$.
  Therefore we obtain
  \begin{align*}
    2\KCurrent{X,q,m}[\psi]
    ={}& 2\KCurrent{X,q,0}[\psi]
         + \Divergence\left(\abs*{\psi}^2\rpOneForm\right)\\
    ={}& 2\KCurrent{X,q,0}[\psi]
         + \Divergence\left(\abs*{\psi}^2\left(\frac{r}{\abs*{q}^2}\partial_rf e_4\right)\right)
         + \Divergence\left(\abs*{\psi}^2\left(\frac{r}{\abs*{q}^2}\rpOneForm_{\Lambda}\right)\right)
    \\
    ={}& 2\KCurrent{X,q,0}[\psi]
         + \frac{2r}{\abs*{q}^2}\partial_rf 
         \nabla_4\psi\cdot \psi
         + \left(\partial_r\left(
         \frac{r}{\abs*{q}^2}
         \right)\partial_rf
         + \frac{r}{\abs*{q}^2}\partial_{rr}f
         + \frac{r}{\abs*{q}^2}\partial_rf \Divergence e_4
         \right)\abs*{\psi}^2\\
       &  + \frac{2-p}{3}\Lambda f \nabla_4\psi\cdot\psi
         + \frac{2-p}{3}\Lambda\left(\partial_rf + 2r^{-1}f\right)\abs*{\psi}^2
         + Err
    \\
    ={}& \left(
         -f\partial_r\left(\frac{\Delta}{\abs*{q}^2}\right)
         + \frac{\Delta}{\abs*{q}^2}\partial_rf \right)\nabla_4\psi\cdot \nabla_4\psi
         + \left(2\frac{fr}{\abs*{q}^2} - \partial_rf\right)\left(\abs*{\NablaAngular\psi}^2
         + V\abs*{\psi}^2
         \right)\\
       & + \left(
         -\left(2\frac{fr}{\abs*{q}^2}V\abs*{\psi}^2 + f\partial_rV\right)
         - r^{-1}\partial_{rr}f
         - \frac{\Lambda}{3}\left(r\partial_{rr}f + 2\partial_rf - 2fr^{-1}\right)
         \right)\abs*{\psi}^2\\
       &  + \frac{2-p}{3}\Lambda f \nabla_4\psi\cdot\psi
         + \frac{2-p}{3}\Lambda\left(\partial_rf + 2r^{-1}f\right)\abs*{\psi}^2
         + Err\\
    ={}& pr^{p-1}\left(
         \nabla_4\psi\cdot\nabla_4\psi
         + \frac{2}{r}\nabla_4\psi\cdot\psi
         + \frac{1}{r^2}\abs*{\psi}^2
         \right)
         + \left(
         2-p
         \right)r^{p-1}\left(\abs*{\nabla \psi}^2 + V\abs*{\psi}^2\right)\\
       & + %
         \frac{2-p}{3}\Lambda r^{p+1}\abs*{\nabla_4\psi}^2
         + \left(- \frac{4\Lambda}{3}r^{p-1}
         + \frac{\Lambda}{3}\left(p^2+p-2\right)r^{p-1}\right)\abs*{\psi}^{2}\\
       & + \frac{2-p}{3}\Lambda r^{p+1}\abs*{\widecheck{\nabla}_4\psi}^2
         - \frac{2-p}{3}\Lambda r^{p+1}\abs*{\nabla_4\psi}^2
         +\frac{\left(2-p\right)\left(p+1\right)}{3}\Lambda r^{p-1}\abs*{\psi}^2
         + Err\\
    ={}& r^{p-1}\abs*{\widecheck{\nabla}_4\psi}^2 + (2-p)r^{p-1}\left(\abs*{\nabla\psi}^2 + V\abs*{\psi}^2\right)\\
       &  + \frac{2-p}{3}r^{p+1}\Lambda\abs*{\widecheck{\nabla}_4\psi}^2
         + \frac{2\Lambda}{3}r^{p-1}\left(2-p\right)\abs*{\psi}^2
         + Err
         .
  \end{align*}
  Then, \zcref[noname]{eq:rp:bulk:specific} follows from the fact that we assume that
  $\delta\le p \le 2-\delta$.
\end{proof}

\paragraph{Boundary terms.}

\begin{lemma}
  \label{lemma:rp-boundary:e4-e3-computations}
  Given $\JCurrent{\rpVF, \rpLagrangianCorr, \rpOneForm}[\psi]$ as in
  \zcref[noname]{eq:J-current:def}, with $\left(\rpVF, \rpLagrangianCorr, \rpOneForm\right)$
  as in \zcref[noname]{eq:rp-Kerr:multiplier-choice:general}, we have that
  \begin{align}
    \begin{split}
      \label{eq:rp-boundary:e4-e3-computations:e4}
          \JCurrent{\rpVF, \rpLagrangianCorr, \rpOneForm}[\psi]\cdot e_4
    ={}& f\abs*{\widecheck{\nabla}_4\psi}^2
    - \frac{1}{2}r^{-2}\nabla_4\left(r f \abs*{\psi}^2\right)
    + \frac{1}{2}r^{-2}\lambda \EMTensor(e_4, e_3)\\
    & + \left(O\left(r^{-3}\right)\partial_rf + O\left(r^{-4}\right)f\right)\abs*{\psi}^2.
    \end{split}        
    \\
    \begin{split}
      \label{eq:rp-boundary:e4-e3-computations:e3}
      \JCurrent{\rpVF, \rpLagrangianCorr, \rpOneForm}[\psi]\cdot e_3
      ={}& \frac{1}{2}r^{-2}\nabla_3\left(r f \abs*{\psi}^2\right)
           - \frac{\Lambda}{3}r^2\partial_rf\abs*{\psi}^2
           + f\EMTensor[\psi]_{34}
           + \frac{1}{2}r^{-2}f\lambda \EMTensor[\psi]_{33} \\
         & + O\left(r^{-3}\right)O\left(1 + \Lambda r^2\right)\partial_rf \abs*{\psi}^2
           + O\left(r^{-4}\right)O\left(1 + \Lambda r^2\right)f\abs*{\psi}^2 \\
         & + a^2O\left(r^{-3}\right)f \left(\abs*{\psi}^2 + \abs*{\nabla_3\psi}^2\right)
           .
    \end{split}    
    \\
    \begin{split}
      \label{eq:rp-boundary:e4-e3-computations:Y}
      \JCurrent{\rpVF, \rpLagrangianCorr, \rpOneForm}[\psi]\cdot Y
      ={}& \abs*{f}\abs*{\widecheck{\nabla}_4\psi}\abs*{Y}\abs*{\nabla\psi}
           + \frac{1}{2}r^{-2}\lambda\abs*{f}\abs*{\nabla_3\psi}\abs*{\nabla_Y\psi}
           + f O(R^{-1}r^{-3})\abs*{\nabla\psi}^2\\
         & + f O\left(R^{-1}r^{-3}\right) \abs*{\psi}^2.         
    \end{split}    
  \end{align}
  It will also be convenient to have the alternative expressions
  \begin{align}
    \begin{split}
      \label{eq:rp-boundary:e4-e3-computations:e4:alt}
      \JCurrent{\rpVF, \rpLagrangianCorr, \rpOneForm}[\psi]\cdot e_4
    ={}& f\abs*{\nabla_4\psi}^2
         + f\frac{r}{\abs*{q}^2}\psi\cdot \nabla_4\psi
         - \frac{1}{2}e_4\left(f r^{-1}\right)\abs*{\psi}^2
         + \frac{1}{2}r^{-2}\lambda \EMTensor[\psi]_{34}\\
       & + \left( O(r^{-3})\partial_rf + O(r^{-4})f \right)\abs*{\psi}^2,
    \end{split}
    \\
    \begin{split}
      \label{eq:rp-boundary:e4-e3-computations:e3:alt}
          \JCurrent{\rpVF, \rpLagrangianCorr, \rpOneForm}[\psi]\cdot e_3
    ={}& f\EMTensor[\psi]_{34}
         + \frac{1}{2}fr^{-2}\lambda\EMTensor[\psi]_{33}
         + r^{-1}f\psi\cdot\nabla_3\psi
         - \frac{1}{2}e_3(r^{-1}f)\abs*{\psi}^2
         - \frac{r}{\abs*{q}^2}f'\abs*{\psi}^2 \\
       & + \left(O(r^{-3})f + O(r^{-2})\partial_rf\right)\abs*{\psi}^2         
         + O(r^{-3})f\abs*{\nabla_3\psi}^2. 
    \end{split}
  \end{align}
\end{lemma}

\begin{proof}
  Since we chose $\rpVF=f\left(e_4 + \frac{1}{2}r^{-2}\lambda
    e_3\right)$, we have that
  \begin{align*}
    \JCurrent{\rpVF, \rpLagrangianCorr, \rpOneForm}[\psi]\cdot e_4
    ={}& \left(
         \EMTensor_{\mu\nu}\rpVF^{\nu}
         +\rpLagrangianCorr\psi\cdot\CovariantDeriv_{\mu}\psi
         - \frac{1}{2}\abs*{\psi}^2\partial_{\mu}\psi
         + \frac{1}{2}\abs*{\psi}^2\rpOneForm_{\mu}
         \right)e_4^{\mu}\\
    ={}& \EMTensor(\rpVF, e_4)
         + f\frac{r}{\abs*{q}^2}\psi\cdot \nabla_4\psi
         - \frac{1}{2}e_4\left(f \frac{r}{\abs*{q}^2}\right)\abs*{\psi}^2\\
    ={}& f\EMTensor(e_4, e_4)
         + \frac{1}{2}r^{-2}\lambda \EMTensor(e_4, e_3)
         + f\frac{r}{\abs*{q}^2}\psi\cdot \nabla_4\psi
         - \frac{1}{2}e_4\left(f \frac{r}{\abs*{q}^2}\right)\abs*{\psi}^2\\
    ={}& \underbrace{f\EMTensor(e_4, e_4)         
         + f\frac{r}{\abs*{q}^2}\psi\cdot \nabla_4\psi
         - \frac{1}{2}e_4\left(f r^{-1}\right)\abs*{\psi}^2}_{=\vcentcolon I}
         + \frac{1}{2}r^{-2}\lambda \EMTensor(e_4, e_3)
        + \frac{1}{2}e_4\left(f\frac{a^2\cos^2\theta}{\abs*{q}^2r}\right)\abs*{\psi}^2.
  \end{align*}
  With $I$ defined as above, we get
  \begin{align*}
    I ={}&  f\left(
         \abs*{\nabla_4\psi}^2
         + \frac{1}{r}\psi\cdot\nabla_4\psi
         \right)
         - \frac{1}{2}e_4\left(r^{-1}f\right)\abs*{\psi}^2\\
    ={}& f\abs*{\widecheck{\nabla}_4\psi}^2
         - \frac{1}{2}r^{-2}\nabla_4\left(r f \abs*{\psi}^2\right)
         + r^{-2}\left(e_4(r) - 1\right)f\abs*{\psi}^2\\
    ={}&f\abs*{\widecheck{\nabla}_4\psi}^2
    - \frac{1}{2}r^{-2}\nabla_4\left(r f \abs*{\psi}^2\right),
  \end{align*}
  where we also used the fact that we work on exact \KdS. This gives us the desired expression for $\JCurrent{\rpVF, \rpLagrangianCorr, \rpOneForm}[\psi]\cdot e_4$. We now consider $\JCurrent{\rpVF, \rpLagrangianCorr, \rpOneForm}[\psi]\cdot e_3$. Observe that
  \begin{align*}
    \JCurrent{\rpVF, \rpLagrangianCorr, \rpOneForm}[\psi]\cdot e_3
    ={}& \left(
         \EMTensor_{\mu\nu}\rpVF^{\nu}
         +\rpLagrangianCorr\psi\cdot\CovariantDeriv_{\mu}\psi
         - \frac{1}{2}\abs*{\psi}^2\partial_{\mu}\psi
         + \frac{1}{2}\abs*{\psi}^2\rpOneForm_{\mu}
         \right)e_3^{\mu}\\
    ={}& \EMTensor\left(\rpVF, e_3\right)
         + f\frac{r}{\abs*{q}^2}\psi\cdot\nabla_3\psi
         - \frac{1}{2}e_3\left(f\frac{r}{\abs*{q}^2}\right)\abs*{\psi}^2
         - \frac{r}{\abs*{q}^2}\partial_rf\abs*{\psi}^2\\
    ={}& f\EMTensor[\psi]_{34}
         + \frac{1}{2}r^{-2}f\lambda \EMTensor[\psi]_{33}
         + \frac{1}{2}f\frac{r}{\abs*{q}^2}\nabla_3\left(\abs*{\psi}^2\right)
         - \frac{1}{2}e_3\left(f \frac{r}{\abs*{q}^2}\right)\abs*{\psi}^2
         - \frac{r}{\abs*{q}^2}\partial_rf \abs*{\psi}^2\\
    ={}& \frac{1}{2}r^{-1}f\nabla_3\left(\abs*{\psi}^2\right)
         - \frac{1}{2}e_3\left(r^{-1}f\right)\abs*{\psi}^2
         - \frac{1}{2}\partial_rf\abs*{\psi}^2
         + f\EMTensor[\psi]_{34}
         + \frac{1}{2}r^{-2}f\lambda \EMTensor[\psi]_{33}\\
    &   - \frac{1}{2}f\frac{a^2\cos^2\theta}{\abs*{q}^2r}\nabla_3\left(\abs*{\psi}^2\right)
         + \frac{1}{2}e_3\left(f\frac{a^2\cos^2\theta}{\abs*{q}^2r}\right)\abs*{\psi}^2.
  \end{align*}
  Observe that directly from the second-to-last line of the previous
  equation, we can use the fact that
  $\frac{r}{\abs*{q}^2} - \frac{1}{r} = O(r^{-3})$ to write 
  \begin{align*}
    \JCurrent{\rpVF, \rpLagrangianCorr, \rpOneForm}[\psi]\cdot e_3
    ={}& f\EMTensor[\psi]_{34}
         + \frac{1}{2}r^{-2}f\lambda \EMTensor[\psi]_{33}
         + r^{-1}f\psi\cdot\nabla_3\psi
         - \frac{1}{2}e_3(r^{-1}f)\abs*{\psi}^2
         - \frac{r}{\abs*{q}^2}\partial_rf\abs*{\psi}^2\\
       & + \left(O(r^{-3})f + O(r^{-2})\partial_rf\right)\abs*{\psi}^2         
         + O(r^{-3})f\abs*{\nabla_3\psi}^2.
  \end{align*}
  Denoting
  \begin{equation*}
    J = \frac{1}{2}r^{-1}f\nabla_3\left(\abs*{\psi}^2\right)
    - \frac{1}{2}e_3\left(r^{-1}f\right)\abs*{\psi}^2
    - \frac{1}{2}\partial_rf\abs*{\psi}^2,
  \end{equation*}
  we have that
  \begin{align*}
    2J ={}& r^{-2}\nabla_3\left(rf\abs*{\psi}^2\right)
            - r^{-2}e_3\left(rf \right)\abs*{\psi}^2
            - e_3\left(r^{-1}f\right)\abs*{\psi}^2
            - 2r^{-1}\partial_rf \abs*{\psi}^2\\
    ={}& r^{-2}\nabla_3\left(rf\abs*{\psi}^2\right)
         - 2\partial_rf\left(r^{-1}e_3(r) + r^{-1}\right)\abs*{\psi}^2.
  \end{align*}
  On exact \KdS, we have that $e_3(r) = -\frac{\Delta}{\abs*{q}^2}$,
  so we have in fact that
  \begin{equation*}
    \frac{\Delta}{\abs*{q}^2} - 1
    = -\frac{2M}{r}
    - \frac{\Lambda}{3}r^2
    + aO\left(1\right)\Lambda
    + O\left(r^{-2}\right).
  \end{equation*}
  As a result, we have that
  \begin{equation*}
    J ={} \frac{1}{2}r^{-2}\nabla_3\left(rf\abs*{\psi}^2\right)
    + \left(
      - \frac{\Lambda}{3}r
      + O\left(r^{-1}\right)\Lambda
      + O\left(r^{-2}\right)\right)\partial_rf\abs*{\psi}^2.
  \end{equation*}
  Therefore, we have that
  \begin{equation*}
    \begin{split}
      \JCurrent{\rpVF, \rpLagrangianCorr, \rpOneForm}[\psi]\cdot e_3
      ={}& \frac{1}{2}\left(
           r^{-2}\nabla_3\left(rf\abs*{\psi}^2\right)
           - 2\partial_rf\left(r^{-1}e_3(r) + r^{-1}\right)\abs*{\psi}^2
           \right)
           + f\EMTensor[\psi]_{34}
           + \frac{1}{2}r^{-2}f\lambda \EMTensor[\psi]_{33}\\
         &   - \frac{1}{2}f\frac{a^2\cos^2\theta}{\abs*{q}^2r}\nabla_3\left(\abs*{\psi}^2\right)
           + \frac{1}{2}e_3\left(f\frac{a^2\cos^2\theta}{\abs*{q}^2r}\right)\abs*{\psi}^2\\
      ={}& \frac{1}{2}r^{-2}\nabla_3\left(r f \abs*{\psi}^2\right)
           - \frac{\Lambda}{3}r\partial_rf\abs*{\psi}^2           
           + f\EMTensor[\psi]_{34}
           + \frac{1}{2}r^{-2}f\lambda \EMTensor[\psi]_{33}\\
         & + O\left(r^{-3}\right)e_3(r)\partial_rf \abs*{\psi}^2
           + e_3(r)O\left(r^{-4}\right)f\abs*{\psi}^2
           + a^2O\left(r^{-3}\right)f \left(\abs*{\psi}^2 + \abs*{\nabla_3\psi}^2\right)\\
      ={}& \frac{1}{2}r^{-2}\nabla_3\left(r f \abs*{\psi}^2\right)
           - \frac{\Lambda}{3}r\partial_rf\abs*{\psi}^2 
           + f\EMTensor[\psi]_{34}
           + \frac{1}{2}r^{-2}f\lambda \EMTensor[\psi]_{33}\\
         & + O\left(r^{-3}\right)O\left(1 + \Lambda r^2\right)\partial_rf \abs*{\psi}^2
           + O\left(r^{-4}\right)O\left(1 + \Lambda r^2\right)f\abs*{\psi}^2\\
         & + a^2O\left(r^{-3}\right)f \left(\abs*{\psi}^2 + \abs*{\nabla_3\psi}^2\right)
           .
    \end{split}
  \end{equation*}

  Finally, we also have that
  \begin{align*}
    \JCurrent{\rpVF, \rpLagrangianCorr, \rpOneForm}[\psi]\cdot Y
    ={}& \left(
         \EMTensor_{\mu\nu}\rpVF^{\nu}
         +\rpLagrangianCorr\psi\cdot\CovariantDeriv_{\mu}\psi
         - \frac{1}{2}\abs*{\psi}^2\partial_{\mu}\psi
         + \frac{1}{2}\abs*{\psi}^2\rpOneForm_{\mu}
         \right)Y^{\mu}\\
    ={}& \EMTensor\left(\rpVF, Y\right)
         + \frac{r}{\abs*{q}^2}f\psi\cdot\nabla_Y\psi
         - \abs*{\psi}^2Y\left(\frac{r}{2\abs*{q}^2}f\right)\\
    ={}& f \EMTensor\left(e_4, Y\right)
         + \frac{1}{2}r^{-2}\lambda f \EMTensor(e_3, Y)
         + \frac{r}{\abs*{q}^2}f\psi\cdot\nabla_Y\psi
         - \abs*{\psi}^2Y\left(\frac{r}{2\abs*{q}^2}f\right)\\
    ={}& f\nabla_4\psi\cdot\nabla_Y\psi
         + \frac{1}{2}r^{-2}f\lambda\nabla_3\psi\cdot\nabla_Y\psi
         + \frac{r}{\abs*{q}^2}f\psi\cdot\nabla_Y\psi
         - \abs*{\psi}^2Y\left(\frac{r}{2\abs*{q}^2}f\right)\\
    ={}& f\left(\widecheck{\nabla}_4\psi - r^{-1}\psi\right)\cdot \nabla_Y\psi
         + \frac{1}{2}r^{-2}f\lambda\nabla_3\psi\cdot\nabla_Y\psi
         + \frac{r}{\abs*{q}^2}f\psi\cdot\nabla_Y\psi
         - \abs*{\psi}^2Y\left(\frac{r}{2\abs*{q}^2}f\right). 
  \end{align*}
  As a result,
  \begin{equation*}
    \JCurrent{\rpVF, \rpLagrangianCorr, \rpOneForm}[\psi]\cdot Y
    = f\widecheck{\nabla}_4\psi\cdot\nabla_Y\psi
    + \left(\frac{r}{\abs*{q}^2}- \frac{1}{r}\right)f\psi\cdot\nabla_Y\psi
    + \frac{1}{2}r^{-2}f\lambda\nabla_3\psi\cdot\nabla_Y\psi
    - \abs*{\psi}^2Y\left(\frac{r}{2\abs*{q}^2}f\right).
  \end{equation*}
  Hence we have the estimate
  \begin{align*}
    \abs*{\JCurrent{\rpVF, \rpLagrangianCorr, \rpOneForm}[\psi]\cdot Y}
    ={}&  \abs*{f}\abs*{\widecheck{\nabla}_4\psi}\abs*{Y}\abs*{\nabla\psi}
         + \frac{1}{2}r^{-2}\lambda\abs*{f}\abs*{\nabla_3\psi}\abs*{\nabla_{Y}\psi}
         + O(r^{-4})f\abs*{\psi}\abs*{\nabla\psi}
        + O(R^{-1}r^{-3})f\abs*{\psi}^2\\
    ={}& \abs*{f}\abs*{\widecheck{\nabla}_4\psi}\abs*{Y}\abs*{\nabla\psi}
         + \frac{1}{2}r^{-2}\lambda\abs*{f}\abs*{\nabla_3\psi}\abs*{\nabla_{Y}\psi}
         + O(R^{-1})r^{-3}f\abs*{\nabla\psi}^2
         + O(R^{-1}r^{-3})f\abs*{\psi}^2,
  \end{align*}
  as stated.  
\end{proof}

We are now ready to prove the main inequalities for the relevant
boundary quantities that appear in the $r^p$-estimates.

\begin{proposition}
  \label{prop:rp:main-boundary-estimates}
  The following bounds hold true for $r_0\ge R$ sufficiently large
  \begin{enumerate}
  \item On $\Sigma=\Sigma(\tau)$, for $0\le p \le 2-\delta$,
    \begin{equation}
      \label{eq:rp:main-boundary:estimates:Sigmatau:general}
      \begin{split}
        \JCurrent{\rpVF, \rpLagrangianCorr, \rpOneForm}[\psi]\cdot N_{\Sigma}
        \gtrsim{}& f\abs*{\widecheck{\nabla}_4\psi}^2
                   + r^{-2}f\abs*{\nabla\psi}^2
                   - \frac{1}{2}\Divergence_{\Sigma}\left(
                   r^{-1}f\abs*{\psi}^2\nu_{\Sigma}
                   \right)\\
                 & - O(R^{-1})\left(r^{-3}f + r^{-2}\partial_rf\right)\abs*{\psi}^2
                   .
      \end{split}      
    \end{equation}
  \item On $\SigmaStar$, for $0\le p\le 2-\delta$,
    \begin{equation}
      \label{eq:rp:main-boundary:estimates:SigmaStar:general}
      \begin{split}
        \JCurrent{\rpVF, \rpLagrangianCorr, \rpOneForm}[\psi]\cdot N_{\SigmaStar}
        \gtrsim{}& \frac{1}{2}\left(
         f\left(\abs*{\nabla\psi}^2+\left(\frac{4}{\abs*{q}^2} + \frac{2\Lambda}{3}\right)\abs*{\psi}^2\right)
         - \frac{\Lambda}{3}r^2\partial_rf\abs*{\psi}^2
         + \frac{1}{2}r^{-2}f\lambda\abs*{\nabla_3\psi}^2
         \right)\\
       & + \frac{-\Delta}{2\abs*{q}^2}f\abs*{\widecheck{\nabla}_4\psi}^2
         + \frac{1}{2}r^{-2}\Divergence_{\SigmaStar}\left(rf\abs*{\psi}^2\nu_{\SigmaStar}\right).
      \end{split}      
    \end{equation}

    In particular, with $f = r^p$, we have that
    \begin{equation}
      \label{eq:rp:main-boundary:estimates:SigmaStar:f-specific:full}
      \begin{split}
      \JCurrent{\rpVF, \rpLagrangianCorr, \rpOneForm}[\psi]\cdot N_{\SigmaStar}
      \ge{}& \frac{1}{2}
             \left(
             r^p\abs*{\nabla\psi}^2
             + \frac{4}{\abs*{q}^2}r^p\abs*{\psi}^2
             + \frac{2-p}{3}\Lambda r^p\abs*{\psi}^2
             \right)
             + \frac{-\Delta}{2\abs*{q}^2}r^p\abs*{\widecheck{\nabla}_4\psi}^2\\
           &  + \frac{1}{4}r^{p-2}\lambda\abs*{\nabla_3\psi}^2
             + \frac{1}{2}\Divergence_{\SigmaStar}\left(r^{p-1}\abs*{\psi}^2\nu_{\SigmaStar}\right),
    \end{split}    
    \end{equation}
  \item On $\Sigma=\Sigma(\tau)$, with $f = r^p$, for
    $\delta\le p\le 1-\delta$, we have that
    \begin{equation}
      \label{eq:rp:main-boundary:estimates:Sigmatau:f-specific:1-delta}
      \begin{split}
        \JCurrent{\rpVF, \rpLagrangianCorr, \rpOneForm}[\psi]\cdot N_{\Sigma}
       \ge{}& \frac{(p-1)^2}{8}r^{p-2}\abs*{\psi}^2
              + \frac{1}{2}r^{p-2}\lambda \abs*{\nabla\psi}^2\\
             & - \frac{p}{2}\Divergence_{\Sigma}\left(r^{p-1}\abs*{\psi}^2\nu_{\Sigma}\right)              
              + O(r^{p-3})\left(\abs*{\nabla_3\psi}^2 
               + \abs*{\psi}^2\right).
      \end{split}
    \end{equation}
  \item On $\SigmaStar$ with $f= r^p$ and for $0\le p \le 1-\delta$,
    we have that
    \begin{equation}
      \label{eq:rp:main-boundary:estimates:SigmaStar:f-specific:1-delta}
      \begin{split}
        \JCurrent{\rpVF, \rpLagrangianCorr, \rpOneForm}[\psi]\cdot N_{\SigmaStar}
        \ge{}& \frac{1}{2}
               \left(
               r^p\abs*{\nabla\psi}^2
               + \frac{4}{\abs*{q}^2}r^p\abs*{\psi}^2
               + \frac{2-p}{3}\Lambda r^p\abs*{\psi}^2
               \right)
               + \frac{-\Delta}{2\abs*{q}^2}r^p\abs*{\widecheck{\nabla}_4\psi}^2\\
             &  + \frac{1}{4}r^{p-2}\lambda\abs*{\nabla_3\psi}^2
               + \frac{1}{2}\Divergence_{\SigmaStar}\left(r^{p-1}\abs*{\psi}^2\nu_{\SigmaStar}\right)\\
             & + O\left(r^{p-3}\abs*{\psi}^2 + r^{p-3}\abs*{\nabla_3\psi}^2\right).
      \end{split}    
    \end{equation}

  \end{enumerate}
\end{proposition}
\begin{proof}
  Recall that $N_{\Sigma} =  e_4 + \frac{1}{2}r^{-2}\lambda e_3 + Y^be_b$ and $\nu_{\Sigma} = e_4 - \frac{1}{2}r^{-2}\lambda e_3$. We have that
  \begin{align*}
    &\JCurrent{\rpVF, \rpLagrangianCorr, \rpOneForm}[\psi]\cdot \left(e_4 + \frac{1}{2}r^{-2}\lambda e_3 \right)\\
    ={}& f\abs*{\widecheck{\nabla}_4\psi}^2
         -\frac{1}{2}r^{-2}\nabla_4\left(rf\abs*{\psi}^2\right)
         + \frac{1}{2}r^{-2}\lambda\EMTensor\left(e_3,e_4\right)
         + \frac{1}{2}e_4\left(
         f\frac{a^2\cos^2\theta}{\abs*{q}^2r}
         \right)\abs*{\psi}^2\\
    & + \frac{1}{2}r^{-2}\lambda\left(\frac{1}{2}\left(
      r^{-2}\nabla_3\left(rf\abs*{\psi}^2\right)
      - 2 \partial_rf\left(r^{-1}e_3(r)+r^{-1}\right)\abs*{\psi}^2
      \right)
      + f\EMTensor(e_3,e_4)
      + \frac{1}{2}r^{-2}f\lambda\EMTensor(e_3,e_3)\right)\\
    & + \frac{1}{2}r^{-2}\lambda\left(
      - \frac{1}{2}f\frac{a^2\cos^2\theta}{\abs*{q}^2r}\nabla_3\left(\abs*{\psi}^2\right)
      + \frac{1}{2}e_3\left(f\frac{a^2\cos^2\theta}{\abs*{q}^2r}\right)\abs*{\psi}^2
      \right)\\
    ={}& f\abs*{\widecheck{\nabla}_4\psi}^2
         + r^{-2}\lambda f \EMTensor_{34}
         + \frac{1}{4}r^{-4}f\lambda^2\EMTensor_{33}
         - \frac{1}{2}r^{-2}\left(\nabla_4 - \frac{1}{2}r^{-2}\lambda\nabla_3\right)\left(rf\abs*{\psi}^2\right)\\
    &- \frac{1}{2}r^{-2}\lambda\partial_rf\left(r^{-1}e_3(r) + r^{-1}\right)\abs*{\psi}^2
      + \frac{1}{2}\left(
      e_4\left(
      f\frac{a^2\cos^2\theta}{\abs*{q}^2r}
      \right)\abs*{\psi}^2\right.\\
    &\left. + r^{-2}\lambda\left(
      - \frac{1}{2}f\frac{a^2\cos^2\theta}{\abs*{q}^2r}\nabla_3\left(\abs*{\psi}^2\right)
      + \frac{1}{2}e_3\left(f\frac{a^2\cos^2\theta}{\abs*{q}^2r}\right)\abs*{\psi}^2
      \right)
      \right)
  \end{align*}
  Denoting
  \begin{equation*}    
    Err^{\bowtie} \vcentcolon=
    \frac{1}{2}\left(
      e_4\left(
        f\frac{a^2\cos^2\theta}{\abs*{q}^2r}
      \right)\abs*{\psi}^2 + r^{-2}\lambda\left(
        - \frac{1}{2}f\frac{a^2\cos^2\theta}{\abs*{q}^2r}\nabla_3\left(\abs*{\psi}^2\right)
        + \frac{1}{2}e_3\left(f\frac{a^2\cos^2\theta}{\abs*{q}^2r}\right)\abs*{\psi}^2
      \right)
    \right),
  \end{equation*}
  and observing $Err^{\bowtie}
    = \left(O\left(r^{-4}\right)f + O\left(r^{-3}\right)\partial_rf\right)\abs*{\psi}^2$ we have that
  \begin{align*}
    \JCurrent{\rpVF, \rpLagrangianCorr, \rpOneForm}[\psi]\cdot \left(e_4 + \frac{1}{2}r^{-2}\lambda e_3 \right)
    ={}& f\abs*{\widecheck{\nabla}_4\psi}^2
         + r^{-2}\lambda f \left(\abs*{\nabla\psi}^2 + V\abs*{\psi}^2\right)
         + \frac{1}{4}r^{-4}f\lambda^2\abs*{\nabla_3\psi}^2\\
       &- \frac{1}{2}r^{-2}\nu_{\Sigma}\left(rf\abs*{\psi}^2\right)
         - \frac{1}{2}r^{-2}\lambda\partial_rf\left(r^{-1}e_3(r) + r^{-1}\right)\abs*{\psi}^2
         + Err^{\bowtie}.
  \end{align*}
  Then observe that in \KdS, we have $e_3(r)+1=\frac{2M}{r} + O(r^{-2}) + \frac{r^2}{3}\Lambda + O(1)\Lambda$ so that
  \begin{equation}
    \label{eq:rp:boundary:Sigma:general:L2-estimate}
    -r^{-2}\lambda\partial_rf\left(r^{-1}e_3(r) + r^{-1}\right)\abs*{\psi}^2
    = - \frac{\Lambda}{3}\lambda r^{-1} \partial_rf \abs*{\psi}^2
    + \left(O(r^{-4}) + O(r^{-3})\Lambda\right)\partial_rf\abs*{\psi}^2.
  \end{equation}  
  Now observe that
  \begin{align*}
    \nu_{\Sigma}(r)
    ={}& e_4(r) - \frac{1}{2}r^{-2}\lambda e_3(r)
    ={} 1 + \frac{1}{2}r^{-2}\lambda \frac{\Delta}{\abs*{q}^2}
    ={} 1 +  O(r^{-2})\lambda + O(\Lambda)\lambda.
  \end{align*}
  We can now write that
  \begin{align*}
    \Divergence_{\Sigma}\left(r^{-1}f\abs*{\psi}^2\nu_{\Sigma}\right)
    ={}& \nu_{\Sigma}\left(r^{-1}f\abs*{\psi}^2\right)
         + r^{-1}f\abs*{\psi}^2\Divergence_{\Sigma}\left(\nu_{\Sigma}\right)\\
    ={}& r^{-2}\nu_{\Sigma}\left(rf\abs*{\psi}^2\right)
         + \nu_{\Sigma}(r^{-2})rf\abs*{\psi}^2
         + r^{-1}f \abs*{\psi}^2\Divergence_{\Sigma}\left(\nu_{\Sigma}\right)\\
    ={}& r^{-2}\nu_{\Sigma}\left(rf\abs*{\psi}^2\right)
         - 2r^{-2}\nu_{\Sigma}\left(rf\abs*{\psi}^2\right)
         + r^{-1}f\abs*{\psi}^2\Divergence_{\Sigma}(\nu_{\Sigma}).
  \end{align*}
  As a result, we have that
  \begin{align*}
    r^{-2}\nu_{\Sigma}\left(rf\abs*{\psi}^2\right)
    ={}& \Divergence_{\Sigma}\left(r^{-1}f\abs*{\psi}^2\nu_{\Sigma}\right)
         - r^{-1}f\abs*{\psi}^2\Divergence_{\Sigma}(\nu_{\Sigma})
         + 2fr^{-2}\abs*{\psi}^2\nu_{\Sigma}(r)\\
    ={}& \Divergence_{\Sigma}\left(r^{-1}f\abs*{\psi}^2\nu_{\Sigma}\right)
         + fr^{-1}\abs*{\psi}^2\left(
         2r^{-1}\nu_{\Sigma}(r)
         - \Divergence_{\Sigma}(\nu_{\Sigma})
         \right).
  \end{align*}
  Then, recalling \zcref{eq:div-Sigma} and the exact values for
  $\Trace\chi$ and $\xi$ from \zcref{lemma:Kerr:outgoing-PG:Ric-and-curvature},
  we have that
  \begin{align*}
    2r^{-1}\nu_{\Sigma}(r)-\Divergence_{\Sigma}(\nu_{\Sigma})
    ={}& \frac{2}{r} - \frac{\Delta}{r^3\abs*{q}^2}
         - \frac{2r}{\abs*{q}^2} + \frac{\Delta }{r \abs*{q}^4}\lambda
         + O(M^2 r^{-3})
    \\
    ={}& O(r^{-3}) + O(r^{-1})\Lambda + O(M^2 r^{-3}).
  \end{align*}
  As a result, we have that
  \begin{equation}
    \label{eq:rp:boundary:Sigma:general:div:aux-nu-Sigma-and-div}
    r^{-2}\nu_{\Sigma}(rf\abs*{\psi}^2)
    = \Divergence_{\Sigma}\left(r^{-1}f\abs*{\psi}^2\nu_{\Sigma}\right)
    + \left(O(r^{-4})+ O(r^{-2})\Lambda\right)f\abs*{\psi}^2.
  \end{equation}
  As a result, we have that
  \begin{align*}
    \JCurrent{\rpVF, \rpLagrangianCorr, \rpOneForm}[\psi]\cdot \left(e_4 + \frac{1}{2}r^{-2}\lambda e_3\right)
    ={}& f\abs*{\widecheck{\nabla}_4\psi}^2
         + \frac{1}{4}r^{-4}f\lambda^2\abs*{\nabla_3\psi}^2
         + r^{-2}\lambda f\left(\abs*{\nabla\psi}^2 + V\abs*{\psi}^2\right)\\
       & - \frac{1}{2}\Divergence_{\Sigma}\left(r^{-1}f\abs*{\psi}^2\nu_{\Sigma}\right)
         - \frac{\Lambda}{6}\lambda r^{-1}\partial_rf\abs*{\psi}^2 + Err^{\bowtie}_1,
  \end{align*}
  where
  \begin{align*}
    Err^{\bowtie}_1 & = Err^{\bowtie}
    + \left(O(r^{-4}) + \Lambda O(r^{-3})\right)\partial_rf\abs*{\psi}^2
    + \left(O(r^{-4})+ O(r^{-2})\Lambda\right)f\abs*{\psi}^2
    \\& = \left(O\left(r^{-4}\right)f + O\left(r^{-3}\right)\partial_rf\right)\abs*{\psi}^2
    + f O(r^{-5})\abs*{\nabla_3\psi}^2.
  \end{align*}
  Recalling from \zcref[noname]{eq:model-problem-gRW} that
  $V= \frac{4}{\abs*{q}^2}\frac{\Delta}{r^2+a^2} + 2\Lambda$, we have
  that
  \begin{align*}
    \JCurrent{\rpVF, \rpLagrangianCorr, \rpOneForm}[\psi]\cdot \left(e_4 + \frac{1}{2}r^{-2}\lambda e_3\right)
    ={}& f\abs*{\widecheck{\nabla}_4\psi}^2
         + \frac{1}{4}r^{-4}f\lambda^2\abs*{\nabla_3\psi}^2
         + r^{-2}\lambda f\left(\abs*{\nabla\psi}^2 + \left(\frac{4}{r^2} + \frac{2\Lambda}{3}\right)\abs*{\psi}^2\right)\\
       & - \frac{1}{2}\Divergence_{\Sigma}\left(r^{-1}f\abs*{\psi}^2\nu_{\Sigma}\right)
         - \frac{\Lambda}{6}\lambda r^{-1}\partial_rf\abs*{\psi}^2
         + Err^{\bowtie}_1.
  \end{align*}  
  Using \zcref[noname]{eq:rp-boundary:e4-e3-computations:Y} and that
  $\delta\le p\le 2-\delta$, we have that
  \begin{align*}
    & \JCurrent{\rpVF, \rpLagrangianCorr, \rpOneForm}[\psi]\cdot
      \left(e_4 + \frac{1}{2}r^{-2}\lambda e_3 + Y^be_b\right)\\
    \ge{}& f\abs*{\widecheck{\nabla}_4\psi}^2
           + \frac{1}{4}r^{-4}f\lambda^2\abs*{\nabla_3\psi}^2
           + r^{-2}\lambda f\left(\abs*{\nabla\psi}^2 + \frac{4}{r^2}\abs*{\psi}^2\right)\\
    & - \frac{1}{2}\Divergence_{\Sigma}\left(r^{-1}f\abs*{\psi}^2\nu_{\Sigma}\right)
      - f\abs*{\widecheck{\nabla}_4\psi}\abs*{Y}\abs*{\nabla\psi}
      + \frac{1}{2}r^{-2}f\lambda\nabla_3\psi\cdot\nabla_Y\psi
      + Err^{\bowtie}_\Sigma,
  \end{align*}
  where $Err^{\bowtie}_\Sigma
    = Err^{\bowtie}_1
    + O(R^{-1}r^{-3})\abs*{\nabla\psi}^2$.
  Observe that $\abs*{Y} = O(r^{-1})$ and moreover using
  \zcref[noname]{eq:n-Sigma-nSigma-metric-value},
  \begin{align*}
    &\abs*{\widecheck{\nabla}_4\psi}^2
      - \abs*{\widecheck{\nabla}_4\psi}\abs*{Y}\abs*{\nabla\psi}
      + \frac{1}{4}r^{-4}\lambda^2\abs*{\nabla_3\psi}^2
      + \frac{1}{2}r^{-2}\lambda\nabla_3\psi\cdot\nabla_Y\psi
      + r^{-2}\lambda \abs*{\nabla\psi}^2\\
    \ge{}& \frac{\delta}{4}r^{-4}\abs*{\lambda\nabla_3\psi}^2
           + \frac{1}{4}\left(\sqrt{1-\delta}r^{-2}\lambda\nabla_3\psi
           + \frac{1}{\sqrt{1-\delta}}\nabla_Y\psi\right)^2           
           - \frac{1}{2(1-\delta)}\abs*{\nabla_Y\psi}^2\\
    & - \frac{4(1-\delta)}{4}\abs*{\widecheck{\nabla}_4\psi}^2
      + \abs*{\widecheck{\nabla}_4\psi}^2
      + r^{-2}\lambda \abs*{\nabla\psi}^2
    \\
    \ge{}& \frac{\delta}{4}r^{-4}\abs*{\lambda\nabla_3\psi}^2
           + \frac{\delta}{4}\abs*{\widecheck{\nabla}_4\psi}^2
           + \frac{1}{2}\left(2r^{-2}\lambda - \abs*{Y}^2\right)\abs*{\nabla\psi}^2
           - \frac{\delta}{2(1-\delta)}\abs*{\nabla_Y\psi}^2\\
    \ge{}& \frac{\delta}{4}r^{-4}\abs*{\lambda\nabla_3\psi}^2
           +\frac{\delta}{4}\abs*{\widecheck{\nabla}_4\psi}^2
           + \frac{c_0M^2}{2r^2}\abs*{\nabla\psi}^2
           - \frac{\delta}{2(1-\delta)}\abs*{\nabla_Y\psi}^2
           ,
  \end{align*}
  where $c_1$ is the implicit constant in
  \zcref[noname]{eq:n-Sigma-nSigma-metric-value}. Since $Y = O(r^{-1})$, we
  can choose $\delta$ sufficiently small so that there exists some
  $\delta_1$ such that $\frac{c_1M^2}{2r^2} - \frac{\delta}{2(1-\delta)}\abs*{Y}^2 > \frac{\delta_1}{r^2}$. This then allows us to conclude that
  \begin{equation*}
    \begin{split}
      \JCurrent{\rpVF, \rpLagrangianCorr, \rpOneForm}[\psi]\cdot \left(e_4 + \frac{1}{2}r^{-2}\lambda e_3 + Y^be_b\right)
      \gtrsim{}& r^{-4}f\abs*{\nabla_3\psi}^2
                 +f\abs*{\widecheck{\nabla}_4\psi}^2
                 + r^{-2}f\abs*{\nabla\psi}^2
                 + r^{-4}\lambda f \abs*{\psi}^2
      \\
               & - \frac{1}{2}\Divergence_{\Sigma}\left(r^{-1}f\abs*{\psi}^2\nu_{\Sigma}\right)
                 + Err^{\bowtie}_\Sigma.
    \end{split}
  \end{equation*}
  Taking $R$ sufficiently large allows us to absorb the
  first-derivative error terms in $Err_{\Sigma}^{\bowtie}$, and we are
  left with \zcref[noname]{eq:rp:main-boundary:estimates:Sigmatau:general}. We now move onto proving
  \zcref[noname]{eq:rp:main-boundary:estimates:SigmaStar:general}.  Recall that since $\SigmaStar$ is an $r$-constant
  hypersurface, we have $N_{\SigmaStar} = \frac{1}{2}e_4(r)e_3
    + \frac{1}{2}e_3(r)e_4$. Then,
  \begin{align}
    &\JCurrent{\rpVF, \rpLagrangianCorr, \rpOneForm}[\psi]\cdot N_{\SigmaStar}\notag\\
    ={}& \frac{1}{2}\left(
         \frac{1}{2}r^{-2}\nabla_3\left(r f \abs*{\psi}^2\right)
         - \frac{\Lambda}{3}r\partial_rf\abs*{\psi}^2
         + f\EMTensor[\psi]_{34}
         + \frac{1}{2}r^{-2}f\lambda \EMTensor[\psi]_{33}\right)\notag\\
    & + \frac{1}{2}\left(O\left(r^{-3}\right)O\left(1 + \Lambda r^2\right)\partial_rf \abs*{\psi}^2
      + O\left(r^{-2}\right)\Lambda f\abs*{\psi}^2
      + a^2O\left(r^{-3}\right)f \left(\abs*{\psi}^2 + \abs*{\nabla_3\psi}^2\right)
      \right)\notag\\
    & + \frac{-\Delta}{2\abs*{q}^2}\left(
      f\abs*{\widecheck{\nabla}_4\psi}^2
      - \frac{1}{2}r^{-2}\nabla_4\left(r f \abs*{\psi}^2\right)
      + \frac{1}{2}r^{-2}\lambda \EMTensor[\psi]_{34}
      + \left(O\left(r^{-3}\right)\partial_rf + O\left(r^{-4}\right)f\right)\abs*{\psi}^2
      \right)\notag\\
    ={}& \frac{1}{2}\left(
         - \frac{\Lambda}{3}r\partial_rf\abs*{\psi}^2
         + f\EMTensor[\psi]_{34}
         + \frac{1}{2}r^{-2}f\lambda \EMTensor[\psi]_{33}\right)\notag\\
    & + \frac{1}{2}\left(O\left(r^{-3}\right)\partial_rf \abs*{\psi}^2
      + O\left(r^{-4}\right)\Lambda f\abs*{\psi}^2
      + a^2O\left(r^{-3}\right)f \left(\abs*{\psi}^2 + \abs*{\nabla_3\psi}^2\right)
      \right)\notag\\
    & + \frac{-\Delta}{2\abs*{q}^2}\left(
      f\abs*{\widecheck{\nabla}_4\psi}^2
      + \frac{1}{2}r^{-2}\lambda \EMTensor[\psi]_{34}
      + \left(O\left(r^{-3}\right)\partial_rf + O\left(r^{-4}\right)f\right)\abs*{\psi}^2
      \right)\notag\\
    & + \frac{1}{2}r^{-2}\nu_{\SigmaStar}\left(rf\abs*{\psi}^2\right)
      . \label{eq:rp:Kerr:boundary:Sigma-star:1}
  \end{align}
  Recalling \zcref{lemma:rp:div-SigmaStar}
  we can compute that
  \begin{align*}
    &\Divergence_{\SigmaStar}\left(r^{-1}f\abs*{\psi}^2\nu_{\SigmaStar}\right)\\
    ={}& r^{-2}\nu_{\SigmaStar}\left(rf \abs*{\psi}^2\right)
         + \nu_{\SigmaStar}(r^{-2})rf\abs*{\psi}^2
         + \nabla_{\SigmaStar}\cdot\nu_{\SigmaStar} r^{-1}f\abs*{\psi}^2\\
    ={}& r^{-2}\nu_{\SigmaStar}\left(rf \abs*{\psi}^2\right)
         + \left(\nabla_{\SigmaStar}\cdot\nu_{\SigmaStar}
         - \frac{2\nu_{\SigmaStar}(r)}{r}\right)r^{-1}f\abs*{\psi}^2\\
    ={}& r^{-2}\nu_{\SigmaStar}\left(rf \abs*{\psi}^2\right)
         + \left(
         \frac{1}{2}e_4(r)\left(
         \Trace\chiBar - \frac{2e_3(r)}{r}\right)
         - \frac{1}{2}e_3(r)\left(
         \Trace\chi - \frac{2e_4(r)}{r}
         \right)
         \right)r^{-1}f\abs*{\psi}^2.
  \end{align*}
  Observe that we have
  \begin{align*}
    \Trace\chiBar - \frac{2e_3(r)}{r}
    ={}& -\frac{\Delta}{\abs*{q}^2}a^2 O\left(r^{-3}\right) ,& \Trace\chi - \frac{2e_4(r)}{r}
    ={}& a^2 O\left(r^{-3}\right).
  \end{align*}
Therefore we find
  \begin{equation}
    \label{eq:rp:kerr:boundary:Sigma-star:div-Sigma-star-aux-eqn}
    \Divergence_{\SigmaStar}\left(r^{-1}f\abs*{\psi}^2\nu_{\SigmaStar}\right)
    = r^{-2}\nu_{\SigmaStar}\left(rf\abs*{\psi}^2\right)
    + a^2\frac{\Delta}{2\abs*{q}^2}O\left(r^{-4}\right)f\abs*{\psi}^2.
  \end{equation}
  Plugging in
  \zcref[noname]{eq:rp:kerr:boundary:Sigma-star:div-Sigma-star-aux-eqn} to
  \zcref[noname]{eq:rp:Kerr:boundary:Sigma-star:1}, we have that
  \begin{align*}
    \JCurrent{\rpVF, \rpLagrangianCorr, \rpOneForm}[\psi]\cdot N_{\SigmaStar}
    ={}& \frac{1}{2}\left(
         f\EMTensor[\psi]_{43}
         - \frac{\Lambda}{3}r\partial_rf\abs*{\psi}^2
         + \frac{1}{2}r^{-2}f\lambda\EMTensor[\psi]_{33}
         \right)\\
       & + \frac{-\Delta}{2\abs*{q}^2}f\abs*{\widecheck{\nabla}_4\psi}^2
         + \frac{1}{2}\Divergence_{\SigmaStar}\left(r^{-1}f\abs*{\psi}^2\nu_{\SigmaStar}\right)
         + Err_{\SigmaStar}^{\bowtie}\\
    ={}& \frac{1}{2}\left(
         f\left(\abs*{\nabla\psi}^2+V\abs*{\psi}^2\right)
         - \frac{\Lambda}{3}r^2\partial_rf\abs*{\psi}^2
         + \frac{1}{2}r^{-2}f\lambda\abs*{\nabla_3\psi}^2
         \right)\\
       & + \frac{-\Delta}{2\abs*{q}^2}f\abs*{\widecheck{\nabla}_4\psi}^2
         + \frac{1}{2}\Divergence_{\SigmaStar}\left(r^{-1}f\abs*{\psi}^2\nu_{\SigmaStar}\right)
         + Err_{\SigmaStar}^{\bowtie},
  \end{align*}
  where
  \begin{equation}
    \label{eq:rp:boundary:Err-SigmaStar:def}
    \begin{split}
      Err_{\SigmaStar}^{\bowtie}
    ={}& \frac{1}{2}\left(O\left(r^{-3}\right)O\left(1 + \Lambda r^2\right)\partial_rf \abs*{\psi}^2
      + O\left(r^{-2}\right)\Lambda f\abs*{\psi}^2
       + a^2O\left(r^{-3}\right)f \left(\abs*{\psi}^2 + \abs*{\nabla_3\psi}^2\right)
    \right)\\
    & + \frac{-\Delta}{2\abs*{q}^2}\left(O\left(r^{-3}\right)\partial_rf + O\left(r^{-4}\right)f\right)\abs*{\psi}^2
    \end{split}    
  \end{equation}
  Recalling that
  $V = \frac{4}{\abs*{q}^2}\frac{\Delta}{r^2+a^2} + 2\Lambda + O(ar^{-4})$, we
  have that
  \begin{equation*}
    \begin{split}
      \JCurrent{\rpVF, \rpLagrangianCorr, \rpOneForm}[\psi]\cdot N_{\SigmaStar}
    ={}& \frac{1}{2}\left(
         f\left(\abs*{\nabla\psi}^2+\left(\frac{4}{\abs*{q}^2} + \frac{2\Lambda}{3}\right)\abs*{\psi}^2\right)
         - \frac{\Lambda}{3}r\partial_rf\abs*{\psi}^2
         + \frac{1}{2}r^{-2}f\lambda\abs*{\nabla_3\psi}^2
         \right)\\
       & + \frac{-\Delta}{2\abs*{q}^2}f\abs*{\widecheck{\nabla}_4\psi}^2
         + \frac{1}{2}\Divergence_{\SigmaStar}\left(r^{-1}f\abs*{\psi}^2\nu_{\SigmaStar}\right)
         + Err_{\SigmaStar}^{\bowtie},
    \end{split}    
  \end{equation*}
  where we observe that for $f=f_p$, we have that $Err_{\SigmaStar}^{\Aux}
    = O\left(r^{p-3}\abs*{\psi}^2 + r^{p-3}\abs*{\nabla_3\psi}^2\right)$. We now move on to showing
  \zcref[noname]{eq:rp:main-boundary:estimates:Sigmatau:f-specific:1-delta}
  and
  \zcref[noname]{eq:rp:main-boundary:estimates:SigmaStar:f-specific:1-delta}.
  Combining the identities in
  \zcref[noname]{eq:rp-boundary:e4-e3-computations:e4:alt} in
  \zcref[noname]{eq:rp-boundary:e4-e3-computations:e3:alt}, we have 
  \begin{align*}
    \JCurrent{\rpVF, \rpLagrangianCorr, \rpOneForm}[\psi]\cdot e_4
    ={}& f\abs*{\nabla_4\psi}^2
         + fr^{-1}\psi\nabla_4\psi
         - \frac{1}{2}e_4(fr^{-1})\abs*{\psi}^2
         + \frac{1}{2}r^{-2}\lambda f \EMTensor[\psi]_{34} + O(R^{-1})r^{-3}f\abs*{\psi}^2,\\
    \JCurrent{\rpVF, \rpLagrangianCorr, \rpOneForm}[\psi]\cdot e_3
    ={}&f\EMTensor[\psi]_{34}
         + \frac{1}{2}fr^{-2}\lambda\EMTensor[\psi]_{33}
         + r^{-1}f\psi\cdot\nabla_3\psi
         - \frac{1}{2}e_3(r^{-1}f)\abs*{\psi}^2\\
       & - \frac{r}{\abs*{q}^2}f'\abs*{\psi}^2
         + O(R^{-1})r^{-3}f\abs*{\psi}^2. 
  \end{align*}
  Then, recalling that
  $N_{\Sigma} = e_4 + \frac{1}{2}r^{-2}\lambda e_3 + Y^be_b$, we have
  that
  \begin{equation*}
    \begin{split}
      \JCurrent{\rpVF, \rpLagrangianCorr, \rpOneForm}[\psi]\cdot N_{\Sigma}
    ={}& f\abs*{\nabla_4\psi}^2
    + fr^{-1}\psi\cdot\nu_{\Sigma}\psi
    - \frac{1}{2}\nu_{\Sigma}(fr^{-1})\abs*{\psi}^2
    + r^{-2}\lambda\EMTensor[\psi]_{34}
    + \frac{1}{4}fr^{-4}\lambda^2\abs*{\nabla_3\psi}^2\\
    & + O(r^{-3})f\left( \abs*{\psi}^2 + \abs*{\nabla_3\psi}^2\right)
    + \JCurrent{\rpVF, \rpLagrangianCorr, \rpOneForm}[\psi]\cdot N_{\Sigma}\cdot Y,
    \end{split}    
  \end{equation*}
  where we recall that
  $\nu_{\Sigma} = e_4- \frac{1}{2}r^{-2}\lambda e_3$. Then observe that
  \begin{equation*}
    \abs*{\nabla_4\psi}^2
    + \frac{1}{4}r^{-4}\lambda^2\abs*{\nabla_3\psi}^2
    = \abs*{\nabla_{\nu_\Sigma}\psi}^2
    + r^{-2}\lambda\nabla_4\psi\cdot\nabla_3\psi.
  \end{equation*}
  Thus, we have that
  \begin{equation*}
    \begin{split}
      \JCurrent{\rpVF, \rpLagrangianCorr, \rpOneForm}[\psi]\cdot N_{\Sigma}
    ={}& f\abs*{\nabla_{\nu_{\Sigma}}\psi}^2
    + fr^{-1}\psi\cdot\nu_{\Sigma}\psi
    - \frac{1}{2}\nu_{\Sigma}(fr^{-1})\abs*{\psi}^2
    + r^{-2}\lambda\EMTensor[\psi]_{34}
    + fr^{-2}\lambda\nabla_4\psi\cdot\nabla_3\psi\\
    & + O(R^{-1})r^{-3}f\abs*{\psi}^2
    + \JCurrent{\rpVF, \rpLagrangianCorr, \rpOneForm}[\psi]\cdot Y.
    \end{split}    
  \end{equation*}
  Next observing that
  \begin{equation*}
    r^{-2}\lambda\nabla_4\psi\cdot\nabla_3\psi
    \lesssim r^{-1}\abs*{\nabla_{\nu_{\Sigma}}\psi}^2
    + r^{-3}\abs*{\nabla_3\psi}^2,
  \end{equation*}
  and using \zcref[noname]{eq:rp-boundary:e4-e3-computations:Y}, we obtain
  that for $0<\delta_1<1$ chosen sufficiently small, 
  \begin{equation*}
    \begin{split}
      \JCurrent{\rpVF, \rpLagrangianCorr, \rpOneForm}[\psi]\cdot N_{\Sigma}
      \ge{}& f(1-O(r^{-1})\abs*{\nabla_{\nu_{\Sigma}}\psi}^2
      + fr^{-1}\psi\cdot\nabla_{\nu_{\Sigma}}\psi
      - \frac{1}{2}\nu_{\Sigma}(fr^{-1})\abs*{\psi}^2\\
      &+ r^{-2}\lambda f\EMTensor[\psi]_{34}
      - \frac{1-\delta_1}{2}f\abs*{\widecheck{\nabla}_4\psi}^2
      - \frac{1}{2(1-\delta_1)}f\abs*{Y}^2\abs*{\nabla\psi}^2\\
           &+ O(R^{-\delta})r^{-2}f\abs*{\nabla\psi}^2
             + O(r^{-3})\left(\abs*{\nabla_3\psi}^2 + \abs*{\psi}^2\right).
    \end{split}
  \end{equation*}
  Then, denote by $c_1>0$ the constant implicit in
  \zcref[noname]{eq:n-Sigma-nSigma-metric-value}, so that $\abs*{Y}^2\le 2r^{-2}\lambda - c_1\frac{M^2}{r^2}$. We have that
  \begin{equation*}
    r^{-2}\lambda - \frac{1}{2(1-\delta_1)}\abs*{Y}^2
    \ge -\frac{\delta_1}{(1-\delta_1)r^2}\lambda
    + \frac{c_1M^2}{2(1-\delta_1)r^2}
    \ge \frac{c_1M^2}{4r^2},
  \end{equation*}
  provided that $0<\delta_1\ll c_1$ is chosen sufficiently small.
  This implies that there exists some $\mathcal{B} =  O(r^{-1})$
  \begin{equation}
    \label{eq:rp:main-boundary:estimates:Sigmatau:f-specific:1-delta:aux1}
    \begin{split}
      \JCurrent{\rpVF, \rpLagrangianCorr, \rpOneForm}[\psi]\cdot N_{\Sigma}
      \ge{}& f\left( 1-\mathcal{B} \right)\abs*{\nabla_{\nu_{\Sigma}}\psi}^2
             + fr^{-1}\psi\cdot \nu_{\Sigma}\psi
             - \frac{1}{2}\nu_{\Sigma}(fr^{-1})\abs*{\psi}^2
             - \frac{1-\delta_1}{2}f\abs*{\nabla_{\nu_{\Sigma}}\psi + \frac{1}{r}\psi}^2\\
           &  + \frac{C_1M^2}{4r^2}f\EMTensor[\psi]_{34}
             + O(R^{-\delta})r^{-2}f\abs*{\nabla\psi}^2
             + O(r^{-3})f(\abs*{\nabla_3\psi}^2 + \abs*{\psi}^2).
    \end{split}
  \end{equation}
  We now consider
  \begin{equation*}
    J = f(1-\mathcal{B})\abs*{\nabla_{\nu_{\Sigma}}\psi}^2
    + f r^{-1}\psi\cdot\nu_{\Sigma}\psi
    - \frac{1}{2}\nu_{\Sigma}(fr^{-1})\abs*{\psi}^2
    - \frac{1-\delta_1}{2}f\abs*{\nabla_{\nu_{\Sigma}}\psi+ \frac{1}{r}\psi}^2
    .
  \end{equation*}
  We can compute that %
  \begin{align*}
    &J + \frac{p}{2}r^{-2}\nu_{\Sigma}(rf\abs*{\psi}^2)\\
    ={}& f\left( 1-\mathcal{B} \right)\abs*{\nabla_{\nu_{\Sigma}}\psi}^2
    + f r^{-1}\psi\cdot\nu_{\Sigma}\psi
    - \frac{1}{2}\nu_{\Sigma}(fr^{-1})\abs*{\psi}^2
         - \frac{1-\delta_1}{2}f\abs*{\nabla_{\nu_{\Sigma}}\psi+ \frac{1}{r}\psi}^2\\
    & + pr^{-1}f\psi\cdot\nabla_{\nu_{\Sigma}}\psi
      + \frac{p}{2}r^{-2}\nu_{\Sigma}(rf)\abs*{\psi}^2\\
    ={}& \left(\frac{1+\delta_1}{2}-\mathcal{B}\right)\abs*{\nabla_{\nu_{\Sigma}}\psi}^2
         + (p+\delta_1)r^{-1}f\psi\cdot\nabla_{\nu_{\Sigma}}\psi
         + \left(-\frac{1-\delta_1}{2} + \frac{p^2+1}{2}\nu_{\Sigma}(r)\right)r^{-2}f\abs*{\psi}^2.
  \end{align*}
  Redefining $\mathcal{B}$, we can write that
  \begin{align*}
    &J + \frac{p}{2}r^{-2}\nu_{\Sigma}(rf\abs*{\psi}^2)\\
    ={}& \frac{1+\delta_1}{2}f\left(1 -\mathcal{B}\right)\abs*{\nabla_{\nu_{\Sigma}}\psi}^2
         + (p+\delta_1)r^{-1}f\psi\cdot\nabla_{\nu_{\Sigma}}\psi
         + \left(-\frac{1-\delta_1}{2} + \frac{p^2+1}{2}\nu_{\Sigma}(r)\right)r^{-2}f\abs*{\psi}^2\\
    ={}& \frac{1+\delta_1}{2}f(1-\mathcal{B})\abs*{\nabla_{\nu_{\Sigma}}\psi + \frac{p+\delta_1}{(1-\mathcal{B})(1+\delta_1)}r^{-1}f\psi}^2\\
    &+ \frac{1}{2}\left(
      -\frac{(p+\delta_1)^2}{(1+\delta_1)(1-\mathcal{B})}
      - (1-\delta_1)
      + (p^2+1)\nu_{\Sigma}(r)
      \right)r^{-2}f\abs*{\psi}^2.
  \end{align*}
  Defining
  \begin{equation*}
    L\vcentcolon= -\frac{(p+\delta_1)^2}{(1+\delta_1)(1-\mathcal{B})}
      - (1-\delta_1)
      + (p^2+1)\nu_{\Sigma}(r)
    ,
  \end{equation*}
  we observe that since $\nu_{\Sigma}(r)=1 + O(r^{-2})$, and $\mathcal{B} = O(r^{-1})$, we deduce that
   \begin{align*}
     L ={}& - \frac{(p+\delta_1)^2}{1+\delta_1} - (1-\delta_1) + p^2+1
            + O(R^{-1})
     ={} \frac{(p-1)^2\delta_1}{1+\delta_1} + O(R^{-1}).
   \end{align*}
   For $R$ sufficiently large, and in particular
   $R\gg \frac{1}{(p-1)^2\delta_1}$, we then have that
   \begin{equation*}
     J + \frac{p}{2}r^{-2}\nu_{\Sigma}(rf\abs*{\psi}^2)\ge \frac{(p-1)^2\delta_1}{4}r^{p-2}\abs*{\psi}^2.
   \end{equation*}
   Combining with
   \zcref[noname]{eq:rp:main-boundary:estimates:Sigmatau:f-specific:1-delta:aux1},
   and recalling that $\EMTensor[\psi]_{34} = \abs*{\nabla\psi}^2 + V\abs*{\psi}^2$,
   we have that
   \begin{equation*}
     \begin{split}
       \JCurrent{\rpVF, \rpLagrangianCorr, \rpOneForm}[\psi]\cdot N_{\Sigma}
       \ge{}& c_0(p-1)^2r^{p-2}\abs*{\psi}^2
              + M^2c_0r^{p-2}\abs*{\nabla\psi}^2                            
              - \frac{p}{2}r^{-2}\nu_{\Sigma}(r^{p+1}\abs*{\psi}^2)
               + O(r^{p-3})\left(\abs*{\nabla_3\psi}^2 + \abs*{\psi}^2\right),
     \end{split}     
   \end{equation*}
   where $4c_0 = \min(c_1,\delta_1)>0$. 
   Next, recall from
   \zcref[noname]{eq:rp:boundary:Sigma:general:div:aux-nu-Sigma-and-div} that
   \begin{equation*}
     r^{-2}\nu_{\Sigma}(rf\abs*{\psi}^2)
    = \Divergence_{\Sigma}\left(r^{-1}f\abs*{\psi}^2\nu_{\Sigma}\right)
    + \left(O(r^{-4})+ O(r^{-2})\Lambda\right)f\abs*{\psi}^2.
  \end{equation*}
  As a result, since $\Lambda\lesssim r^{-2}$ we can conclude that for $r$ sufficiently large,
  \begin{equation*}
    \begin{split}
      \JCurrent{\rpVF, \rpLagrangianCorr, \rpOneForm}[\psi]\cdot N_{\Sigma}
      \ge{}& c_0(p-1)^2r^{p-2}\abs*{\psi}^2
              + M^2c_0r^{p-2}\abs*{\nabla\psi}^2              
              - \frac{p}{2}\Divergence_{\Sigma}\left(r^{-1}f\abs*{\psi}^2\nu_{\Sigma}\right) 
              + O(r^{p-3})\left(\abs*{\nabla_3\psi}^2 + \abs*{\psi}^2\right).
    \end{split}
  \end{equation*}
  We now move onto showing
  \zcref[noname]{eq:rp:main-boundary:estimates:SigmaStar:f-specific:1-delta}.
  To this end, recalling the expression for $V$, we have that on $\SigmaStar$
  \begin{equation*}
    V - \frac{\Lambda}{3}p
    = \frac{4}{\abs*{q}^2}\frac{r^2+a^2-2Mr}{r^2}
    + \frac{2-p}{3}\Lambda
    + a^2O\left(r^{-4}\right). 
  \end{equation*}
  In particular, as a result, we have that for $r$ sufficiently large
  with respect to $M$\footnote{Observe that since $r$ is constant on
    $\SigmaStar$, this is equivalent to asking for $\Lambda$
    sufficiently small.}
  \begin{equation*}
    \begin{split}
      \JCurrent{\rpVF, \rpLagrangianCorr, \rpOneForm}[\psi]\cdot N_{\SigmaStar}
      \ge{}& \frac{1}{2}
             \left(
             r^p\abs*{\nabla\psi}^2
             + \frac{4}{\abs*{q}^2}r^p\abs*{\psi}^2
             + \frac{2-p}{3}\Lambda r^p\abs*{\psi}^2
             \right)
             + \frac{-\Delta}{2\abs*{q}^2}r^p\abs*{\widecheck{\nabla}_4\psi}^2\\
           &  + \frac{1}{4}r^{p-2}\lambda\abs*{\nabla_3\psi}^2
             + \frac{1}{2}\Divergence_{\SigmaStar}\left(r^{p-1}\abs*{\psi}^2\nu_{\SigmaStar}\right)
             + Err_{\SigmaStar}^{\bowtie}.
    \end{split}    
  \end{equation*}
  Then, observe that for $r_{\CosmologicalHorizon}> R$, for $R$
  sufficiently large but independent of $\Lambda$, for $f = r^p$,
  $\delta\le p\le 2-\delta$, we can absorb
  $Err_{\SigmaStar}^{\bowtie}$ into the other terms on the right-hand
  side, as desired.
\end{proof}

\paragraph{Proof of \zcref[cap]{prop:rp:Kerr}}

We will now prove \zcref[cap]{prop:rp:Kerr} in the case $s=0$ for
the model equation given by \zcref[noname]{eq:model-problem-gRW}. The proof is similar to the proof of Proposition 10.1.2 in
\cite{giorgiWaveEquationsEstimates2024} and the proof of Theorem 10.37 in
\cite{klainermanGlobalNonlinearStability2020}. The main difference between the present
proof and the proofs in \cite{giorgiWaveEquationsEstimates2024} and
\cite{klainermanGlobalNonlinearStability2020} lies in the more delicate analysis near
the cosmological region in our present setting. 

\paragraph{Step 0: Reducing to the region $r\ge R$.}

We begin by integrating the application of the divergence theorem in
\zcref[noname]{eq:div-thm:general} and write
\begin{equation}
  \label{eq:rp:div-thm-with-rp-quantities:integrated}
  \begin{split}
    &\int_{\Sigma(\tau_2)}\JCurrent{\rpVF, \rpLagrangianCorr, \rpOneForm}[\psi]\cdot N_{\Sigma}
  + \int_{\SigmaStar(\tau_1,\tau_2)}\JCurrent{\rpVF, \rpLagrangianCorr, \rpOneForm}[\psi]\cdot N_{\SigmaStar}
  + \int_{\Manifold(\tau_1,\tau_2)}\KCurrent{\rpVF, \rpLagrangianCorr, \rpOneForm}[\psi]\\
  ={}& \int_{\Sigma(\tau_1)}\JCurrent{\rpVF, \rpLagrangianCorr, \rpOneForm}[\psi]\cdot N_{\Sigma}
  + \mathcal{F}(\tau_1,\tau_2),
  \end{split}
\end{equation}
where
\begin{equation}
  \label{eq:rp:div-thm:extra-terms}
  \begin{split}
    \mathcal{F}(\tau_1,\tau_2)
    \vcentcolon={}&  -\int_{\Manifold(\tau_1,\tau_2)}f_p\nabla_4\psi\cdot N
  - \frac{1}{2}\int_{\Manifold(\tau_1,\tau_2)}r^{-2}f_p\lambda \nabla_3\psi\cdot N
  - \int_{\Manifold(\tau_1,\tau_2)}f_p\nabla_4\psi \cdot \frac{4a\cos\theta}{\abs*{q}^2}\LeftDual{\nabla}_{\KillT}\psi\\
        & - \frac{1}{2}\int_{\Manifold(\tau_1,\tau_2)}r^{-2}f_p\lambda \nabla_3\psi\cdot\frac{4a\cos\theta}{\abs*{q}^2}\LeftDual{\nabla}_{\KillT}\psi\\
        & - \int_{\Manifold(\tau_1,\tau_2)}\frac{r}{\abs*{q}^2}f_p\psi\cdot N
          - \int_{\Manifold(\tau_1,\tau_2)}\frac{r}{\abs*{q}^2}f_p\psi \cdot \frac{4a\cos\theta}{\abs*{q}^2}\LeftDual{\nabla}_{\KillT}\psi
          ,
  \end{split}  
\end{equation}
We will denote the boundary terms,
\begin{equation}
  \label{eq:rp:div-thm:boundary-terms}
  \begin{split}
    J_{\ge R}(\tau_1, \tau_2)
    \vcentcolon={}& \int_{\Sigma_{\ge R}(\tau_2)}\JCurrent{\rpVF, \rpLagrangianCorr, \rpOneForm}[\psi]\cdot N_{\Sigma}
          + \int_{\SigmaStar(\tau_1,\tau_2)}\JCurrent{\rpVF, \rpLagrangianCorr, \rpOneForm}[\psi]\cdot N_{\SigmaStar}
          - \int_{\Sigma_{\ge R}(\tau_1)}\JCurrent{\rpVF, \rpLagrangianCorr, \rpOneForm}[\psi]\cdot N_{\Sigma}\\
    J_{\le R}(\tau_1, \tau_2)
    \vcentcolon={}& \int_{\Sigma_{\le R}(\tau_1)}\JCurrent{\rpVF, \rpLagrangianCorr, \rpOneForm}[\psi]\cdot N_{\Sigma}
          - \int_{\Sigma_{\le R}(\tau_2)}\JCurrent{\rpVF, \rpLagrangianCorr, \rpOneForm}[\psi]\cdot N_{\Sigma},
  \end{split}  
\end{equation}
we then rewrite \zcref[noname]{eq:rp:div-thm-with-rp-quantities:integrated}
\begin{equation}
  \label{eq:rp:div-thm-with-rp-quantities:integrated:alt-form}
  J_{\ge R}(\tau_1,\tau_2)
  + \int_{\Manifold_{\ge R}(\tau_1,\tau_2)}\KCurrent{\rpVF, \rpLagrangianCorr, \rpOneForm}[\psi]
  = J_{\le R}(\tau_1,\tau_2)
  - \int_{\Manifold_{\le R}(\tau_1,\tau_2)}\KCurrent{\rpVF, \rpLagrangianCorr, \rpOneForm}[\psi]
  + \mathcal{F}(\tau_1,\tau_2), 
\end{equation}
which follows from the fact that for
$\left(\rpVF, \rpLagrangianCorr, \rpOneForm\right)$ as chosen in
\zcref[noname]{eq:rp-Kerr:multiplier-choice:general} we have $\rpVF + \rpLagrangianCorr = f_p\widecheck{\nabla}_4 + \frac{1}{2}r^{-2}f_p\lambda \nabla_3$. We first estimate the term $J_{\le R}(\tau_1,\tau_2)
- \int_{\Manifold_{\le R}(\tau_1,\tau_2)}\KCurrent{\rpVF, \rpLagrangianCorr, \rpOneForm}[\psi]$
on the \RHS.

\begin{lemma}
  \label{lemma:rp:final:step0}
  Let $J_{\ge R}$ be as defined in \zcref[noname]{eq:rp:div-thm:boundary-terms}. Then, 
  \begin{equation}
    \label{eq:rp:final:step0}
    \begin{split}
      & J_{\ge R}(\tau_1,\tau_2)
    + \int_{\Manifold_{\ge R}(\tau_1,\tau_2)}\KCurrent{\rpVF, \rpLagrangianCorr, \rpOneForm}[\psi]    
    \lesssim{} R^{p+2}\left(
      E_{\frac{R}{2}\le r\le R}[\psi](\tau_1)
      + \MorNorm_{\frac{R}{2}\le r\le R}[\psi](\tau_1,\tau_2)
    \right)
    + \mathcal{F}(\tau_1,\tau_2).
    \end{split}   
  \end{equation}
\end{lemma}
\begin{proof}
  Using the support of $f_p$, we have that
  \begin{align*}
    \abs*{\int_{\Sigma_{\le R}(\tau_1)}\JCurrent{\rpVF, \rpLagrangianCorr, \rpOneForm}[\psi]\cdot N_{\Sigma}}
    &\lesssim  R^p\EnergyFlux_{\frac{R}{2}\le r\le R}[\psi](\tau_1),\\
    \abs*{\int_{\Sigma_{\le R}(\tau_2)}\JCurrent{\rpVF, \rpLagrangianCorr, \rpOneForm}[\psi]\cdot N_{\Sigma}}
    &\lesssim R^p\EnergyFlux_{\frac{R}{2}\le r\le R}[\psi](\tau_2),\\
    \abs*{\int_{\Manifold_{\le R}(\tau_1,\tau_2)}\rpKCurrent}
    &\lesssim R^{p+2}\MorNorm_{\frac{R}{2}\le r\le R}[\psi](\tau_1,\tau_2).
  \end{align*}
  Hence, using the support of $f_p$, we have that 
  \begin{equation*}
    \begin{split}
      & J_{\le R}(\tau_1, \tau_2)
        - \int_{\Manifold_{\le R}(\tau_1,\tau_2)}\rpKCurrent\\
      \lesssim{}& R^{p+2}\left(
                  \EnergyFlux_{\frac{R}{2}\le r\le R}[\psi](\tau_1)
                  + \EnergyFlux_{\frac{R}{2}\le r\le R}[\psi](\tau_2)
                  + \MorNorm_{\frac{R}{2}\le r\le R}[\psi](\tau_1,\tau_2)
                  \right).
    \end{split}    
  \end{equation*}
  Plugging this into
  \zcref[noname]{eq:rp:div-thm-with-rp-quantities:integrated:alt-form}
  directly yields the conclusion of \zcref[cap]{lemma:rp:final:step0}.
\end{proof}

\paragraph{Step 1: Boundary terms.}

We prove the following lower bound for the boundary terms
\begin{equation*}
  J_{\ge R}(\tau_1,\tau_2)
  = \int_{\Sigma_{\ge R(\tau_2)}} \JCurrent{\rpVF, \rpLagrangianCorr, \rpOneForm}[\psi]\cdot N_{\Sigma}
  + \int_{\SigmaStar(\tau_1,\tau_2)}\JCurrent{\rpVF,\rpLagrangianCorr, \rpOneForm}[\psi]\cdot N_{\SigmaStar}
  - \int_{\Sigma_{\ge R(\tau_1)}} \JCurrent{\rpVF, \rpLagrangianCorr, \rpOneForm}[\psi]\cdot N_{\Sigma},
\end{equation*}
by making use of
\zcref[noname]{eq:rp:main-boundary:estimates:Sigmatau:general}.

\begin{lemma}
  \label{lemma:rp:main-proof:boundary-terms}
  Let $J_{\ge R}(\tau_1,\tau_2)$ be as defined in
  \zcref[noname]{eq:rp:div-thm:boundary-terms}.  Then for $R$ sufficiently
  large,
  \begin{equation}
    \label{eq:rp:boundary:first-lower-bound}
    \begin{split}
      J_{\ge R}(\tau_1,\tau_2)
      \gtrsim{}&  \int_{\SigmaStar(\tau_1,\tau_2)}
                 r^p\abs*{\nabla\psi}^2
                 + \frac{-\Delta}{\abs*{q}^2}r^p\abs*{\widecheck{\nabla}_4\psi}^2
                 + r^{p-2}\abs*{\nabla_3\psi}^2\\
               & + \int_{\Sigma_{\ge R}(\tau_2)}r^p\left(\abs*{\widecheck{\nabla}_4\psi}^2 + r^{-2}\abs*{\nabla\psi}^2\right)
                 - \EnergyFluxFar_{p,R}[\psi](\tau_1).
    \end{split}
  \end{equation}
\end{lemma}

\begin{proof}
  Integrating \zcref[noname]{eq:rp:main-boundary:estimates:Sigmatau:general},
  we see that for $R$ sufficiently large,
  \begin{align*}
    \int_{\Sigma_{\ge R}(\tau_2)}\JCurrent{\rpVF,\rpLagrangianCorr,\rpOneForm}[\psi]\cdot N_{\Sigma}
    \gtrsim{}& \int_{\Sigma_{\ge R}(\tau_2)}r^p\left(\abs*{\widecheck{\nabla}_4\psi}^2 + r^{-2}\abs*{\nabla\psi}^2\right)
               + \int_{\Sigma_{\ge R}(\tau_2)}r^{p-4}\lambda^2\abs*{\nabla_3\psi}^2
    \\
             & + \frac{\Lambda}{6}\int_{\Sigma_{\ge R}(\tau_2)}p\lambda r^{p-1}\abs*{\psi}^2
               - \frac{1}{2}\int_{\Sigma_{\ge R}(\tau_2)}\Divergence_{\Sigma}\left(r^{p-1}\abs*{\psi}^2\nu_{\Sigma}\right).
  \end{align*}
  Similarly, integrating
  \zcref[noname]{eq:rp:main-boundary:estimates:SigmaStar:f-specific:1-delta},
  we have that for $R$ sufficiently large,
  \begin{align*}
    \int_{\SigmaStar}\rpJCurrent\cdot N_{\SigmaStar}
    \ge{}& \frac{1}{2}\int_{\SigmaStar(\tau_1,\tau_2)}r^p\left(
           \abs*{\nabla\psi}^2
           + \frac{\lambda}{2r^{-2}}\abs*{\nabla_3\psi}^2
           + \frac{p-2}{3}\Lambda\abs*{\psi}^2
           + \frac{-\Delta}{\abs*{q}^2}\abs*{\widecheck{\nabla}_4\psi}^2
           \right)\\
         & + \frac{1}{2}\int_{\SigmaStar(\tau_1,\tau_2)} \Divergence_{\SigmaStar}\left(r^{p-1}\abs*{\psi}^2\nu_{\SigmaStar}\right)
          .
  \end{align*}
  Applying the divergence theorem on
  $\Sigma(\tau_2)\bigcup \SigmaStar(\tau_1,\tau_2)$ to
  $\Divergence_{\Sigma}(r^{p-1}\abs*{\psi}^2\nu_{\Sigma})$, we see that
  \begin{equation}
    \label{eq:rp:main-proof:boundary-terms:Sigma-tau-sphere-divergence}
    \int_{\Sigma(\tau)}\Divergence_{\Sigma}\left(r^{p-1}\abs*{\psi}^2\nu_{\Sigma}\right)
    = \int_{S_{*}(\tau)}r^{p-1}\abs*{\psi}^2 - \int_{S_R(\tau)}r^{p-1}\abs*{\psi}^2.
  \end{equation}
  Similarly, we have that 
  \begin{equation}
    \label{eq:rp:main-proof:boundary-terms:SigmaStar-sphere-divergence}
    \int_{\SigmaStar(\tau_1,\tau_2)}\Divergence_{\SigmaStar}\left(r^{p-1}\abs*{\psi}^2\nu_{\SigmaStar}\right)  
    = \int_{S_{*}(\tau_2)}r^{p-1}\abs*{\psi}^2
    - \int_{S_{*}(\tau_1)}r^{p-1}\abs*{\psi}^2.
  \end{equation}
  As a result, notice that the integrals over $S_{*}(\tau_1)$ and
  $S_{*}(\tau_2)$ cancel out. Thus, we have that
  \begin{equation*}
    \begin{split}
      & \int_{\Sigma_{\ge R}(\tau_2)}\JCurrent{\rpVF,\rpLagrangianCorr,\rpOneForm}[\psi]\cdot N_{\Sigma}
      + \int_{\SigmaStar}\JCurrent{\rpVF,\rpLagrangianCorr,\rpOneForm}[\psi]\cdot N_{\SigmaStar}\\
    \gtrsim{}& \int_{\Sigma_{\ge R}(\tau_2)}r^p\left(\abs*{\widecheck{\nabla}_4\psi}^2 + r^{-2}\abs*{\nabla\psi}^2\right)
               + \int_{\Sigma_{\ge R}(\tau_2)}r^{p-4}\lambda^2\abs*{\nabla_3\psi}^2
    \\
             & + \frac{\Lambda}{6}\int_{\Sigma_{\ge R}(\tau_2)}p\lambda r^{p-1}\abs*{\psi}^2
               + \frac{1}{2}\int_{S_R(\tau_2)}r^{p-1}\abs*{\psi}^2
               - \frac{1}{2}\int_{S_*(\tau_1)}r^{p-1}\abs*{\psi}^2\\
      & + \frac{1}{2}\int_{\SigmaStar(\tau_1,\tau_2)}r^p\left(
           \abs*{\nabla\psi}^2
           + \frac{\lambda}{2r^{-2}}\abs*{\nabla_3\psi}^2
           + \frac{p-2}{3}\Lambda\abs*{\psi}^2
           + \frac{-\Delta}{\abs*{q}^2}\abs*{\widecheck{\nabla}_4\psi}^2
           \right).
    \end{split}
  \end{equation*}
  Thus for $R$ sufficiently large, and
  $\delta\le p\le 2-\delta$,
  \begin{equation*}
    \begin{split}
      J_{\ge R}(\tau_1,\tau_2)
      \gtrsim{}& \int_{\SigmaStar(\tau_1,\tau_2)}r^p
                 \abs*{\nabla\psi}^2
                 + \frac{-\Delta}{\abs*{q}^2}r^p\abs*{\widecheck{\nabla}_4\psi}^2
                 + \lambda r^{p-2}\abs*{\nabla_3\psi}^2\\
               & + \int_{\Sigma_{\ge R}(\tau_2)}r^p\left(\abs*{\widecheck{\nabla}_4\psi}^2 + r^{-2}\abs*{\nabla\psi}^2\right)
                 - \EnergyFluxFar_{p,R}[\psi](\tau_1)
                 .
    \end{split}  
  \end{equation*}
  as desired, concluding the proof of \zcref{lemma:rp:main-proof:boundary-terms}.
\end{proof}

Observe that in the estimate in \zcref[cap]{lemma:rp:main-proof:boundary-terms}, $J_{\ge R}(\tau_1,\tau_2)$
does not control the $L^2$ norm of $\psi$ on the boundary. To recover
the $L^2$ norm, we will have to treat separately the case where
$p\le 1-\delta$, and the case where $1-\delta\le p \le 2-\delta$. 
\begin{lemma}
  \label{lemma:rp:main-proof:boundary-terms:1-delta}
  For $\delta\le p \le 1-\delta$, we have that
  \begin{equation}
    \label{eq:rp:main-proof:boundary-terms:1-delta}
    \begin{split}
      J_{\ge R}(\tau_1,\tau_2)
    \gtrsim{}& \int_{\Sigma_{\ge R}(\tau_2)}r^p\left(
      \abs{\widecheck{\nabla}_4\psi}^2
    + r^{-2}\abs*{\nabla\psi}^2
    + r^{-2}\abs*{\psi}^2\right)\\
  & + \int_{\SigmaStar(\tau_1,\tau_2)}r^p\left(
    \abs*{\nabla\psi}^2
    + r^{-2}\abs*{\psi}^2
    + r^{-2}\abs*{\nabla_3\psi}^2
  \right)
  - \EnergyFluxFar_{p,R}[\psi](\tau_1).
    \end{split}    
  \end{equation}
\end{lemma}
\begin{proof}
  Integrating, \zcref[noname]{eq:rp:main-boundary:estimates:Sigmatau:f-specific:1-delta} we have that
  \begin{equation}
    \label{eq:rp:main-proof:boundary-terms:1-delta:Sigma-integrated}
    \begin{split}
      \int_{\Sigma_{\ge R}(\tau_2)}\rpJCurrent\cdot N_{\Sigma}
      \ge{}& \int_{\Sigma_{\ge R}(\tau_2)}\left( \frac{\delta^2}{8}r^{p-2}\abs*{\psi}^2
             + \frac{1}{2}r^{p-2}\lambda\abs*{\nabla\psi}^2
             - \frac{p}{2}\Divergence_{\Sigma}\left(r^{p-1}\abs*{\psi}^2\nu_{\Sigma}\right) \right)\\
           & + \int_{\Sigma_{\ge R}(\tau_2)}O(r^{p-3})\left(\abs*{\nabla_3\psi}^2 + \abs*{\psi}^2\right)
             .
    \end{split}    
  \end{equation}
  Similarly, integrating
  \zcref[noname]{eq:rp:main-boundary:estimates:SigmaStar:f-specific:1-delta},
  we have that
  \begin{equation}
    \label{eq:rp:main-proof:boundary-terms:1-delta:SigmaStar-integrated}
    \begin{split}
      \int_{\SigmaStar}\rpJCurrent\cdot N_{\SigmaStar}
      \ge{}& \int_{\SigmaStar} \frac{1}{2}\left(
             r^p\abs*{\nabla\psi}^2
             + 4 r^{p-2}\abs*{\psi}^2
             \right)
             + \int_{\SigmaStar}\frac{-\Delta}{2\abs*{q}^2}r^p\abs*{\widecheck{\nabla}_4\psi}^2\\
           &  + \int_{\SigmaStar}\frac{1}{4}r^{p-2}\lambda\abs*{\nabla_3\psi}^2
             + \frac{1}{2}\int_{\SigmaStar}\Divergence_{\SigmaStar}\left(r^{p-1}\abs*{\psi}^2\nu_{\SigmaStar}\right)\\
      & + \int_{\SigmaStar}O(r^{p-3})\left(\abs*{\nabla_3\psi}^2 + \abs*{\psi}^2\right).
    \end{split}
  \end{equation}
  Applying again the divergence theorem on
  $\Sigma(\tau_2)\bigcup \SigmaStar(\tau_1,\tau_2)$ to
  $\Divergence_{\Sigma}(r^{p-1}\abs*{\psi}^2\nu_{\Sigma})$, we again see that
  \begin{align*}
    \int_{\Sigma_{r\ge R}(\tau)}\Divergence_{\Sigma}\left(r^{p-1}\abs*{\psi}^2\nu_{\Sigma}\right)
    ={}& \int_{S_{*}(\tau)}r^{p-1}\abs*{\psi}^2 - \int_{S_R(\tau)}r^{p-1}\abs*{\psi}^2
         + \int_{\Sigma_{r\ge R}(\tau)}O(r^{p-3})\abs*{\psi}^2
         ,\\
    \int_{\SigmaStar(\tau_1,\tau_2)}\Divergence_{\SigmaStar}\left(r^{p-1}\abs*{\psi}^2\nu_{\SigmaStar}\right)  
    ={}& \int_{S_{*}(\tau_2)}r^{p-1}\abs*{\psi}^2
         - \int_{S_{*}(\tau_1)}r^{p-1}\abs*{\psi}^2
         + \int_{\Sigma_{*}(\tau_1, \tau_2)}O(r^{p-3})\abs*{\psi}^2
         .
  \end{align*}
  Combining
  \zcref[noname]{eq:rp:main-proof:boundary-terms:1-delta:Sigma-integrated},
  \zcref[noname]{eq:rp:main-proof:boundary-terms:1-delta:SigmaStar-integrated},
  and the applications of the divergence theorem on
  $\Sigma(\tau_2)\bigcup \SigmaStar(\tau_1,\tau_2)$ above, we have
  that
  \begin{align*}
    &\int_{\Sigma_{\ge R(\tau_2)}} \JCurrent{\rpVF, \rpLagrangianCorr, \rpOneForm}[\psi]\cdot N_{\Sigma}
      + \int_{\SigmaStar(\tau_1,\tau_2)}\JCurrent{\rpVF,\rpLagrangianCorr, \rpOneForm}[\psi]\cdot N_{\SigmaStar}\\
    \ge{}& \int_{\Sigma_{\ge R}(\tau_2)}\left( \frac{\delta^2}{8}r^{p-2}\abs*{\psi}^2
           + \frac{1}{2}r^{p-2}\lambda\abs*{\nabla\psi}^2\right)
           +  \int_{\SigmaStar} \frac{1}{2}\left(
           r^p\abs*{\nabla\psi}^2
           + 4r^{p-2}\abs*{\psi}^2
           \right)\\
    & + \int_{\SigmaStar}\frac{1}{4}r^{p-2}\lambda\abs*{\nabla_3\psi}^2
      + \frac{p}{2}\int_{S_R(\tau_2)}r^{p-1}\abs*{\psi}^2
      - \frac{1}{2}\int_{S_{*}(\tau_1)}r^{p-1}\abs*{\psi}^2\\
    & + \frac{1-p}{2}\int_{S_{*}(\tau_2)}r^{p-1}\abs*{\psi}^2
      +  \int_{\Sigma_{\ge R}(\tau_2)}O(r^{p-3})\left(\abs*{\nabla_3\psi}^2 + \abs*{\psi}^2\right)
      + \int_{\SigmaStar}O(r^{p-3})\left(\abs*{\nabla_3\psi}^2 + \abs*{\psi}^2\right)
      .
  \end{align*}
  Since we are only considering $\delta\le p\le 1-\delta$, we then
  have that the first term in the last line in the above equation is actually positive,
  and can be excluded, so we in fact have that for $\delta\le p\le 1-\delta$, and $R$ sufficiently large that
  \begin{equation}
    \label{eq:rp:main-proof:boundary-terms:1-delta:aux:L2-conclusion}
    \begin{split}
      &\int_{\Sigma_{\ge R(\tau_2)}} \JCurrent{\rpVF, \rpLagrangianCorr, \rpOneForm}[\psi]\cdot N_{\Sigma}
        + \int_{\SigmaStar(\tau_1,\tau_2)}\JCurrent{\rpVF,\rpLagrangianCorr, \rpOneForm}[\psi]\cdot N_{\SigmaStar}\\
      \gtrsim{}& \int_{\Sigma_{\ge R}(\tau_2)}\left( \frac{\delta^2}{8}r^{p-2}\abs*{\psi}^2
             + \frac{1}{2}r^{p-2}\lambda\abs*{\nabla\psi}^2\right)
             +  \int_{\SigmaStar} \frac{1}{2}\left(
             r^p\abs*{\nabla\psi}^2
             + 4r^{p-2}\abs*{\psi}^2
             \right)\\
      & + \int_{\SigmaStar}\frac{1}{4}r^{p-2}\lambda\abs*{\nabla_3\psi}^2       
        - \EnergyFluxFar_{p,R}[\psi](\tau_1)
        .
    \end{split}
  \end{equation}
  Then combining
  \zcref[noname]{eq:rp:main-proof:boundary-terms:1-delta:aux:L2-conclusion} 
  with \zcref[noname]{eq:rp:boundary:first-lower-bound} then yields
  \zcref[noname]{eq:rp:main-proof:boundary-terms:1-delta} for
  $\delta\le p\le 1-\delta$.
\end{proof}

\paragraph{Step 2: Bulk terms for $r\ge R$.}

We prove the following lower bound for
$\int_{\Manifold_{\ge R}(\tau_1,\tau_2)}\rpKCurrent$.

\begin{lemma}
  \label{lemma:rp:main-proof:bulk:lower-bound}
  Given a fixed $\delta>0$, we have that for all
  $\delta\le p\le 2-\delta$ and $R\gg \frac{M}{\delta}$,
  $\epsilon\ll \delta$,
  \begin{equation}
    \label{eq:rp:main-proof:bulk:lower-bound}
    \begin{split}
      \int_{\Manifold_{\ge R}(\tau_1,\tau_2)}\rpKCurrent
      \ge{}& \frac{1}{4}\int_{\Manifold_{\ge R}(\tau_1,\tau_2)}
             r^{p-1}\left(p\abs*{
             \widecheck{\nabla}_4(\psi)}^2
             + (2-p)\left(
             \abs*{\nabla \psi}^2
             + r^{-2}\abs*{\psi}^2
             \right)
             \right)\\
           & + O\left(R^{-1}\right)\int_{\Manifold_{\ge R}(\tau_1,\tau_2)}r^{p-3}\abs*{\nabla_3\psi}^2
             .
    \end{split}    
  \end{equation}
\end{lemma}

\begin{proof}
  Integrating \zcref[noname]{eq:rp:bulk:specific} over
  $\Manifold_{\ge R}(\tau_1,\tau_2)$ directly yields the conclusion.
\end{proof}

\paragraph{Step 3: Preliminary control using estimates on $\nabla_3\psi$.}

We prove the following Lemma using \zcref[cap]{prop:rp:preliminary-nabla-3-estimate}.
\begin{lemma}
  \label{lemma:rp:main-proof:combined-with-nabla-3}
  For $\delta\le 1-p$ and $R$ sufficiently large,
  \begin{equation}
    \label{eq:rp:main-proof:combined-with-nabla-3}
    \begin{split}
      \WeightedBEFNorm{p,\ge R}[\psi](\tau_1,\tau_2)
      \lesssim{}& \EnergyFluxWeighted[\psi](\tau_1)
                  + \ForcingTermCombinedNorm{}[\psi, N](\tau_1,\tau_2)
      \\
                &  + R^{p+2}\left(
                  \EnergyFlux_{\frac{R}{2}\le r \le R}[\psi](\tau_1)
                  + \MorNorm_{\frac{R}{2}\le r\le R}[\psi](\tau_1,\tau_2)
                  \right)
                  + \mathcal{F}(\tau_1,\tau_2).
    \end{split}    
  \end{equation}
  On the other hand, for $1-\delta<p\le 2-\delta$, we have that
  \begin{equation}
    \label{eq:rp:main-proof:combined-with-nabla-3:2-delta}
    \begin{split}
      \WeightedBEFNorm{p,\ge R}[\psi](\tau_1,\tau_2)
      \lesssim{}& \EnergyFluxWeighted[\psi](\tau_1)
                  + \ForcingTermCombinedNorm{}[\psi, N](\tau_1,\tau_2)
      \\
                &  + R^{p+2}\left(
                  \EnergyFlux_{\frac{R}{2}\le r \le R}[\psi](\tau_1)
                  + \MorNorm_{\frac{R}{2}\le r\le R}[\psi](\tau_1,\tau_2)
                  \right)\\
               &  + \mathcal{F}(\tau_1,\tau_2)
                  + \EnergyFluxFar_{1-\delta}[\psi](\tau_2)
                  .
    \end{split}
  \end{equation}
\end{lemma}
\begin{proof}
  We begin by proving \zcref[noname]{eq:rp:main-proof:combined-with-nabla-3}
  and considering only $\delta\le p\le 1-\delta$. To this end,
  consider adding \zcref[noname]{eq:rp:main-proof:boundary-terms:1-delta},
  \zcref[noname]{eq:rp:boundary:first-lower-bound} and
  \zcref[noname]{eq:rp:main-proof:bulk:lower-bound}. This yields
  \begin{equation}
    \label{eq:rp:main-proof:combined-with-nabla-3:aux-1}
    \begin{split}
        \rpBulkWeighted{p}{R}[\psi](\tau_1,\tau_2)
        + \EnergyFluxWeighted[\psi](\tau_2)
        + \SpacelikeFluxWeighted[\psi](\tau_1,\tau_2)
      &\lesssim{} \int_{\Manifold_{\ge R}(\tau_1,\tau_2)}\rpKCurrent
                  + J_{\ge R}(\tau_1,\tau_2)\\
                & + O(R^{-1})\int_{\Manifold_{\ge R}(\tau_1,\tau_2)}r^{p-3}\abs*{\nabla_3\psi}^2
                  + \EnergyFluxWeighted[\psi](\tau_1).
    \end{split}    
  \end{equation}
  Then, applying 
  \zcref{prop:rp:preliminary-nabla-3-estimate}, we see that
  \begin{equation}
    \label{eq:rp:step-2:aux-1}
    \begin{split}
      \int_{\Manifold(\tau_1,\tau_2)}r^{-1-\delta}\abs*{\nabla_3\psi}^2
      \lesssim{}& \int_{\Manifold(\tau_1,\tau_2)}r^{-1-\delta}\left(
                  \abs*{\nabla_4\psi}^2 + \abs*{\nabla\psi}^2 + r^{-2}\abs*{\psi}^2
                  \right)\\
                & + \int_{\Sigma(\tau_2)}a^2r^{-2-\delta}\abs*{\nabla_4\psi}^2
                 + \int_{\Sigma(\tau_1)}\abs*{\frac{-\Delta}{\abs*{q}^2}}r^{-\delta}\left(\abs*{\widecheck{\nabla}_4\psi}^2 + ar^{-2}\abs*{\nabla\psi}^2 \right)\\
                & + \int_{\SigmaStar(\tau_1,\tau_2)}\frac{a}{\abs*{q}^2}r^{-\delta}\abs*{\nabla\psi}^2
                  + \frac{\Delta^2}{\abs*{q}^4}r^{-\delta}\abs*{\widecheck{\nabla}_4\psi}^2\\
                & + \EnergyHorizonDeg[\psi](\tau_2)
                  + \MorNorm_{\frac{R}{2}\le r\le R}[\psi](\tau_1,\tau_2)
                  + \ForcingTermCombinedNorm{}[\psi, N](\tau_1,\tau_2). 
    \end{split}
  \end{equation}
  In particular, we can rewrite
  \begin{equation}
    \label{eq:rp:step-2:aux-2}
    \begin{split}
      \int_{\Manifold}r^{-1-\delta}\abs*{\nabla_3\psi}^2
      \lesssim{}& R^{-2\delta}\rpBulkWeighted{p}{R}[\psi](\tau_1,\tau_2)
                  + R^{-2}\EnergyFluxWeighted[\psi](\tau_2)
                  + R^{-2\delta}\EnergyFluxWeighted[\psi](\tau_1)
                  + R^{-2\delta}\SpacelikeFluxWeighted[\psi](\tau_1,\tau_2)\\
                & + \EnergyHorizonDeg[\psi](\tau_2)
                  + \MorNorm_{\frac{R}{2}\le r\le R}[\psi](\tau_1,\tau_2)
                  + \ForcingTermCombinedNorm{}[\psi, N](\tau_1,\tau_2). 
    \end{split}
  \end{equation}
  Combining
  \zcref[noname]{eq:rp:div-thm-with-rp-quantities:integrated:alt-form},
  \zcref[noname]{eq:rp:final:step0},
  \zcref[noname]{eq:rp:main-proof:combined-with-nabla-3:aux-1} and
  \zcref[noname]{eq:rp:step-2:aux-2}, we see that for $R$ sufficiently large,
  we have that
  \begin{equation*}
    \begin{split}
      \WeightedBEFNorm{p,\ge R}[\psi](\tau_1,\tau_2)
      \lesssim{}&  R^{-2\delta}\rpBulkWeighted{p}{R}[\psi](\tau_1,\tau_2)
                  + R^{-2\delta}\EnergyFluxWeighted[\psi](\tau_2)
                  + R^{-2\delta}\SpacelikeFluxWeighted[\psi](\tau_1,\tau_2)\\
                & + \EnergyHorizonDeg[\psi](\tau_2)                  
                  + \EnergyFluxWeighted[\psi](\tau_1)
                                    + \mathcal{F}(\tau_1,\tau_2)
                  + \ForcingTermCombinedNorm{}[\psi, N](\tau_1,\tau_2)
      \\
                &  + R^{p+2}\left(
                  \EnergyFlux_{\frac{R}{2}\le r \le R}[\psi](\tau_1)
                  + \MorNorm_{\frac{R}{2}\le r\le R}[\psi](\tau_1,\tau_2)
                  \right)
      . 
    \end{split}
  \end{equation*}
  For $R$ sufficiently large, the first line of the right-hand side
  can be absorbed into
  $\WeightedBEFNorm{p,\ge R}[\psi](\tau_1,\tau_2)$, and as a result,
  \begin{equation*}
    \begin{split}
      \WeightedBEFNorm{p,\ge R}[\psi](\tau_1,\tau_2)
      \lesssim{}& \EnergyHorizonDeg[\psi](\tau_2)
                  + \EnergyFluxWeighted[\psi](\tau_1)
                  + \ForcingTermCombinedNorm{}[\psi, N](\tau_1,\tau_2)
      \\
                &  + R^{p+2}\left(
                  \EnergyFlux_{\frac{R}{2}\le r \le R}[\psi](\tau_1)
                  + \MorNorm_{\frac{R}{2}\le r\le R}[\psi](\tau_1,\tau_2)
                  \right)
                  + \mathcal{F}(\tau_1,\tau_2)
                  ,
    \end{split}
  \end{equation*}
  as desired. This proves \zcref[noname]{eq:rp:main-proof:combined-with-nabla-3}. To prove \zcref[noname]{eq:rp:main-proof:combined-with-nabla-3:2-delta}, it
  suffices to observe that the only control missing in the case of
  $1-\delta<p\le 2-\delta$ is control over the $L^2$ terms on the
  boundary. But from the definition of
  $\EnergyFluxFar_{p,R}[\psi](\tau)$ in
  \zcref[noname]{eq:energy-flux:weighted:def}
  for $1-\delta<p\le 2-\delta$ we have $\EnergyFluxFar_{p,R}[\psi](\tau)
    \le \EnergyFluxFar_{1-\delta,R}[\psi](\tau)$ which yields the desired estimate.
\end{proof}

\paragraph{Controlling the forcing terms.}
We now show how to control the $\mathcal{F}(\tau_1,\tau_2)$ term in
\zcref{eq:rp:div-thm:extra-terms}. 
\begin{lemma}
  \label{lemma:rp:extra-forcing-terms}
  Let $\mathcal{F}(\tau_1,\tau_2)$ be as defined in
  \zcref[noname]{eq:rp:div-thm:extra-terms}. Then
  \begin{equation}
    \label{eq:rp:extra-forcing-terms}
    \begin{split}
      \abs*{\mathcal{F}(\tau_1,\tau_2)}
    \lesssim{}& \ForcingTermWeightedNorm{p}{\ge \frac{R}{2}}[\psi,N](\tau_1,\tau_2)
    + R^{p+2}\MorNorm_{\frac{R}{2}\le r \le R}[\psi](\tau_1,\tau_2)
    + R^{-\delta}\BulkNormWeighted{p}[\psi](\tau_1,\tau_2)\\
    &+ R^{-\delta}\EnergyFluxWeighted[\psi](\tau_2)
    + R^{-\delta}\EnergyFluxWeighted[\psi](\tau_1)
    + R^{-\delta}\SpacelikeFluxWeighted[\psi](\tau_1,\tau_2).
    \end{split}    
  \end{equation}
\end{lemma}
\begin{proof}
  Recalling
  \zcref{eq:rp:div-thm:extra-terms}, we have $\mathcal{F}(\tau_1,\tau_2)
    = \mathcal{F}_1(\tau_1,\tau_2) + \mathcal{F}_2(\tau_1,\tau_2) + \mathcal{F}_3(\tau_1,\tau_2)$ with
  \begin{align*}
    \mathcal{F}_1(\tau_1,\tau_2)
    \vcentcolon={}& -\int_{\Manifold(\tau_1,\tau_2)}f_p\nabla_4\psi\cdot N
                    - \frac{1}{2}\int_{\Manifold(\tau_1,\tau_2)}r^{-2}f_p\lambda \nabla_3\psi\cdot N
                    - \int_{\Manifold(\tau_1,\tau_2)}\frac{r}{\abs*{q}^2}f_p\psi\cdot N
                    ,
    \\
    \mathcal{F}_2(\tau_1,\tau_2)
    \vcentcolon={}& - f_p\int_{\Manifold(\tau_1,\tau_2)}\left( \nabla_4\psi + \frac{r}{\abs*{q}^2}\psi \right) \cdot \frac{4a\cos\theta}{\abs*{q}^2}\LeftDual{\nabla}_{\KillT}\psi,\\
    \mathcal{F}_3(\tau_1,\tau_2)
    \vcentcolon={}&  - \frac{1}{2}\int_{\Manifold(\tau_1,\tau_2)}r^{-2}f_p\lambda \nabla_3\psi \cdot \frac{4a\cos\theta}{\abs*{q}^2}\LeftDual{\nabla}_{\KillT}\psi.
  \end{align*}
  Recalling \zcref[noname]{eq:forcing-term-weighted-norms:def}, we directly get $\abs*{\mathcal{F}_1(\tau_1,\tau_2)} \lesssim \ForcingTermWeightedNorm{p}{\frac{R}{2}}[\psi,N]\left(\tau_1,\tau_2\right)$. To bound $\mathcal{F}_3(\tau_1,\tau_2)$, recall
  \zcref[noname]{eq:T-R-Z:principal-outgoing:expression:T-Z}. Then,
  \begin{align*}
    &-\frac{2a\cos\theta}{\abs*{q}^2}r^{-2}f_p\lambda\nabla_3\psi\cdot\LeftDual{\nabla}_{\KillT}\psi\\
    ={}& -\frac{a\cos\theta}{\abs*{q}^2}r^{-2}f_p\lambda\nabla_3\psi\cdot\LeftDual{\nabla}_3\psi\\
    & + O\left(r^{-4}\right)f_p\abs*{\widecheck{\nabla}_4\psi}\cdot\left(\left(1- \frac{\Lambda}{3}r^2 + O\left(r^{-1}\right) + O\left(r\right)\Lambda\right)\abs*{\nabla_4\psi} + O\left(r^{-1}\right)\abs*{\nabla\psi}\right).
  \end{align*}
  Then, recalling the definition of
  $\rpBulkWeighted{p}{R}[\psi](\tau_1,\tau_2)$, we get $\abs*{\mathcal{F}_3(\tau_1,\tau_2)} \lesssim R^{-1}\rpBulkWeighted{p}{R}[\psi](\tau_1,\tau_2)$. To bound $\mathcal{F}_2(\tau_1,\tau_2)$, recall 
  \zcref[noname]{eq:T-R-Z:principal-outgoing:expression:T-Z} again to obtain 
  \begin{align*}
    &-\frac{4a\cos\theta}{\abs*{q}^2}f_p\left( \nabla_4\psi + \frac{r}{\abs*{q}^2}\psi \right)
      \cdot\LeftDual{\nabla}_{\KillT}\psi
      \\
    ={}& -\frac{2a\cos\theta}{\abs*{q}^2}f_p\widecheck{\nabla}_4\psi \cdot\LeftDual{\nabla}_3\psi\\
       & + O\left(r^{-2}\right)f_p\left( \abs*{\nabla_4\psi} + r^{-1}\abs*{\psi} \right)\cdot\left(\left(1- \frac{\Lambda}{3}r^2 + O\left((1+\Lambda r^2)r^{-1}\right) \Lambda\right)\abs*{\nabla_4\psi} + O\left(r^{-1}\right)\abs*{\nabla\psi}\right).
  \end{align*}
  As a result, we have that
  \begin{align*}
    \mathcal{F}_2\left(\tau_1,\tau_2\right)
    ={}& -\int_{\Manifold(\tau_1,\tau_2)}
         f_p\widecheck{\nabla}_4\psi \cdot\left(\frac{2a\cos\theta}{\abs*{q}^2}\LeftDual{\nabla}_3\psi\right)
          + R^{-1}\int_{\Manifold_{r\ge R}(\tau_1,\tau_2)}r^{p-1}\left(
         \abs*{\widecheck{\nabla}_4\psi}^2 + r^{-2}\abs*{\psi}^2 + \abs*{\nabla\psi}^2 
         \right)\\
       & +
         R^{-1}\Lambda\int_{\Manifold_{r\ge R}(\tau_1,\tau_2)}r^{p+1}\left(
         \abs*{\widecheck{\nabla}_4\psi}^2 + r^{-2}\abs*{\psi}^2 + \abs*{\nabla\psi}^2
         \right)
         \\ \lesssim{}& \underbrace{-\int_{\Manifold(\tau_1,\tau_2)}
         f_p\widecheck{\nabla}_4\psi \cdot\left(\frac{2a\cos\theta}{\abs*{q}^2}\LeftDual{\nabla}_3\psi\right)}_{=\vcentcolon I} + R^{-1}\BulkNormWeighted{p}[\psi](\tau_1,\tau_2),
  \end{align*}
  where we used $\Lambda\lesssim r^{-2}$. To estimate $|I|$ we observe that for $p\le 2-\delta$,
  \begin{align*}
    I\lesssim{}& \int_{\Manifold_{r\ge \frac{R}{2}}(\tau_1,\tau_2)}r^{p-2}\abs*{\widecheck{\nabla}_4\psi}\abs*{\nabla_3\psi}\\
    \lesssim{}&\left(\int_{\Manifold_{r\ge \frac{R}{2}}(\tau_1,\tau_2)}r^{p-1}\abs*{\widecheck{\nabla}_4\psi}^2\right)^{\frac{1}{2}}
                \left(\int_{\Manifold_{r\ge \frac{R}{2}}(\tau_1,\tau_2)}r^{p-3}\abs*{\nabla_3\psi}^2\right)^{\frac{1}{2}}\\
    \lesssim{}&R^{-\frac{2-\delta-p}{2}}\left(\int_{\Manifold_{r\ge \frac{R}{2}}(\tau_1,\tau_2)}r^{p-1}\abs*{\widecheck{\nabla}_4\psi}^2\right)^{\frac{1}{2}}
                \left(\int_{\Manifold_{r\ge \frac{R}{2}}(\tau_1,\tau_2)}r^{-1-\delta}\abs*{\nabla_3\psi}^2\right)^{\frac{1}{2}}.
  \end{align*}
  Then, using \zcref[cap]{prop:rp:preliminary-nabla-3-estimate}
  with $s=\delta$, we have that
  \begin{align*}
    I\lesssim{}& R^{-\frac{2-\delta-p}{2}}\int_{\Manifold_{\ge \frac{R}{2}}(\tau_1,\tau_2)}r^{p-1}\abs*{\widecheck{\nabla}_4\psi}^2
                 +  R^{-\frac{2-\delta-p}{2}-2\delta}\rpBulkWeighted{p}{R}[\psi](\tau_1,\tau_2)\\
               & + R^{-\frac{2-\delta-p}{2}-2\delta}\EnergyFluxWeighted[\psi](\tau_2)
                 + R^{-\frac{2-\delta-p}{2}-2\delta}\EnergyFluxWeighted[\psi](\tau_1)
                 + R^{-\frac{2-\delta-p}{2}-2\delta}\SpacelikeFluxWeighted[\psi](\tau_1,\tau_2).
  \end{align*}
  This concludes the proof of \zcref[cap]{lemma:rp:extra-forcing-terms}. 
\end{proof}

\paragraph{Conclusion of the proof of \zcref[cap]{prop:rp:Kerr}.}

Combining \zcref[cap]{lemma:rp:main-proof:combined-with-nabla-3} and
\zcref[cap]{lemma:rp:extra-forcing-terms}, we have that
\begin{equation*}
  \begin{split}
    \WeightedBEFNorm{p,\ge R}[\psi](\tau_1,\tau_2)
      \lesssim{}& \EnergyFluxWeighted[\psi](\tau_1)
                  + \ForcingTermCombinedNorm{p}[\psi, N](\tau_1,\tau_2)
                  + R^{p+2}\left(
                  \EnergyFlux_{\frac{R}{2}\le r \le R}[\psi](\tau_1)
                  + \MorNorm_{\frac{R}{2}\le r\le R}[\psi](\tau_1,\tau_2)
                  \right)\\
    &+ R^{-\delta}\BulkNormWeighted{p}[\psi](\tau_1,\tau_2)
    + R^{-\delta}\EnergyFluxWeighted[\psi](\tau_2)%
      + R^{-\delta}\EnergyFluxWeighted[\psi](\tau_1)
    + R^{-\delta}\SpacelikeFluxWeighted[\psi](\tau_1,\tau_2).
  \end{split}
\end{equation*}
Then, for $R$ sufficiently large, we can absorb the last four terms on
the right-hand side to the left-hand side.
As a result, we have that
\begin{equation*}
  \WeightedBEFNorm{p, \ge R}[\psi](\tau_1,\tau_2)
  \lesssim \EnergyFluxWeighted[\psi](\tau_1)
  + \ForcingTermCombinedNorm{p}[\psi, N](\tau_1,\tau_2)
  + R^{p+2}\left(
    \EnergyFlux_{\frac{R}{2}\le r \le R}[\psi](\tau_1)
    + \MorNorm_{\frac{R}{2}\le r\le R}[\psi](\tau_1,\tau_2)
  \right),
\end{equation*}
as desired. This concludes the proof of \zcref[cap]{prop:rp:Kerr}
for $\delta\le p\le 1-\delta$. To prove \zcref[cap]{prop:rp:Kerr} for $1-\delta< p\le 2-\delta$,
it suffices to use the case for $\delta\le p\le 1-\delta$ to bound the
extra $\EnergyFluxFar_{1-\delta,R}[\psi](\tau_2)$ term on the right-hand side of
\zcref[noname]{eq:rp:main-proof:combined-with-nabla-3:2-delta}, and conclude.

\subsection{Proof of \zcref[cap]{prop:rp:Kerr} for \texorpdfstring{$s\ge 1$}{higher-order derivatives}}
\label{sec:rp:higher-order}

In this section, we will prove \zcref[cap]{prop:rp:Kerr} for
$s\ge 1$. To this end, we introduce the following elliptic estimates
proven in \cite{giorgiWaveEquationsEstimates2024} for the horizontal
Hodge operators.
\begin{lemma}[Proposition 9.3.2 in
  \cite{giorgiWaveEquationsEstimates2024}]
  \label{lemma:hodge:elliptic}
  Consider a sphere $S\subset \Manifold$ of the type $S(\tau,r)$. For
  $a$ sufficiently small, the following estimates hold
  \begin{enumerate}
    \item For $f\in \realHorkTensor{p},p=1,2$ we have 
    \begin{equation}
      \label{eq:hodge:elliptic:p-1-2:lower-order}
      \begin{split}
        \int_S\left(r^2\abs*{\nabla f}^2 + \abs*{f}^2\right)
        \lesssim \int_Sr^2\abs*{\HodgeOp{p}f}^2
        + O(a^2)\int_S\abs*{\nabla_{\KillT}f}^2.
      \end{split}    
    \end{equation}
    Moreover, for higher derivatives, in a simplified form, we have
    \begin{equation}
      \label{eq:hodge:elliptic:p-1-2:higher-order}
      \int_S\left(r^2\abs*{\nabla \frakWeightedDeriv^{\le k}f}^2
        + \abs*{\frakWeightedDeriv^{\le k}f}^2\right)
        \lesssim \int_Sr^2\abs*{\frakWeightedDeriv^{\le k}\HodgeOp{p}f}^2
        + O(a^2)\int_S\abs*{\frakWeightedDeriv^{\le k+1}f}^2.
    \end{equation}
  \item For $f\in \realHorkTensor{p},p=0,1$ we have
    \begin{equation}
      \label{eq:hodge:elliptic:p-0-1:lower-order}
      \int_Sr^2\abs*{\nabla f}^2
      \lesssim \int_S\left(r^2\abs*{\HodgeOpDual{p}f}^2 + \abs*{f}^2\right)
      + O(a)^2\int_S\abs*{\nabla_{\KillT}f}^2
      .
    \end{equation}
    Moreover, for higher derivatives, in a simplified form, we have
    \begin{equation}
      \label{eq:hodge:elliptic:p-0-1:higher-order}
      \int_Sr^2\abs*{\nabla\frakWeightedDeriv^{\le k}f}^2
      \lesssim \int_S\left(r^2\abs*{\frakWeightedDeriv^{\le k}\HodgeOpDual{p}f}^2 + \abs*{\frakWeightedDeriv^{\le k}f}^2\right)
      + O(a^2)\int_S\abs*{\frakWeightedDeriv^{\le k+1}f}^2.
    \end{equation}
  \end{enumerate}
\end{lemma}
\begin{proof}
  See the proof of  Proposition 9.3.2 in
  \cite{giorgiWaveEquationsEstimates2024}.
\end{proof}

\paragraph{Step 0: Auxiliary lemma.}
We start by proving the following lemma which will be used to deal
with some terms generated by commutations with $\WaveOpHork{2}$.
\begin{lemma}
  \label{lemma:rp:higher-order-estimates:auxiliary-lemma}
  Let $F$ be a given tensor, and let $\widetilde{N}$ be a tensor with the following schematic structure
  \begin{equation*}
    \widetilde{N} = O(r^{-2})\frakWeightedDeriv F
    + O(r^{-2})\frakWeightedDeriv^{\le1}\psi
    + O(r^{-3})\frakWeightedDeriv^{\le 2}\psi.
  \end{equation*}
  Then, for a given $\widetilde{\psi}$ satisfying
  \begin{equation*}
    \abs*{\widetilde{\psi}}\lesssim \abs*{\frakWeightedDeriv^{\le 1 }\psi}, 
  \end{equation*}
  we have that for $\delta\le p\le 2-\delta$,
  \begin{equation*}
    \begin{split}
      \ForcingTermWeightedNorm{p}{\ge\frac{R}{2}}[\widetilde{\psi}, \widetilde{N}](\tau_1,\tau_2)
      \lesssim{}& \sqrt{\BulkNormWeighted{p,\ge\frac{R}{2}}[\widetilde{\psi}](\tau_1,\tau_2)}
      \left(
      \sqrt{\BulkNormWeighted{p,\ge\frac{R}{2}}[F](\tau_1,\tau_2)}
      + \sqrt{\BulkNormWeighted{p,\ge\frac{R}{2}}[\psi](\tau_1,\tau_2)}
                  \right)\\
      &+ R^{-1}\BulkNormWeighted{p, \ge \frac{R}{2}}[\psi](\tau_1,\tau_2).
    \end{split}
  \end{equation*}
\end{lemma}
\begin{proof}
  In view of the definition of $\ForcingTermWeightedNorm{p}{\frac{R}{2}}$, we have
  that for $\delta\le p \le 2-\delta$,
  \begin{align}
    \label{eq:higher-order-estimates:auxiliary-lemma:aux1}
    \begin{split}         
      \ForcingTermWeightedNorm{p}{\frac{R}{2}}[\widetilde{\psi}, \widetilde{N}](\tau_1,\tau_2)
      \lesssim{}& \int_{\Manifold_{r\ge \frac{R}{2}}(\tau_1,\tau_2)}\abs*{\widetilde{N}}\left(
                  \abs*{\nabla_3\widetilde{\psi}}
                  + r^{p-1}\abs*{\frakWeightedDeriv^{\le 1}\widetilde{\psi}}
                  \right)\\
      \lesssim{}& \sqrt{\BulkNormWeighted{p,\ge\frac{R}{2}}[\widetilde{\psi}](\tau_1,\tau_2)}\left( \int_{\Manifold_{r\ge \frac{R}{2}}(\tau_1,\tau_2)}r^{p+1}\abs*{\widetilde{N}}^2 \right)^{\frac{1}{2}}.
    \end{split}
  \end{align}
  In view of the assumed structure of $\widetilde{N}$, we also have that
  \begin{align}
    \int_{\Manifold_{r\ge \frac{R}{2}}(\tau_1,\tau_2)}r^{p+1}\abs*{\widetilde{N}}^2
    \lesssim{}& \int_{\Manifold_{r\ge \frac{R}{2}}(\tau_1,\tau_2)}r^{p-3}\left(
                \abs*{\frakWeightedDeriv F}^2
                + \abs*{\frakWeightedDeriv^{\le 1}\psi}^2
                + O(R^{-2})\abs*{\frakWeightedDeriv^{\le 2}\psi}^2
                \right)\notag\\
    \lesssim{}& \BulkNormWeighted{p, \ge \frac{R}{2}}[F](\tau_1,\tau_2)
                + \BulkNormWeighted{p, \ge \frac{R}{2}}[\psi](\tau_1,\tau_2)
                + R^{-2}\BulkNormWeighted{p, \ge \frac{R}{2}}^1[\psi](\tau_1,\tau_2).\label{eq:higher-order-estimates:auxiliary-lemma:aux2}
  \end{align}
  Combining \zcref[noname]{eq:higher-order-estimates:auxiliary-lemma:aux1} and
  \zcref[noname]{eq:higher-order-estimates:auxiliary-lemma:aux2} we then have that
  \begin{equation*}
\begin{split}
      \ForcingTermWeightedNorm{p}{\ge\frac{R}{2}}[\widetilde{\psi}, \widetilde{N}](\tau_1,\tau_2)
      \lesssim{}& \sqrt{\BulkNormWeighted{p,\ge\frac{R}{2}}[\widetilde{\psi}](\tau_1,\tau_2)}
      \left(
      \sqrt{\BulkNormWeighted{p,\ge\frac{R}{2}}[F](\tau_1,\tau_2)}
      + \sqrt{\BulkNormWeighted{p,\ge\frac{R}{2}}[\psi](\tau_1,\tau_2)}
                  \right)\\
      &+ R^{-1}\BulkNormWeighted{p, \ge \frac{R}{2}}^1[\psi](\tau_1,\tau_2),
    \end{split}
  \end{equation*}
  as desired.
\end{proof}

\paragraph{Step 1: Commutation with $\HorLieDeriv_{\KillT}$ and $\HorLieDeriv_{\KillPhi}$.}

We first derive an estimate for $\HorLieDeriv_{\KillT}\psi$. To this
end, we first commute \zcref[noname]{eq:model-problem-gRW} with
$\HorLieDeriv_{\KillT}$. Since $\KillT$ and $\KillPhi$ are Killing vectorfields, we have
\begin{align*}
  \WaveOpHork{2}\HorLieDeriv_{\KillT}\psi
  - V\HorLieDeriv_{\KillT}\psi
  ={}& - \frac{4a\cos\theta}{\abs*{q}^2}\LeftDual{\nabla_{\KillT}}\HorLieDeriv_{\KillT}\psi
  + \HorLieDeriv_{\KillT}N,&  \WaveOpHork{2}\HorLieDeriv_{\KillPhi}\psi
  - V\HorLieDeriv_{\KillPhi}\psi
  ={}& - \frac{4a\cos\theta}{\abs*{q}^2}\LeftDual{\nabla_{\KillPhi}}\HorLieDeriv_{\KillPhi}\psi
  + \HorLieDeriv_{\KillPhi}N.
\end{align*}
We can then immediately apply the $r^p$-weighted estimates in \zcref[noname]{eq:rp:Kerr}
with $s=0$ to the commuted equations, which yields
\begin{equation*}
  \begin{split}
    \WeightedBEFNorm{p,\ge R}[(\HorLieDeriv_{\KillT}, \HorLieDeriv_{\KillPhi})\psi](\tau_1,\tau_2)
    \lesssim{}& \EnergyFluxWeighted^1[\psi](\tau_1)
    + R^{p+2}\MorNorm_{\frac{R}{2}\le r \le R }^1[\psi](\tau_1,\tau_2)
    + \ForcingTermWeightedNorm{p}{\frac{R}{2}}^1(\tau_1,\tau_2). 
    \end{split}    
\end{equation*}
Using the relationship between the horizontal Lie derivative and the
horizontal covariant derivative in \zcref[cap]{lemma:LieT-LieZ-to-nabT-nabZ-comparison}, we have that 
\begin{equation}
  \label{eq:rp:higher-order-estimates:commutation-with-T-Z-final}
  \begin{split}
    \WeightedBEFNorm{p,\ge R}[(\HorLieDeriv_{\KillT}, \HorLieDeriv_{\KillPhi})\psi](\tau_1,\tau_2)
    \lesssim{}& \EnergyFluxWeighted^1[\psi](\tau_1)
    + R^{p+2}\MorNorm_{\frac{R}{2}\le r \le R }^1[\psi](\tau_1,\tau_2)
     + \ForcingTermWeightedNorm{p}{\frac{R}{2}}^1(\tau_1,\tau_2).
  \end{split}
\end{equation}

\paragraph{Step 2: Commute with $\abs*{q}\HodgeOp{2}$.}

Next, we commute \zcref[noname]{eq:model-problem-gRW} by
$\abs*{q}\HodgeOp{2}$. Using 
\zcref[noname]{eq:wave-operator-basic-commutation:Hodge-op}, this gives
\begin{equation*}
  \WaveOpHork{1}\left(\abs*{q}\HodgeOp{2}\psi\right)
  - V\abs*{q}\HodgeOp{2}\psi
  = - \frac{4a\cos\theta}{\abs*{q}^2}\LeftDual{\nabla}_{\KillT}\left(\abs*{q}\HodgeOp{2}\psi\right)
  + N_{\abs*{q}\HodgeOp{2}},
\end{equation*}
where, since $\psi$ in fact solves \zcref[noname]{eq:model-problem-gRW}, we have
\begin{align*}
  N_{\abs*{q}\HodgeOp{2}}
  & = q\HodgeOp{2}N
  - \frac{3}{r^2}\abs*{q}\HodgeOp{2}\psi
  - \frac{2a\cos\theta}{\abs*{q}}\LeftDual{\HodgeOp{2}}\HorLieDeriv_{\KillT}\psi
  + O(ar^{-2})\frakWeightedDeriv^{\le 1}\psi
  + O(ar^{-3})\frakWeightedDeriv^{\le 2}\psi
\\& = q\HodgeOp{2}N
  + O(r^{-2})\frakWeightedDeriv\HorLieDeriv_{\KillT}\psi
  + O(ar^{-2})\frakWeightedDeriv^{\le 1}\psi
  + O(ar^{-3})\frakWeightedDeriv^{\le 2}\psi.
\end{align*}
We can then apply \zcref[noname]{eq:rp:Kerr} with $s=0$ again so that
\begin{equation*}
  \begin{split}
    \WeightedBEFNorm{p, \ge R}[\abs*{q}\HodgeOp{2}\psi](\tau_1,\tau_2)
    \lesssim{}&  \EnergyFluxWeighted[\abs*{q}\HodgeOp{2}\psi](\tau_1)
    + R^{p+2}\MorNorm_{\frac{R}{2}\le r \le R }[\abs*{q}\HodgeOp{2}\psi](\tau_1,\tau_2)\\
    & + \ForcingTermWeightedNorm{p}{\frac{R}{2}}[\abs*{q}\HodgeOp{2}\psi, N_{\abs*{q}\HodgeOp{2}}](\tau_1,\tau_2). 
  \end{split}  
\end{equation*}
Next, using \zcref[cap]{lemma:rp:higher-order-estimates:auxiliary-lemma} with
$(\widetilde{\psi}, F) = (\abs*{q}\HodgeOp{2}\psi,
\HorLieDeriv_{\KillT}\psi)$, we can estimate that
\begin{align*}
  \ForcingTermWeightedNorm{p}{\frac{R}{2}}[\abs*{q}\HodgeOp{2}\psi, N_{\abs*{q}\HodgeOp{2}}](\tau_1,\tau_2)
  \lesssim{}& \sqrt{\BulkNormWeighted{p, \frac{R}{2}}[\abs*{q}\HodgeOp{2}\psi](\tau_1,\tau_2)}
              \left(
              \sqrt{\BulkNormWeighted{p, \frac{R}{2}}[\HorLieDeriv_{\KillT}\psi](\tau_1,\tau_2)}
              + \sqrt{\BulkNormWeighted{p,\frac{R}{2}}[\psi](\tau_1,\tau_2)}
              \right)\\
            & + R^{-1}\BulkNormWeighted{p, \frac{R}{2}}[\psi](\tau_1,\tau_2)
              + \ForcingTermWeightedNorm{p}{\frac{R}{2}}^1[\psi, N](\tau_1,\tau_2)
              .
\end{align*}
As a result, we have that
\begin{align*}
  \WeightedBEFNorm{p, \ge R}[\abs*{q}\HodgeOp{2}\psi](\tau_1,\tau_2)
  \lesssim{}&  \EnergyFluxWeighted[\abs*{q}\HodgeOp{2}\psi](\tau_1)
              + R^{p+2}\MorNorm_{\frac{R}{2}\le r \le R }[\abs*{q}\HodgeOp{2}\psi](\tau_1,\tau_2)\\
  & + \sqrt{\BulkNormWeighted{p, \frac{R}{2}}[\abs*{q}\HodgeOp{2}\psi](\tau_1,\tau_2)}
    \left(
    \sqrt{\BulkNormWeighted{p, \frac{R}{2}}[\HorLieDeriv_{\KillT}\psi](\tau_1,\tau_2)}
    + \sqrt{\BulkNormWeighted{p,\frac{R}{2}}[\psi](\tau_1,\tau_2)}
    \right)\\
  & + R^{-1}\BulkNormWeighted{p, \frac{R}{2}}^1[\psi](\tau_1,\tau_2)
    + \ForcingTermWeightedNorm{p}{\frac{R}{2}}^1[\psi, N](\tau_1,\tau_2).
\end{align*}
Then, using \zcref[noname]{eq:rp:Kerr} with $s=0$ and the commuted estimate in
\zcref[noname]{eq:rp:higher-order-estimates:commutation-with-T-Z-final}, we have that
\begin{align*}
  & \WeightedBEFNorm{p, \ge R}[\abs*{q}\HodgeOp{2}\psi](\tau_1,\tau_2)
  + \WeightedBEFNorm{p, \ge R}[\HorLieDeriv_{\KillT}\psi](\tau_1,\tau_2)
    + \WeightedBEFNorm{p, \ge R}[\psi](\tau_1,\tau_2)\\
  \lesssim{}& \EnergyFluxWeighted^1[\psi](\tau_1)
              + R^{-1}\BulkNormWeighted{p, \frac{R}{2}}^1[\psi](\tau_1,\tau_2)
              + \ForcingTermWeightedNorm{p}{\frac{R}{2}}^1(\tau_1,\tau_2)
               + R^{p+2}\left(
              \MorNorm_{\frac{R}{2}\le r \le R }^1[\psi](\tau_1,\tau_2)
              + \EnergyFlux_{\frac{R}{2}\le r\le R}^1[\psi](\tau_1)
              \right)
    ,
\end{align*}
Then using the Hodge estimates from \zcref[cap]{lemma:hodge:elliptic} and
the comparison of $\HorLieDeriv_{\KillT}$ and $\nabla_{\KillT}$ in
\zcref[cap]{lemma:LieT-LieZ-to-nabT-nabZ-comparison}, we get
\begin{equation}
  \label{eq:rp:higher-order-estimates:commutation-with-qD2-final}
  \begin{split}
    \WeightedBEFNorm{p, \ge R}[(\nabla_{\KillT}, \frakWeightedDerivAngular)\psi](\tau_1,\tau_2)
    \lesssim{}& \EnergyFluxWeighted^1[\psi](\tau_1)
                + R^{-1}\BulkNormWeighted{p, \frac{R}{2}}^1[\psi](\tau_1,\tau_2)
                + \ForcingTermWeightedNorm{p}{\frac{R}{2}}^1(\tau_1,\tau_2)\\
              & + R^{p+2}\left(
                \MorNorm_{\frac{R}{2}\le r \le R }^1[\psi](\tau_1,\tau_2)
                + \EnergyFlux^1_{\frac{R}{2}\le r\le R}[\psi](\tau_1)
                \right)
                .
  \end{split}  
\end{equation}

\paragraph{Step 3: Commuting with $r\nabla_4$.}

We now derive an estimate for $r\nabla_4\psi$. To this end, we commute
\zcref[noname]{eq:model-problem-gRW} with $r\nabla_4$, and obtain from \zcref[cap]{lemma:wave-operator-basic-commutation},
\begin{equation*}
  \label{eq:rp:higher-order-estimates:rnab4-commuted-equation}
  \WaveOpHork{2}(r\nabla_4\psi)
  - Vr\nabla_4\psi
  = -\frac{4a\cos\theta}{\abs*{q}^2}\LeftDual{\nabla}_{\KillT}(r\nabla_4\psi)
  + N_{r\nabla_4},
\end{equation*}
where $N_{r\nabla_4}
    = \frac{1}{r}\nabla_4(r\nabla_4\psi) + \widetilde{N}_{r\nabla_4}$ with
\begin{equation*}
  \begin{split}
    \widetilde{N}_{r\nabla_4}={}& O(1)\frakWeightedDeriv^{\le 1}N
                                  + O(r^{-2})\frakWeightedDerivAngular^2\psi
                                  + O(r^{-2})\frakWeightedDeriv^{\le 1}\psi
                                  + O(r^{-3})\frakWeightedDeriv^{\le 2}\psi.
  \end{split}  
\end{equation*}
Next, we apply \zcref[noname]{eq:rp:Kerr} with $s=0$ to the commuted equation
\zcref[noname]{eq:rp:higher-order-estimates:rnab4-commuted-equation}. Since the
$\frac{1}{r}\nabla_4(r\nabla_4\psi)$ term in $N_{r\nabla_4}$ has a
good sign, we can drop it and then we immediately obtain
\begin{equation}
  \label{eq:rp:higher-order-estimates:rnab4:aux1}
  \begin{split}
    \WeightedBEFNorm{p,R}[r\nabla_4\psi](\tau_1,\tau_2)
    \lesssim{}& \EnergyFluxFar_{p, R}^1[\psi](\tau_1)
                + \ForcingTermWeightedNorm{p}{\frac{R}{2}}[r\nabla_4\psi, \widetilde{N}_{r\nabla_4}](\tau_1,\tau_2)\\
              & + R^{p+2}\left(\EnergyFlux_{\frac{R}{2}\le r\le R}^1[\psi](\tau_1) + \MorNorm_{\frac{R}{2}\le r\le R}^1[\psi](\tau_1,\tau_2)\right).
  \end{split}
\end{equation}
Next, applying \zcref[cap]{lemma:rp:higher-order-estimates:auxiliary-lemma} with
$F=\frakWeightedDerivAngular\psi$, we have that
\begin{equation}
  \label{eq:rp:higher-order-estimates:rnab4:aux2}
  \begin{split}
    &\ForcingTermWeightedNorm{p}{\frac{R}{2}}[r\nabla_4,\widetilde{N}_{r\nabla_4}](\tau_1,\tau_2)\\
      \lesssim{}&\sqrt{\BulkNormWeighted{p,\frac{R}{2}}[r\nabla_4\psi](\tau_1,\tau_2)}
                  \left(
                  \sqrt{\BulkNormWeighted{p,\frac{R}{2}}[\frakWeightedDerivAngular\psi](\tau_1,\tau_2)}
                  + \sqrt{\BulkNormWeighted{p,\frac{R}{2}}[\psi](\tau_1,\tau_2)}
                  \right)\\
                 & + R^{-1}\BulkNormWeighted{p,\frac{R}{2}}^1[\psi](\tau_1,\tau_2)
                  + \ForcingTermWeightedNorm{p}{R}^1[\upsilon, N](\tau_1,\tau_2).
  \end{split}
\end{equation}
Plugging \zcref[noname]{eq:rp:higher-order-estimates:rnab4:aux2} into
\zcref[noname]{eq:rp:higher-order-estimates:rnab4:aux1}, we have that
\begin{equation*}
  \begin{split}
    \WeightedBEFNorm{p,R}[r\nabla_4\psi](\tau_1,\tau_2)
    \lesssim{}& \EnergyFluxFar_{p, R}^1[\psi](\tau_1)
                + \ForcingTermWeightedNorm{p}{\frac{R}{2}}^1[\psi, N](\tau_1,\tau_2)
                + R^{-1}\BulkNormWeighted{p,\frac{R}{2}}^1[\psi](\tau_1,\tau_2)\\
              & \sqrt{\BulkNormWeighted{p,\frac{R}{2}}[r\nabla_4\psi](\tau_1,\tau_2)}
                  \left(
                  \sqrt{\BulkNormWeighted{p,\frac{R}{2}}[\frakWeightedDerivAngular\psi](\tau_1,\tau_2)}
                  + \sqrt{\BulkNormWeighted{p,\frac{R}{2}}[\psi](\tau_1,\tau_2)}
                  \right)\\
              & + R^{p+2}\left(\EnergyFlux_{\frac{R}{2}\le r\le R}^1[\psi](\tau_1) + \MorNorm_{\frac{R}{2}\le r\le R}^1[\psi](\tau_1,\tau_2)\right).
  \end{split}
\end{equation*}
Applying
\zcref[noname]{eq:rp:higher-order-estimates:commutation-with-qD2-final} and
\zcref[noname]{eq:rp:higher-order-estimates:commutation-with-T-Z-final}, we then have that
\begin{align*}
  \WeightedBEFNorm{p, R}[(\nabla_{\KillT}, \frakWeightedDerivAngular, r\nabla_4)\psi](\tau_1,\tau_2)
  +\WeightedBEFNorm{p, R}[\psi](\tau_1,\tau_2)
  \lesssim{}& \EnergyFluxFar_{p, R}^1[\psi](\tau_1)
                + \ForcingTermWeightedNorm{p}{\frac{R}{2}}^1[\psi, N](\tau_1,\tau_2)
                + R^{-1}\BulkNormWeighted{p,\frac{R}{2}}^1[\psi](\tau_1,\tau_2)\\
              & + R^{p+2}\left(\EnergyFlux_{\frac{R}{2}\le r\le R}^1[\psi](\tau_1) + \MorNorm_{\frac{R}{2}\le r\le R}^1[\psi](\tau_1,\tau_2)\right).
\end{align*}
Since $(\nabla_{\KillT}, \frakWeightedDerivAngular, r\nabla_4)$ span
$\frakWeightedDeriv$ away from the event horizon, we have that
\begin{align*}
  \WeightedBEFNorm{p, R}^1[\psi](\tau_1,\tau_2)
  \lesssim{}& \EnergyFluxFar_{p, R}^1[\psi](\tau_1)
                + \ForcingTermWeightedNorm{p}{\frac{R}{2}}^1[\psi, N](\tau_1,\tau_2)
                + R^{-1}\BulkNormWeighted{p,\frac{R}{2}}^1[\psi](\tau_1,\tau_2)\\
              & + R^{p+2}\left(\EnergyFlux_{\frac{R}{2}\le r\le R}^1[\psi](\tau_1) + \MorNorm_{\frac{R}{2}\le r\le R}^1[\psi](\tau_1,\tau_2)\right).
\end{align*}
Taking $R$ sufficiently large then allows us to absorb the
$R^{-1}\BulkNormWeighted{p,\frac{R}{2}}^1[\psi](\tau_1,\tau_2)$ term
on the \RHS{} to the \LHS, establishing the desired estimate for $s=1$.


\section{Redshift estimates}
\label{sec:redshift}

In this section, we prove the redshift estimates near the event
horizon. Recall that the redshift effect is also present at the
cosmological horizon for a \KdS{} black hole. However, in the limit of
$\Lambda\to 0$, the redshift effect at the cosmological horizon
vanishes\footnote{This is related to the cosmological horizon turning
  into an extremal horizon in the limit.}. Thus, in this section, we
will only prove the redshift estimates near the event horizon.  For a
discussion and proof of the redshift estimates near the cosmological
horizon for slowly-rotating \KdS, see for instance
\cite{fangLinearStabilitySlowlyRotating2026}. Consistent with the fact
that we are interested in the redshift effect near the event horizon,
the principal global null frame $(e_3,e_4)$ coincides with the
principal ingoing null frame in this section. The main proposition of this section is then as follows.

\begin{proposition}[Redshift estimates]
  \label{prop:redshift-main}
  Let $\psi\in \realHorkTensor{2}$ be a solution to the model problem in
  \zcref[noname]{eq:model-problem-gRW}. Then for $\frac{\abs*{a}}{M}$
  sufficiently small, there exists a sufficiently small constant
  $\delta_{\operatorname{red}}>0$ such that
  $\delta_{\operatorname{red}} =
  \delta_{\operatorname{red}}(M-\abs*{a})$ with
  $\delta_{\operatorname{red}} > \delta_{\Horizon}$ and a small
  constant $c_0>0$ with $c_0=c_0(M-\abs*{a})$ such that the following estimate holds true in $\Manifold(\tau_1,\tau_2)$.
  \begin{equation}
    \label{eq:redshift-main}
    \begin{split}
      &c_0\EnergyFlux_{r\le r_{\EventHorizon}(1+\delta_{\operatorname{red}})}^s[\psi](\tau_2)
    + c_0\MorNorm_{\operatorname{red}}^s[\psi](\tau_1,\tau_2)
    + c_0\SpacelikeFlux_{\mathcal{A}}^s[\psi](\tau_1,\tau_2)\\
    \le{}& \EnergyFlux_{r\le r_{\EventHorizon}(1+2\delta_{\operatorname{red}})}^s[\psi](\tau_1)
      + \delta_{\operatorname{red}}^{-3}\MorNorm_{r_{\EventHorizon}(1+\delta_{\operatorname{red}})\le r\le r_{\EventHorizon}(1+2\delta_{\operatorname{red}})}^s[\psi](\tau_1,\tau_2)
       + \int_{\Manifold(\tau_1,\tau_2)\bigcap\curlyBrace*{\frac{r}{r_{\EventHorizon}}\le 1+2\delta_{\operatorname{red}}}}\abs*{\frakWeightedDeriv^{\le s}N}^2.
    \end{split}    
  \end{equation}
\end{proposition}

\begin{lemma}
  \label{lemma:redshift:Y:event-horizon:bulk-and-boundary}
  Given the vectorfield $Y = \dBar(r)e_3 + d(r)e_4$ and assuming
  \begin{equation*}
    \sup_{r\in \left((1-\delta_{\Horizon})r_{\EventHorizonFuture}, r_0\right)}\left(\abs*{d(r)} + \abs*{\partial_rd(r)}+\abs*{\dBar(r)} + \abs*{\partial_r\dBar(r)}\right) \lesssim 1,
  \end{equation*}
  we have for
  $r\in \left((1-\delta_{\Horizon})r_{\EventHorizonFuture}, r_0\right)$,
  \begin{equation}
    \label{eq:redshift:Y:event-horizon:bulk:EMT-DeformTensor}
    \begin{split}
      2\EMTensor[\psi]\cdot\DeformationTensor[]{Y}
      ={}& \left(
           \partial_r\left(\frac{\Delta}{\abs*{q}^2}\right)\dBar(r)      
           - \frac{\Delta}{\abs*{q}^2}\partial_r\dBar(r)
           \right)\abs*{\nabla_3\psi}^2
           + \partial_rd(r)\abs*{\nabla_4\psi}^2\\
         & + \left(
           \partial_r\dBar(r)
           - \partial_r\left(\frac{\Delta}{\abs*{q}^2}\right)d(r)
           - \frac{\Delta}{\abs*{q}^2}\partial_rd(r)
           \right)\left(\abs*{\nabla\psi}^2 + V\abs*{\psi}^2\right)\\
         & + \frac{2ar}{\abs*{q}^2}\dBar(r) \Re\CCOneFormJ\cdot\nabla\psi\cdot\nabla_3\psi
           - \frac{2r}{\abs*{q}^2}\left(\dBar(r) - \frac{\Delta}{\abs*{q}^2}d(r)\right)\left(
           \nabla_3\psi\cdot\nabla_4\psi
           - V\abs*{\psi}^2
           \right),
    \end{split}    
  \end{equation}
  and
  \begin{equation}
    \label{eq:redshift:Y:event-horizon:bulk:full}
    \KCurrent{Y, 0, 0}[\psi]
    ={} \EMTensor\cdot\DeformationTensor[]{Y}
    + 2\left(\dBar(r) - \frac{\Delta}{\abs*{q}^2}d(r)\right)\partial_r\left(\frac{\Delta}{(r^2+a^2)\abs*{q}^2}\right)\abs*{\psi}^2.
  \end{equation}
\end{lemma}

\begin{proof}
  We can calculate that
  \begin{align*}
    \DeformationTensor[_{\mu\nu}]{Y}
    ={}& \dBar(r)\DeformationTensor[_{\mu\nu}]{e_3}
         + e_{\mu}(r)\partial_r\dBar(r)\Metric(e_3, e_\nu)
         + e_{\nu}(r)\partial_r\dBar(r)\Metric(e_3, e_{\mu})\\
       & + d(r)\DeformationTensor[_{\mu\nu}]{e_4}
         + e_{\mu}(r)\partial_rd(r)\Metric(e_4,e_{\nu})
         + e_{\nu}(r)\partial_rd(r)\Metric(e_4,e_{\mu}).
  \end{align*}
  Using the deformation tensor in the
  principal ingoing frame from
  \zcref[noname]{eq:e3-deformation-tensors:principal-ingoing} and
  \zcref[noname]{eq:e4-deformation-tensors:principal-ingoing}, we get
  \begin{align*}
    \DeformationTensor[_{44}]{Y}
    ={}& 2\partial_r\left(\frac{\Delta}{\abs*{q}^2}\right)\dBar(r)
         - 2\frac{\Delta}{\abs*{q}^2}\partial_r\dBar(r)
         ,&
    \DeformationTensor[_{34}]{Y}
    ={}& \partial_r\dBar(r)
         - \frac{1}{2}\partial_r\left(\frac{\Delta}{\abs*{q}^2}\right)d(r)
         - \frac{\Delta}{\abs*{q}^2}\partial_rd(r),\\
    \DeformationTensor[_{33}]{Y}
    ={}& 2\partial_rd(r),&
    \DeformationTensor[_{4a}]{Y}
    ={}& -\dBar(r)\frac{ar}{\abs*{q}^2}\Re\CCOneFormJ_a,\\
    \DeformationTensor[_{3a}]{Y}
    ={}& 0,&
    \DeformationTensor[_{ab}]{Y}
    ={}& -\frac{r}{\abs*{q}^2}\left(\dBar(r) -  d(r)\frac{\Delta}{\abs*{q}^2}\right)\delta_{ab}.
  \end{align*}
  Recalling the null components of $\EMTensor$ this gives
  \begin{equation}
    \begin{split}
      2\EMTensor[\psi]\cdot\DeformationTensor[]{Y}
      ={}& \left(\dBar(r)\partial_r\left(\frac{\Delta}{\abs*{q}^2}\right) - \frac{\Delta}{\abs*{q}^2}\partial_r\dBar(r)\right)\abs*{\nabla_3\psi}^2
      - \partial_rd(r)\abs*{\nabla_4\psi}^2\\
         &+ \left(
           \partial_r\dBar(r)
           -\partial_r\left(\frac{\Delta}{\abs*{q}^2}\right)d(r)
           - \frac{\Delta}{\abs*{q}^2}\partial_rd(r)
           \right)\left(\abs*{\nabla\psi}^2 + V\abs*{\psi}^2\right)\\
         &+ 2\dBar(r)\left(\frac{ar}{\abs*{q}^2}\Re\CCOneFormJ\right)\cdot\nabla\psi\cdot\nabla_3\psi
          - 2\frac{r}{\abs*{q}^2}\left(\dBar(r)-\frac{\Delta}{\abs*{q}^2}d(r)\right)\left(\nabla_3\psi\cdot\nabla_4\psi - V\abs*{\psi}^2\right).
    \end{split}
  \end{equation}
  Furthermore, we have that $Y(V) = 4\left(\frac{\Delta}{\abs*{q}^2}d(r) - \dBar(r)\right)\partial_r\left(\frac{\Delta}{(r^2+a^2)\abs*{q}^2}\right)$. Thus, we deduce \zcref[noname]{eq:redshift:Y:event-horizon:bulk:full}
  directly from the definition of $\KCurrent{Y, 0, 0}[\psi]$.
\end{proof}

\begin{corollary}
  \label{coro:redshift:estimate:at-event-horizon}
  Let
  \begin{equation}
    \label{eq:redshift:redshift-Y-def}
    \RedShiftY_0\vcentcolon= Y + 2\KillT, \qquad
    Y = \dBar(r)e_3 + d(r) e_4. 
  \end{equation}
  If we choose
  \begin{equation*}
    d(r_{\EventHorizonFuture}) = 0,\qquad
    \dBar(r_{\EventHorizonFuture}) = 1, \qquad
    \partial_rd(r_{\EventHorizonFuture})
    \ge \frac{c_1}{\partial_r\Delta} + \evalAt*{\frac{\partial_r\Delta}{32 r^2}}_{r=r_{\EventHorizonFuture}},\qquad
    \partial_r\dBar(r_{\EventHorizonFuture})
    \ge \frac{c_1}{\partial_r\Delta} + \evalAt*{\frac{\partial_r\Delta}{32 r^2}}_{r=r_{\EventHorizonFuture}}, 
  \end{equation*}
  for some sufficiently large universal constant $c_1$, then at
  $r=r_{\EventHorizonFuture}$, we have that for any sphere
  $S = S(\tau, r_{\EventHorizonFuture})$,
  \begin{equation}
    \label{eq:redshift:estimate:at-event-horizon}
    \begin{split}
      \int_S\CovariantDeriv\cdot\JCurrent{\RedShiftY_{0}, 0, 0}[\psi]
      \ge{}& \frac{\partial_r\Delta(r_{\EventHorizonFuture})}{64 r_{\EventHorizonFuture}^2} \int_S\left(\abs*{\nabla_3\psi}^2 + \abs*{\nabla_4\psi}^2 + \abs*{\nabla \psi}^2 + r^{-2}\abs*{\psi}^2\right)\\
           & + \int_S\RedShiftY_0\psi\cdot\left(\WaveOpHork{2}\psi - V\psi + \frac{4a\cos\theta}{\abs*{q}^2}\LeftDual{\nabla}_{\KillT}\psi\right).
    \end{split}
  \end{equation}
\end{corollary}
\begin{proof}
  Observe that
  \begin{equation*}
    \Delta(r_{\EventHorizonFuture})=0, \qquad
    \evalAt*{\partial_r\frac{\Delta}{\abs*{q}^2}}_{r=r_{\EventHorizonFuture}} = \evalAt*{\frac{\partial_r\Delta}{\abs*{q}^2}}_{r=r_{\EventHorizonFuture}},\qquad
    V(r_{\EventHorizonFuture}) = 2\Lambda.  %
  \end{equation*}
  At $r = r_{\EventHorizonFuture}$ we have the assumptions that
  $d(r_{\EventHorizonFuture})=0$, and
  $\dBar(r_{\EventHorizonFuture}) = 1$. Thus from 
  \zcref{lemma:redshift:Y:event-horizon:bulk-and-boundary}, we have that
  \begin{align*}
    2\KCurrent{\RedShiftY_0, 0, 0}[\psi]
    = {}& \partial_rd(r_{\EventHorizonFuture})\abs*{\nabla_4\psi}^2
          + \evalAt*{\frac{\partial_r\Delta}{\abs*{q}^2}}_{r=r_{\EventHorizonFuture}}\abs*{\nabla_3\psi}^2
          + \partial_r\dBar(r_{\EventHorizonFuture}) \left(\abs*{\nabla\psi}^2 + V\abs*{\psi}^2\right)\\
        & + 2\evalAt*{\frac{ar}{\abs*{q}^2}\Re\CCOneFormJ}_{r=r_{\EventHorizonFuture}}\cdot\nabla\psi\cdot\nabla_3\psi
          -  \evalAt*{\frac{2r}{\abs*{q}^2}}_{r=r_{\EventHorizonFuture}}\left(\nabla_3\psi\cdot\nabla_4\psi - V\abs*{\psi}^2\right)\\
        & + 4\evalAt*{\partial_r\left(\frac{\Delta}{\left(r^2+a^2\right)\abs*{q^2}}\right)}_{r=r_{\EventHorizonFuture}}
          \abs*{\psi}^2.
  \end{align*}
  We can then make the following computation exactly at the event horizon
  \begin{align*}
    \evalAt*{\KCurrent{\RedShiftY_0, 0, 0}[\psi]}_{r=r_{\EventHorizonFuture}}
    \ge{}& \evalAt*{\frac{\partial_r\Delta}{8\abs*{q}^2}}_{r=r_{\EventHorizonFuture}}\abs*{\nabla_3\psi}^2
           + \left(\frac{\partial_rd(r_{\EventHorizonFuture})}{2} - \evalAt*{\frac{2r^2}{\partial_r\Delta\abs*{q}^2}}_{r=r_{\EventHorizonFuture}}\right)\abs*{\nabla_4\psi}^2\\           
         & + \left(
           \frac{\partial_r\dBar(r_{\EventHorizonFuture})}{2}
           - \evalAt*{\frac{2a^2r^2}{\abs*{q}^2\partial_r\Delta}\abs*{\Re\CCOneFormJ}^2}_{r=r_{\EventHorizonFuture}}
           \right)\abs*{\nabla\psi}^2
           + 
           2\evalAt*{\frac{\partial_r\Delta}{(r^2+a^2)\abs*{q}^2}}_{r=r_{\EventHorizonFuture}}    
           \abs*{\psi}^2\\
         & + \left(\frac{\partial_r\dBar(r_{\EventHorizonFuture})}{2} + \evalAt*{\frac{r}{\abs*{q}^2}}_{r=r_{\EventHorizonFuture}}\right)2\Lambda\abs*{\psi}^2.%
  \end{align*}

  Now, furthermore, observe that
  \begin{align*}
    \RedShiftY_0^{\mu}\HorCovDeriv^{\nu}\psi^a\HorRiem_{ab\nu\mu}\psi^b
    ={}& Y^{\mu}\HorCovDeriv^{\nu}\psi^a\HorRiem_{ab\nu\mu}\psi^b
         + 2\KillT^{\mu}\HorCovDeriv^{\nu}\psi^a\HorRiem_{ab\nu\mu}\psi^b\\
    ={}& -\LeftDual{\rho}\in_{AB}\left(d(r)\nabla_4\psi^A\psi^B - \dBar(r)\nabla_3\psi^A\psi^B\right)
         - \LeftDual{\rho}\in_{AB}\left(\frac{\Delta}{\abs*{q}^2}\nabla_4\psi^A\psi^B - \nabla_3\psi^A\psi^B\right)\\
       & - 2a\LeftDual{\Re \CCOneFormJ}^d\nabla_d\psi_A\rho\LeftDual{\psi}^A.
  \end{align*}
  As a result, at $r=r_{\EventHorizonFuture}$, we have that
  \begin{align*}
    \RedShiftY_0^{\mu}\HorCovDeriv^{\nu}\psi^a\HorRiem_{ab\nu\mu}\psi^b
    ={}& 2\LeftDual{\rho}\in_{AB}\nabla_3\psi^A\psi^B
         - 2a\LeftDual{\Re \CCOneFormJ}^d\nabla_d\psi_A\rho\LeftDual{\psi}^A.
  \end{align*}
  As a result, we can write that
  \begin{align*}
    \evalAt*{\CovariantDeriv\cdot\JCurrent{\RedShiftY_0, 0, 0}[\psi]}_{r=r_{\EventHorizonFuture}}
    \ge{}&  \evalAt*{\frac{\partial_r\Delta}{8\abs*{q}^2}}_{r=r_{\EventHorizonFuture}}\abs*{\nabla_3\psi}^2
           + \left(\frac{\partial_rd(r_{\EventHorizonFuture})}{2} - \evalAt*{\frac{2r^2}{\partial_r\Delta\abs*{q}^2}}_{r=r_{\EventHorizonFuture}}\right)\abs*{\nabla_4\psi}^2\\           
         & + \left(
           \frac{\partial_r\dBar(r_{\EventHorizonFuture})}{2}
           + \evalAt*{\frac{2a^2r^2}{\abs*{q}^2\partial_r\Delta}\abs*{\Re\CCOneFormJ}_{r=r_{\EventHorizonFuture}}}^2
           - a\abs*{\Re\CCOneFormJ}\abs*{\rho}r_{\EventHorizonFuture}
           \right)\abs*{\nabla\psi}^2\\
         & +
           \left(
           \evalAt*{2\frac{\partial_r\Delta}{(r^2+a^2)\abs*{q}^2}}_{r=r_{\EventHorizonFuture}}
           -\evalAt*{\frac{8\abs*{q}^2\LeftDual{\rho}^2}{\partial_r\Delta}}_{r=r_{\EventHorizonFuture}}
           -\evalAt*{\frac{a\abs*{\Re\CCOneFormJ}\abs*{\rho}}{r_{\EventHorizonFuture}}}_{r=r_{\EventHorizonFuture}}
           \right)
           \abs*{\psi}^2\\
         & + \left(\frac{\partial_r\dBar(r_{\EventHorizonFuture})}{2} + \evalAt*{\frac{r}{\abs*{q}^2}}_{r=r_{\EventHorizonFuture}}\right)2\Lambda\abs*{\psi}^2%
           +\RedShiftY_0\psi\cdot\left(\WaveOpHork{2}\psi-V\psi\right)
           .
  \end{align*}
  Observe also that $-\frac{4a\cos\theta}{\abs*{q}^2}\RedShiftY_0\psi\cdot\LeftDual{\nabla}_{\KillT}\psi=-\frac{4a\cos\theta}{\abs*{q}^2}\nabla_{Y}\psi\cdot\LeftDual{\nabla}_{\KillT}\psi$. Thus, at $r=r_{\EventHorizonFuture}$, we have that
  \begin{align*}
    -\frac{4a\cos\theta}{\abs*{q}^2}\RedShiftY_0\psi\cdot\LeftDual{\nabla}_{\KillT}\psi
    ={}& -\frac{2a\cos\theta}{\abs*{q}^2(1+\gamma)}\nabla_3\psi\cdot\nabla_{e_4-2a\Re\CCOneFormJ^b e_b}\LeftDual{\psi}\\
    \ge{}& -\frac{2a}{\abs*{q}^2}\abs*{\nabla_3\psi}\abs*{\nabla_4\psi}
           -\frac{4a^2}{\abs*{q}^2} \abs*{\Re\CCOneFormJ}\abs*{\nabla \psi}\abs*{\nabla_3\psi}\\
    \ge{}& - \frac{\partial_r\Delta}{16\abs*{q}^2}\abs*{\nabla_3\psi}^2
           - \frac{32a^2}{\partial_r\Delta\abs*{q}^2}\abs*{\nabla_4\psi}^2
           - \frac{128a^4}{\partial_r\Delta\abs*{q}^2}\abs*{\Re\CCOneFormJ}^2\abs*{\nabla\psi}^2.
  \end{align*}
  As a result, we have that
  \begin{align*}
    \evalAt*{\CovariantDeriv\cdot\JCurrent{\RedShiftY_0, 0, 0}[\psi]}_{r=r_{\EventHorizonFuture}}
    \ge{}& \evalAt*{\frac{\partial_r\Delta}{16\abs*{q}^2}}_{r=r_{\CosmologicalHorizonFuture}}\abs*{\nabla_3\psi}^2
           + \left(
           \frac{\partial_rd(r_{\EventHorizonFuture})}{2}
           - \evalAt*{\frac{r^2}{\partial_r\Delta\abs*{q}^2}}_{r=r_{\EventHorizonFuture}}
           - \evalAt*{\frac{32a^2}{\partial_r\Delta\abs*{q}^2}}_{r=r_{\EventHorizonFuture}}
           \right)\abs*{\nabla_4\psi}^2\\           
         &  + \left(
           \frac{\partial_r\dBar(r_{\EventHorizonFuture})}{2} %
           - \frac{2a^2r^2}{\abs*{q}^2\partial_r\Delta}\abs*{\Re\CCOneFormJ}^2
           - a\evalAt*{\abs*{\Re \CCOneFormJ}\abs*{\rho}r}_{r=r_{\EventHorizonFuture}}
           - \frac{128a^4}{\partial_r\Delta\abs*{q}^2}\abs*{\Re\CCOneFormJ}^2
           \right)\abs*{\nabla\psi}^2\\
         & + \left(
           \evalAt*{\frac{2\partial_r\Delta}{(r^2+a^2)\abs*{q}^2}}_{r=r_{\EventHorizonFuture}}
           - \evalAt*{\frac{8\abs*{q}^2\LeftDual{\rho}^2}{\partial_r\Delta}}_{r=r_{\EventHorizonFuture}}
           - \evalAt*{\frac{a\abs*{\Re\CCOneFormJ}\abs*{\rho}}{r}}_{r=r_{\EventHorizonFuture}}
           \right)\abs*{\psi}^2\\
         & + \left(\partial_r\dBar(r_{\EventHorizonFuture}) + \evalAt*{\frac{r}{\abs*{q}^2}}_{r=r_{\EventHorizonFuture}}\right)2\Lambda\abs*{\psi}^2%
           +\RedShiftY_0\psi\cdot\left(\WaveOpHork{2}\psi-V\psi + \frac{4a\cos\theta}{\abs*{q}^2}\LeftDual{\nabla}_{\KillT}\psi\right)
           .
  \end{align*}
  Then, recall that
  \begin{equation*}
    r\le \abs*{q}\le 2r,\qquad
    \abs*{\rho} \le \frac{2 M}{r^{-3}} + O(a r^{-4})
    \qquad
    \abs*{\LeftDual{\rho}} \le  \frac{6aM}{r^4} + O(ar^{-5}),\qquad
    \abs*{\CCOneFormJ} \le \frac{\sqrt{\kappa}}{r} + O(ar^{-2}).
  \end{equation*}
  We can then infer that at $r=r_{\EventHorizonFuture}$, for
  $\partial_r\dBar(r_{\EventHorizonFuture})\ge 0$,
  \begin{align*}
    \evalAt*{\CovariantDeriv\cdot\JCurrent{\RedShiftY_0, 0, 0}[\psi]}_{r=r_{\EventHorizonFuture}}
    \ge{}& \evalAt*{\frac{\partial_r\Delta}{32r^2}}_{r=r_{\EventHorizonFuture}}\abs*{\nabla_3\psi}^2
           + \left(
           \frac{\partial_rd(r_{\EventHorizonFuture})}{2}
           - \evalAt*{\frac{33}{\partial_r\Delta}}_{r=r_{\CosmologicalHorizonFuture}}
           \right)\abs*{\nabla_4\psi}^2\\
         &  + \left(
           \frac{\partial_r\dBar(r_{\CosmologicalHorizonFuture})}{2}
           - \evalAt*{\frac{131Ma}{r^2\partial_r\Delta}}_{r=r_{\CosmologicalHorizonFuture}}
           \right)\abs*{\nabla\psi}^2
           + \left(
           \evalAt*{\frac{\partial_r\Delta}{8r^4}}_{r=r_{\CosmologicalHorizonFuture}}
           + \frac{1154 a M}{r^4\partial_r\Delta}
           \right)\abs*{\psi}^2\\
         & + \RedShiftY_0\psi\cdot\left(\WaveOpHork{2}\psi - V\psi + \frac{4a\cos\theta}{\abs*{q}^2}\LeftDual{\nabla}_{\KillT}\psi\right).
  \end{align*}
  Then, recall from the Poincar\'{e} inequality in \zcref[cap]{lemma:Poincare}, we
  have that for $a$ sufficiently small, there exists some universal
  constant $c_0>0$
  \begin{equation*}
    \int_S\left(
      \abs*{\nabla_4\psi}^2
      + \frac{\Delta}{r^2}\abs*{\nabla_3\psi}^2
      + \abs*{\nabla\psi}^2
    \right)
    \ge \frac{c_0}{r^2}\int_S\abs*{\psi}^2. 
  \end{equation*}
  In particular, at $r=r_{\EventHorizonFuture}$, we have that
  \begin{equation*}
    \int_S\left(
      \abs*{\nabla_4\psi}^2
      + \abs*{\nabla\psi}^2
    \right)
    \ge \frac{c_0}{r^2}\int_S\abs*{\psi}^2. 
  \end{equation*}
  As a result, for any sphere $S = S(\tau, r_{\EventHorizonFuture})$,
  \begin{equation}
    \label{eq:redshift:main-bulk:before-simplification}
    \begin{split}
      \CovariantDeriv\cdot\JCurrent{\RedShiftY_0,0,0}[\psi]
      \ge{}& \evalAt*{\frac{\partial_r\Delta}{32r^2}}_{r=r_{\EventHorizonFuture}}\abs*{\nabla_3\psi}^2
             + \left(
             \frac{\partial_rd(r_{\EventHorizonFuture})}{2}
             - \evalAt*{\frac{33}{\partial_r\Delta}}_{r=r_{\EventHorizonFuture}}
             -\evalAt*{\frac{1}{c_0}\frac{1154 a M}{r^2\partial_r\Delta}}_{r=r_{\EventHorizonFuture}}
             \right)\abs*{\nabla_4\psi}^2\\
           &  + \left(
             \frac{\partial_r\dBar(r_{\EventHorizonFuture})}{2}
             - \evalAt*{\frac{131Ma}{r^2\partial_r\Delta}}_{r=r_{\EventHorizonFuture}}
             -\evalAt*{\frac{1}{c_0}\frac{1154 a M}{r^2\partial_r\Delta}}_{r=r_{\EventHorizonFuture}}
             \right)\abs*{\nabla\psi}^2\\
           & + \RedShiftY_0\psi\cdot\left(\WaveOpHork{2}\psi - V\psi + \frac{4a\cos\theta}{\abs*{q}^2}\LeftDual{\nabla}_{\KillT}\psi\right).
    \end{split}    
  \end{equation}
  Then, there exists some sufficiently large $c_1$ such that if
  $d(r), \dBar(r)$ is chosen so that
  \begin{equation*}
    \partial_rd(r_{\EventHorizonFuture})
    \ge \frac{c_1}{\partial_r\Delta} + \evalAt*{\frac{\partial_r\Delta}{32 r^2}}_{r=r_{\EventHorizonFuture}},\qquad
    \partial_r\dBar(r_{\EventHorizonFuture})
    \ge \frac{c_1}{\partial_r\Delta} + \evalAt*{\frac{\partial_r\Delta}{32 r^2}}_{r=r_{\EventHorizonFuture}}, 
  \end{equation*}
  then,
  \begin{align*}
    \int_S\CovariantDeriv\cdot\JCurrent{\RedShiftY_0,0,0}[\psi]
    \ge{}& \evalAt*{\frac{\partial_r\Delta}{64 r^2}}_{r=r_{\EventHorizonFuture}}
           \int_S\left(\abs*{\nabla_3\psi}^2 + \abs*{\nabla_4\psi}^2 + \abs*{\nabla\psi}^2 + r^{-2}\abs*{\psi}^2\right)\\
         &+ \int_S\RedShiftY_0\psi\cdot\left(\WaveOpHork{2}\psi - V\psi + \frac{4a\cos\theta}{\abs*{q}^2}\LeftDual{\nabla}_{\KillT}\psi\right),
  \end{align*}
  as stated. This concludes the proof of 
  \zcref[cap]{coro:redshift:estimate:at-event-horizon}. 
\end{proof}

\begin{proposition}
  \label{prop:redshift:main-props}
  Let $\dot{\chi}(r)$ be a positive bump function such that $\supp \dot{\chi}\subset [-2,2]$ and $\dot{\chi}= 1$ on $[-1,1]$. Also, for $a<M$, let $\delta_{red}>0$ be a sufficiently small
  constant such that $\delta_{red} = \delta_{red}(M-a)$ with
  $\delta_{red}>\delta_{\Horizon}$. Furthermore, let
  $\RedShiftY_{\EventHorizon}\vcentcolon= \chi_{\EventHorizon}\RedShiftY_0$ with $\chi_{\EventHorizon}= \dot{\chi}\left(\frac{\frac{r}{r_{\EventHorizon}} - 1}{\delta_{red}}\right)$ and where $\RedShiftY_0$ is the vectorfield of 
  \zcref[cap]{coro:redshift:estimate:at-event-horizon} defined in
  \zcref[noname]{eq:redshift:redshift-Y-def}. Then, the following estimate
  holds for any sphere $S=S(\tau,r)$
  \begin{equation}
    \label{eq:redshift:main-props:bulk}
    \begin{split}
      \int_S\CovariantDeriv\cdot\JCurrent{\RedShiftY_{\EventHorizon},0,0}[\psi]
    \ge{}& -\evalAt*{\frac{\partial_r\Delta}{128 r^2}}_{r=r_{\EventHorizon}}
    \int_S\left(\abs*{\nabla_3\psi}^2 + \abs*{\nabla_4\psi}^2 + \abs*{\nabla\psi}^2 + r^{-2}\abs*{\psi}^2\right)1_{\abs*{\frac{r}{r_{\EventHorizon}}-1}\le \delta_{red}}\\
    &- O\left(r_{\EventHorizon}^{-1}\delta_{red}^{-1}\right)
    \int_S\left(\abs*{\nabla_3\psi}^2 + \abs*{\nabla_4\psi}^2 + \abs*{\nabla\psi}^2 + r^{-2}\abs*{\psi}^2\right)1_{\delta_{red}\le \abs*{\frac{r}{r_{\EventHorizon}}-1}\le 2\delta_{red}}\\
   & + \int_S\RedShiftY_{\EventHorizon}\psi\cdot\left(\WaveOpHork{2}\psi- V\psi + \frac{4a\cos\theta}{\abs*{q}^2}\LeftDual{\nabla}_{\KillT}\psi\right).
    \end{split}
  \end{equation}
  In addition, $\JCurrent{\RedShiftY_{\EventHorizon},0, 0}[\psi]\cdot N_{\Sigma} \ge 0$ where $N_{\Sigma}$ is the  normal to $\Sigma$ given by
  $N_{\Sigma}= - \Metric^{\alpha\beta}\partial_{\alpha}\tau\partial_{\beta}$. In
  addition, there exists a constant $c_0=c_0(M-a)$ such that for any
  sphere $S=S(\tau,r)$ with
  $\abs*{\frac{r}{r_{\EventHorizon}}-1}\le \delta_{red}$,
  \begin{equation}
    \label{eq:redshift:main-props:boundary:Sigma}
    \int_S\JCurrent{\RedShiftY_{\EventHorizon},0,0}[\psi]\cdot N_{\Sigma}
    \ge c_0\int_S\left(
      \abs*{\nabla_3\psi}^2
      + \abs*{\nabla_4\psi}^2
      + \abs*{\nabla\psi}^2
      + r^{-2}\abs*{\psi}^2
    \right).
  \end{equation}
  Furthermore, there exists a constant $c_0=c_0(M-a)$ such that for
  any sphere $S=S(\tau,r_{*})$,
  \begin{equation}
    \label{eq:redshift:main-props:boundary:A}
    \int_S\JCurrent{\RedShiftY_{\EventHorizon},0,0}[\psi]\cdot N_{\mathcal{A}}
    \ge c_0\int_S\left(
      \delta_{\Horizon}\abs*{\nabla_3\psi}^2
      + \abs*{\nabla_4\psi}^2
      + \abs*{\nabla\psi}^2
      + r^{-2}\abs*{\psi}^2
    \right).
  \end{equation}
\end{proposition}
\begin{proof}
  We first prove \zcref[noname]{eq:redshift:main-props:bulk}. In view of the
  properties of $\RedShiftY_{\EventHorizon}$ and the properties of
  $\chi_{\EventHorizon}$, we have that
  \begin{align}
    \CovariantDeriv\cdot\JCurrent{\RedShiftY_{\EventHorizon}, 0, 0}[\psi]
    ={}& \chi_{\EventHorizon}\CovariantDeriv\cdot\JCurrent{\RedShiftY_0, 0, 0}[\psi]
         + \EMTensor\left(\RedShiftY_0,d\chi_{\EventHorizon}\right)\notag\\
    ={}&\chi_{\EventHorizon}\CovariantDeriv\cdot\JCurrent{\RedShiftY_0, 0, 0}[\psi]
         - O\left(r_{\EventHorizon}^{-1}\delta_{\operatorname{red}}^{-1}\right)\bOne_{\delta_{\operatorname{red}}\le \abs*{\frac{r}{r_{\EventHorizon}}-1}\le 2\delta_{\operatorname{red}}}
         \left(
         \abs*{\nabla_3\psi}^2
         + \abs*{\nabla_4\psi}^2
         + \abs*{\nabla\psi}^2
         + r^{-2}\abs*{\psi}^2
         \right). \label{eq:redshift:Div-J-Yh:cutoff-expansion}
  \end{align}
  Then, using \zcref[cap]{coro:redshift:estimate:at-event-horizon},
  we have that for $\delta_{\operatorname{red}}$ chosen sufficiently
  small, for any sphere $S=S(\tau,r)$,
  \begin{equation*}
    \begin{split}
      \int_S\CovariantDeriv\cdot\JCurrent{\RedShiftY_{0}, 0, 0}[\psi]
      \ge{}& \frac{r_{\EventHorizonFuture}-M}{256r_{\EventHorizonFuture}^2}
             \int_S\left(\abs*{\nabla_3\psi}^2 + \abs*{\nabla_4\psi}^2 + \abs*{\nabla \psi}^2 + r^{-2}\abs*{\psi}^2\right)
             \bOne_{\abs*{\frac{r}{r_{\EventHorizon}-1}}\le \delta_{\operatorname{red}}}\\
           & - O\left(r_{\EventHorizon}^{-1}\delta_{\operatorname{red}}^{-1}\right)
             \int_{S}\bOne_{\delta_{\operatorname{red}}\le \abs*{\frac{r}{r_{\EventHorizon}}-1}\le 2\delta_{\operatorname{red}}}
             \left(
             \abs*{\nabla_3\psi}^2
             + \abs*{\nabla_4\psi}^2
             + \abs*{\nabla\psi}^2
             + r^{-2}\abs*{\psi}^2
             \right)
      \\
           & + \int_S\RedShiftY_0\psi\cdot\left(\WaveOpHork{2}\psi - V\psi + \frac{4a\cos\theta}{\abs*{q}^2}\LeftDual{\nabla}_{\KillT}\psi\right),
    \end{split}
  \end{equation*}
  as stated. We now turn to proving \zcref[noname]{eq:redshift:main-props:boundary:Sigma}.
  We can compute that
  \begin{align*}
    \JCurrent{\RedShiftY_{\EventHorizon},0,0}[\psi]\cdot N_{\Sigma}
    ={}&\EMTensor\left(\RedShiftY_{\EventHorizon}, N_{\Sigma}\right)\\
    ={}& \chi_{\EventHorizon}\EMTensor\left(\RedShiftY_0, N_{\Sigma}\right)\\
    ={}& \chi_{\EventHorizon}\EMTensor\left(
         \dBar(r)e_3
         + d(r)e_4
         + \frac{1}{1+\gamma}e_4
         + \frac{\Delta}{(1+\gamma)\abs*{q}^2}e_3
         - 2a\Re\CCOneFormJ^be_b,
         N_{\Sigma}
         \right)\\
    ={}& \chi_{\EventHorizon}\EMTensor\left(
         e_3
         + e_4
         + O(r-r_{\EventHorizon})e_4
         + O(r-r_{\EventHorizon})e_3
         - 2a\Re\CCOneFormJ^be_b,
         N_{\Sigma}
         \right).
  \end{align*}
  Using the support of $\chi_{\EventHorizon}$, we then have that
  \begin{align*}
    \JCurrent{\RedShiftY_{\EventHorizon},0,0}[\psi]\cdot N_{\Sigma}
    ={}& \chi_{\EventHorizon}\left(
         \EMTensor\left(
         e_3+e_4-2a\Re\CCOneFormJ^be_b,N_{\Sigma}
         \right)
         + O\left(\delta_{\operatorname{red}}\right)\left(
         \abs*{\nabla_3\psi}^2
         + \abs*{\nabla_4\psi}^2
         + \abs*{\nabla\psi}^2
         + r^{-2}\abs*{\psi}^2
         \right)
         \right).
  \end{align*}
  Then, we recall that we chose $\tau$ in \zcref[cap]{definition de k et de h} such that $\Metric(N_{\Sigma},N_{\Sigma}) \le \frac{M^2}{-8r^2} < 0$. Moreover, we have that
  \begin{align*}
    \Metric\left(e_3+e_4+2a\Re\CCOneFormJ^be_b,e_3+e_4+2a\Re\CCOneFormJ^be_b\right)
    ={}& -4\left(1-\frac{a^2\kappa\sin^2\theta}{\abs*{q}^2}\right)
    \le{} -4\left(1-\frac{a^2}{r_{\EventHorizon}^2} + O(\delta_{\operatorname{red}})\right) <0
  \end{align*}
  for $\abs*{a}<M$, and $\delta_{\operatorname{red}}$ sufficiently
  small.  Thus we have shown that both $N_{\Sigma}$ and
  $e_3+e_4+2a\Re\CCOneFormJ^be_b$ are uniformly timelike on the support
  of $\chi_{\EventHorizon}$. Next, recall from the explicit form of $V$ in
  \zcref[noname]{eq:model-problem-gRW} that
  $V(r_{\EventHorizon})=2\Lambda$. As a result, we have that on the
  support of $\chi_{\EventHorizon}$,
  \begin{gather*}
    \EMTensor_{33}[\psi] = \abs*{\nabla_3\psi}^2,\qquad
    \EMTensor_{44}[\psi] = \abs*{\nabla_4\psi}^2,\qquad
    \EMTensor_{34}[\psi] = \abs*{\nabla\psi}^2 + \left(2\Lambda+O(\delta_{\operatorname{red}})\right)\abs*{\psi}^2,\\
    \EMTensor_{4a}[\psi] = \nabla_4\psi\cdot\nabla_a\psi,\qquad
    \EMTensor_{3a}[\psi] = \nabla_3\psi\cdot\nabla_a\psi,\\
    \Trace \EMTensor[\psi] = \nabla_4\psi\cdot\nabla_3\psi + \left(-2\Lambda + O(\delta_{\operatorname{red}})\right)\abs*{\psi}^2,\qquad
    \widehat{\EMTensor}_{ab}[\psi] = \frac{1}{2}\left(\nabla\psi\SymTracelessTensorProd\nabla\psi\right)_{ab}.
  \end{gather*}
  As a result, since both $N_{\Sigma}$ and $e_3+e_4+2a\Re\CCOneFormJ$
  are timelike, we have that
  \begin{align*}
    \EMTensor[\psi]\left(e_3+e_4-2a\Re\CCOneFormJ^be_b,N_{\Sigma}\right)
    \ge{}& \EMTensor_0[\psi]\left(e_3+e_4-2a\Re\CCOneFormJ^be_b,N_{\Sigma}\right)\\
    & + O\left(\delta_{\operatorname{red}}\right)\left(
      \abs*{\nabla_3\psi}^2
      + \abs*{\nabla_4\psi}^2
      + \abs*{\nabla\psi}^2
      + r^{-2}\abs*{\psi}^2,
    \right)
  \end{align*}
  uniformly in $\Lambda$, where $\EMTensor_0[\psi]$ is the symmetric 2-tensor defined by $\EMTensor_0[\psi]\vcentcolon= \HorCovDeriv_{\mu}\psi\cdot\HorCovDeriv_{\nu}\psi
    - \frac{1}{2}\Metric_{\mu\nu}\HorCovDeriv^{\lambda}\psi\cdot\HorCovDeriv_{\lambda}\psi$. As a result, on the support of $\chi_{\EventHorizon}$,
  \begin{equation*}
    \JCurrent{\RedShiftY_{\EventHorizon},0,0}[\psi]\cdot N_{\Sigma}
    \ge \chi_{\EventHorizon}\left(
      \EMTensor_0[\psi]\left(e_3+e_4-2a\Re\CCOneFormJ^be_b,N_{\Sigma}\right)
      + O\left(\delta_{\operatorname{red}}\right)\left(
        \abs*{\nabla_3\psi}^2
        + \abs*{\nabla_4\psi}^2
        + \abs*{\nabla\psi}^2
        + r^{-2}\abs*{\psi}^2
      \right)
    \right).
  \end{equation*}
  Combined with the
  Poincar\'{e} inequality of \zcref[cap]{lemma:Poincare}, we then have
  that there exists some constant $c_0>0$ such that for
  $\delta_{\operatorname{red}}$ sufficiently small, for any sphere
  $S=S(\tau,r)$,
  \begin{equation*}
    \int_S\JCurrent{\RedShiftY_{\EventHorizon},0,0}[\psi]\cdot N_{\Sigma}
    \ge c_0\chi_{\EventHorizon}\int_S\left(
      \abs*{\nabla_3\psi}^2
      + \abs*{\nabla_4\psi}^2
      + \abs*{\nabla\psi}^2
      + r^{-2}\abs*{\psi}^2
    \right). 
  \end{equation*}
  In particular, since $\chi_{\EventHorizon}$ is a non-negative bump
  function, we have that on $\Manifold$, $\int_S\JCurrent{\RedShiftY_{\EventHorizon},0,0}[\psi]\cdot N_{\Sigma}\ge 0$. Moreover, given the definition of $\chi_{\EventHorizon}$, we have that
  for any sphere $S(\tau,r)$ with
  $\abs*{\frac{r}{r_{\EventHorizon}}-1}\le \delta_{\operatorname{red}}$
  \begin{equation*}
    \int_S\JCurrent{\RedShiftY_{\EventHorizon},0,0}[\psi]\cdot N_{\Sigma}
    \ge c_0\int_S\left(
      \abs*{\nabla_3\psi}^2
      + \abs*{\nabla_4\psi}^2
      + \abs*{\nabla\psi}^2
      + r^{-2}\abs*{\psi}^2
    \right), 
  \end{equation*}
  as stated. Finally, we prove
  \zcref[noname]{eq:redshift:main-props:boundary:A}. To this end, observe that
  on $\mathcal{A} = \curlyBrace*{r = (1-\delta_{\Horizon})r_{\EventHorizon}}$ we have $N_{\mathcal{A}}=\frac{1}{2}e_4%
         - \frac{1}{2}\left(\frac{\Delta}{(1+\gamma)\abs*{q}^2}%
         \right)e_3$ and also
  \begin{align*}
    \RedShiftY_0
    ={}& \left(1+ O(\delta_{\Horizon})\right)e_3
         + \left(1+ O(\delta_{\Horizon})\right)e_4
         - 2a\Re\CCOneFormJ^be_b.
  \end{align*}
  Recall that we defined
  $\delta_{\operatorname{red}}\ge \delta_{\Horizon}$ so that
  $\chi_{\EventHorizon}=1$ on $\mathcal{A}$. Thus,
  \begin{align*}
    \JCurrent{\RedShiftY_{\EventHorizon},0,0}[\psi]\cdot N_{\mathcal{A}}
    ={}& \EMTensor[\psi](\RedShiftY_{\EventHorizon}, N_{\mathcal{A}})\\
    ={}& \EMTensor[\psi](\RedShiftY_0, N_{\mathcal{A}})\\
    ={}& \frac{1}{2}\left(1+O(\delta_{\Horizon})\right)\EMTensor_{43}[\psi]
         + \frac{1}{2}\left(1+O(\delta_{\Horizon})\right)\EMTensor_{44}[\psi]
         - a\Re\CCOneFormJ^b\EMTensor_{4b}\\
       &+ \frac{1}{2}\left(\frac{\abs*{\Delta}}{(1+\gamma)\abs*{q}^2}\right)\left(1+O(\delta_{\Horizon})\EMTensor_{33}[\psi] - 2a\Re\CCOneFormJ^b\EMTensor_{3b}[\psi]\right)
         .
  \end{align*}
  Next, recall from the explicit form of $V$ in
  \zcref[noname]{eq:model-problem-gRW} that
  $V(r_{\EventHorizon})=2\Lambda$. As a result, we have that on
  $\mathcal{A}$,
  \begin{gather*}
    \EMTensor_{33}[\psi] = \abs*{\nabla_3\psi}^2,\qquad
    \EMTensor_{44}[\psi] = \abs*{\nabla_4\psi}^2,\qquad
    \EMTensor_{34}[\psi] = \abs*{\nabla\psi}^2 + \left(2\Lambda+O(\delta_{\Horizon})\right)\abs*{\psi}^2,\\
    \EMTensor_{4a}[\psi] = \nabla_4\psi\cdot\nabla_a\psi,\qquad
    \EMTensor_{3a}[\psi] = \nabla_3\psi\cdot\nabla_a\psi,\\
    \delta^{ab} \EMTensor_{ab}[\psi] = \nabla_4\psi\cdot\nabla_3\psi + \left(2\Lambda + O(\delta_{\Horizon})\right)\abs*{\psi}^2,\qquad
    \widehat{\EMTensor}_{ab}[\psi] = \frac{1}{2}\left(\nabla\psi\SymTracelessTensorProd\nabla\psi\right)_{ab}
    + \left(2\Lambda + O(\delta_{\Horizon})\right)\abs*{\psi}^2.
  \end{gather*}
  As a result, we have that on $\mathcal{A}$,
  \begin{align*}
    \JCurrent{\RedShiftY_{\EventHorizon},0,0}[\psi]\cdot N_{\mathcal{A}}
    ={}& \frac{1}{2}\abs*{\nabla\psi}^2
         + \frac{1}{2}\abs*{\nabla_4\psi}^2
         - a \Re\CCOneFormJ^b\nabla\psi\cdot\nabla_b\psi\\
       & + \frac{1}{2}\frac{\abs*{\Delta}}{(1+\gamma)\abs*{q}^2}\left(
         \abs*{\nabla_3\psi}^2
         - 2a\Re\CCOneFormJ^b\nabla_3\psi\cdot\nabla_b\psi
         \right)\\
       & + O\left(\delta_{\Horizon}^2\right)\abs*{\nabla_3\psi}^2
         + O(\delta_{\Horizon})\left(
         \abs*{\nabla_4\psi}^2
         + \abs*{\nabla\psi}^2
         + r^{-2}\abs*{\psi}^2
         \right).
  \end{align*}
  The same arguments as above show that for $\abs*{a}<M$ and
  $\delta_{\Horizon}$, there exists
  a constant $c_0>0$ depending on $a,M$ such that
  \begin{equation*}
    \JCurrent{\RedShiftY_{\EventHorizon},0,0}[\psi]\cdot N_{\mathcal{A}}
    \ge c_0\left(\delta_{\Horizon}\abs*{\nabla_3\psi}^2 + \abs*{\nabla_4\psi}^2 + \abs*{\nabla\psi}^2\right) + O(\delta_{\Horizon})\abs*{\psi}^2.
  \end{equation*}
  Again using the Poincar\'{e} inequality, we get that for any sphere $S=S(\tau, (1-\delta_{\Horizon})r_{\EventHorizon})$ 
  \begin{equation*}
    \int_S\JCurrent{\RedShiftY_{\EventHorizon},0,0}[\psi]\cdot N_{\mathcal{A}}
    \ge c_0\int_S\left(\delta_{\Horizon}\abs*{\nabla_3\psi}^2 + \abs*{\nabla_4\psi}^2 + \abs*{\nabla\psi}^2 + r^{-2}\abs*{\psi}^2\right),
  \end{equation*}
  for a possibly smaller $c_0>0$, as desired. 
\end{proof}

We now present the proof of \zcref[cap]{prop:redshift-main}.
\begin{proof}[Proof of Proposition \ref{prop:redshift-main}]
  Let $\RedShiftY_{\EventHorizon}$ be the vectorfield of
  \zcref[cap]{prop:redshift:main-props}. Integrating
  $\CovariantDeriv\cdot\JCurrent{\RedShiftY_{\EventHorizon},0,0}[\psi]$
  on $\Manifold(\tau_1,\tau_2)$ and applying the divergence theorem,
  the proof of \zcref{prop:redshift-main} follows
  immediately from the lower bounds
  \zcref[noname]{eq:redshift:main-props:bulk},
  \zcref[noname]{eq:redshift:main-props:boundary:Sigma}, and
  \zcref[noname]{eq:redshift:main-props:boundary:A}.
\end{proof}

We can now prove in fact the so-called enhanced redshift
estimates. The main idea here is to commute with the transverse null
vectorfield at the horizon. In particular, since we are only
interested in the redshift near the event horizon, we will commute
with with $\ein_3$, although we remark that in proving the
higher-order $r^p$-estimates in \zcref{sec:rp:higher-order}, we are essentially
taking advantage of the redshift effect near the cosmological horizon
by commuting with $\eout_4$. 
\begin{lemma}
  \label{lemma:redshift:enhanced-redshift:nab3-wave-commutation}
  If $r\le 4M$, then
  \begin{equation*}
    [\nabla_3,\WaveOpHork{2}]
    = - \partial_r\left(\frac{\Delta}{\abs*{q}^2}\right)\nabla_3^2\psi
    + O(1)\nabla\nabla_3\psi
    + O(1)\nabla_4\nabla_3\psi
    + O(1)\WaveOpHork{2}\psi
    + O(1)\frakWeightedDeriv^{\le 1}\psi.
  \end{equation*}
\end{lemma}
\begin{proof}
  We recall the expression for $\WaveOpHork{2}$ from
  \zcref[noname]{eq:wave-operator-null-frame-decomp:wave-2}
  \begin{equation}
    \begin{split}
      \WaveOpHork{2}\psi
      ={}& - \nabla_4\nabla_3\psi - \frac{1}{2}\Trace\chiBar\nabla_4\psi
           + \left(2\omega - \frac{1}{2}\Trace\chi\right)\nabla_3\psi
           + \Laplace_2\psi
           + 2\etaBar\cdot\nabla\psi
          + 2\ImagUnit \left(\LeftDual{\rho} - \eta\wedge \etaBar\right)\psi.
    \end{split}    
  \end{equation}
  As a result, we infer that for $r\le 4M$,
  \begin{align*}
    \left[\nabla_3, \WaveOpHork{2}\right]
    ={}& \left[
         \nabla_3,
         -\nabla_4\nabla_3
         - \frac{1}{2}\Trace \chiBar\nabla_4
         + \left(2\omega - \frac{1}{2}\Trace\chi\right)\nabla_3
         + \Laplace_2
         + 2\etaBar\cdot\nabla
         \right]\psi
       + \left[\nabla_3, 2\ImagUnit \left(\LeftDual{\rho} - \eta\wedge \etaBar\right)
         \right]\psi\\
    ={}& - \left[\nabla_3,\nabla_4\right]\nabla_3\psi
         + \left[\nabla_3,\Laplace_2\right]\psi
         + O(1)\frakWeightedDeriv^{\le 1}\psi\\
    ={}& 2\omega \nabla_3^2\psi
         + 2(\eta-\etaBar)\cdot\nabla\nabla_3\psi
         + 2\omegaBar\nabla_4\nabla_3\psi
         - 2\chiBar\cdot\nabla\cdot\nabla\psi\\
       &+ 2\left(\eta-\zeta\right)\cdot\nabla\nabla_3\psi
         + 2\xiBar\cdot\nabla_4\nabla\psi
         + O(1)\frakWeightedDeriv^{\le 1}\psi \\
    ={}& 2\omega\nabla_3^2\psi
         + 2(\eta-\etaBar)\cdot\nabla\nabla_3\psi
         - \Trace \chiBar\LaplaceHor_2\psi
         + O(1)\frakWeightedDeriv^{\le 1}\psi, 
  \end{align*}
  where we used the commutators from
  \zcref[cap]{corollary:commutation-formula:B-applied} and in the last step,
  we used the fact that in the principal ingoing frame,
  $\xiBar=\omegaBar=0$. Then, using \zcref[noname]{eq:wave-operator-null-frame-decomp:wave-2} again
  to rewrite the $\LaplaceHor_2$ term on the \RHS{} yields that
  \begin{align*}
    \left[\nabla_3, \WaveOpHork{2}\right]
    ={} 2\omega\nabla_3^2\psi
         + 2(\eta-\etaBar)\cdot\nabla\nabla_3\psi
         - \Trace\chiBar\nabla_4\nabla_3\psi
         - \Trace \chiBar\WaveOpHork{2}\psi
         + O(1)\frakWeightedDeriv^{\le 1}\psi.
  \end{align*}
  Plugging in the values for the Ricci coefficients from 
  \zcref[cap]{lemma:Kerr:ingoing-PG:Ric-and-curvature} then yields the result.
\end{proof}

We then have the following corollary of \zcref[cap]{prop:redshift-main}.
\begin{corollary}
  \label{coro:redshift:e3-commuted-estimate}
  Let $\psi$ be a solution in $\Manifold\cap \{r\le 4M\}$ to
  \begin{equation}
    \label{eq:redshift:e3-commuted-estimate:commuted-eqn}
    \WaveOpHork{2}\psi - V\psi
    = - \frac{4a\cos\theta}{\abs*{q}^2}\LeftDual{\nabla}_{\KillT}\psi
    + N
    + \partial_r\left(\frac{\Delta}{\abs*{q}^2}\right)\nabla_3\psi
    + O(1)\nabla\psi
    + O(1)\nabla_4\psi.
  \end{equation}
  Then, for $\abs*{a}< M$, there exists a small enough constant
  $\delta_{\operatorname{red}}>0$ such that $\delta_{\operatorname{red}} = \delta_{\operatorname{red}}(M - \abs*{a})$ and $\delta_{\operatorname{red}}\ge \delta_{\Horizon}$ and a small constant $c_0=c_0(M,a) > 0$ such that the following
  holds true in $\Manifold(\tau_1,\tau_2)$
  \begin{equation}
    \label{eq:redshift:e3-commuted-estimate}
    \begin{split}
      &c_0\EnergyFlux_{r\le r_{\EventHorizon}(1+\delta_{\operatorname{red}})}[\psi](\tau_2)
    + c_0\MorNorm_{\operatorname{red}}[\psi](\tau_1,\tau_2)
    + c_0\SpacelikeFlux_{\mathcal{A}}[\psi](\tau_1,\tau_2)\\
    \le{}& \EnergyFlux_{r\le r_{\EventHorizon}(1+2\delta_{\operatorname{red}})}[\psi](\tau_1)
      + \delta_{\operatorname{red}}^{-1}\MorNorm_{r_{\EventHorizon}(1+\delta_{\operatorname{red}})\le r\le r_{\EventHorizon}(1+2\delta_{\operatorname{red}})}[\psi](\tau_1,\tau_2)
       + \int_{\Manifold(\tau_1,\tau_2)\bigcap\curlyBrace*{\frac{r}{r_{\EventHorizon}}\le 1+2\delta_{\operatorname{red}}}}\abs*{N}^2.
    \end{split}    
  \end{equation}
\end{corollary}

\begin{proof}
  Let $\RedShiftY_{\Horizon}$ be the redshift vectorfield as
  constructed in \zcref[cap]{prop:redshift:main-props}. As shown
  in \zcref[noname]{eq:redshift:Div-J-Yh:cutoff-expansion}, we have that
  \begin{equation*}
    \CovariantDeriv\cdot\JCurrent{\RedShiftY_{\EventHorizon}, 0, 0}[\psi]
    ={}\chi_{\EventHorizon}\CovariantDeriv\cdot\JCurrent{\RedShiftY_0, 0, 0}[\psi]
    - O\left(r_{\EventHorizon}^{-1}\delta_{\operatorname{red}}^{-1}\right)\bOne_{\delta_{\operatorname{red}}\le \abs*{\frac{r}{r_{\EventHorizon}}-1}\le 2\delta_{\operatorname{red}}}
    \left(
      \abs*{\nabla_3\psi}^2
      + \abs*{\nabla_4\psi}^2
      + \abs*{\nabla\psi}^2
      + r^{-2}\abs*{\psi}^2
    \right).
  \end{equation*}
  Then, using the support properties of $\chi_{\Horizon}$, and
  \zcref[noname]{eq:redshift:main-bulk:before-simplification}, we have that
  for any sphere $S=S(\tau,r)$, 
  \begin{equation*}    
    \begin{split}
      \int_S\CovariantDeriv\cdot\JCurrent{\RedShiftY_\Horizon,0,0}[\psi]
      \ge{}& \evalAt*{\frac{\partial_r\Delta}{128r^2}}_{r=r_{\EventHorizonFuture}}
             \int_S\left( \abs*{\nabla_3\psi}^2 + \abs*{\nabla_4\psi}^2 + \abs*{\nabla\psi}^2 + r^{-2}\abs*{\psi}^2 \right)\bOne_{\abs*{\frac{r}{r_{\EventHorizonFuture}}}\le\delta_{\operatorname{red}}}\\
           &  + \frac{1}{4}\left(
             \partial_rd(r_{\EventHorizonFuture})
             - \evalAt*{\frac{c_1}{\partial_r\Delta}}_{r=r_{\EventHorizonFuture}}             
             \right)\int_S\abs*{\nabla_4\psi}^2\bOne_{\abs*{\frac{r}{r_{\EventHorizonFuture}}}\le\delta_{\operatorname{red}}}
      \\
           & + \frac{1}{4}\left(
             \partial_r\dBar(r_{\EventHorizonFuture})
             - \evalAt*{\frac{c_1}{\partial_r\Delta}}_{r=r_{\EventHorizonFuture}}
             \right)\int_S\abs*{\nabla\psi}^2\bOne_{\abs*{\frac{r}{r_{\EventHorizonFuture}}}\le\delta_{\operatorname{red}}}\\
           & + \int_S\RedShiftY_\Horizon\psi\cdot\left(\WaveOpHork{2}\psi - V\psi + \frac{4a\cos\theta}{\abs*{q}^2}\LeftDual{\nabla}_{\KillT}\psi\right)\\
      & - O(r_{\EventHorizonFuture}^{-1}\delta_{\operatorname{red}}^{-1})\int_S\left( \abs*{\nabla_3\psi}^2 + \abs*{\nabla_4\psi}^2 + \abs*{\nabla\psi}^2 + r^{-2}\abs*{\psi}^2 \right)\bOne_{\delta_{\operatorname{red}}\le \abs*{\frac{r}{r_{\EventHorizonFuture}}-1}\le2\delta_{\operatorname{red}}}
             .
    \end{split}    
  \end{equation*}
  Using the fact that $\psi$ solves
  \zcref[noname]{eq:redshift:e3-commuted-estimate:commuted-eqn}, we have that
  \begin{equation*}
    \begin{split}
      \int_S\CovariantDeriv\cdot\JCurrent{\RedShiftY_\Horizon,0,0}[\psi]
      \ge{}& \evalAt*{\frac{\partial_r\Delta}{128r^2}}_{r=r_{\EventHorizonFuture}}
             \int_S\left( \abs*{\nabla_3\psi}^2 + \abs*{\nabla_4\psi}^2 + \abs*{\nabla\psi}^2 + r^{-2}\abs*{\psi}^2 \right)\bOne_{\abs*{\frac{r}{r_{\EventHorizonFuture}}}\le\delta_{\operatorname{red}}}\\
           &  + \frac{1}{4}\left(
             \partial_rd(r_{\EventHorizonFuture})
             - \evalAt*{\frac{c_1}{\partial_r\Delta}}_{r=r_{\EventHorizonFuture}}             
             \right)\int_S\abs*{\nabla_4\psi}^2\bOne_{\abs*{\frac{r}{r_{\EventHorizonFuture}}}\le\delta_{\operatorname{red}}}
      \\
           & + \frac{1}{4}\left(
             \partial_r\dBar(r_{\EventHorizonFuture})
             - \evalAt*{\frac{c_1}{\partial_r\Delta}}_{r=r_{\EventHorizonFuture}}
             \right)\int_S\abs*{\nabla\psi}^2\bOne_{\abs*{\frac{r}{r_{\EventHorizonFuture}}}\le\delta_{\operatorname{red}}}\\
           & + \int_S\RedShiftY_\Horizon\psi\cdot\left(N
             + \partial_r\left(\frac{\Delta}{\abs*{q}^2}\right)\nabla_3\psi
             + O(1)\nabla\psi
             + O(1)\nabla_4\psi\right)\\
           & - O(r_{\EventHorizonFuture}^{-1}\delta_{\operatorname{red}}^{-1})\int_S\left( \abs*{\nabla_3\psi}^2 + \abs*{\nabla_4\psi}^2 + \abs*{\nabla\psi}^2 + r^{-2}\abs*{\psi}^2 \right)\bOne_{\delta_{\operatorname{red}}\le \abs*{\frac{r}{r_{\EventHorizonFuture}}-1}\le2\delta_{\operatorname{red}}}
             .
    \end{split}    
  \end{equation*}
  Then, we observe that 
  \begin{align*}
    &\int_S\RedShiftY_\Horizon\psi\cdot\left(\partial_r\left(\frac{\Delta}{\abs*{q}^2}\right)\nabla_3\psi
    + O(1)\nabla\psi
    + O(1)\nabla_4\psi\right)\\
    ={}& \int_S\chi_{\Horizon}\left(
         (1+O(\delta_{\operatorname{red}}))\nabla_3
         + O(1)\nabla_4
         + O(1)\nabla
         \right)\psi
         \cdot
         \left(
         \partial_r\left(\frac{\Delta}{\abs*{q}^2}\right)\nabla_3\psi
         + O(1)\nabla\psi
         + O(1)\nabla_4\psi
         \right).
  \end{align*}
  Recall that $\delta_{\operatorname{red}}$ can be chosen so that
  $\partial_r\left(\frac{\Delta}{\abs*{q}^2}\right)>0$ on the support
  of $\chi_{\Horizon}$. We thus get
  \begin{align*}
    &\int_S\RedShiftY_\Horizon\psi\cdot\left(\partial_r\left(\frac{\Delta}{\abs*{q}^2}\right)\nabla_3\psi
    + O(1)(\nabla\psi
      + \nabla_4\psi)\right)
    \ge{} -O(1)\int_S\chi_{\Horizon}\left(
           \abs*{\nabla_4\psi}^2
           + \abs*{\nabla\psi}^2
           + \abs*{\nabla_3\psi}\left(\abs*{\nabla_4\psi} + \abs*{\nabla\psi}\right)
           \right).
  \end{align*}
  Thus, choosing $\partial_rd, \partial_r\dBar$ sufficiently large, we have that
  \begin{align*}
    \int_S\CovariantDeriv\cdot\JCurrent{\RedShiftY_{\Horizon},0,0}[\psi]
    \ge{}& \evalAt*{\frac{\partial_r\Delta}{256r^2}}_{r=r_{\EventHorizon}}\int_S\left( \abs*{\nabla_3\psi}^2 + \abs*{\nabla_4\psi}^2 + \abs*{\nabla\psi}^2 + r^{-2}\abs*{\psi}^2 \right)\bOne_{\abs*{\frac{r}{r_{\EventHorizonFuture}}}\le\delta_{\operatorname{red}}}\\
         & - O(r_{\EventHorizonFuture}^{-1}\delta_{\operatorname{red}}^{-1})\int_S\left( \abs*{\nabla_3\psi}^2 + \abs*{\nabla_4\psi}^2 + \abs*{\nabla\psi}^2 + r^{-2}\abs*{\psi}^2 \right)\bOne_{\delta_{\operatorname{red}}\le \abs*{\frac{r}{r_{\EventHorizonFuture}}}\le2\delta_{\operatorname{red}}}\\
         & - O(1)\int_S\abs*{N}^2\bOne_{\abs*{\frac{r}{r_{\EventHorizonFuture}}-1}\le\delta_{\operatorname{red}}}.
  \end{align*}
  Integrating
  $\CovariantDeriv\cdot\JCurrent{\RedShiftY_{\Horizon},0,0}[\psi]$ on
  $\Manifold(\tau_1,\tau_2)$ and applying the divergence theorem, the
  proof of \zcref[cap]{coro:redshift:e3-commuted-estimate} follows
  then from the above lower bounds for
  $\CovariantDeriv\cdot\JCurrent{\RedShiftY_{\Horizon},0,0}[\psi]$ and
  from the lower bound for
  $\JCurrent{\RedShiftY_{\Horizon},0,0}[\psi]\cdot N_{\Sigma}$ and
  $\JCurrent{\RedShiftY_{\Horizon},0,0}[\psi]\cdot N_{\mathcal{A}}$
  derived in \zcref[cap]{prop:redshift:main-props}.
\end{proof}

\paragraph{Step 1: Commute with $\HorLieDeriv_{\KillT}, \HorLieDeriv_{\KillPhi}$.}

Since $\KillT$ and $\KillPhi$ are Killing vectorfields, we have
\begin{align*}
  \WaveOpHork{2}\HorLieDeriv_{\KillT}\psi
  - V\HorLieDeriv_{\KillT}\psi
  ={}& - \frac{4a\cos\theta}{\abs*{q}^2}\LeftDual{\nabla_{\KillT}}\HorLieDeriv_{\KillT}\psi
  + \HorLieDeriv_{\KillT}N ,& 
  \WaveOpHork{2}\HorLieDeriv_{\KillPhi}\psi
  - V\HorLieDeriv_{\KillPhi}\psi
  ={}& - \frac{4a\cos\theta}{\abs*{q}^2}\LeftDual{\nabla_{\KillPhi}}\HorLieDeriv_{\KillPhi}\psi
  + \HorLieDeriv_{\KillPhi}N.
\end{align*}
We can then immediately apply the estimate in \zcref[noname]{eq:redshift-main}
with $s=0$ to the commuted equations, which yields
\begin{equation*}
  \begin{split}
      &c_0\EnergyFlux_{r\le r_{\EventHorizon}(1+\delta_{\operatorname{red}})}[\left( \HorLieDeriv_{\KillT}, \HorLieDeriv_{\KillPhi} \right)\psi](\tau_2)
    + c_0\MorNorm_{\operatorname{red}}[\left( \HorLieDeriv_{\KillT}, \HorLieDeriv_{\KillPhi} \right)\psi](\tau_1,\tau_2)
    + c_0\SpacelikeFlux_{\mathcal{A}}[\left( \HorLieDeriv_{\KillT}, \HorLieDeriv_{\KillPhi} \right)\psi](\tau_1,\tau_2)\\
    \le{}& \EnergyFlux_{r\le r_{\EventHorizon}(1+2\delta_{\operatorname{red}})}[\left( \HorLieDeriv_{\KillT}, \HorLieDeriv_{\KillPhi} \right)\psi](\tau_1)
      + \delta_{\operatorname{red}}^{-3}\MorNorm_{r_{\EventHorizon}(1+\delta_{\operatorname{red}})\le r\le r_{\EventHorizon}(1+2\delta_{\operatorname{red}})}[\left( \HorLieDeriv_{\KillT}, \HorLieDeriv_{\KillPhi} \right)\psi](\tau_1,\tau_2)\\
      & + \int_{\Manifold(\tau_1,\tau_2)\bigcap\curlyBrace*{\frac{r}{r_{\EventHorizon}}\le 1+2\delta_{\operatorname{red}}}}\abs*{\frakWeightedDeriv^{\le 1}N}^2.
    \end{split}    
\end{equation*}
Using the relationship between the horizontal Lie derivative and the
horizontal covariant derivative in \zcref[cap]{lemma:LieT-LieZ-to-nabT-nabZ-comparison}, we have that 
\begin{equation}
  \label{eq:redshift:higher-order:nab-T-nab-Z-commutation}
  \begin{split}
    &c_0\EnergyFlux_{r\le r_{\EventHorizon}(1+\delta_{\operatorname{red}})}[\left( \nabla_{\KillT}, \nabla_{\KillPhi} \right)\psi](\tau_2)
      + c_0\MorNorm_{\operatorname{red}}[\left( \nabla_{\KillT}, \nabla_{\KillPhi} \right)\psi](\tau_1,\tau_2)
      + c_0\SpacelikeFlux_{\mathcal{A}}[\left( \nabla_{\KillT}, \nabla_{\KillPhi} \right)\psi](\tau_1,\tau_2)\\
    \le{}& \EnergyFlux_{r\le r_{\EventHorizon}(1+2\delta_{\operatorname{red}})}[\left( \nabla_{\KillT}, \nabla_{\KillPhi} \right)\psi](\tau_1)
           + \delta_{\operatorname{red}}^{-3}\MorNorm_{r_{\EventHorizon}(1+\delta_{\operatorname{red}})\le r\le r_{\EventHorizon}(1+2\delta_{\operatorname{red}})}[\left( \nabla_{\KillT}, \nabla_{\KillPhi} \right)\psi](\tau_1,\tau_2)\\
    & + \int_{\Manifold(\tau_1,\tau_2)\bigcap\curlyBrace*{\frac{r}{r_{\EventHorizon}}\le 1+2\delta_{\operatorname{red}}}}\abs*{\frakWeightedDeriv^{\le 1}N}^2.
  \end{split}    
\end{equation}

\paragraph{Step 2: Commute with $\abs*{q}\HodgeOp{2}$.}

Next, we commute the model equation \zcref[noname]{eq:model-problem-gRW} by
$\abs*{q}\HodgeOp{2}$. Using 
\zcref[noname]{eq:wave-operator-basic-commutation:Hodge-op}, we have that
\begin{equation*}
  \WaveOpHork{1}\left(\abs*{q}\HodgeOp{2}\psi\right)
  - V\abs*{q}\HodgeOp{2}\psi
  = - \frac{4a\cos\theta}{\abs*{q}^2}\LeftDual{\nabla}_{\KillT}\left(\abs*{q}\HodgeOp{2}\psi\right)
  + N_{\abs*{q}\HodgeOp{2}},
\end{equation*}
where, since $\psi$ solves \zcref[noname]{eq:model-problem-gRW}, we have
\begin{equation*}
  N_{\abs*{q}\HodgeOp{2}}
  = q\HodgeOp{2}N
  + O(r^{-2})\frakWeightedDeriv\HorLieDeriv_{\KillT}\psi
  + O(ar^{-2})\frakWeightedDeriv^{\le 1}\psi
  + O(ar^{-3})\frakWeightedDeriv^{\le 2}\psi.
\end{equation*}
We can then apply \zcref[noname]{eq:redshift-main} with $s=0$ again so
that
\begin{equation*}
  \begin{split}
    &c_0\EnergyFlux_{r\le r_{\EventHorizon}(1+\delta_{\operatorname{red}})}[\abs*{q}\HodgeOp{2}\psi](\tau_2)
      + c_0\MorNorm_{\operatorname{red}}[\abs*{q}\HodgeOp{2}\psi](\tau_1,\tau_2)
      + c_0\SpacelikeFlux_{\mathcal{A}}[\abs*{q}\HodgeOp{2}\psi](\tau_1,\tau_2)\\
    \le{}& \EnergyFlux_{r\le r_{\EventHorizon}(1+2\delta_{\operatorname{red}})}[\abs*{q}\HodgeOp{2}\psi](\tau_1)
           + \delta_{\operatorname{red}}^{-3}\MorNorm_{r_{\EventHorizon}(1+\delta_{\operatorname{red}})\le r\le r_{\EventHorizon}(1+2\delta_{\operatorname{red}})}[\abs*{q}\HodgeOp{2}\psi](\tau_1,\tau_2)\\
    & + \int_{\Manifold(\tau_1,\tau_2)\bigcap\curlyBrace*{\frac{r}{r_{\EventHorizon}}\le 1+2\delta_{\operatorname{red}}}}\abs*{N_{\abs*{q}\HodgeOp{2}}}^2.
  \end{split}    
\end{equation*}
Next, observe that
\begin{equation*}
  \int_{\Manifold(\tau_1,\tau_2)}\abs*{N_{\abs*{q}\HodgeOp{2}}}^2
  \lesssim{} \int_{\Manifold(\tau_1,\tau_2)}\abs*{q\HodgeOp{2}N}^2
  + \MorNorm_{\operatorname{red}}[\HorLieDeriv_{\KillT}\psi](\tau_1,\tau_2)
  + a \MorNorm_{\operatorname{red}}^1[\psi](\tau_1,\tau_2).
\end{equation*}
Then, using \zcref[noname]{eq:redshift:higher-order:nab-T-nab-Z-commutation},
and the Hodge estimates in \zcref[cap]{lemma:hodge:elliptic}, we have
that
\begin{equation*}
  \begin{split}
    &c_0\EnergyFlux_{r\le r_{\EventHorizon}(1+\delta_{\operatorname{red}})}[(\nabla_{\KillT},\frakWeightedDerivAngular)\psi](\tau_2)
      + c_0\MorNorm_{\operatorname{red}}[(\nabla_{\KillT},\frakWeightedDerivAngular)\psi](\tau_1,\tau_2)
      + c_0\SpacelikeFlux_{\mathcal{A}}[(\nabla_{\KillT},\frakWeightedDerivAngular)\psi](\tau_1,\tau_2)\\
    \le{}& \EnergyFlux_{r\le r_{\EventHorizon}(1+2\delta_{\operatorname{red}})}^1[\psi](\tau_1)
           + \delta_{\operatorname{red}}^{-3}\MorNorm_{r_{\EventHorizon}(1+\delta_{\operatorname{red}})\le r\le r_{\EventHorizon}(1+2\delta_{\operatorname{red}})}^1[\psi](\tau_1,\tau_2)\\
    & + \int_{\Manifold(\tau_1,\tau_2)}\abs*{\frakWeightedDeriv^{\le 1}N}^2
     + a \MorNorm_{\operatorname{red}}^1[\psi](\tau_1,\tau_2).
  \end{split}    
\end{equation*}

\paragraph{Step 3: Commute with $e_3$.}  

Next, we commute \zcref[noname]{eq:model-problem-gRW} with
$\nabla_3$. Then, from
\zcref[cap]{lemma:redshift:enhanced-redshift:nab3-wave-commutation}
and \zcref[cap]{coro:redshift:e3-commuted-estimate}, we have that
\begin{equation*}
  \begin{split}
    &c_0\EnergyFlux_{r\le r_{\EventHorizon}(1+\delta_{\operatorname{red}})}[\nabla_3\psi](\tau_2)
      + c_0\MorNorm_{\operatorname{red}}[\nabla_3\psi](\tau_1,\tau_2)
      + c_0\SpacelikeFlux_{\mathcal{A}}[\nabla_3\psi](\tau_1,\tau_2)\\
    \le{}& \EnergyFlux_{r\le r_{\EventHorizon}(1+2\delta_{\operatorname{red}})}[\nabla_3\psi](\tau_1)
           + \delta_{\operatorname{red}}^{-3}\MorNorm_{r_{\EventHorizon}(1+\delta_{\operatorname{red}})\le r\le r_{\EventHorizon}(1+2\delta_{\operatorname{red}})}[\nabla_3\psi](\tau_1,\tau_2)
     + \int_{\Manifold(\tau_1,\tau_2)\bigcap\curlyBrace*{\frac{r}{r_{\EventHorizon}}\le 1+2\delta_{\operatorname{red}}}}\abs*{N}^2.
  \end{split}    
\end{equation*}
Since $(\nabla_{\KillT}, \frakWeightedDerivAngular, \nabla_3)$ span
$\frakWeightedDeriv$ away from the cosmological horizon, we have that
\begin{equation*}
  \begin{split}
    &c_0\EnergyFlux_{r\le r_{\EventHorizon}(1+\delta_{\operatorname{red}})}^1[\psi](\tau_2)
      + c_0\MorNorm_{\operatorname{red}}^1[\psi](\tau_1,\tau_2)
      + c_0\SpacelikeFlux_{\mathcal{A}}^1[\psi](\tau_1,\tau_2)\\
    \le{}& \EnergyFlux_{r\le r_{\EventHorizon}(1+2\delta_{\operatorname{red}})}^1[\psi](\tau_1)
           + \delta_{\operatorname{red}}^{-3}\MorNorm_{r_{\EventHorizon}(1+\delta_{\operatorname{red}})\le r\le r_{\EventHorizon}(1+2\delta_{\operatorname{red}})}^1[\psi](\tau_1,\tau_2)\\
    & + \int_{\Manifold(\tau_1,\tau_2)\bigcap\curlyBrace*{\frac{r}{r_{\EventHorizon}}\le 1+2\delta_{\operatorname{red}}}}\abs*{\frakWeightedDeriv^{\le 1}N}^2
      + a\MorNorm_{\operatorname{red}}^1[\psi](\tau_1,\tau_2)
      .
  \end{split}    
\end{equation*}
Then choosing $a$ sufficiently small concludes the proof of
\zcref[cap]{prop:redshift-main}.


\section{Combined estimates for the model Regge-Wheeler equation}
\label{sec:combined-estimates-on-KdS}

The main goal of this section will be to prove the following proposition.
\begin{proposition}
  \label{prop:combined-estimate:main}
  Let $\psi\in\realHorkTensor{2}$ solve \zcref[noname]{eq:model-problem-gRW}, $s\ge 0$, and
  let $\delta\le p\le 2-\delta$, then we have that
  \begin{equation}
    \label{eq:combined-estimate:main}
    \CombinedBEFNorm{p}^s[\psi](\tau_1,\tau_2)
    \lesssim  \EnergyFluxCombined^s[\psi](\tau_1)
    + \ForcingTermCombinedNorm{p}^s[\psi, N](\tau_1,\tau_2).
  \end{equation}
\end{proposition}
The main idea will be to combine the estimates we have already proven
in \zcref[cap]{prop:redshift-main}, \zcref[cap]{prop:rp:Kerr}, \zcref[cap]{prop:Morawetz:KdS:main}, and
\zcref[cap]{prop:Killing-energy-estimate}. There are two main
sets of estimates that need to be combined. In a neighborhood of the
event horizon (in particular, a neighborhood containing the ergoregion
associated to the event horizon), we need to combine the redshift
estimates, the Killing energy estimate, and the Morawetz estimate.  On
the other hand, in a neighborhood of the cosmological region (again,
in particular, in a neighborhood of the cosmological horizon
encompassing its associated ergoregion), we need to combined the
$r^p$-estimates, the Killing energy estimate, and the Morawetz
estimate.  Because we are working in the slowly-rotating regime, these
two combinations can done separately in physical space and do not
intersect. In particular, the relevant regions of the spacetime where
we are combining estimates is separated by the trapping region of the
spacetime.

\subsection{Proof of  \zcref[cap]{prop:combined-estimate:main} for $s=0$}

  We recall the estimates we have proven so far which we seek to combine.
  We have from \zcref[noname]{eq:redshift-main} that
  \begin{equation}
    \label{eq:combined-estimate:redshift}
    \begin{split}
      &c_0\EnergyFlux_{r\le r_{\EventHorizon}(1+\delta_{\operatorname{red}})}[\psi](\tau_2)
    + c_0\MorNorm_{\operatorname{red}}[\psi](\tau_1,\tau_2)
    + c_0\SpacelikeFlux_{\mathcal{A}}[\psi](\tau_1,\tau_2)\\
    \le{}& \EnergyFlux_{r\le r_{\EventHorizon}(1+2\delta_{\operatorname{red}})}[\psi](\tau_1)
      + \delta_{\operatorname{red}}^{-3}\MorNorm_{r_{\EventHorizon}(1+\delta_{\operatorname{red}})\le r\le r_{\EventHorizon}(1+2\delta_{\operatorname{red}})}[\psi](\tau_1,\tau_2)\\
      & + \int_{\Manifold(\tau_1,\tau_2)\bigcap\curlyBrace*{\frac{r}{r_{\EventHorizon}}\le 1+2\delta_{\operatorname{red}}}}\abs*{N}^2.
    \end{split}    
  \end{equation}
  We also have from the $r^p$-weighted estimate \zcref[noname]{eq:rp:Kerr}
  that there exists some\footnote{This is a slight abuse of notation,
    but we will take $c_0$ to be the minimum of the implicit constants
    from the redshift and the $r^p$-weighted estimates.} $c_0$ such that
  \begin{equation}
    \label{eq:combined estimate:rp}
    \begin{split}
      &c_0\left(\BulkNormWeighted{p, r\geq R}[\psi](\tau_1,\tau_2)
    + \EnergyFluxFar_{r\ge R}[\psi](\tau_2)
    + \SpacelikeFluxFar_{\SigmaStar}[\psi](\tau_1,\tau_2)\right)\\
      \le{}& \EnergyHorizonDeg[\psi](\tau_2)
             + \EnergyFluxWeightedOpt{p}{\frac{R}{2}}[\psi](\tau_1)
    + R^{p+2}\MorNorm_{\frac{R}{2}\le r \le R }[\psi](\tau_1,\tau_2)
    + \ForcingTermWeightedNorm{p}{R}(\tau_1,\tau_2). 
    \end{split}    
  \end{equation}
  We also have from the Morawetz estimate in
  \zcref[noname]{eq:Morawetz:KdS:main} that
  \begin{equation}
    \label{eq:combined-estimate:Morawetz}
    \begin{split}
      \MorNorm[\psi](\tau_1,\tau_2)
      \lesssim{}& \sup_{\tau\in[\tau_1,\tau_2]}\EnergyHorizonDeg[\psi](\tau)
                  + \int_{\mathcal{A}(\tau_1,\tau_2)}\abs*{\nabla_{4}\psi}^2 + \int_{\SigmaStar(\tau_1,\tau_2)}\abs*{\nabla_{3}\psi}^2  \\
                  &+ \delta_{\Horizon}\SpacelikeFlux_{\mathcal{A}}[\psi](\tau_1,\tau_2)
                  + \delta_{\Horizon}\SpacelikeFlux_{\SigmaStar}[\psi](\tau_1,\tau_2)
                + \int_{\Manifold(\tau_1,\tau_2)}\left(\abs*{\nabla_{\HprVF}\psi}+r^{-1}\abs*{\psi}\right)\abs*{N},
    \end{split}
  \end{equation}
  In addition, we have the following energy estimate from
  \zcref[noname]{eq:Killing-energy-estimate},
  \begin{equation}
    \label{eq:combined-estimate:Killing}
    \begin{split}
      & \EnergyFlux_{\operatorname{deg}}[\psi](\tau_2) 
        + \int_{\mathcal{A}(\tau_1,\tau_2)}\abs*{\nabla_{4}\psi}^2 + \int_{\SigmaStar(\tau_1,\tau_2)}\abs*{\nabla_{3}\psi}^2\\      
      \lesssim{}& \EnergyFlux_{\deg}[\psi](\tau_1)
                  + \delta_{\Horizon}\left(
                  \EnergyFlux_{r\le r_{\EventHorizon}}[\psi](\tau_2)
                  + \EnergyFluxFar_{r\ge r_{\CosmologicalHorizon}}[\psi](\tau_2)
                  + \SpacelikeFlux_{\mathcal{A}}[\psi](\tau_1, \tau_2)
                  + \SpacelikeFluxFar_{\SigmaStar}[\psi](\tau_1,\tau_2)
                  \right)\\
                &+ \frac{a}{M}\MorNorm[\psi](\tau_1,\tau_2)
                  + \abs*{\int_{\Manifold(\tau_1,\tau_2)}\nabla_{\TAlmostKilling}\psi\cdot N}
                  + \int_{\Manifold(\tau_1,\tau_2)}\abs*{N}^2.
    \end{split}    
  \end{equation}
  We also observe that 
  \begin{align*}
    \EnergyFluxCombined[\psi](\tau_2)
    \lesssim{}& \EnergyFlux_{r\le r_{\EventHorizon}(1+\delta_{\operatorname{red}})}[\psi](\tau_2)
    + \EnergyFlux_{\operatorname{deg}}[\psi](\tau_2)
    + \EnergyFluxFar_{r\ge R}[\psi](\tau_2),\\
    \rpBulkCombined{p}[\psi](\tau_1,\tau_2)
    \lesssim{}&
    \MorNorm_{\operatorname{red}}[\psi](\tau_1,\tau_2)
    + \MorNorm[\psi](\tau_1,\tau_2)
    +\rpBulkWeighted{p}{R}[\psi](\tau_1,\tau_2).
  \end{align*}
  We now proceed by first summing
  \zcref[noname]{eq:combined-estimate:Morawetz} and $C\gg1$ of
  \zcref[noname]{eq:combined-estimate:Killing} to obtain that for $a$
  sufficiently small,
  \begin{equation}
    \label{eq:combined-estimate:Morawetz-and-Killing}
    \begin{split}
      &\MorNorm[\psi](\tau_1,\tau_2)
        + C\EnergyFlux_{\operatorname{deg}}[\psi](\tau_2) 
        + \int_{\mathcal{A}(\tau_1,\tau_2)}\abs*{\nabla_{4}\psi}^2 + \int_{\SigmaStar(\tau_1,\tau_2)}\abs*{\nabla_{3}\psi}^2\\
      \lesssim{}& C\EnergyFlux_{\deg}[\psi](\tau_1)
                  + C\delta_{\Horizon}\left(
                  \EnergyFlux_{r\le r_{\EventHorizon}}[\psi](\tau_2)
                  + \EnergyFluxFar_{r\ge r_{\CosmologicalHorizon}}[\psi](\tau_2)
                  + \SpacelikeFlux_{\mathcal{A}}[\psi](\tau_1, \tau_2)
                  + \SpacelikeFluxFar_{\SigmaStar}[\psi](\tau_1,\tau_2)
                  \right)\\
                & + C\abs*{\int_{\Manifold(\tau_1,\tau_2)}\nabla_{\TAlmostKilling}\psi\cdot N}
                  + C\int_{\Manifold(\tau_1,\tau_2)}\abs*{N}^2
                  + \delta_{\Horizon}\SpacelikeFlux_{\mathcal{A}}[\psi](\tau_1,\tau_2)
                  + \delta_{\Horizon}\SpacelikeFlux_{\SigmaStar}[\psi](\tau_1,\tau_2)\\
                &+ \int_{\Manifold(\tau_1,\tau_2)}\left(\abs*{\nabla_{\HprVF}\psi}+r^{-1}\abs*{\psi}\right)\abs*{N}
                  .
    \end{split}    
  \end{equation}
  Next, let
  $\dot{\delta}^{-3} :=
  \max\left(\delta_{\operatorname{red}}^{-3},R^{\delta+2}\right)$ and
  multiply \zcref[noname]{eq:combined-estimate:redshift} with $s=0$ and
  \zcref[noname]{eq:combined estimate:rp} with $p=\delta$ by
  $\dot{\delta}^{4}$ to obtain
  \begin{equation}
    \label{eq:combined-estimate:rp-and-redshift}
    \begin{split}
      &c_0\dot{\delta}^4\left(\BulkNormWeighted{p, r\geq R}[\psi](\tau_1,\tau_2)
        + \EnergyFluxFar_{r\ge R}[\psi](\tau_2)
        + \SpacelikeFluxFar_{\SigmaStar}[\psi](\tau_1,\tau_2)\right)\\
      &c_0\dot{\delta}^4\left( \EnergyFlux_{r\le r_{\EventHorizon}(1+\delta_{\operatorname{red}})}[\psi](\tau_2)
        + \MorNorm_{\operatorname{red}}[\psi](\tau_1,\tau_2)
        + \SpacelikeFlux_{\mathcal{A}}[\psi](\tau_1,\tau_2) \right)\\
      \le{}& \dot{\delta}^4\EnergyFluxWeightedOpt\delta{\frac{R}{2}}[\psi](\tau_1)
             + \dot{\delta}\MorNorm_{\frac{R}{2}\le r \le R }[\psi](\tau_1,\tau_2)
             + \dot{\delta}^4\EnergyHorizonDeg[\psi](\tau_2)
      \\
      &+  \dot{\delta}^4\EnergyFlux_{r\le r_{\EventHorizon}(1+2\delta_{\operatorname{red}})}[\psi](\tau_1)
        + \dot{\delta}\MorNorm_{r_{\EventHorizon}(1+\delta_{\operatorname{red}})\le r\le r_{\EventHorizon}(1+2\delta_{\operatorname{red}})}[\psi](\tau_1,\tau_2)\\
      & + \dot{\delta}^4\int_{\Manifold(\tau_1,\tau_2)\bigcap\curlyBrace*{\frac{r}{r_{\EventHorizon}}\le 1+2\delta_{\operatorname{red}}}}\abs*{N}^2
        + \dot{\delta}^4\ForcingTermWeightedNorm\delta{R}(\tau_1,\tau_2)
        .
    \end{split}
  \end{equation}
  Now combining
  \zcref[noname]{eq:combined-estimate:Morawetz-and-Killing} and
  \zcref[noname]{eq:combined-estimate:rp-and-redshift}, with
  $\dot{\delta}$ sufficiently small, we have that
  \begin{equation*}
    \begin{split}
      &\MorNorm[\psi](\tau_1,\tau_2)
        + C\EnergyFlux_{\operatorname{deg}}[\psi](\tau_2) 
        + \int_{\mathcal{A}(\tau_1,\tau_2)}\abs*{\nabla_{4}\psi}^2 + \int_{\SigmaStar(\tau_1,\tau_2)}\abs*{\nabla_{3}\psi}^2
      \\
     & +  c_0\dot{\delta}^4\left(\BulkNormWeighted{p, r\geq R}[\psi](\tau_1,\tau_2)
        + \EnergyFluxFar_{r\ge R}[\psi](\tau_2)
        + \SpacelikeFluxFar_{\SigmaStar}[\psi](\tau_1,\tau_2)\right)\\
      &c_0\dot{\delta}^4\left( \EnergyFlux_{r\le r_{\EventHorizon}(1+\delta_{\operatorname{red}})}[\psi](\tau_2)
        + \MorNorm_{\operatorname{red}}[\psi](\tau_1,\tau_2)
        + \SpacelikeFlux_{\mathcal{A}}[\psi](\tau_1,\tau_2) \right)\\
\lesssim{}& C\EnergyFlux_{\deg}[\psi](\tau_1)
                  + C\delta_{\Horizon}\left(
                  \EnergyFlux_{r\le r_{\EventHorizon}}[\psi](\tau_2)
                  + \EnergyFluxFar_{r\ge r_{\CosmologicalHorizon}}[\psi](\tau_2)
                  + \SpacelikeFlux_{\mathcal{A}}[\psi](\tau_1, \tau_2)
                  + \SpacelikeFluxFar_{\SigmaStar}[\psi](\tau_1,\tau_2)
                  \right)\\
                & + C\abs*{\int_{\Manifold(\tau_1,\tau_2)}\nabla_{\TAlmostKilling}\psi\cdot N}
                  + C\int_{\Manifold(\tau_1,\tau_2)}\abs*{N}^2
                  + \delta_{\Horizon}\SpacelikeFlux_{\mathcal{A}}[\psi](\tau_1,\tau_2)
                  + \delta_{\Horizon}\SpacelikeFlux_{\SigmaStar}[\psi](\tau_1,\tau_2)\\
      &+ \int_{\Manifold(\tau_1,\tau_2)}\left(\abs*{\nabla_{\HprVF}\psi}+r^{-1}\abs*{\psi}\right)\abs*{N}
        + \dot{\delta}^4\int_{\Manifold(\tau_1,\tau_2)\bigcap\curlyBrace*{\frac{r}{r_{\EventHorizon}}\le 1+2\delta_{\operatorname{red}}}}\abs*{N}^2
        + \dot{\delta}^4\ForcingTermWeightedNorm\delta{R}(\tau_1,\tau_2).
    \end{split}
  \end{equation*}
  Now taking $\delta_{\Horizon}$ sufficiently small so that $ C\delta_{\Horizon}\ll \dot{\delta}^4$ we find that
  \begin{equation*}
    \begin{split}
      &\MorNorm[\psi](\tau_1,\tau_2)
       + C\EnergyFlux_{\operatorname{deg}}[\psi](\tau_2) + \int_{\mathcal{A}(\tau_1,\tau_2)}\abs*{\nabla_{4}\psi}^2 + \int_{\SigmaStar(\tau_1,\tau_2)}\abs*{\nabla_{3}\psi}^2
      \\
     & +  c_0\dot{\delta}^4\left(\BulkNormWeighted{p, r\geq R}[\psi](\tau_1,\tau_2)
        + \EnergyFluxFar_{r\ge R}[\psi](\tau_2)
        + \SpacelikeFluxFar_{\SigmaStar}[\psi](\tau_1,\tau_2)\right)\\
      &c_0\dot{\delta}^4\left( \EnergyFlux_{r\le r_{\EventHorizon}(1+\delta_{\operatorname{red}})}[\psi](\tau_2)
        + \MorNorm_{\operatorname{red}}[\psi](\tau_1,\tau_2)
        + \SpacelikeFlux_{\mathcal{A}}[\psi](\tau_1,\tau_2) \right)\\
      \lesssim{}& C\EnergyFlux_{\deg}[\psi](\tau_1)
                  + \abs*{\int_{\Manifold(\tau_1,\tau_2)}\nabla_{\TAlmostKilling}\psi\cdot N}
                  + \int_{\Manifold(\tau_1,\tau_2)}\abs*{N}^2
                  + \int_{\Manifold(\tau_1,\tau_2)}\left(\abs*{\nabla_{\HprVF}\psi}+r^{-1}\abs*{\psi}\right)\abs*{N}\\
      &  + \int_{\Manifold(\tau_1,\tau_2)\bigcap\curlyBrace*{\frac{r}{r_{\EventHorizon}}\le 1+2\delta_{\operatorname{red}}}}\abs*{N}^2
        + \ForcingTermWeightedNorm\delta{R}(\tau_1,\tau_2).
    \end{split}
  \end{equation*}
  This shows that
  \begin{equation*}
    \begin{split}
      \CombinedBEFNorm{\delta}[\psi](\tau_1,\tau_2)
      \lesssim{}& \EnergyFluxCombined[\psi](\tau_1)
                  + \ForcingTermCombinedNorm{\delta}[\psi, N](\tau_1,\tau_2)
                  + C \abs*{\int_{\Manifold(\tau_1,\tau_2)}\nabla_{\TAlmostKilling}\psi\cdot N}
                  ,
    \end{split}
  \end{equation*}
  which combined with \zcref[noname]{eq:rp:Kerr} (and $\dot{\delta}$
  potentially even smaller), shows that for $\delta\le p\le 2-\delta$,
  \begin{equation*}
    \begin{split}
      \CombinedBEFNorm{p}[\psi](\tau_1,\tau_2)
      \lesssim{}& \EnergyFluxCombined[\psi](\tau_1)
                  + \ForcingTermCombinedNorm{p}[\psi, N](\tau_1,\tau_2)
                  + C \abs*{\int_{\Manifold(\tau_1,\tau_2)}\nabla_{\TAlmostKilling}\psi\cdot N}
                  .
    \end{split}
  \end{equation*}
  It remains to bound the last term on the \RHS. To this end,
  \begin{align*}
    \abs*{\int_{\Manifold(\tau_1,\tau_2)}\nabla_{\TAlmostKilling}\psi\cdot N}
  \lesssim{}& \abs*{\int_{\Manifold_{\operatorname{trap}}(\tau_1,\tau_2)}\nabla_{\TAlmostKilling}\psi\cdot N}
              + \abs*{\int_{\Manifold_{\cancel{\operatorname{trap}}}(\tau_1,\tau_2)}\nabla_{\TAlmostKilling}\psi\cdot N}
    \\
    \lesssim{}& \abs*{\int_{\Manifold_{\operatorname{trap}}(\tau_1,\tau_2)}\nabla_{\TAlmostKilling}\psi\cdot N}
                + \varepsilon\BulkNormWeighted{p}[\psi](\tau_1,\tau_2)
                + \varepsilon^{-1}\int_{\Manifold_{\cancel{\operatorname{trap}}(\tau_1,\tau_2)}}\abs*{N}^2\\
    \lesssim{}& \varepsilon\BulkNormWeighted{p}[\psi](\tau_1,\tau_2)
                +  \varepsilon^{-1}\ForcingTermCombinedNorm{p}[\psi, N](\tau_1,\tau_2),
  \end{align*}
  for any $\varepsilon>0$. Then, for $\varepsilon$ sufficiently small,
  we can absorb the
  $\varepsilon\BulkNormWeighted{p}[\psi](\tau_1,\tau_2)$ term onto the
  right-hand side, resulting in the desired estimate.

\subsection{Basic computations involving \texorpdfstring{$\KillT$, $\KillPhi$}{the Killing vectorfields}}

In this section, we set up a few preliminary computations that involve
$\KillT$, $\KillPhi$.

\begin{lemma}
  \label{lemma:LieT-LieZ-to-nabT-nabZ-comparison}
  Let $U$ be a horizontal covariant $k$-tensor. Then,
  \begin{equation}
    \label{eq:LieT-LieZ-to-nabT-nabZ-comparison}
    \begin{split}
      \nabla_{\KillT}U_{b_1\cdots b_k}
    ={}& \HorLieDeriv_{\KillT}U_{b_1\cdots b_k}
    + \frac{2a\cos\theta}{\abs*{q}^4}O(r,r^4\Lambda)\sum_{j=1}^k\volFormHor_{b_jc}U_{b_1\cdots c \cdots b_k}\\
    \nabla_{\KillPhi}U_{b_1\cdots b_k}
      ={}& \HorLieDeriv_{\KillPhi}U_{b_1\cdots b_k}    
           +O(1)\sum_{j=1}^k\volFormHor_{b_jc}U_{b_1\cdots c\cdots b_k}    
    \end{split}    
  \end{equation}
\end{lemma}
\begin{proof}
  For simplicity, we use the principal ingoing null frame for the
  subsequent computations, but we remark that the computations are
  equivalent in the principal outgoing frame. We have
  \begin{align*}
    2(1+\gamma)\Metric(\CovariantDeriv_b\KillT, e_c)
    ={}& \Metric\left(
      \CovariantDeriv_b\left(\ein_4 + \frac{\Delta}{\abs*{q}^2}\ein_3 - 2a(1+\gamma)\Re\CCOneFormJ_de_d, e_c\right)
    \right)\\
    ={}& \chi_{bc}^{(\operatorname{\in})}
         + \frac{\Delta}{\abs*{q}^2}\chiBar^{(\operatorname{\in})}
         - 2a(1+\gamma)\nabla_b\Re\CCOneFormJ_c\\
    ={}& \frac{1}{2}\left(\Trace\chi^{(\operatorname{in})} + \frac{\Delta}{\abs*{q}^2}\Trace\chiBar^{(\operatorname{in})}\right)\delta_{bc}
         + \frac{1}{2}\left(\aTrace\chi^{(\operatorname{in})} + \frac{\Delta}{\abs*{q}^2}\aTrace\chiBar^{(\operatorname{in})}\right)\volForm_{bc}\\
         & - a(1+\gamma)\Divergence\Re\CCOneFormJ\delta_{bc}
         - a(1+\gamma)\Curl\Re\CCOneFormJ\delta_{bc}\\
    ={}& \left(
         \frac{2a\Delta\cos\theta}{\abs*{q}^4}
         - \frac{2a(1+\gamma)\cos\theta\left(r^2+a^2+\gamma\left(\abs*{q}^2\cos^2\theta - r^2\sin^2\theta\right)\right)}{\abs*{q}^4}
         \right).
  \end{align*}
  Then, since
  \begin{equation*}
    \HorLieDeriv_{\KillT}U_{b_1\cdots b_k}
    = \nabla_{\KillT}U_{b_1\cdots b_k}
    + \sum_{j=1}^k\Metric(\CovariantDeriv_{b_j}\KillT, e_c)U_{b_1\cdots c \cdots b_k},
  \end{equation*}
  the equation for $\nabla_{\KillT}$ in
  \zcref[noname]{eq:LieT-LieZ-to-nabT-nabZ-comparison} follows as desired. Moreover, we have
  \begin{align*}
    &2(1+\gamma)\Metric(\CovariantDeriv_b\KillPhi, e_c)\\
    ={}& \Metric\left(
         \CovariantDeriv_b\left(
         \frac{2(r^2+a^2)\sqrt{\kappa}\sin\theta}{\abs*{q}}e_2
         - a\sin^2\theta\left(
         \ein_4
         + \frac{\Delta}{\abs*{q}^2}\ein_3
         \right)
         \right), e_c
         \right),\\
    ={}& 2(r^2+a^2)\nabla_b\Re\CCOneFormJ_c
         + 4re_b(r)\Re\CCOneFormJ_c
         - a\sin^2\theta\left(\chi_{bc}^{(\operatorname{in})} + \frac{\Delta}{\abs*{q}^2}\chiBar_{bc}^{(\operatorname{in})}\right)\\
    ={}& \left(r^2+a^2\right)\Divergence\Re\CCOneFormJ\delta_{bc}
         +\left(r^2+a^2\right)\Curl\Re\CCOneFormJ\volFormHor_{bc}\\
       &  - \frac{a\sin^2\theta}{2}\left(\Trace\chi^{(\operatorname{in})} + \frac{\Delta}{\abs*{q}^2}\Trace\chiBar^{(\operatorname{in})}\right)\delta_{bc}
         - \frac{a\sin^2\theta}{2}\left(\aTrace\chi^{(\operatorname{in})} + \frac{\Delta}{\abs*{q}^2}\aTrace\chiBar^{(\operatorname{in})}\right)\volFormHor_{bc}\\
    ={}& \left( -\frac{2a^2\sin^2\theta\cos\theta\Delta}{\abs*{q}^4}
         + \frac{2\cos\theta(r^2+a^2)\left( r^2+a^2+\gamma\left(\abs*{q}^2\cos^2\theta - r^2\sin^2\theta\right) \right)}{\abs*{q}^4} \right)\volFormHor_{bc}.
  \end{align*}
  Then, since
  \begin{equation*}
    \HorLieDeriv_{\KillPhi}U_{b_1\cdots b_k}
    = \nabla_{\KillPhi}U_{b_1\cdots b_k}
    + \sum_{j=1}^k\Metric(\CovariantDeriv_{b_j}\KillPhi, e_c)U_{b_1\cdots c \cdots b_k},
  \end{equation*}
  the conclusion in \zcref[noname]{eq:LieT-LieZ-to-nabT-nabZ-comparison} follows as desired.
\end{proof}

\subsection{Higher-order estimates}

We now prove \zcref[cap]{prop:combined-estimate:main} for $s\ge
1$. Below, we illustrate how to prove
\zcref[cap]{prop:combined-estimate:main} with $s=1$ given that we
already have \zcref[cap]{prop:combined-estimate:main} for $s=0$. The
general proof follows a similar inductive argument.

\paragraph{Step 1: Commutation with $\HorLieDeriv_{\KillT}$ and
  $\HorLieDeriv_{\KillPhi}$.} Since $\KillT$ and $\KillPhi$ are Killing vectorfields, we have
\begin{align*}
  \WaveOpHork{2}\HorLieDeriv_{\KillT}\psi
  - V\HorLieDeriv_{\KillT}\psi
  ={}& - \frac{4a\cos\theta}{\abs*{q}^2}\LeftDual{\nabla_{\KillT}}\HorLieDeriv_{\KillT}\psi
  + \HorLieDeriv_{\KillT}N ,& 
  \WaveOpHork{2}\HorLieDeriv_{\KillPhi}\psi
  - V\HorLieDeriv_{\KillPhi}\psi
  ={}& - \frac{4a\cos\theta}{\abs*{q}^2}\LeftDual{\nabla_{\KillPhi}}\HorLieDeriv_{\KillPhi}\psi
  + \HorLieDeriv_{\KillPhi}N.
\end{align*}
We can then immediately apply the energy estimates in
\zcref[noname]{eq:combined-estimate:main} to the commuted equations, which
yield
\begin{equation*}
  \begin{split}
    \CombinedBEFNorm{p}[(\HorLieDeriv_{\KillT}, \HorLieDeriv_{\KillPhi})\psi](\tau_1,\tau_2)
    \lesssim{}& \EnergyFluxCombined[(\HorLieDeriv_{\KillT}, \HorLieDeriv_{\KillPhi})\psi](\tau_1)
                + \ForcingTermCombinedNorm{p}[(\HorLieDeriv_{\KillT}, \HorLieDeriv_{\KillPhi})\psi, N](\tau_1,\tau_2).
  \end{split}
\end{equation*}
Using the relationship between the horizontal Lie derivative and the
horizontal covariant derivative in \zcref[cap]{lemma:LieT-LieZ-to-nabT-nabZ-comparison}, we have that 
\begin{equation}
  \label{eq:higher-order:T-Z-commutation}
  \begin{split}
    \CombinedBEFNorm{p}[(\nabla_{\KillT}, \nabla_{\KillPhi})\psi](\tau_1,\tau_2)
    \lesssim{}& \EnergyFluxCombined[(\nabla_{\KillT}, \nabla_{\KillPhi})\psi](\tau_1)
                + \ForcingTermCombinedNorm{p}[(\nabla_{\KillT}, \nabla_{\KillPhi})\psi, N](\tau_1,\tau_2),
  \end{split}
\end{equation}

\paragraph{Step 2: Recover higher-order bulk estimates.}
In this step we prove the following lemma.
\begin{lemma}
  \label{lemma:higher-order:bulk-RHat-Angular}
  We have that
  \begin{equation}
    \label{eq:higher-order:bulk-RHat-Angular}
    \MorNorm^1[\psi](\tau_1,\tau_2)
    \lesssim \EnergyFluxCombined^1[\psi](\tau_1)
    + \ForcingTermCombinedNorm{p}^1[\psi, N](\tau_1,\tau_2)
    + \sqrt{\operatorname{EF}_p[\psi](\tau_1,\tau_2)}\sqrt{\operatorname{EF}_p^1[\psi](\tau_1,\tau_2)}
    .
  \end{equation}
\end{lemma}
\begin{proof}
  From \zcref[noname]{eq:higher-order:T-Z-commutation}, we
  already have $\MorNorm^1[(\nabla_{\KillT}, \nabla_{\KillPhi})\psi](\tau_1,\tau_2)
    \lesssim \EnergyFluxCombined^1[\psi](\tau_1)
    + \ForcingTermCombinedNorm{p}^1[\psi, N](\tau_1,\tau_2)$ so it suffices to show that
  \begin{equation}
    \label{eq:higher-order:bulk:sufficient-condition}
    \begin{split}
      \MorNorm^1[(\nabla_{\HprVF}, \nabla)\psi](\tau_1,\tau_2)
    \lesssim{}& \EnergyFluxCombined^1[\psi](\tau_1)
    + \ForcingTermCombinedNorm{p}^1[\psi, N](\tau_1,\tau_2)
     + \sqrt{\operatorname{EF}_p[\psi](\tau_1,\tau_2)}\sqrt{\operatorname{EF}_p^1[\psi](\tau_1,\tau_2)}
    .
    \end{split}    
  \end{equation}
  Using the explicit representation of the wave operator
  in \zcref[noname]{eq:wave-using-THat-RHat}, we have
  that 
  \begin{equation}
    \label{eq:higher-order:wave-THat-RHat-decomp}
    \begin{split}
     \frac{(r^2+a^2)^2}{\Delta}\nabla_{\HprVF}\nabla_{\HprVF}\psi
      + \abs*{q}^2\LaplaceHor_2\psi
      ={}& \frac{(r^2+a^2)^2}{\Delta}\nabla_{\HawkingVF}\nabla_{\HawkingVF}\psi
           - 2r\nabla_{\HprVF}\psi
           - \abs*{q}^2\WaveOpHork{2}\psi
           - \abs*{q}^2\left(\eta + \etaBar\right)\cdot\nabla\psi.
    \end{split}    
  \end{equation}
  Multiplying both sides by $\frac{\Delta}{(r^2+a^2)^2\abs*{q}^2}\nabla_{\HprVF}^2\psi$, we have that
  \begin{equation*}
    \begin{split}
      \frac{1}{\abs*{q}^2}\nabla_{\HprVF}^2\psi\cdot\nabla_{\HprVF}^2\psi
      + \frac{\Delta}{(r^2+a^2)^2}\LaplaceHor_2\psi\cdot\nabla_{\HprVF}^2\psi
      ={}& \frac{1}{\abs*{q}^2}\nabla_{\HawkingVF}\nabla_{\HawkingVF}\psi\cdot \nabla_{\HprVF}^2\psi
           - \frac{2r\Delta}{\abs*{q}^2(r^2+a^2)^2}\nabla_{\HprVF}\psi\cdot \nabla_{\HprVF}^2\psi\\
         &   - \frac{\Delta}{(r^2+a^2)^2}\WaveOpHork{2}\psi \cdot \nabla_{\HprVF}^2\psi
           - \frac{\Delta}{(r^2+a^2)^2}\left(\eta + \etaBar\right)\cdot\nabla\psi \cdot \nabla_{\HprVF}^2\psi.
    \end{split}    
  \end{equation*}
  Integrating by parts over $\Manifold(\tau_1,\tau_2)$, we have that
  \begin{equation}
    \label{eq:higher-order:bulk:HprVF2}
    \begin{split}
      &\int_{\Manifold(\tau_1,\tau_2)}\left( \frac{1}{r^2}\abs*{\nabla_{\HprVF}^2\psi}^2
      + \abs*{\frac{\Delta}{(r^2+a^2)^2}}\abs*{\nabla \nabla_{\HprVF}\psi}^2 \right)\\
      \lesssim{}& \MorNorm[\nabla_{\KillT}\psi, \nabla_{\KillPhi}\psi](\tau_1,\tau_2)
                  + \MorNorm[\psi](\tau_1,\tau_2)
                  + \int_{\mathcal{M}(\tau_1,\tau_2)}\abs*{\frac{\Delta}{(r^2+a^2)^2}\WaveOpHork{2}\psi\cdot\nabla_{\HprVF}^2\psi}\\
      & + \sqrt{\operatorname{EF}_p[\psi](\tau_1,\tau_2)}\sqrt{\operatorname{EF}_p^1[\psi](\tau_1,\tau_2)}
        + \delta_{\Horizon}\BulkNormWeighted{\delta}^1[\psi](\tau_1,\tau_2)
        ,
    \end{split}    
  \end{equation}
  where we remark that the
  $\delta_{\Horizon}\BulkNormWeighted{\delta}^1[\psi](\tau_1,\tau_2)$
  error term in the last line comes from the fact that
  $\mathcal{M}(\tau_1,\tau_2)$ contains a small region beyond the event
  and cosmological horizons.  On the other hand, multiplying
  \zcref[noname]{eq:higher-order:wave-THat-RHat-decomp} by
  $\frac{\Delta}{(r^2+a^2)^2\abs*{q}^2}\psi$, we have that
  \begin{equation*}
    \begin{split}
      \frac{1}{\abs*{q}^2}\nabla_{\HprVF}^2\psi\cdot\psi
      + \frac{\Delta}{(r^2+a^2)^2}\LaplaceHor_2\psi\cdot\psi
      ={}& \frac{1}{\abs*{q}^2}\nabla_{\HawkingVF}\nabla_{\HawkingVF}\psi\cdot \psi
           - \frac{2r\Delta}{\abs*{q}^2(r^2+a^2)^2}\nabla_{\HprVF}\psi\cdot \psi\\
         &   - \frac{\Delta}{(r^2+a^2)^2}\WaveOpHork{2}\psi \cdot \psi
           - \frac{\Delta}{(r^2+a^2)^2}\left(\eta + \etaBar\right)\cdot\nabla\psi \cdot \psi.
    \end{split}    
  \end{equation*}
  Integrating over $\Manifold(\tau_1,\tau_2)$ and integrating
  by parts, we then have that
  \begin{equation}
    \label{eq:higher-order:bulk:nabla}
    \begin{split}
      &\int_{\Manifold(\tau_1,\tau_2)}
         \left( \frac{1}{r^2}\abs*{\nabla_{\HprVF}\psi}^2+
        \abs*{\frac{\Delta}{(r^2+a^2)^2}}\abs*{\nabla\psi}^2 \right)\\
      \lesssim{}& \MorNorm[\nabla_{\KillT}\psi, \nabla_{\KillPhi}\psi](\tau_1,\tau_2)
                  + \MorNorm[\psi](\tau_1,\tau_2)  + \int_{\mathcal{M}(\tau_1,\tau_2)}\abs*{\frac{\Delta}{(r^2+a^2)^2}\WaveOpHork{2}\psi\cdot\psi}\\
      & + \sqrt{\operatorname{EF}_p[\psi](\tau_1,\tau_2)}\sqrt{\operatorname{EF}_p^1[\psi](\tau_1,\tau_2)}
        + \delta_{\Horizon}\BulkNormWeighted{\delta}^1[\psi](\tau_1,\tau_2)
        .
    \end{split}
  \end{equation}  
  On the other hand, multiplying
  \zcref[noname]{eq:higher-order:wave-THat-RHat-decomp} by
  $\frac{\Delta^2}{(r^2+a^2)^3\abs*{q}^2}\LaplaceHor_2\psi$, we have
  that
  \begin{equation*}
    \begin{split}
      &\frac{\Delta}{\abs*{q}^2(r^2+a^2)}\nabla_{\HprVF}^2\psi\cdot\LaplaceHor_2\psi
      + \frac{\Delta^2}{(r^2+a^2)^3}\LaplaceHor_2\psi\cdot\LaplaceHor_2\psi\\
      ={}& \frac{\Delta}{\abs*{q}^2(r^2+a^2)}\nabla_{\HawkingVF}\nabla_{\HawkingVF}\psi\cdot \LaplaceHor_2\psi
           - \frac{2r\Delta^2}{\abs*{q}^2(r^2+a^2)^3}\nabla_{\HprVF}\psi\cdot \LaplaceHor_2\psi\\
         &   - \frac{\Delta^2}{(r^2+a^2)^3}\WaveOpHork{2}\psi \cdot \LaplaceHor_2\psi
           - \frac{\Delta^2}{(r^2+a^2)^3}\left(\eta + \etaBar\right)\cdot\nabla\psi \cdot \LaplaceHor_2\psi.
    \end{split}    
  \end{equation*}
  Again integrating by parts over
  $\Manifold_{\cancel{\operatorname{trap}}}(\tau_1,\tau_2)$, we have
  that
  \begin{equation*}
    \begin{split}
      \int_{\Manifold_{\cancel{\operatorname{trap}}}(\tau_1,\tau_2)}\frac{\Delta^2}{(r^2+a^2)^3}\abs*{\LaplaceHor_2\psi}^2
      \lesssim{}& \MorNorm[(\nabla_{\KillT}, \nabla_{\KillPhi}, \nabla_{\HprVF})\psi](\tau_1,\tau_2)
                  + \MorNorm[\psi](\tau_1,\tau_2)
                  + \ForcingTermCombinedNorm{p}^1[\psi, N](\tau_1,\tau_2)\\
                & + \sqrt{\operatorname{EF}_p[\psi](\tau_1,\tau_2)}\sqrt{\operatorname{EF}_p^1[\psi](\tau_1,\tau_2)}
                  + \delta_{\Horizon}\BulkNormWeighted{\delta}^1[\psi](\tau_1,\tau_2)
                  ,
    \end{split}
  \end{equation*}
  where we observe that the boundary terms on
  $\abs*{r-3M}=\delta_{\operatorname{trap}}$ vanish since
  $\nabla(r)=0$. Then using \zcref[noname]{eq:laplacian-bochner:s2}, we have
  that actually,
  \begin{equation}
    \label{eq:higher-order:bulk:nabla2}
    \begin{split}
      \int_{\Manifold_{\cancel{\operatorname{trap}}}(\tau_1,\tau_2)}\frac{\Delta^2}{(r^2+a^2)^3}\abs*{\nabla^2\psi}^2
      \lesssim{}& \MorNorm[(\nabla_{\KillT}, \nabla_{\KillPhi}, \nabla_{\HprVF})\psi](\tau_1,\tau_2)
                  + \MorNorm[\psi](\tau_1,\tau_2)
                  + \ForcingTermCombinedNorm{p}^1[\psi, N](\tau_1,\tau_2)\\
                & + \sqrt{\operatorname{EF}_p[\psi](\tau_1,\tau_2)}\sqrt{\operatorname{EF}_p^1[\psi](\tau_1,\tau_2)}
                  + \delta_{\Horizon}\BulkNormWeighted{\delta}^1[\psi](\tau_1,\tau_2)
                  .
    \end{split}
  \end{equation}
  Combining \zcref[noname]{eq:higher-order:bulk:HprVF2},
  \zcref[noname]{eq:higher-order:bulk:nabla}, and
  \zcref[noname]{eq:higher-order:bulk:nabla2} with
  \zcref{prop:rp:Kerr,prop:redshift-main} for $s=1$ then shows
  \zcref[noname]{eq:higher-order:bulk:sufficient-condition}, which as
  mentioned is sufficient to conclude the proof of
  \zcref[cap]{lemma:higher-order:bulk-RHat-Angular}.
\end{proof}

\paragraph{Step 3: Recover higher-order energy estimates.}
In this step we prove the following lemma.
\begin{lemma}
  \label{lemma:higher-order:boundary}
  We have that
  \begin{equation*}
    \EnergyHorizonDegAux[\dk\psi](\tau)
    \lesssim \EnergyFluxCombined^1[\psi](\tau_1)
    + \ForcingTermCombinedNorm{p}^1[\psi, N](\tau_1,\tau_2)
    + \sqrt{\operatorname{EF}_p[\psi](\tau_1,\tau_2)}\sqrt{\operatorname{EF}_p^1[\psi](\tau_1,\tau_2)},
  \end{equation*}
  where
    \begin{equation*}
      \begin{split}
        \EnergyHorizonDegAux[\psi](\tau)
        ={}& \int_{(1-\delta_{\Horizon})r_{\EventHorizon}}^{(1+\delta_{\Horizon})r_{\CosmologicalHorizon}}
             \left(
             \frac{1}{r^2}\abs*{\HawkingVF\psi}^2
             + \frac{1}{r^2}\abs{\HprVF\psi}^2
             + \frac{\Delta}{r^2+a^2}\abs*{\nabla\psi}^2
             + r^{-2}\abs*{\psi}^2\right)(\tau, r,\omega)\,r^2drd\mathring{\gamma}.
      \end{split}      
  \end{equation*}
\end{lemma}

\begin{proof}
  From \zcref[noname]{eq:higher-order:T-Z-commutation}, we have $\EnergyHorizonDegAux[\psi](\tau)[(\nabla_{\KillT}, \nabla_{\KillPhi})\psi](\tau_1,\tau_2)
    \lesssim \EnergyFluxCombined^1[\psi](\tau_1)
    + \ForcingTermCombinedNorm{p}^1[\psi, N](\tau_1,\tau_2)$ so it suffices to show that
  \begin{equation}
    \label{eq:higher-order:boundary:sufficient-condition}
    \begin{split}
      \EnergyHorizonDegAux[\psi](\tau)[(\nabla_{\HprVF}, \nabla)\psi](\tau_1,\tau_2)
      \lesssim{}& \EnergyFluxCombined^1[\psi](\tau_1)
                  + \ForcingTermCombinedNorm{p}^1[\psi, N](\tau_1,\tau_2)
                 + \sqrt{\operatorname{EF}_p[\psi](\tau_1,\tau_2)}\sqrt{\operatorname{EF}_p^1[\psi](\tau_1,\tau_2)}
                  .
    \end{split}    
  \end{equation}
  To this end, we observe that combining
  \zcref[noname]{eq:higher-order:bulk:HprVF2,eq:higher-order:bulk:nabla2}
  and integrating across $\Manifold(\tau_0,\tau_0+1)$ for some
  $\tau_0\in (\tau_1,\tau_2-1)$ instead of integrating across
  $\Manifold(\tau_1,\tau_2)$, we have that
  \begin{equation*}
    \begin{split}
      &\int_{\Manifold(\tau_0,\tau_0+1)}\frac{\Delta}{\left( r^2+a^2 \right)^2}\abs*{\nabla^2\psi}^2
    + \frac{1}{r}\abs*{\nabla_{\HprVF}^2\psi}^2
        + \frac{1}{r^2}\abs*{\nabla_{\HprVF}\nabla\psi}^2\\
      \lesssim{}& \MorNorm^1[\psi](\tau_0,\tau_0+1)
                  + \ForcingTermCombinedNorm{p}^1[\psi, N](\tau_0,\tau_0+1)
                  + \sqrt{\operatorname{EF}_p[\psi](\tau_0,\tau_0+1)}\sqrt{\operatorname{EF}_p^1[\psi](\tau_0,\tau_0+1)},
    \end{split}        
  \end{equation*}
  Using \zcref[cap]{lemma:higher-order:bulk-RHat-Angular}, we then have that
  \begin{equation*}
    \begin{split}
      &\int_{\Manifold(\tau_0,\tau_0+1)}\frac{\Delta}{\left( r^2+a^2 \right)^2}\abs*{\nabla^2\psi}^2
    + \frac{1}{r}\abs*{\nabla_{\HprVF}^2\psi}^2
        + \frac{1}{r^2}\abs*{\nabla_{\HprVF}\nabla\psi}^2\\
      \lesssim{}& \EnergyFluxCombined^1[\psi](\tau_0)
                  +\ForcingTermCombinedNorm{p}^1[\psi, N](\tau_0,\tau_0+1)
                  + \sqrt{\operatorname{EF}_p[\psi](\tau_0,\tau_0+1)}\sqrt{\operatorname{EF}_p^1[\psi](\tau_0,\tau_0+1)}.
    \end{split}        
  \end{equation*}
  Observe that given the definition of $\EnergyFluxCombined[\psi]$ in
  \zcref[noname]{eq:energy-flux:combined:def}, we thus have that
  \begin{equation*}
    \begin{split}
      \int_{\tau_0}^{\tau_0+1}\EnergyHorizonDegAux[\psi](\tau)\,d\tau
      \lesssim{}& \EnergyFluxCombined^1[\psi](\tau_0)
                  +\ForcingTermCombinedNorm{p}^1[\psi, N](\tau_0,\tau_0+1)
                 + \sqrt{\operatorname{EF}_p[\psi](\tau_0,\tau_0+1)}\sqrt{\operatorname{EF}_p^1[\psi](\tau_0,\tau_0+1)},
    \end{split}
  \end{equation*}
  We can thus control
  \begin{equation*}
    \begin{split}
      \inf_{\tau\in [\tau_0,\tau_0+1]}\EnergyHorizonDegAux[\psi](\tau)\,d\tau
      \lesssim{}& \EnergyFluxCombined^1[\psi](\tau_0)
                  +\ForcingTermCombinedNorm{p}^1[\psi, N](\tau_0,\tau_0+1)
                + \sqrt{\operatorname{EF}_p[\psi](\tau_0,\tau_0+1)}\sqrt{\operatorname{EF}_p^1[\psi](\tau_0,\tau_0+1)},
    \end{split}    
  \end{equation*}
  as desired.
\end{proof}

\paragraph{Step 4: Closing the higher-order estimates.}

We can now prove \zcref[noname]{eq:combined-estimate:main} with $s=1$. Let $\tau_{*}\in [\tau_0, \tau_0+1]$ be such that $\EnergyFluxCombined^1[\psi](\tau_{*})
    = \inf_{\tau\in [\tau_0, \tau_0+1]}\EnergyFluxCombined^1[\psi](\tau)$. Then combining
  \zcref[cap]{lemma:higher-order:bulk-RHat-Angular,lemma:higher-order:boundary},
  we have that
  \begin{equation*}
    \begin{split}
      \EnergyHorizonDegAux^1[\psi](\tau_{*})
      + \MorNorm^1[\psi](\tau_0,\tau_{*})
      \lesssim  \EnergyFluxCombined^{1}[\psi](\tau_0)
      + \ForcingTermCombinedNorm{p}^1[\psi, N](\tau_0,\tau_*)
      + \sqrt{\operatorname{EF}_p[\psi](\tau_0,\tau_*)}\sqrt{\operatorname{EF}_p^1[\psi](\tau_0,\tau_*)}.
    \end{split}
  \end{equation*}
  Since $\tau_0$ was arbitrary, we thus in fact have that
  \begin{equation*}
    \begin{split}
      \EnergyHorizonDegAux^1[\psi](\tau_{2})
      + \MorNorm^1[\psi](\tau_1,\tau_{2})
      \lesssim  \EnergyFluxCombined^{1}[\psi](\tau_1)
      + \ForcingTermCombinedNorm{p}^1[\psi, N](\tau_1,\tau_2)
      + \sqrt{\operatorname{EF}_p[\psi](\tau_1,\tau_2)}\sqrt{\operatorname{EF}_p^1[\psi](\tau_1,\tau_2)}.
    \end{split}
  \end{equation*}  
  Adding \zcref[noname]{eq:rp:Kerr} and
  \zcref[noname]{eq:redshift-main} to the above like we did in the
  proof of \zcref[cap]{prop:combined-estimate:main} for $s=0$, we get
  \begin{equation*}
    \CombinedBEFNorm{p}^1[\psi](\tau_1, \tau_2)
    \lesssim \EnergyFluxCombined^{1}[\psi](\tau_1)
    + \ForcingTermCombinedNorm{p}^1[\psi, N](\tau_1,\tau_2)
    + \sqrt{\operatorname{EF}_p[\psi](\tau_1,\tau_2)}\sqrt{\operatorname{EF}_p^1[\psi](\tau_1,\tau_2)}.
  \end{equation*}
  Using Cauchy-Schwarz to control the last error term on the
  \RHS{} of the equation (and potentially adding a large amount of
  \zcref[noname]{eq:combined-estimate:main} with $s=0$) concludes the proof of
  \zcref[noname]{eq:combined-estimate:main} with $s=1$.



\section{Proof of the main estimate \texorpdfstring{\zcref{eq:main-q-A-combined-estimate}}{in the internal region}}
\label{sec:transport}

In this section, we finally prove the main theorem in
\zcref[cap]{MAINTHEOREM} and recover the weighted estimates for the
full Teukolsky wave-transport system in
\zcref[noname]{eq:full-RW-system:wave} and
\zcref[noname]{eq:full-RW-system:transport}.

\subsection{Preliminaries for transport equations}

\subsubsection{Divergence theorem for transport equations}

We first introduce the main application of the divergence theorem that
we will use to control the transport equations.

\begin{lemma}[Divergence lemma]
  \label{lemma:transport:div-lemma}
  Consider a vectorfield $X\in T\mathcal{M}$. We have
  \begin{equation*}
    -\int_{\mathcal{A}(\tau_1,\tau_2)}\Metric(X, N)
    - \int_{\Sigma(\tau_2)}\Metric(X, N)
    - \int_{\SigmaStar(\tau_1,\tau_2)}\Metric(X, N)
    + \int_{\Sigma(\tau_1)}\Metric(X, N)
    = \int_{\mathcal{M}(\tau_1,\tau_2)}\Divergence(X),
  \end{equation*}
  where $N$ is the normal to the boundary such that $\Metric(N,e_3)=-1$. This can be rewritten as
  \begin{equation*}
    -\int_{\partial^+\mathcal{M}(\tau_1,\tau_2)}\Metric(X, N)
    + \int_{\partial^-\mathcal{M}(\tau_1,\tau_2)}\Metric(X, N)
    = \int_{\mathcal{M}(\tau_1,\tau_2)}\Divergence(X),
  \end{equation*}
  where $\partial^+\mathcal{M}(\tau_1,\tau_2)
    = \mathcal{A}(\tau_1,\tau_2)\cup \Sigma(\tau_2)\cup \SigmaStar(\tau_1,\tau_2)$ and $\partial^-\mathcal{M}(\tau_1,\tau_2)=\Sigma(\tau_1)$.
\end{lemma}
\begin{proof}
  Basic application of the standard divergence lemma. 
\end{proof}

We will also make use of the following properties of the boundary of $\mathcal{M}(\tau_1,\tau_2)$.
\begin{lemma}
  \label{lemma:transport:boundary}
  We have the following properties on the boundary of $\mathcal{M}(\tau_1,\tau_2)$.
  \begin{enumerate}
  \item On the boundary $\mathcal{A}$,
    \begin{equation*}
      \Metric(N_{\mathcal{A}}, e_3)=-1,\qquad
      \Metric(N_{\mathcal{A}}, e_4)\le -2\delta_{\Horizon}, \qquad
      \Metric(N_{\mathcal{A}},e_a) = O(\delta_{\Horizon}).
    \end{equation*}
  \item On the boundary $\SigmaStar$ we have
    \begin{equation*}
      \Metric(N_{\SigmaStar}, e_3) = 2\frac{\Delta}{\abs*{q}^2},\qquad
      \Metric(N_{\SigmaStar},e_4)=-2 \qquad
      \Metric(N_{\SigmaStar}, e_a) = 0.
    \end{equation*}
  \item On the boundary $\Sigma(\tau)$, we have
    \begin{equation*}
      \Metric(N_{\Sigma},N_{\Sigma}) \le -\frac{M^2}{8r^2},\qquad
      \Metric(N_{\Sigma}, e_4) = - \frac{2\lambda}{r^2},\qquad
      \Metric(N_{\Sigma}, e_3) = -2.
    \end{equation*}
  \end{enumerate}
\end{lemma}
\begin{proof}
  Direct computation.
\end{proof}

We also have the following direct computation.
\begin{lemma}
  \label{lemma:transport:div-e3}
  For any function $f$, we have that
  \begin{equation*}
    \Divergence(fe_3) = e_3(f) + (-2\underline{\omega} + \Trace \underline{\chi})f.
  \end{equation*}
\end{lemma}
\begin{proof}
  Direct computation.
\end{proof}

\subsubsection{General transport lemmas}

In this section, we prove a few general transport equations that we
will apply repeatedly in what follows.

\begin{lemma}
  \label{lemma:transport:bulk}
  Suppose $\Phi_1,\Phi_2\in \realHorkTensor{2}(\Complex)$
  satisfy the differential
  relation
  \begin{equation}
    \label{eq:transport:main-lemma:Phi1-Phi2-relation}
    \ConformalInvDeriv_3\Phi_1 = \Phi_2,    
  \end{equation}
  where we also assume that $\Phi_1$ is $s$-conformally invariant. Then for all $p\ge \delta$, there exist constants
  $c_0(p), C_0(p) > 0$ such that
  \begin{equation}
    \label{eq:transport:bulk}
    c_0\abs*{q}^{p-3}\abs*{\Phi_1}^2
    \lesssim C_0(p)\abs*{q}^{p-1}\abs*{\Phi_2}^2
    - \Divergence\left(\abs*{q}^{p-2}\abs*{\Phi_1}^2e_3\right),
  \end{equation}
  and the integrated version
  \begin{equation}
    \label{eq:transport:bulk:integrated}
    \int_{\mathcal{M}(\tau_1,\tau_2)}r^{p-3}\abs*{\Phi_1}^2
    + \int_{\partial \mathcal{M}(\tau_1,\tau_2)}r^{p-2}\abs*{\Phi_1}^2
    \lesssim \int_{\mathcal{M}(\tau_1,\tau_2)}r^{p-1}\abs*{\Phi_2}^2
    + \int_{\Sigma(\tau_1)}r^{p-2}\abs*{\Phi_1}^2.
  \end{equation}
\end{lemma}
\begin{proof}
  We first observe that multiplying \zcref[noname]{eq:transport:main-lemma:Phi1-Phi2-relation} by $\overline{\Phi_1}$, we have that
  \begin{align*}
    e_3\left(\abs*{\Phi_1}^2\right)
    ={}& 2\Re \left((\Phi_2 + 2s\underline{\omega}\Phi_1)\cdot\overline{\Phi_1}\right)
    ={} 2\Re\left(\Phi_2\cdot\overline{\Phi_1}\right)
         + 4s\underline{\omega}\abs*{\Phi_1}^2
         .      
  \end{align*}
  Then, using the fact that $\frac{e_3(\abs*{q}^2)}{{\abs*{q}^2}} = \Trace \underline{\chi}$, we have that
  \begin{align*} 
    2\abs*{q}^{p-2}\Re\left(\Phi_2\cdot\overline{\Phi_1}\right)
    ={}& e_3(\abs*{q}^{p-2}\abs*{\Phi_1}^2)
         - \frac{p-2}{2}\abs*{q}^{p-2}\abs*{\Phi_1}^2\Trace \underline{\chi}
         - 4s\underline{\omega}\abs*{q}^{p-2}\abs*{\Phi_1}^2
         .
  \end{align*}
  Using \zcref[cap]{lemma:transport:div-e3}, we can write that
  \begin{align*}
    \Divergence\left(\abs*{q}^{p-2}\abs*{\Phi_1}^2e_3\right)
    ={}& e_3\left(\abs*{q}^{p-2}\abs*{\Phi_1^2}\right)
         + \abs*{q}^{p-2}\abs*{\Phi_1}^2\Divergence(e_3)\\
    ={}& 2\abs*{q}^{p-2}\Re\left(\Phi_2\cdot\overline{\Phi_1}\right)
         + \frac{p-2}{2}\abs*{q}^{p-2}\abs*{\Phi_1}^2\Trace \underline{\chi}
         + 4s\underline{\omega}\abs*{q}^{p-2}\abs*{\Phi_1}^2\\
       & + \abs*{q}^{p-2}\abs*{\Phi_1}^2\left(-2\underline{\omega} + \Trace\underline{\chi}\right)\\
    ={}& \abs*{q}^{p-2}\left(
         \frac{p}{2}\Trace \underline{\chi}
         + 2(2s-1)\underline{\omega}
         \right)\abs*{\Phi_1}^2
         + 2\abs*{q}^{p-2}\Re\left(\Phi_2\cdot\overline{\Phi_1}\right).
  \end{align*}
  We can then easily compute that
  \begin{equation*}
    \begin{split}
      &\frac{p}{2}\Trace \underline{\chi}
      + 2(2s-1)\underline{\omega}\\
    ={}&
    \begin{cases}
      -\frac{pr}{\abs*{q}^2}& r\le \frac{r_0}{2},\\
      -\frac{pr\Delta}{\abs*{q}^4} + (2s-1)\partial_r\left(\frac{\Delta}{\abs*{q}^2}\right)
       = - \frac{1}{r}\left(p\frac{\Delta}{\abs*{q}^2} + \frac{2}{3}(2s-1)\Lambda r^2\right)\left(1+O(r^{-1})\right)
                            & r\ge r_0,      \\
      -\frac{1}{r}\left(p + \frac{2}{3}(2s-1)\Lambda r^2\right)\left(1+O(r_0^{-1})\right) & r\in \left[\frac{r_0}{2}, r_0\right].
    \end{cases}
    \end{split}    
  \end{equation*}
  As a result, we see that for all $s>\frac{1}{2}$, for $\delta_{\Horizon}$ sufficiently small and $r_0$ sufficiently large (but independent of $\Lambda$), and for all $p\ge \delta$, there exists some $c_0(p)>0$ such that
  \begin{equation*}
    -\Divergence\left(\abs*{q}^{p-2}\abs*{\Phi_1}^2e_3\right)
    \ge c_0(p)r^{p-3}\abs*{\Phi_1}^2
    - 2\abs*{q}^{p-2}\Phi_2\cdot\overline{\Phi_1}.
  \end{equation*}
Cauchy-Schwarz then gives a constant $C_0(p)>0$ such that 
  \begin{equation*}
    c_0(p)r^{p-3}\abs*{\Phi_1}^2
    \le C_0\abs*{q}^{p-1}\abs*{\Phi_2}^2
    -\Divergence\left(\abs*{q}^{p-2}\abs*{\Phi_1}^2e_3\right)
    ,
  \end{equation*}
  as desired, proving \zcref[noname]{eq:transport:bulk}. The
  integrated estimate in \zcref[noname]{eq:transport:bulk:integrated}
  follows immediately by applying the divergence theorem in
  \zcref{lemma:transport:div-lemma}. 
\end{proof}

In what follows, we will also need to control $\nabla_{\HprVF}\Phi_1$,
for $\Phi_1$ solving
\zcref[noname]{eq:transport:main-lemma:Phi1-Phi2-relation}. We observe
that $\HprVF$ plays a special role as it is the only derivative of
$\mathfrak{q}$ that we have control of globally from the hyperbolic
estimates in \zcref[cap]{prop:combined-estimate:main} due to the loss
of the other derivatives at the trapped set. 
\begin{lemma}
  \label{lemma:transport:nab-3-RHat-rescaled-commutation}
  Let $U\in \realHorkTensor{k}(\Manifold)$. Then we have that
  \begin{equation*}
    \begin{split}
      \left[\nabla_3, \frac{r^2+a^2}{\abs*{q}^2}\HprVF\right]U
    ={}& O(ar^{-1})\nabla U
    + O(r^{-2})\nabla_3 U
         + O(r^{-3})U 
      + O(r^{-1})\nabla_4U
    ,
    \end{split}    
  \end{equation*}
  where we recall that $2(1+\gamma)\HprVF= \frac{\Delta}{r^2+a^2}\lambdaglo^{-1}e_4 - \frac{\abs*{q}^2}{r^2+a^2}\lambdaglo e_3$.
\end{lemma}

\begin{proof}
  We can directly compute that
  \begin{align*}
    \left[\nabla_3, \frac{\Delta}{\abs*{q}^2}\lambdaglo^{-1}\nabla_4 - \lambdaglo\nabla_3\right]
    ={}&\frac{\Delta}{\abs*{q}^2\lambdaglo}[\nabla_3,\nabla_4]
         + e_3\left(\frac{\Delta}{\abs*{q}^2\lambdaglo}\right)\nabla_4
         - e_3(\lambdaglo)\nabla_3.
  \end{align*}
  Then, observing that
  $\eta^{\glo}, \underline{\eta}^{\glo}\in O(ar^{-2})$, using
  \zcref[noname]{eq:commutation-formula:B-applied:e3e4}, and noting that the
  asymptotic behavior of the commutator comes from when
  $(e_4,e_3) = (\eout_4,\eout_3)$, we have that
  \begin{align*}
    \left[\nabla_3,\nabla_4\right]U
    ={}& 2\omega\nabla_{3}U
         + O(ar^{-2})\nabla U
         + O(r_0^{-1}\bOne_{r<r_0})\nabla_{4} U                  
         + O(r^{-3})U\\
    ={}& -\partial_r\left(\frac{\Delta}{\abs*{q}^2}\right)\nabla_{3}U
         + O(r_0^{-1}\bOne_{r<r_0})\nabla_{3}U
         + O(ar^{-2})\nabla U
         + O(r_0^{-1}\bOne_{r<r_0})\nabla_{4} U                  
         + O(r^{-3})U
         .
  \end{align*}
  Realizing also that
  \begin{equation*}
    \frac{\Delta}{\abs*{q}^2\lambdaglo} = O(1), \qquad
    e_3\left(\frac{\Delta}{\abs*{q}^2\lambdaglo}\right) = O(r^{-1} (1+ r^2\Lambda)),
  \end{equation*}
  concludes the proof of \zcref[cap]{lemma:transport:nab-3-RHat-rescaled-commutation}, using also $r^2\La \lesssim 1$.
\end{proof}

We now prove the following commuted transport lemma.
\begin{lemma}
  \label{lemma:transport:RHat-control}
  Suppose that $\Phi_1,\Phi_2\in \realHorkTensor{2}(\mathbb{C})$
  satisfy the relation
  \zcref[noname]{eq:transport:main-lemma:Phi1-Phi2-relation}. Then for
  every $p\ge \delta$, we have that
  \begin{equation}
    \label{eq:transport:RHat-control}
    \begin{split}
      &\int_{\Manifold(\tau_1,\tau_2)}r^{p-3}\left(
      \abs*{\nabla_{\HprVF}\Phi_1}^2
      + \abs*{\Phi_1}^2
      \right)
        + \int_{\partial\Manifold^+(\tau_1,\tau_2)}r^{p-2}\left(\abs*{\nabla_{\HprVF}\Phi_1}^2 + \abs*{\Phi_1}^2\right)\\
      \lesssim{}& \int_{\Manifold(\tau_1,\tau_2)}r^{p-1}\left(\abs*{\nabla_{\HprVF}\Phi_2}^2 + \abs*{\Phi_2}^2\right)
                  + \int_{\Sigma(\tau_1)}r^{p-2}\left(\abs*{\nabla_{\HprVF}\Phi_1}^2 + \abs*{\Phi_1}^2\right)\\
      &  + a^2\int_{\Manifold(\tau_1,\tau_2)}r^{p-3}\abs*{\nabla\Phi_1}^2
        + \int_{\Manifold(\tau_1,\tau_2)}r^{p-3}\abs*{\nabla_4\Phi_1}^2
        .
    \end{split}
  \end{equation}
\end{lemma}


\begin{proof}
  We first commute the transport equation for $\Phi_1$
  with $\HprVF$. Using 
  \zcref[cap]{lemma:transport:nab-3-RHat-rescaled-commutation}, we have that
  \begin{align*}
    \ConformalInvDeriv_{3}\left(\nabla_{\HprVF}\Phi_1\right)
    ={}& \left[\nabla_{3}-2s\underline{\omega}, \nabla_{\HprVF}\right]\Phi_1
         + \nabla_{\HprVF}\Phi_2\\
    ={}& O(r^{-2})\nabla_{3}\Phi_1
         + O(ar^{-1})\nabla \Phi_1
         + O(r^{-1})\Phi_1 
        + O(r^{-1})\nabla_{4} \Phi_1
         + O(1)\nabla_{\HprVF}\Phi_2.
  \end{align*}
  We then apply the transport estimate in \zcref[noname]{eq:transport:bulk:integrated}, deducing that
  \begin{equation*}
    \begin{split}
      &\int_{\Manifold(\tau_1,\tau_2)}r^{p-3}\abs*{\nabla_{\HprVF}\Phi_1}^2
        + \int_{\partial\Manifold^{+}(\tau_1,\tau_2)}r^{p-2}\abs*{\nabla_{\HprVF}\Phi_1}^2\\
      \lesssim{}& \int_{\Manifold(\tau_1,\tau_2)}r^{p-1}\left(         
                  a^2r^{-2}\abs*{\nabla\Phi_1}^2
                  + r^{-2}\abs*{\Phi_1}^2
                  + r^{-2}\abs*{\nabla_{4} \Phi_1}^2
                  + r^{-4}\abs*{\nabla_{3}\Phi_1}^2
                  + \abs*{\nabla_{\HprVF}\Phi_2}^2
                  \right) + \int_{\Sigma(\tau_1)}r^{p-2}\abs*{\nabla_{\HprVF}\Phi_1}^2.
    \end{split}
  \end{equation*}
  Since $\Phi_1,\Phi_2$ satisfy
  \zcref[noname]{eq:transport:main-lemma:Phi1-Phi2-relation}, we can
  also combine the above with a large multiple of 
  \zcref[noname]{eq:transport:bulk:integrated} to see that
  \begin{equation*}
    \begin{split}
      &\int_{\Manifold(\tau_1,\tau_2)}r^{p-3}\left( \abs*{\nabla_{\HprVF}\Phi_1}^2  + C\abs*{\Phi_1}^2\right)
        + \int_{\partial\Manifold^{+}(\tau_1,\tau_2)}r^{p-2}\left( \abs*{\nabla_{\HprVF}\Phi_1}^2 + C\abs*{\Phi_1}^2 \right)
      \\
      \lesssim{}& \int_{\Manifold(\tau_1,\tau_2)}r^{p-1}\left(         
                  a^2r^{-2}\abs*{\nabla\Phi_1}^2
                  + r^{-2}\abs*{\Phi_1}^2
                  + r^{-2}\abs*{\nabla_{4} \Phi_1}^2
                  + r^{-4}\abs*{\Phi_2}^2
                  + \abs*{\nabla_{\HprVF}\Phi_2}^2
                  \right)\\
      & + \int_{\Sigma(\tau_1)}r^{p-2}\abs*{\nabla_{\HprVF}\Phi_1}^2
        + C\int_{\mathcal{M}(\tau_1,\tau_2)}r^{p-1}\abs*{\Phi_2}^2
        + C\int_{\Sigma(\tau_1)}r^{p-2}\abs*{\Phi_1}^2
        .
    \end{split}
  \end{equation*}
  Thus we see that for $C, r_0$ sufficiently large,
  \begin{equation*}
    \begin{split}
      &\int_{\Manifold(\tau_1,\tau_2)}r^{p-3}\left( \abs*{\nabla_{\HprVF}\Phi_1}^2  + \abs*{\Phi_1}^2\right)
        + \int_{\partial\Manifold^{+}(\tau_1,\tau_2)}r^{p-2}\left( \abs*{\nabla_{\HprVF}\Phi_1}^2 + \abs*{\Phi_1}^2 \right)
      \\
      \lesssim{}& \int_{\Manifold(\tau_1,\tau_2)}r^{p-1}\left(         
                  a^2r^{-2}\abs*{\nabla\Phi_1}^2
                  + r^{-2}\abs*{\nabla_4 \Phi_1}^2
                  + \abs*{\Phi_2}^2
                  + \abs*{\nabla_{\HprVF}\Phi_2}^2
                  \right) + \int_{\Sigma(\tau_1)}r^{p-2}\left( \abs*{\nabla_{\HprVF}\Phi_1}^2 + \abs*{\Phi_1}^2 \right)
        ,
    \end{split}
  \end{equation*}
  as desired, concluding the proof of \zcref{lemma:transport:RHat-control}.
\end{proof}

In practice, we would also like to recover the $\nabla_4$ derivative
of $\Phi_1$. In the particular scenarios that we are considering, it
will not be possible to achieve full control of $\nabla_4$ due to the
loss of derivatives at trapping. However, observe that away from
trapping there is no loss of derivatives, and thus we can prove the
following general lemma.
\begin{lemma}
  \label{lemma:general-transport-estimate}
  Let $\Phi_1,\Phi_2\in \realHorkTensor{2}(\Complex)$ satisfy the
  differential relation
  \zcref[noname]{eq:transport:main-lemma:Phi1-Phi2-relation}. Moreover,
  let $\breve{\chi}=\breve{\chi}(r)\in C^{\infty}$ be a smooth cutoff
  function such that
  \begin{equation*}
    \breve{\chi}(r) =
    \begin{cases}
      0 & r-3M < \delta_{\operatorname{trap}},\\
      1 & r > r_0.
    \end{cases}
  \end{equation*}
  Then for all $p\ge \delta$, we have
  \begin{equation*}
    \begin{split}
      &\int_{\mathcal{M}(\tau_1,\tau_2)}r^{p-3}\left(
      r^2\abs*{\nabla_{3}\Phi_1} + r^2\abs*{\nabla_{4}\Phi_1}^2 + \abs*{\Phi_1}^2
      \right)\\
      & + \int_{\partial \mathcal{M}^+(\tau_1,\tau_2)}r^{p-2}\left(
        r^2\breve{\chi}^2\abs*{\nabla_{4}\Phi_1}^2
        + \abs*{\nabla_{\HprVF}\Phi_1}^2
        + \abs*{\Phi_1}^2
        \right)\\
      \lesssim{}& \int_{\mathcal{M}(\tau_1,\tau_2)}r^{p-1}\left(
                  r^2\breve{\chi}^2\abs*{\nabla_{4}\Phi_2}^2
                  + \abs*{\nabla_{\HprVF}\Phi_2}^2
                  + \abs*{\Phi_2}^2
                  \right)\\
      & + \int_{\Sigma(\tau_1)} r^{p-2}\left(
        r^2\breve{\chi}^2\abs*{\nabla_{4}\Phi_1}^2
        + \abs*{\nabla_{\HprVF}\Phi_1}^2
                  + \abs*{\Phi_1}^2
        \right)
       + a^2\int_{\mathcal{M}(\tau_1,\tau_2)}r^{p-1}\abs*{\nabla\Phi_1}^2,
    \end{split}
  \end{equation*}
  where
  \begin{equation*}
    \partial^+\mathcal{M}^+(\tau_1,\tau_2)
    = \mathcal{A}(\tau_1,\tau_2)\cup \Sigma(\tau_2)\cup \SigmaStar(\tau_1,\tau_2)
  \end{equation*}
  denotes the future boundary of $\mathcal{M}(\tau_1,\tau_2)$. 
\end{lemma}
\begin{proof}
  The main goal of the lemma is to control
  $\breve{\chi}\nabla_4\Phi_1$. To this end, we recall from \zcref[cap]{lemma:commutation-formula:conformally-invariant-derivaties} that
  \begin{align*}
    [\ConformalInvDeriv_4, \ConformalInvDeriv_3]
    ={}& 2\left(\eta-\etaBar\right)\cdot \ConformalInvDeriv U
    + 2 s\left(
      \rho- \frac{\Lambda}{3} -\eta\cdot \etaBar
    \right)U
         + 4\ImagUnit \left(-\LeftDual{\rho} + \eta\wedge \etaBar\right)U\\
    ={}& O(a r^{-1})\ConformalInvDeriv U
         + O(\Lambda + r^{-3})U
         .
  \end{align*}
  Thus we have that
  \begin{align*}
    \ConformalInvDeriv_{3}\left( \breve{\chi}\ConformalInvDeriv_{4}\Phi_1 \right)
    ={}& \breve{\chi}\ConformalInvDeriv_{3}\ConformalInvDeriv_{4}\Phi_1
         + \partial_r\breve{\chi}e_3(r)\ConformalInvDeriv_{4}\Phi_1\\
    ={}& \breve{\chi}\ConformalInvDeriv_{4}\Phi_2
         + \breve{\chi}\left[  \ConformalInvDeriv_{3}, \ConformalInvDeriv_{4} \right]\Phi_1
         + \partial_r\breve{\chi}e_3(r)\ConformalInvDeriv_{4}\Phi_1\\
    ={}& \breve{\chi}\ConformalInvDeriv_{4}\Phi_2
         + O(\Lambda + r^{-3})\Phi_1
         + O(ar^{-1})\ConformalInvDeriv\Phi_1
         + \partial_r\breve{\chi}e_3(r)\ConformalInvDeriv_{4}\Phi_1.
  \end{align*}
  We can now use \zcref[cap]{lemma:transport:bulk} to write that
  \begin{equation*}
    \begin{split}
      &\int_{\mathcal{M}(\tau_1,\tau_2)}r^{p-3}\left(r^2\left( \breve{\chi}\ConformalInvDeriv_{4}\Phi_1 \right)^2\right)
      + \int_{\partial \mathcal{M}(\tau_1,\tau_2)}r^{p-2}\left(r^2\left( \breve{\chi}\ConformalInvDeriv_{4}\Phi_1 \right)^2\right)\\
      \lesssim{}& \int_{\mathcal{M}(\tau_1,\tau_2)}r^{p-1}\left(r^2\abs*{\widetilde{\Phi}_2}^2 \right)
                  + \int_{\Sigma(\tau_1)}r^{p-2}\left(r^2\left( \breve{\chi}\ConformalInvDeriv_{4}\Phi_1 \right)^2\right),
    \end{split}        
  \end{equation*}
  where
  \begin{equation*}
    \widetilde{\Phi}_2 = \breve{\chi}\ConformalInvDeriv_{4}\Phi_2
    + O(\Lambda + r^{-3})\Phi_1
    + O(ar^{-1})\ConformalInvDeriv\Phi_1
    + \partial_r\breve{\chi}e_3(r)\ConformalInvDeriv_{4}\Phi_1.
  \end{equation*}
  We can use the null decomposition of $\HprVF$, the fact that
  $\partial_r\breve{\chi}$ is compactly supported away from the
  horizons, and 
  \zcref[noname]{eq:transport:main-lemma:Phi1-Phi2-relation} to write
  that $\partial_r\breve{\chi}e_3(r)\nabla_4\Phi_1
    = \partial_r\breve{\chi}O(1)\nabla_{\HprVF}\Phi_1 + \partial_r\breve{\chi}O(1)\Phi_2$. As a result, we have that 
  \begin{equation}
    \label{eq:transport:RHat-control:Phi2-tilde-aux}
    \widetilde{\Phi}_2
    = \breve{\chi}\ConformalInvDeriv_{4}\Phi_2
    + O(\Lambda + r^{-3})\Phi_1
    + O(ar^{-1})\ConformalInvDeriv\Phi_1
    + \partial_r\breve{\chi}O(1)\ConformalInvDeriv_{\HprVF}\Phi_1
    + \partial_r\breve{\chi}O(1)\Phi_2
    .
  \end{equation}
  We can therefore write that 
  \begin{equation*}
    \begin{split}
      \abs*{\widetilde{\Phi}_2}^2
      \lesssim \abs*{\ConformalInvDeriv_{4}\Phi_2}^2
      + r^{-4}\abs*{\Phi_1}^2
      + ar^{-2}\abs*{\ConformalInvDeriv\Phi_1}^2
      + \abs*{\ConformalInvDeriv_{\HprVF}\Phi_1}^2
      + \abs*{\Phi_2}^2,
    \end{split}
  \end{equation*}
  where we used the fact that $\Lambda = O(r^{-2})$.
  As a result, we now have that
  \begin{equation*}
    \begin{split}
      &\int_{\mathcal{M}(\tau_1,\tau_2)}r^{p-3}\left(r^2 \breve{\chi}^2\abs*{\ConformalInvDeriv_{4}\Phi_1}^2\right)
      + \int_{\partial \mathcal{M}(\tau_1,\tau_2)}r^{p-2}\left(r^2\breve{\chi}^2\abs*{\ConformalInvDeriv_{4}\Phi_1}^2\right)\\
      \lesssim{}& \int_{\mathcal{M}(\tau_1,\tau_2)}r^{p-1}\left(r^2\abs*{\ConformalInvDeriv_{4}\Phi_2}^2                  
                  + \abs*{\Phi_2}^2
                  \right)
                  + \int_{\mathcal{M}(\tau_1,\tau_2)}r^{p-3}\left( \abs*{\nabla_{\HprVF}\Phi_1}^2 + \abs*{\Phi_1}^2 \right)\\
                 & + \int_{\mathcal{M}(\tau_1,\tau_2)}ar^{p-1}\abs*{\ConformalInvDeriv\Phi_1}^2
                  + \int_{\Sigma(\tau_1)}r^{p-2}\left(r^2 \breve{\chi}^2\abs*{\ConformalInvDeriv_{4}\Phi_1}^2\right),
    \end{split}        
  \end{equation*}
  where we remark that we used the fact that $\partial_r\breve{\chi}$
  has compact support to control the $\partial_r\breve{\chi}\Phi_2$
  term and the $\nabla_{\HprVF}\Phi_1$ term on the \RHS{} coming from
  \zcref[noname]{eq:transport:RHat-control:Phi2-tilde-aux}.  We can
  now apply \zcref[cap]{lemma:transport:RHat-control}. This gives us
  extra control on the left-hand side as well as allowing us to
  control the
  $\int_{\mathcal{M}(\tau_1,\tau_2)}r^{p-3}\abs*{\Phi_1}^2$ term on
  the \RHS, resulting in the following estimate.
  \begin{equation*}
    \begin{split}
      &\int_{\mathcal{M}(\tau_1,\tau_2)}r^{p-3}\left(
        r^2\breve{\chi}^2\abs*{\nabla_{4}\Phi_1}^2
        + \abs*{\nabla_{\HprVF}\Phi_1}^2
        + \abs*{\Phi_1}^2
        \right)\\
      &  + \int_{\partial \mathcal{M}(\tau_1,\tau_2)}r^{p-2}\left(
        r^2\breve{\chi}^2\abs*{\nabla_{4}\Phi_1}^2
        + \abs*{\nabla_{\HprVF}\Phi_1}^2
        + \abs*{\Phi_1}^2
        \right)\\
      \lesssim{}& \int_{\mathcal{M}(\tau_1,\tau_2)}r^{p-1}\left(
                  r^2\abs*{\ConformalInvDeriv_{4}\Phi_2}^2
                  + \abs*{\nabla_{\HprVF}\Phi_2}^2
                  + \abs*{\Phi_2}^2
                  \right)
      \\
      & + \int_{\Sigma(\tau_1)}r^{p-2}\left(
        r^2 \breve{\chi}^2\abs*{\ConformalInvDeriv_{4}\Phi_1}^2
        + \abs*{\nabla_{\HprVF}\Phi_1}^2
        + \abs*{\Phi_1}^2
        \right)\\
      & + a^2\int_{\mathcal{M}(\tau_1,\tau_2)}r^{p-1}\abs*{\nabla\Phi_1}^2
      + \int_{\mathcal{M}(\tau_1,\tau_2)}r^{p-3}\abs*{\nabla_4\Phi_1}^2
        .
    \end{split}
  \end{equation*}
  Taking $r_0$ sufficiently large then allows us to absorb the last
  term on the \RHS{} into the \LHS{} and close the proof of
  \zcref[cap]{lemma:general-transport-estimate}.
\end{proof}

\subsection{Estimate for \texorpdfstring{$\BEFNormAux{p}[A](\tau_1,\tau_2)$}{the auxiliary BEF norm of A}}

We now recover an estimate for the
$\BEFNormAux{p}[A](\tau_1,\tau_2)$ norm of $A$. Recall from its
definition in \zcref[noname]{eq:BEF-A-norm:def} that the
$\BEFNormAux{p}[A](\tau_1,\tau_2)$ norm does not control all
second derivatives of $A$ (lacking notably in the angular
directions). We will recover the angular derivatives in the next section. 
\begin{proposition}
  \label{prop:BEF-A:tranport-estimates}
  The following estimates hold true for all $\delta\le p \le 2-\delta$,
  \begin{equation}
    \label{eq:BEF-A:tranport-estimates}
    \BEFNormAux{p}[A](\tau_1,\tau_2)
    \lesssim \BulkNormWeighted{p}[\mathfrak{q}](\tau_1,\tau_2)
    + \EnergyFlux_p[A](\tau_1)
    + a^2\int_{\mathcal{M}(\tau_1,\tau_2)}r^{p+3}\left(
      r^2\abs*{\nabla\nabla_3A}^2
      + \abs*{\nabla A}^2
    \right).
  \end{equation}
\end{proposition}
\begin{proof}
  Recall from \zcref[cap]{coro:Teukolsky-wave-transport} that $\ConformalInvDeriv_3\Psi = O(r^{-2})\mathfrak{q}$. Applying \zcref[cap]{lemma:general-transport-estimate} with $(\Phi_1, \Phi_2) = (\Psi, O(r^{-2})\mathfrak{q})$, we have then that for $p'\ge \delta$,
  \begin{equation*}
    \begin{split}
      &\int_{\mathcal{M}(\tau_1,\tau_2)}r^{p'-3}\left(
        r^2\abs*{\nabla_{3}\Psi}^2
        + r^2\abs*{\nabla_{4}\Psi}^2
        + \abs*{\Psi}^2
        \right)
        + \int_{\partial \mathcal{M}^+(\tau_1,\tau_2)}r^{p'-2}\left(r^2\breve{\chi}^2\abs*{\nabla_4\Psi}^2 + \abs*{\nabla_{\HprVF}\Psi}^2 + \abs*{\Psi}^2\right)\\
      \lesssim{}& \int_{\mathcal{M}(\tau_1,\tau_2)}r^{p'-5}\left(
                  r^2\breve{\chi}^2\abs*{\nabla_4\mathfrak{q}}^2
                  + \abs*{\nabla_{\HprVF}\mathfrak{q}}^2
                  + \abs*{\mathfrak{q}}^2
                  \right)
                  + \int_{\Sigma(\tau_1)}r^{p'-2}\left(
                  r^2\breve{\chi}^2\abs*{\nabla_4\Psi}^2
                  + \abs*{\nabla_{\HprVF}\Psi}^2
                  + \abs*{\Psi}^2
                  \right)\\
      & + a^2\int_{\mathcal{M}(\tau_1,\tau_2)}r^{p'-1}\abs*{\nabla\Psi}^2.
    \end{split}
  \end{equation*}
  Choosing $p'=p+2$, we obtain for $\delta\le p\le 2-\delta$ that
  \begin{equation*}
    \begin{split}
      &\int_{\mathcal{M}(\tau_1,\tau_2)}r^{p-1}\left(
        r^2\abs*{\nabla_{3}\Psi}^2
        + r^2\abs*{\nabla_{4}\Psi}^2
        + \abs*{\Psi}^2
        \right)
        + \int_{\partial \mathcal{M}^+(\tau_1,\tau_2)}r^p\left(r^2\breve{\chi}^2\abs*{\nabla_4\Psi}^2 + \abs*{\nabla_{\HprVF}\Psi}^2 + \abs*{\Psi}^2\right)\\
      \lesssim{}& \int_{\mathcal{M}(\tau_1,\tau_2)}r^{p-3}\left(
                  r^2\breve{\chi}^2\abs*{\nabla_4\mathfrak{q}}^2
                  + \abs*{\nabla_{\HprVF}\mathfrak{q}}^2
                  + \abs*{\mathfrak{q}}^2
                  \right)
                  + \int_{\Sigma(\tau_1)}r^p\left(
                  r^2\breve{\chi}^2\abs*{\nabla_4\Psi}^2
                  + \abs*{\nabla_{\HprVF}\Psi}^2
                  + \abs*{\Psi}^2
                  \right)\\
      & + a^2\int_{\mathcal{M}(\tau_1,\tau_2)}r^{p+1}\abs*{\nabla\Psi}^2.
    \end{split}
  \end{equation*}
  In view of the definitions of the norms $\BulkNormWeighted{p}[A]$ in
  \zcref[noname]{eq:A-norms} and $\BulkNormWeighted{p}[\mathfrak{q}]$ in
  \zcref[noname]{eq:rp:weighted-bulk:def}, we infer that for $\delta\le p\le 2-\delta$,
  \begin{equation}
    \label{eq:transport:BEF-A:aux1}
    \begin{split}
      &\int_{\mathcal{M}(\tau_1,\tau_2)}r^{p-1}\left(
        r^2\abs*{\nabla_{3}\Psi}^2
        + r^2\abs*{\nabla_{4}\Psi}^2
        + \abs*{\Psi}^2
        \right)
        + \int_{\partial \mathcal{M}^+(\tau_1,\tau_2)}r^p\left(r^2\breve{\chi}^2\abs*{\nabla_4\Psi}^2 + \abs*{\nabla_{\HprVF}\Psi}^2 + \abs*{\Psi}^2\right)\\
      \lesssim{}& \BulkNormWeighted{p}[\mathfrak{q}](\tau_1,\tau_2)
                  + \int_{\Sigma(\tau_1)}r^p\left(
                  r^2\breve{\chi}^2\abs*{\nabla_4\Psi}^2
                  + \abs*{\nabla_{\HprVF}\Psi}^2
                  + \abs*{\Psi}^2
                  \right)
                  + a^2\int_{\mathcal{M}(\tau_1,\tau_2)}r^{p+1}\abs*{\nabla\Psi}^2.
    \end{split}
  \end{equation}
  Then applying \zcref[cap]{lemma:general-transport-estimate} with $(\Phi_1,\Phi_2)=\left(
      \frac{\overline{q}^4}{r^2}A,
      \Psi
    \right)$ we obtain that for $p\ge \delta$,
  \begin{equation}
    \label{eq:transport:BEF-A:aux2}
    \begin{split}
      &\int_{\mathcal{M}(\tau_1,\tau_2)}r^{p+1}\left(
        r^2\abs*{\nabla_{3}A}^2
        + r^2\abs*{\nabla_{4}A}^2
        + \abs*{A}^2
        \right)
        + \int_{\partial \mathcal{M}^+(\tau_1,\tau_2)}r^{p+2}\left(r^2\breve{\chi}^2\abs*{\nabla_4A}^2 + \abs*{\nabla_{\HprVF}A}^2 + \abs*{A}^2\right)\\
      \lesssim{}& \int_{\Sigma(\tau_1)}r^{p-1}\left(
                  r^2\breve{\chi}^2\abs*{\nabla_4\Psi}^2
                  + \abs*{\nabla_{\HprVF}\Psi}^2
                  + \abs*{\Psi}^2
                  \right)
                  + \int_{\Sigma(\tau_1)}r^{p+2}\left(
                  r^2\breve{\chi}^2\abs*{\nabla_4A}^2
                  + \abs*{\nabla_{\HprVF}A}^2
                  + \abs*{A}^2
                  \right)\\
      & + a^2\int_{\mathcal{M}(\tau_1,\tau_2)}r^{p+3}\abs*{\nabla A}^2.        
    \end{split}
  \end{equation}
  Combining \zcref[noname]{eq:transport:BEF-A:aux1} and
  \zcref[noname]{eq:transport:BEF-A:aux2}, we have that
  \begin{equation}
    \label{eq:transport:BEF-A:aux3}
    \begin{split}
      & \int_{\mathcal{M}(\tau_1,\tau_2)}r^{p+1}\left(
        r^2\abs*{\nabla_{3}A}^2
        + r^2\abs*{\nabla_{4}A}^2
        + \abs*{A}^2
        \right)\\
      &  + \int_{\mathcal{M}(\tau_1,\tau_2)}r^{p-1}\left(
        r^2\abs*{\nabla_{3}\Psi}^2
        + r^2\abs*{\nabla_{4}\Psi}^2
        + \abs*{\Psi}^2
        \right)\\
      &  + \int_{\partial \mathcal{M}^+(\tau_1,\tau_2)}r^{p+2}\left(r^2\breve{\chi}^2\abs*{\nabla_4A}^2 + \abs*{\nabla_{\HprVF}A}^2 + \abs*{A}^2\right)\\
      & + \int_{\partial \mathcal{M}^+(\tau_1,\tau_2)}r^p\left(r^2\breve{\chi}^2\abs*{\nabla_4\Psi}^2 + \abs*{\nabla_{\HprVF}\Psi}^2 + \abs*{\Psi}^2\right)\\
      \lesssim{}& \BulkNormWeighted{p}[\mathfrak{q}](\tau_1,\tau_2)
                  + \int_{\Sigma(\tau_1)}r^p\left(
                  r^2\breve{\chi}^2\abs*{\nabla_4\Psi}^2
                  + \abs*{\nabla_{\HprVF}\Psi}^2
                  + \abs*{\Psi}^2
                  \right)\\
      &  + \int_{\Sigma(\tau_1)}r^{p+2}\left(
        r^2\breve{\chi}^2\abs*{\nabla_4A}^2
        + \abs*{\nabla_{\HprVF}A}^2
        + \abs*{A}^2
        \right)
        + a^2\int_{\mathcal{M}(\tau_1,\tau_2)}\left( 
        r^{p+1}\abs*{\nabla\Psi}^2
        + r^{p+3}\abs*{\nabla A}^2
        \right).
    \end{split}
  \end{equation}
  Next, observe that from \zcref[noname]{eq:Teukolsky:Psi:def} that $ \Psi = \frac{\overline{q}^4}{r^2}\nabla_3A + O(r)A$. This directly implies that
  \begin{equation}
    \label{eq:transport:BEF-A:nab-A-in-terms-of-Psi:1}
    \begin{split}
      \abs*{\nabla_3A} &\lesssim r^{-2}\abs*{\Psi} + r^{-1}\abs*{A},\\
      \abs*{\nabla_4\nabla_3A} &\lesssim r^{-2}\abs*{\nabla_4\Psi}
                                 + r^{-3}\abs*{\Psi}
                                 + r^{-1}\abs*{\nabla_4A}
                                 + r^{-2}\abs*{A},\\
      \abs*{\nabla_{\HprVF}\nabla_3A}
                       & \lesssim r^{-2}\abs*{\nabla_{\HprVF}\Psi}
                         + r^{-2}\abs*{\Psi}
                         + r^{-1}\abs*{\nabla_{\HprVF}A}
                         + r^{-1}\abs*{A}.
    \end{split}
  \end{equation}
  Similarly, the definition of $\mathfrak{q}$
  from~\zcref[noname]{definition of q frak} gives $\ConformalInvDeriv_3^2A
    = \frac{1}{q\overline{q}^3}\mathfrak{q}
    + O(r^{-1})\ConformalInvDeriv_3A
    + O(r^{-2})A$, which directly implies that
  \begin{equation}
    \label{eq:transport:BEF-A:nab-A-in-terms-of-Psi:2}
    \abs*{\nabla_3\ConformalInvDeriv_3A}
    \lesssim r^{-4}\abs*{\mathfrak{q}}
    + r^{-1}\abs*{\nabla_3A}
    + r^{-2}\abs*{A}. 
  \end{equation}
  We can use \zcref[noname]{eq:transport:BEF-A:nab-A-in-terms-of-Psi:1}, and
  \zcref[noname]{eq:transport:BEF-A:nab-A-in-terms-of-Psi:2} to rewrite
  \zcref[noname]{eq:transport:BEF-A:aux3} as
  \begin{equation*}
    \begin{split}
      &\int_{\mathcal{M}(\tau_1,\tau_2)}
    r^{p+1}\left(
      r^4\abs*{\nabla_3\ConformalInvDeriv_3A}^2
      + r^4\abs*{\nabla_4\nabla_3A}^2
      +r^2\abs*{\nabla_3A}^2
      + r^2\abs*{\nabla_4A}^2
      + \abs*{A}^2
    \right)\\
    & + \int_{\partial\Manifold^+(\tau_1,\tau_2)}r^{p+2}\left(
      r^4\breve{\chi}^2\abs*{\nabla_4\nabla_3A}^2
      + r^2\abs*{\nabla_{\HprVF}\nabla_3A}^2
      +r^2\abs*{\nabla_3A}^2
      + r^2\abs*{\nabla_4A}^2
      + \abs*{A}^2
    \right)\\
    \lesssim{}& \BulkNormWeighted{p}[\mathfrak{q}](\tau_1,\tau_2)
    + \int_{\Sigma(\tau_1)}r^{p+2}\left(
      r^4\breve{\chi}^2\abs*{\nabla_4\nabla_3A}^2
      + r^2\abs*{\nabla_{\HprVF}\nabla_3A}^2
      + r^2\abs*{\nabla_3A}^2
      + r^2\abs*{\nabla_4A}^2
      + \abs*{A}^2
    \right)\\
    & + a^2\int_{\mathcal{M}(\tau_1,\tau_2)}\left( r^{p+1}\abs*{\nabla\Psi}^2 + r^{p+3}\abs*{\nabla A}^2 \right).
    \end{split}    
  \end{equation*}
  Then, the definition of the $\WeightedBEFNorm{p}[A]$ norm in
   \zcref[noname]{eq:BEF-A-norm:def} and of $\EnergyFlux_p[A](\tau_1)$ in \zcref{eq:bulk-A-aux-norm:def} imply that
  \begin{equation*}
    \BEFNormAux{p}[A]\lesssim \BulkNormWeighted{p}[\mathfrak{q}](\tau_1,\tau_2)
    + \EnergyFlux_p[A](\tau_1)
    + a^2\int_{\mathcal{M}(\tau_1,\tau_2)}r^{p+3}\left(r^2\abs*{\nabla\nabla_3A}^2 + \abs*{\nabla A}^2\right),
  \end{equation*}
  as desired, concluding the proof of \zcref[noname]{eq:BEF-A:tranport-estimates}.
\end{proof}

\subsection{Estimates for the angular derivatives of \texorpdfstring{$A$}{A}}

We now recover the angular derivatives of $A$ which are missing in the
$\WeightedBEFNorm{p}[A]$ norm. We will use the following lemma to control the angular derivatives.
\begin{lemma}
  \label{lemma:transport:angular-derivatives-of-A:main-aux-lemma}
  For any sphere $S\subset \mathcal{M}$, and
  $U\in \realHorkTensor{2}(\Complex)$
  \begin{equation}
    \label{eq:transport:angular-derivatives-of-A:main-aux-lemma}
    \int_S\left(
      \abs*{\nabla U}^2
      + r^{-2}\abs*{U}^2
    \right)
    \lesssim \abs*{\int_S\overline{U}\cdot\ConformalComplexDeriv\SymTracelessTensorProd\left( \overline{\ConformalComplexDeriv}\cdot U \right)}
    + O(a^2)\int_Sr^{-2}\abs*{(\nabla_3,\nabla_4)U}^2.
  \end{equation}
\end{lemma}
\begin{proof}
  See the proof of Lemma 11.4.12 in
  \cite{giorgiWaveEquationsEstimates2024}. The proof is identical
  after the application of $\Lambda=O(r^{-2})$, using in
  particular \zcref{eq:K:asymptotic-behavior}.
\end{proof}

We first control
$\ConformalInvDeriv\SymTracelessTensorProd\left(
  \overline{\ConformalComplexDeriv}\cdot A \right)$ and
$\ConformalInvDeriv\SymTracelessTensorProd\left(
  \overline{\ConformalInvDeriv}\cdot\nabla_3A \right)$ by using the
Teukolsky equation itself in \zcref[noname]{eq:teuk}.
\begin{lemma}
  \label{lemma:transport:D-hot-D-A-control}
  Let $A\in \realHorkTensor{2}$ be a 2-conformally invariant tensor
  solving $\TeukOp[A]=0$, with $\TeukOp$ as defined in \zcref[noname]{eq:Teuk:A:def}.
  For $\delta\le p \le 2-\delta$, we have the following estimates 
  \begin{equation}
    \label{eq:transport:D-hot-D-A-control}
    \begin{split}
      \int_{\mathcal{M}(\tau_1,\tau_2)}r^{p+5}\abs*{\ConformalComplexDeriv\SymTracelessTensorProd\left(
      \overline{\ConformalComplexDeriv}\cdot A \right)}^2
      \lesssim{}& a^2\int_{\mathcal{M}(\tau_1,\tau_2)}r^{p+1}\abs*{\nabla A}^2
                  + \BulkNormWeighted{p}[A](\tau_1,\tau_2),\\
      \int_{\Sigma(\tau)}r^{p+6}\abs*{\ConformalComplexDeriv\SymTracelessTensorProd\left(
      \overline{\ConformalComplexDeriv}\cdot A \right)}^2
      \lesssim{}& \int_{\Sigma(\tau)}r^{p+2}\abs*{\nabla A}^2
                  + \EnergyFlux_p[A](\tau),
    \end{split}
  \end{equation}
  as well as 
  \begin{equation}
    \label{eq:transport:D-hot-D-A-control:D-hot-D-nab3-A}
    \begin{split}
            \int_{\mathcal{M}(\tau_1,\tau_2)}r^{p+7}\abs*{\ConformalComplexDeriv\SymTracelessTensorProd\left( 
  \overline{\ConformalComplexDeriv}\cdot\nabla_3A \right)}^2
      \lesssim{}& a^2\int_{\mathcal{M}(\tau_1,\tau_2)}r^{p+3}\left( \abs*{\nabla\nabla_3A}+\abs*{\nabla A}^2 \right)\\
                & + \BulkNormWeighted{p}[\mathfrak{q}](\tau_1,\tau_2)
                  + \BulkNormWeighted{p}[A](\tau_1,\tau_2),\\
      \int_{\Sigma(\tau)}r^{p+8}\abs*{\ConformalComplexDeriv\SymTracelessTensorProd\left(\overline{\ConformalComplexDeriv}\cdot\nabla_3A \right)}^2
      \lesssim{}& \int_{\Sigma(\tau)}r^{p+4}\left( \abs*{\nabla\nabla_3A}+\abs*{\nabla A}^2 \right)\\
                & + \EnergyFlux_p[\mathfrak{q}](\tau)
                  + \EnergyFlux_p[A](\tau),
    \end{split}
  \end{equation}
\end{lemma}
\begin{proof}
  From \zcref[cap]{prop:Teuk:Teuk-eq-for-A}, we have that
  \begin{equation}
    \label{eq:transport:D-hot-D-A-control:reduced-Teuk}
    \begin{split}
      \frac{1}{4}\ConformalComplexDeriv\SymTracelessTensorProd\left(\overline{\ConformalComplexDeriv}\cdot A\right)
      ={}& \ConformalInvDeriv_4\ConformalInvDeriv_3A
           + \left(\frac{5}{r} + O(ar^{-2})\right)\ConformalInvDeriv_3A
           + O(r^{-1})\ConformalInvDeriv_4A\\
         & + O(ar^{-2})\ConformalInvDeriv A
           + O(r^{-2}) A,
    \end{split}
  \end{equation}
  where we used the fact that $\Lambda = O(r^{-2})$. We can then
  integrate directly to see that
  \begin{equation*}
    \begin{split}
      \int_{\mathcal{M}(\tau_1,\tau_2)}r^{p+5}\abs*{\ConformalComplexDeriv\SymTracelessTensorProd\left(
      \overline{\ConformalComplexDeriv}\cdot A \right)}^2
      \lesssim{}& a^2\int_{\mathcal{M}(\tau_1,\tau_2)}r^{p+1}\abs*{\nabla A}^2
                  + \BulkNormWeighted{p}[A](\tau_1,\tau_2),\\
      \int_{\Sigma(\tau)}r^{p+6}\abs*{\ConformalComplexDeriv\SymTracelessTensorProd\left(
      \overline{\ConformalComplexDeriv}\cdot A \right)}^2
      \lesssim{}& \int_{\Sigma(\tau)}r^{p+2}\abs*{\nabla A}^2
                  + \EnergyFlux_p[A](\tau),
    \end{split}
  \end{equation*}
  as desired. To prove the control in
  \zcref[noname]{eq:transport:D-hot-D-A-control:D-hot-D-nab3-A}, we
  can compute that
  \begin{align*}
    &\nabla_3\left(
      \left(\frac{5}{r} + O(ar^{-2})\right)\ConformalInvDeriv_3A
      + O(r^{-1})\ConformalInvDeriv_4A
      + O(ar^{-2})\ConformalInvDeriv A
      + O(r^{-2}) A
      \right)\\
    ={}& \left(\frac{5}{r} + O(ar^{-2})\right)\ConformalInvDeriv_3\ConformalInvDeriv_3A
         + O(r^{-1})\ConformalInvDeriv_4\ConformalInvDeriv_3A
         + O(r^{-1})\left[\ConformalInvDeriv_3,\ConformalInvDeriv_4\right]\\
    & + O(ar^{-2})\ConformalInvDeriv\ConformalInvDeriv_3 A
      + O(ar^{-2})\left[ \ConformalInvDeriv_3,\ConformalInvDeriv \right] A\\
    &  + O(r^{-2})\ConformalInvDeriv_3A
      + O(r^{-2})\ConformalInvDeriv_4A
      + O(ar^{-3})\ConformalInvDeriv A
      + O(r^{-3}) A.
  \end{align*}
  Then using
  \zcref[cap]{lemma:commutation-formula:conformally-invariant-derivaties},
  we see that
  \begin{align*}
    &\nabla_3\left(
      \left(\frac{5}{r} + O(ar^{-2})\right)\ConformalInvDeriv_3A
      + O(r^{-1})\ConformalInvDeriv_4A
      + O(ar^{-2})\ConformalInvDeriv A
      + O(r^{-2}) A
      \right)\\
    ={}& \left(\frac{5}{r} + O(ar^{-2})\right)\ConformalInvDeriv_3\ConformalInvDeriv_3A
         + O(r^{-1})\ConformalInvDeriv_4\ConformalInvDeriv_3A
     + O(ar^{-2})\ConformalInvDeriv\ConformalInvDeriv_3 A\\
    &  + O(r^{-2})\ConformalInvDeriv_3A
      + O(r^{-2})\ConformalInvDeriv_4A
      + O(ar^{-3})\ConformalInvDeriv A
      + O(r^{-3}) A.
  \end{align*}
  Using the definition of $\mathfrak{q}$ from
  \zcref[noname]{definition of q frak}, we then have that
  \begin{equation}
    \label{eq:transport:D-hot-D-A-control:nab3-2-A-and-qfrak-relation}
    \ConformalInvDeriv_3^2A = \left( \frac{1}{r^4} + O(ar^{-5}) \right)\mathfrak{q}
    + O(r^{-1})\ConformalInvDeriv_3 A
    + O(r^{-2})A. 
  \end{equation}
  Thus, we have that
  \begin{equation}
    \label{eq:transport:D-hot-D-A-control:D-hot-D-nab3-A:aux1}
    \begin{split}
      &\nabla_3\left(
      \left(\frac{5}{r} + O(ar^{-2})\right)\ConformalInvDeriv_3A
      + O(r^{-1})\ConformalInvDeriv_4A
      + O(ar^{-2})\ConformalInvDeriv A
      + O(r^{-2}) A
      \right)\\
    ={}& \left(\frac{5}{r^5} + O(ar^{-6})\right)\mathfrak{q}
         + O(r^{-1})\ConformalInvDeriv_4\ConformalInvDeriv_3A
     + O(ar^{-2})\ConformalInvDeriv\ConformalInvDeriv_3 A\\
    &  + O(r^{-2})\ConformalInvDeriv_3A
      + O(r^{-2})\ConformalInvDeriv_4A
      + O(ar^{-3})\ConformalInvDeriv A
      + O(r^{-3}) A.
    \end{split}        
  \end{equation}
  Now, taking a $\ConformalInvDeriv_3$ derivative of
  \zcref[noname]{eq:transport:D-hot-D-A-control:reduced-Teuk} and
  applying
  \zcref[noname]{eq:transport:D-hot-D-A-control:D-hot-D-nab3-A:aux1},
  we see that
  \begin{equation*}
    \begin{split}
      &\frac{1}{4}\ConformalInvDeriv_3\ConformalComplexDeriv\SymTracelessTensorProd\left(\overline{\ConformalComplexDeriv}\cdot A\right)\\
        ={}& \ConformalInvDeriv_3\ConformalInvDeriv_4\ConformalInvDeriv_3A
             + \left(\frac{5}{r^5} + O(ar^{-6})\right)\mathfrak{q}
         + O(r^{-1})\ConformalInvDeriv_4\ConformalInvDeriv_3A
     + O(ar^{-2})\ConformalInvDeriv\ConformalInvDeriv_3 A\\
    &  + O(r^{-2})\ConformalInvDeriv_3A
      + O(r^{-2})\ConformalInvDeriv_4A
      + O(ar^{-3})\ConformalInvDeriv A
      + O(r^{-3}) A.
    \end{split}
  \end{equation*}
  Using again the commutation formulas for the conformally invariant
  derivatives in
  \zcref[cap]{lemma:commutation-formula:conformally-invariant-derivaties}
  and using
  \zcref[noname]{eq:transport:D-hot-D-A-control:nab3-2-A-and-qfrak-relation}
  we have that
  \begin{align*}
      \frac{1}{4}\ConformalInvDeriv_3\ConformalComplexDeriv\SymTracelessTensorProd\left(\overline{\ConformalComplexDeriv}\cdot A\right)
      ={}& \ConformalInvDeriv_4\left( \left(\frac{1}{r^4} + O(ar^{-5})\right)\mathfrak{q} \right)
           + \left(\frac{5}{r^5} + O(ar^{-6})\right)\mathfrak{q}
           + O(r^{-1})\ConformalInvDeriv_4\ConformalInvDeriv_3A\\
      & + O(ar^{-2})\ConformalInvDeriv\ConformalInvDeriv_3 A
        + O(r^{-2})\ConformalInvDeriv_3A
        + O(r^{-2})\ConformalInvDeriv_4A\\
      & + O(ar^{-3})\ConformalInvDeriv A
        + O(r^{-3}) A\\
        ={}& O(r^{-5})\ConformalInvDeriv_4\left( r\mathfrak{q} \right)
           + O(r^{-6})\mathfrak{q}
         + O(r^{-1})\ConformalInvDeriv_4\ConformalInvDeriv_3A
     + O(ar^{-2})\ConformalInvDeriv\ConformalInvDeriv_3 A\\
    &  + O(r^{-2})\ConformalInvDeriv_3A
      + O(r^{-2})\ConformalInvDeriv_4A
      + O(ar^{-3})\ConformalInvDeriv A
      + O(r^{-3}) A.
  \end{align*}
  Applying the commutator formulas in
  \zcref[cap]{lemma:commutation-formula:conformally-invariant-derivaties}
  once more, we see that
  \begin{equation}
    \label{eq:transport:D-hot-D-A-control:aux2}
    \begin{split}
      &\ConformalInvDeriv_3\ConformalComplexDeriv\SymTracelessTensorProd\left(\overline{\ConformalComplexDeriv}\cdot A\right)
        - \ConformalComplexDeriv\SymTracelessTensorProd\left(\overline{\ConformalComplexDeriv}\cdot \ConformalInvDeriv_3A\right)\\
      ={}& O(r^{-1})\ConformalComplexDeriv\SymTracelessTensorProd\left(\overline{\ConformalComplexDeriv}\cdot A\right)
           + O(ar^{-2})\nabla\nabla_3A
           + O(ar^{-3})\nabla_3A
            + O(ar^{-3})\nabla A
           + O(ar^{-4})A.
    \end{split}    
  \end{equation}
  Using
  \zcref[noname]{eq:transport:D-hot-D-A-control:D-hot-D-nab3-A:aux1,eq:transport:D-hot-D-A-control:reduced-Teuk}
  we then have that  for $\delta\le p\le 2-\delta$, we have that
  \begin{equation}  
    \begin{split}
      \int_{\mathcal{M}(\tau_1,\tau_2)}r^{p+7}\abs*{\ConformalComplexDeriv\SymTracelessTensorProd\left( 
      \overline{\ConformalComplexDeriv}\cdot\nabla_3A \right)}^2
      \lesssim{}& a^2\int_{\mathcal{M}(\tau_1,\tau_2)}r^{p+3}\left( \abs*{\nabla\nabla_3A}+\abs*{\nabla A}^2 \right)\\
                & + \BulkNormWeighted{p}[\mathfrak{q}](\tau_1,\tau_2)
                  + \BulkNormWeighted{p}[A](\tau_1,\tau_2),\\
      \int_{\Sigma(\tau)}r^{p+8}\abs*{\ConformalComplexDeriv\SymTracelessTensorProd\left(\overline{\ConformalComplexDeriv}\cdot\nabla_3A \right)}^2
      \lesssim{}& \int_{\Sigma(\tau)}r^{p+4}\left( \abs*{\nabla\nabla_3A}+\abs*{\nabla A}^2 \right)\\
                & + \EnergyFlux_p[\mathfrak{q}](\tau)
                  + \EnergyFlux_p[A](\tau),
    \end{split}
  \end{equation}
  as desired, concluding the proof of \zcref[cap]{lemma:transport:D-hot-D-A-control}.
\end{proof}

Combining \zcref[cap]{lemma:transport:angular-derivatives-of-A:main-aux-lemma} and
\zcref[cap]{lemma:transport:D-hot-D-A-control}, it is a simple
application to recover estimates for $\nabla A$ and $\nabla\nabla_3A$.
\begin{corollary}
  \label{coro:transport:angular-derivatives-of-A}
  Let $\delta\le p\le 2-\delta$. Then,
  \begin{align*}
    \int_{\mathcal{M}(\tau_1,\tau_2)}r^{p+3}\abs*{\nabla A}^2
    \lesssim{}& \BulkNormWeighted{p}[A](\tau_1,\tau_2),&
    \int_{\Sigma(\tau)}r^{p+4}\abs*{\nabla A}^2
    \lesssim{}& \EnergyFlux_p[A](\tau),\\
    \int_{\mathcal{M}(\tau_1,\tau_2)}r^{p+5}\abs*{\nabla\nabla_3A}^2
    \lesssim{}& \BulkNormWeighted{p}[\mathfrak{q}](\tau_1,\tau_2)
                + \BulkNormWeighted{p}[A](\tau_1,\tau_2),&
    \int_{\Sigma(\tau)}r^{p+6}\abs*{\nabla\nabla_3 A}^2
    \lesssim{}& \EnergyFlux_p[A](\tau)
                + \EnergyFlux_p[\mathfrak{q}](\tau).
  \end{align*}
\end{corollary}
\begin{proof}
   The proof follows identically to that of Corollary 11.4.13 in \cite{giorgiWaveEquationsEstimates2024}.
\end{proof}

\subsection{Reduction of the source terms}

Recall that in each of the estimates in
\zcref{prop:rp:Kerr,prop:Morawetz:KdS:main,prop:Killing-energy-estimate,prop:redshift-main},
we accrued terms related to the forcing terms that we have left on the
\RHS{} contained in $\mathcal{N}$-norms. In this section, we now
address these terms in the case where the forcing term $N$ is
precisely that in the actual Teukolsky wave-transport system in \zcref{eq:full-RW-system:wave}. 
\begin{lemma}
  \label{lemma:transport:control-of-source-terms}
  For $\delta\le p\le 2-\delta$, and $N=N_0+N_L$ defined by
  \zcref[noname]{eq:full-RW-system:N0:def,eq:full-RW-system:NL:def},
  we have that
  \begin{equation}
    \label{eq:control-of-source-terms}
    \ForcingTermCombinedNorm{p}^s[\psi, N](\tau_1,\tau_2)
    \lesssim \abs*{a}\CombinedBEFNorm{\delta}^s[\psi, A](\tau_1,\tau_2).
  \end{equation}
\end{lemma}

To prove \zcref[cap]{lemma:transport:control-of-source-terms}, we define
\begin{equation}
  \label{eq:transport:forcing-terms:aux-forcing-norms-def}
  \begin{split}
    \ForcingTermNormReg{Mor}^s[\psi, N](\tau_1,\tau_2)
    &=\sum_{k\le s}\int_{\mathcal{M}(\tau_1,\tau_2)}\left(
      \abs*{\nabla_{\HprVF}(\frakWeightedDeriv^k\psi)}
      + r^{-1}\abs*{\frakWeightedDeriv^k\psi}
      \right)\abs*{\frakWeightedDeriv^kN},\\
    \ForcingTermNormReg{ext}_p^s[\psi, N](\tau_1,\tau_2)
    &= \sum_{k\le s}\int_{\mathcal{M}_{r\ge 4M}(\tau_1,\tau_2)}
      r^{p-1}\left( \abs*{\nabla_4(r\frakWeightedDeriv^k\psi)} + r^{-2}\abs*{\frakWeightedDeriv^k\psi} \right)\cdot\abs*{\frakWeightedDeriv^kN},\\
    \ForcingTermNormReg{En}[\psi, N](\tau_1,\tau_2)
    &= \sum_{k\le s}\abs*{\int_{\mathcal{M}(\tau_1,\tau_2)}\nabla_{\TAlmostKilling}(\frakWeightedDeriv^k\psi) \cdot\frakWeightedDeriv^kN},
  \end{split}
\end{equation}
to be the contributions to $\ForcingTermNorm$ coming from the Morawetz
estimate, the $r^p$-weighted estimates, and the energy estimate
respectively. We will now control each of the contributions to
$\ForcingTermCombinedNorm{p}^s[\psi, N](\tau_1,\tau_2)$ independently.

\begin{lemma}  
  \label{lemma:transport:nonlinear-terms:aux}
  For $\delta\le p\le 2-\delta$, and $N=N_0+N_L$ defined by
  \zcref[noname]{eq:full-RW-system:N0:def} and
  \zcref[noname]{eq:full-RW-system:NL:def},
  we have that
  \begin{equation}
    \label{eq:transport:nonlinear-terms:aux}
    \int_{\mathcal{M}(\tau_1,\tau_2)}r^{\delta+3}\abs*{N_0+N_L}^2
    \lesssim \abs*{a}^2\BulkNormWeighted{\delta}^s[\psi, A](\tau_1,\tau_2).
  \end{equation}
\end{lemma}
\begin{proof}
  We will make use of the fact that $N_0 = O(ar^{-4})\psi$ and $N_L = O(a)\frakWeightedDeriv^{\le 1}\nabla_3 A
    + O(ar^{-1})\frakWeightedDeriv^{\le 1}A$. Thus, we have that
  \begin{align*}
    \int_{\mathcal{M}(\tau_1,\tau_2)}r^{\delta+3}\abs*{N_0+N_L}^2
    \lesssim{}& \abs*{a}^2\int_{\mathcal{M}(\tau_1,\tau_2)}r^{\delta-5}\abs*{\frakWeightedDeriv^{\le s}\psi}^2
                + \abs*{a}^2\int_{\mathcal{M}(\tau_1,\tau_2)}r^{\delta+3}\left( \abs*{\frakWeightedDeriv^{\le s+1}\nabla_3A}^2
                + r^{-2}\abs*{\frakWeightedDeriv^{\le s+1}A}^2 \right)\\
    \lesssim{}& \abs*{a}^2\BulkNormWeighted{\delta}^s[\psi, A],
  \end{align*}
  as desired.
\end{proof}

We now control the forcing term norms that arise from doing the
Morawetz estimate and the $r^p$-weighted estimates.
\begin{lemma}
  \label{lemma:transport:nonlinear-terms:Mor+ext}
  We have the following estimate
  \begin{equation}
    \label{eq:transport:nonlinear-terms:Mor+ext}
    \ForcingTermNormReg{Mor}^s[\psi, N_0 + N_L](\tau_1,\tau_2)
    +\ForcingTermNormReg{ext}_p^s[\psi, N_0+N_L]
    \lesssim \abs*{a}\BulkNormWeighted{\delta}^s[\psi](\tau_1,\tau_2).
  \end{equation}
\end{lemma}
\begin{proof}
  Using the definition of $\ForcingTermNormReg{Mor}^s$ in \zcref[noname]{eq:transport:forcing-terms:aux-forcing-norms-def}, a basic Cauchy-Schwarz, and the fact that
  \begin{equation*}
    \abs*{\nabla_{\HprVF}\frakWeightedDeriv^{\le s}\psi}
    \lesssim \abs*{\nabla_{\HprVF}\frakWeightedDeriv^{\le s}\psi}\bOne_{r\le r_0}
    + \left( \abs*{\nabla_{4}\frakWeightedDeriv^{\le s}\psi}^2
      + \abs*{\frakWeightedDeriv^{\le s+1}\psi}^2      
    \right)\bOne_{r\ge r_0}
    + r^{-2}\abs*{\frakWeightedDeriv^{\le s}\psi}^2
    ,
  \end{equation*}
  we have that 
  \begin{align*}
    &\ForcingTermNormReg{Mor}^s[\psi, N_0+ N_L](\tau_1,\tau_2)\\
    ={}&\int_{\Manifold(\tau_1,\tau_2)}\left(
         \abs*{\nabla_{\HprVF}\frakWeightedDeriv^{\le s}\psi}
         + r^{-1}\abs*{\frakWeightedDeriv^{\le s}\psi}
         \right)\abs*{\frakWeightedDeriv^{\le s}(N_0+N_L)}\\
    \lesssim{}& \left(
                \int_{\Manifold(\tau_1,\tau_2)}r^{-1-\delta}\left(
                \abs*{\nabla_{\HprVF}\frakWeightedDeriv^{\le s}\psi}^2\bOne_{r\le r_0}
                + \left( \abs*{\nabla_4\frakWeightedDeriv^{\le s}\psi}^2
                + \abs*{\frakWeightedDeriv^{\le s+1}\psi}^2                
                \right)\bOne_{r\ge r_0}
                \right)
                + r^{-2}\abs*{\frakWeightedDeriv^{\le s}\psi}^2
                \right)^{\frac{1}{2}}\\
    &\times\left(\int_{\Manifold(\tau_1,\tau_2)}r^{1+\delta}\abs*{\frakWeightedDeriv^{\le s}\left( N_0+N_L \right)}^2\right)^{\frac{1}{2}}\\
    \lesssim{}& \left(\BulkNormWeighted{\delta}^s[\psi](\tau_1,\tau_2)\right)^{\frac{1}{2}}\left(\int_{\Manifold(\tau_1,\tau_2)}r^{1+\delta}\abs*{\frakWeightedDeriv^{\le s}\left( N_0+N_L \right)}^2\right)^{\frac{1}{2}}.
  \end{align*}
  Then, we have directly from \zcref[cap]{lemma:transport:nonlinear-terms:aux}
  that $\ForcingTermNormReg{Mor}^s[\psi, N_0+ N_L](\tau_1,\tau_2)
    \lesssim{} \abs*{a}\BulkNormWeighted{\delta}^s[\psi,A](\tau_1,\tau_2)$. Using a similar argument, we can obtain that
  \begin{align*}
    \ForcingTermNormReg{ext}^s_p[\psi, N_0+ N_L](\tau_1,\tau_2)
    ={}& \int_{\Manifold_{r \ge R}(\tau_1,\tau_2)}
         r^{p-1}\left( \abs*{\nabla_4(r\frakWeightedDeriv^k\psi)} + r^{-2}\abs*{\frakWeightedDeriv^k\psi} \right)\cdot\abs*{\frakWeightedDeriv^{\le s}\left( N_0+N_L \right)}\\
    \lesssim{}& \left( \int_{\Manifold_{r \ge R}(\tau_1,\tau_2)}r^{p-3}\abs*{\frakWeightedDeriv^{\le s+1}\psi}^2 \right)^{\frac{1}{2}}
                \left( \int_{\Manifold_{r \ge R}(\tau_1,\tau_2)}r^{p+1}\abs*{\frakWeightedDeriv^{\le s}\left( N_0+N_L \right)} \right)^{\frac{1}{2}}\\
    \lesssim{}& \left(\BulkNormWeighted{\delta}^s[\psi]\right)^{\frac{1}{2}}
                \left( \int_{\Manifold_{r \ge R}(\tau_1,\tau_2)}r^{p+1}\abs*{\frakWeightedDeriv^{\le s}\left( N_0+N_L \right)}^2\right)^{\frac{1}{2}},
  \end{align*}
  where the last equality follows from the fact that
  $\delta\le p \le 2-\delta$. Another application of \zcref[cap]{lemma:transport:nonlinear-terms:aux} then yields that $\ForcingTermNormReg{ext}_p^s[\psi, N_0+N_L]
    \lesssim\abs*{a}\BulkNormWeighted{\delta}^s[\psi, A](\tau_1,\tau_2)$. We can thus conclude that
  \begin{equation*}
    \ForcingTermNormReg{Mor}^s[\psi, N_0+N_L]
    +\ForcingTermNormReg{ext}_p^s[\psi, N_0+N_L]
    \lesssim\abs*{a}\BulkNormWeighted{\delta}^s[\psi, A](\tau_1,\tau_2),
  \end{equation*}
  proving \zcref[noname]{eq:transport:nonlinear-terms:Mor+ext}.
\end{proof}

\begin{lemma}
  \label{lemma:transport:nonlinear-terms:En}
  We have that
  \begin{align}
    \label{eq:transport:nonlinear-terms:En:N0}
    \ForcingTermNormReg{En}[\psi, N_0](\tau_1,\tau_2)
    \lesssim{}&\abs*{a} \left( \sup_{\tau\in [\tau_1,\tau_2]}\EnergyFlux^s[\psi](\tau)
    + \BulkNormWeighted{\delta}^s[\psi](\tau_1,\tau_2) \right),\\
    \label{eq:transport:nonlinear-terms:En:NL}
    \ForcingTermNormReg{En}[\psi, N_L](\tau_1,\tau_2)
    \lesssim{}&\abs*{a} \left( \sup_{\tau\in [\tau_1,\tau_2]}\EnergyFlux^s[\psi, A](\tau)
    + \BulkNormWeighted{\delta}^s[\psi, A](\tau_1,\tau_2) \right).
  \end{align}
\end{lemma}
\begin{proof}
  Observe that
  \begin{align*}
    \ForcingTermNormReg{En}[\psi, N_0](\tau_1,\tau_2)
    ={}& \abs*{\int_{\Manifold(\tau_1,\tau_2)}\nabla_{\TAlmostKilling}(\frakWeightedDeriv^{\le s}\psi)\cdot\frakWeightedDeriv^{\le s}N_0}\\
    \lesssim{}& \abs*{\int_{\Manifold(\tau_1,\tau_2)}\chi\nabla_{\TAlmostKilling}(\frakWeightedDeriv^{\le s}\psi)\cdot\frakWeightedDeriv^{\le s}N_0}
                + \int_{\Manifold(\tau_1,\tau_2)}(1-\chi)\abs*{\nabla_{\TAlmostKilling}(\frakWeightedDeriv^{\le s}\psi)}\abs*{\frakWeightedDeriv^{\le s}N_0},
  \end{align*}
  where $\chi=\chi(r)$ is some smooth cut-off function which is 1 on $\Manifold_{\trap}$ and is supported in $(1+\delta_{\operatorname{red}})r_{\EventHorizon}< r< r_0$.
  Then we observe that since $1-\chi$ is thus supported on $r\ge r_0$, so that in particular,
  \begin{equation*}
    \begin{split}
      &\int_{\Manifold(\tau_1,\tau_2)}(1-\chi)\abs*{\nabla_{\TAlmostKilling}(\frakWeightedDeriv^{\le s}\psi)}\abs*{\frakWeightedDeriv^{\le s}N_0}\\
    \lesssim{}& \left( \int_{\Manifold(\tau_1,\tau_2)}(1-\chi)r^{-1-\delta}\abs*{\nabla_{\TAlmostKilling}\left( \frakWeightedDeriv^{\le s}\psi \right)}^2 \right)^{\frac{1}{2}}
    \left( \int_{\Manifold(\tau_1,\tau_2)}(1-\chi)r^{1+\delta}\abs*{\frakWeightedDeriv^{\le s}N_0}^2 \right)^{\frac{1}{2}},
    \end{split}    
  \end{equation*}
  which we can bound as before using 
  \zcref{lemma:transport:nonlinear-terms:aux}, giving
  \begin{equation*}
    \int_{\mathcal{M}(\tau_1,\tau_2)}(1-\chi)\abs*{\nabla_{\TAlmostKilling}(\frakWeightedDeriv^{\le s}\psi)}\abs*{\frakWeightedDeriv^{\le s}N_0}
    \lesssim a \BulkNormWeighted{\delta}^s[\psi](\tau_1,\tau_2).
  \end{equation*}
  We now move to estimating $\abs*{\int_{\Manifold(\tau_1,\tau_2)}\chi\nabla_{\TAlmostKilling}(\frakWeightedDeriv^{\le s}\psi)\cdot\frakWeightedDeriv^{\le s}N_0}$.
  Using that $N_0=O(ar^{-4})\psi$, we  have that
  \begin{align*}
    \chi\nabla_{\TAlmostKilling}(\frakWeightedDeriv^{\le s}\psi)\cdot\frakWeightedDeriv^{\le s}N_0={}& O(ar^{-4})\chi\nabla_{\TAlmostKilling}\left( \frakWeightedDeriv^{\le s}\psi \right)\cdot\frakWeightedDeriv^{\le s}\psi\\
    ={}& \Divergence_{\Metric_{M,a,\Lambda}}\left(O(ar^{-4})\chi\abs*{\frakWeightedDeriv^{\le s}\psi}^2\TAlmostKilling\right)
         - O(ar^{-4})\chi\abs*{\frakWeightedDeriv^{\le s}\psi}^2\Divergence_{\Metric_{M,a,\Lambda}}\TAlmostKilling.
  \end{align*}
  Applying the divergence theorem on the first term, we then have that
  \begin{equation*}
    \abs*{\int_{\mathcal{M}(\tau_1,\tau_2)}\chi\nabla_{\TAlmostKilling}(\frakWeightedDeriv^{\le s}\psi)\cdot\frakWeightedDeriv^{\le s}N_0}
    \lesssim \abs*{a}\left(
      \sup_{\tau\in [\tau_1,\tau_2]}\EnergyFlux^s[\psi](\tau)
      + \BulkNormWeighted{\delta}^s[\psi](\tau_1,\tau_2)
    \right),
  \end{equation*}
  where we remark that here, powers of $r$ do not matter since $\chi$
  has compact support. This concludes the proof of
  \zcref[noname]{eq:transport:nonlinear-terms:En:N0}. We now consider the proof of
  \zcref[noname]{eq:transport:nonlinear-terms:En:NL}. We once again divide
  $\ForcingTermNormReg{En}$ into an integral over a function with
  compact support away from the cosmological horizon, which we will
  now control, and an integral over a function supported in a region
  containing the cosmological horizon, which we control again using
  \zcref[cap]{lemma:transport:nonlinear-terms:aux}.  That is, we
  observe that
  \begin{align*}
    \ForcingTermNormReg{En}^s[\psi, N_L](\tau_1,\tau_2)
    \lesssim{}& \abs*{\int_{\mathcal{M}(\tau_1,\tau_2)}\chi \nabla_{\TAlmostKilling}\left( \frakWeightedDeriv^{\le s}\psi \right)\cdot\frakWeightedDeriv^{\le s}N_L}
                +\abs*{\int_{\mathcal{M}(\tau_1,\tau_2)}\left( 1-\chi \right) \nabla_{\TAlmostKilling}\left( \frakWeightedDeriv^{\le s}\psi \right)\cdot\frakWeightedDeriv^{\le s}N_L}\\
    \lesssim{}& \abs*{\int_{\mathcal{M}(\tau_1,\tau_2)}\chi \nabla_{\TAlmostKilling}\left( \frakWeightedDeriv^{\le s}\psi \right)\cdot\frakWeightedDeriv^{\le s}N_L}
                + \abs*{a}\BulkNormWeighted{\delta}^s[\psi, A](\tau_1,\tau_2).
  \end{align*}
  We can then decompose
  \begin{align*}
    \int_{\mathcal{M}(\tau_1,\tau_2)}\chi \nabla_{\TAlmostKilling}\left( \frakWeightedDeriv^{\le s}\psi \right)\cdot\frakWeightedDeriv^{\le s}N_L
    ={}& \int_{\mathcal{M}(\tau_1,\tau_2)}\chi \nabla_{\KillT}\left( \frakWeightedDeriv^{\le s}\psi \right)\cdot\frakWeightedDeriv^{\le s}N_L\\
    & + \int_{\mathcal{M}(\tau_1,\tau_2)}\chi \frac{a}{r^2+a^2}\chi_0\left( \delta^{-1}\frac{3M}{r} \right)\nabla_{\KillPhi}\left( \frakWeightedDeriv^{\le s}\psi \right)\cdot\frakWeightedDeriv^{\le s}N_L.
  \end{align*}
  Recalling that $\chi_0\left(\delta^{-1}\frac{3M}{r}\right)=0$ on $\mathcal{M}_{\trap}(\tau_1,\tau_2)$, we see that
  \begin{equation*}
    \left( \int_{\mathcal{M}(\tau_1,\tau_2)}\abs*{\chi \frac{a}{r^2+a^2}\chi_0\left( \delta^{-1}\frac{3M}{r} \right)\nabla_{\KillPhi}\left( \frakWeightedDeriv^{\le s}\psi \right)} \right)^2
    \lesssim \abs*{a}\BulkNormWeighted{\delta}^s[\psi](\tau_1,\tau_2).
  \end{equation*}
  Applying \zcref[cap]{lemma:transport:nonlinear-terms:aux} as before,
  we have that
  \begin{equation*}
    \abs*{\int_{\mathcal{M}(\tau_1,\tau_2)}\chi \frac{a}{r^2+a^2}\chi_0\left( \delta^{-1}\frac{3M}{r} \right)\nabla_{\KillPhi}\left( \frakWeightedDeriv^{\le s}\psi \right)\cdot\frakWeightedDeriv^{\le s}N_L}\lesssim \abs*{a}\BulkNormWeighted{\delta}^s[\psi](\tau_1,\tau_2).
  \end{equation*}
  Now, recalling the definition of $\mathfrak{q}$ in
  \zcref[noname]{definition of q frak}, we have that $\psi = \psi_0 + \frakWeightedDeriv^{\le 1}A$ where we denote
  \begin{equation}
    \label{eq:transport:forcing-term-control:psi-0:def}
    \psi_0 \vcentcolon= \Re(q \overline{q} \nabla_{3}\nabla_{3}A).
  \end{equation}
  Then, we observe that again using \zcref[cap]{lemma:transport:nonlinear-terms:aux},
  \begin{align*}
    \ForcingTermNormReg{En}^s[\psi, N_L](\tau_1,\tau_2)
    \lesssim{}& \abs*{\int_{\mathcal{M}(\tau_1,\tau_2)}\chi \nabla_{\KillT}(\frakWeightedDeriv^{\le s}\psi_0)\cdot\frakWeightedDeriv^{\le s}N_L}
                + \abs*{a}\BulkNormWeighted{\delta}^s[\psi, A](\tau_1,\tau_2).
  \end{align*}
  To handle the $\abs*{\int_{\mathcal{M}(\tau_1,\tau_2)}\chi \nabla_{\KillT}(\frakWeightedDeriv^{\le s}\psi_0)\cdot\frakWeightedDeriv^{\le s}N_L}$ term, we observe that
  \begin{align*}
    \chi\nabla_{\KillT}(\frakWeightedDeriv^{\le s}\psi_0)\cdot \frakWeightedDeriv^{\le s}N_L
    ={}& -\chi(\frakWeightedDeriv^{\le s}\psi_0)\cdot \nabla_{\KillT}\frakWeightedDeriv^{\le s}N_L
         + \Divergence_{\Metric_{M,a,\Lambda}}(\chi \left( \frakWeightedDeriv^{\le s}\psi_0 \right)\cdot \left( \frakWeightedDeriv^{\le s}N_L \right)\KillT)\\
       &- \chi(\frakWeightedDeriv^{\le s}\psi_0)\cdot(\frakWeightedDeriv^{\le s}N_L)\Divergence_{\Metric_{M,a,\Lambda}}\KillT\\
    ={}& -\chi(\frakWeightedDeriv^{\le s}\psi_0)\cdot \frakWeightedDeriv^{\le s}\nabla_{\KillT}N_L
         + \Divergence_{\Metric_{M,a,\Lambda}}(\chi \left( \frakWeightedDeriv^{\le s}\psi_0 \right)\cdot \left( \frakWeightedDeriv^{\le s}N_L \right)\KillT)\\
        &- \chi(\frakWeightedDeriv^{\le s}\psi_0)\cdot(\frakWeightedDeriv^{\le s}N_L)\Divergence_{\Metric_{M,a,\Lambda}}\KillT.
  \end{align*}
  We now recall from \zcref[noname]{eq:full-RW-system:NL:def} that
  $N_L = O(a)\frakWeightedDeriv^{\le 1}\nabla_3^{\le 1}A$, so we have that
  \begin{equation*}
    \begin{split}
      &\bigg\vert
      \int_{\Manifold(\tau_1,\tau_2)}\left[
      \chi\nabla_{\KillT}(\frakWeightedDeriv^{\le s}\psi_0)\cdot \frakWeightedDeriv^{\le s}N_L
      + \CovariantDeriv_{\alpha}\left(
      \chi\left(\frakWeightedDeriv^{\le s}\psi_0\right)\cdot\left( \frakWeightedDeriv^{\le s}N_L \right)\KillT^{\alpha}
      \right)\right.\\
      & \left.-\chi\left(\frakWeightedDeriv^{\le s}\psi_0\right)\cdot\left( \frakWeightedDeriv^{\le s}N_L \right)\Divergence\KillT
        - \partial_r\chi\KillT(r)\left(\frakWeightedDeriv^{\le s}\psi_0\right)\cdot\left( \frakWeightedDeriv^{\le s}N_L \right)
      \right]
        \bigg\vert\\
      \lesssim{}&\abs*{a}\left(
                  \sup_{\tau\in [\tau_1,\tau_2]}\EnergyFluxCombinedOpt{\delta}^s[\psi, A]
                  + \BulkNormWeighted{\delta}^s[\psi, A](\tau_1,\tau_2)
                  \right).
    \end{split}
  \end{equation*}
  To handle the
  $\chi(\frakWeightedDeriv^{\le s}\psi_0)\cdot \frakWeightedDeriv^{\le
    s}\nabla_{\KillT}N_L$ term, we recall from
  \zcref[noname]{eq:full-RW-system:NL:def} that
  \begin{equation*}
    N_L = -\Re\left(
      \frac{(1+\gamma)8 q\overline{q}^3\Delta}{\lambdaglo r^2\abs*{q}^4} \left( a^2\nabla_{\KillT} + a\nabla_{\KillPhi} \right) \nabla_3A 
    \right)
    + O(a)\frakWeightedDeriv^{\le 1}A. 
  \end{equation*}
  Using $\KillT
    ={} \HprVF + \frac{\abs*{q}^2}{\lambdaglo(1+\gamma)\left( r^2+a^2 \right)}e_3 - \frac{a}{r^2+a^2}\KillPhi$ we then have that
  \begin{align*}
    \nabla_{\KillT} N_L
    ={}& - \Re\left(  \frac{(1+\gamma)8 q\overline{q}^3\Delta}{\lambdaglo r^2\abs*{q}^4}\left(
         a^2 \nabla_{\KillT}^2
         + a\nabla_{\HprVF}\nabla_{\KillPhi}
         + \frac{a\abs*{q}^2}{\lambdaglo(1+\gamma)\left( r^2+a^2 \right)}\nabla_{3}\nabla_{\KillPhi}
         - \frac{a^2}{r^2+a^2}\nabla_{\KillPhi}^2
         \right)\nabla_{3}A \right)\\
    & + O(a)\frakWeightedDeriv^{\le s}A\\
    ={}& - \Re\left(  \frac{(1+\gamma)8 q\overline{q}^3\Delta}{\lambdaglo r^2\abs*{q}^4}\left(
         a^2 \nabla_{\KillT}^2
         + \frac{a\abs*{q}^2}{\lambdaglo(1+\gamma)\left( r^2+a^2 \right)}\nabla_{3}\nabla_{\KillPhi}
         - \frac{a^2}{r^2+a^2}\nabla_{\KillPhi}^2
         \right)\nabla_{3}A \right)\\
    & + O(a)\nabla_{\HprVF}\frakWeightedDeriv^{\le 2}A
      + O(a)\frakWeightedDeriv^{\le s}A.
  \end{align*}
  We can then use the product rule to write that
  \begin{align*}
    &\chi \left( \frakWeightedDeriv^{\le s}\psi_0 \right)\cdot\frakWeightedDeriv^{\le s}\nabla_{\KillT}N_L\\
    ={}& -\chi \left( \frakWeightedDeriv^{\le s}\psi_0 \right)
         \cdot \Re\left(  \frac{(1+\gamma)8 q\overline{q}^3\Delta}{\lambdaglo r^2\abs*{q}^4}\left(
         a^2 \nabla_{\KillT}^2
         + \frac{a\abs*{q}^2}{\lambdaglo(1+\gamma)\left( r^2+a^2 \right)}\nabla_{3}\nabla_{\KillPhi}
         - \frac{a^2}{r^2+a^2}\nabla_{\KillPhi}^2
         \right)\nabla_{3}A \right)\\
    & + O(a)\chi \left( \frakWeightedDeriv^{\le s}\psi_0 \right)\cdot\left( \nabla_{\HprVF}\frakWeightedDeriv^{\le 2}A \right)
      + O(a)\chi \left( \frakWeightedDeriv^{\le s}\psi_0 \right)\cdot\left( \frakWeightedDeriv^{\le s}A \right)\\
    ={}& \CovariantDeriv_{\mu}\left(
         O(a) \chi \frakWeightedDeriv^{\le s}\psi_0\cdot\frakWeightedDeriv^{\le s+2}A(\HprVF^{\mu}, \KillT^{\mu}, \KillPhi^{\mu})
         \right)
         + O(a)\chi \left( \frakWeightedDeriv^{\le s}\psi_0 \right)\cdot\left( \nabla_{\HprVF}\frakWeightedDeriv^{\le 2}A \right)
         + O(a)\chi \left( \frakWeightedDeriv^{\le s}\psi_0 \right)\cdot\left( \frakWeightedDeriv^{\le s}A \right)\\
       &+ \chi(\nabla_{\KillT}\frakWeightedDeriv^{\le s}\psi_0)\cdot\frakWeightedDeriv^{\le s}\Re \left(\frac{(1+\gamma)8 q\overline{q}^3\Delta}{\lambdaglo r^2\abs*{q}^4}a^2\nabla_{\KillT}\nabla_{3}A\right)\\
    &- \chi(\frakWeightedDeriv^{\le s}\psi_0)\cdot\frakWeightedDeriv^{\le s}\Re \left(
      \frac{8q\overline{q}^3}{\abs*{q}^2r^2}\frac{a\Delta}{r^2+a^2}\nabla_{\KillPhi}\nabla_{3}^2A\right)\\
    &- \chi(\nabla_{\KillPhi}\frakWeightedDeriv^{\le s}\psi_0)\cdot\frakWeightedDeriv^{\le s}\Re \left(\frac{(1+\gamma)8 q\overline{q}^3\Delta}{\lambdaglo r^2\abs*{q}^4}\frac{a^2}{r^2+a^2}\nabla_{\KillPhi}\nabla_{3}A\right).
  \end{align*}
  We now rewrite the \RHS{} in terms of $A$, using  \zcref[noname]{eq:transport:forcing-term-control:psi-0:def}
  to see that
  \begin{align*}
    \chi \left( \frakWeightedDeriv^{\le s}\psi_0 \right)\cdot\frakWeightedDeriv^{\le s}\nabla_{\KillT}N_L    
    ={}& \CovariantDeriv_{\mu}\left(
         O(a) \chi \frakWeightedDeriv^{\le s}\psi_0\cdot\frakWeightedDeriv^{\le s+2}A(\HprVF^{\mu}, \KillT^{\mu}, \KillPhi^{\mu})
         \right)\\
       & + O(a)\chi \left(\abs*{\partial_r\chi} + \abs{\chi}\right)\left(\frakWeightedDeriv^{\le s+2}A\right)\cdot\left(\frakWeightedDeriv^{\le s+2}A\right)
         + O(a)\chi\nabla_{\HprVF}\frakWeightedDeriv^{\le s}\psi \cdot\frakWeightedDeriv^{\le s+1}A\\      
       &+ \chi(\nabla_{\KillT}\frakWeightedDeriv^{\le s}\Re(q \overline{q}^3\nabla_{3}^2A))\cdot\frakWeightedDeriv^{\le s}\Re \left(\frac{(1+\gamma)8 q\overline{q}^3\Delta}{\lambdaglo r^2\abs*{q}^4}a^2\nabla_{\KillT}\nabla_{3}A\right)\\
    &- \chi(\frakWeightedDeriv^{\le s}\Re(q \overline{q}^3\nabla_{3}^2A))\cdot\frakWeightedDeriv^{\le s}\Re \left(\frac{8q\overline{q}^3}{\abs*{q}^2r^2}\frac{a\Delta}{r^2+a^2}\nabla_{\KillPhi}\nabla_{3}^2A\right)\\
    &- \chi(\nabla_{\KillPhi}\frakWeightedDeriv^{\le s}\Re(q \overline{q}^3\nabla_{3}^2A))\cdot\frakWeightedDeriv^{\le s}\Re \left(\frac{(1+\gamma)8 q\overline{q}^3\Delta}{\lambdaglo r^2\abs*{q}^4}\frac{a^2}{r^2+a^2}\nabla_{\KillPhi}\nabla_{3}A\right).
  \end{align*}
  We can now rewrite the last three terms on the \RHS{} as divergences
  of squares, up to some error terms,
  \begin{align*}
    \chi \left( \frakWeightedDeriv^{\le s}\psi_0 \right)\cdot\frakWeightedDeriv^{\le s}\nabla_{\KillT}N_L    
    & = \CovariantDeriv_{\mu}\left(
         O(a) \chi \frakWeightedDeriv^{\le s}\psi_0\cdot\frakWeightedDeriv^{\le s+2}A(\HprVF^{\mu}, \KillT^{\mu}, \KillPhi^{\mu})
         \right)\\
       &\quad + O(a)\chi \left(\abs*{\partial_r\chi} + \abs{\chi}\right)\left(\frakWeightedDeriv^{\le s+2}A\right)\cdot\left(\frakWeightedDeriv^{\le s+2}A\right)
         + O(a)\chi\nabla_{\HprVF}\frakWeightedDeriv^{\le s}\psi \cdot\frakWeightedDeriv^{\le s+1}A\\      
       &\quad+ \frac{8a^2\Delta(1+\gamma)}{\lambdaglo r^2\abs*{q}^4}\chi\nabla_{3}(\frakWeightedDeriv^{\le s}\Re(q \overline{q}^3\nabla_{\KillT}\nabla_{3}A))\cdot\frakWeightedDeriv^{\le s}\Re \left(q \overline{q}^3\nabla_{\KillT}\nabla_{3}A\right)\\
    &\quad- \frac{8}{r^2\abs*{q}^2}\frac{a\Delta}{r^2+a^2}\chi\nabla_{\KillPhi}(\frakWeightedDeriv^{\le s}\Re(q \overline{q}^3\nabla_{3}^2A))\cdot\left(\frakWeightedDeriv^{\le s}\Re q \overline{q}^3\nabla_{3}^2A\right)\\
    &\quad-\frac{8(1+\gamma)\Delta}{\lambdaglo r^2\abs*{q}^4}\frac{a^2}{r^2+a^2} \chi\nabla_{3}(\frakWeightedDeriv^{\le s}\Re(q \overline{q}^3\nabla_{\KillPhi}\nabla_{3}A))\cdot\left(\frakWeightedDeriv^{\le s}\Re q \overline{q}^3\nabla_{\KillPhi}\nabla_{3}A\right)
    \\& = \CovariantDeriv_{\mu}\left(
           O(a)\chi\frakWeightedDeriv^{\le s+2}A \cdot\frakWeightedDeriv^{\le s+2}A\left((e_3)^{\mu}, \KillPhi^{\mu}\right)
           \right)\\
         &\quad  + \CovariantDeriv_{\mu}\left(
           O(a)\chi\frakWeightedDeriv^{\le s}\psi_0\cdot \frakWeightedDeriv^{\le s+2} A(\HprVF^{\mu}, \KillT^{\mu}, \KillPhi^{\mu})
           \right)\\
      &\quad + O(a)\chi \left(\abs*{\partial_r\chi} + \abs{\chi}\right)\left(\frakWeightedDeriv^{\le s+2}A\right)\cdot\left(\frakWeightedDeriv^{\le s+2}A\right)
         + O(a)\chi\nabla_{\HprVF}\frakWeightedDeriv^{\le s}\psi \cdot\frakWeightedDeriv^{\le s+1}A.
  \end{align*}
  Applying the divergence theorem and recalling the definition of the
  $\EnergyFlux_{\delta}[\psi, A](\tau)$ norm in
  \zcref[noname]{eq:BEF-A-norm:def}, we see that
  \begin{equation*}
    \abs*{\int_{\mathcal{M}(\tau_1,\tau_2)}\chi(\frakWeightedDeriv^{\le s}\psi_0)\cdot \frakWeightedDeriv^{\le s}\nabla_{\KillT} N_L}
    \lesssim{} \abs*{a}\left(
      \sup_{\tau\in [\tau_1,\tau_2]}\EnergyFlux_{\delta}^s[\psi, A](\tau)
      + \BulkNormWeighted{\delta}^s[\psi, A](\tau_1,\tau_2)
    \right),
  \end{equation*}
  as desired, concluding the proof of \zcref[cap]{lemma:transport:nonlinear-terms:En}.
\end{proof}

\begin{proof}[Proof of Lemma \ref{lemma:transport:control-of-source-terms}]
  Combining
  \zcref[cap]{lemma:transport:nonlinear-terms:En,lemma:transport:nonlinear-terms:aux,lemma:transport:nonlinear-terms:Mor+ext}
  directly yields \zcref[noname]{eq:control-of-source-terms}.
\end{proof}

\subsection{Proof of \zcref{eq:main-q-A-combined-estimate}}

We first prove an auxiliary proposition, from which the proof of
\zcref{eq:main-q-A-combined-estimate} will follow immediately. 
\begin{proposition}
  \label{prop:BEF-controlled-by-B-E}
  Let $(\mathfrak{q}, A)$ be as in \zcref[cap]{MAINTHEOREM}, and let
  $\psi = \Re \mathfrak{q}$.  For $\delta\le p\le 2-\delta$, we have
  that
  \begin{equation}
    \label{eq:BEF-controlled-by-B-E}
    \CombinedBEFNorm{p}^s[A](\tau_1,\tau_2)
    \lesssim \BulkNormWeighted{p}^s[\psi](\tau_1,\tau_2)
    + \EnergyFlux_p^s[A](\tau_1).
  \end{equation}
\end{proposition}
\begin{proof}
  Given \zcref[cap]{coro:transport:angular-derivatives-of-A,prop:BEF-A:tranport-estimates}, we
  have that for all $\delta\le p \le 2-\delta$ and for $a$ sufficiently small,
  \begin{equation*}
    \BEFNormAux{p}[A](\tau_1,\tau_2)
    + \int_{\mathcal{M}(\tau_1,\tau_2)}r^{p+3}\left(r^2\abs*{\nabla\nabla_3A}^2 +\abs*{\nabla A}^2\right)
    \lesssim\BulkNormWeighted{p}[\mathfrak{q}](\tau_1,\tau_2)
    + \EnergyFlux_p[A](\tau_1).
  \end{equation*}
  Then, observe from the definitions of
  $\BulkNormWeightedAux{p}[A](\tau_1,\tau_2)$ and
  $\BulkNormWeighted{p}[A](\tau_1,\tau_2)$ in
  \zcref[noname]{eq:bulk-A-aux-norm:def} and \zcref[noname]{eq:A-norms} that 
  \begin{equation*}
    \BulkNormWeighted{p}[A](\tau_1,\tau_2)
    = \BulkNormWeightedAux{p}[A](\tau_1,\tau_2)
    + \int_{\mathcal{M}(\tau_1,\tau_2)}r^{p+3}\left(r^2\abs*{\nabla\nabla_3A}^2 + \abs*{\nabla A}^2\right),
  \end{equation*}
  so that $\CombinedBEFNorm{p}[A](\tau_1,\tau_2)
    \lesssim\BulkNormWeighted{p}[\mathfrak{q}](\tau_1,\tau_2)
    + \EnergyFlux_p[A](\tau_1)$. Then, since $\mathfrak{q} = \psi + \ImagUnit \LeftDual{\psi}$, we have that $\CombinedBEFNorm{p}[A](\tau_1,\tau_2)
    \lesssim\BulkNormWeighted{p}[\psi](\tau_1,\tau_2)
    + \EnergyFlux_p[A](\tau_1)$. This proves \zcref[noname]{eq:BEF-controlled-by-B-E} when $s=0$. It remains to recover \zcref[noname]{eq:BEF-controlled-by-B-E} for
  $s>0$. We show how to inductively go from $s=0$ to $s=1$. The
  general case is similar. Observe that
  $\left( \nabla_{\KillT}, r\nabla_{4}, r\nabla, \dot{\chi}\nabla_{3}
  \right)$ span $\frakWeightedDeriv$. We can now commute the transport
  equations in \zcref[noname]{eq:full-RW-system:transport} using the
  commutation relations
  \zcref[cap]{lemma:commutation-formula:conformally-invariant-derivaties}
  to see that
  \begin{align*}
    \ConformalInvDeriv_{3}\HorLieDeriv_{\KillT}\Psi
    ={}& O(r^{-2})\HorLieDeriv_{\KillT}\mathfrak{q}
         + O(r^{-3})\mathfrak{q},\\
    \ConformalInvDeriv_{3}\HorLieDeriv_{\KillT}\left(\frac{\overline{q}^4}{r^2}A\right)
    ={}& \HorLieDeriv_{\KillT}\Psi,\\
    \ConformalInvDeriv_{3}\left( q\overline{\ConformalComplexDeriv}\cdot\Psi \right)
    ={}& O(r^{-1})\nabla \mathfrak{q}+ O(r^{-2})\mathfrak{q},\\
    \ConformalInvDeriv_{3}\left( q\overline{\ConformalComplexDeriv}\cdot\left(\frac{\overline{q}^4}{r^2}A\right) \right)
    ={}& q\overline{\ConformalComplexDeriv}\cdot\Psi,\\
    \ConformalInvDeriv_{3}\ConformalInvDeriv_{4}\Psi
    ={}& O(r^{-2})\ConformalInvDeriv_{4}\mathfrak{q}
         + O(r^{-3})\mathfrak{q}
         + O(r^{-3})\Psi,\\
    \ConformalInvDeriv_{3}\ConformalInvDeriv_{4}\left(\frac{\overline{q}^4}{r^2}A\right)
    ={}& \ConformalInvDeriv_4\Psi + O(a)\nabla A+ O(r^{-1})A.
  \end{align*}
  Using the fact that \zcref[noname]{eq:combined-estimate:main} holds
  for $s=0$, we thus have that
  \begin{equation*}
    \WeightedBEFNorm{p}[(\HorLieDeriv_{\KillT},\ConformalInvDeriv_4, q\overline{\ConformalComplexDeriv}\cdot)(\psi, A)](\tau_1,\tau_2)
    \lesssim \EnergyFlux_p[(\HorLieDeriv_{\KillT},\ConformalInvDeriv_4, q\overline{\ConformalComplexDeriv}\cdot)(\psi, A)](\tau_1)
    + \ForcingTermCombinedNorm{p}^1[\psi, \Nonlinearity](\tau_1,\tau_2).
  \end{equation*}
  Using the transport system for $(\Psi, O(r^{-2})\mathfrak{q})$ to
  recover control of the $\nabla_{3}$ derivative of $A$ close to the
  horizon, the relationship between $\HorLieDeriv_{\KillT}$ and
  $\nabla_{\KillT}$ in
  \zcref[noname]{eq:LieT-LieZ-to-nabT-nabZ-comparison} and the Hodge
  elliptic estimates of \zcref[cap]{lemma:hodge:elliptic}, and the
  fact that we assumed that \zcref[noname]{eq:BEF-controlled-by-B-E}
  is true for $s=0$ to absorb lower order terms, we can conclude
  directly.
\end{proof}

\begin{proof}[Proof of \eqref{eq:main-q-A-combined-estimate}]
  By combining \zcref[cap]{prop:combined-estimate:main} and
  \zcref{lemma:transport:control-of-source-terms}, we see that
  \begin{equation*}
    \CombinedBEFNorm{p}^s[\psi](\tau_1,\tau_2)
    \lesssim  \EnergyFluxCombined^s[\psi](\tau_1)
    + \abs*{a}\CombinedBEFNorm{\delta}^s[\psi, A](\tau_1,\tau_2).
  \end{equation*}
  But then using \zcref[cap]{prop:BEF-controlled-by-B-E}, we see that
  for $\frac{\abs*{a}}{M}$ sufficiently small we have $\CombinedBEFNorm{p}^s[\psi, A](\tau_1,\tau_2)
    \lesssim  \EnergyFluxCombined^s[\psi, A](\tau_1)$. Since $\mathfrak{q} = \psi + \ImagUnit\LeftDual{\psi}$, we then have that $\CombinedBEFNorm{p}^s[\mathfrak{q}, A](\tau_1,\tau_2)
    \lesssim  \EnergyFluxCombined^s[\mathfrak{q}, A](\tau_1)$ as desired.
\end{proof}

\section{Estimates in the external region \texorpdfstring{ $\Manifold_e$}{} }
\label{sec:external-estimates}

In this section, we give a proof of \zcref{eq:main-external-estimate},
the main estimate in the external region of the \KdS{} spacetime,
$\Manifold_e$. To prove \zcref{eq:main-external-estimate}, we mainly
repeat the $r^p$ estimates we conducted in $\Manifold$ in
\zcref{sec:transport,sec:rp} in $\Manifold_e$. We first repeat the
$r^p$ estimates in \zcref{sec:rp} concerning the wave equation in
\zcref{eq:full-RW-system:wave}, before combining these with $r^p$
estimates for the transport equations in
\zcref{eq:full-RW-system:transport}.

\subsection{Wave estimate in \texorpdfstring{$\mathcal{M}_{e}$}{the external region}}
We first prove an estimate for $\mathfrak{q}$ by repeating the
proof of \zcref{prop:rp:Kerr} in the region $\Manifold_e$. 
\begin{lemma}
  \label{lemma:external-estimate:q-frak}
  Let $\psi\in \realHorkTensor{2}$ be a solution to
  \zcref[noname]{eq:model-problem-gRW}, and $R\gg M$ be sufficiently
  large. Then we have that for $\delta\le p \le 2-\delta$,
  \begin{equation}
    \label{eq:external-estimate:q-frak}
    \begin{split}
      \ExternalBEFNorm{p}^s[\psi](\tau_{\CosmologicalHorizon},\tau)
      \lesssim{}& \EnergyInit{p}^{s-1}[\psi](\tau)
                + \ForcingTermExternalNorm{p}^s[\psi, N](\tau_{\CosmologicalHorizon},\tau). 
    \end{split}    
  \end{equation}
\end{lemma}
\begin{proof}
  We first prove \zcref{lemma:external-estimate:q-frak} for $s=0$. To
  this end, we choose \zcref{eq:rp-Kerr:multiplier-choice:general},
  where $f_p$ is a non-negative function that is equal to $r^p$ for
  $r\ge 2r_0$, and zero for $r\le r_0$, i.e. with $R = 2r_0$. Applying the divergence theorem in \zcref{eq:div-thm:general}, we
  then have that
  \begin{equation*}
    \begin{split}
      & \int_{\Sigma(\tau)}\JCurrent{\rpVF, \rpLagrangianCorr, \rpOneForm}[\psi]\cdot N_{\Sigma}
    + \int_{\SigmaStar(\tau_{\CosmologicalHorizon},\tau)}\JCurrent{\rpVF, \rpLagrangianCorr, \rpOneForm}[\psi]\cdot N_{\SigmaStar}
        + \int_{\Manifold_e(\tau_{\CosmologicalHorizon},\tau)}\KCurrent{\rpVF, \rpLagrangianCorr, \rpOneForm}[\psi]\\
      \lesssim{}& \int_{\SigmaInit(\tau)}\JCurrent{\rpVF, \rpLagrangianCorr, \rpOneForm}[\psi] \cdot N_{\SigmaInit}
                   + \mathcal{F}(\tau_{\CosmologicalHorizon},\tau),
    \end{split}    
  \end{equation*}
  where
  \begin{equation}
    \begin{split}
      \mathcal{F}(\tau_{\CosmologicalHorizon},\tau)
      :={}&  -\int_{\Manifold_e(\tau_{\CosmologicalHorizon},\tau)}f_p\nabla_4\psi\cdot N
            - \frac{1}{2}\int_{\Manifold_e(\tau_{\CosmologicalHorizon},\tau)}r^{-2}f_p\lambda \nabla_3\psi\cdot N\\
          & - \int_{\Manifold_e(\tau_{\CosmologicalHorizon},\tau)}f_p\nabla_4\psi \cdot \frac{4a\cos\theta}{\abs*{q}^2}\LeftDual{\nabla}_{\KillT}\psi\\
          & - \frac{1}{2}\int_{\Manifold_e(\tau_{\CosmologicalHorizon},\tau)}r^{-2}f_p\lambda \nabla_3\psi\cdot\frac{4a\cos\theta}{\abs*{q}^2}\LeftDual{\nabla}_{\KillT}\psi\\
          & - \int_{\Manifold_e(\tau_{\CosmologicalHorizon},\tau)}\frac{r}{\abs*{q}^2}f_p\psi\cdot N
            - \int_{\Manifold_e(\tau_{\CosmologicalHorizon},\tau)}\frac{r}{\abs*{q}^2}f_p\psi \cdot \frac{4a\cos\theta}{\abs*{q}^2}\LeftDual{\nabla}_{\KillT}\psi
            .
    \end{split}  
  \end{equation}
  Observe that we can trivially bound $\int_{\SigmaInit(\tau)}\JCurrent{\rpVF, \rpLagrangianCorr, \rpOneForm}[\psi] \cdot N_{\SigmaInit}
    \lesssim \EnergyInit{p}[\psi](\tau)$. We can now apply
  \zcref{prop:rp:bulk,prop:rp:main-boundary-estimates}, and a direct
  adaptation of \zcref{lemma:rp:extra-forcing-terms} to conclude the
  proof of \zcref{lemma:external-estimate:q-frak} for $s=0$. To prove
  the higher-order estimates, it suffices to perform a simple
  adaptation of the proof in \zcref{sec:rp:higher-order}.
\end{proof}

\begin{lemma}
  \label{lemma:external:control-of-source-terms}
  For $\delta\le p\le 2-\delta$, and $N=N_0+N_L$ defined by
  \zcref[noname]{eq:full-RW-system:N0:def,eq:full-RW-system:NL:def},
  we have that
  \begin{equation}
    \label{eq:external:control-of-source-terms}
    \ForcingTermExternalNorm{p}^s[\psi, N](\tau_{\CosmologicalHorizon},\tau)
    \lesssim \abs*{a}\CombinedBEFNorm{\delta}^s[\psi, A](\tau_{\CosmologicalHorizon},\tau).
  \end{equation}
\end{lemma}
\begin{proof}
  The proof follows exactly as the control for
  $\ForcingTermNormReg{ext}_p^s[\psi, N](\tau_1,\tau_2)$ in the proof
  of \zcref{lemma:transport:control-of-source-terms}, which also arise
  from the $r^p$-weighted estimates.
\end{proof}

\subsection{Transport estimates in \texorpdfstring{$\mathcal{M}_{e}$}{the external region}}

We first observe that in $\Manifold_e$, we have a simplified version
of the divergence lemma in \zcref{lemma:transport:div-lemma}.
\begin{lemma}[Divergence lemma]
  \label{lemma:external:transport:div-lemma}
  Consider a vectorfield $X\in T\Manifold_e$. We have that
  \begin{equation*}
    - \int_{\Sigma(\tau)}\Metric(X, N)
    - \int_{\SigmaStar(\tau_{\CosmologicalHorizon},\tau)}\Metric(X, N)
    + \int_{\SigmaInit(\tau)}\Metric(X, N)
    = \int_{\Manifold_e(\tau_{\CosmologicalHorizon},\tau)}\Divergence(X),
  \end{equation*}
  where $N$ is the normal to the boundary such that $\Metric(N,e_3)=-1$. This can be rewritten as
  \begin{equation*}
    -\int_{\partial^+\Manifold_e(\tau_{\CosmologicalHorizon},\tau)}\Metric(X, N)
    + \int_{\partial^-\Manifold_e(\tau_{\CosmologicalHorizon},\tau)}\Metric(X, N)
    = \int_{\Manifold_e(\tau_{\CosmologicalHorizon},\tau)}\Divergence(X),
  \end{equation*}
  where $\partial^+\Manifold_e(\tau_{\CosmologicalHorizon},\tau)
    = \Sigma(\tau)\cup \SigmaStar(\tau_{\CosmologicalHorizon},\tau)$ and $\partial^-\Manifold_e(\tau_{\CosmologicalHorizon},\tau)=\SigmaInit(\tau)$.
\end{lemma}
\begin{proof}
  Basic application of the standard divergence lemma. 
\end{proof}
Correspondingly, on $\Manifold_e$, we also have the following
adaptation of \zcref{lemma:general-transport-estimate}, where we note
that since $r\ge r_0$ on $\Manifold_e$, there are no trapping dynamics
in $\Manifold_e$.

\begin{lemma}
  \label{lemma:external:general-transport-estimate}
  Let $\Phi_1,\Phi_2\in \realHorkTensor{2}(\Complex)$ satisfy the
  differential relation
  \zcref{eq:transport:main-lemma:Phi1-Phi2-relation}. 
  Then for all $p\ge \delta$, we have
  \begin{equation*}
    \begin{split}
      &\int_{\Manifold_e(\tau_{\CosmologicalHorizon},\tau)}r^{p-3}\left(
      r^2\abs*{\nabla_{3}\Phi_1} + r^2\abs*{\nabla_{4}\Phi_1}^2 + \abs*{\Phi_1}^2
      \right)\\
      & + \int_{\partial \Manifold^+_e(\tau_{\CosmologicalHorizon},\tau)}r^{p-2}\left(
        r^2\abs*{\nabla_{4}\Phi_1}^2
        + \abs*{\nabla_{\HprVF}\Phi_1}^2
        + \abs*{\Phi_1}^2
        \right)\\
      \lesssim{}& \int_{\Manifold_e(\tau_{\CosmologicalHorizon},\tau)}r^{p-1}\left(
                  r^2\abs*{\nabla_{4}\Phi_2}^2
                  + \abs*{\nabla_{\HprVF}\Phi_2}^2
                  + \abs*{\Phi_2}^2
                  \right)\\
      & + \int_{\SigmaInit(\tau)} r^{p-2}\left(
        r^2\abs*{\nabla_{4}\Phi_1}^2
        + \abs*{\nabla_{\HprVF}\Phi_1}^2
                  + \abs*{\Phi_1}^2
        \right)
       + a^2\int_{\Manifold_e(\tau_{\CosmologicalHorizon},\tau)}r^{p-1}\abs*{\nabla\Phi_1}^2,
    \end{split}
  \end{equation*}
  where $\partial^+\mathcal{M}^+(\tau_{\CosmologicalHorizon},\tau)
    = \Sigma(\tau)\cup \SigmaStar(\tau_{\CosmologicalHorizon},\tau)$ denotes the future boundary of $\mathcal{M}_e(\tau_{\CosmologicalHorizon},\tau)$. 
\end{lemma}
\begin{proof}
  The proof follows in the same way as that of
  \zcref{lemma:general-transport-estimate}.
\end{proof}

Using the divergence theorem in the external region in
\zcref{lemma:external:general-transport-estimate}, we can prove the
analogue of \zcref{prop:BEF-controlled-by-B-E} in the external region. 
\begin{proposition}
  \label{prop:external:BEF-controlled-by-B-E}
  Let $(\mathfrak{q}, A)$ be as in \zcref[cap]{MAINTHEOREM}, and let
  $\psi = \Re \mathfrak{q}$.  For $\delta\le p\le 2-\delta$, we have
  that
  \begin{equation}
    \label{eq:external:BEF-controlled-by-B-E}
    \CombinedBEFNorm{p}^s[A](\tau_{\CosmologicalHorizon},\tau)
    \lesssim \BulkNormWeighted{p}^s[\psi](\tau_{\CosmologicalHorizon},\tau)
    + \EnergyFluxInitAux{p}^s[A](\tau).
  \end{equation}
\end{proposition}
\begin{proof}
  The proof of \zcref{prop:external:BEF-controlled-by-B-E} follows
  identically to the proof of \zcref{prop:BEF-controlled-by-B-E} with
  $\EnergyFluxInitAux{p}^s[A](\tau)$ in the place of
  $\EnergyFlux_p^s[A](\tau_1)$, and integrating and applying the
  divergence theorem in
  \zcref{lemma:external:general-transport-estimate} in the place of
  \zcref{lemma:general-transport-estimate}. We remark that the absence of
  $\breve{\chi}$ terms in the conclusions in the external region are a
  result of $\Manifold_e$ being far from $r=3M$. 
\end{proof}

\subsection{Proof of \zcref{eq:main-external-estimate}}
\label{sec:proof-of-external-estimate}

We are now ready to give a proof of \zcref{eq:main-external-estimate}.
By combining \zcref{lemma:external-estimate:q-frak} and
\zcref{lemma:external:control-of-source-terms}, we see that
\begin{equation*}
  \CombinedBEFNorm{p}^s[\psi](\tau_{\CosmologicalHorizon},\tau)
  \lesssim  \EnergyInit{p}^{s-1}[\psi](\tau)
  + \abs*{a}\CombinedBEFNorm{\delta}^s[\psi, A](\tau_{\CosmologicalHorizon},\tau).
\end{equation*}
Using \zcref{prop:external:BEF-controlled-by-B-E}, we see that
for $\frac{\abs*{a}}{M}$ sufficiently small, $\CombinedBEFNorm{p}^s[\psi, A](\tau_{\CosmologicalHorizon},\tau)
  \lesssim  \EnergyInit{p}^{s-1}[\psi](\tau) + \EnergyFluxInitAux{p}^s[A](\tau)$. Since $\mathfrak{q} = \psi + \ImagUnit\LeftDual{\psi}$, we then have that
\begin{equation}
  \label{eq:proof-of-main-ext-estimate:aux1}
  \CombinedBEFNorm{p}^s[\mathfrak{q}, A](\tau_{\CosmologicalHorizon},\tau)
  \lesssim \EnergyInit{p}^{s-1}[\mathfrak{q}](\tau) + \EnergyFluxInitAux{p}^s[A](\tau).
\end{equation}
Finally, we observe from \zcref{definition of q frak} that $\mathfrak{q} = O(r^4)\widehat{\dk}^2A$. As a result, we have that
\begin{equation}
  \label{eq:proof-of-main-ext-estimate:aux2}
  \EnergyInit{p}^{s-1}[\mathfrak{q}](\tau) + \EnergyFluxInitAux{p}^s[A](\tau)
  \lesssim \EnergyInit{p}^{s+1}[A](\tau).
\end{equation}
Combining
\zcref{eq:proof-of-main-ext-estimate:aux1,eq:proof-of-main-ext-estimate:aux2}
directly yields \zcref{eq:main-external-estimate}.

\section{Proof of \zcref[cap]{theo comparaison}}
\label{section preuve theorem comparaison}

This section is devoted to the proof of \zcref[cap]{theo
  comparaison}. The general idea is as follows: based on the given
solution $(A^0,\q^0)$ of the Teukolsky-wave system on Kerr, we define
a corresponding $(A^\La,\q^\La)$ solving the Teukolsky-wave system on
Kerr-de Sitter and converging on compact sets towards $(A^0,\q^0)$; we
then apply \zcref[cap]{MAINTHEOREM} to $(A^\La,\q^\La)$ and pass to
the limit $\La\to0$ in the estimates, crucially benefiting from the
fact that the estimates are uniform with respect to $\La$. Before putting into practice this strategy, we introduce some useful
notations to distinguish and/or relate Kerr-de Sitter to Kerr. The main spacetime regions are directly defined in the coordinates introduced in \zcref[cap]{sec:adaptedglobalcoordinates}:
\begin{align*}
  \MM_{\mathrm{tot},0} &  =  \pth{\mathbb{R}_\tau \times \left( r_{\EventHorizon,0}(1-\delta_{\Horizon}), +\infty \right)_r\times \Sphere^2}\cap \{\underline{t}\geq 0\},
  \\ \MM_{\mathrm{tot},\La} &  =  \pth{\mathbb{R}_\tau \times \left( r_{\EventHorizon,\La}(1-\delta_{\Horizon}), r_{\CosmologicalHorizon,\La}(1+\delta_{\Horizon}) \right)_r\times \Sphere^2}\cap \{\underline{t}\geq 0\},
\end{align*}
where $r_{\EventHorizon,0}$, $r_{\EventHorizon,\La}$ and $r_{\CosmologicalHorizon,\La}$ are the positions of the event horizon in Kerr and Kerr-de Sitter and of the cosmological horizon in Kerr-de Sitter. We also define
\begin{align*}
  \MM'_{\mathrm{tot},0} &  =  \pth{\mathbb{R}_\tau \times \left( r_{\EventHorizon,0}\pth{1-\frac{\delta_{\Horizon}}{2}}, +\infty \right)_r\times \Sphere^2}\cap \{\underline{t}\geq 0\}.
\end{align*}
Since $\lim_{\Lambda\to 0}r_{\CosmologicalHorizon,\La} = +\infty$, $\lim_{\Lambda\to 0}r_{\EventHorizon, \Lambda} = r_{\EventHorizon, 0}$ and $r_{\EventHorizon, 0} < r_{\EventHorizon, \La}$ for $0<\Lambda\ll 1$, we have $\bigcup_{0<\La\ll 1} \MM_{\mathrm{tot},\La} = \MM_{\mathrm{tot},0}$, and for every compact $K\subset \MM'_{\mathrm{tot},0} $ we have $K\subset \MM_{\mathrm{tot},\La}$ for $\La$ small enough. The notation for the initial hypersurfaces is straightforward: $\widehat{\Si}_{\init,\La}$ denotes $\MM_{\mathrm{tot},\La}\cap\{\underline{t}=0\}$ for both $\La=0$ and $\La>0$.

We let $\Metric_0 = \Metric_{M,a,0}$,
$\Metric_{\Lambda} = \Metric_{M,a,\Lambda}$ be the Kerr and \KdS{}
metric, respectively defined on $\MM_{\mathrm{tot},0}$ and
$\MM_{\mathrm{tot},\La}$. We also consider the two global
principal null frames $\pth{e^\La_\mu}_{\mu=1,2,3,4}$ for $\La>0$
and $\La=0$, as defined in \zcref[cap]{sec:global null
  frame}. With respect to these null frames, let
$\horMetric_0, \horMetric_{\Lambda}$ be the horizontal metric of
Kerr and \KdS{} respectively. Crucially, observe that with respect
to global coordinates, both horizontal structures for Kerr and
\KdS{} are spanned by $\dr_\th$ and
$\dr_\ffi + a\sin^2\th\dr_\tau$. Therefore the horizontal
structures coincide, in the sense that $\HorVFSpace_k(\MM_{\mathrm{tot},\La}) =
\evalAt*{\HorVFSpace_k(\MM_{\mathrm{tot},0})}_{\MM_{\mathrm{tot},\La}}$ for all $\Lambda>0$. However, since the horizontal metrics $g_\La$
and $g_0$ do not coincide, the traceless property of horizontal
tensors is not preserved and the sets $\mathfrak{s}_{2,\La}(\CCC)$
and $\mathfrak{s}_{2,0}(\CCC)$ (defined as in
\zcref[cap]{def:hor-struc:hor-tensor-field} in Kerr or \KdS{})
thus do not coincide. To go from one to the other, we consider the
map
\begin{equation}\label{def PiLa}
  \begin{aligned}
    \Pi_{\Lambda}: \realHorkTensor{2,0}(\Complex)&\longrightarrow \realHorkTensor{2, \Lambda}(\Complex),
    \\ \psi & \longmapsto \Pi_\La \psi   \vcentcolon= \evalAt*{\psi}_{\MM_{\mathrm{tot},\La}} - \half \horMetric_\La^{cd} \pth{\evalAt*{\psi}_{\MM_{\mathrm{tot},\La}}}_{cd}\;\horMetric_\La.
  \end{aligned}
\end{equation}
Note also that the conformal invariance of a family of tensors (in the
sense of \zcref[cap]{def:s-conformally-invariants}) is defined
with respect to the global principal null frames
$\pth{e^\La_\mu}_{\mu=1,2,3,4}$ and $\pth{e^0_\mu}_{\mu=1,2,3,4}$
in the Kerr-de Sitter and Kerr case respectively. Since horizontal
metrics are $0$-conformally invariant, \zcref[noname]{def PiLa}
implies that if $\psi\in \realHorkTensor{2,0}(\Complex)$ is
$s$-conformally invariant in Kerr, then so is $\Pi_\La\psi\in \realHorkTensor{2, \Lambda}(\Complex)$ in
Kerr-de Sitter.

Horizontal covariant derivatives are denoted $\nab^0$ and $\nab^\La$,
and weighted derivatives are denoted $\dk^0$ and $\dk^\La$ or
$\widehat{\dk}^0$ and $\widehat{\dk}^\La$. The norms from
\zcref[cap]{sec:main-quantities}, first introduced on Kerr-de Sitter,
can all be considered on Kerr since the integration domains are
coordinate-based and we can replace every covariant derivatives and
every vector fields by their Kerr counterparts, as well as the
horizontal tensor norms (denoted $|\cdot|_\La$ and $|\cdot|_0$). In
terms of notation, we differentiate between Kerr and \KdS{} quantities
with a $0$ or $\La$ subscript or superscript: for instance,
$\E_{p}^{s,0}[\psi](\tau)$ and $\E_{p}^{s,\La}[\psi](\tau)$ denote the
weighted energy for higher derivatives of $\psi$ on Kerr and
Kerr-de Sitter respectively.

\subsection{Definition of $A^\La$ and initial boundedness}

Let $A^0$ be as in \zcref[cap]{theo comparaison}, and define $A^\La$
to be the solution on $\MM_{\mathrm{tot},\La}$ of
$\LL^\La(A^\La)=0$\footnote{$\LL^\La$ here denotes the Teukolsky
  operator on Kerr-de Sitter.} with data
\begin{align}\label{eq:data for A0}
    \pth{A^\La, \nab^\La_{\dr_{\underline{t}}}A^\La } & \vcentcolon = \pth{\Pi_\La A^0, \nab^\La_{\dr_{\underline{t}}}\Pi_\La A^0  }_{|\widehat{\Si}_{\init,\La}},
\end{align}
on $\widehat{\Si}_{\init,\La}$. In this section we prove a uniform
initial bound for $A^\La$.

\begin{proposition}\label{prop:uniform-bound-data}
We have
\begin{align}\label{eq:initial-boundedness}
    \int_{\widehat{\Si}_{\init,\La}} r^{4-\de} \left| \pth{ \widehat{\dk}^\La}^{\leq s+2} A^\La \right|^2_\La \lesssim C_0,
\end{align}
uniformly in $\La$, where $C_0 \vcentcolon = \int_{\widehat{\Si}_{\init,0}} r^{4-\de} \left| \pth{ \widehat{\dk}^0}^{\leq s+2} A^0 \right|^2_0$.
\end{proposition}

The proof of \zcref[cap]{prop:uniform-bound-data} is based on the
following relatively straightforward lemmas. First we show that we can consider
the following coordinate-based weighted derivatives\footnote{Recall that a coordinate chart on the spheres can be defined to include the poles. For simplicity, we just use the $(\th,\ffi)$ coordinates here.}
\begin{align*}
  \widehat{\dk}^\La_c \vcentcolon = \left\{ r\nab^\La_{\dr_{\underline{t}}} , r \nab^\La_{\dr_r} , \nab^\La_{\dr_\ffi} , \nab^\La_{\dr_\th} \right\},
\end{align*}
instead of $\widehat{\dk}^\La$ and $\widehat{\dk}^0$, allowing a
simpler comparison between Kerr-de Sitter and Kerr. Then we use the
equation satisfied by $A^\La$ to express any of its derivatives on
$\widehat{\Si}_{\init,\La}$ in terms of $A^0$ (which requires a uniformly
spacelike hypersurface, unlike the asymptotically null levels of $\tau$). \zcref[cap]{prop:uniform-bound-data} then follows directly.

\begin{lemma}\label{lem:equiv-coord-derivatives}
    For $\La\geq 0$ and for all $\psi\in\O_2(\MM_{\mathrm{tot},\La})$ and $s\geq 0$ we have $\pth{ \widehat{\dk}^\La}^{\leq s} \psi   = \GO{1}  \pth{ \widehat{\dk}^\La_c}^{\leq s} \psi $ and $\pth{ \widehat{\dk}^\La_c}^{\leq s} \psi   = \GO{1}  \pth{ \widehat{\dk}^\La}^{\leq s} \psi$ uniformly over $\La\in[0,\La_0]$ and on $\MM_{\mathrm{tot},\La}$.
\end{lemma}

\begin{proof}
    From the expressions of $e^\La_\mu(\underline{t})$, $e^\La_\mu(r)$, $e^\La_\mu(\ffi)$ and $e^\La_\mu(\th)$ for $\mu=1,2,3,4$ (see \zcref[cap]{lem:action-global-frame} and \zcref[noname]{eq:action on underline t appendix}), we obtain that the transition matrix $\CC_\La$ between the frame $\pth{re_4,re_3,re_2,re_1}$ defining $\widehat{\dk}^\La$ and the frame $\pth{ r\dr_{\underline{t}} , r\dr_r, \dr_\ffi,\dr_\th  }$ defining $\widehat{\dk}^\La_c$  satisfies
    \begin{align*}
        \CC_\La & = \begin{pmatrix}
                    1   & 1 & 0 & 0
                    \\ 1 + \frac{\La r^2}{3} &  - 1 +  \frac{\La r^2}{3} &  0 & 0
                    \\ 0 & 0 &1 & 0
                    \\ 0 & 0 & 0 & 1
                    \end{pmatrix} 
        + \GO{r^{-1}},&
     \CC_\La^{-1} & = -\frac{1}{2} 
                    \begin{pmatrix}
                     - 1 +  \frac{\La r^2}{3} & -1 & 0 & 0
                    \\ -1 - \frac{\La r^2}{3} &  1 &  0 & 0
                    \\ 0 & 0 &1 & 0
                    \\ 0 & 0 & 0 & 1
                    \end{pmatrix}
        + \GO{r^{-1}},
    \end{align*}
    where the $\GO{r^{-1}}$ are uniform over $\La\in[0,\La_0]$ and on
    $\MM_{\mathrm{tot},\La}$. This shows that both $\CC_\La$ and
    $\CC_\La^{-1}$ are uniformly $\GO{1}$, as well as their
    $\widehat{\dk}^\La$ or $\widehat{\dk}^\La_c$ derivatives, which
    concludes the proof of the lemma.
\end{proof}

We recall that $\O(\MM_{\mathrm{tot},\La})=\O(\MM_{\mathrm{tot},0})_{|\MM_{\mathrm{tot},\La}}$, so that a tensor $\psi\in\O_2(\MM_{\mathrm{tot},\La})$ can be differentiated with both $\nab^\La$ and $\nab^0$. 

\begin{lemma}\label{lem:comparaison-derivatives}
    If $\psi\in\O_2(\MM_{\mathrm{tot},\La})$ is a symmetric horizontal 2-tensor then $ \pth{\widehat{\dk}^\La_c}^s\psi  = \GO{1}\pth{\widehat{\dk}^0_c}^{\leq s}\psi$.
\end{lemma}

\begin{proof}
By definition of the horizontal covariant derivative we have 
\begin{align}\label{eq:diff-derivee-covariante}
    \pth{\nab^\La_{X} \psi  - \nab^0_{X} \psi}\pth{e^\La_a,e^\La_b} & = \psi\pth{ \pth{\nab^0_{X}-\nab^\La_{X}}e^\La_a ,e^\La_b } + \psi\pth{ \pth{\nab^0_{X}-\nab^\La_{X}}e^\La_b ,e^\La_a } ,
\end{align}
where $X\in\{ r\dr_{\underline{t}} , r\dr_r, \dr_\ffi,\dr_\th\}$. Following the notations in the proof of \zcref[cap]{lem:equiv-coord-derivatives}, we can schematically write $X=r (\CC_\La^{-1})^\mu e^\La_\mu$ for $\La\geq 0$ and where $(\CC_\La^{-1})^\mu$ denotes coefficient of $\CC_\La^{-1}$. Therefore
\begin{align*}
    \pth{\nab^0_{X}-\nab^\La_{X}}e^\La_a & = r (\CC_0^{-1})^\mu\frac{e^0_\mu(f^\La_a )}{f^\La_a} e^\La_a + r (\CC_0^{-1})^\mu f^\La_a \nab^0_{e^0_\mu} e^0_a  -  r (\CC_\La^{-1})^\mu\nab^\La_{e^\La_\mu}e^\La_a,
\end{align*}
where we also used the notation $e^\La_a = f^\La_a e^0_a$ (i.e $f^\La_1=\sqrt{\ka}$ and $f^\La_2=\frac{1+\ga}{\sqrt{\ka}}$). Recalling the notation $(\Lambda_{\alpha})_{\beta\gamma} = \Metric\left(\CovariantDeriv_{\alpha}e_{\gamma},e_{\beta}\right)$, we get $\pth{\nab^0_{X}-\nab^\La_{X}}e^\La_a = (\DD^\La_X)^d_a e^\La_d$ with 
\begin{align*}
    (\DD^\La_X)^d_a & =  r (\CC_0^{-1})^\mu \pth{\frac{e^0_\mu(f^\La_a )}{f^\La_a} \de^d_a  +  \pth{\frac{f^\La_a}{f^\La_d} - 1 } \de^{cd} \pth{ \La^0_\mu}_{ca} }
    \\&\quad +  \de^{cd} r \pth{ \pth{ (\CC_0^{-1})^\mu - (\CC_\La^{-1})^\mu }  \pth{ \La^0_\mu}_{ca}    +  (\CC_\La^{-1})^\mu \pth{ \pth{ \La^0_\mu}_{ca}  -  \pth{ \La^\La_\mu}_{ca} } }.
\end{align*}
Using \zcref[noname]{eq:Kerr:outgoing:Lambda-Ricci},
$\CC_0^{-1} - \CC_\La^{-1}=\GO{1}$, and $\CC_\La^{-1}=\GO{1}$ (see the
proof of \zcref[cap]{lem:equiv-coord-derivatives}), the second line in the above expression of $(\DD^\La_X)^d_a $ is thus $O(1)$. Similarly, the fact that $f_a^\La = 1 + \GO{\La a^2}$ finally implies that $(\DD^\La_X)^d_a=\GO{1}$. Coming back to \zcref[noname]{eq:diff-derivee-covariante} this implies $\widehat{\dk}^\La_c \psi  = \GO{1} \pth{\widehat{\dk}^0_c}^{\leq 1} \psi$. We conclude the proof of the lemma by induction.
\end{proof}

\begin{lemma}\label{lem:expressing-derivatives}
    For all $s\geq 0$ we have
    \begin{align*}
        \pth{\pth{\widehat{\dk}^\La_c}^{\leq s} A^\La}_{|\widehat{\Si}_{\init,\La}} & = \GO{1} \pth{\pth{\widehat{\dk}^\La_c}^{\leq s} \Pi_\La A^0}_{|\widehat{\Si}_{\init,\La}}.
    \end{align*}
\end{lemma}

\begin{proof}
    We introduce the notation $ \widehat{\dk}^\La_{c,tan} \vcentcolon = \left\{ r \nab^\La_{\dr_r} , \nab^\La_{\dr_\ffi} , \nab^\La_{\dr_\th} \right\}$ for derivatives tangent to $\widehat{\Si}_{\init,\La}$. Thanks to \zcref[noname]{eq:data for A0} we already have 
    \begin{align}\label{initialization}
        \pth{\pth{\widehat{\dk}^\La_{c,tan}}^{\leq s}\pth{r\nab^\La_{ \dr_{\underline{t}} }}^{\leq 1} A^\La}_{|\widehat{\Si}_{\init,\La}} & =  \pth{\pth{\widehat{\dk}^\La_{c,tan}}^{\leq s} \pth{r\nab^\La_{ \dr_{\underline{t}} }}^{\leq 1} \Pi_\La A^0}_{|\widehat{\Si}_{\init,\La}},
    \end{align}
    for all $s\geq 0$. For higher order normal derivatives of $A^\La$, we use the equation it solves, namely $\LL^\La(A^\La)=0$. By expanding all the conformal derivatives appearing in \zcref[noname]{eq:Teuk:A:def} this equation becomes 
    \begin{align*}
        -\nab^\La_{e^\La_4}  \nab^\La_{e^\La_3} A^\La     + \de^{ab} \nab^\La_{e^\La_a}\nab^\La_{e^\La_b}A^\La  & = - \pth{\mathcal{T}^\La_3 \nab^\La_{e^\La_3} + \mathcal{T}^\La_4 \nab^\La_{e^\La_4} + \de^{ab}\mathcal{T}^\La_a \nab^\La_{e^\La_b}  + \mathcal{T}^\La_0 } A^\La ,
    \end{align*}
    where the lower order terms $\mathcal{T}_{\mu}^{\Lambda}$ can be estimated
    with \zcref[cap]{lemma:Kerr:outgoing-PG:Ric-and-curvature} and
    satisfy
    \begin{align*}
        \mathcal{T}^\La_3 & = \GO{r^{-1}} ,& \mathcal{T}^\La_4 & = \GO{r^{-1}} ,&  \mathcal{T}^\La_b & = \GO{a r^{-2}} ,&  \mathcal{T}^\La_0 & = \GO{r^{-2}}.
    \end{align*}
Therefore we have $-\nab^\La_{e^\La_4}  \nab^\La_{e^\La_3} A^\La    +  \de^{ab} \nab^\La_{e^\La_a}\nab^\La_{e^\La_b}A^\La  = \GO{r^{-2}} \pth{\widehat{\dk}^\La_c}^{\leq 1}A^\La  $, where we also used \zcref[cap]{lem:equiv-coord-derivatives} to replace $\pth{\widehat{\dk}^\La}^{\leq 1}$ by $\pth{\widehat{\dk}^\La_c}^{\leq 1}$. By writing $\nab^\La_{e^\La_\mu} = e^\La_\mu(\underline{t})\nab^\La_{\dr_{\underline{t}}} + \GO{r^{-1}}\widehat{\dk}^\La_{c,tan}$ we can expand the \LHS{} of this equation, and obtain 
    \begin{align*}
        \g^\La\pth{\D^\La\underline{t},\D^\La\underline{t}} \pth{r\nab^\La_{\dr_{\underline{t}} }}^2 A^\La = \GO{1} \pth{ \widehat{\dk}^\La_{c,tan} \pth{r\nab^\La_{\dr_{\underline{t}} }} + \pth{\widehat{\dk}^\La_{c,tan}}^2 + \pth{\widehat{\dk}^\La_c}^{\leq 1} }A^\La.    
    \end{align*}
    Using \zcref[noname]{eq:asympt-gtt} and applying an arbitrary number of $r\nab^\La_{\dr_{\underline{t}}}$ this becomes
    \begin{align*}
        \pth{r\nab^\La_{\dr_{\underline{t}} }}^{s} A^\La = \GO{1}\sum_{\substack{p+q\leq s\\q\leq s-1}}  \pth{\widehat{\dk}^\La_{c,tan}}^{p} \pth{r\nab^\La_{\dr_{\underline{t}} }}^q A^\La.
    \end{align*}
    Together with \zcref[noname]{initialization} this can be used to prove by induction on $s'$ that 
    \begin{align*}
        \pth{\pth{\widehat{\dk}^\La_{c,tan}}^{\leq s}\pth{r\nab^\La_{ \dr_{\underline{t}} }}^{\leq s'} A^\La}_{|\widehat{\Si}_{\init,\La}} & = \GO{1} \pth{\pth{\widehat{\dk}^\La_{c}}^{\leq s+s'}  \Pi_\La A^0}_{|\widehat{\Si}_{\init,\La}},
    \end{align*}
    for $s,s'\geq 0$, which concludes the proof.
\end{proof}

We are now ready to prove \zcref[cap]{prop:uniform-bound-data}.

\begin{proof}[Proof of Proposition \ref{prop:uniform-bound-data}]
  By first applying \zcref[cap]{lem:equiv-coord-derivatives} and then \zcref[cap]{lem:expressing-derivatives} we obtain
  \begin{align*}
    \left| \pth{ \widehat{\dk}^\La}^{\leq s} A^\La \right|_\La & \lesssim \left| \pth{ \widehat{\dk}^\La_c}^{\leq s} A^\La \right|_\La
     \lesssim \left| \pth{ \widehat{\dk}^\La_c}^{\leq s} \Pi_\La A^0 \right|_\La.
  \end{align*}
  Recalling \zcref[noname]{def PiLa} we note that $\Pi_\La$ commutes with the $\widehat{\dk}^\La_c$ derivatives and also satisfies $|\Pi_\La\psi|_\La \lesssim |\psi|_\La$, giving $\left| \pth{ \widehat{\dk}^\La}^{\leq s} A^\La \right|_\La  \lesssim \left| \pth{ \widehat{\dk}^\La_c}^{\leq s} A^0 \right|_\La$. Using now \zcref[cap]{lem:comparaison-derivatives} and the uniform closeness of the horizontal metrics $g_\La$ and $g_0$ (which can be seen from $e^\La_b=\pth{1+\GO{\La a^2}}e^0_b$) this becomes $\left| \pth{ \widehat{\dk}^\La}^{\leq s} A^\La \right|_\La  \lesssim \left| \pth{ \widehat{\dk}^0_c}^{\leq s} A^0 \right|_0$. Using again \zcref[cap]{lem:equiv-coord-derivatives} in the $\La=0$ case, we finally obtain
  \begin{align*}
    \left| \pth{ \widehat{\dk}^\La}^{\leq s} A^\La \right|_\La & \lesssim \left| \pth{ \widehat{\dk}^0}^{\leq s} A^0 \right|_0,
  \end{align*}
  from which it is straightforward to deduce \zcref[cap]{prop:uniform-bound-data} since $\widehat{\Si}_{\init,\La}\subset\widehat{\Si}_{\init,0}$.
\end{proof}

\subsection{Convergence on compact sets}

To describe the convergence properties of $(A^\La,\q^\La)$\footnote{$\q^\La$ is defined with respect to $A^\La$ according to \zcref[noname]{definition of q frak}.} when $\La\to0$, we introduce the following notation: for $\La>0$ or $\La=0$ and $\Om\subset\MM_{\mathrm{tot},\La}$ or $\Om\subset\MM'_{\mathrm{tot},0}$, denote by $\mathbf{BE}^{s,\La}_{p,\Om}\left[\q^\La,A^\La\right]$ the (higher-order) bulk-energy norm of $(A^\La,\q^\La)$ defined as the sum of the corresponding quantities from \zcref[cap]{sec:main-quantities} where the integration domains are all intersected with $\Om$. The particular cases of $\Om=\MM_{\mathrm{tot},\La}$ or $\Om=\MM'_{\mathrm{tot},0}$ are denoted $\mathbf{BE}^{s,\La}_{p,\mathrm{tot}}\left[\q^\La,A^\La\right]$. We start by deducing uniform global boundedness from \zcref[cap]{MAINTHEOREM} and \zcref[cap]{prop:uniform-bound-data}.

\begin{lemma}
  For $C_0$ as defined in \zcref{prop:uniform-bound-data} and uniformly over $\La\in(0,\La_0]$ we have
  \begin{align}\label{eq:totale estimée La}
    \mathbf{BE}^{s-1,\La}_{2-\de,\mathrm{tot}} \left[\q^\La,A^\La\right] \lesssim C_0,
  \end{align}
\end{lemma}

\begin{proof}
  We apply \zcref[cap]{MAINTHEOREM} to the sequence $(\q^\La,A^\La)$,
  starting with the external region. Thanks to
  \zcref[noname]{eq:initial-boundedness,eq:main-external-estimate} we have
  \begin{align}\label{estimée totale externe La}
    \ExternalBEFNorm{2-\de}^{s-1,\La}\left[\mathfrak{q}^\La, A^\La\right](\tau_{\overline{\HH},\La},0) \lesssim C_0.
  \end{align}
  Turning now to the region $\tau\geq 0$, we decompose
  \begin{equation}
    \label{eq:convergence:external-E-decomposition}
    \mathbf{E}^{s-1,\La}_{2-\de}\left[\q^\La,A^\La\right](0)
    = \mathbf{E}^{s-1,\La}_{2-\de,r\leq 2r_0}\left[\q^\La,A^\La\right](0)
    + \mathbf{E}^{s-1,\La}_{2-\de,r\geq 2r_0}\left[\q^\La,A^\La\right](0).
  \end{equation}
  The second term on the \RHS{} is trivially bounded by
  $\ExternalBEFNorm{2-\de}^{s-1,\La}\left[\mathfrak{q}^\La,
    A^\La\right](\tau_{\overline{\HH},\La},0)$ and thus
  \zcref[noname]{estimée totale externe La} implies
  $\mathbf{E}^{s-1,\La}_{2-\de,r\geq
    2r_0}\left[\q^\La,A^\La\right](0)\lesssim C_0$. For the first term
  on the \RHS{} of \zcref{eq:convergence:external-E-decomposition}, we
  can adjust the weights at the cost of additional powers of $r_0$
  which we absorb in the universal constants, and show that
  \begin{align*}
    \mathbf{E}^{s-1,\La}_{2-\de,r\leq 2r_0}\left[\q^\La,A^\La\right](0) \lesssim \int_{\widehat{\Si}_{\init,\La}\cap\{r\leq 2r_0\}} r^{4-\de} \left| \pth{ \widehat{\dk}^\La}^{\leq s+2} A^\La \right|^2_\La \lesssim C_0,
  \end{align*}
  where the second inequality follows from 
  \zcref[noname]{eq:initial-boundedness} again. %
  Our uniform estimate
  \zcref[noname]{eq:main-q-A-combined-estimate} then implies that $\mathbf{BE}^{s-1,\La}_{2-\de}\left[\q^\La,A^\La\right](0,+\infty)\lesssim C_0$. Adding this to \zcref[noname]{estimée totale externe La} concludes the proof of \zcref[noname]{eq:totale estimée La}.
\end{proof}

We now focus on particular compact subsets of $\MM'_{\mathrm{tot},0}$ defined for $n$ large enough:
\begin{align*}
    K_n \vcentcolon = \MM'_{\mathrm{tot},0} \cap \{-n\leq \tau \leq n\} \cap \left\{ \tau + \frac{3}{2}r\leq 2n \right\}.
\end{align*}
The sets $K_n$ are clearly compact and cover $\MM'_{\mathrm{tot},0}$,
i.e. $\MM'_{\mathrm{tot},0}=\bigcup_{n}K_n$. Moreover, one can easily
check that $K_n$ are causal for the Kerr metric, in the sense that
$K_n\subset J^+_0\pth{ K_n \cap \widehat{\Si}_{\init,0}}$, where
$J^+_0$ is the causal future defined with respect to the Kerr
metric. Finally, for all $n$, we have
$K_n\subset \MM_{\mathrm{tot},\La}$ for $\La$ small enough depending
on $n$. In particular, $(\q^\La,A^\La)$ is defined on $K_n$ for $\La$ small enough depending
on $n$.

\begin{proposition}
  For all $n$ we have
  \begin{align}\label{limite sur Kn}
    \lim_{\La\to0}\mathbf{BE}^{s-5,0}_{2-\de,K_n}\left[ \q^\La - \q^0, A^\La - A^0 \right] & = 0.
  \end{align}
\end{proposition}

\begin{proof}
  On a fixed compact subset of both $\MM'_{\mathrm{tot},0}$ and $\MM_{\mathrm{tot},\La}$, the Kerr and Kerr-de Sitter metric uniformly converge to each other when $\La\to0$. In particular, we can show that
  \begin{align}\label{changement de norme}
    \mathbf{BE}^{s-4,0}_{2-\de,K_n}\left[\q^\La,A^\La\right] = \mathbf{BE}^{s-4,\La}_{2-\de,K_n}\left[\q^\La,A^\La\right] + \GO{\La}\mathbf{BE}^{s-1,\La}_{2-\de,K_n}\left[\q^\La,A^\La\right],
  \end{align}
  where the $\GO{\La}$ depends on $n$ and the loss of derivatives in
  the second term is due to the off-diagonal terms in the transition
  matrix between the Kerr and Kerr-de Sitter null frames (which
  appears when going from $\mathbf{BE}^{s-4,0}_{2-\de,K_n}$ to
  $\mathbf{BE}^{s-4,\La}_{2-\de,K_n}$). Since
  $K_n\subset \MM_{\mathrm{tot},\La}$ for $\La$ small enough depending
  on $n$ we deduce from \zcref[noname]{eq:totale estimée La} the
  uniform bound
  \begin{align}\label{eq:uniform bound sur Kn}
    \mathbf{BE}^{s-4,0}_{2-\de,K_n}\left[\q^\La,A^\La\right] \lesssim C_0.
  \end{align}
  Therefore, by compact embedding there exists a subsequence
  $(\q^{\La_m},A^{\La_m})_m$ with $\lim_{m\to \infty}\La_m= 0$ and
  some $(\tilde{\q}^0,\tilde{A}^0)$ defined on $K_n$ such that
  \begin{align}\label{limite subsequence}
    \lim_{m\to+\infty}\mathbf{BE}^{s-5,0}_{2-\de,K_n}\left[\q^{\La_m} - \tilde{\q}^0,A^{\La_m} - \tilde{A}^0\right] & = 0.
  \end{align}
  Since $s\geq 5$, \zcref[noname]{limite subsequence} allows us to pass to the limit in the Teukolsky-wave system satisfied by $(\q^{\La_m},A^{\La_m})$ and get that $(\tilde{\q}^0,\tilde{A}^0)$ is a solution to the Teukolsky-wave system of Kerr on $K_n$. However, according to \zcref[cap]{eq:data for A0}, the limit \zcref[noname]{limite subsequence} implies
  \begin{align*}
    \pth{ \pth{\dk^0}^{\leq 1}\tilde{\q}^0 , \pth{\dk^0}^{\leq 1}\tilde{A}^0 }_{|\widehat{\Si}_{\init,0}\cap K_n}  & = \pth{ \pth{\dk^0}^{\leq 1}\q^0 , \pth{\dk^0}^{\leq 1}A^0 }_{|\widehat{\Si}_{\init,0}\cap K_n} .
  \end{align*}
  Therefore, $(\tilde{\q}^0,\tilde{A}^0)$ and $(\q^0,A^0)$ solve the
  same system on $K_n$ and have the same data on
  $\widehat{\Si}_{\init,0}\cap K_n$. Since
  $K_n\subset J^+_0\pth{\widehat{\Si}_{\init,0}\cap K_n}$, uniqueness
  for the Teukolsky wave system on Kerr implies
  $(\tilde{\q}^0,\tilde{A}^0)=(\q^0,A^0)$ on $K_n$. We have thus proved that $\lim_{m\to+\infty}\mathbf{BE}^{s-5,0}_{2-\de,K_n}\left[\q^{\La_m} - \q^0,A^{\La_m} - A^0\right]  = 0$. Since the same reasoning could have been applied to any subsequence
  of $(\q^\La,A^\La)_\La$, this shows that this limit actually holds
  for the whole sequence and not just a subsequence, which concludes
  the proof.
\end{proof}

We can already deduce a global estimate for $(\q^0,A^0)$.

\begin{corollary}
  We have
  \begin{align}\label{eq:totale estimée 0}
    \mathbf{BE}^{s-5,0}_{2-\de,\mathrm{tot}}\left[ \q^0,A^0\right] \lesssim C_0.
  \end{align}
\end{corollary}

\begin{proof}
  Since $\MM'_{\mathrm{tot},0}=\bigcup_n K_n$, the Beppo-Levi theorem implies $ \mathbf{BE}^{s-5,0}_{2-\de,\mathrm{tot}}\left[ \q^0,A^0\right]  = \lim_{n\to +\infty}\mathbf{BE}^{s-5,0}_{2-\de,K_n}\left[ \q^0,A^0\right]$. However, \zcref[noname]{limite sur Kn} and \zcref[noname]{eq:uniform bound sur Kn} imply that $\mathbf{BE}^{s-5,0}_{2-\de,K_n}\left[ \q^0,A^0\right]  = \lim_{\La\to0}\mathbf{BE}^{s-5,0}_{2-\de,K_n}\left[ \q^\La,A^\La\right] \lesssim C_0$ for all $n$, which concludes the proof. 
\end{proof}

By losing a bit of decay, we can also get a convergence statement for
the energy on a $\tau$-slice in Kerr-de Sitter.

\begin{corollary}
    For $p$ as in \zcref[cap]{theo comparaison}, i.e. $p<2-\de$, we have
    \begin{align}\label{convergence de l'énergie quitte à perdre du poids}
        \lim_{\La\to0} \mathbf{E}^{s-6,\La}_{p}\left[\q^\La,A^\La \right](\tau) & = \mathbf{E}^{s-6,0}_{p}\left[\q^0,A^0 \right](\tau).
    \end{align}
\end{corollary}

\begin{proof}
  We first note that due to the scaling of the $r^p$ energy terms, we have
  \begin{align*}
    \dot{\mathbf{E}}^{s-5,\La}_{p,r\geq R}\left[\q^\La\right](\tau) + \mathbf{E}^{s-5,\La}_{p,r\geq R} \left[A^\La\right](\tau) & \lesssim \frac{1}{R^{2-\de-p}}\pth{ \dot{\mathbf{E}}^{s-5,\La}_{2-\de,r\geq R}\left[\q^\La\right](\tau) + \mathbf{E}^{s-5,\La}_{2-\de,r\geq R} \left[A^\La\right](\tau) }
     \lesssim \frac{C_0}{R^{2-\de-p}},
  \end{align*}
  where we also used \zcref[noname]{eq:totale estimée La}. The same
  argument in the Kerr case, but with \zcref[noname]{eq:totale estimée
    0} instead, shows that
  \begin{align*}
    \dot{\mathbf{E}}^{s-5,0}_{p,r\geq R}\left[\q^0\right](\tau) + \mathbf{E}^{s-5,0}_{p,r\geq R} \left[A^0\right](\tau) & \lesssim \frac{C_0}{R^{2-\de-p}}.
  \end{align*}
  To bound the energy terms that do not scale with $p$,
  i.e. $\mathbf{E}^{s-5,\La}_{r\geq R}\left[\q^\La\right](\tau)$ and
  $\mathbf{E}^{s-5,0}_{r\geq R}\left[\q^0\right](\tau)$, we need to
  lose a derivative and use the zero-th order term in the $r^p$
  energy, which gives (after using \zcref[noname]{eq:totale estimée
    La} and \zcref[noname]{eq:totale estimée 0})
  \begin{align*}
    \mathbf{E}^{s-6,\La}_{r\geq R}\left[\q^\La\right](\tau) + \mathbf{E}^{s-6,0}_{r\geq R}\left[\q^0\right](\tau) & \lesssim \frac{1}{R^{1-\de}}\pth{ \dot{\mathbf{E}}^{s-5,\La}_{1-\de,r\geq R}\left[ \q^\La\right] + \dot{\mathbf{E}}^{s-5,0}_{1-\de,r\geq R}\left[ \q^0\right]  }
     \lesssim \frac{C_0}{R^{1-\de}}.
  \end{align*}
  Therefore we have
  \begin{align}\label{tail estimate}
    \mathbf{E}^{s-6,\La}_{p,r\geq R}\left[\q^\La,A^\La \right](\tau) + \mathbf{E}^{s-6,0}_{p,r\geq R}\left[\q^0,A^0 \right](\tau) \lesssim \frac{C_0}{R^\nu},
  \end{align}
  for $\nu=\min(2-\de-p,1-\de)>0$. Now, for $\e>0$ fixed, let $R$ be sufficiently
  large so that $\frac{C_0}{R^\nu}<\e$. Thanks to
  \zcref[noname]{limite sur Kn} and \zcref[noname]{changement de
    norme} we have $\left| \mathbf{E}^{s-6,\La}_{p,r\leq R}\left[\q^\La,A^\La \right](\tau) - \mathbf{E}^{s-6,0}_{p,r\leq R}\left[\q^0,A^0 \right](\tau) \right| \leq \e$
  for $\La$ sufficiently small depending only on $\e$. Thanks to \zcref[noname]{tail estimate} we thus have $\left| \mathbf{E}^{s-6,\La}_{p}\left[\q^\La,A^\La \right](\tau) - \mathbf{E}^{s-6,0}_{p}\left[\q^0,A^0 \right](\tau) \right| \lesssim \e$
  for $\La$ sufficiently small depending only on $\e$, which concludes the proof.
\end{proof}

\subsection{Conclusion of the proof of \zcref[cap]{theo comparaison}}

We conclude here the proof of \zcref[cap]{theo comparaison}. Again
from the Beppo-Levi theorem and $p< 2-\de$ we get
\begin{align*}
    \mathbf{BE}^{s-6,0}_p\left[\q^0,A^0 \right](\tau_1,\tau_2) & = \lim_{n\to+\infty}\mathbf{BE}^{s-6,0}_{p,K_n}\left[\q^0,A^0 \right](\tau_1,\tau_2).
\end{align*}
Therefore, in order to prove \zcref[noname]{internal estimate pour Kerr}, it is
enough to prove
\begin{align}\label{intermediaire internal estimate pour Kerr}
    \mathbf{BE}^{s-6,0}_{p,K_n}\left[\q^0,A^0 \right](\tau_1,\tau_2) \lesssim \mathbf{E}^{s-6,0}_p\left[\q^0,A^0 \right](\tau_1),
\end{align}
uniformly in $n$. According to \zcref[noname]{limite sur Kn} and
\zcref[noname]{changement de norme} we have
\begin{align*}
    \mathbf{BE}^{s-6,0}_{p,K_n}\left[\q^0,A^0 \right](\tau_1,\tau_2) & = \lim_{\La\to 0} \mathbf{BE}^{s-6,0}_{p,K_n}\left[\q^\La,A^\La \right](\tau_1,\tau_2)
    \\& = \lim_{\La\to 0} \mathbf{BE}^{s-6,\La}_{p,K_n}\left[\q^\La,A^\La \right](\tau_1,\tau_2)
    \\& \leq \limsup_{\La\to 0} \mathbf{BE}^{s-6,\La}_{p}\left[\q^\La,A^\La \right](\tau_1,\tau_2).
\end{align*}
Applying \zcref[noname]{eq:main-q-A-combined-estimate} we get $\mathbf{BE}^{s-6,0}_{p,K_n}\left[\q^0,A^0 \right](\tau_1,\tau_2) \lesssim \limsup_{\La\to 0} \mathbf{E}^{s-6,\La}_{p}\left[\q^\La,A^\La \right](\tau_1)$ uniformly in $n$. Finally, the limit \zcref[noname]{convergence de
  l'énergie quitte à perdre du poids} concludes the proof of
\zcref[noname]{intermediaire internal estimate pour Kerr} and thus of
\zcref[noname]{internal estimate pour Kerr}. The proof of the external
counterpart \zcref[noname]{external estimate pour Kerr} is
identical. This concludes the proof of \zcref[cap]{theo comparaison}.

\appendix

\section{Proofs of \zcref[cap]{sec:adaptedglobalcoordinates}}

\subsection{Proof of \zcref[cap]{prop:properties-of-global-coordinates}}
\label{appendix:properties-of-global-coordinates}

Based on the formulas of \zcref[cap]{lem:action-global-frame} we get 
\begin{align*}
\g(\D\tau,\D\tau) & = - e_3(\tau)e_4(\tau) + (e_2(\tau))^2
 = \frac{\De}{|q|^2} (k')^2  -  \frac{(1+\ga)^2(r^2+a^2)^2}{|q|^2\De}       +  \frac{a^2(1+\ga)^2\sin^2\th}{|q|^2 \ka} .
\end{align*}
Plugging the expression of $k'$ from \zcref[cap]{definition de k et de h} this becomes $\frac{\g(\D\tau,\D\tau)}{(1+\ga)^2} +  \frac{M^2}{8r^2}   =  \pth{ 1 - \pth{1-2\chi_{\mathrm{glo}}}^2 }  F_0 + F_1 +  \frac{M^2}{8r^2}$
where
\begin{align*}
F_0 & = -\frac{(r^2+a^2)^2}{|q|^2\De}   - \frac{\De}{|q|^2}  \frac{M^4}{r^4}  + 2  \frac{r^2+a^2}{|q|^2}  \frac{M^2}{r^2} ,
& F_1 & = \frac{\De}{|q|^2}  \frac{M^4}{r^4}  -2   \frac{r^2+a^2}{|q|^2}  \frac{M^2}{r^2}     +  \frac{a^2\sin^2\th}{|q|^2 \ka} .
\end{align*}
Since $1 - \pth{1-2\chi_{\mathrm{glo}}}^2 \in [0,1]$, a sufficient condition for $\frac{\g(\D\tau,\D\tau)}{(1+\ga)^2} +  \frac{M^2}{8r^2} $ to be nonpositive is $F_1+  \frac{M^2}{8r^2}$ and $F_0 + F_1+  \frac{M^2}{8r^2} \leq 0$. We have
\begin{align*}
\frac{r^2}{M^2}\pth{F_1+  \frac{M^2}{8r^2}}  & =   \frac{1}{8} - 2   \frac{r^2+a^2}{|q|^2}  +  \frac{\De}{|q|^2}  \frac{M^2}{r^2}      + \pth{ \frac{a}{M}}^2  \frac{ r^2 \sin^2\th}{|q|^2 \ka}     ,
\\ \frac{r^2}{M^2}\pth{F_0 + F_1+  \frac{8M^2}{r^2}} & =  \frac{1}{8} - \frac{r^2}{|q|^2}\frac{(r^2+a^2)^2}{M^2\De}   +  \pth{ \frac{a}{M}}^2   \frac{r^2\sin^2\th}{|q|^2 \ka}  .
\end{align*}
Now, on $r_\HH(1-\de_\HH) \leq r \leq r_{\overline{\HH}}(1+\de_\HH)$ we have
\begin{align*}
\frac{1}{1+ \pth{\frac{a}{r_\HH(1-\de_\HH)}}^2} \leq \frac{r^2+a^2}{|q|^2} \leq 1 + \pth{ \frac{a}{r_\HH(1-\de_\HH)}}^2.
\end{align*}
Using this and also $\De\leq r^2+a^2$ we thus obtain the uniform bounds
\begin{align}
\frac{r^2}{M^2}\pth{F_1+  \frac{M^2}{8r^2}} & \leq \frac{1}{8}-\frac{2}{1+ \pth{\frac{a}{r_\HH(1-\de_\HH)}}^2} + \pth{ \frac{a}{M}}^2  +  \pth{1 + \pth{ \frac{a}{r_\HH(1-\de_\HH)}}^2} \frac{M^2}{\pth{r_\HH(1-\de_\HH)}^2}  , \label{premiere inegalite}
\\ \frac{r^2}{M^2}\pth{F_0 + F_1+  \frac{M^2}{8r^2}} & \leq \frac{1}{8} -  \frac{r_\HH^2(1-\de_\HH)^2}{M^2\pth{1+\pth{\frac{a}{r_\HH(1-\de_\HH)}}^2}}  +  \pth{ \frac{a}{M}}^2  .  \label{deuxieme inegalite}
\end{align}
Remembering the expression of $r_\HH$ we find that $\lim_{a,\de_\HH,\La\to 0} \pth{\text{RHS of both \eqref{premiere inegalite} and \eqref{deuxieme inegalite}}}  < 0$. Therefore, by choosing $a$, $\de_\HH$ and $\La$ small enough compared to $M$, we have proved that $\frac{\g(\D\tau,\D\tau)}{(1+\ga)^2} \leq -  \frac{M^2}{8r^2} $, which implies \zcref[noname]{eq:tau-foliation:properties:NSigma-uniformly-timelike} since $\ga>0$.

Now, the inequality \zcref[noname]{eq:tau-foliation:properties:NSigma-uniformly-timelike} implies that $- 8e_3(\tau)e_4(\tau) +  9(e_2(\tau))^2  \leq -\frac{M^2}{r^2} +  (e_2(\tau))^2 $. Plugging the value of $e_2(\tau)$ from \zcref[cap]{lem:action-global-frame} gives
\begin{align*}
- 8e_3(\tau)e_4(\tau) +  9(e_2(\tau))^2 & \leq  -  M^2\frac{r^2\pth{ 1+\ga\cos^2\th -  \frac{a^2}{M^2}(1+\ga)^2\sin^2\th }  + a^2\cos^2\th(1+\ga\cos^2\th) }{r^2|q|^2\ka } 
\\& \leq   \frac{ a^2(1+\ga)^2  - M^2 }{|q|^2\ka } ,
\end{align*}
which is indeed nonpositive if $a$ and $\La$ are small enough. We have proved that $| \nab \tau|^2 \leq \frac{8}{9}e_4(\tau) e_3(\tau)$. To get both $e_4(\tau)>0$ and $e_3(\tau)>0$, it is thus sufficient to prove $e_4(\tau)>0$. We have
\begin{align*}
\frac{e_4(\tau)}{ \la_{\mathrm{glo}}(1+\ga)} & =  (1-\chi_{\mathrm{glo}}) \pth{ \frac{2(r^2+a^2)}{\De} - \frac{ M^2}{r^2} } + \chi_{\mathrm{glo}}\frac{ M^2}{r^2}.
\end{align*}
If $r\leq r_0-M$, then $\De\leq r^2+a^2$ implies $\frac{e_4(\tau)}{ (1+\ga)} \geq  \frac{2(r^2+a^2)}{|q|^2} \pth{1 -  \frac{ M^2}{ 2r^2}} $ which is positive on $r_\HH(1-\de_\HH) \leq r \leq r_{\overline{\HH}}(1+\de_\HH)$. If $r\geq r_0+M$ then $\frac{e_4(\tau)}{ 1+\ga}  =  \frac{ M^2}{r^2}$ which is obviously positive. In the transition region $r_0-M\leq r\leq r_0+M$ we have $\frac{2(r^2+a^2)}{\De} - \frac{ M^2}{r^2} \geq 2 - \frac{ M^2}{r^2}>0$, implying that $\frac{e_4(\tau)}{ \la_{\mathrm{glo}}(1+\ga)} $ is an interpolation between two positive quantities, and is thus positive. We have proved $e_4(\tau)>0$ and $e_3(\tau)>0$, and thus \zcref[noname]{eq:tau-foliation:properties:ingoing}. Finally, when $r\geq r_0+M$, we have $\frac{e_4(\tau)}{1+\ga}  =  \frac{ M^2}{r^2}$ and $ \frac{e_3(\tau)}{1+\ga} = 2\frac{(r^2+a^2)}{|q|^2}  - \frac{\De}{|q|^2} \frac{M^2}{r^2}$ which indeed imply \zcref[noname]{eq:tau-foliation:properties:asymptotic-behavior}.

To conclude the proof of \zcref[cap]{prop:properties-of-global-coordinates} we only need to compute all the new metric components $\g(\D x^\mu,\D x^\nu)$ in the coordinate system $(\tau,r,\th,\ffi)$ (with the exception of $\g(\D\tau,\D\tau)$ that has already been estimated). More precisely, we only need to compute $-\half e_4(x^\mu)e_3(x^\nu)-\half e_4(x^\nu)e_3(x^\mu)$ since the horizontal part can be seen to be regular on the formulas of \zcref[cap]{lem:action-global-frame} (we can also neglect $\th$ since $e_3(\th)=e_4(\th)=0$). From the same lemma and after plugging the values of $k'$ and $h'$ from \zcref[cap]{definition de k et de h} we find
\begin{align*}
  -\half e_4(\tau)e_3(r)- \half e_4(r)e_3(\tau) & =    (1+\ga) \pth{1-2\chi_{\mathrm{glo}}}\pth{ \frac{r^2+a^2}{|q|^2} - \frac{\De M^2}{|q|^2r^2}}    ,
  \\ -\half e_4(\tau)e_3(\ffi)- \half e_4(\ffi)e_3(\tau) & = - (1-2\chi_{\mathrm{glo}})^2\frac{a(1+\ga)^2}{|q|^2}  \frac{M^2}{r^2}
 + \pth{(1-2\chi_{\mathrm{glo}})^2   - 1 }\frac{a(1+\ga)^2(r^2+a^2)}{\De|q|^2}  ,
  \\ -\half e_4(r)e_3(\ffi)-\half e_4(\ffi)e_3(r) & =     (1-2\chi_{\mathrm{glo}})\frac{a(1+\ga)}{|q|^2}  , 
  \\ -e_4(\ffi)e_3(\ffi) & = \pth{(1-2\chi_{\mathrm{glo}})^2 - 1}\frac{a^2(1+\ga)^2}{\De|q|^2}    ,
\end{align*}
which are all regular expressions on $r_\HH(1-\de_\HH)\leq r \leq r_{\overline{\HH}}(1+\de_\HH)$ since $(1-2\chi_{\mathrm{glo}})^2 - 1$ vanishes outside of $r\in[r_0-M,r_0+M]$. This concludes the proof of \zcref[cap]{prop:properties-of-global-coordinates}.

\subsection{Proof of \zcref[cap]{prop:properties-underline-t}}
\label{appendix:proof lem underline t}

From $\underline{t}=\tau + \ell(r)$ and the expressions of \zcref[cap]{lem:action-global-frame} we have
\begin{equation}\label{eq:action on underline t appendix}
\begin{aligned}
  e_4(\underline{t}) & = \lambdaglo \pth{  2\pth{1-\chi_{\mathrm{glo}}} \frac{(1+\ga)(r^2+a^2)}{\De} - (1+\ga) \pth{1-2\chi_{\mathrm{glo}}}\frac{M^2}{r^2} + \ell' },
  \\ e_3(\underline{t}) & = \lambdaglo^{-1} \pth{ \frac{(1+\ga)(r^2+a^2)}{|q|^2} - \frac{\De}{|q|^2} (1+\ga) \pth{1-2\chi_{\mathrm{glo}}}\pth{ \frac{r^2+a^2}{\De} - \frac{M^2}{r^2}} - \frac{\De}{|q|^2}\ell'},
  \\ e_1(\underline{t}) & = 0,
  \\ e_2(\underline{t}) & = \frac{a(1+\ga)\sin\th}{|q|\sqrt{\ka}},
\end{aligned}
\end{equation}
where we already used the expression of $k'$ from \zcref[cap]{definition de k et de h}. In the region $r\geq 3r_0$ where $\chiglo=\lambdaglo=\ell'=1$ we indeed get the asymptotics \zcref[noname]{eq:action-frame-underline-t}. Since $\g(\D\underline{t},\D\underline{t})=-e_4(\underline{t})e_3(\underline{t}) + \pth{e_2(\underline{t})}^2$, the asymptotics \zcref[noname]{eq:action-frame-underline-t} imply \zcref[noname]{eq:asympt-gtt}. To prove \zcref[noname]{eq:underline-t-timelike}, we deduce from the above expression of the action of the global frame on $\underline{t}$ that
\begin{equation}
  \label{eq:underline-t-timelike:aux0}
  |q|^2\pth{\g(\D\underline{t},\D\underline{t}) - \g(\D\tau,\D\tau)} =  \De\ell' \pth{ \ell' + 2(1+\ga) \pth{1-2\chi_{\mathrm{glo}}}\pth{ \frac{r^2+a^2}{\De} - \frac{M^2}{r^2}}} .
\end{equation}
If $r\leq 2r_0$ then $\underline{t}=\tau$ and we have already prove $\g(\D\tau,\D\tau)<0$. If $r>2r_0$ then $\ell'>0$ and $\chiglo=1$, so that \zcref[noname]{eq:underline-t-timelike:aux0} becomes
\begin{equation*}
  \frac{|q|^2}{\ell'(r^2+a^2)}\pth{\g(\D\underline{t},\D\underline{t}) - \g(\D\tau,\D\tau)} =     - 2(1+\ga) + \frac{\De}{r^2+a^2}\pth{\ell' + 2(1+\ga) \frac{ M^2}{r^2} }.
\end{equation*}
Using $\De\leq r^2+a^2$, $\ga>0$, $0\leq \ell'\leq 1$ and $r>2r_0$ we get that the right-hand side is less than $ - 1 + (1+\ga) \frac{ M^2}{2r_0^2} $. By taking $r_0$ large enough compared to $M$ the RHS of this inequality can be made less than $-\half$, thus proving that $\g(\D\underline{t},\D\underline{t})<0$ (since $\g(\D\tau,\D\tau)$ has already been proved to be negative). This concludes the proof of \zcref[noname]{eq:underline-t-timelike}. Finally, the above expression of the action of the global frame on $\underline{t}$ also imply that
\begin{align*}
  \g(\D\underline{t},\D r) & =  \frac{\De}{ |q|^2}\ell' +  (1+\ga) \pth{1-2\chi_{\mathrm{glo}}}\pth{ \frac{r^2+a^2}{|q|^2} - \frac{\De M^2}{|q|^2r^2}}   ,
  \\ \g(\D\underline{t},\D \th) & = 0,
  \\ \g(\D\underline{t},\D\ffi ) & = \frac{a(1+\ga)^2}{|q|^2}\pth{ - (1-2\chi_{\mathrm{glo}})^2   \frac{M^2}{r^2} + \pth{(1-2\chi_{\mathrm{glo}})^2   - 1 }\frac{(r^2+a^2)}{\De}   +  \frac{1}{ \ka}  +  \frac{1-2\chi_{\mathrm{glo}}}{1+\ga} \ell' } .
\end{align*}
Again, these are all regular expressions on $r_\HH(1-\de_\HH)\leq r \leq r_{\overline{\HH}}(1+\de_\HH)$ since $(1-2\chi_{\mathrm{glo}})^2 - 1$ vanishes outside of $r\in[r_0-M,r_0+M]$. This concludes the proof of \zcref[cap]{prop:properties-underline-t}.

\section{Proofs of \zcref[cap]{sec:derivation-of-RW}}
\label{appendix:sec:derivation-of-RW}

\subsection{Proof of \zcref[cap]{prop:Teuk:Teuk-eq-for-A}}
\label{appendix:prop:Teuk:Teuk-eq-for-A}

This appendix is devoted to the proof of \zcref[cap]{prop:Teuk:Teuk-eq-for-A}, which follows closely the one of Proposition 5.1.1 in \cite{giorgiWaveEquationsEstimates2024}. 

\subsubsection{Material for the derivation of Teukolsky}
\label{appendix:material-derivation-teuk}

Here, we consider any solution of the Einstein equation with $\La>0$ and any null frame.

\paragraph{Some null structure equations.}

\begin{lemma}
We have
\begin{align}
\,^{(c)}\nab_4 \tr\underline{X} & = - \half \tr X \tr\underline{X} + \DD\cdot \overline{\underline{H}} + \underline{H}\cdot \overline{\underline{H}} + 2\overline{P} + \frac{4\La}{3} + \underline{\Xi}\cdot\overline{\Xi} - \half \widehat{\underline{X}}\cdot \overline{\widehat{X}}, \label{cnab4trXb}
\\ \,^{(c)}\nab_4 H & = \,^{(c)}\nab_3\Xi - \half \overline{\tr X} (H-\underline{H}) - \half \widehat{X}\cdot (\overline{H} - \overline{\underline{H}}) - B, \label{cnab4H}
\\ \,^{(c)}\nab_4 \widehat{X} & = - \Re(\tr X) \widehat{X} + \half \,^{(c)}\DD\widehat{\otimes}\Xi + \half \Xi \otimes (\underline{H} + H) - A, \label{cnab4hatX}
\\ \half \overline{\,^{(c)}\DD}\cdot \widehat{X} & =  \half \,^{(c)}\DD\overline{\tr X} - i \Im(\tr X) H - i \Im(\tr\underline{X})\Xi - B. \label{codazzi}
\end{align}
\end{lemma}

\begin{proof}
By definition of the Riemann curvature tensor and the Ricci coefficients we have
\begin{align*}
\nab_4 \underline{\chi}_{ab} & = 2 \nab_b \underline{\eta}_a + 2\om \underline{\chi}_{ab} - \chi_{b}^{\; c} \underline{\chi}_{ac} + 2 \pth{ \xi_b \underline{\xi}_a + \underline{\eta}_b \underline{\eta}_a } + \R_{b34a}.
\end{align*}
We consider the trace and the anti-trace of these equations and obtain
\begin{align*}
\nab_4 \tr \underline{\chi} & = -\widehat{\chi}\c\widehat{\underline{\chi}} -\frac 1 2 \tr \chi\tr \underline{\chi}+\frac 1 2 \atrch\atrchb    +  2   \div \etab+ 2 \om \trchb + 2\big( \xi\c \xib +|\etab|^2\big)+ g^{ab}\R_{b34a},
\\ \nab_4 \atr\underline{\chi} & = -\widehat{\chi}\wedge\widehat{\underline{\chi}}-\frac 1 2(\atrch \trchb+\tr\chi\atrchb)+ 2 \curl \etab + 2 \om \atrchb + 2 \xi\wedge\xib+ \in^{ab}\R_{b34a}.
\end{align*}
Using $\R_{\mu\nu}=\La\g_{\mu\nu}$ and \zcref[noname]{link Riemann Weyl} we have $g^{ab}\R_{b34a} = 2 \rho + \frac{4\La}{3}$ and $\in^{ab}\R_{b34a} = 2 \,^*\rho$ so that
\begin{align*}
\nab_4 \tr \underline{\chi} & = -\widehat{\chi}\c\widehat{\underline{\chi}} -\frac 1 2 \tr \chi\tr \underline{\chi}+\frac 1 2 \atrch\atrchb    +  2   \div \etab+ 2 \om \trchb + 2\big( \xi\c \xib +|\etab|^2\big)+ 2 \rho + \frac{4\La}{3},
\\ \nab_4 \atr\underline{\chi} & = -\widehat{\chi}\wedge\widehat{\underline{\chi}}-\frac 1 2(\atrch \trchb+\tr\chi\atrchb)+ 2 \curl \etab + 2 \om \atrchb + 2 \xi\wedge\xib+ 2 \,^*\rho.
\end{align*}
We then pass to conformal derivatives and use complex notations to obtain \zcref[noname]{cnab4trXb}. We proceed similarly for \zcref[noname]{cnab4H}, \zcref[noname]{cnab4hatX} and \zcref[noname]{codazzi}, noting that the Riemann components appearing, namely $\R_{a443}$, $g^{ab}\R_{a44b}$ and $g^{ac}\R_{a4cb}$, coincide with the Weyl components (see \zcref[noname]{link Riemann Weyl}), so that these null structure equations can be directly taken from Proposition 2.4.13 in \cite{giorgiWaveEquationsEstimates2024}. 
\end{proof}

\paragraph{Some Bianchi equations.} 
Since the Bianchi equations are derived directly for the Weyl tensor $\W$, the Bianchi equations from \cite{giorgiWaveEquationsEstimates2024} holding for the null components of $\R$ (which coincides with $\W$ in the $\La=0$ case) also hold for the null components of $\W$ in the $\La\neq 0$ case. Therefore, we obtain from Proposition 2.4.14 in \cite{giorgiWaveEquationsEstimates2024} the following equations for $A$, $B$ and $P$:
\begin{align}
\,^{(c)}\nab_3 A - \half \,^{(c)}\DD\widehat{\otimes} B & = - \half \tr\underline{X} A + 2 H \widehat{\otimes} B - 3 \overline{P} \widehat{X}, \label{cnab3A}
\\ \,^{(c)}\nab_4 B - \half \overline{\,^{(c)}\DD}\cdot A & = - 2 \overline{\tr X} B + \half A \cdot \overline{\underline{H}} + 3 \overline{P}\Xi, \label{cnab4B}
\\ \,^{(c)}\nab_3 B - \,^{(c)}\DD \overline{P} & = - \tr\underline{X} B + \overline{\underline{B}}\cdot \widehat{X} + 3 \overline{P}H + \half A \cdot \overline{\underline{\Xi}}, \label{cnab3B}
\\ \,^{(c)}\nab_4P - \half \,^{(c)}\DD\cdot \overline{B} & = - \frac{3}{2} \tr X P + \underline{H} \cdot \overline{B} - \overline{\Xi} \cdot \underline{B} - \frac{1}{4} \widehat{\underline{X}}\cdot \overline{A}. \label{cnab4P}
\end{align}

\paragraph{One commutation formula.} For the derivation of Teukolsky we only need $\left[ \,^{(c)}\nab_4,\,^{(c)}\DD\widehat{\otimes} \right] F$ for $F\in\mathfrak{s}_1(\mathbb{C})$. According to the discussion before \zcref[cap]{lemma:commutation-formula:conformally-invariant-derivaties}, this commutator is unchanged when going to $\La>0$. Therefore, from \cite{giorgiWaveEquationsEstimates2024} we find
\begin{align}
\left[ \,^{(c)}\nab_4,\,^{(c)}\DD\widehat{\otimes} \right] F & =    -\half \tr X \pth{ \,^{(c)}\DD\widehat{\otimes} F  + (1-s)\underline{H}\widehat{\otimes} F} + \underline{H} \widehat{\otimes} \,^{(c)}\nab_4 F + \mathrm{L}(F), \label{cnab4nDDotimesF}
\end{align}
where $\mathrm{L}(\cdot)$ is a notation used to denote any linear first order operator with coefficients being at least linear in the quantities $\Xi$, $\underline{\Xi}$, $\widehat{X}$, $\widehat{\underline{X}}$, $A$, $\underline{A}$, $B$ and $\underline{B}$.

\subsubsection{Derivation of Teukolsky}

The first step of the derivation is to apply $\,^{(c)}\nab_4$ to \zcref[noname]{cnab3A}:
\begin{align*}
\,^{(c)}\nab_4\,^{(c)}\nab_3A  & = \frac{1}{2}  \,^{(c)}\nab_4\pth{ \,^{(c)}\DD\widehat{\otimes} B }  - \frac{1}{2}\,^{(c)}\nab_4\pth{\tr\underline{X} A} + 2\,^{(c)}\nab_4\pth{H\widehat{\otimes} B} - 3 \,^{(c)}\nab_4\pth{\overline{P}\widehat{X}}
\\& = \vcentcolon \mathbf{I} + \mathbf{II} + \mathbf{III} + \mathbf{IV}.
\end{align*}      
First, note that \zcref[noname]{cnab4nDDotimesF} and the fact that $B$ is $1$-conformally invariant imply
\begin{align*}
\left[ \,^{(c)}\nab_4,\,^{(c)}\DD\widehat{\otimes} \right] B & =    -\half \tr X  \,^{(c)}\DD\widehat{\otimes} B   + \underline{H} \widehat{\otimes} \,^{(c)}\nab_4 B    + \mathrm{L}(B) .
\end{align*}
Together with \zcref[noname]{cnab4B} this allows us to obtain
\begin{align*}
\mathbf{I} & =  -\frac{1}{4}  \tr X  \,^{(c)}\DD\widehat{\otimes} B +\frac{1}{2}  \underline{H} \widehat{\otimes} \,^{(c)}\nab_4 B    + \mathrm{L}(B)  +  \half  \,^{(c)}\DD\widehat{\otimes} \pth{ \half \overline{\,^{(c)}\DD}\cdot A - 2 \overline{\tr X}B + \half A \cdot \overline{\underline{H}} + 3 \overline{P}\Xi } .
\end{align*} 
Using \zcref[noname]{cnab4B} again this becomes
\begin{align*}
\mathbf{I} & =   \frac{1}{4} \,^{(c)}\DD\widehat{\otimes}\pth{ \overline{\,^{(c)}\DD}\cdot A} + \pth{\underline{H} +  \overline{\underline{H}}}\cdot \,^{(c)}\nab A   + \pth{ \half \DD\cdot \overline{\underline{H}} + \frac{1}{2} \underline{H}\cdot \overline{\underline{H}}   }A +  \pth{ -\frac{1}{4}  \tr X  -  \overline{\tr X} } \,^{(c)}\DD\widehat{\otimes} B 
\\&\quad + \pth{  -  \overline{\tr X} \underline{H}   -  \,^{(c)}\DD \overline{\tr X} } \widehat{\otimes}B  + \frac{3}{2} \overline{P}  \underline{H} \widehat{\otimes} \Xi + \frac{3}{2} \overline{P} \,^{(c)}\DD\widehat{\otimes} \Xi  + \frac{3}{2} \,^{(c)}\DD \overline{P}\widehat{\otimes} \Xi  + \mathrm{L}(B) ,
\end{align*} 
where we also used $ \underline{H} \widehat{\otimes} \pth{ \overline{\underline{H}} \cdot A }=2 \pth{\underline{H}\cdot \overline{\underline{H}} } A$, $\,^{(c)}\DD\widehat{\otimes}\pth{ \overline{\underline{H}} \cdot A} =2 \pth{ \DD\cdot \overline{\underline{H}}}A + 4 \overline{\underline{H}}\cdot \,^{(c)}\nab A$ and $\underline{H}\widehat{\otimes}\pth{ \overline{\,^{(c)}\DD}\cdot A } = 4 \underline{H}\cdot \,^{(c)}\nab A$. Using \zcref[noname]{cnab4trXb} we obtain
\begin{align*}
\mathbf{II} & =  - \frac{1}{2}\tr\underline{X}\,^{(c)}\nab_4 A   - \half \pth{ - \frac{1}{2}\tr X\tr\underline{X} +  \DD\cdot\overline{\underline{H}} +  \underline{H} \cdot \overline{\underline{H}} + 2\overline{P}  + \frac{4\Lambda}{3} +  \underline{\Xi}\cdot \overline{\Xi} - \frac{1}{2}\widehat{\underline{X}}\cdot\overline{\widehat{X}} } A.
\end{align*}
Note that this computation is the only one differing from the derivation of Teukoslky in Kerr, coming from the null structure equations for $\,^{(c)}\nab_4 \tr\underline{X}$. With the notation $\mathrm{L}$ this rewrites as
\begin{align*}
\mathbf{II} & =  - \frac{1}{2}\tr\underline{X}\,^{(c)}\nab_4 A   +  \pth{  \frac{1}{4}\tr X\tr\underline{X} - \half \DD\cdot\overline{\underline{H}} - \half  \underline{H} \cdot \overline{\underline{H}} - \overline{P}  - \frac{2\Lambda}{3} } A + \mathrm{L}(A).
\end{align*}
Using \zcref[noname]{cnab4B} and \zcref[noname]{cnab4H} we obtain
\begin{align*}
\mathbf{III} & =  H\widehat{\otimes}\pth{  \overline{\,^{(c)}\DD}\cdot A - 4 \overline{\tr X}B +  A \cdot \overline{\underline{H}} + 6 \overline{P}\Xi } +  \pth{2\,^{(c)}\nab_3\Xi  - \overline{\tr X}(H-\underline{H}) - \widehat{X}\cdot\left(\overline{H} - 2\overline{\underline{H}}\right) - B } \widehat{\otimes} B 
\\& =  4 H \cdot \,^{(c)}\nab A  +  2\pth{H\cdot \overline{\underline{H}}}A + 6\overline{P} H\widehat{\otimes}\Xi +\overline{\tr X} \pth{ \underline{H} - 5 H } \widehat{\otimes} B + \mathrm{L}(B).
\end{align*}
Finally, using \zcref[noname]{cnab4hatX} and \zcref[noname]{cnab4P} we obtain
\begin{align*}
\mathbf{IV} & = - 3  \overline{P} \pth{ -  \Re(\tr X)\widehat{X} + \frac{1}{2}\,^{(c)}\DD\widehat{\otimes}\Xi + \frac{1}{2}\Xi\widehat{\otimes}(\underline{H}  + H ) - A }
\\&\quad  - 3 \pth{  \frac{1}{2}\overline{\,^{(c)}\DD }\cdot {B} - \frac{3}{2}\overline{\tr X} \overline{P} + \overline{\underline{H}}\cdot {B} - {\Xi}\cdot\overline{\underline{B}} - \frac{1}{4}\overline{\widehat{\underline{X}}} \cdot A }  \widehat{X}
\\& = 3  \overline{P} \pth{ \frac{3}{2}\overline{\tr X}  \widehat{X}  + \Re(\tr X)\widehat{X} - \frac{1}{2}\,^{(c)}\DD\widehat{\otimes}\Xi - \frac{1}{2}\Xi\widehat{\otimes}(\underline{H}  + H ) + A } + \mathrm{L}( A ) + \mathrm{L}(B).
\end{align*}
Putting everything together, we obtain
\begin{align*}
\,^{(c)}\nab_4\,^{(c)}\nab_3A  & =  \frac{1}{4} \,^{(c)}\DD\widehat{\otimes}\pth{ \overline{\,^{(c)}\DD}\cdot A} - \frac{1}{2}\tr\underline{X}\,^{(c)}\nab_4 A + \pth{4H +\underline{H} +  \overline{\underline{H}}}\cdot \,^{(c)}\nab A 
\\&\quad  + \pth{  \frac{1}{4}\tr X\tr\underline{X}   + 2 \overline{P}  - \frac{2\Lambda}{3}  + 2\pth{H\cdot \overline{\underline{H}}} }A
\\&\quad +  \pth{ -\frac{1}{4}  \tr X  -  \overline{\tr X} } \,^{(c)}\DD\widehat{\otimes} B + \pth{ - 5 \overline{\tr X}  H    -  \,^{(c)}\DD \overline{\tr X}  } \widehat{\otimes}B 
\\&\quad  + \frac{3}{2} \,^{(c)}\DD \overline{P}\widehat{\otimes} \Xi   + \overline{P} \pth{   \frac{9}{2}\overline{\tr X}  \widehat{X}  + 3\Re(\tr X)\widehat{X}   + \frac{9}{2} H\widehat{\otimes}\Xi   }   + \mathrm{L}( A ) + \mathrm{L}(B).
\end{align*}
We now use \zcref[noname]{codazzi} to replace $\,^{(c)}\DD\overline{\tr X}$, \zcref[noname]{cnab3A} to replace $ \,^{(c)}\DD\widehat{\otimes} B$ and \zcref[noname]{cnab3B} to replace $\,^{(c)}\DD \overline{P}$:
\begin{align*}
\,^{(c)}\nab_4\,^{(c)}\nab_3A  & =  \frac{1}{4} \,^{(c)}\DD\widehat{\otimes}\pth{ \overline{\,^{(c)}\DD}\cdot A}  - \frac{1}{2}\tr\underline{X}\,^{(c)}\nab_4 A +  \pth{ -\frac{1}{2}  \tr X  - 2 \overline{\tr X} }\,^{(c)}\nab_3 A   + \pth{4H +\underline{H} +  \overline{\underline{H}}}\cdot \,^{(c)}\nab A 
\\&\quad  + \pth{    -  \overline{\tr X}\tr\underline{X}  + 2 \overline{P}  - \frac{2\Lambda}{3}  + 2 H\cdot \overline{\underline{H}} }A+ \mathrm{L}( A ) + \mathrm{L}(B) + \mathrm{L}(\underline{B}).
\end{align*}
This concludes the proof of \zcref[cap]{prop:Teuk:Teuk-eq-for-A}, since $ \mathrm{L}( A ) + \mathrm{L}(B) + \mathrm{L}(\underline{B})$ is of the form of the right hand side in \zcref[noname]{eq:teuk-with-Q}.

\subsection{Proof of \zcref[cap]{proposition T to RW}}
\label{appendix:proposition T to RW}

This appendix is devoted to the proof of \zcref[cap]{proposition T to RW}, which follows closely the one of Theorem 5.2.9 in \cite{giorgiWaveEquationsEstimates2024}. The quantity $\q$ defined by \zcref[noname]{definition of q frak} and ultimately satisfying the generalized Regge-Wheeler equation \zcref[noname]{gRW} will be first sought under the form 
\begin{align}\label{q first expression}
\q & = f \pth{   \,^{(c)}\nab_3\,^{(c)}\nab_3 A  + C_1 \,^{(c)}\nab_3 A + C_2 A },
\end{align}
for well-chosen complex-valued functions $f$, $C_1$ and $C_2$. They will be chosen in order to achieve the structural properties of \zcref[noname]{gRW}:
\begin{itemize}
\item The function $f$ will be such that the first order derivatives of $\q$ don't contain $\dr_\th$ derivatives and have a vanishing real part.
\item The function $C_1$ will be such that the potential $V$ is real and positive and such that the first order derivatives of $\,^{(c)}\nab_3A$ don't contain $\dr_\th$ or $\dr_r$ derivatives.
\item The function $C_2$ will ensure the fall-off of $W^{[4]}$, $W^{[3]}$, $W^{[h]}$ and $W^{[0]}$, as well as the factorization of \zcref[noname]{q first expression} into \zcref[noname]{definition of q frak}.
\end{itemize}

\subsubsection{Material for the derivation of Regge-Wheeler}

As opposed to the derivation of the Teukolski equation for $A$, the derivation of the Regge-Wheeler equation for $\mathfrak{q}$ is directly performed in the global principal null frame introduced in \zcref[cap]{sec:global null frame}. We can thus benefit from the properties of this frame stated in \zcref[noname]{eq:global-frame:vanishing-qtys} which simplify the null structure equations, Bianchi equations and commutation formulas required for the derivation compared to the ones stated in \zcref[cap]{appendix:material-derivation-teuk}.

\paragraph{Some null structure equations.} 

\begin{lemma}
In the global principal null frame of \zcref[cap]{def:global null frame}, the following null structure equations hold
\begin{align}
\,^{(c)}\nab_3 \tr X   & = - \half \tr\underline{X}\tr X + \DD\cdot \overline{H} + H\cdot \overline{H} + 2P + \frac{4\La}{3} , \label{cnab3trX}
\\ \,^{(c)}\nab_3\tr\underline{X} & = - \half (\tr\underline{X})^2 , \label{cnab3trXb}
\\ \,^{(c)}\nab_4 {\tr\underline{X}}  & = - \half {\tr X} {\tr\underline{X}} + {\DD} \cdot \overline{\underline{H}} + \underline{H}\cdot \overline{\underline{H}} +  2\overline{P} + \frac{4\La}{3} , \label{cnab4trXbbis}
\\ \,^{(c)}\nab_3 \underline{H} & = - \half \overline{\tr\underline{X}} \pth{ \underline{H} - H } . \label{cnab3Hb}
\end{align}
Moreover, by considering the real and imaginary parts of \zcref[noname]{cnab3trX}, \zcref[noname]{cnab3trXb} and \zcref[noname]{cnab4trXbbis} we also obtain
\begin{align}
\,^{(c)}\nab_3\tr\underline{\chi} & = - \half  \tr\underline{\chi}^2 + \half  \,^{(a)}\tr\underline{\chi}^2 , \label{cnab3trchib}
\\ \,^{(c)}\nab_3\,^{(a)}\tr\underline{\chi} & = - \tr\underline{\chi} \,^{(a)}\tr\underline{\chi} ,\label{cnab3atrchib}
\\ \,^{(c)}\nab_4 \tr\underline{\chi}  & = - \half {\tr X} {\tr\underline{X}} + 2 \div\underline{\eta}  + 2 |\eta|^2 +  2\rho  + \frac{4\La}{3} , \label{cnab4trchib}
\\ \,^{(c)}\nab_4 \atr\underline{\chi}  & =  2 \curl\underline{\eta}  + 2  \,^*\rho  , \label{cnab4atrchi}
\\ \,^{(c)}\nab_3 \tr\chi &  = - \half \tr\underline{X}\tr X + 2\div\eta + 2 |\eta|^2 + 2\rho + \frac{4\La}{3} , \label{cnab3trchi}
\end{align}
where we note that ${\tr X} {\tr\underline{X}}$ is a real-valued function\footnote{This follows from ${\tr X} {\tr\underline{X}}$ being $0$-conformally invariant and thus equal to $-\frac{4\De}{|q|^4}$.}.
\end{lemma}

\begin{proof}
By definition of the Riemann curvature tensor and the Ricci coefficients we have
\begin{align*}
\nab_3 {\chi}_{ab} & = 2 \nab_b {\eta}_a + 2\underline{\om} {\chi}_{ab} - \underline{\chi}_{b}^{\; c} {\chi}_{ac} + 2 \pth{ \underline{\xi_b} {\xi}_a + {\eta}_b {\eta}_a } + \R_{b43a}.
\end{align*}
We consider the trace and the anti-trace of these equations and obtain
\begin{align*}
\nab_3 \tr{\chi} & = -\widehat{\underline{\chi}}\c\widehat{\chi} -\frac 1 2 \trchb\tr\chi+\frac 1 2 \atrchb\atrch    +   2   \div \eta+ 2 \omb \trch + 2 \big(\xi\c \xib +|\eta|^2\big)+ g^{ab}\R_{b43a},
\\ \nab_3 \atr{\chi} & = -\widehat{\underline{\chi}}\wedge\widehat{\chi}-\frac 1 2(\atrchb \tr\chi+\trchb\atrch)+ 2 \curl \eta + 2 \omb \atrch + 2 \xib\wedge\xi  + \in^{ab} \R_{b43a}.
\end{align*}
Using $\R_{\mu\nu}=\La\g_{\mu\nu}$ and \zcref[noname]{link Riemann Weyl} we have $g^{ab}\R_{b43a} = 2 \rho + \frac{4\La}{3}$ and $\in^{ab} \R_{b43a} = - 2\,^*\rho$. Plugging this into the above null structure equations and using also the simplifications from \zcref[noname]{eq:global-frame:vanishing-qtys} we thus obtain
\begin{align*}
\nab_3 \tr{\chi} & =  -\frac 1 2 \trchb\tr\chi+\frac 1 2 \atrchb\atrch    +   2   \div \eta+ 2 \omb \trch + 2 |\eta|^2 + 2 \rho + \frac{4\La}{3},
\\ \nab_3 \atr{\chi} & = -\frac 1 2(\atrchb \tr\chi+\trchb\atrch)+ 2 \curl \eta + 2 \omb \atrch  - 2\,^*\rho.
\end{align*}
We then pass to conformal derivatives and use complex notations to obtain \zcref[noname]{cnab3trX}. We proceed similarly for \zcref[noname]{cnab3trXb} and \zcref[noname]{cnab3Hb}, noting that the Riemann components appearing, namely $g^{ab}\R_{a33b}$ and $\R_{a334}$, coincide with the Weyl components (see \zcref[noname]{link Riemann Weyl}), so that these null structure equations can be directly taken from Proposition 2.4.13 in \cite{giorgiWaveEquationsEstimates2024}. Note also that \zcref[noname]{cnab4trXbbis} follows from \zcref[noname]{cnab4trXb} and \zcref[noname]{eq:global-frame:vanishing-qtys}.
\end{proof}

We now state relations that are not consequence of the null structure equations but rather exact computations on Kerr-de Sitter. They are all based on the following properties of $q$:
\begin{equation}\label{properties of q appendix}
\begin{aligned}
\nab_4q & = \half \tr X q, & \nab_3q & = \half \overline{\tr\underline{X}}q,
& \nab q & = \half \pth{\underline{H} + \overline{H}}q, & \DD q & = q \underline{H}.
\end{aligned}
\end{equation}
Using the third identity in \zcref[noname]{properties of q appendix} and the value of $\tr\underline{X}$ in either the outgoing or ingoing principal null frames we can also obtain the following relation, valid in the global principal null frame:
\begin{align}
\,^{(c)}\nab\tr \underline{X} & = - \pth{ H + \half \overline{\underline{H}} + \half \overline{H} } \tr\underline{X}, 
\\ 2\,^{(c)}\nab\tr \underline{\chi} & = -  \tr\underline{\chi} \pth{   \underline{\eta} + 3 \eta }  -  \atr\underline{\chi} \pth{  \,^*\eta -  \,^*\underline{\eta} } , \label{cnabtrchib}
\\ 2 \,^{(c)}\nab\atr \underline{\chi} & = -  \atr\underline{\chi} \pth{  \underline{\eta} + 3 \eta }    -  \tr\underline{\chi} \pth{ \,^*\underline{\eta} - \,^*\eta    }.  \label{cnabatrchib}
\end{align}
From the fourth property of \zcref[noname]{properties of q appendix} we also obtain
\begin{align}
\,^{(c)}\DD \overline{\tr\underline{X}} & = -  \overline{\tr\underline{X}}(H + \underline{H}).     \label{cDtrXb}
\end{align}
Moreover, from $H = \frac{a}{\bar{q}}\mathfrak{J}$, the second property of \zcref[noname]{properties of q appendix} and the equation satisfied by $\mathfrak{J}$ from \zcref[cap]{lemma:canonical-complex-one-form-j-basic-identities} we obtain
\begin{align}
\,^{(c)} \nab_3H & = - \tr\underline{X} H. \label{cnab3H}
\end{align}
We can also obtain the decay
\begin{align}
\,^{(c)}\DD\cdot \overline{H} = \,^{(c)}\DD\cdot \overline{\underline{H}} = \GO{ a r^{-3}}. \label{cnabHHb} 
\end{align}

\paragraph{One Bianchi equation.}
As explained above, the Bianchi identities match the ones from \cite{giorgiWaveEquationsEstimates2024}. Therefore, from Proposition 2.4.14 in \cite{giorgiWaveEquationsEstimates2024} and \zcref[noname]{eq:global-frame:vanishing-qtys} we get the following equation in the global principal null frame:
\begin{align}
\,^{(c)}\nab_3 P & = - \frac{3}{2} \overline{\tr\underline{X}} P, \label{cnab3P}
\\ \,^{(c)}\nab_3 \rho & = - \frac{3}{2} \pth{ \tr\underline{\chi}\rho - \,^{(a)}\tr\chi \,^*\rho} . \label{cnab3rho}
\end{align}

\paragraph{Some commutation formulas.} We also need some commutation formulas for conformal derivatives applied to $F\in\mathfrak{s}_1(\mathbb{C})$ and $U\in\mathfrak{s}_2(\mathbb{C})$ (all assumed to be $s$-conformally invariant). As explained in \zcref[cap]{lemma:commutation-formula:conformally-invariant-derivaties}, the following commutation formulas are the same as in the $\La=0$ case, therefore, Lemma 4.2.2 in \cite{giorgiWaveEquationsEstimates2024} gives
\begin{align}
\left[  \,^{(c)}\nab_3 , \,^{(c)}\DD \widehat{\otimes}   \right] F & = - \half \tr\underline{X} \pth{ \,^{(c)}\DD \widehat{\otimes} F + (s+1)H\widehat{\otimes}F } + H \widehat{\otimes}\,^{(c)}\nab_3 F , \label{cDDotimescnab_3 F} 
\\ \left[ \,^{(c)}\nab_3 , F\cdot \,^{(c)}\nab \right] U & = \pth{ \,^{(c)}\nab_3 F  - \half \tr\underline{\chi}  F  + \half \,^{(a)}\tr\underline{\chi} \,^*F }\cdot  \,^{(c)}\nab U + (F\cdot  \eta) \,^{(c)}\nab_3 U  \label{cnab_3Fcdotcnab U}
\\&\quad + \GO{ F \tr\underline{X} H U } . \nonumber
\end{align}

\subsubsection{Commutators involving the Teukolsky operator}\label{section teukolsky operator}

To ease the computations, we decompose the Teukolsky operator as
\begin{align}\label{L(A)=loool}
\mathcal{L}(A) & =  - \,^{(c)}\nab_4 \,^{(c)}\nab_3 A + \frac{1}{4} \,^{(c)}\DD \widehat{\otimes} \pth{ \overline{ \,^{(c)}\DD}\cdot A } + \mathcal{L}^{(1)}(A) + \mathcal{L}^{(0)}(A),
\end{align}
where $\mathcal{L}^{(1)}(A)  = T^{[3]}  \,^{(c)}\nab_3 A + T^{[4]} \,^{(c)}\nab_4 A + T^{[h]}  \cdot \,^{(c)}\nab A$ and $\mathcal{L}^{(0)}(A)  = T^{[0]} A$ with coefficients
\begin{align*}
T^{[3]} & =   -\half \tr X - 2 \overline{\tr X}, & T^{[4]} & = - \half \tr\underline{X} ,
\\ T^{[h]} & =  4H + \underline{H} + \overline{\underline{H}} , & T^{[0]} & =  - \overline{\tr X} \tr \underline{X} + 2\overline{P} - \frac{2\La}{3} + 2 H\cdot \overline{\underline{H}}.
\end{align*}
Since the global principal null frame coincides with the outgoing principal null frame for large radius, the fall-off properties that will be stated in this section can all be deduced from the various expressions of the Ricci and curvature coefficients in the outgoing principal null frame given in \zcref[cap]{lemma:Kerr:outgoing-PG:Ric-and-curvature}. In particular, we get that
\begin{align*}
T^{[3]} & = - \frac{5}{2}\tr\chi + \GO{a r^{-2}},& T^{[4]} & = - \half \tr\underline{\chi} + \GO{a r^{-2}},
\\ T^{[h]} & = \GO{a r^{-2}},& T^{[0]} & = - \tr\chi \tr \underline{\chi} + 2\rho - \frac{2\La}{3} + \GO{a r^{-3}}.
\end{align*}

\begin{lemma}\label{lemma Lf}
Let $U\in\mathfrak{s}_2(\CCC)$ and $f$ a scalar function. We have
\begin{align*}
[\mathcal{L},f](U) & = E^{[3]}(f)  \,^{(c)}\nab_3 U  + E^{[4]}(f)  \,^{(c)}\nab_4 U + E^{[h]}(f)\cdot \,^{(c)}\nab U+ E^{[0]}(f) U,
\end{align*}
where
\begin{align*}
E^{[3]}(f) & = - \,^{(c)}\nab_4 f , & E^{[4]}(f) & = -  \,^{(c)}\nab_3 f ,
\\ E^{[h]}(f) & = 2  \,^{(c)}\nab f , & E^{[0]}(f) & = T^{[3]}  \,^{(c)}\nab_3 f  + T^{[4]} \,^{(c)}\nab_4 f + T^{[h]}  \cdot \,^{(c)}\nab f - \,^{(c)}\nab_4   \,^{(c)}\nab_3 f +  \half \,^{(c)}\DD \cdot \overline{ \,^{(c)}\DD} f.
\end{align*}
\end{lemma}

\begin{proof}
Straightforward computations using the Leibniz rule.
\end{proof}

\begin{lemma}\label{lemma Lnab3}
Let $U\in\mathfrak{s}_2(\CCC)$ be $s$-conformally invariant. We have
\begin{align*}
\left[ \mathcal{L},\,^{(c)}\nab_3\right](U) & = F^{[h,h]} \,^{(c)}\DD \widehat{\otimes} \pth{  \overline{\,^{(c)}\DD}\cdot U   }  + F^{[h,3]}\cdot \,^{(c)} \nab \,^{(c)}\nab_3 U  
\\&\quad + F^{[h,s]}\cdot \,^{(c)}\nab   U  +  F^{[3,s]} \,^{(c)}\nab_3 U  + F^{[4]} \,^{(c)}\nab_4 U + F^{[0,s]} U ,
\end{align*}
where
\begin{align*}
F^{[h,h]} & =  \frac{1}{4}  \tr\underline{\chi} ,
\\ F^{[h,3]} & = - 2\underline{\eta} ,
\\ F^{[3,s]} & =  - (T^{[h]}\cdot  \eta)-  \frac{1}{2}    \pth{  \,^{(c)}\DD \cdot  \overline{H} }  -   |\eta|^2- \,^{(c)}\nab_3( T^{[3]}) + 2(s-1) \pth{ \rho - \frac{\La}{3} - \eta\cdot \underline{\eta} }  + 4 \ImagUnit \pth{ - \,^*\rho + \eta\wedge\underline{\eta}} ,
\\ F^{[h,s]} & = -\,^{(c)}\nab_3 T^{[h]}  + \half \tr\underline{\chi}  T^{[h]}  - \half \,^{(a)}\tr\underline{\chi} \,^*T^{[h]} -2T^{[4]}  (\eta-\underline{\eta}) 
\\&\quad + \frac{1}{2}   \,^{(c)}\DD\pth{ \overline{\tr \underline{X}} } + \frac{1}{2}  \overline{\tr \underline{X}}   H+ \frac{s+1}{2} \tr\underline{X}H  + \frac{ s-2}{2}    \overline{\tr \underline{X}}\overline{H} ,
\\ F^{[4]} & = - \,^{(c)}\nab_3(T^{[4]}) ,
\\ F^{[0,s]} & =  - \,^{(c)}\nab_3(T^{[0]})+ \frac{ s-2}{4}     \,^{(c)}\DD  \cdot \pth{ \overline{\tr \underline{X}}\overline{H}}     + \frac{s-2}{2}  \overline{\tr \underline{X}} |\eta|^2 
\\&\quad -T^{[4]} \pth{  2s \pth{ \rho - \frac{\La}{3} - \eta\cdot \underline{\eta} }  + 4 \ImagUnit \pth{ - \,^*\rho + \eta\wedge\underline{\eta}}}   + \GO{ T^{[h]} \tr\underline{X} H }  .
\end{align*}
\end{lemma}

\begin{proof}
We have
\begin{align*}
\left[ \mathcal{L},\,^{(c)}\nab_3\right](U)  & = \left[ \,^{(c)}\nab_3 ,\,^{(c)}\nab_4\right](\,^{(c)}\nab_3 U)   + \frac{1}{4} \left[ \,^{(c)}\DD \widehat{\otimes} \pth{ \overline{ \,^{(c)}\DD}\cdot  },\,^{(c)}\nab_3\right](U) 
\\&\quad +  T^{[4]} \left[  \,^{(c)}\nab_4,\,^{(c)}\nab_3\right](U)   + \left[ T^{[h]}  \cdot \,^{(c)}\nab,\,^{(c)}\nab_3\right](U) 
\\&\quad - \,^{(c)}\nab_3( T^{[3]}) \,^{(c)}\nab_3 U - \,^{(c)}\nab_3(T^{[4]}) \,^{(c)}\nab_4 U - \,^{(c)}\nab_3(T^{[0]}) U.
\end{align*}
We compute each terms one by one, starting with the two terms involving $\left[ \,^{(c)}\nab_3 ,\,^{(c)}\nab_4\right]$. Since $\,^{(c)}\nab_3 U$ and $U$ are respectively $(s-1)$- and $s$-conformally invariant, we obtain from \zcref[noname]{eq:commutation-formula:conformally-invariant-derivaties:tensor:3-4} that
\begin{align*}
\left[ \,^{(c)}\nab_3 ,\,^{(c)}\nab_4\right](\,^{(c)}\nab_3 U)  & = 2(\eta-\underline{\eta})\cdot \,^{(c)} \nab \,^{(c)}\nab_3 U 
\\&\quad + \pth{ 2(s-1) \pth{ \rho - \frac{\La}{3} - \eta\cdot \underline{\eta} }  + 4 \ImagUnit\pth{ - \,^*\rho + \eta\wedge\underline{\eta}}}\,^{(c)}\nab_3 U,
\\ T^{[4]} \left[  \,^{(c)}\nab_4,\,^{(c)}\nab_3\right](U)  & = -2T^{[4]}  (\eta-\underline{\eta})\cdot \,^{(c)} \nab U  
 -T^{[4]} \pth{  2s \pth{ \rho - \frac{\La}{3} - \eta\cdot \underline{\eta} }  + 4 \ImagUnit \pth{ - \,^*\rho + \eta\wedge\underline{\eta}}}U.
\end{align*}
We now look at the terms involving horizontal derivatives. First, using that $\overline{ \,^{(c)}\DD}\cdot U$ is $s$-conformally invariant we obtain from \zcref{cDDotimescnab_3 F} and \zcref[noname]{eq:eq:commutation-formula:conformally-invariant-derivaties:tensor:3-div} that
\begin{align*}
&\frac{1}{4} \left[ \,^{(c)}\DD \widehat{\otimes} \pth{ \overline{ \,^{(c)}\DD}\cdot  },\,^{(c)}\nab_3\right](U) 
\\&\qquad = \frac{1}{4} \,^{(c)}\DD \widehat{\otimes} \pth{ \left[   \overline{ \,^{(c)}\DD}\cdot  ,\,^{(c)}\nab_3\right](U) } + \frac{1}{4} \left[ \,^{(c)}\DD \widehat{\otimes} ,\,^{(c)}\nab_3\right]\pth{\overline{ \,^{(c)}\DD}\cdot U}
\\&\qquad  = - \frac{1}{4} \,^{(c)}\DD \widehat{\otimes} \pth{ - \half \overline{\tr \underline{X}} \pth{ \overline{\,^{(c)}\DD}\cdot U + (s-2) \overline{H}\cdot U } + \overline{H} \cdot \,^{(c)} \nab_3 U  } 
\\&\qquad\quad + \frac{1}{4} \pth{  \half \tr\underline{X} \pth{ \,^{(c)}\DD \widehat{\otimes} \pth{\overline{ \,^{(c)}\DD}\cdot U} + (s+1)H\widehat{\otimes}\pth{\overline{ \,^{(c)}\DD}\cdot U} } - H \widehat{\otimes}\,^{(c)}\nab_3 \pth{\overline{ \,^{(c)}\DD}\cdot U} }.
\end{align*}
Expanding using the Leibniz rule we obtain
\begin{equation}\label{gros commutateur}
\begin{aligned}
\frac{1}{4} \left[ \,^{(c)}\DD \widehat{\otimes} \pth{ \overline{ \,^{(c)}\DD}\cdot  },\,^{(c)}\nab_3\right](U) 
& =  \frac{1}{4}  \tr\underline{\chi}  \,^{(c)}\DD \widehat{\otimes} \pth{  \overline{\,^{(c)}\DD}\cdot U   }   -  2\eta \cdot  \,^{(c)}\nab   \,^{(c)} \nab_3 U
\\&\quad + \pth{ \frac{1}{2}   \,^{(c)}\DD\pth{ \overline{\tr \underline{X}} } + \frac{1}{2}  \overline{\tr \underline{X}}   H+ \frac{s+1}{2} \tr\underline{X}H  + \frac{ s-2}{2}    \overline{\tr \underline{X}}\overline{H} } \cdot \,^{(c)}\nab   U 
\\&\quad + \pth{  -  \frac{1}{2}    \pth{  \,^{(c)}\DD \cdot  \overline{H} }  -   |\eta|^2 } \,^{(c)} \nab_3 U 
\\&\quad + \pth{  \frac{ s-2}{4}    \pth{ \,^{(c)}\DD  \cdot \pth{ \overline{\tr \underline{X}}\overline{H}} }    + \frac{s-2}{2}  \overline{\tr \underline{X}} |\eta|^2 } U.
\end{aligned}
\end{equation}
Finally, we use \zcref[noname]{cnab_3Fcdotcnab U} to get
\begin{align*}
\left[ T^{[h]}  \cdot \,^{(c)}\nab,\,^{(c)}\nab_3\right](U)  & = \pth{ -\,^{(c)}\nab_3 T^{[h]}  + \half \tr\underline{\chi}  T^{[h]}  - \half \,^{(a)}\tr\underline{\chi} \,^*T^{[h]} }\cdot  \,^{(c)}\nab U
\\&\quad - (T^{[h]}\cdot  \eta) \,^{(c)}\nab_3 U  + \GO{ T^{[h]} \tr\underline{X} H U }  .
\end{align*}
Putting everything together and collecting each derivatives of $U$ concludes the proof.
\end{proof}

From \zcref[cap]{lemma:Kerr:outgoing-PG:Ric-and-curvature}, \zcref[noname]{cnab3trX}, \zcref[noname]{cnab3trXb}, \zcref[noname]{cnab3Hb}, \zcref[noname]{cDtrXb}, \zcref[noname]{cnab3H}, \zcref[noname]{cnabHHb} and \zcref[noname]{cnab3P} we obtain that
\begin{align*}
F^{[h,3]} & = \GO{ar^{-2}},
\\ F^{[3,s]} & =  - \frac{5}{4}  \tr\underline{\chi}\tr \chi  +   5 \pth{  \rho + \frac{2\La}{3} } + 2(s-1) \pth{ \rho - \frac{\La}{3} }   + \GO{a r^{-3}},
\\ F^{[h,s]} & =  \GO{ a r^{-3}},
\\ F^{[4]} & = - \frac{1}{4} \tr\underline{\chi}^2 + \GO{a r^{-3}}, 
\\ F^{[0,s]} & = -  \tr \chi  \tr\underline{\chi}^2   +\tr \underline{\chi} \pth{  5\rho + \frac{4\La}{3}  + s    \pth{ \rho - \frac{\La}{3} }  }    + \GO{ a r^{-4}} .
\end{align*}

\subsubsection{Link with the wave operator}

The following lemma connects the Teukolsky operator $\mathcal{L}$ and the wave operator $\dot{\Box}_2$.

\begin{lemma}
Let $U\in\mathfrak{s}_2(\CCC)$ which is $0$-conformally invariant. We have
\begin{equation}\label{Box in terms of L}
\begin{aligned}
\dot{\Box}_2 U & = \mathcal{L}(U)  + \pth{ - T^{[4]}  - \half \tr\underline{X}}\,^{(c)}\nab_4 U +\pth{ -  T^{[3]} - \half \tr X }  \,^{(c)}\nab_3 U + \pth{2 \underline{\eta} - T^{[h]} }\cdot \,^{(c)}\nab U  
\\&\quad  - \pth{ \frac{1}{4} \tr X \overline{\tr\underline{X}} + \frac{1}{4} \tr\underline{X}\overline{\tr X} +   2 \overline{P}  - \frac{2\La}{3} + 2\ImagUnit \eta\wedge\underline{\eta} + T^{[0]} }U .
\end{aligned}
\end{equation}
\end{lemma}

\begin{proof}
Recall from \zcref[noname]{eq:wave-operator-null-frame-decomp:wave-2} that
\begin{align*}
\dot{\Box}_2 U & = - \,^{(c)}\nab_4 \,^{(c)}\nab_3 U  - \half \tr\underline{\chi} \,^{(c)}\nab_4 U  - \half \tr \chi \,^{(c)}\nab_3U + \Laplace_2U + 2 \underline{\eta}\cdot \,^{(c)}\nab U + 2 \ImagUnit (\,^*\rho - \eta\wedge\underline{\eta})U,
\end{align*}
where we used  $ \,^{(c)}\nab_4 \,^{(c)}\nab_3 U  = \nab_4 \nab_3 U  - 2\om \,^{(c)}\nab_3 U$. Now, as in Appendix C.8 of \cite{giorgiWaveEquationsEstimates2024} we find
\begin{align*}
\Laplace_2 U  = \frac{1}{4} \DD\widehat{\otimes} \pth{ \overline{\DD}\cdot U} +  \ImagUnit[\nab_1,\nab_2] U
& = \frac{1}{4} \,^{(c)}\DD\widehat{\otimes} \pth{ \overline{\,^{(c)}\DD}\cdot U} +   \half \pth{ \ImagUnit\,^{(a)}\tr\chi\,^{(c)}\nab_3 U + \ImagUnit\,^{(a)}\tr\underline{\chi} \,^{(c)}\nab_4 U} 
\\&\quad +  \pth{ - \frac{1}{4}  \tr X \overline{\tr\underline{X}} - \frac{1}{4} \tr\underline{X}\overline{\tr X}   - 2\rho + \frac{2\La}{3} } U ,
\end{align*}
where we also used \zcref[noname]{eq:gauss-eqn:k-1-2}, \zcref[noname]{eq:h-K:def} and \zcref[noname]{link Riemann Weyl}. This gives
\begin{align*}
\dot{\Box}_2 U & = - \,^{(c)}\nab_4 \,^{(c)}\nab_3 U +   \frac{1}{4} \,^{(c)}\DD\widehat{\otimes} \pth{ \overline{\,^{(c)}\DD}\cdot U}   - \half \tr\underline{X} \,^{(c)}\nab_4 U  - \half \tr X \,^{(c)}\nab_3U  + 2 \underline{\eta}\cdot \,^{(c)}\nab U 
\\&\quad -  \pth{  \frac{1}{4}  \tr X \overline{\tr\underline{X}} + \frac{1}{4} \tr\underline{X}\overline{\tr X}   + 2\overline{P} - \frac{2\La}{3} + 2 \ImagUnit  \eta\wedge\underline{\eta} } U .
\end{align*} 
Using \zcref[noname]{L(A)=loool} concludes the proof of the lemma.
\end{proof}

\subsubsection{The intermediate quantity $\tilde{\q}$}

As stated in \zcref[cap]{proposition T to RW}, we fix $A\in \mathfrak{s}_2(\mathbb{C})$ $2$-conformally invariant and satisfying $\mathcal{L}(A)=0$ in the global frame. To ease the presentation, we define the intermediate quantity
\begin{align}\label{definition of tilde q frak}
\tilde{\q} & = f \pth{\,^{(c)}\nab_3\,^{(c)}\nab_3 A  + C_1 \,^{(c)}\nab_3 A},
\end{align}
for $f$ and $C_1$ two functions to be determined that are respectively $0$- and $-1$-conformally invariant, thus making $\tilde{\q} $ $0$-conformally invariant. In the following proposition, we derive a wave equation for $\tilde{\q}$, making use of the commutators previously computed. Note that the absence of $\,^{(c)}\nab_a \,^{(c)}\nab_4 A$ and $ \,^{(c)}\nab_4 \,^{(c)}\nab_4 A$ in \zcref[noname]{eq tilde q frak} is due to the vanishing of  $\underline{\xi}$ in the global frame. 

\begin{proposition}\label{prop q tilde}
The quantity $\tilde{\q}$ defined by \zcref[noname]{definition of tilde q frak} satisfies
\begin{equation}\label{eq tilde q frak}
\begin{aligned}
\dot{\Box}_2\tilde{\q} & = \nab_{Y}\tilde{\q} + \widehat{V} \tilde{\q} + f \RR^{[4,3]} \,^{(c)}\nab_4\,^{(c)}\nab_3 A  + f \RR^{[h,3]} \cdot \,^{(c)}\nab\,^{(c)}\nab_3 A
\\&\quad + \RR^{[3]} \,^{(c)}\nab_3 A  + \RR^{[4]} \,^{(c)}\nab_4 A  + \RR^{[h]} \cdot \,^{(c)}\nab A  + \RR^{[0]}  A,
\end{aligned}
\end{equation}
where
\begin{itemize}
\item the vector field $Y$ is given by
\begin{align}\label{vector field Y}
Y & =  \pth{f^{-1} E^{[4]}(f) +8F^{[h,h]} - T^{[4]}  - \half \tr\underline{X}}e_4  +\pth{f^{-1}E^{[3]}(f)  -  T^{[3]} - \half \tr X }e_3 
\\&\quad + \pth{2 \underline{\eta} - T^{[h]} +  f^{-1}E^{[h]}(f) + 2F^{[h,3]}}^a e_a   , \nonumber
\end{align}
\item the potential $\widehat{V}$ is given by
\begin{equation}\label{fonction V tilde}
\begin{aligned}
\widehat{V}  & =  -  \frac{1}{4} \tr X \overline{\tr\underline{X}} - \frac{1}{4} \tr\underline{X}\overline{\tr X} -   2 \overline{P}  + \frac{2\La}{3} - 2i  \eta\wedge\underline{\eta} - T^{[0]}  + F^{[3,s=2]} +F^{[3,s=1]}  + (F^{[h,3]}\cdot  \eta) 
\\&\quad -8 F^{[h,h]} T^{[3]} - f^{-2} E^{[3]}(f)  \nab_3 \pth{ f } - f^{-2} E^{[4]}(f)  \nab_4 \pth{  f  }   - f^{-2} E^{[h]}(f)\cdot \nab \pth{  f } 
\\&\quad   - 2 f^{-1} F^{[h,3]}\cdot  \nab \pth{f} -  8f^{-1}F^{[h,h]} \nab_4 \pth{  f } + f^{-1} E^{[0]}(f) + f^{-1} E^{[3]}(f  C_1)  -  f^{-1} C_1 E^{[3]}(f) ,
\end{aligned}
\end{equation}
\item the second order coefficients $\RR^{[4,3]}$ and $\RR^{[h,3]}$ are given by
\begin{align}
\RR^{[4,3]} & =  E^{[4]}(C_1)  + 2 F^{[4]}   + 4 F^{[h,h]} \pth{ \tr\underline{\chi}    -2T^{[4]} - C_1}   + 4  \,^{(c)}\nab_3\pth{ F^{[h,h]} } ,  \label{coeff R 43}
\\ \RR^{[h,3]} & = E^{[h]}(C_1) +   F^{[h,s=1]} + F^{[h,s=2]}   +  8  F^{[h,h]} (\eta-2\underline{\eta} - T^{[h]} )   \label{coeff R h3}
\\&\quad    +  \,^{(c)}\nab_3 F^{[h,3]}  - \pth{  C_1  + \half \tr\underline{\chi}  }F^{[h,3]}  + \half \,^{(a)}\tr\underline{\chi} \,^*F^{[h,3]} , \nonumber 
\end{align}
\item the first and zero-th order coefficients $ \RR^{[3]}$, $ \RR^{[4]}$, $ \RR^{[h]}$ and $\RR^{[0]}$ are given by
\begin{align}
f^{-1}\RR^{[3]} & =   - f^{-1}E^{[3]}(f)  \,^{(c)}\nab_3 \pth{  C_1}   - C_1  E^{[3]}( C_1 )   - f^{-1} E^{[4]}(f)  \,^{(c)}\nab_4 \pth{   C_1} - f^{-1} E^{[h]}(f)\cdot \,^{(c)}\nab \pth{  C_1}   \nonumber
\\&\quad    - 4  F^{[h,h]}  \pth{ 2\,^{(c)}\nab_3 T^{[3]} + \tr\underline{\chi} T^{[3]}  + 2T^{[0]} - T^{[3]}   C_1 + 2 \,^{(c)}\nab_4 \pth{  C_1 }  - 4 \pth{ \rho - \frac{\La}{3}  } }  \label{RR3}
\\&\quad  +  F^{[0,s=1]}    +   F^{[0,s=2]}     +   \,^{(c)}\nab_3\pth{    F^{[3,s=2]} }      - 4  \,^{(c)}\nab_3\pth{     F^{[h,h]}  } T^{[3]}      \nonumber
\\&\quad    - 2  F^{[h,3]}\cdot \,^{(c)} \nab \pth{  C_1}   -  C_1 (F^{[h,3]}\cdot  \eta)   \nonumber
\\&\quad  -  F^{[3,s=1]}   C_1  - f^{-1}E^{[0]}(f) C_1 + f^{-1}E^{[0]}(f C_1 )  +    \GO{a r^{-4}}  , \nonumber
\\ f^{-1}\RR^{[4]} & =    2 F^{[h,h]}  \tr\underline{X} \pth{    \tr\underline{\chi}   - \half \tr\underline{X} }  +   \,^{(c)}\nab_3\pth{   F^{[4]} }  - 4    \,^{(c)}\nab_3\pth{   F^{[h,h]} T^{[4]} }  \label{RR4}
\\&\quad +   C_1 \pth{ F^{[4]}  - 4   F^{[h,h]} T^{[4]} } , \nonumber
\\ f^{-1} \RR^{[h]} & =     C_1  \GO{ a r^{-3}} + \GO{a r^{-4}}, \label{RRh}
\\ f^{-1}\RR^{[0]} & =  - 4  F^{[h,h]}   \pth{  \tr\underline{\chi} T^{[0]}   + \,^{(c)}\nab_3\pth{ T^{[0]}} + 8T^{[4]}   \pth{ \rho - \frac{\La}{3}  }   }  \nonumber
\\&\quad  +   \,^{(c)}\nab_3\pth{ F^{[0,s=2]}}  - 4    \,^{(c)}\nab_3\pth{      F^{[h,h]} T^{[0]} } + 4F^{[4]}   \pth{ \rho - \frac{\La}{3} } \label{RR0}
\\&\quad + C_1 \pth{ F^{[0,s=2]} - 4   F^{[h,h]} T^{[0]} }  +  \GO{ar^{-6}} . \nonumber
\end{align}
\end{itemize}
\end{proposition}

\begin{proof}
Since $\tilde{\q}$ is $0$-conformally invariant we have from \zcref[noname]{Box in terms of L} that
\begin{align*}
\dot{\Box}_2 \tilde{\q} & = \mathcal{L}(\tilde{\q})  + \pth{ - T^{[4]}  - \half \tr\underline{X}}\nab_4 \tilde{\q} +\pth{ -  T^{[3]} - \half \tr X }  \nab_3 \tilde{\q} + \pth{2 \underline{\eta} - T^{[h]} }\cdot \nab \tilde{\q}  
\\&\quad  - \pth{ \frac{1}{4} \tr X \overline{\tr\underline{X}} + \frac{1}{4} \tr\underline{X}\overline{\tr X} +   2 \overline{P}  - \frac{2\La}{3} + 2i  \eta\wedge\underline{\eta} + T^{[0]} }\tilde{\q} .
\end{align*}
Looking at the definition of $\tilde{\q}$ we have
\begin{align*}
\mathcal{L}(\tilde{\q})  & = \left[\mathcal{L}, f   \right] (\,^{(c)}\nab_3\,^{(c)}\nab_3A)   + f \left[\mathcal{L},  \,^{(c)}\nab_3  \right] (\,^{(c)}\nab_3A)  + f  \,^{(c)}\nab_3\pth{\left[\mathcal{L}, \,^{(c)}\nab_3  \right] (A)  }
\\&\quad  + f  C_1  \left[ \mathcal{L},  \,^{(c)}\nab_3 \right](A)   + \left[ \mathcal{L},  f  C_1  \right](\,^{(c)}\nab_3 A) 
\\& = \vcentcolon \mathbf{I} + \mathbf{II} + \mathbf{III} + \mathbf{IV} + \mathbf{V}.
\end{align*}
We compute each terms separately using the decompositions of \zcref[cap]{lemma Lf} and \zcref[cap]{lemma Lnab3}. We will also use the following relations, consequence of the definition of $\tilde{\q}$ and of $\mathcal{L}(A)=0$:
\begin{align}
\,^{(c)}\nab_3\,^{(c)}\nab_3 A   & = f^{-1}\tilde{\q} - C_1 \,^{(c)}\nab_3 A , \label{consequence def q tilde}
\\    \,^{(c)}\DD \widehat{\otimes} \pth{ \overline{ \,^{(c)}\DD}\cdot A } & = 4 \,^{(c)}\nab_4 \,^{(c)}\nab_3 A  - 4  T^{[3]}  \,^{(c)}\nab_3 A - 4 T^{[4]} \,^{(c)}\nab_4 A  \label{consequence LA=0}
  - 4 T^{[h]}  \cdot \,^{(c)}\nab A   - 4 T^{[0]} A . 
\end{align}
For $\mathbf{I}$ we directly use \zcref[cap]{lemma Lf} and \zcref[noname]{consequence def q tilde}:
\begin{align*}
\mathbf{I} & =  f^{-1}E^{[3]}(f)  \nab_3 \tilde{\q}  + f^{-1} E^{[4]}(f)  \nab_4 \tilde{\q} +  f^{-1}E^{[h]}(f)\cdot \nab \tilde{\q} 
\\&\quad + \pth{ - f^{-2} E^{[3]}(f)  \nab_3 \pth{ f } -  f^{-1} C_1 E^{[3]}(f) - f^{-2} E^{[4]}(f)  \nab_4 \pth{  f  }   - f^{-2} E^{[h]}(f)\cdot \nab \pth{  f } + f^{-1} E^{[0]}(f)  } \tilde{\q}
\\&\quad - C_1 E^{[4]}(f)  \,^{(c)}\nab_4 \,^{(c)}\nab_3 A - C_1 E^{[h]}(f)\cdot \,^{(c)}\nab  \,^{(c)}\nab_3 A
\\&\quad +\pth{  - E^{[3]}(f)  \,^{(c)}\nab_3 \pth{  C_1}   + C_1^2 E^{[3]}(f)  - E^{[4]}(f)  \,^{(c)}\nab_4 \pth{   C_1} - E^{[h]}(f)\cdot \,^{(c)}\nab \pth{  C_1}  - E^{[0]}(f) C_1 } \,^{(c)}\nab_3 A. 
\end{align*}
For $\mathbf{II}$, we apply \zcref[cap]{lemma Lnab3} with $s=1$ (since $\,^{(c)}\nab_3A$ is $1$-conformally invariant) and use \zcref[noname]{consequence def q tilde}:
\begin{align*}
\mathbf{II} & = fF^{[h,h]} \,^{(c)}\DD \widehat{\otimes} \pth{  \overline{\,^{(c)}\DD}\cdot \,^{(c)}\nab_3A   }  + F^{[h,3]}\cdot  \nab \tilde{\q}     + \pth{  F^{[3,s=1]}  - f^{-1} F^{[h,3]}\cdot  \nab \pth{f} } \tilde{\q}  
\\&\quad + \pth{ fF^{[h,s=1]}   - fC_1F^{[h,3]} }\cdot \,^{(c)} \nab  \,^{(c)}\nab_3 A   + fF^{[4]} \,^{(c)}\nab_4 \,^{(c)}\nab_3A 
\\&\quad  + \pth{ fF^{[0,s=1]}  - fF^{[h,3]}\cdot \,^{(c)} \nab \pth{  C_1 }  -  fF^{[3,s=1]}   C_1 } \,^{(c)}\nab_3 A  .
\end{align*}
We now expand the first line using \zcref[noname]{gros commutateur} and the various fall-off properties listed above to obtain
\begin{align*}
\mathbf{II} & =  fF^{[h,h]} \,^{(c)}\nab_3\pth{ \,^{(c)}\DD \widehat{\otimes} \pth{  \overline{\,^{(c)}\DD}\cdot A   }  }  +  fF^{[h,h]}   \tr\underline{\chi}  \,^{(c)}\DD \widehat{\otimes} \pth{  \overline{\,^{(c)}\DD}\cdot A   } 
\\&\quad  + F^{[h,3]}\cdot  \nab \tilde{\q}     + \pth{  F^{[3,s=1]}  - f^{-1} F^{[h,3]}\cdot  \nab \pth{f} } \tilde{\q}  
\\&\quad + \pth{ fF^{[h,s=1]}   - fC_1F^{[h,3]} -  fF^{[h,h]}  8\eta }\cdot \,^{(c)} \nab  \,^{(c)}\nab_3 A    + fF^{[4]} \,^{(c)}\nab_4 \,^{(c)}\nab_3A 
\\&\quad  + f \pth{ F^{[0,s=1]}  - F^{[h,3]}\cdot \,^{(c)} \nab \pth{  C_1 }  -  F^{[3,s=1]}   C_1   +   \GO{a r^{-4}} } \,^{(c)} \nab_3 A  +  f \GO{a r^{-4}} \cdot \,^{(c)}\nab   A .
\end{align*}
We now use \zcref[noname]{consequence LA=0} to replace the two $  \,^{(c)}\DD \widehat{\otimes} \pth{ \overline{ \,^{(c)}\DD}\cdot A }$ in the first line:
\begin{align*}
\mathbf{II} & = 4  fF^{[h,h]} \,^{(c)}\nab_3 \,^{(c)}\nab_4 \,^{(c)}\nab_3 A   -4  fF^{[h,h]} \,^{(c)}\nab_3\pth{ T^{[3]}  \,^{(c)}\nab_3 A  }  -4 fF^{[h,h]} \,^{(c)}\nab_3\pth{  T^{[4]} \,^{(c)}\nab_4 A  } 
\\&\quad -4 fF^{[h,h]} \,^{(c)}\nab_3\pth{  T^{[h]}  \cdot \,^{(c)}\nab A   }  -4 fF^{[h,h]} \,^{(c)}\nab_3\pth{ T^{[0]} A } 
\\&\quad  + F^{[h,3]}\cdot  \nab \tilde{\q}     + \pth{  F^{[3,s=1]}  - f^{-1} F^{[h,3]}\cdot  \nab \pth{f} } \tilde{\q}  
\\&\quad + f \pth{  4  F^{[h,h]}   \tr\underline{\chi} + F^{[4]} } \,^{(c)}\nab_4 \,^{(c)}\nab_3A  + f \pth{ F^{[h,s=1]}   - C_1F^{[h,3]} -  F^{[h,h]}  8\eta }\cdot \,^{(c)} \nab  \,^{(c)}\nab_3 A   
\\&\quad  + f \pth{ F^{[0,s=1]}  - F^{[h,3]}\cdot \,^{(c)} \nab \pth{  C_1 }  -  F^{[3,s=1]}   C_1 - 4  F^{[h,h]}   \tr\underline{\chi} T^{[3]}    +   \GO{a r^{-4}} } \,^{(c)} \nab_3 A   
\\&\quad - 4  fF^{[h,h]}   \tr\underline{\chi} T^{[4]} \,^{(c)}\nab_4 A +  f \GO{a r^{-4}} \cdot \,^{(c)}\nab   A     - 4  fF^{[h,h]}   \tr\underline{\chi} T^{[0]} A .
\end{align*}
Using now \zcref[noname]{eq:commutation-formula:conformally-invariant-derivaties:tensor:3-4} and \zcref[noname]{cnab_3Fcdotcnab U} we finally obtain
\begin{align*}
\mathbf{II} & =   4  F^{[h,h]} \nab_4 \tilde{\q}   + F^{[h,3]}\cdot  \nab \tilde{\q}    
\\&\quad + \pth{  F^{[3,s=1]}  - f^{-1} F^{[h,3]}\cdot  \nab \pth{f}  -4  F^{[h,h]} T^{[3]}  -  4  f^{-1}F^{[h,h]} \nab_4 ( f )} \tilde{\q}  
\\&\quad + f \pth{  4  F^{[h,h]}   \tr\underline{\chi} + F^{[4]}  -4 F^{[h,h]}  T^{[4]}  - 4  F^{[h,h]} C_1 } \,^{(c)}\nab_4 \,^{(c)}\nab_3A     
\\&\quad + f \pth{ F^{[h,s=1]}   - C_1F^{[h,3]} -  8 F^{[h,h]}  \underline{\eta}  -4 F^{[h,h]} T^{[h]}  }\cdot \,^{(c)} \nab  \,^{(c)}\nab_3 A   
\\&\quad  + f \Bigg(  - 4  F^{[h,h]}  \pth{\,^{(c)}\nab_3 T^{[3]} + \tr\underline{\chi} T^{[3]}  + T^{[0]} - T^{[3]}   C_1 + \,^{(c)}\nab_4 \pth{  C_1 }  - 2 \pth{ \rho - \frac{\La}{3}  } } 
\\&\hspace{5cm}  +  F^{[0,s=1]}  - F^{[h,3]}\cdot \,^{(c)} \nab \pth{  C_1 }  -  F^{[3,s=1]}   C_1   +   \GO{a r^{-4}}  \Bigg)  \,^{(c)}\nab_3 A
\\&\quad + 2 fF^{[h,h]}  \tr\underline{X} \pth{    \tr\underline{\chi}   - \half \tr\underline{X} } \,^{(c)}\nab_4 A  +  f \GO{a r^{-4}} \cdot \,^{(c)}\nab   A  
\\&\quad - 4  fF^{[h,h]}   \pth{  \tr\underline{\chi} T^{[0]}   + \,^{(c)}\nab_3\pth{ T^{[0]}} + 4T^{[4]}   \pth{ \rho - \frac{\La}{3}  }  + \GO{a r^{-5} } } A  .
\end{align*}
For $\mathbf{III}$, we first use \zcref[cap]{lemma Lnab3} and \zcref[noname]{consequence LA=0} to obtain:
\begin{align*}
\mathbf{III} & = 4 f   F^{[h,h]} \,^{(c)}\nab_3  \,^{(c)}\nab_4 \,^{(c)}\nab_3 A     + 4 f  \,^{(c)}\nab_3\pth{    F^{[h,h]} } \,^{(c)}\nab_4 \,^{(c)}\nab_3 A    
\\&\quad  - 4  f  \,^{(c)}\nab_3\pth{    F^{[h,h]} T^{[h]}  \cdot \,^{(c)}\nab A   }  + f  \,^{(c)}\nab_3\pth{ F^{[h,3]}\cdot \,^{(c)} \nab \,^{(c)}\nab_3 A }
\\&\quad + f  \,^{(c)}\nab_3\pth{  F^{[h,s=2]}\cdot \,^{(c)}\nab   A   } + f  \pth{  F^{[3,s=2]}    - 4    F^{[h,h]} T^{[3]} } \,^{(c)}\nab_3  \,^{(c)}\nab_3 A      
\\&\quad + f  \pth{ F^{[4]}   - 4   F^{[h,h]} T^{[4]} } \,^{(c)}\nab_3 \,^{(c)}\nab_4 A    
\\&\quad  + f \pth{ \,^{(c)}\nab_3\pth{   F^{[4]} }  - 4    \,^{(c)}\nab_3\pth{   F^{[h,h]} T^{[4]} } } \,^{(c)}\nab_4 A  
\\&\quad + f \pth{ F^{[0,s=2]}   +   \,^{(c)}\nab_3\pth{    F^{[3,s=2]} }   - 4    F^{[h,h]} T^{[0]}   - 4  \,^{(c)}\nab_3\pth{     F^{[h,h]} T^{[3]} }  } \,^{(c)}\nab_3 A    
\\&\quad + f \pth{  \,^{(c)}\nab_3\pth{ F^{[0,s=2]}}  - 4    \,^{(c)}\nab_3\pth{      F^{[h,h]} T^{[0]} } } A    .
\end{align*}
We now use \zcref[noname]{eq:commutation-formula:conformally-invariant-derivaties:tensor:3-4}, \zcref[noname]{cnab_3Fcdotcnab U} and \zcref[noname]{consequence def q tilde} and obtain
\begin{align*}
\mathbf{III} & =  4    F^{[h,h]} \nab_4 \tilde{\q} +   F^{[h,3]}\cdot  \nab \tilde{\q}
\\&\quad +  \pth{  F^{[3,s=2]}    - 4    F^{[h,h]} T^{[3]}    -  4 f^{-1}   F^{[h,h]} \,^{(c)}\nab_4 (f) - f^{-1}   F^{[h,3]}\cdot \,^{(c)} \nab \pth{ f } + (F^{[h,3]}\cdot  \eta)  } \tilde{\q}   
\\&\quad + f \pth{  F^{[4]}   - 4   F^{[h,h]} T^{[4]}    + 4   \,^{(c)}\nab_3\pth{    F^{[h,h]} }   - 4 C_1 F^{[h,h]} } \,^{(c)}\nab_4 \,^{(c)}\nab_3 A   
\\&\quad + f \Bigg( 8  F^{[h,h]} (\eta-\underline{\eta})  +  F^{[h,s=2]} -  C_1 F^{[h,3]}  -  4  F^{[h,h]} T^{[h]} 
\\&\quad\hspace{3cm}  +  \,^{(c)}\nab_3 F^{[h,3]}  - \half \tr\underline{\chi}  F^{[h,3]}  + \half \,^{(a)}\tr\underline{\chi} \,^*F^{[h,3]} \Bigg)\cdot  \,^{(c)}\nab \,^{(c)}\nab_3 A 
\\&\quad  + f \pth{ \,^{(c)}\nab_3\pth{   F^{[4]} }  - 4    \,^{(c)}\nab_3\pth{   F^{[h,h]} T^{[4]} } } \,^{(c)}\nab_4 A  +  f \GO{a r^{-4}}\cdot \,^{(c)} \nab A 
\\&\quad + f \Bigg(  F^{[0,s=2]}   +   \,^{(c)}\nab_3\pth{    F^{[3,s=2]} }   - 4    F^{[h,h]} T^{[0]}   - 4  \,^{(c)}\nab_3\pth{     F^{[h,h]} T^{[3]} }   
\\&\quad \hspace{2cm} -  C_1  \pth{  F^{[3,s=2]}    - 4    F^{[h,h]} T^{[3]} }       +  8   F^{[h,h]}   \pth{ \rho - \frac{\La}{3} }   
\\&\quad \hspace{2cm} - 4   F^{[h,h]} \,^{(c)}\nab_4 \pth{ C_1 }  -   F^{[h,3]}\cdot \,^{(c)} \nab \pth{  C_1}   -  C_1 (F^{[h,3]}\cdot  \eta)  +    \GO{a r^{-5}}  \Bigg) \,^{(c)}\nab_3 A 
\\&\quad + f \pth{  \,^{(c)}\nab_3\pth{ F^{[0,s=2]}}  - 4    \,^{(c)}\nab_3\pth{      F^{[h,h]} T^{[0]} } + \pth{ F^{[4]}   - 4   F^{[h,h]} T^{[4]} }   \pth{ 4\rho - \frac{4\La}{3} } +  \GO{ar^{-6}} } A.    
\end{align*}
For $\mathbf{IV}$, we use \zcref[cap]{lemma Lnab3} and \zcref[noname]{consequence LA=0} to get
\begin{align*}
\mathbf{IV} & =   4 f  C_1  F^{[h,h]} \,^{(c)}\nab_4 \,^{(c)}\nab_3 A   +  f  C_1 F^{[h,3]}\cdot \,^{(c)} \nab \,^{(c)}\nab_3 A 
\\&\quad   +  f  C_1  \GO{ a r^{-3}}  \cdot \,^{(c)}\nab A    +  f  C_1 \pth{ F^{[3,s=2]}   - 4   F^{[h,h]} T^{[3]} }  \,^{(c)}\nab_3 A
\\&\quad +  f  C_1 \pth{ F^{[4]}  - 4   F^{[h,h]} T^{[4]} } \,^{(c)}\nab_4 A + fC_1 \pth{ F^{[0,s=2]} - 4   F^{[h,h]} T^{[0]} } A.
\end{align*}
Finally, for $\mathbf{V}$, we use \zcref[cap]{lemma Lf} and \zcref[noname]{consequence def q tilde} to get
\begin{align*}
\mathbf{V} & =  f^{-1} E^{[3]}(f C_1 )  \tilde{\q}  
\\&\quad  + E^{[4]}(f C_1 )  \,^{(c)}\nab_4 \,^{(c)}\nab_3 A  + E^{[h]}(f C_1 )\cdot \,^{(c)}\nab \,^{(c)}\nab_3 A + \pth{ E^{[0]}(f C_1 )  - C_1 E^{[3]}(f C_1 )  } \,^{(c)}\nab_3 A .
\end{align*}
Putting everything together concludes the proof.
\end{proof}

\subsubsection{Choice of $f$}

In this section, we choose $f$ in \zcref[noname]{definition of tilde q frak} in order to cancel the $\dr_\th$ derivative in the vector field $Y$ in \zcref[noname]{eq tilde q frak}, as well as to cancel its real part.

\begin{lemma}\label{lem choice of f}
If we set $f=q\bar{q}^3$ then the vector field $Y   =  \ImagUnit\frac{4a(1+\ga)\cos\th}{|q|^2} \T$ in \zcref[noname]{eq tilde q frak}.
\end{lemma}

\begin{proof}
From \zcref[noname]{vector field Y} we see that $Y = Y^4 e_4 + Y^3 e_3 + \,^{(h)}Y^a e_a$ with
\begin{align*}
Y^4 & =  f^{-1} E^{[4]}(f)   - \half \tr\underline{X} +8F^{[h,h]} - T^{[4]},& Y^3 & = f^{-1}E^{[3]}(f)   - \half \tr X   -  T^{[3]},
\\ \,^{(h)}Y & =   f^{-1}E^{[h]}(f) +  2 \underline{\eta}   + 2F^{[h,3]}    - T^{[h]}.
\end{align*}
Plugging the values of $F^{[h,h]}$, $F^{[h,3]}$, $T^{[4]}$, $T^{[3]}$, $T^{[h]}$, $E^{[4]}(f)$, $E^{[3]}(f)$ and $E^{[h]}(f)$ from \zcref[cap]{section teukolsky operator}, \zcref[cap]{lemma Lf} and \zcref[cap]{lemma Lnab3} this becomes
\begin{align*}
Y^4 & = - f^{-1} \nab_3f + 2 \tr\underline{\chi},  & Y^3 & = - f^{-1}\nab_4f  + 2 \overline{\tr X}  ,& \,^{(h)}Y & = 2 f^{-1}\nab f   -4 \underline{\eta} - 4H .
\end{align*}
We search $f$ under the form $f=q^n\bar{q}^m$ for $n,m\in\mathbb{Z}$ and with \zcref[noname]{properties of q appendix} we obtain
\begin{equation}\label{equations for f}
\begin{aligned}
f^{-1}\nab_3 f & = \frac{n}{2} \overline{\tr\underline{X}} + \frac{m}{2} \tr\underline{X},& f^{-1}\nab_4 f & = \frac{n}{2} \tr X + \frac{m}{2} \overline{\tr X},
& f^{-1} \nab f & = \frac{m}{2} H + \frac{n}{2} \overline{H} + \frac{n}{2} \underline{H} + \frac{m}{2} \overline{\underline{H}},
\end{aligned}
\end{equation}
and thus
\begin{align*}
Y^4 & = - \half (m+n-4) \tr\underline{\chi} + \ImagUnit \frac{m-n}{2}  \,^{(a)}\tr\underline{\chi} ,
\\ Y^3 & = -\half (m+n-4) \tr\chi + \ImagUnit\frac{4-m+n}{2} \,^{(a)}\tr\chi,
\\ \,^{(h)}Y & = (m+n-4) (\eta + \underline{\eta})    + \ImagUnit \pth{ (m-n-4)\,^*\eta    + ( n - m) \,^*\underline{\eta}   }.
\end{align*}
Choosing $m=4-n$ implies that the real part of these quantities vanish. We are left with
\begin{align*}
Y^4 & =  \ImagUnit (2-n)  \,^{(a)}\tr\underline{\chi} ,
& Y^3 & = \ImagUnit n \,^{(a)}\tr\chi ,
& \,^{(h)}Y & =  \ImagUnit \pth{ -2n\,^*\eta    + ( 2n - 4) \,^*\underline{\eta}   }.
\end{align*} 
As it can be seen in \zcref[noname]{eq:global-frame:non-vanishing-qtys}, we have $(\,^*\eta +\,^*\underline{\eta})_1=(\eta +\underline{\eta})_2=0$. To benefit from this cancellation in $\,^{(h)}Y$ we want $n$ to satisfy $-2n=2n-4$ i.e $n=1$. We thus obtain $Y =  \ImagUnit\pth{ \,^{(a)}\tr\underline{\chi} e_4 +  \,^{(a)}\tr\chi e_3 -2 \pth{ \,^*\eta_2 + \,^*\underline{\eta}_2   } e_2 }$. This vector field is 0-conformally invariant so that we can compute it using the outgoing principal null frame. From \zcref[cap]{lemma:Kerr:outgoing-PG:Ric-and-curvature} we obtain $Y =  \ImagUnit\frac{2 a \cos\th}{|q|^2}  \pth{ \frac{\De}{|q|^2} e_4^{\operatorname{(out)}} +  e_3^{\operatorname{(out)}} - 2 \frac{a\sqrt{\ka}\sin\th}{|q|} e_2 }$ which concludes the proof of the lemma after having used the first relation of \zcref[noname]{eq:T-R-Z:principal-outgoing:expression:T-Z}.
\end{proof}

\subsubsection{Choice of $C_1$}

In terms of $A$, the potential term $\widehat{V}$ and the second order derivatives coefficients $\mathcal{R}^{[4,3]}$ and $\mathcal{R}^{[h,3]}$ in \zcref[noname]{eq tilde q frak} are at the same level since $\tilde{\q}$ contains $\,^{(c)}\nab_3\,^{(c)}\nab_3$. The function $C_1$ will be chosen to achieve structural properties of these terms. First, we compute them more explicitly, plugging in particular our previous choice of $f=q\bar{q}^3$.

\begin{lemma}\label{lem intermédiaire V R43 Rh3}
If $f=q\bar{q}^3$ then 
\begin{align}
\widehat{V} & =   -  \,^{(c)}\nab_4( C_1)    + 4\rho  - \half \pth{ 3  +  \Re \pth{ \frac{\bar{q}}{q} } } \tr X \tr\underline{X} +  8 \div\eta+2 \eta\cdot\underline{\eta}  +  5 |\eta|^2 \label{inter Vtilde}
 - 8\ImagUnit\,^*\rho+ \frac{14\La}{3}   + 8 \ImagUnit \div\,^*\eta , 
\\ \RR^{[4,3]} & =    -\,^{(c)}\nab_3 (C_1)    -  \tr\underline{\chi} C_1 +  \left| \tr\underline{X} \right|^2 , \label{inter R43}
\\ \RR^{[h,3]} & = 2\,^{(c)}\nab(C_1)  +  (\underline{H} + \overline{\underline{H}}) C_1   + 2 \overline{\tr\underline{X}} H    + 2  \tr\underline{X} \overline{H}   + 8  \tr\underline{X}H   . \label{inter Rh3}
\end{align}
\end{lemma}

\begin{proof}
We start by $\widehat{V}$. By plugging the values of $F^{[h,3]}$ and $F^{[h,h]}$ from \zcref[cap]{lemma Lnab3} in \zcref[noname]{fonction V tilde} we obtain
\begin{align*}
\widehat{V} & =  -  \frac{1}{4} \tr X \overline{\tr\underline{X}} - \frac{1}{4} \tr\underline{X}\overline{\tr X} -   2 \overline{P}  + \frac{2\La}{3} - 2\ImagUnit \eta\wedge\underline{\eta} - T^{[0]} 
 + F^{[3,s=2]} +F^{[3,s=1]}   - 2 \eta\cdot\underline{\eta}    - 2 \tr\underline{\chi} T^{[3]}
\\&\quad  - f^{-2} E^{[3]}(f)  \nab_3 \pth{ f } - f^{-2} E^{[4]}(f)  \nab_4 \pth{  f  }   - f^{-2} E^{[h]}(f)\cdot \nab \pth{  f } 
\\&\quad + f^{-1} E^{[3]}(f  C_1)  -  f^{-1} C_1 E^{[3]}(f)  + 4 f^{-1} \underline{\eta}\cdot  \nab \pth{f} -  2 \tr\underline{\chi} f^{-1} \nab_4 \pth{  f } + f^{-1} E^{[0]}(f) .
\end{align*}
We now use the expression of $E$'s operators from \zcref[cap]{lemma Lf} together with the formulas \zcref[noname]{equations for f} with $n=1$ and $m=3$ and their consequence $ f^{-1} \overline{\DD}(f)  =   \overline{H}  + 3 \overline{\underline{H}} $ as well as \zcref[noname]{cnab4trXbbis}. We obtain
\begin{align*}
\widehat{V}  & =   -  \,^{(c)}\nab_4( C_1)  - 2{P} + 6  \overline{P} + \frac{14\La}{3}   -   \frac{3}{2} {\tr\underline{X}}   {\tr X}   - \frac{1}{4} \overline{\tr X}  \tr\underline{X}  - \frac{1}{4} \tr X \overline{\tr\underline{X}}  
\\&\quad     + \half H   \cdot    \overline{\underline{H}}    + \half  \underline{H} \cdot    \overline{H}     +  \frac{5}{2} {H} \cdot   \overline{{H}} + 4\overline{\DD}\cdot {H} +  \half  \DD \cdot  \overline{H}      - \half \overline{\DD} \cdot {\underline{H}} ,
\end{align*}
where we also used the values of $F^{[3,s]}$ from \zcref[cap]{lemma Lnab3}. We then use $\overline{\DD}\cdot\underline{H}=\DD\cdot \overline{H}$ and $\frac{\overline{\tr\underline{X}}}{\tr\underline{X}}=\frac{\bar{q}}{q}$ to obtain the desired expression of $\widehat{V}$. For $ \RR^{[4,3]}$, we start by plugging the expressions of $E^{[4]}$ and $T^{[4]}$ in \zcref[noname]{coeff R 43} to obtain 
\begin{align*}
\RR^{[4,3]} & = -\,^{(c)}\nab_3(C_1)  + 2 F^{[4]}   + 4 F^{[h,h]} \pth{ \frac{3}{2} \tr\underline{X} + \half \overline{\tr\underline{X}}   - C_1}   + 4  \,^{(c)}\nab_3\pth{ F^{[h,h]} } .
\end{align*}
We now plug the expressions of $F^{[4]}$ and $F^{[h,h]}$ and use \zcref[noname]{cnab3trXb} to obtain the desired expression of $\RR^{[4,3]}$. We proceed similarly for $\RR^{[h,3]}$ given by \zcref[noname]{coeff R h3}, using \zcref[noname]{cnab3Hb} as well as $\tr\underline{X}\underline{H} = - \overline{\tr\underline{X}}H$.
\end{proof}

By writing $\widehat{V} = \widetilde{V} + \widehat{V} - \widetilde{V}$ in \zcref[noname]{eq tilde q frak} for $\widetilde{V}$ yet to be defined, we obtain 
\begin{align*}
\dot{\Box}_2\tilde{\q} -  \ImagUnit\frac{4a(1+\ga)\cos\th}{|q|^2} \nab_{\T}\tilde{\q} & =  \widetilde{V} \tilde{\q}  + f  \,^{(c)}\nab_\VFArthur  \,^{(c)}\nab_3 A + \pth{  (\widehat{V}-\widetilde{V} ) f C_1 + \RR^{[3]} } \,^{(c)}\nab_3 A  
\\&\quad + \RR^{[4]} \,^{(c)}\nab_4 A  + \RR^{[h]} \cdot \,^{(c)}\nab A  + \RR^{[0]}  A ,
\end{align*}
where $\VFArthur = (\widehat{V}-\widetilde{V} ) e_3  +  \RR^{[4,3]} e_4  +  \RR^{[h,3]}_1 e_1 + \RR^{[h,3]}_2 e_2$. By using the expression of the global frame we get
\begin{align*}
\VFArthur &  = \pth{ (\widehat{V}-\widetilde{V} ) \lambdaglo^{-1}  \frac{(1+\ga)(r^2+a^2)}{|q|^2}  +  \RR^{[4,3]} \lambdaglo   \frac{(1+\ga)(r^2+a^2)}{\De} + \RR^{[h,3]}_2 \frac{(1+\ga)}{\sqrt{\ka}}  \frac{a\sin\th}{|q|}} \dr_t  
\\&\quad  + \pth{ - (\widehat{V}-\widetilde{V} ) \lambdaglo^{-1}  \frac{\De}{|q|^2}  +  \RR^{[4,3]} \lambdaglo } \dr_r 
\\&\quad  + \pth{ (\widehat{V}-\widetilde{V} ) \lambdaglo^{-1}  \frac{(1+\ga)a}{|q|^2} +  \RR^{[4,3]} \lambdaglo  \frac{(1+\ga)a}{\De} +  \RR^{[h,3]}_2 \frac{(1+\ga)}{\sqrt{\ka}} \frac{1}{|q|\sin\th} } \dr_\phi  +  \RR^{[h,3]}_1 \frac{\sqrt{\ka}}{|q|}\dr_\th .
\end{align*}
Cancelling the $\dr_r$ derivative foreces us to define $\widetilde{V}$ by $ \widetilde{V}  \vcentcolon = \widehat{V}- \lambdaglo^2 \frac{|q|^2}{\De}  \RR^{[4,3]} $. This quantity appears as the potential in the above equation for $\tilde{\q}$ and we will choose $\Im(C_1)$ to ensure $\Im(\widetilde{V})=0$. We will then choose $\Re(C_1)$ to cancel the $\dr_\th$ derivative in $\VFArthur$, i.e. to ensure $ \RR^{[h,3]}_1=0$.

\begin{lemma}\label{lem potentiel réel}
If $C_1=2\tr\underline{\chi} + \tilde{C}_1$ with $\Im(\tilde{C}_1)=-4\,^{(a)}\tr\underline{\chi}$ then
\begin{align*}
\widetilde{V} & =     \,^{(c)}\nab_{\R}( \Re(\tilde{C}_1))  -\tr\chi \Re(\tilde{C}_1) +   \frac{4 \De r^2 }{|q|^6} + 2\La   +  4 \div\eta +   |\eta|^2  + 2\eta\cdot\underline{\eta},
\end{align*}
where $\R\vcentcolon= -e_4 + \lambdaglo^2 \frac{|q|^2}{\De} e_3$. In particular we have $\Im\pth{ \widetilde{V} }=0$.
\end{lemma}

\begin{proof}
From \zcref[noname]{inter Vtilde} and \zcref[noname]{inter R43} we obtain
\begin{align*}
\widetilde{V}  & =   \pth{\,^{(c)}\nab_\R  - \tr\chi} C_1   + 4\rho  - \half \pth{ 1  +  \Re \pth{ \frac{\bar{q}}{q} } } \tr X \tr\underline{X}   +  8 \div\eta+2 \eta\cdot\underline{\eta}  +  5 |\eta|^2 - 8\ImagUnit\,^*\rho+ \frac{14\La}{3}   + 8 \ImagUnit \div\,^*\eta,
\end{align*}
with $\R$ defined in the lemma and where we used $ \overline{\tr\underline{X}} = -\lambdaglo^{-2} \frac{\De}{|q|^2}\tr X$. Plugging now $C_1=2\tr\underline{\chi} + \tilde{C}_1$ and using \zcref[noname]{cnab3trXb} and \zcref[noname]{cnab4trXbbis} this becomes
\begin{align*}
\widetilde{V}   & =     \,^{(c)}\nab_{\R}( \Re(\tilde{C}_1))  -\tr\chi \Re(\tilde{C}_1) - \half \pth{ 1 +  \Re \pth{  \frac{\bar{q}}{q}  } } \tr X \tr\underline{X} + 2\La  
  +  4 \div\eta +   |\eta|^2  + 2\eta\cdot\underline{\eta}
\\&\quad  + \ImagUnit \pth{    \,^{(c)}\nab_{\R}( \Im(\tilde{C}_1))  -\tr\chi \Im(\tilde{C}_1)     - 8\,^*\rho  + 8\div\,^*\eta  } .
\end{align*}
Using again the structure equations for $\tr\underline{X}$ we see that choosing $\Im(\tilde{C}_1)=-4\,^{(a)}\tr\underline{\chi}$ cancels the imaginary part of the above expression. Using then $\Re\pth{\frac{\bar{q}}{q}}=\frac{r^2-a^2\cos^2\th}{|q|^2}$ and $\tr X \tr\underline{X} = - \frac{4 \De }{|q|^4}$ concludes the proof of the lemma.
\end{proof}

In the next lemma, we choose $\Re(\tilde{C}_1)$ to cancel the $e_1$ component of $\mathcal{R}^{[h,3]}$. We also give the values of the $e_2$ component, $\widetilde{V} $ and $\RR^{[4,3]}$ with this choice of $\Re(\tilde{C}_1)$.

\begin{lemma}
If $C_1=2\tr\underline{\chi} + \tilde{C}_1$ with $\Im(\tilde{C}_1)=-4\,^{(a)}\tr\underline{\chi}$ and $\Re(\tilde{C}_1)= -2 \frac{\,^{(a)}\tr\underline{\chi}^2}{\tr\underline{\chi}}$ then
\begin{align}
\mathcal{R}^{[h,3]}_1 & = 0.  \label{Rh3_1=0}
\end{align}
Moreover we have
\begin{align}
\mathcal{R}^{[h,3]}_2 & = -  \lambdaglo^{-1} \frac{8a\De \sqrt{\ka}\sin\th}{|q|^5}  ,  \label{Rh3_2}
\\ \widehat{V} - \lambdaglo^2 \frac{|q|^2}{\De}  \RR^{[4,3]}  & =      \frac{4 \De r^2 }{|q|^6} + 2\La   +  4 \div\eta +   |\eta|^2  + 2\eta\cdot\underline{\eta}  +   \frac{8\,^{(a)}\tr\underline{\chi}}{\tr\underline{\chi}} \pth{  \curl\underline{\eta}  +   \,^*\rho  } \nonumber
\\&\quad  +  \frac{\,^{(a)}\tr\underline{\chi}^2}{\tr\underline{\chi}^2} \pth{   {\tr X} {\tr\underline{X}} - 4 \div\underline{\eta}  - 4 |\eta|^2 -  4\rho  - \frac{8\La}{3} } \label{Vtilde-R43 final}
\\&\quad + \lambdaglo^2 \frac{|q|^2}{\De}  \,^{(a)}\tr\underline{\chi}^2  \pth{  3 + \frac{  \,^{(a)}\tr\underline{\chi}^2}{\tr\underline{\chi}^2} }  +2\tr\chi  \frac{\,^{(a)}\tr\underline{\chi}^2}{\tr\underline{\chi}}, \nonumber
\\ \RR^{[4,3]} & = - \lambdaglo^{-2} \frac{4a^2 \De^2 \cos^2\th}{r^2|q|^6}     . \label{R43 final}
\end{align}
\end{lemma}

\begin{proof}
We use \zcref[noname]{cnabtrchib} and $C_1=2\tr\underline{\chi} + \tilde{C}_1$ to deduce from \zcref[noname]{inter Rh3} that
\begin{align*}
\RR^{[h,3]} & =    2\,^{(c)}\nab(\tilde{C}_1)  +  2 \tilde{C}_1 \underline{\eta}    + 6  \tr\underline{\chi}  \eta  +  2 \atr\underline{\chi}   \,^*\eta + 2   \tr\underline{\chi}   \underline{\eta} +2  \atr\underline{\chi} \,^*\underline{\eta}   + 8\ImagUnit \pth{ \tr\underline{\chi}\,^*\eta - \atr\underline{\chi} \eta }.
\end{align*}
Using now \zcref[noname]{cnabatrchib} and $\Im(\tilde{C}_1)=-4\,^{(a)}\tr\underline{\chi}$ we can compute the terms depending on $\tilde{C}_1$:
\begin{align*}
2\,^{(c)}\nab(\tilde{C}_1)  + 2  \tilde{C}_1 \underline{\eta}  & = 2\,^{(c)}\nab(\Re(\tilde{C}_1))  + 2 \Re( \tilde{C}_1)\underline{\eta} + 4  \ImagUnit \pth{   \atr\underline{\chi} \pth{ - \underline{\eta}  + 3 \eta }    +  \tr\underline{\chi} \pth{ \,^*\underline{\eta} - \,^*\eta    } },
\end{align*}
and thus
\begin{align*}
\RR^{[h,3]} & =   2\,^{(c)}\nab(\Re(\tilde{C}_1))  + 2 \Re( \tilde{C}_1)\underline{\eta}
\\&\quad   + 6  \tr\underline{\chi}  \eta  +  2 \atr\underline{\chi}   \,^*\eta + 2   \tr\underline{\chi}   \underline{\eta} +2  \atr\underline{\chi} \,^*\underline{\eta} + 4\ImagUnit \pth{    \atr\underline{\chi} \pth{   \eta - \underline{\eta}}    +  \tr\underline{\chi} \pth{ \,^*\underline{\eta} + \,^*\eta    }  }.
\end{align*}
However, we can deduce from \zcref[cap]{lemma:Kerr:outgoing-PG:Ric-and-curvature} that the quantity $  \atr\underline{\chi} \pth{   \eta - \underline{\eta}}    +  \tr\underline{\chi} \pth{ \,^*\underline{\eta} + \,^*\eta    }$ vanishes in the outgoing null frame and is $-1$-conformally invariant so that it also vanishes in the global frame. We thus have
\begin{align*}
\RR^{[h,3]} & =   2\,^{(c)}\nab(\Re(\tilde{C}_1))  + 2 \Re( \tilde{C}_1)\underline{\eta}  + 6  \tr\underline{\chi}  \eta  +  2 \atr\underline{\chi}   \,^*\eta + 2   \tr\underline{\chi}   \underline{\eta} +2  \atr\underline{\chi} \,^*\underline{\eta} .
\end{align*}
Using the same structure equations as above and $\Re(\tilde{C}_1)= -2 \frac{\,^{(a)}\tr\underline{\chi}^2}{\tr\underline{\chi}}$ we get that
\begin{align*}
2\,^{(c)}\nab(\Re(\tilde{C}_1))  + 2 \Re( \tilde{C}_1)\underline{\eta} & =  6 \frac{\,^{(a)}\tr\underline{\chi}^2}{\tr\underline{\chi}}   \eta  - 2 \frac{\,^{(a)}\tr\underline{\chi}^2}{\tr\underline{\chi}}   \underline{\eta}  + 2 \atr\underline{\chi}  \pth{ 2  +  \frac{\,^{(a)}\tr\underline{\chi}^2}{\tr\underline{\chi}^2} } \pth{ \,^*\underline{\eta} -  \,^*\eta   }  ,
\end{align*}
so that
\begin{align*}
\RR^{[h,3]} & =  \pth{ 1 + \frac{\,^{(a)}\tr\underline{\chi}^2}{\tr\underline{\chi}^2}   }\pth{ - 2  \atr\underline{\chi}   \,^*\eta + 6 \tr\underline{\chi}   \eta } - 2  \tr\underline{\chi} \pth{ \frac{\,^{(a)}\tr\underline{\chi}^2}{\tr\underline{\chi}^2} -  1   }  \underline{\eta} 
  + 2  \atr\underline{\chi} \pth{  3 +  \frac{\,^{(a)}\tr\underline{\chi}^2}{\tr\underline{\chi}^2}   } \,^*\underline{\eta}     .
\end{align*}
Together with $\eta_1=\underline{\eta}_1$, $\underline{\eta}_2 =-\eta_2$ and the $e_2$ component of $ \atr\underline{\chi} \pth{   \eta - \underline{\eta}}    +  \tr\underline{\chi} \pth{ \,^*\underline{\eta} + \,^*\eta    }=0$, this indeed implies that $\RR^{[h,3]}_1  = \pth{     8 + 4 \frac{\,^{(a)}\tr\underline{\chi}^2}{\tr\underline{\chi}^2}      }  \pth{ \tr\underline{\chi}  \eta_1   -  \atr\underline{\chi}  \eta_2  } = 0$ and
\begin{align*}
\RR^{[h,3]}_2 & =    4  \tr\underline{\chi} \eta_2 + 8\frac{\,^{(a)}\tr\underline{\chi}^2}{\tr\underline{\chi}}     \eta_2        -4 \atr\underline{\chi}   \eta_1  =   4  \tr\underline{\chi} \eta_2 + 4\atr\underline{\chi}  \eta_1     ,
\end{align*}
which indeed implies \zcref[noname]{Rh3_2} after having plugged the values of $\tr\underline{\chi}$, $\atr\underline{\chi}$, $\eta_1$ and $\eta_2$ in the global principal null frame.  For \zcref[noname]{Vtilde-R43 final}, we use \zcref[noname]{cnab3trchib}-\zcref[noname]{cnab4atrchi} to compute
\begin{align*}
\pth{\,^{(c)}\nab_{\R}  -\tr\chi } \Re(\tilde{C}_1) & =     \frac{8\,^{(a)}\tr\underline{\chi}}{\tr\underline{\chi}} \pth{ \curl\underline{\eta}  +   \,^*\rho  } +  \frac{\,^{(a)}\tr\underline{\chi}^2}{\tr\underline{\chi}^2} \pth{   {\tr X} {\tr\underline{X}} - 4 \div\underline{\eta}  - 4 |\eta|^2 -  4\rho  - \frac{8\La}{3} }
\\&\quad + \lambdaglo^2 \frac{|q|^2}{\De}  \,^{(a)}\tr\underline{\chi}^2  \pth{  3 + \frac{  \,^{(a)}\tr\underline{\chi}^2}{\tr\underline{\chi}^2} }  +2\tr\chi  \frac{\,^{(a)}\tr\underline{\chi}^2}{\tr\underline{\chi}}.
\end{align*}
Plugging this into the expression from \zcref[cap]{lem potentiel réel} concludes the proof of \zcref[noname]{Vtilde-R43 final}. For \zcref[noname]{R43 final}, we start from \zcref[noname]{inter R43} and use \zcref[noname]{cnab3trXb} to obtain $\RR^{[4,3]}  = - \pth{ 1  + \frac{\atr\underline{\chi}^2}{\tr\underline{\chi}^{2}}     } \,^{(a)}\tr\underline{\chi}^2$, which indeed implies \zcref[noname]{R43 final} after having plugged the values of $\tr\underline{\chi}$ and $\atr\underline{\chi}$ in the global principal null frame. This concludes the proof of the lemma.
\end{proof}

By collecting the various computations performed in the previous lemmas we get the equation satisfied by the intermediate quantity $\tilde{\q}$.

\begin{corollary}\label{coro intermediate conclusion}
If $f$ and $C_1$ are chosen to be $f  = q \bar{q}^3$ and $C_1  = 2\tr\underline{\chi} -2 \frac{\,^{(a)}\tr\underline{\chi}^2}{\tr\underline{\chi}}   -4\ImagUnit \,^{(a)}\tr\underline{\chi}$
then $\tilde{\q}$ defined by \zcref[noname]{definition of tilde q frak} satisfies
\begin{align*}
\dot{\Box}_2\tilde{\q} -  \ImagUnit\frac{4a(1+\ga)\cos\th}{|q|^2} \nab_{\T}\tilde{\q} & =  \widetilde{V} \tilde{\q}  + f  \,^{(c)}\nab_\VFArthur  \,^{(c)}\nab_3 A+ \pth{  \lambdaglo^2 \frac{|q|^2}{\De}  \RR^{[4,3]}  f C_1 + \RR^{[3]} } \,^{(c)}\nab_3 A  
\\&\quad + \RR^{[4]} \,^{(c)}\nab_4 A  + \RR^{[h]} \cdot \,^{(c)}\nab A  + \RR^{[0]}  A ,
\end{align*}
with 
\begin{align*}
\widetilde{V} & =    \frac{4\De}{(r^2+a^2)|q|^2} + 2\La + \GO{ a r^{-4}},& \VFArthur & = - \lambdaglo^{-1} (1+\ga)   \frac{8a\De }{r^2|q|^4}  \pth{   a  \T   +  \mathbf{\Phi} } ,
\end{align*}
and where $ \RR^{[4,3]}$ is given in \zcref[noname]{R43 final}, and $ \RR^{[3]}$, $ \RR^{[4]}$, $ \RR^{[h]}$ and $\RR^{[0]} $ are given in \zcref[cap]{prop q tilde}.
\end{corollary}

\begin{proof}
We recall from the discussion above \zcref[cap]{lem potentiel réel} that $\tilde{\q}$ satisfies
\begin{align*}
\dot{\Box}_2\tilde{\q} -  \ImagUnit\frac{4a(1+\ga)\cos\th}{|q|^2} \nab_{\T}\tilde{\q} & =  \widetilde{V} \tilde{\q}  + f  \,^{(c)}\nab_\VFArthur  \,^{(c)}\nab_3 A + \pth{  \lambdaglo^2 \frac{|q|^2}{\De}  \RR^{[4,3]}  f C_1 + \RR^{[3]} } \,^{(c)}\nab_3 A  + \RR^{[4]} \,^{(c)}\nab_4 A  
\\&\quad + \RR^{[h]} \cdot \,^{(c)}\nab A  + \RR^{[0]} A ,
\end{align*}
where the vector field $\VFArthur$ is given by
\begin{align*}
\VFArthur &  = (1+\ga) \pth{ \pth{ 2 \lambdaglo   \RR^{[4,3]}    \frac{(r^2+a^2)}{\De}  + \RR^{[h,3]}_2   \frac{a\sin\th}{\sqrt{\ka}|q|}} \T  +  \pth{  2 \lambdaglo  \RR^{[4,3]}    \frac{a}{\De}  +   \frac{\RR^{[h,3]}_2 }{\sqrt{\ka}|q|\sin\th} } \mathbf{\Phi} } ,
\end{align*}
where we already used \zcref[noname]{Rh3_1=0} to remove the $\dr_\th$ derivative, and where the potential $\widetilde{V}$ is given by \zcref[noname]{Vtilde-R43 final}. The expansion of $\widetilde{V}$ stated in the corollary follows from the fall-off properties of the various quantities involved, see for instance \zcref[noname]{cnabHHb}. It remains to plug \zcref[noname]{Rh3_2} and \zcref[noname]{R43 final} inside the above expression for $\VFArthur$ to conclude the proof.
\end{proof}

\subsubsection{The quantity $\q$ and the choice of the leading part of $C_2$}

We now consider the quantity $\q = \tilde{\q} + f C_2 A$, with $C_2$ a $-2$-conformally invariant function yet to be determined. From the previous corollary we first obtain
\begin{align*}
\dot{\Box}_2 \q & =  \ImagUnit\frac{4a(1+\ga)\cos\th}{|q|^2} \nab_{\T} \q +  \widetilde{V} \q    + f  \,^{(c)}\nab_\VFArthur  \,^{(c)}\nab_3 A
\\&\quad + \pth{  \la^2 \frac{|q|^2}{\De}  \RR^{[4,3]}  f C_1  + \RR^{[3]}  -   f C_2  \ImagUnit    \atr\chi } \,^{(c)}\nab_3 A  + \pth{ \RR^{[4]} -   f C_2  \ImagUnit    \atr\underline{\chi} } \,^{(c)}\nab_4 A 
\\&\quad  + \pth{ \RR^{[h]}  +2  f C_2  \ImagUnit    (\,^*\eta + \,^*\underline{\eta}) }\cdot \,^{(c)}\nab A  + \pth{ \RR^{[0]}  - \ImagUnit\frac{4a(1+\ga)\cos\th}{|q|^2} \,^{(c)}\nab_{\T} ( f C_2 ) -  \widetilde{V}  f C_2 } A
\\&\quad + \dot{\Box}_2 (f C_2 A),
\end{align*}
where we also used $ \frac{4a(1+\ga)\cos\th}{|q|^2} \T  =  \,^{(a)}\tr\underline{\chi} e_4 +  \,^{(a)}\tr\chi e_3 -2 \pth{ \,^*\eta_1 + \,^*\underline{\eta}_1   } e_1 -2 \pth{ \,^*\eta_2 + \,^*\underline{\eta}_2   } e_2$ which follows from the first relation in \zcref[noname]{eq:T-R-Z:principal-outgoing:expression:T-Z}. Using now \zcref[noname]{Box in terms of L} and the commutator from \zcref[cap]{lemma Lf} we can compute the last term in the above equation and obtain
\begin{equation}\label{equation finale pour q}
\begin{aligned}
\dot{\Box}_2 \q & =  \ImagUnit\frac{4a(1+\ga)\cos\th}{|q|^2} \nab_{\T} \q +  \widetilde{V} \q    + f  \,^{(c)}\nab_\VFArthur  \,^{(c)}\nab_3 A
 + \tilde{\RR}^{[4]}  \,^{(c)}\nab_4 A  + \tilde{\RR}^{[3]} \,^{(c)}\nab_3 A   + \tilde{\RR}^{[h]}  \cdot \,^{(c)}\nab A + \tilde{\RR}^{[0]} A,
\end{aligned}
\end{equation}
where 
\begin{align*}
\tilde{\RR}^{[4]} & =  \RR^{[4]} -   f C_2  \ImagUnit   \atr\underline{\chi} - \,^{(c)}\nab_3(fC_2)   ,
\\ \tilde{\RR}^{[3]} & =   \la^2 \frac{|q|^2}{\De}  \RR^{[4,3]}  f C_1  + \RR^{[3]}  -   f C_2  \ImagUnit  \atr\chi   - \,^{(c)}\nab_4(fC_2)    + 2 \overline{\tr X} fC_2 ,
\\ \tilde{\RR}^{[h]} & =  \RR^{[h]}   + 2   f C_2  \ImagUnit (\,^*\eta + \,^*\underline{\eta}) +  2  \,^{(c)}\nab(fC_2) -4 fC_2 H,
\\ \tilde{\RR}^{[0]} & =  \RR^{[0]}  - \pth{ \frac{1}{4} \tr X \overline{\tr\underline{X}} + \frac{1}{4} \tr\underline{X}\overline{\tr X} +   2 \overline{P}  - \frac{2\La}{3} + 2\ImagUnit \eta\wedge\underline{\eta} + T^{[0]} + \widetilde{V} }f C_2 
\\&\quad   - \ImagUnit\frac{4a(1+\ga)\cos\th}{|q|^2} \,^{(c)}\nab_{\T} ( f C_2 )  -\half \tr X   \,^{(c)}\nab_3 (f C_2)  - \half \tr\underline{X}  \,^{(c)}\nab_4 (f C_2) 
\\&\quad + 2\underline{\eta} \cdot \,^{(c)}\nab (f C_2)  - \,^{(c)}\nab_4   \,^{(c)}\nab_3 (f C_2) +  \half \,^{(c)}\DD \cdot \overline{ \,^{(c)}\DD} (f C_2).
\end{align*}
The decay properties of these coefficients which are required by \zcref[cap]{proposition T to RW} are stated in the following lemma and they are obtained by a careful choice of the leading part of $C_2$. More precisely, we can see in the following proof that having $C_2$ proportional to $\tr\underline{\chi}^2$ (at leading order) is enough for $\tilde{\RR}^{[h]} $, while for $\tilde{\RR}^{[4]} $, $\tilde{\RR}^{[3]} $ and $\tilde{\RR}^{[0]} $ we need exact cancellations.

\begin{lemma}\label{lem C2}
If $C_2= \half \tr\underline{\chi}^2 + \tilde{C}_2$ with $\tilde{C}_2=\GO{a r^{-3}}$ and $\,^{(c)}\nab_\mu (\tilde{C}_2)=\GO{a r^{-4}}$ for $\mu=1,2,3,4$ then $\tilde{\RR}^{[h]}  = \tilde{\RR}^{[4]} = \tilde{\RR}^{[3]} =   \GO{a }$ and $\tilde{\RR}^{[0]}  =  \GO{ar^{-1}}$.
\end{lemma}

\begin{proof}
We start with $\tilde{\RR}^{[h]}$. From \zcref[noname]{RRh} and $C_1 = \GO{r^{-1}}$ we find $\RR^{[h]} =    f  \GO{a r^{-4}}$ and thus
\begin{align*}
f^{-1}\tilde{\RR}^{[h]} & = 2  \,^{(c)}\nab(C_2) + \pth{  4 \underline{\eta}    - \ImagUnit\,^*\eta } C_2   +   \GO{a r^{-4}},
\end{align*}
where we also used \zcref[noname]{equations for f} to compute $\nab f$. Using now $C_2= \half \tr\underline{\chi}^2 + \tilde{C}_2$ and \zcref[noname]{cnabtrchib} (which implies $\,^{(c)}\nab\tr\underline{\chi}=\GO{a r^{-3}}$) we obtain $\tilde{\RR}^{[h]}=f   \GO{a r^{-4}}$, as stated in the lemma. We now look at $\tilde{\RR}^{[4]}$. From \zcref[noname]{cnab3trXb}, \zcref[noname]{RR4} and $C_1 = 2\tr\underline{\chi} + \GO{a r^{-2}}$ we find $ \RR^{[4]}   = f\pth{ \half \tr\underline{\chi}^3  + \GO{a r^{-4}} }$ so that
\begin{align*}
f^{-1}\tilde{\RR}^{[4]} & = \half \tr\underline{\chi}^3   - \pth{ \,^{(c)}\nab_3(C_2) + 2 \tr\underline{\chi} C_2 } + \GO{a r^{-4}},
\end{align*}
where we also used \zcref[noname]{equations for f} to compute $\nab_3 f$. Using $C_2= \half \tr\underline{\chi}^2 + \tilde{C}_2$ and \zcref[noname]{cnab3trchib} we find that 
\begin{align*}
  \,^{(c)}\nab_3(C_2) +  2 \tr\underline{\chi} C_2 & = \half \tr\underline{\chi}^3 +   \,^{(c)}\nab_3(\tilde{C}_2) +  2 \tr\underline{\chi} \tilde{C}_2 + \GO{a r^{-4}}.
\end{align*}
We therefore obtain a cancellation leading to $\tilde{\RR}^{[4]}=f   \GO{a r^{-4}}$. We now look at $\tilde{\RR}^{[3]}$. From \zcref[noname]{equations for f}, \zcref[noname]{R43 final}, $C_1=\GO{r^{-1}}$ and $C_2=\GO{r^{-2}}$ we get
\begin{align*}
f^{-1}\tilde{\RR}^{[3]} & = f^{-1}\RR^{[3]}    - \,^{(c)}\nab_4(C_2)  + \GO{a r^{-4}}.
\end{align*}
From \zcref[noname]{RR3}, \zcref[noname]{equations for f}, \zcref[cap]{lemma Lf}, \zcref[noname]{cnab3trchib}, \zcref[noname]{cnab4trchib}, \zcref[noname]{cnabtrchib}, \zcref[noname]{cnab3rho}, the expression and the fall-off of the $F$'s, the expression of $C_1$ we get
\begin{align*}
f^{-1}\RR^{[3]} & = \tr\underline{\chi} \pth{ - \half {\tr \chi} {\tr\underline{\chi}}      +  2\rho  + \frac{4\La}{3}} +  \GO{a r^{-4}} .
\end{align*}
However using \zcref[noname]{cnab3trchib}, \zcref[noname]{cnab4trchib} and \zcref[noname]{cnabtrchib} we also obtain
\begin{align*}
 \,^{(c)}\nab_4(C_2) & = \tr\underline{\chi} \pth{ - \half {\tr \chi} {\tr\underline{\chi}}  +  2\rho  + \frac{4\La}{3}} +  \,^{(c)}\nab_4(\tilde{C}_2) +  \GO{a r^{-4}} .
\end{align*}
Using the assumptions on $\tilde{C}_2$ we then indeed obtain $\tilde{\RR}^{[3]}=f\GO{a r^{-4}}$. We finally compute $\tilde{\RR}^{[0]}$. From \zcref[noname]{RR0}, \zcref[noname]{cnab3trX}, \zcref[noname]{cnab3trXb}, \zcref[noname]{cnab3rho}, $C_1=2\tr\underline{\chi} + \GO{ar^{-2}}$ and  $\La\lesssim r^{-2}$ we find
\begin{align*}
f^{-1}\RR^{[0]} & =     \tr \underline{\chi}^2 \pth{   6 \rho  + 3 \La}  +  \GO{ar^{-5}} .
\end{align*}
Therefore, using $\widetilde{V}=-\tr\chi\tr\underline{\chi} + 2\La + \GO{a r^{-4}}$ we obtain
\begin{align*}
f^{-1}\tilde{\RR}^{[0]} & =    \tr \underline{\chi}^2 \pth{   6 \rho  + 3 \La}  + \pth{  \frac{3}{2} \tr \chi {\tr\underline{\chi}}  -   4 \rho  - \frac{2\La}{3}   }C_2   +  \GO{ar^{-5}}
\\&\quad - \half \overline{\tr X} f^{-1}\,^{(c)}\nab_3 ( f C_2 )    - \half \overline{\tr\underline{X}} f^{-1} \,^{(c)}\nab_4 (f C_2)  +2 f^{-1} \pth{ \underline{\eta}   +  \ImagUnit  \pth{ \,^*\eta + \,^*\underline{\eta}   } }\cdot \,^{(c)}\nab ( f C_2 )
\\&\quad  - f^{-1} \,^{(c)}\nab_4   \,^{(c)}\nab_3 (f C_2) + f^{-1} \half \,^{(c)}\DD \cdot \overline{ \,^{(c)}\DD} (f C_2)
\end{align*}
We now use \zcref[noname]{equations for f} to compute the various derivatives of $f$, \zcref[noname]{cnab4trXbbis} and $C_2=\GO{r^{-2}}$ to get
\begin{align*}
f^{-1}\tilde{\RR}^{[0]} & =    \tr \underline{\chi}^2 \pth{   6 \rho  + 3 \La}  + \pth{    - \frac{7}{2} \tr \chi {\tr\underline{\chi}}   -   8 \rho  - \frac{10\La}{3}   }C_2          +  \GO{ar^{-5}}
\\&\quad - \half \pth{    \tr X + 4 \overline{\tr X} } \,^{(c)}\nab_3 ( C_2 )   - \half  \pth{ 2  \overline{\tr\underline{X}}   + 3 \tr\underline{X}   } \,^{(c)}\nab_4 (C_2)  +  \pth{ 6\underline{\eta}    + 4 H   }\cdot \,^{(c)}\nab (  C_2 )
\\&\quad  -  \,^{(c)}\nab_4  \,^{(c)}\nab_3 ( C_2)  +  \half \,^{(c)}\DD \cdot   \overline{ \,^{(c)}\DD} (C_2) .
\end{align*}
Using now $C_2= \half \tr\underline{\chi}^2 + \tilde{C}_2$ and  \zcref[noname]{cnab3trchib}, \zcref[noname]{cnab4trchib}, \zcref[noname]{cnabtrchib}, we finally obtain a cancellation at leading order and we get
\begin{align*}
f^{-1}\tilde{\RR}^{[0]} & =   - \half \pth{    \tr X + 4 \overline{\tr X} } \,^{(c)}\nab_3 (  \tilde{C}_2 )   - \half  \pth{ 2  \overline{\tr\underline{X}}   + 3 \tr\underline{X}   } \,^{(c)}\nab_4 ( \tilde{C}_2)  +  \pth{ 6\underline{\eta}    + 4 H   }\cdot \,^{(c)}\nab (   \tilde{C}_2 )
\\&\quad  -  \,^{(c)}\nab_4  \,^{(c)}\nab_3 (  \tilde{C}_2)  +  \half \,^{(c)}\DD \cdot   \overline{ \,^{(c)}\DD} ( \tilde{C}_2) +  \GO{ar^{-5}},
\end{align*}
which indeed gives $\tilde{\RR}^{[0]}=f   \GO{a r^{-5}}$ thanks to the assumptions on $\tilde{C}_2$. 
\end{proof}

Note that with the notations of \zcref[cap]{proposition T to RW} and in view of \zcref[noname]{equation finale pour q}, we have $f W^{[i]} = \tilde{\RR}^{[i]}$ for $i=3,4,h,0$. Thus, \zcref[cap]{lem C2} states the correct fall-off properties for the lower order coefficients in \zcref[noname]{gRW}.

\subsubsection{Factorization of $\q$ and the choice of $\tilde{C}_2$}

To conclude the proof of \zcref[cap]{proposition T to RW}, it remains to prove that $\q  = f \pth{   \,^{(c)}\nab_3\,^{(c)}\nab_3 A  + C_1 \,^{(c)}\nab_3 A + C_2 A } $ with 
\begin{align*}
f & = q \bar{q}^3 , & C_1 & = 2\tr\underline{\chi} -2 \frac{\,^{(a)}\tr\underline{\chi}^2}{\tr\underline{\chi}}   -4\ImagUnit \,^{(a)}\tr\underline{\chi} ,&  C_2 & = \half \tr\underline{\chi}^2 + \tilde{C}_2,
\end{align*}
can be factorized as \zcref[noname]{definition of q frak}. This requirement is actually what defines $\tilde{C}_2$, and one can check in its expression below that it indeed satisfies the fall-off properties stated in \zcref[cap]{lem C2}.

\begin{lemma}
If 
\begin{align*}
\tilde{C}_2 & = - 4 \atr\underline{\chi}^2 + \frac{3}{2} \frac{\atr\underline{\chi}^4}{\tr\underline{\chi}^2} + \ImagUnit\pth{ - 2\tr\underline{\chi} \atr\underline{\chi} + 4 \frac{\atr\underline{\chi}^3}{\tr\underline{\chi}} },
\end{align*}
then $\q  = \frac{q}{\bar{q}} r^2 \,^{(c)}\nab_3 \,^{(c)} \nab_3 \pth{ \frac{\bar{q}^4}{r^2}A} $.
\end{lemma}

\begin{proof}
First, using the second property of \zcref[noname]{properties of q appendix} we compute
\begin{align*}
2\frac{\nab_3 \pth{ \frac{ \bar{q}^4}{r^2} } }{ \frac{ \bar{q}^4}{r^2}} & = 2\tr\underline{\chi} - 4 \ImagUnit \atr\underline{\chi} +   \frac{ 2a\cos\th}{r} \atr\underline{\chi} .
\end{align*}
Using $-2 \frac{\atr\underline{\chi}}{\tr\underline{\chi}}  =   \frac{ 2a\cos\th}{r}$ in the global principal null frame and the fact that $\frac{ \bar{q}^4}{r^2} $ is $0$-conformally invariant we have thus obtained
\begin{align}\label{derivee loga 1}
2\frac{\,^{(c)}\nab_3 \pth{ \frac{ \bar{q}^4}{r^2} } }{ \frac{ \bar{q}^4}{r^2}} & = C_1.
\end{align} 
Using this identity twice we also obtain
\begin{align*}
\frac{\,^{(c)}\nab_3 \,^{(c)} \nab_3 \pth{ \frac{ \bar{q}^4}{r^2} }}{\frac{ \bar{q}^4}{r^2} } & = \frac{1}{\frac{ \bar{q}^4}{r^2} } \,^{(c)}\nab_3 \pth{ \half \frac{ \bar{q}^4}{r^2} C_1  } =  \half \,^{(c)}\nab_3 C_1   + \frac{1}{4} C_1^2  .
\end{align*}
We now choose $\tilde{C}_2$ so that $ \half \,^{(c)}\nab_3 C_1   + \frac{1}{4} C_1^2 = C_2 $ where we recall that $C_2 = \half \tr\underline{\chi}^2 + \tilde{C}_2$. This forces us to choose $\tilde{C}_2  = \half \,^{(c)}\nab_3 C_1   + \frac{1}{4} C_1^2 - \half \tr\underline{\chi}^2 $ which indeed gives the expression stated in the lemma. With this choice of $\tilde{C}_2$, we thus have
\begin{align}\label{derivee loga 2}
\frac{\,^{(c)}\nab_3 \nab_3 \pth{ \frac{ \bar{q}^4}{r^2} }}{\frac{ \bar{q}^4}{r^2} }  & = C_2.
\end{align}
Now, \zcref[noname]{derivee loga 1} and \zcref[noname]{derivee loga 2} imply that
\begin{align*}
\frac{ \,^{(c)}\nab_3 \,^{(c)} \nab_3 \pth{ \frac{ \bar{q}^4}{r^2} A} }{\frac{ \bar{q}^4}{r^2}} & =   \,^{(c)}\nab_3   \,^{(c)} \nab_3  A   + 2 \frac{ \nab_3  \pth{\frac{ \bar{q}^4}{r^2}}   }{\frac{ \bar{q}^4}{r^2}} \,^{(c)} \nab_3  A +  \frac{ \,^{(c)}\nab_3   \nab_3  \pth{\frac{ \bar{q}^4}{r^2}} }{\frac{ \bar{q}^4}{r^2}} A
\\& = \,^{(c)}\nab_3\,^{(c)}\nab_3 A  + C_1 \,^{(c)}\nab_3 A + C_2 A ,
\end{align*}
which concludes the proof of the lemma.
\end{proof}

\printbibliography

\end{document}


%% file: Images/Intro-Penrose-KdS.tex
\begin{tikzpicture}[scale=0.7,every node/.style={scale=0.7}]


  \def \s{3} 
  \def \exts{0.2} 
  \def \t{0.4}
  \def \Tlen{.5}

  \coordinate (tInf) at (0,\s); 
  \coordinate (EventZero) at (-\s,0); 
  \coordinate (CosmoZero) at (\s,0); 
  \coordinate (tNegInf) at (0,-\s);
  \coordinate (SigmaZeroEvent) at (-\s + 0.15*\s, 0.15*\s);
  \coordinate (SigmaZeroCosmo) at (\s - 0.15*\s, 0.15*\s);

  \draw[shorten >= -10,name path=EventFuture] (tInf) --
  node[pos=0.5,left]{$\EventHorizonFuture$} (EventZero) ;  
  \draw[shorten >= -10,name path=CosmoFuture] (tInf) --
  node[inner sep=10pt,scale=1.0,pos=0.5,right]{$\CosmologicalHorizonFuture$} (CosmoZero) ;

  \draw[decoration={complete sines},decorate] (\exts * \s,\s) --  (-\exts * \s,\s);


  \draw[dashed] (tInf)    -- (-\s - \exts*\s,-\exts*\s);
  \draw[dashed] (tInf)    -- ( \s + \exts*\s,-\exts*\s);

  \draw[dashed, decoration={complete sines},decorate] (\exts * \s,\s)
  --node[pos=0.3,above]{} (2 * \exts * \s,\s);
  \draw[dashed, decoration={complete sines},decorate] (-\exts * \s,\s)
  --node[pos=0.3,left]{}  (-2 * \exts * \s,\s);

  \path[name path=EventLowerBound] (-1.6*\s, 0.6*\s)-- (-\s, 0)--
  (tInf);
  \path[name path=CosmoLowerBound] (tInf) -- (\s, 0) -- (1.6*\s, 0.6*\s);


  \path[name path=SigmaZero, blue,thick,out=-39,in=-141,shorten >= -20, shorten <= -20]
  (SigmaZeroEvent) edge node[pos=0.5,below]{$\Sigma_0$} node[pos=0.8](TInit){} (SigmaZeroCosmo) ;

  \path[fill = black, fill opacity=0.1, draw=none](SigmaZeroEvent) to [out=-35,in=-145] (SigmaZeroCosmo) to (tInf) --cycle;



  
  \node[scale=0.5,fill=white,draw,circle,label=above:$\timelikeInf$]at(tInf){};

  \node[label=left:${r=r_{\mathcal{H},\Lambda}}$]at(SigmaZeroEvent){};
  \node[label=right:${r=r_{\overline{\mathcal{H}},\Lambda}}\sim \Lambda^{-\frac{1}{2}}$]at(SigmaZeroCosmo){};

\end{tikzpicture}


%% file: Images/Intro-Penrose-Kerr.tex
\begin{tikzpicture}[scale=0.7,every node/.style={scale=0.7}]


  \def \s{3} 
  \def \exts{0.2} 
  \def \t{0.4}
  \def \Tlen{.5}

  \coordinate (tInf) at (0,\s); 
  \coordinate (EventZero) at (-\s,0); 
  \coordinate (CosmoZero) at (\s,0); 
  \coordinate (tNegInf) at (0,-\s);
  \coordinate (SigmaZeroEvent) at (-\s + 0.15*\s, 0.15*\s);
  \coordinate (SigmaZeroCosmo) at (\s - 0.15*\s, 0.15*\s);

  \draw[shorten >= -10,name path=EventFuture] (tInf) --
  node[pos=0.5,left]{$\EventHorizonFuture$} (EventZero) ;  
  \draw[shorten >= -10,name path=CosmoFuture,dashed] (tInf) --
  node[inner sep=10pt,scale=1.0,pos=0.5,right]{$\mathcal{I}^{+}$} (CosmoZero) ; 


  \draw[decoration={complete sines},decorate] (-\exts * \s,\s) --  (0,\s);
  

  \draw[dashed] (tInf)    -- (-\s - \exts*\s,-\exts*\s);
  \draw[dashed] (tInf)    -- ( \s + \exts*\s,-\exts*\s);

  \draw[dashed, decoration={complete sines},decorate] (-\exts * \s,\s)
  --node[pos=0.3,left]{}  (-2 * \exts * \s,\s);

  \path[name path=EventLowerBound] (-1.6*\s, 0.6*\s)-- (-\s, 0)--
  (tInf);
  \path[name path=CosmoLowerBound] (tInf) -- (\s, 0) -- (1.6*\s, 0.6*\s);


  \path[name path=SigmaZero, blue,thick,out=-39,in=-141,shorten <= -20]
  (SigmaZeroEvent) edge node[pos=0.5,below]{$\Sigma_0$} node[pos=0.8](TInit){} (SigmaZeroCosmo) ;

  \path[fill = black, fill opacity=0.1, draw=none](SigmaZeroEvent) to [out=-35,in=-145] (SigmaZeroCosmo) to (tInf) --cycle;



  
  \node[scale=0.5,fill=white,draw,circle,label=above:$\timelikeInf$]at(tInf){};

  \node[label=left:${r=r_{\mathcal{H},0}}$]at(SigmaZeroEvent){};
  \node[label=right:${r=\infty}$]at(SigmaZeroCosmo){};

\end{tikzpicture}


%% file: Images/intro-penrose-3-obstacles.tex
\begin{tikzpicture}[scale=1,every node/.style={scale=1}]


  \def \s{3} 
  \def \exts{0.2} 
  \def \t{0.4}
  \def \Tlen{.5}

  \coordinate (tInf) at (0,\s); 
  \coordinate (EventZero) at (-\s,0); 
  \coordinate (CosmoZero) at (\s,0); 
  \coordinate (tNegInf) at (0,-\s);
  \coordinate (SigmaZeroEvent) at (-\s + 0.15*\s, 0.15*\s);
  \coordinate (SigmaZeroCosmo) at (\s - 0.15*\s, 0.15*\s);

  \draw[shorten >= -10,name path=EventFuture] (tInf) --
  node[pos=0.5,left]{$\EventHorizonFuture$} (EventZero) ;  
  \draw[shorten >= -10,name path=CosmoFuture] (tInf) --
  node[inner sep=10pt,scale=1.0,pos=0.5,right]{$\CosmologicalHorizonFuture$} (CosmoZero) ;

  \draw[decoration={complete sines},decorate] (\exts * \s,\s) --  (-\exts * \s,\s);


  \draw[dashed] (tInf)    -- (-\s - \exts*\s,-\exts*\s);
  \draw[dashed] (tInf)    -- ( \s + \exts*\s,-\exts*\s);

  \draw[dashed, decoration={complete sines},decorate] (\exts * \s,\s)
  --node[pos=0.3,above]{} (2 * \exts * \s,\s);
  \draw[dashed, decoration={complete sines},decorate] (-\exts * \s,\s)
  --node[pos=0.3,left]{}  (-2 * \exts * \s,\s);

  \path[name path=EventLowerBound] (-1.6*\s, 0.6*\s)-- (-\s, 0)--
  (tInf);
  \path[name path=CosmoLowerBound] (tInf) -- (\s, 0) -- (1.6*\s, 0.6*\s);



  \path[fill = black, fill opacity=0.1, draw=none](SigmaZeroEvent) to [out=-35,in=-145] (SigmaZeroCosmo) to (tInf) --cycle;



  
  

  \node[label=left:${r=r_{\mathcal{H},\Lambda}}$]at(SigmaZeroEvent){};
  \node[label=right:${r=r_{\overline{\mathcal{H}},\Lambda}}\sim \Lambda^{-\frac{1}{2}}$]at(SigmaZeroCosmo){};

\coordinate (ergo) at (- 0.7*\s, 0.0*\s+0.17);
  \path[fill = black, fill opacity=0.3, draw=none](SigmaZeroEvent) to [out=-35,in=150] (ergo) to (tInf) --cycle;
    \node[text=white] at (-2.2, 0.5) {$A$};

\coordinate (T1) at (-1.3, -0.19);
\coordinate (T2) at (-0.8, -0.32);
  \path[fill = black, fill opacity=0.3, draw=none](T1) to [out=-15,in=172] (T2) to (tInf) --cycle;
\node[text=white] at (-1, 0) {$B$};
\node at (-1.2, -0.6) {$r\sim 3M$};

\coordinate (ergobis) at (0.7*\s, 0.0*\s+0.17);
  \path[fill = black, fill opacity=0.3, draw=none](ergobis) to [out=35,in=210] (SigmaZeroCosmo) to (tInf) --cycle;
      \node[text=white] at (2.2, 0.5) {$C$};

 \node[scale=0.5,fill=white,draw,circle,label=above:$\timelikeInf$]at(tInf){}; 
   \path[name path=SigmaZero, blue,thick,out=-39,in=-141,shorten >= -20, shorten <= -20] (SigmaZeroEvent) edge node[pos=0.5,below]{$\Sigma_0$} node[pos=0.8](TInit){} (SigmaZeroCosmo) ;
\end{tikzpicture}


%% file: Images/intro-penrose-3.tex
\begin{tikzpicture}[scale=1,every node/.style={scale=1}]


  \def \s{3} 
  \def \exts{0.2} 
  \def \t{0.4}
  \def \Tlen{.5}

  \coordinate (tInf) at (0,\s); 
  \coordinate (EventZero) at (-\s,0); 
  \coordinate (CosmoZero) at (\s,0); 
  \coordinate (tNegInf) at (0,-\s);
  \coordinate (SigmaZeroEvent) at (-\s + 0.15*\s, 0.15*\s);
  \coordinate (SigmaZeroCosmo) at (\s - 0.15*\s, 0.15*\s);

  \draw[shorten >= -10,name path=EventFuture] (tInf) --
  node[pos=0.5,left]{$\EventHorizonFuture$} (EventZero) ;  
  \draw[shorten >= -10,name path=CosmoFuture,dashed] (tInf) --
  node[inner sep=10pt,scale=1.0,pos=0.5,right]{$\mathcal{I}^{+}$} (CosmoZero) ; 


  \draw[decoration={complete sines},decorate] (-\exts * \s,\s) --  (0,\s);
  

  \draw[dashed] (tInf)    -- (-\s - \exts*\s,-\exts*\s);
  \draw[dashed] (tInf)    -- ( \s + \exts*\s,-\exts*\s);

  \draw[dashed, decoration={complete sines},decorate] (-\exts * \s,\s)
  --node[pos=0.3,left]{}  (-2 * \exts * \s,\s);

  
  \path[name path=EventLowerBound] (-1.6*\s, 0.6*\s)-- (-\s, 0)-- (tInf);
  \path[name path=CosmoLowerBound] (tInf) -- (\s, 0) -- (1.6*\s, 0.6*\s);
  


  \path[name path=SigmaZero, blue,thick,out=-39,in=-141,shorten <= -20]
  (SigmaZeroEvent) edge node[pos=0.5,below left]{$\Sigma_0$} node[pos=0.8](TInit){} (SigmaZeroCosmo) ;

  \path[fill = black, fill opacity=0.1, draw=none](SigmaZeroEvent) to [out=-35,in=-145] (SigmaZeroCosmo) to (tInf) --cycle;




  \coordinate (rLambda) at (\s - 0.7*\s, -0.11*\s);
  \node[label=below:${r=r_{\overline{\mathcal{H}},\La}}$] at (rLambda) {};
  \node[label=below:${\mathcal{M}_{\mathrm{tot},\La}}$] at (-0.3\s,0.99\s) {};
  \path[fill = black, fill opacity=0.3, draw=none](SigmaZeroEvent) to [out=-27,in=-165] (rLambda) to (tInf) --cycle;
  \draw (tInf)    -- (rLambda);
  
  \node[scale=0.5,fill=white,draw,circle,label=above:$\timelikeInf$]at(tInf){};

  \node[label=left:${r=r_{\mathcal{H},0}}$]at(SigmaZeroEvent){};
  \node[label=right:${r=\infty}$]at(SigmaZeroCosmo){};

  \node[label=below:${\mathcal{M}_{\mathrm{tot},0}}$] at (0.5,2) {};

  
\end{tikzpicture}


%% file: Images/Global-foliation.tex
\begin{tikzpicture}[scale=1.2,every node/.style={scale=1.0}]


  \def \s{3} 
  \def \exts{0.2} 
  \def \t{0.4}
  \def \Tlen{.5}
  \def \endloc{1.3} 
  \def \tPlus{0.2} 

  \coordinate (tInf) at (0,\s); 
  \coordinate (EventZero) at (-\s,0); 
  \coordinate (CosmoZero) at (\s,0); 
  \coordinate (tNegInf) at (0,-\s);
  \coordinate (SigmaZeroEvent) at (-\s + 0.15*\s, 0.15*\s);
  \coordinate (SigmaZeroCosmo) at (\s - 0.15*\s, 0.15*\s);
  \coordinate (SigmaZeroEventExt) at ($(SigmaZeroEvent) + (-\exts, \exts)$);
  \coordinate (SigmaZeroCosmoExt) at ($(SigmaZeroCosmo) + (\exts, \exts)$);
  \path[name path=SigmaZero,in=-35,out=-145]
   (SigmaZeroCosmoExt) to node[pos=0.5,below]{$\SigmaInit$} node[pos=0.5](SigmaZeroHalf){} (SigmaZeroEventExt) ;
   \coordinate (SigmaZeroCosmoEnd) at ($(SigmaZeroHalf) + (\endloc*\s, 0)$);
   \def \negLeft{0.65} 
   \def \negRight{0.5} 
   \def \negLeftAngle{30} 
   \coordinate (SigmaNegLeftEnd) at ($(SigmaZeroHalf)!\negLeft!(SigmaZeroCosmoEnd)$);

   \path[draw=none, fill=olive!10] (SigmaZeroCosmoExt) to [out=-40, in = 135]
   coordinate[pos=\negRight](SigmaNegRightEnd)
   (SigmaZeroCosmoEnd) to (SigmaZeroHalf) --cycle;
   \draw[purple!50, thick, dashed] (SigmaZeroCosmoExt) to [out=-40, in = 136] node[pos=0.5, right]{$\Sigma_{{*,e}}$} (SigmaZeroCosmoEnd);

   \draw[blue,dotted,thick] (SigmaNegLeftEnd) to [out=\negLeftAngle, in = -135] node[pos=0.6, left]{$\Sigma(-1)$} (SigmaNegRightEnd);

  \path[name path = MBoundary, draw=none, fill=black!5]
  (SigmaZeroEventExt) to [out = 40, in = -135] coordinate[pos=\tPlus](SigmaPlusEventExt)
  (tInf) to [out = -45, in=140] coordinate[pos=1-\tPlus](SigmaPlusCosmoExt)
  (SigmaZeroCosmoExt) to [out=-145,in=-35] (SigmaZeroEventExt)
  ;  
  \draw[purple!90, thick,dashed] (SigmaZeroEventExt) to [out = 40, in = -135]
  node[pos=0.5,left]{$\mathcal{A}$}
  (tInf) to [out = -45, in=140]
  node[pos=0.5,right]{$\SigmaStar$}
  (SigmaZeroCosmoExt);
  \draw[blue,thick,dotted] (SigmaZeroCosmoExt) to [out=-145,in=0] node[pos=0.3, left]{$\Sigma(0)$} (SigmaZeroHalf);
  \draw[blue,thick,dashed] (SigmaZeroEventExt) to [out=-35,in=180] (SigmaZeroHalf) to (SigmaZeroCosmoEnd);
  \draw[blue, thick,dotted] (SigmaPlusEventExt) to [out=-35, in=-145] node[pos=0.5,above]{$\Sigma(1)$}(SigmaPlusCosmoExt);

  \draw[shorten >= -10,name path=EventFuture] (tInf) --
  node[inner sep=10pt, pos=1.1,below]{$\EventHorizonFuture$} (EventZero) ;  
  \draw[shorten >= -10,name path=CosmoFuture] (tInf) --
  node[inner sep=10pt,scale=1.0,pos=1.1,below]{$\CosmologicalHorizonFuture$} (CosmoZero) ;



  \draw[dashed] (tInf)    -- (-\s - \exts*\s,-\exts*\s);
  \draw[dashed] (tInf)    -- ( \s + \exts*\s,-\exts*\s);


  \path[name path=EventLowerBound] (-1.6*\s, 0.6*\s)-- (-\s, 0)--
  (tInf);
  \path[name path=CosmoLowerBound] (tInf) -- (\s, 0) -- (1.6*\s, 0.6*\s);

  \node[scale=0.5,fill=white,draw,circle,label=above:$\timelikeInf$]at(tInf){};


\end{tikzpicture}


%% file: Images/Internal-Boundary.tex
\begin{tikzpicture}[scale=1.2,every node/.style={scale=1.0}]


  \def \s{3} 
  \def \exts{0.2} 
  \def \t{0.4}
  \def \Tlen{.5}
  \def \endloc{1.3} 
  \def \tPlus{0.2} 
  \def \tWidth{0.2} 

  \coordinate (tInf) at (0,\s); 
  \coordinate (EventZero) at (-\s,0); 
  \coordinate (CosmoZero) at (\s,0); 
  \coordinate (tNegInf) at (0,-\s);
  \coordinate (SigmaZeroEvent) at (-\s + 0.15*\s, 0.15*\s);
  \coordinate (SigmaZeroCosmo) at (\s - 0.15*\s, 0.15*\s);
  \coordinate (SigmaZeroEventExt) at ($(SigmaZeroEvent) + (-\exts, \exts)$);
  \coordinate (SigmaZeroCosmoExt) at ($(SigmaZeroCosmo) + (\exts, \exts)$);
  \path[name path=SigmaZero,in=-35,out=-145]
   (SigmaZeroCosmoExt) to node[pos=0.5,below]{$\SigmaInit$} node[pos=0.5](SigmaZeroHalf){} (SigmaZeroEventExt) ;
   \coordinate (SigmaZeroCosmoEnd) at ($(SigmaZeroHalf) + (\endloc*\s, 0)$);
   \def \negLeft{0.65} 
   \def \negRight{0.5} 
   \def \negLeftAngle{30} 
   \coordinate (SigmaNegLeftEnd) at ($(SigmaZeroHalf)!\negLeft!(SigmaZeroCosmoEnd)$);

   \path[draw=none, fill=black!10] (SigmaZeroCosmoExt) to [out=-40, in = 135]
   coordinate[pos=\negRight](SigmaNegRightEnd)
   (SigmaZeroCosmoEnd) to (SigmaZeroHalf) --cycle;
   \draw[black!50, thick, dashed] (SigmaZeroCosmoExt) to [out=-40, in = 136] node[pos=0.5, right]{$\Sigma_{{*,e}}$} (SigmaZeroCosmoEnd);


  \path[name path = MBoundary, draw=none, fill=black!5]
  (SigmaZeroEventExt) to [out = 40, in = -135]
  coordinate[pos=\tPlus](SigmaPlusEventExt)
  coordinate[pos=\tPlus+\tWidth](SigmaPlus2EventExt)
  (tInf) to [out = -45, in=140]
  coordinate[pos=1-\tPlus](SigmaPlusCosmoExt)
  coordinate[pos=1-\tPlus-\tWidth](SigmaPlus2CosmoExt)
  (SigmaZeroCosmoExt) to [out=-145,in=-35] (SigmaZeroEventExt)
  ;  
  \draw[black!90, thick,dashed] (SigmaZeroEventExt) to [out = 40, in = -135]
  (tInf) to [out = -45, in=140]
  (SigmaZeroCosmoExt);
  \draw[black,thick,dotted] (SigmaZeroCosmoExt) to [out=-145,in=0] (SigmaZeroHalf);
  \draw[black,thick,dashed] (SigmaZeroEventExt) to [out=-35,in=180] (SigmaZeroHalf) to (SigmaZeroCosmoEnd);

  \draw[black!70,shorten >= -10,name path=EventFuture] (tInf) --
  (EventZero) ;  
  \draw[black!70,shorten >= -10,name path=CosmoFuture] (tInf) --
  (CosmoZero) ;



  \draw[dashed] (tInf)    -- (-\s - \exts*\s,-\exts*\s);
  \draw[dashed] (tInf)    -- ( \s + \exts*\s,-\exts*\s);
  
  \draw[purple, ultra thick] (SigmaPlus2EventExt) to [out=-35, in= -145]
  node[pos=0.5,above]{$\Sigma(\tau_1)$}
  (SigmaPlus2CosmoExt)
  to   node[pos=0.5,right]{$\SigmaStar(\tau_1,\tau_2)$} (SigmaPlusCosmoExt) to [out = -145, in=-35]
  node[pos=0.5,below]{$\Sigma(\tau_2)$}
  (SigmaPlusEventExt) --   node[pos=0.5,left]{$\mathcal{A}(\tau_1,\tau_2)$} cycle ;


  \path[name path=EventLowerBound] (-1.6*\s, 0.6*\s)-- (-\s, 0)--
  (tInf);
  \path[name path=CosmoLowerBound] (tInf) -- (\s, 0) -- (1.6*\s, 0.6*\s);

  \node[scale=0.5,fill=white,draw,circle,label=above:$\timelikeInf$]at(tInf){};


\end{tikzpicture}


%% file: Images/External-boundary.tex
\begin{tikzpicture}[scale=1.2,every node/.style={scale=1.0}]


  \def \s{3} 
  \def \exts{0.2} 
  \def \t{0.4}
  \def \Tlen{.5}
  \def \endloc{1.3} 
  \def \tPlus{0.2} 
  \def \tWidth{0.2} 

  \coordinate (tInf) at (0,\s); 
  \coordinate (EventZero) at (-\s,0); 
  \coordinate (CosmoZero) at (\s,0); 
  \coordinate (tNegInf) at (0,-\s);
  \coordinate (SigmaZeroEvent) at (-\s + 0.15*\s, 0.15*\s);
  \coordinate (SigmaZeroCosmo) at (\s - 0.15*\s, 0.15*\s);
  \coordinate (SigmaZeroEventExt) at ($(SigmaZeroEvent) + (-\exts, \exts)$);
  \coordinate (SigmaZeroCosmoExt) at ($(SigmaZeroCosmo) + (\exts, \exts)$);
  \path[name path=SigmaZero,in=-35,out=-145]
   (SigmaZeroCosmoExt) to node[pos=0.5,below]{$\SigmaInit$} node[pos=0.5](SigmaZeroHalf){} (SigmaZeroEventExt) ;
   \coordinate (SigmaZeroCosmoEnd) at ($(SigmaZeroHalf) + (\endloc*\s, 0)$);
   \def \negLeft{0.65} 
   \def \negRight{0.5} 
   \def \negLeftAngle{30} 
   \coordinate (SigmaNegLeftEnd) at ($(SigmaZeroHalf)!\negLeft!(SigmaZeroCosmoEnd)$);

   \path[draw=none, fill=black!10] (SigmaZeroCosmoExt) to [out=-40, in = 135]
   coordinate[pos=\negRight](SigmaNegRightEnd)
   (SigmaZeroCosmoEnd) to (SigmaZeroHalf) --cycle;
   \draw[black!50, thick, dashed] (SigmaZeroCosmoExt) to [out=-40, in = 136]
   (SigmaZeroCosmoEnd);

  \path[name path = MBoundary, draw=none, fill=black!5]
  (SigmaZeroEventExt) to [out = 40, in = -135]
  coordinate[pos=\tPlus](SigmaPlusEventExt)
  coordinate[pos=\tPlus+\tWidth](SigmaPlus2EventExt)
  (tInf) to [out = -45, in=140]
  coordinate[pos=1-\tPlus](SigmaPlusCosmoExt)
  coordinate[pos=1-\tPlus-\tWidth](SigmaPlus2CosmoExt)
  (SigmaZeroCosmoExt) to [out=-145,in=-35] (SigmaZeroEventExt)
  ;  
  \draw[black!90, thick,dashed] (SigmaZeroEventExt) to [out = 40, in = -135]
  (tInf) to [out = -45, in=140]
  (SigmaZeroCosmoExt);
  \draw[black,thick,dotted] (SigmaZeroCosmoExt) to [out=-145,in=0] (SigmaZeroHalf);
  \draw[black,thick,dashed] (SigmaZeroEventExt) to [out=-35,in=180] (SigmaZeroHalf) to (SigmaZeroCosmoEnd);

  \draw[black!80, shorten >= -10,name path=EventFuture] (tInf) --
  (EventZero) ;  
  \draw[black!80, shorten >= -10,name path=CosmoFuture] (tInf) --
  (CosmoZero) ;


  \draw[black!80, dashed] (tInf)    -- (-\s - \exts*\s,-\exts*\s);
  \draw[black!80, dashed] (tInf)    -- ( \s + \exts*\s,-\exts*\s);

     \draw[purple,ultra thick]
   (SigmaNegLeftEnd) to [out=\negLeftAngle, in = -135]
   node[pos=0.8, left]{$\Sigma(\tau)$}
   (SigmaNegRightEnd) to
   node[pos=0.5, right]{$\SigmaStar(\tau_{\CosmologicalHorizon}, \tau)$}
   (SigmaZeroCosmoEnd) --
   node[pos=0.5, below]{$\SigmaInit(\tau)$}
   cycle
   ;

  \path[name path=EventLowerBound] (-1.6*\s, 0.6*\s)-- (-\s, 0)--
  (tInf);
  \path[name path=CosmoLowerBound] (tInf) -- (\s, 0) -- (1.6*\s, 0.6*\s);

  \node[scale=0.5,fill=white,draw,circle,label=above:$\timelikeInf$]at(tInf){};

\end{tikzpicture}
